\documentclass[a4paper,10pt]{article}
\usepackage{hyperref}

\usepackage[utf8]{inputenc}

\usepackage[english ]{babel}
 
\usepackage{xcolor}
\usepackage{bbding}

\usepackage[utf8]{inputenc}

\usepackage{tikz-cd}

\usepackage{amsfonts}
\usepackage{amssymb}
\usepackage{graphicx}
\usepackage{mathtools}
\usepackage{skak} 
\usepackage{foekfont}

\usepackage[T1]{fontenc}

\usepackage{amsmath,mathdots}

\usepackage[utf8]{inputenc}
\usepackage{devanagari}

\usepackage{mathtools}
\usepackage{mathdots}
\usepackage{tikz}

\definecolor{amber(sae/ece)}{rgb}{1.0, 0.49, 0.0}
\newfont{\rsfsten}{rsfs10 scaled 1200}

\usepackage{MnSymbol}
\usepackage{tikz}
\usepackage{amscd}
\usepackage{MnSymbol}
\usepackage{wasysym}

\usepackage{mathrsfs}

\newcommand*{\rom}[1]{\expandafter\@slowromancap\romannumeral #1@}

\usepackage{xcolor}

\usepackage{wrapfig}
\newcommand{\tightunderset}[2]{%
  \mathop{#2}\limits_{\vbox to .3ex{\kern-0.95ex\hbox{$#1$}\vss}}}

\usepackage[utf8]{inputenc}

\usepackage{CJKutf8}

\usepackage{pmboxdraw}

\usepackage{stmaryrd}
\usepackage{url}
\usepackage{harmony}
\usepackage{ifsym}

\usepackage{eufrak}

\usepackage{rotating}

\usepackage{dingbat}

\title {Four Lectures on Scalar Curvature}

\author{Misha Gromov}

\begin{document}

\maketitle

 Unlike  manifolds with controlled sectional and Ricci  curvatures, those with their  {\sl scalar curvatures bounded from below} are not configured in  specific rigid   forms   but display  
 an uncertain  variety of  flexible shapes similar to what one  sees  in geometric topology.
 
 Yet, there are definite limits to this flexibility, where determination of such  limits  crucially depends, at least in the known cases,  on two seemingly unrelated   analytic means:  {\it \color {blue} index theory of Dirac operators} and the {\it \color {blue}geometric measure theory,}\footnote{Spaces of metrics with $Sc\geq \sigma$ on 3-manifolds  are amenable to the global  study with the  {\it Hamilton's Ricci flow},  which also applies, at the present moment only  $C^0$-locally,   in higher dimensions.
Also,    much    topological and geometrical information on 4-manifolds  with $Sc\geq \sigma$, for positive as well as negative $\sigma$,   is obtained, exclusively,  with     the {\it Seiberg-Witten equations.} }

 \vspace {1mm}

  The emergent   picture  of spaces with $Sc.curv\geq 0$,    where     topology  and geometry are intimately  intertwined, is reminiscent of  
 {\it the  symplectic geometry},\footnote{Geometric  invariants associated with the scalar curvature, such as the {\it $K$-area},
are  linked with the  symplectic invariants (see [G(positive)  1996], [Polterovich(rigidity) 1996], [Entov(Hofer metric) 2001], [Savelyev(jumping) 2012[), but this link is still poorly understood.}
    but the former has not reached   yet  the maturity of the latter.  \vspace{1mm}
    
 \hspace {10mm}   {\large \sf\color {red!40!black} The mystery of  the  scalar curvature   remains unsolved.}

 \vspace{2mm}
 
 What follows is an extended account  of my lectures, delivered  during the  Spring 2019 at 
 IHES.
 
  In \S1, we give an outline of     results,  techniques  and problems  in  scalar curvature.
  
 In  \S2,  we  spend   a few  dozen pages   on  background Riemannian geometry, with 
   another dozen in section  \ref{Clifford3}  on Clifford algebras and Dirac operators.

 In \S 3,   we   overview  main topics   in geometry and topology of  manifolds  with their  scalar curvatures bounded from below,  state  theorems, explain the ideas of their  proofs and formulate a variety of problems and conjectures. 
 
In \S\S 4 and 5, we reformulate,  in a more precise and general form, what was  stated  in the earlier sections 
and  expose   technical aspects of    the proofs.

 In \S6, we  describe  connective links between different facets of the  scalar curvature 
presented  in  the earlier sections with an emphasis on open  problems.

Finally,   in \S7  we  overview      metric invariants  that are influenced by and/or going along  with the  scalar curvature.
   \vspace {1mm}
 
 I have made a maximal  effort to  lighten the  burden on the reader of  locating the place 
 where a  certain  notation or definition  was introduced.

  Our terminology is displayed in the table of contents. 
  
  When     returning to  the same topic   -- this happens again and again  --
    we, besides    recalling   definitions and  formulas,    explain   what is needed   for the matter at hand,  rather than referring to  earlier   sections  in the  text.
      Everything  needed  for  understanding    a statement   on page "x" can be    found  on  a couple of  preceding pages.

  \tableofcontents

\section {Preliminaries}\label {1}

%%%%%%%%%%%%%%%%%%%%%%%%%%%%%%%%%%%%%%

\subsection {\color{blue} Geometrically Deceptive  Definition}\label {deceptive1}

  %%%%%%%%%%%%%%%%%%%%%%%%%%%%%%%%%%%%%% 
 The scalar curvature of a $C^2$-smooth  Riemannian manifold $X=(X,g)$, denoted 
 $$Sc=Sc(X,x)=Sc(X,g)=Sc(g)=Sc_g(x)$$
  is a continuous  function on $X$, which is  traditionally defined as

 {\sl  the sum of the values  of the   sectional curvatures  at the $n(n-1)$  ordered bivectors  of  an orthonormal frame in $X$, 
{\color {blue!40!black} $$Sc(X,x)=Sc(X)(x)=\sum_{i, j} \kappa_{ij}(x),\mbox {  }  i\neq j=1,...,n,$$}}\vspace {1mm}
 where this sum doesn't depend on the choice of this frame by the Pythagorean theorem.

{\it Algebraically},  this formula defines  a {\it second order differential }   $$g\mapsto Sc(g)$$  
from the space $G_+$ of positive
 definite quadratic differential forms on $X$ to the space  $S$ of functions on $X$, that is 
 characterised uniquely, up to a scalar multiple,  by  two properties.
 
{\sf {\Large $ \star $} the   $g\mapsto Sc(g)$ is {\it equivariant} under  the natural actions of   diffeomorphisms of $X$ in the spaces    $G_+$ and $S$.

{\sf  \Large $ \star $} the   $g\mapsto Sc(g)$ is {\it linear in the second derivatives} of $g$.} \vspace {1mm}

To make geometric sense of this, let us summarize basic    properties of $Sc(X)$.

 $\bullet_1$  {\it  Additivity under Cartesian-Riemannian Products.} 
 $$Sc(X_1\times X_2, g_1\plus g_2)=Sc(X_1,g_1)+ Sc(X_2, g_2).$$

   $\bullet_2$  {\it Quadratic Scaling}. 
   $$Sc(\lambda\cdot X)  =\lambda^{-2} Sc(X), \mbox { for all  } \lambda >0,$$
  where 
  $$\lambda\cdot X=\lambda\cdot (X, dist_X) =_{def} (X, dist_{\lambda \cdot X})\mbox { for }  dist_{\lambda \cdot X}= \lambda\cdot  dist(X)$$
   for all metric spaces $X=(X, dist_X)$  and where $dist\mapsto   \lambda\cdot  dist(X)$ corresponds to $g\mapsto   \lambda^2\cdot g$ for the Riemannian quadratic form  $g$.  
  
{\it Example.} The Euclidean spaces  are scalar-flat, $Sc(\mathbb R^n)=0$, since $\lambda \cdot  \mathbb R^n$ 
is {\it isometric} to $
\mathbb R^n$.

\vspace {1mm}
  
     $\bullet_3$  {\it Volume Comparison.} If the scalar curvatures of $n$-dimensional manifolds $X$ and $X'$ at some points $x\in X$  and  $x'\in X'$ are related by the strict inequality
  $$Sc(X)(x)< Sc(X')(x'),$$
then the Riemannian volumes of  the   $\varepsilon$-balls around these points   satisfy
  $$vol(B_{x}(X, \varepsilon)) > vol(B_{x'}(X', \varepsilon))$$
  for all sufficiently small $\varepsilon>0$.

  Observe that this volume inequality is {\it additive under Riemannian products}:
 if 

 $$\mbox {$vol(B_{x_i}(X, \varepsilon)) > vol(B_{x'_i}(X_i', \varepsilon))$, for $\varepsilon \leq \varepsilon_0$,}$$

and {\it for  all  points}   $x_i\in  X_i$ and $x'_l\in X'_i$,  $i=1,2$,
then 
 $$ vol_n(B_{(x_1, x_2)}(X_1\times X_2, \varepsilon_0)) > vol_n(B_{(x'_1,x_2')}(X'_1\times X'_2, \varepsilon_0)$$
   {\it for  all}   $(x_1,x_2) \in  X_i\times X_2 $ and  $(x'_1,x'_2) \in  X'_1\times X'_2 $.
 
 This follows from the Pythagorean  formula 
 $$dist_{X_1\times X_2} =\sqrt{dist_{X_1}^2+dist_{X_2}^2}.$$ 
 and  the Fubini theorem applied to the  "fibrations" of balls over balls:
 $$B_{(x_1, x_2)}(X_1\times X_2, \varepsilon_0)) \to B_{x_1}(X_1,\varepsilon_0)\mbox  { and } 
B_{(x'_1, x'_2)}(X'_1\times X'_2, \varepsilon_0))  \to B_{x_1}(X'_1,\varepsilon_0),$$
 where the fibers are balls of radii $\varepsilon \in [0,\varepsilon_0]$ in $X_2$ and $X'_2$.
 \vspace {1mm}
 
  $\bullet_4$  {\it Normalisation/Convention  for Surfaces with  Constant Sectional Curvatures.} The unit spheres $S^2(1)$ have constant scalar curvature  $2$ and the hyperbolic plane $H^2(-1)$ with the sectional curvature $-1$ has  scalar curvature $-2$
 \footnote {The equality  $Sc(H^2)=-2$ follows from $Sc(S^2)=2$  by comparing the volumes of  small balls in $S^2\times H^2$ and in $\mathbb R^4$.}
 \vspace {1mm}
 
It is an elementary exercise  to prove the following. \vspace {1mm}

{\Large \color {blue}$\mathbf \star$}$_1$  {\sl The function $Sc(X,g)(x)$ which satisfies $\bullet_1$-$\bullet_4$   exists  and unique};\vspace {1mm}
  
{\Large  \color {blue} $\mathbf \star$}$_2$  {\sl  The unit spheres d  the hyperbolic spaces with $sect.curv=-1$ satisfy
  $$Sc(S^n(1))=n(n-1)\mbox  {  and } Sc(H^n(-1))=-n(n-1).$$}
  Thus, 
  $$Sc(S^n(1)\times  H^n(-1))=0=Sc(\mathbb R^{2n}),$$
 which implies  that 
 
 {\sf the volumes of the  small $\varepsilon$-balls in $S^n(1)\times \mathbf H^n(-1)$ are "very close" to the
 
  volumes 
  of  the  $\varepsilon$-balls  in the Euclidean space $\mathbb R^{2n}$.}
 \vspace {1mm}

Also it is  elementary   to show that the  definition of the scalar curvature via volumes of balls  agrees with the traditional  $Sc= \sum\kappa_{ij}$,
where the definition via volumes  seem  to have an advantage of being   geometrically more  usable.

{\sf \large But this  is an illusion}: \vspace {1mm}

{\color {magenta!90 !blue}\hspace {20mm} {\sf there is no single known (are there unknown?)

\hspace {12mm}geometric 
argument, which would make use  of this  definition.}}

The immediate reason for this is {\it the  infinitesimal} nature of the volume comparison property: it {\it doesn't integrate}  to the corresponding property of balls of specified, let them be small,  radii $r\leq \varepsilon>0.$ \footnote{An attractive conjecture to the contrary  appears in  [Guth(volumes of balls-large) 2011], also 
see  [Guth(volumes of balls-width)) 2011].}

\vspace {1mm}

The following {\it alternative, let it be also only infinitesimal, property} of the scalar curvature seems more promising:

\vspace {1mm}
{\Large \color {blue} $\circledast$} {\sf the inequality $Sc(X, x) < Sc(X{'}, x{'})$ is equivalent to the following
relation between the average mean curvatures  of  the (very)  small $\varepsilon$-spheres $S_x^{n-1}(\varepsilon)\subset X$  and  
$S_{x{'}}^{n-1}(\varepsilon)\subset X{'}$:
$$ \frac {\int_{S_x^{n-1}(\varepsilon)} mean.curv(S_x^{n-1}(\varepsilon),s)ds}{vol_{n-1}(S_x^{n-1}(\varepsilon))} > \frac {\int_{S_{x{'}}^{n-1}(\varepsilon)} mean.curv(S_{x{'}}^{n-1}(\varepsilon),s{'})ds{'}}{vol_{n-1}(S_{x{'}}^{n-1}(\varepsilon))}.$$}

There are  also  several  {\it non-local inequalities} for the mean curvatures of manifolds $B$  with
boundaries $S$, in terms of  the scalar curvatures of $B$ (and sometimes of  sizes of $B$) that we shall see in these lectures,  e.g.
  {\color {blue} \EllipseSolid} \hspace {1mm} and    {\large$ \color {blue}\blacksquare$}  in section \ref {Sc-criteria3},  
but  we are still far from  the  ultimate  inequality of this kind.\vspace{1mm}

%%%%%%%%%%%%%%%%%%%%%%%%%%%%%%%%%%%%%%
 [{\Large $\ast$}]  {\it Exercise: Spherical  Suspension.} Compute the scalar curvature   of the {\it spherical join} of two Riemannian manifolds $X_1$ and $X_2$, that is  the unit sphere in the product of the Euclidean cones over these manifolds:
$$X_1\ast X_2\subset {\sf C} X_1\times {\sf   C}X_2,$$
where ${ \sf C}X=(X\times \mathbb R^\times _+, r^2dx^2+dr^2)$, accordingly 
$${\sf C} X_1\times {\sf   C}X_2=(X_1\times X_2\times  \mathbb R_+^\times \times \mathbb R^\times _+,
r_1^2dx_1^2+r_2^2dx_2^2 +  dr_1^2+dr_2^2)$$
and where the hypersurface  $X_1\ast X_2\subset {\sf C} X_1\times {\sf   C}X_2$ is defined by the equation 
$$r_1^2+r_2^2=1. $$
(The manifold $X_1\ast X_2$ with this metric, which is  defined  for $r_1,r_2>0$, is incomplete; if completed, it becomes singular, unless
 $X_1$ and $X_2$ are isometric to the  unit spheres $S^{n_1}$ and $S^{n_2}$.) 

Show,  in particular, that if $Sc(X_i)\geq n_i(n_i-1)=Sc(S^{n_i})$,  $n_i=dim(X_i)$, $i=1,2$, then 
$$Sc(X_1\ast X_2)\geq (n_1+n_2)(n_1+n_2-1).$$

{\it Hint.}  Use the  formula  for the curvature  of warped products from section \ref {warped2}.
%%%%%%%%%%%%%%%%%%%%%%%%%%%%%%%%%%%%%%

\subsection{\color{blue}Fundamental Examples of Manifolds with $Sc\geq 0$}\label {SYS1}

%%%%%%%%%%%%%%%%%%%%%%%%%%%%%%%%%%%%%%

{\it   \color {blue!50!black} Symmetric and homogeneous spaces.}  Since compact
 symmetric spaces $X$ have non-negative sectional curvatures $\kappa$, they   satisfy $Sc(X)\geq 0,  $
where the equality holds only for flat tori.  

Since the bi-variant metrics on Lie  groups have $\kappa\geq 0$  and since the  inequality $\kappa\geq 0$ is preserved under dividing spaces by isometry groups,   
 all compact homogeneous spaces $G/H$ carry such  metrics, .\footnote{This is also true for   non-compact homogeneous spaces the isometry groups of which contain compact semisimple factors.} \vspace {1mm}

 Furthermore,
 
 {\it   quotients of compact homogeneous spaces by compact  freely acting isometry groups  carry metrics with $Sc\geq 0$,}
 
 where prominent examples of these are

 \hspace{19mm} {\it spheres divided by finite free isometry groups.}

 Thus,   in particular,  
 
 {\it all homology classes in the classifying spaces ${\sf B}(G)$   of finite cyclic groups $G$ are representable by  compact  manifolds
with $Sc>0$  mapped to these spaces.}
 
 But, at the present moment,  it  is  {\color {red!40!black}unknown} if this remains true for all finite groups $G$.\footnote {This was pointed out to me by Bernhard Hanke.}
  
 On the other  extreme, {\sf \color {red!40!black}  there are no known  examples} of ''$Sc>0$  {\it representable}''  {\it non-torsion} homology classes in the classifying spaces of   {\it infinite  countable} groups or of (possibly torsion) homology classes  in the classifying spaces of   groups {\it without torsion.}
 \vspace {1mm}

 (We shall see in the following sections  that  majority of  known topological obstructions to   metrics with $Sc\geq 0$  come from {\it the  rational homology} and {\it K-theory} of classifying spaces of infinite groups.

 Also we  shall meet   examples -- we call these  {\it \textbf Schoen-\textbf Yau-\textbf Schick -manifolds} --  where non-trivial  obstructions to  $Sc\geq 0$,  which  
 reside in  {\it  the integer}  homology  classes in  ${\sf B} (\mathbb Z^n\times  \mathbb Z/p  \mathbb Z)$,
 {\it vanish for non-zero multiples} of these classes.)

 \vspace {1mm}

 {\it \color {blue!50!black} Fibrations.} Since the scalar curvature is additive,  
  {\color  {blue!50!black} fibered spaces $X\to Y$  with compact non-flat  homogeneous fiberes carry metrics  
  with $Sc>0$.}
 
 (This is seen by scaling  metrics  in  $Y$ by  large constants.)

 \vspace {1mm}

{\it    \color {blue!50!black}Convex Hypersurfaces.} Since  convex hypersurfaces in $\mathbb R^n$ as well  as    in general spaces with sectional curvatures $\kappa\geq 0$,  their  scalar curvatures are also non-negative.

\vspace {1mm}

 {\it \color {blue!50!black}Fano, Uniruled  and Calabi-Yau Manifolds}. Smooth  Fano varieties\footnote{A smooth  algebraic  variety $X$ is {\it Fano} if the {\it anticanonical line  bundle}  $L$,
 that is the top exterior power of the tangent bundle, $L= \wedge^nT(X)$, $n=dim X$, is {\it ample}, that is the subset  $Z_x$ of sections in the space $S$ of all sections of   some  power $L^{\otimes m}$ that vanish at $x\in X$ has codimension $n$ for all $x\in X$ and the resulting map $x \mapsto Z_x$ from $X$ to the 
  the space of  codimension $n$ subspaces  in the space $S$  is  a  {\it smooth embedding}.}
   e.g. complex projective  hypersurfaces $X\subset  \mathbb CP^n$ 
 of degree $\leq n$ 
 admit 
 K\"ahler  metrics $g$ with $Sc>0$. 
 
 In fact,  by Yau's solution of the {\it Calabi conjecture},  Fano varieties   carry  K\"ahler metrics with  {\it positive Ricci}  curvatures, while hypersurfaces of 
  if degree n+1 carry   {\it Calabi-Yau} K\"ahler metrics, i.e. with {\it zero Ricci} curvature. 
 
A distinctive  geometric feature of  Fano varieties  is that they are {\it uniruled},  i.e.    covered by rational curves\footnote{The proof  of this relies   on Mori's argument of reduction of  the general  
 case to that of varieties over   {\it finite fields}.}  and  it is {\it conjectured}  that, in general,   
\vspace {1mm}

{\it Uniruled Varieties admit K\"ahler   metrics with $Sc>0$}.\footnote {See    [Debarre(lectures] 2003), [Ballmann(lectures) 2006], [Yang-complex(2017)] and references therein.} 
 
\vspace {1mm}

Conversely, one   knows that

\hspace {-0mm} ({\Large $\star$}) {\sl compact K\"ahler manifolds with $Sc>0$ are
 uniruled.}
 
 \vspace {1mm}

 In fact, this is proven in [Heier-Wong(uniruled) 2012] under the weaker assumption  of positivity  of the   {\it  integral} of the scalar curvature, where, observe, this integral depends only on the first Chern class of $X$ and the   cohomology class of 
  (the symplectic part $\omega$  of) the K\"ahler  metric: 
 $\int_XSc(X,x)dx   =  4\pi/(n-1)! (c_1\smile [\omega^{n-1}](X))$, for $n=dim_{\mathbb C}(X).$ 
 
 There is also a non-trivial geometric constraint on  $\int_XSc(X,x)dx$
 for general  compact  Riemannian manifolds $X$:  \vspace {1mm}
 
 {\sl this integral
 can be bounded from above in terms of  dimension  $dim (X)$, 
diameter, and a lower bound on the sectional curvature of $X$}, see [Petrunin(upper bound) 2008)].
 
 \vspace{1mm}

Yet, it is unclear if there is a {\it true  Riemannian counterpart  of   ({\Large $\star$}):}

{\sf  the literal topological  translation      of  ({\Large $\star$}) may be deceptive:    the connected sum  $X=S^2\times S^2\#S^2\times S^2$, which, as we explain below,   carries  metrics with $Sc>0$,  
admits, however,  no map of non-zero degree from the total space of  any 2-sphere bundle over a surface. }

But if one allows 

\hspace{-3mm} {\sf multiparametric families of maps  $S^2\to X$ and/or suitably controlled 

\hspace{-3mm} discontinuities/singularities}, 

\hspace {-6mm}then this $X$, and apparently all known manifolds $X$ which admits metrics with $Sc>0$, start looking "topologically unirational".

%%%%%%%%%%%%%%%%%%%%%%%%%%%%%%%%%%%%%%
 
\subsection{\color {blue} Thin   Surgery with  $Sc> \sigma$} \label {thin1}

%%%%%%%%%%%%%%%%%%%%%%%%%%%%%%%%%%%%%%

{\it \textbf {Assumptions.}}  {\sf Let an $n$-dimensional manifolds  $X$ bounds a Riemannian manifold $X_+$, i.e. $X_+$ is an  $(n+1)$-manifold with boundary $\partial X_+=X$
and let $Z\subset X_+$ be a  submanifold,  which  meets  $X$ transversally along its boundary  denoted $Y=\partial Z =Z\cap  X=\partial X_+$.\footnote{Here   "boundary"  and  the equality  $Y=\partial Z$ mean  that {\it $Z$ is a manifold-with-boundary}, where this boundary is equal to the intersection of $Z$ with $X$; this is different from  the boundary of  $Z$ as a subset
   in $X_+$, which is in our examples,  where $codim(Z)>0$, coincides with all of $Z$.}}

If $U\subset X_+$  is a tubular neighbourhood of $Z$, then the boundary   $X'=\partial (X\cup U)\subset X_+$ is a smooth (never mind the corner along $\partial U \cap X$) manifold   that implements     surgery of $X$ along $Y=\partial Z\subset X$.\vspace {1mm}

{\it Connected Sum Example.} If $X$ consists of two connected  components,
$X=X_1\sqcup X_2$ and $Z$ is a smooth segments with the ends $y_1\in X_1$ and $y_2\in X_2$, then $X'$ is topological connected sum of $X_1$ and $X_2$  that is performed   by the "tube" $T=\partial U\subset X_+$ joining $X_1$ and $X_2$.

 Observe  that the connected sums and all surgeries performed over  $X$ in general can be realized 
 as above by embedding $X$ as a boundary in a larger (non-compact)   manifold $X_+$,  where one may assume, if one wishes so,  that $X_+$ {\it metrically  splits} near  the boundary $\partial X_+=X$, i.e. 
 it is  isometric to $X\times \mathbb  R_+$ near $X$  and that $Z\subset  X_+$ {\it agrees with this splitting}
 by being equal to  $Y\times \mathbb R_+$ near  $X$.

\vspace{1mm}

If  $Z$ is compact and $\delta>0$ is small, then   the $\delta$-neighbourhood $U_\delta(Z)\subset X_+$ can be taken for $U$. It is also clear that  if the codimension  of $Z$ in $X_+$  satisfies  $k=codim (Y)\geq 3$, e.g. 
if $Y$ is a  curve in  a Riemannian  4-manifold, 
then

{\sl $T_\delta$ with the Riemannian metric induced from 
$X_+$
 has  {\it large positive}  scalar curvature  $\delta$-away  from $X$.} Namely,   by   {\it Gauss' Theorema Egregium},
$$\mbox {$Sc(T_\delta)\sim \frac {(k-1)(k-2)} {\delta^2}$ for small $\delta\to 0$}.$$.

What is more interesting is  that the  submanifold   $X'_\delta=\partial (X\cup U_\delta(Y))$ can be smoothed  by  slightly perturbing it in the $\varepsilon$-neighbourhood of $Y=X\cap T_\delta$, 
for   $\varepsilon=\varepsilon(\delta)\to 0$ for $\delta\to 0$,
such that 

{\sl the scalar curvature of the resulting submanifold, call it $X'_{\delta, \varepsilon}=\partial (X\cup  T'_\delta)$, where $T'_\delta$ denotes the smoothed $T_\delta$, becomes almost as positive as that of $X$.}

This is achieved by a  local    "staircase" construction,\footnote{See [GL(classification) 1980] and [BaDoSo(sewing Riemannian manifolds) 2018].  This construction also applies to hypersurfaces with $mean.curv> \mu$, see [G(mean) 2019] and it
 extends to families of metrics, see  [Ebert-Williams(infinite loop spaces) 2017]   and references therein.

  Besides there is  a  {\it non-local construction} with a similar effect  on the scalar curvature that  was suggested  by Schoen and Yao in [SY(structure) 1979]. }   that  makes  $U_\delta$   thinner and thinner   as you move away from $X$   in the $\varepsilon$-vicinity of $Y$.   

Here is a  precise  statement. \vspace {1mm}

%%%%%%%%%%%%%%%%%%%%%%%%%%%&&&&&&&&&&&&&
{\color {blue}$\CIRCLE$\textbf{--}$\CIRCLE$}   \textbf {Proposition: Thin Surgery by Controlled Thickening. }  {\sf Let $\sigma(x_+)$,  $x_+\in X_+$, be a continuous function, such that its restriction to $X$ satisfies
$\sigma(x_+)< Sc(X, x_+)$, for all $x_+ \in X$, where the scalar curvature of $X$ is evaluated ted with  the   Riemannian metric on $X$ induced from  $X_+$.

Let $\delta (z)>0$ and $\varepsilon(y) >0$  be   continuous positive function on $Z$ and on  $Y$.}

Then there exists a family of tubular neighbourhoods  $U_{\delta, \varepsilon}(Y), \subset X$ with the following  four properties.

$\bullet_1$ {\sf the boundary $T_{\delta, \varepsilon} =\partial U_{\delta, \varepsilon}(Y)$ meets $X=\partial X_+$ {\it \color {blue!40!black}  tangentially}
(rather than transversally) and such that  the submanifold  
 $$X'_{\delta, \varepsilon} = \partial (X\cup  U_{\delta, \varepsilon})\subset X_+,$$
  which  is, a priori,  $ C^1$-smooth, actually  is  {\it \color {blue!40!black} $C^\infty$}-smooth.\footnote{Unless stated  otherwise,  all  our manifolds,  submanifolds etc.  are assumed 
 {\it  smooth} meaning $C^\infty$-smooth.}

$\bullet_2$ The scalar curvature of $X'_{\delta, \varepsilon}$ with the metric induced from $X_+\supset X'_{\delta, \varepsilon}$ satisfies
$$Sc(X'_{\delta, \varepsilon}, x_+) \geq \sigma(x_+) \mbox  { for all } x_+ \in X'_{\delta, \varepsilon}.$$}

Furthermore, 

{\sf $\bullet_3$  $U_{\delta, \varepsilon}$ is contained in the $\delta$-neighbourhood  $U_\delta(Z)\subset X_+$ of
 $Z\subset X_+$, that is the union of all $\delta(z)$-balls, 
 $$U_{\delta, \varepsilon}\subset \bigcup_{z\in }B_z(\delta(z)).$$
 
 {\sf $\bullet_4$ There exists a positive   continuous  function  $0< \delta'(z) =\delta'_{\delta, \varepsilon}(z)< \delta(z)$, such that  the neighbourhood  $U_{\delta, \varepsilon}$ within distance $ >\varepsilon$ from $Y$  is equal to the $\delta'$-neighbourhood of $Z$, that is
 $$ U_{\delta, \varepsilon}\setminus 
 U_\varepsilon (Y)=  U_{\delta'}(Z).$$}}

\vspace {1mm}

The domain $U_{\delta, \varepsilon } \subset X_+$ admits a  more concrete description if 
the  $X_+$ and $Z$ metrically split near $X$, that is  $(X_+,Z)=(X,Y)\times \mathbb R_+$ near $X$.
Namely, one can take 

{\sl the $\delta''$-neighbourhood of $Z$ for $U_{\delta, \varepsilon}$,  where  
 $\delta''(z)=\delta''_{\delta, \varepsilon}(z)$  is a  smooth function on $Z$.}

The  most transparent case here is where  $Y$ is compact; here  one can take 
$ \delta''(z)=\rho (dist(z, X))$  for a suitable function $\rho(d)  =\rho_{\delta,\varepsilon} (d)$, where the nature of this $\rho$ is well represented   by the following.\vspace {1mm}

{\it Halfspace Example/Exercise.} Let $X_+  = \mathbb R^n\times \mathbb R_+$, where 
$X=\partial (X_+\times  \mathbb R_+)=  \mathbb R^n\times \{0\}$   and let  $Z$ be the  half line $\{0\}\times \mathbb R_+\subset \mathbb R^n\times \mathbb R_+$.

Find the above mentioned function $\rho$ for this pair $(X_+,Z)$ and then derive the general case from this example.

{\it Hint}. The scalar curvature of the tube can be calculated either with {\it Gauss'  Theorema Egregium} (section \ref {Gauss2}) or with the {Second Main Formula}  \ref {Weyl2}.A.  (A more general statement is formulated in the Cylindrical Extension  Exercise in \ref {warped2}.)

\vspace {1mm}

{\it\color {red!40!black}  Why not   "$\geq $" instead of "$ >$"?}  One  can't replace the strict inequalities $Sc(X)>\sigma$
and $Sc(X''_{\delta, \varepsilon})>\sigma$ by $Sc(X)\geq\sigma$ and  $Sc(X''_{\delta, \varepsilon})\geq\sigma$, not even in the case of $(X_+,Z)=(\mathbb R^n\times \mathbb R_+, \mathbb R_+)$.

In fact 
the flat metric on $\mathbb R^n$ minus a ball  $B\subset \mathbb R^n$ admits no extension to a complete metric with $Sc\geq 0$ as it 
follows from the solution of the positive mass conjecture (section \ref{asymptotic3}) and/or from non-existence of a  complete metric with $Sc\geq 0$ on the punctured torus  (sections \ref {obstructions4},  \ref {obstructions5}).

\vspace {1mm}

{\it  Exercise: Extension of Families of Metrics with $Sc\geq \sigma$.}  Let  $X$ be a smooth manifold, $\sigma(x)$  a continuous function on $X$ 
and let $\mathcal S_{>\sigma}=\mathcal) S_{>\sigma}(X)$ be the space of Riemannian metrics $g$ on $X$ with $Sc(g, x)>\sigma(x).$

Given an open subset $U\subset X$, let  $\mathcal S_{>\sigma}(U)$ denote the space   metrics on $U$ with $Sc(g, x)>\sigma(x)$ and, if $Y\subset X$ is a closed subset, let  $\mathcal S_{>\sigma}(op(Y))$ be the the space of {\it germs of metrics} with $Sc(g, x)>\sigma(x)$ defined in (arbitrarily small)  neighbourhoods $U\supset Y$.

Show that if $codim (Y)\geq 3$, then the natural (restriction) map 
$$\mathcal S_{>\sigma} \to \mathcal S_{>\sigma}(op(Y))$$
 is a {\it Serre fibration}. (Compare with  {\it Chernysh's theorem} as stated in 2.2.3 in  [Ebert-Williams(cobordism category) 2019].)

\vspace {1mm}

{\it Exercise+Question.} generalize the above to the case where $Z$ is a {\it piecewise smooth} polyhedral subset of codimension $\geq 3$ in $  X_+$.

Then try to generalize this to more general closed subsets $Z$. (See [G(mean) 2019]  for discussion on the corresponding problem for hypersurfaces with $mean.curv>\mu$.)

\vspace {1mm}

{\it\color {red!40!black} Discouraging  Remark.} Despite impressive applications of the above {\color {blue}$\CIRCLE$\textbf{--}$\CIRCLE$} and its variations  to the {\it topology} of  manifolds with $Sc>0$, e.g.
\vspace{1mm} 
 
 {\sf \color {blue!50!black} the existence of metrics with positive scalar curvatures on {\it simply connected} 
 
 manifolds 
 of dimension $n\neq 0,1,2,4 \mod 8$,} \vspace{1mm}

\hspace {-6mm} and of spaces of metrics with $Sc>0$, e.g.\vspace{1mm}
 
 {\it \sf \color {blue!50!black}  {\it infiniteness} of the $k$th homotopy groups of the spaces of such metrics on 
 
 the spheres $S^{4m-k-1 }$ for $m>>k$,}
 
 \hspace {-6mm} the actual geometry behind  "thin construction(s)"   
is skin-deep: positivity of the scalar curvatures of the $n$-spheres for $n\geq 2$ and
nothing else.

In fact, besides homogeneous spaces,  the only  known general  source  of "thickness"   with  $Sc>0$ comes from  solutions of Monge-Ampere equations on  K\"ahler manifolds.
 %%%%%%%%%%%%%%%%%%%%%%%%%%%%%%%%%%%%%%

\subsection {\color {blue} Scalar Curvature and Mean Curvature} \label {mean1}

 %%%%%%%%%%%%%%%%%%%%%%%%%%%%%%%%%%%%%%
 A  simple  link between the two notions is provided by the following observation.\footnote {Look at fig  8 in [GL(spin) 1980].}\vspace {1mm}
 
 Let  $X=(X, g) $ be a  Riemannian $n$-manifold  with    boundary  represented by  
 a domain in a slightly larger manifold $X_+\supset X$ and then embedded 
 to the cylinder
$ X_+\times \mathbb R$ for 
$$X= X_0=X\times \{0\}\subset X  \times \mathbb R \subset X_+ \times \mathbb R$$
and let $U_\varepsilon=U_\varepsilon(X_0)\subset  X_+ \times \mathbb R$ be the 
 $\varepsilon$-neighbourhood  of $X_0\subset X_+ \times \mathbb R$.

The boundary $\partial U_\varepsilon$ consists of two parts:
 two 
 "$\varepsilon$-copies" 
 $$X_{\pm\varepsilon}=X\times\{\pm\varepsilon\}\subset  \partial U_\varepsilon$$ 
 of $X$ and the complementary semicircular band, 
 $$\partial X_0\times S^1_+(\varepsilon)\subset \partial X_0\times S^1_+(\varepsilon),$$
  that is 
 one half of the boundary of the $\varepsilon$-neighbourhood of the boundary $\partial X_0\subset  X_+ \times \mathbb R$.

Both parts of the hypersurface  $\partial U_\varepsilon$ are $C^\infty$-smooth,\footnote {Our Riemannian manifolds  are $C^\infty$-smooth unless stated otherwise.} and  $\partial U_\varepsilon$ is also $C^1$-smooth at
the common boundary of these parts. But  the  curvature of  the band  
$\partial X_0\times S^1_+(\varepsilon)\subset \partial U_\varepsilon $  
along  the  semicircles $ \{x\}\times S^1_+(\varepsilon) $, $x\in \partial X_0$,   jumps down  from   $\varepsilon^{-1}$ to  $0$, where this band meats the "flat horizontal" part 
  of $\partial U_\varepsilon$.
  that is  the union of the two  "$\varepsilon$-copies" of $X$,  
 $$X_{-\varepsilon}\cup X_{-\varepsilon}\subset  \partial U_\varepsilon.$$
,

The scalar curvature of this band, computed with the Gauss formula (theorema egregium),  interpolates  between, roughly,   $\varepsilon^{-1} \times  mean.curv(\partial  X_0)$ at the points not too close to the flat part of  $\partial U_\varepsilon$,
where it  becomes equal to the scalar curvature of $X$.  \vspace {1mm}

  {\sf if $Sc (g)> \sigma$ and the boundary    $\partial X\subset X$ is  {\it strictly  mean convex}, i.e. 
 
 \hspace {-6mm}$mean.curv(Y)> 0$,\footnote{Our coorientation convention is such that convex domains are mean convex according to it.} 
 then the boundary $\partial U_\varepsilon$ can be  $C^\infty$-smoothed by  interpolating the curvatures on the two sides of the jump between  $\varepsilon^{-1}$ and  $0$,
    such that}\vspace {1mm} 
    
  {\color {blue!20!black}  {\it  the scalar curvature of the smoothed boundary becomes bounded from below
    
     by the scalar curvature of the original metric  $g$ on $X$.}}\vspace {1mm}
 
(To see this, look at the $(n-2)$-ball in the $n$-space, $X_0=B^{n-2}\subset \mathbb R^n$, where the boundary of its $\varepsilon$-neighbourhood can be $O(n-2)$-invariantly smoothed by $C^\infty$-flattening the  semicircle  
$S^1_+(\varepsilon)$ at the ends, while keeping it convex.)\vspace {1mm}

Since the  boundary $\partial U_\varepsilon $ is naturally  diffeomorphic to the  
  {\it double}  \DD$(X)$    obtained by gluing two copies of $X$ along $\partial X$,
this delivers the following \vspace {1mm}

%%%%%%%%%%%%%%%%%%%%%&&&&&&&&&&&&&&&&&
\textbf {Proposition: Smoothing \DD-Corner}. { \it  \color {blue!40!black} There  exists an approximation of the natural continuous  metric $G_0$ on   the double \DD$(X)=X\cup_{\partial X} X$ by smooth 
metrics $G_\varepsilon$ with 
scalar curvatures bounded from below by $Sc(G_\varepsilon) \geq Sc(X).$}

Moreover,   

{\sf  \color {blue!40!black} {\it strictness}  of positivity of  the mean curvature, can be propagated\footnote {See section 11.2   in [G(inequalities) 2018].}  by a small $C^\infty$-perturbation to such a "strictness" for  the scalar curvature  all over  \DD$(X)$,} thus making

\hspace {15mm} {\it $Sc( G_\varepsilon)$ everywhere strictly greater than $Sc(X)$.}

\vspace {1mm}

 For instance, \vspace {1mm}\vspace {1mm}

{\sf   \color {blue!40!black}  the   doubles of  compact mean convex  bounded Euclidean domains carry 

metrics 
with 
positive scalar curvatures,} \vspace {1mm}

where  the necessary  strictness of mean convexity is  achieved by small perturbations of the boundaries of these domains.\vspace {1mm}

If  you think about this (excessively geometric) construction  in intrinsic terms of $X$, you  will   realize that  the metric $G_\varepsilon$  was actually obtained by stretching the original Riemannian metric  $g$ of $X$  near the boundary $\partial X\subset X$ {\it along geodesic segments  normal to  $\partial X$}. Then  you write   down  everything in the  normal coordinates in a neighbourhood of the boundary $\partial X\subset X$ \footnote {[Almeida(minimal) 1985], [Miao(corners) 2002], [Bre-Mar-Nev(hemisphere) 2011], [G(billiards) 2014].}
  and  arrive at   the following proposition.
   \vspace {1mm}
 
 %%%%%%%%%%%%%%%%%%%%%%%%%&&&&&&&&&&&
   {\it \textbf {Miao's Gluing Lemma.}}  {\sf Let $X_\circlearrowleft$ be   obtained by 
   identifying   pairs of points in the boundary   of a Riemannian manifold $X=(X,g)$  by an isometric involution $I:\partial (X)\to \partial X$ without fixed points.}\footnote  {This $I$ may  be more interesting than interchanging  two isometric components of the boundary, such as the   involution  on the boundary of a centrally symmetric  $X\subset  \mathbb R^n$.}

 {\it If  the sums of the mean curvatures at the identified   points satisfy
$$\mbox { $mean.curv(\partial X, x)+ mean.curv(\partial X,  I(x)))>0$  for all $x\in \partial X$},$$
 then the natural continuous Riemannian  metric $G$ on $X_\circlearrowleft$ can be approximated by smooth metrics $G_\varepsilon$  with their   scalar curvatures strictly bounded from below by the   scalar curvature of $g$.}\footnote {This is similar to   preservation of  lower bounds on (Alexandrov's) {\it sectional curvature}  under   gluing,   where  the second fundamental form II of the boundary satisfies II$_x$+II$_{I(x)}\geq 0$.}
 \vspace {1mm}

{\it The main step  of the poof} is  stretching  $g$  normally to $\partial X$   in a small neighbourhood of  $\partial X$ with no decrease of the scalar curvature and without changing the restriction $g|_{{\partial X}}$, such that  the second fundamental form $A$  for the new metric $g_{new}$ on $X$ will match one another at the $I$-corresponding points, 
i.e.    
$$A_x+ A_{I(x)}=0\mbox { for all $x\in \partial X$}.$$ 

We   implement  such a stretching by extending $X$ with   the $\varepsilon$-cylinder $\partial X\times [0, \varepsilon]$ attached to $X$ by the tautological map $\partial X\times \{0\}\to \partial X$ and 
we endow  this cylinder  with a family of metrics $g_\varepsilon$  defined with the following metrics 
$h_\varepsilon(t)$ on  $\partial X$.

 Let  $A_{old},  A_{new}$ be  quadratic differential forms on $\partial X$, where $A_{old}$ is equal to  the   second fundamental form of $\partial_0=\partial X\times \{0\}=\partial X$ in $X$  and $A_{new}$
is another (desired) quadratic differential form on $\partial X=\partial_\varepsilon=
\partial  X\times \{\varepsilon\}$.

Let $$h_\varepsilon(t)=h + t A_{old} + \frac {t^2}{2\varepsilon}( A_{new}-A_{old}), \mbox { } 0\leq t\leq \varepsilon.\leqno {(++)}$$
and let 
$$g_\varepsilon = h_\varepsilon(t)+dt^2.$$

Then:

(i) {\sl the 
 second fundamental forms of the two boundary parts  $\partial_0$ and $\partial_\varepsilon$ for   the metric $ g_\varepsilon$  are equal to $A_{old }$ and $A_{new}$  correspondingly by the Riemann variation formula   in \ref {equidistant2};} \footnote {These forms are  evaluated on the (same unit) vector field $\frac {d}{dt}$.}
   
(ii)   the scalar curvature of $ g_\varepsilon$ satisfies,
   $$ Sc(g_\varepsilon)=\frac {1}{\varepsilon}  trace(A_{old}-A_{new}) +O(1)$$
by {\it Hermann Weyl's tube formula}  and  Gauss's formula  (see \ref {Weyl2}, \ref {Gauss2}). 

It is also clear that $(g_\varepsilon)_{|\partial_0}=h$ and $(g_\varepsilon)_{|\partial_\varepsilon}=h+o(\varepsilon)$, which allows a small perturbation of $g_\varepsilon$ that makes it equal to $h$
on $\partial_\varepsilon$, while   keeping  (i) and (ii).\footnote{Details can be found in section   11.5 in [G(inequalities 2018].}\vspace {1mm}

{\it Second Step.}  Because of the match of the second quadratic forms,   the  metric  $G_{new} $  on $X_\circlearrowleft$ is now $C^1$-smooth, which allows its  painless smoothing,  while   metric keeping the  scalar curvature almost as 
as positive as that of $g$,  and, due to the strictness condition, even {\it more positive} than $Sc(g)$.\footnote{ This  trivially follows from a  general "{\it local h-principle}", see  section 11.1 in [G(inequalities 2018] and [Baer-Hanke(local flexibility) 2020]. }
\vspace {2mm}

Besides the above   "infinitesimal realtions",  there is   an amusing   similarity between global geometries of n-dimensional  {\it Riemannian manifolds $X$ with positive scalar curvatures} and  {\it mean convex convex hypersurfaces} 
in the Euclidean space $\mathbb R^n$  and in similar spaces.

Although,  in many respects mean convex  hypersurfaces $Y\subset \mathbb R^n$  are better understood then manifolds $X$ with $Sc(X)>0$,  
    essential geometric properties of    $Y$ with $mean.curv \geq \mu$   can be proved at the present moment only in
    the light of the scalar curvature  by means of twisted  Dirac operators or minimal hypersurfaces  and where the transition from mean curvature to the scalar curvature is most clearly seen in the doubling construction. \footnote{See   [G(mean) 2019], [Lott(boundary) 2020], [Cecchini-Zeidler(scalar\&mean) 2021  and 
    section \ref {mean convex3}   for  more about it.}

\vspace {1mm}

{\it Exercises}.  Let $X$ be a Riemannian $n$-manifold with a non-empty mean convex boundary.
Show the following.

(a) If  $X$ has non-negative scalar curvature, then the double of $X$ admits a metric with $Sc>0$, unless $X$ is Riemannian  flat with flat boundary. 

For instance, doubles of mean convex domains in $\mathbb R^n$ carry metrics with positive scalar curvatures.

(b) If $X$ has non-negative Ricci curvature then  either it is diffeomorphic to to a regular neighbourhood of a $(n-2)$-dimensional curve-linear polyhedral  subset  $P^{n-2}\subset X$,  or  it is Riemannian  flat with flat boundary. 

For instance, if $X$ is  connected orientable of  dimension  $n=3$,  then it is either diffeomorphic to a handle body, or it is 
isometric to a flat  torus times a segment $[-d,d]$, or to a flat bundle over a flat Klein bottle with the fiber $[-d,d]$.

(c) If  $X$ admits an equidimensional  isometric immersion to a complete simply connected manifold $\hat X^n$ with non-positive sectional  curvature, then   it is also diffeomorphic to to a regular neighbourhood of  an  $P^{n-2}\subset X$.

Moreover, if $\hat X^n$ is equal to  the hyperbolic space $\mathbf H^n$ with the sectional curvature $-1$, then 
the condition $mean.curv(\partial X)\geq 0$ can be relaxed to  $mean.curv(\partial X)\geq -(n-1)$ and if  $\hat X^n=\mathbb R^n$ then one  needs only the following integral  bound on the negative part $M_-$  of the mean curvature
of $Y=\partial X$, 
$$\int_Y |M_-(y)|^{n-1}dy\leq (n-1)^{(n-1)}\gamma_{n-1},$$
where 
$$M_-(y)=\min(0, mean.curv (Y, y))$$
and
$\gamma_{n-1}$  denotes the volume of the unit sphere $S^{n-1}$. 

{\it Remark/Question.}  The above integral inequality is sharp, where the equality holds for bands between concentric spheres. 

But it is unclear what is the sharp inequality for domains $X\subset \mathbb R^n$ with {\it connected} boundaries $Y$.

For instance 

{\sf what is the infimum of $\int_Y |M_-(y)|^2dy$ for {\it torical} $Y=\partial X\subset \mathbb R^3$, where $X$ is {\it not diffeomorphic to the solid torus}}?

{\sf Is there a lower bound on  $\int_Y |M_-(y)|^2dy$  by the topology of $X$, e.g. by positive  $const_n$  times  the {\it simplicial volume} of $X$?}
(Compare with the {\it simplicial volume conjecture} in section \ref {negative3}.)

%%%%%%%%%%%%%%%%%%%%%%%%%%%

\subsection {\color {blue} Topological and Geometric Domination by Compact  and non-Compact Manifolds with  positive Scalar Curvatures}\label {domination1}

%%%%%%%%%%%%%%%%%%%%%%%%%%
The global effect of positivity of the curvature of a   Riemannian manifold $X$ is  a bound 
on the overall size of $X$. 
This, in the case  the sectional and Ricci curvatures, can be expressed in terms of  simple    geometric  characteristics of $X$, e.g.   the diameter and the volume, which are defined in purely metric  terms with no direct reference to the topology of $X$.

Positivity of scalar curvature  also limits the size of $X$, geometrically as  well as topologically, but here  the  bounds on geometry in terms of $\inf Sc(X)$ can't be even {\it properly}  formulated  
without explicit use of the  underlying topology of $X$. \vspace {1mm}

 {\color {teal} \textbf {1. \it Prelude to Example}.} {\sf Let $g$ be a Riemannian metric on the Euclidean space  $\mathbb R^n$ with uniformly positive scalar curvature, i.e. $Sc(g)\geq \sigma>0$.} Then  
 
 {\it this $g$ can't be greater than the Euclidean metric in two respects.} \vspace {1mm}
   
   {\color {teal} (a) }  {\it For all $D>0$, there exist  points $y_1,y_2\in \mathbb R^n$ with $dist_{Eucl}\geq D$,
     
      such that 
      $$  dist_g(y_1, y_2)\leq const=const_{n,\sigma}, \mbox {  to be specific, say, for } const= \frac {2\pi\sqrt {n(n-1)}}{n\sqrt \sigma}. $$}
(Recall that $n(n-1)$ is the scalar curvature of the unit sphere $S^n.$)\vspace {1mm}

 {\color {teal}(b)} {\it For all $\varepsilon>0$, there  exists a smooth surface $S\subset \mathbb R^n$, such  that
$$area_g(S)\leq \varepsilon\cdot area_{Eucl}(S).$$}\vspace {1mm}

On the  surface of things, there is nothing particularly topological about these (a) and  (b),  but the {\it true 
comparison relations} between metrics  with $Sc>0$  and the Euclidean ones,  which are expressed by means of, in general {\it non-diffeomorphic,  maps} from $X$ to $\mathbb R^n$   are inherently
 topological. \vspace{1mm}

{\color {blue} \textbf {2.  Example}: {\it Euclidean non-Domination  with  $Sc(X)\geq \sigma>0$.} } {\sf Let $X$ be  an orientable  Riemannian manifold of 
 dimension $n$ with uniformly positive scalar curvature, $Sc(X)\geq \sigma>0$,
 and let $f:X\to \mathbb R^n$ be   a smooth proper  map \footnote{A  map is {\it proper} if "infinity goes to  to infinity". Formally: 
  {\it the pullbacks  of compact subsets are compact.}} with {\it non-zero} degree.\footnote{A sufficient geometric condition for this "{\it non-zero}" reads: {\sf there is a  non empty open subset $ U\subset \mathbb R^n$,  such that the pullbacks $f^{-1}(u)\subset X$, $u\in U$,  are finite and contain {\it odd}  numbers of points.}} }Then,
  
  {\it this $f$ can't be uniformly Lipschitz on the large scale, nor can it be uniformly 
  
  area non-expanding.} 
  
  This means   the following.\vspace {1mm}

{\color {blue}\textbf{(a*)}}   {\it For all $D>0$, there exist  points $y_1,y_2\in \mathbb R^n$ with $dist_{Eucl}\geq D$,
     
      such that the distance between their pullbacks $f^{-1}(y_1)\subset X$ and  $f^{-1}(y_2)\subset X$
      is uniformly  bounded,
      $$  dist_X(f^{-1}(x_1), f^{-1}(x_2))<const. \footnote {If $n=dim(X)\leq 9$, a 
       Schoen-Yau kind of argument with minimal  hypersurfaces reduces the problem to  an auxiliary {\it spin} manifold 
        with $Sc \geq \sigma$ to which the Dirac theoretic argument applies, see  section 5.3 in [G(billiards) 2014]. But if  $X$ is
       {\it non-spin}  of dimension $n\geq 10$, I can’t vouch for the proof, since it 
 depends  on "desingularization" of minimal varieties  from the papers
      [SY(singularities) 2017] and/or [Lohkamp(smoothing) 2018],  which I have not studied in depth.}
      $$}

{\color {blue}\textbf{(b*)}}  {\it If $X$ is spin,\footnote {This a somewhat tricky topological condition, which we shall explain later on. It suffices to say at this point that manifolds homeomorphic  to $\mathbb R^n$, and  more generally, those with vanishing cohomology group $H^2(X;\mathbb Z_2)$ are spin. 

But, for instance, the connected sum of $\mathbb R^4$ with the complex projective plane $\mathbb CP^2$ is non-spin.} then for all $\varepsilon>0$, there  exist smooth surfaces $S\subset  X$,  and $\underline S\subset \mathbb R^n$, such  that 
$$area_X(S)\leq \varepsilon\cdot area_{\mathbb R^n}(\underline S)$$
 and such that 
the map $f$ sends $S$ diffeomorphically onto $\underline  S$.}
\vspace {1mm}

{\it Remarks  and Corollaries.} (i) The above  {\color {blue} \textbf {1}} follows from    {\color {blue}  \textbf {2}} applied to the identity map $id: (\mathbb R^n, g)\to  (\mathbb R^n, g_{Eucl})$.\vspace {1mm}

(ii) It follows from  {\color {blue}  \textbf {2}} that   

{\it no compact orientable $n$-manifold $X$ with $Sc(X)>0$ admits a map $f$ with non-zero degree to 
the $n$-torus $\mathbb T^n$, \footnote {This was proved in  [SY(structure 1979] for $n\leq 7$ and in [GL(spin) 1980] for spin manifolds $X$ and all $n$. Nowadays,  (yet unpublished in an academic journal) analysis of singularities of minimal hypersurfaces in dimensions $n\geq 8$ in  [Lohkamp(smoothing) 2018]  and in [SY(singularities) 2017] yields this result for all, not necessarily spin, compact  manifolds  $X$ and all $n$.} while}  {\color {blue}\textbf {1}} {\it yields this (only) for diffeomorphisms $X\to \mathbb T^n$.}

\vspace {1mm}

{\it Proof.} Given a  smooth map $f: X\to \mathbb T^n$, let  $\tilde f: \tilde X\to\mathbb R^n$ be its lift to the $\mathbb Z^n$-coverings   of both manifolds and apply either {\color {blue}\textbf{(a*)}}  or {\color {blue}\textbf{(b*)}} 
to the $\varepsilon$-scaled map $\varepsilon \tilde f: \tilde X\to\mathbb R^n$  for $\varepsilon\to 0$.

\vspace {1mm}

(iii) The proof of {\color {blue}\textbf{(a*)}},  mainly  depends  the {\it geometric measure theory}, (see  sections \ref {SY1}, \ref  {FCS1}
while  {\color {blue}\textbf{(b*)}} relies on an  {\it index theorem for "twisted"  Dirac  
s} (see sections \ref {Dirac1}, \ref {twisted1}).  At the present day,  there  is no alternative proof 
of {\color {blue}\textbf{(b*)}} (not even of  {\color {teal}\textbf{(b)}})  and   {\color {blue}\textbf{(b*)}} remains unknown for general  (non-spin) manifolds $X$of dimension $n\geq 4$.\footnote {All 3-manifolds are spin.}
 \vspace {1mm}

 Motivated by the above example we make the following definition.\vspace {1mm}

 {\it \textbf { Domination by $\mathbf {Sc>0}$.}}  Let  $\underline X$ be a "nice", say locally contractible topological space, e.g. a  cellular or    polyhedral one,  and let 
  $\underline h\in H_n(\underline X)$ be a homology class.
  Say that a, possibly open, oriented connected $n$-manifold $X$ {\it dominates $\underline h$}, if there
   exists a continuous map $f: X\to\underline X$ {\it locally constant at infinity}, called {\it $\underline h$-dominating map}, which   sends the fundamental homology class $[X]$ to $\underline h$.
 \vspace{1mm}

{\large \sf  Quasi-Proper Maps.}  Similarly, if $X$ is a locally compact and countably compact, one defines domination of homology classes  with {\it infinite supports}  in $\underline X$, where the relevant maps $f:X \to \underline X$ are {\it quasi-proper},   i.e. they extends to  continuous maps between the compactified spaces, 
from $ X^{+ends}\supset X$ to $ \underline X^{+ends} \supset \underline X$,  obtained by a attaching the sets of ends to these  spaces.

In simple words, $f$ is quasi-proper if 
 
 {\sf for all proper maps   $\phi :  \mathbb R_+\to X$ (i.e. $\phi(t) \to \infty$ for $ t\to \infty$),  the composed 
 map $f\circ\phi:  \mathbb R_+\to  \underline X$  is either proper or converges to a point in $\underline X$  for $t\to \infty$.}
 
 \vspace {1mm}
 
{\large \sf Domination of Manifolds}. For instance, if $\underline X$ is an oriented   $n$-manifold  or a pseudomanifold,\footnote {An {\it $n$-pseudomanifold}  is a triangulated space, where the singular locus, where  this  space  is {\it not} locally $\mathbb R^n$, has codimension (at least) 2.}
  and $\underline h=[\underline X]$,
 then  these "dominations" also called {\it dominations  with degree 1} are just  quasi-proper maps $X\to \underline X$ of degree 1.  \vspace{1mm} 
 
 More generally, {\it \color {blue} domination with $degree\neq 0$}  -- we  shall meet  these  many time in these lectures  ---   refers to equividimensional
 maps of non-zero degrees  between orientable  manifolds or pseudomanifolds.
 
 \vspace{1mm} 
 
 Next, if  $X$ and  $\underline X$ are  a metric spaces, say that $\underline h$ is {\it $\lambda$-Lipschitz dominated} or {\it distance-wise  $\lambda$-dominated   
 by $X$} if the  map $f$ is $\lambda$-Lipschitz, i.e. $dist_{\underline X} (f(x_1), f(x_2) \leq\lambda \cdot dist_X(x_1,x_2)$.
 
 Similarly, define {\it area-wise $\lambda$-domination}, by the inequality 
 $$area_{\underline X}(f(S))\leq \lambda \cdot area_X(S),$$ 
  provided that areas of (suitable) surfaces $S\subset X$ and of their images $f(S)\subset \underline X$  their are suitably  defined in $X$ and $\underline X$,  e.g. where  these are smooth surfaces in Riemannian manifolds. \vspace{1mm}

  \textbf {3. Positive Scalar Curvature Domination Problems.}  {\sl What are spaces $\underline X$ and classes 
   $\underline h\in H_n(\underline X)$, which can and which can't be  dominated by complete  Riemannian manifolds $X$ with $Sc(X)>0$?
  
  How much does  the answer depend on additional conditions on  topology  and geometry of a dominating manifold $X$?
   
   When  can   such a domination  be implemented with $\lambda$-Lipschitz or with  area $\lambda$-contracting maps?}

     \vspace {1mm} 
 
  Notice that   {\color {blue}\textbf{(a*)}}  says  in this regard  that 
  
  {\it for no $\lambda>0$, a  non-zero multiple of the fundamental homology  class $[\mathbb R^n]$  can be  (large scale) distance-wise  $\lambda$-dominated   by a manifold $X$ with $Sc(X)\geq \lambda>0$},
  
  Similarly,   {\color {blue}\textbf{(b*)}}  can be stated as  {\it non-existence of  {\it area-wise  spin}   $\lambda$-domination.}\vspace{1mm} 
  
 %%%%%%%%%% this ``'' 4`` is '' new, added on  Jul  8 !!!  add in the volume version}
   \textbf {4. From Algebraic Topology to Asymptotic Geometry:  Topological}  
   \textbf{versus Lipschitz 
   Domination.} However trivial, it should be emphasized  that the existence of

   {\sf a positive scalar curvature domination  of  a compact  orientable manifold $\underline X$ 
   
  (or a pseudomanifold)  with degree $d$}

 \hspace{-6mm}    implies 
   
   {\sf  positive scalar curvature 1-{\it Lipschitz} domination  of  {\it all covering of  $\underline X$,} 
   
   $\underline X$, in particular, of   the {\it universal 
   covering} $\tilde {\underline X}$, with degrees $d$.} 
    
    \vspace{0.6mm} 
   
  (Continous maps $f: X\to \underline X$ can be approximated by $\lambda$-Lipschitz maps; these lift to the coverings 
   and can be made 1-Lipschitz by scaling $X=(X,g)\mapsto \lambda\cdot X=(X, \lambda^2\cdot g)$.) 
 
Also it must be noted that the Lipschitz domination between  open  manifolds is  more  general and versatile  relation than the   topological domination for  compact manifolds.\footnote {This is well demonstrated by aspherical 4- and 5-dimensional manifolds in  section \ref  {5D.3}.}

For instance, only exceptional  compact  aspherical  $n$-manifolds $X$ dominate the $n$-torus $\mathbb T^n$, but

{\sf \color {red!20!black} there is no single example (so far), where the  universal  covering $\tilde X$  of a  
compact  

aspherical $X$ wouldn't  
1-Lipschitz dominate $\mathbb  R^n=\tilde {\mathbb T}^{n}$.}

 In view of this, the topological   $Sc$>$0$-domination  problem is shifted to a more fruitful geometric one of the {\color {red!50!black}(non)}existence of  \vspace{1mm} 
 
 {\sf
  
  \hspace{0mm}  {\color {blue!60!black} a 1-Lipschitz domination of an open   $\underline X$ by Riemannian  manifolds $X$ with 
  
  $Sc(X)\geq \sigma>0$, or  -- this is most relevant if $\underline X$ is complete  --  by $X$ with $Sc(X)>0$.}}

 \vspace{1mm} 
 
  \textbf {5.  Domination Equivalence  {\color {red!0!black} Conjecture}.} {\sf  If a homology class $\underline h $ in $\underline X$ (here it is an ordinary one, with compact supports)   is dominated by a
  {\it complete}  manifold $X$ with $Sc>0$, then it  also admits    a {\it compact spin domination}, i.e. by  a   {\it compact spin} manifold $X_o$ with 
  $Sc(X_o)>0$.}
  
  (It may be safer to assume $n\neq 4$; also, to avoid irrelevant purely topological obstructions to dominability,  one  should  replace     "{\sf domination of}  $\underline h$" by  
   "{\sf domination of  a {\it non-zero multiple} of $\underline h$"} in some cases.)
  
  \vspace{1mm}
  
 Let us stress  out that  the  most essential  cases of   this conjecture concern homology  classes in {\it aspherical spaces} 
 that are 
  classifying spaces of 
 of discrete 
 groups, 
 $$\underline X={\sf B}(\Pi)=K(\Pi,1)\mbox {   for } \Pi=\pi_1(\underline X),$$
and that  the {\sf \textbf {main} (topological and  naive)  \textbf {$Sc\ngtr 0$-conjecture} -- the scalar curvature counterpart of {\it Novikov's higher  signatures conjecture} --  formulated in the present terms  reads: \vspace{1mm}

{\color {blue} $\mathbf {[Sc \ngtr \mathbf 0]} $}  {\it \large No non-torsion homology class in the classifying space  of a countable group
 can be dominated by a compact manifold with $Sc>0$}.} \vspace{1mm}

\textbf {6.   From $Sc\ngtr 0$ to $Sc \ngeqq 0$.}
As far as the topology of a complete manifold  $X$ is concerned, there is little difference between  the conditions  $Sc(X)\ngtr 0$ and $Sc(X)\ngeqq 0$,  

\hspace {-6mm} where, observe,  the former corresponds  to the bound  $\inf Sc\leq 0$  and the   latter   to  $\inf Sc<0$. 
\vspace {1mm}

Indeed, according to 
{\it  \textbf  {Kazdan's deformation theorem,} }\vspace {1mm}
 
  {\sf  non-existence of a deformation of   metric on a complete Riemannian manifold $X$ with $Sc\geq 0$ to  a complete metric with 
$Sc>0$ implies that $X$ is {\it Ricci flat,}} [Kazdan(complete) 1982].\vspace {1mm}

If $dim(X) =3$ then then  "Ricci flat"  implies {\it Riemannian flat}; and if $n\geq 4$, the   {\it \textbf {Cheeger-Gromoll splitting
 theorem}} shows in most (all?) of  our cases that $X$ is {\it Riemannian flat}, i.e. {\sf isometric  to the Euclidean space divided by a discrete isometry group. }

Thus, as we shall see in several  examples later on, 
 \vspace {1mm}

 \hspace {4mm}{\it non-existence theorems}  for  $Sc>0$ yield {\it rigidity results  for $Sc\geq 0$.}

\vspace{1mm}

{ \it  Spin Domination Problem.} Non-domination result proven with a use of Dirac operators (these are many, require the dominating manifolds to be {\it spin}\footnote {See   section \ref {spin index3} for the definition of spin and recall that manifolds with $w_2=0$, e.g.  {\it stably  parallelizable} ones,  are spin.}

This could be removed in majority of cases if the following were true.

{\it\color {red!50!black} \sf Unrealistic Conjecture}.  {\sf Compact Riemannian  orientable manifolds $\underline X$ with  
positive scalar curvatures can be dominated with $degree\neq 0$ by compact Riemannin  manifolds $X$ with
 $Sc \geq 0$  and with {\it  universal coverings  spin.}\vspace{1mm}}

{\it Exercise.} Prove these conjecture for manifolds $\underline X$ of dimensions $n\geq 5$   with finite fundamental groups.

{\it Hint. }  Use Thom's theorem on domination  of multiples of homology classes by stably  parallelizable manifolds and classification of simply connected manifolds  with $  Sc>0$ of dimension $\geq 5 $
as in \HandRightUp of \ref {spin index3}.

{\it Remark.}  Proving (maybe  disproving?) this conjecture  seems  possible by the present day tecnniqies   for manifolds with Abelian fundamental groups.
%%%%%%%%%%%%%

\subsection {\color {blue}  Analytic Techniques} \label{analytic1}

%%%%%%%%%%%%%%%%%

The logic of most (all?) arguments  concerning the global geometry of  manifolds $X$ with {\it scalar curvatures bounded from below}  is, in general terms, as follows.

Firstly, one uses  (or proves) the {\it existence theorems for solutions $\Phi $ of certain partial differential equations}, where the existence of these $\Phi $ and their properties depend on global,  topological and/or  geometric assumptions $\cal A$ on $X$, which  are,  a priori, unrelated to  the scalar curvature.

Secondly, one concocts  some {\it algebraic-differential 
 expressions}    ${\cal E}(\Phi, Sc(X))$, where the crucial role is played by certain {\it algebraic formulae} and issuing  inequalities satisfied by ${\cal E}(\Phi, Sc(X))$ under  assumptions $\cal A$.
 
 Then one  arrives at  {\it a contradiction}, by showing that\vspace{0.8mm}
 
 \hspace {7mm} {\it if  
 $Sc(X)\geq \sigma$, then the implied  properties, e.g. the sign,  of  ${\cal E}(\Phi, Sc(X))$   are
 
  \hspace {7mm}    opposite to those satisfied under assumption(s) $\cal A$.}\vspace{1mm}
 
 %%%%%%%%%%%%%%%%%%%%%%
\subsubsection{\color {blue} Spin Manifolds, Dirac Operators $\cal D$,  Atiyah-Singer Index Theorem  and   
S-L-W-(B)  Formula}\label {Dirac1}

%%%%%%%%%%%%%%%%%%%%%%%%%%%

{\color {blue}[I]}  Historically  the  first  $\Phi$ in this story  were {\it harmonic spinors} on a Riemannian manifold $X=(X,g)$, that are solutions $s$ of ${\cal D}(s)=0$, where $\mathcal D =\mathcal D_g$ is the {\it (Atiyah-Singer)-Dirac } on $X$.\footnote { All you have to know at this stage  about $\cal D$   is that  $\cal D$ is a certain first order differential  on sections of some bundle over $X$ associated with the tangent bundle $T(X)$.
 Basics  on  $\cal D$  are  presented  in [Min-Oo(K-Area) 2002]
 and, comprehensively, in  [Lawson\&Michelsohn(spin geometry) 1989].
  Also see sections \ref{Clifford3},\ref {Dirac4}.}
\vspace {1mm}

[{\color {blue}I}{\color {magenta}$_{yes}$}].  {\color {blue!50!magenta} The existence} of {\it non-zero  harmonic} spinors $s$ on certain  smooth manifolds $X$   follows  from {\it non-vanishing of the index of $\mathcal D$}, where this index, which is {\it independent of $g$}, 
identifies,  by the the {\it \color {blue}  Atiyah-Singer  theorem of 1963  }, with a certain (smooth) topological invariant, denoted   {\color {blue}$\hat\alpha (X)$}   (see section \ref {spin index3}).

Then    the relevant  formula involving $Sc(X)$ is  the following  {\sf algebraic identity between the squared  {\it Dirac operator} and the (coarse) {\it Bochner-Laplace operator} $\nabla^\ast \nabla$ also denoted $ \nabla^2$,} \vspace{2mm}

[{{\color {blue}I}{\color {magenta}$_{ no }$}].  \hspace {1mm} {\it  \color {blue!50!black}  \textbf {Schroedinger-Lichnerowicz-Weitzenboeck-}({\color {blue!50!black}{\it \textbf{Bochner})  {\it  \color {blue!50!black}
  \textbf{Formula}}}\footnote{ All natural selfadjoint geometric second order operators differ from the Bochner  Laplacians by {\it zero order} terms, i.e. (curvature related) endomorphisms  of the corresponding vector bundles, but it is { \sf \color {magenta} remarkable} that this  in the case of  $  \mathcal D^2$ reduces to {\it multiplication by a scalar function},  which happens to be equal to $ \frac {1}{4}Sc_X(x)$.  From  a certain perspective,  the existence of  such an   with   a wonderful combination of properties is the most amazing  aspect of the  Atiyah-Singer index theory.} 
 {\color {blue!50!black}$$\mathcal D^2=\nabla^2+\frac{1}{4}Sc,$$}}}
 shows that if $Sc>0$, then  $\mathcal D^2s=0$ implies that $s=0$, since
 $$0=\int \langle  \mathcal D^2 s, s \rangle= \int \langle  \nabla^2 s, s \rangle+  \frac {Sc}{4}||s||^2= \int ||\nabla s||^2 + \frac {Sc}{4}||s||^2,$$
where the latter identity follows by integration by parts  (Green's formula).

 By confronting these   {\color {magenta} \small \it yes }  and {\small \color {magenta}  \it no},  Andr\'e Lichnerowicz\footnote {See  [Lichnerowitz(spineurs harmoniques) 1963]} showed in 1963 that 

\hspace {26mm} {\sl $Sc(g)>0$ $\Rightarrow$  $\hat \alpha(X)=0$. }
 
\hspace {-6mm}and proved the following.

 \vspace{1mm}
 %%%%%%%%%%%%%&&&&&&&&&&&&
 
 {{\sf \large Non-Existence  Theorem Number One:} \textbf{Topological Obstruction to $Sc>0$ for $n=4k.$}  {\color {blue}\it There exists smooth  closed  $4k$-dimensional  manifolds  $X$, for all $k=1,2,...$,  which admit no metrics with $Sc>0$.} \vspace{1mm}
  
  A decade later, empowered by a general  Atiyah-Singer index theorem, Nigel Hitchin  extended Lichnerowitz'  result to manifolds of dimensions $n=8k+1$ and $8k+2$ and showed, in particular, that

  {\sf \color {blue}  the class  of manifolds $X$ with  $\hat \alpha(X)\neq 0$, that  support non-zero $g$-harmonic  spinors all metrics $g$ on $X$ by the  Atiyah-Singer theorem,  hence  no $g$ with $Sc(g)>0$ by \textbf {S-L-W-B} formula,  includes certain {\it homotopy spheres.}}   \footnote {See [AS(index) 1971],  [Hitchin(spinors)1974].}
   \footnote{Prior to 1963, one didn't  even  know if therere were   {\it simply connected} manifold that would admit {\it no metric with positive sectional  curvature} was known.  But   Lichnerowicz' theorem, saying, in fact,  that  
  
  \hspace {20mm} {\it if   $X$ is spin, then   $Sc(X)>0
  \Rightarrow \hat A[X]=0$ }
  
\hspace {-3mm}delivered     lots of {\it simply connected manifolds} $X$  that  admitted {\it no metrics with positive scalar curvatures, (see section \ref {spin index3}}).

Most of these $X$ have {\it large Betti numbers},  that, as we know nowadays, is {\it incompatible} 
with $sect.curv(X)\geq 0$, but 
   one still doesn't know if  there are   homotopy spheres not covered by Hitchin's theorem  which admit no metrics with positive sectional curvatures.  }\vspace{1mm}
 
%%%%%%%%%%%%%%%%%%%%%

\subsubsection{\color {blue} Inductive Descent  with Minimal Hypersurfaces  and Conformal Metrics}  \label {SY1}
%%%%%%%%%%%%%%%%%%%%%%
{\color {blue}[II]}  Another class of solutions $\Phi$  of geometric 
PDE, that  are  essential  for understanding  scalar curvature and that are quite different from harmonic spinors, are {\it solutions to the Plateau problem}.

 More specifically, these are   
   {\it smooth} {\sf stable minimal hypersurfaces} $Y\subset X$ that represent non-zero integer homology classes from  $H_{n-1}(X)$, $n=dim(X)$.
   
     The existence of  minimal $Y,$ possibly singular ones,  was established  by  Herbert Federer and Wendell Fleming in 1960, while the smoothness of these $Y$, that is crucial for our applications, was proven  by
      Federer in 1970 who relied on   {\it regularity}  of volume  minimizing cones of dimensions $\leq 6$ proved  by Jim  Simons in 1968.

The relevance of these {\sf  \color {blue}minimal $Y$ of codimension 1} to the  scalar curvature problems was discovered by    Schoen and Yau who  proved in 1979     that\vspace {1mm}

%%%%%%%%%%%%%%&&&&&&&&&&&

\hspace {-2mm}  {\it \color {blue} $\bigstar^{codim1}_{\sf mini}$ \hspace {2mm} } {\sf  \color {blue!40!black}  if $Sc(X)>0$  and $Y\subset X$ is a smooth  stable minimal hypersurface,}  then {\it  \color {blue!40!black}
   $Y$ 
admits a Riemannian metric $h$ with $Sc(h)>0$.}\footnote{See
   [SY(structure) 1979]:    {\sl On the structure of manifolds with positive scalar
curvature}.} \vspace {1mm}

In fact, if $dim(Y)=n-1=2$, the  stability of $Y$, that is  {\it positivity} of the second variation of the area of $Y$,  implies that (see sections \ref {2nd variation2}, \ref{warped+2})
$$\int_Y (Sc(Y, y)-Sc(X, y))dy \geq 0$$
where the scalar curvature $Sc(Y)$ refers to the metric $h_0$ in $Y$ induced from the Riemannian metric $g$ of $X$. 

Therefore, positivity of $Sc(X)$  implies positivity of the Euler characteristic of $Y$, for 
 $$4\pi \chi(Y)=\int_YSc(Y, y)dy\geq \int_YSc(X, y)dy>0.$$

If $m=n-1\geq 3$, then 
 $h$ is obtained by a {\it conformal modification} of  the metric  $h_0$ on $Y$,
 $$h_0\mapsto h=(f^2)^\frac {2}{m-2} h_0,$$
where, as  in the 1975 "conformal paper"   by Jerry Kazdan  and Frank Warner  
   $f=f(y)$ is the first eigenfunction of the {\it conformal Laplacian $L$} on 

\hspace {-6mm}$Y=(Y,h_0)$,  that is 
$$L_{conf}(f)=-\Delta(f) +\frac {m-2}{4(m-1)}f,$$
where derivation of positivity of the  $L$ from positivity of the second variation of $vol_{n-1}(Y)$  relies  on   the  {\it  \color {blue!50!black} \large Gauss formula} suitably rewritten for this purpose by 
 Schoen and Yau and   where the issuing positivity of $Sc( f^\frac {4}{m-2} h_0)$ follows, as in 
  [Kazdan-Warner(conformal)],\footnote{There is more to  this paper, than the implication 
  $L_{conf}>0\leadsto \exists g$ with  $Sc(g)>0$ on $X$.
For instance, Kazdan  and Warner prove

{\it the existence of metrics $g$ on connected   manifolds $X$, $dim(X)\geq 3$, with prescribed scalar curvatures $Sc(g, x)=\sigma(x)$, for smooth functions $\sigma(x)$, which are {\sf negative somewhere on}  $X$}  and 

{\it the existence of metrics
  with $Sc=0$ on manifolds $X$, which admits metrics with $Sc\geq 0$.}}  by a simple (for those who knows how to do this kind of things) computation.
\footnote{This computation,   probably,  going back at least hundred years,  was  brought from the field of infinitesimal geometry  to the context 
 of non-linear PDE and  {\it global analysis} by Hidehiko Yamabe in his 1960-paper {\sl On a deformation of Riemannian structures on compact manifolds}.}

 Consecutively applied   implication $Sc(X,g)>0\Rightarrow Sc(Y,h)>0$
 delivers a descending chain of closed oriented  submanifolds  
 $$X\supset Y=Y_1\supset Y_2\supset...\supset Y_i...\supset Y_{n-2}$$
of dimensions $n-i$ which support  Riemannian metrics $h_i$  with $Sc(h_i)>0$;  thus,  all   connected components of   $Y_{n-2}$ must be a spherical.

Thus, Schoen  and Yau inductively define a topological class of manifolds  ($\cal C$ in their terms) and  prove, in particular,
 the following.  \vspace {1mm}
 
%%%%%%%%%%%%%%%%%%%%%%%&&&&&&&&&&
  {\sf \large Non-Existence  Theorem Number Two     Accompanied by Rigidity Theorem.}
  {\it \color {blue!40!black} Let a compact oriented manifold $X$ of dimension $n$ dominate (a non-zero multiple of the fundamental class of) the $n$-torus, i.e, $X$  
  admits a map of non-zero 
degree to the $n$-torus $\mathbb T^n$, 
$$f:X\to\mathbb T^n.$$
If $n\leq 7$,\footnote{The dimension restriction was removed in   [Lohkamp(smoothing) 2018] and in [SY(singularities) 2017].} {\it $X$ admits no metric with $Sc>0.$}  then $X$ support no metric $g$ with $Sc(g)>0$.

Moreover, the inequality  $Sc(g)\geq 0$  for a metric $g$ on $X$, implies that $g$ is  Riemannian flat and the universal covering of $(X,g)$  is isometric to the Euclidean space 
$\mathbb R^n$.}\vspace {1mm}

 (The submanifolds $Y_i$ in this case are taken in the homology classes of transversal  $f$-pullbacks of subtori in   
 $\mathbb T^n \supset \mathbb T^{n-1}\supset... \supset\mathbb T^{n-i}\supset...\supset \mathbb T^2$.)   \vspace {1mm}

 {\it Remark.} The authors of   [SY(structure) 1979]  say in  their  paper that it was motivated by problems in general 
 relativity  communicated to one of the authors by Stephen Hawking,
 \footnote {It is shown in [Hawking (black holes) 1972], by an argument   elaborating on ideas from 
 [Penrose(gravitational collapse) 1965] and  resembling those in [SY(structure) 1979],
   that surface of  the event horizon has {\it spherical topology.} (See  [Bengtsson(trapped surfaces) 2011] for more about it.)}  
   but I  as  haven't studied this field I can't judge how  much of the current development in geometry of the  scalar curvature is rooted in  ideas originated in physics. %%%%%%%%%%%%%%%%%%%%%%%%
 
\subsubsection {\color {blue}Twisted Dirac Operators,  Large Manifolds  and Dirac  with Potentials}\label {twisted1} 

%%%%%%%%%%%%%%%%%%%%%%%%%

 The index theorem also applies to    Dirac operators ${\cal D}_{\otimes L}$ that act on spinors with values in Hermitian vector bundles $L\to X$, called {\it $L$-twisted spinors},  where {\it \color {blue!50!magenta}non-vanishing} of the index of  ${\cal D}_{\otimes L}$ and, thus the existence of non-zero $L$-twisted harmonic spinors, is {\it \color {blue!50!magenta}    ensured}  for bundles $L$  with {\it sufficiently large top dimensional Chern numbers,} essentially  regardless of the  topology of the  underlying manifold $X$ itself.

On the other hand, the {\it twisted} S-L-W-(B) formula, which now reads
$$ \mbox { ${\cal D}_{\otimes L}^2=\nabla^2_{\otimes L} + \frac{1}{4}Sc(X) + {\cal R}_{\otimes L}$},$$
shows that such spinors {\color {magenta} don't exist} if the {\it $g$-norm} of the curvature of $L$ is  {\it small} compare with the  scalar curvature of $X=(X,g)$. 
Since this norm is  inverse proportional to the size of $g$,   {\color{blue!50!magenta} large} Riemannian manifolds admit {\it topologically complicated} bundles $L$ with {\it {\color{blue!50!magenta}  small}} curvatures, which, by the above, shows, as it was observed in  [GL(spin) 1980],   that, similarly how it is with the sectional and Ricci curvatures,  

\hspace {12mm} {\color {blue} \it scalar curvatures of large manifolds  must be small.}   \vspace {1mm}

 \hspace {-5mm} This   delivers    confirmation of  the {\sf main $[Sc \ngtr 0]$ conjecture}  from the previous section for certain  
  compact manifolds $X$, with large  fundamental groups, e.g.  for $X$, which support metrics with 
  {\it  non-positive sectional curvatures}:

\vspace {1mm}
%%%%%%%%%%%%%%%%&&&&&&&&&&&&
\textbf {Spin-non-Domination theorem  of  $\kappa\leq 0$ by $Sc>0$}.   {\sf \color {blue!60!black} Non-torsion homology classes of complete manifolds  $\underline X$,   with non-positive sectional curvatures  can't be dominated  by  compact  (and also by complete)  orientable spin manifolds with $Sc>0$.
}\footnote{See [GL(spin) 1980] and  sections \ref {spin index3} and  \ref{obstructions4} for  more  specific statements and proofs.}}

In standard terms, 

{\it If a compact orientable spin  Riemannian manifold $X$ has $Sc>0$ and $\underline X$ is complete with $sect.curv(\underline X)\leq 0$, and if 
$$f:X\to \underline X$$
is a continuous map, then the image of the fundamental class $[X]\in H_n(X)$ is {\sf torsion}: 
some {\sf non-zero} multiple  $i\cdot f_\ast[X]\in H_n(\underline X)$ {\sf vanishes.}}\footnote {It's unclear if   $f_\ast[X]\in H_n(\underline X$ can be non-zero, yet  (odd?) torsion.} 

 For instance, 
\vspace {1mm}
 
 {\sf if $\underline X$  is compact of dimension $n=dim(X)$, then all continuous maps $f:X\to \underline X$ 

 have {\it  zero} degrees. }

\vspace {2mm}

{\sf \textbf {Homotopy Invariance   of  Obstructions to $Sc>0$ that  Issues from {\color {magenta}  $\otimes$} in $\mathcal D_ {\color {magenta}}$.} {\it Non-vanishing of topological invariants  } delivered by   the  twist in $\mathcal D_ {\color {magenta}\otimes L}$} that prevent the existence of metrics with $Sc>0$
 are 
 {\it stable under toplogical domination} that   is,  recall,  a map $X\to \underline X$ of degree $\pm 1$ between orientable  manifolds, such that    \vspace {0.6mm}
 
  \hspace {-4mm} {\it if such an invariant  doesn't vanish for 
 $\underline X$, then it  doesn't vanish for $X$ either.}
  \vspace {0.6mm}

(An  instance of such an invariant is the  $\smile$-product  homomorphism
$\bigwedge^n H^1(X)\to H^n(X)$, $n=dim(X)$ behind the Schoen-Yau $[Sc>0]$-non-existence theorem  in section \ref{SY1} for manifolds mapped  to  the  n-tori)
   \vspace {0.6mm}

This is similar to what happens to invariants issuing by the geometric measure theory but very much unlike to  
 those  coming from the untwisted index theorem, namely to non-vanishing of $\hat \alpha(X)$: the connected sum of two copies of an  $X$ with opposite  orientations satisfies: $\hat \alpha(X\#(-X))=0$.
 
 In fact, {\sf if $X$ is simply connected of dimension $n\geq 5$,} then  {\it  $\hat \alpha(X\#(-X))$ does  admit a metric with $Sc>0$. }\footnote {I am uncertain about $n=4$.}\vspace {1mm}

{\it \color {blue}  Dirac with Potentials.} The contribution of the connection of $L$  to the Dirac operator can be seen as a vector potential added to  $\cal D$ twisted with a the trivial bundle of rank $=rank(L)$. 

Besides, this   there are other kinds of {\it   --} zero order terms  -- that can significantly influence geometric effects of $\cal D$.

As far as the scalar curvature is concerned, the first (to the best of my knowledge)  potential of this kind  ({\it Cartan connection})  was introduced  by Min-Oo in his proof of the  positive mass theorem for hyperbolic spaces, [Min-Oo(hyperbolic) 1989],  and, 
 recently, applications of  {\it Callias-type} potentials in the work by Checcini, Zeidler  and Zhang have  significantly extended  the  range of the Dirac-theoretic applications to the scalar curvature  
 problems.\footnote {Exposition  of Dirac operators with potentials, especially of their recent applications to manifolds with boundaries,  are, regretfully, missing from our lectures. The reader has to turn to the original papers by   Checcini, Zeidler,  Zhang and
  [Guo-Xie-Yu(quantitative K-theory) 2020]. Also  we say very little about the  mass/energy theorems for hyperbolic spaces extending that in [Min-Oo(hyperbolic) 1989]; we refer  for this subject matter to [Chrusciel-Herzlich [asymptotically hyperbolic)  2003], [Chrusciel-Delay(hyperbolic positive energy) 2019],  [Huang-Jang-Martin(hyperbolic mass rigidity) 2019]  and  [Jang-Miao(hyperbolic mass) 2021] where one can find further references.}

%%%%%%%%%%%%%%%%%
\subsubsection {\color{blue} Stable  $\mu$-Bubbles} \label{bubbles1}
%%%%%%%%%%%%%%%%%%%%%%%%

  In general, $\mu$-bubbles} $Y\subset X$,  are solutions of   the "non-homogeneous  Plateau equation" 
  $$mean.curv(Y,y)=\mu(y)$$
  for a given function $\mu(x)$ on $X$.
    
    What we  deal with in this paper are {\it stable $\mu$-bubbles} 
     that are  {\it  local minima} of the functional 
  $$Y\mapsto vol_{n-1}(Y) - \mu(Y_<)$$
 where  $\mu$ is a Borel measure on $X$   and $Y_<\subset X$ is  a   region in $X$ with boundary $ \partial Y_<=Y$ (see section  \ref {bubbles5}). 
 
  Often our measure is "continuous", i.e. representable  as    $\mu(x)dx$, for a continuous function $\mu(x) $ on $X$,and 
  all basic existence and regularity  properties of minimal hypersurfaces  automatically   extend to $\mu$-bubbles in this case.
 
 And what is especially useful for our purposes, is that the Schoen-Yau form of the   the second variation  formula neatly extends to  $\mu$-bubbles with  continuous (and some discontinuous)  $\mu\neq 0$.

 \vspace {1mm}
 
  {\it Example/non-Example.} The unit sphere $S^{n-1}\subset \mathbb R^n$ (with the mean curvature $n-1$)  around the  origin is a {\it stable $\mu$-bubble} for the measure $\mu(x)=(n-1)||x||^{-1}dx$ in $\mathbb R^n$  and  the same sphere also 
   is the $\mu$-bubble for $\mu(x) =(n-1)dx$; but this $\mu $-bubble is an {\it unstable} one.

  \vspace {1mm}

 A significant gain achieved  with $\mu$-bubbles  compared with   the "plain" minimal hypersurfaces  is  due to   the {\it flexibility in the choice of}   $\mu$, which can be  {\sf \color {black!40!blue} adapted to the geometry} of $X$,   similarly  to  how  one uses     {\it twisted} Dirac operators 
 $\mathcal D_{\otimes L}$  on $X$ with " {\sf \color {black!40!blue}adaptable" unitary bundles}   $L\to X$.
 
 For example, one obtains this way the following version of   Schoen-Yau theorem  {\color {blue} $\bigstar$  }  from section \ref{SY1}.
 \vspace {1mm} 
%%%%%%%%%%%%%&&&&&&&&&&&&&

{\color {blue} \FiveStarShadow$_{bbl}^{codim1}$} {\sf Let $X$ be a complete Riemannian $n$-manifold with {\it uniformly positive} scalar curvature, i.e, $Sc(X)\geq \sigma>0$.} If $n\leq 7$, then \vspace {1mm} 
 
 {\it \color {blue!39!black}  $X$  can be exhausted by compact domains with smooth  boundaries,
$$V_1\subset V_2\subset  ...\subset V_i\subset ... X,\mbox { }  \bigcup_iV_i=X,$$
where the boundaries $\partial V_i$, for all $i=1,2,...$,}  {\sf  \color {blue!59!black}   admit metrics with  positive scalar curvatures}.

(Here, as in section  \ref {SY1}, this needs  additional  analytical work to be extended to $n\geq7$.)

%%%%%%%%%%%%%%%%%%%%%%%%% %%%%%%%%%%

\subsubsection {\color{blue}Warped FCS-Symmetrization of Stable Minimal  Hypersurfaces and $\mu$-Bubbles.}\label{FCS1}

%%%%%%%%%%%%%%%%%%%%%%%%%%%%%%%%

Positivity of the conformal Laplacian 
$-\Delta +\frac {m-2}{4(m-1)}Sc$  doesn't fully reflect the positivity of the second variation of the volume $vol_{n-1}(Y)$,
where the former actually yields positivity of the  $-\Delta +\frac {1}{2}Sc$,  which is, 
a priori,  smaller then  $-\Delta +\frac {m-2}{4(m-1)}Sc$, since $-\Delta\geq 0$ and $\frac {1}{2}>\frac {m-2}{4(m-1)}.$

{\sf \color {magenta} Remarkably},  positivity of the  $-\Delta +\frac {1}{2}Sc$ on $Y=(Y, h_0)$ {\it neatly implies} positivity of the scalar curvature of the (warped product) metric $h^\rtimes=h_0(y)+\phi^2(y)dt^2$  for the first eigenfunction $\phi$ of $-\Delta +\frac {1}{2}Sc$,  where this metric is defined  on 
the products of $Y$ with the real  line $\mathbb R$  and with the unit circle $S^1(1)=\mathbb T= \mathbb R)/\mathbb Z$, and where the resulting Riemannian manifolds are denoted 
$$ \bar Y^\rtimes = Y\rtimes \mathbb R=(Y\times \mathbb R, h^\rtimes)\mbox {  and }  Y^\rtimes= Y\rtimes \mathbb T=  \bar Y^\rtimes/\mathbb Z.$$

  In fact, 
if $(-\Delta +\frac {1}{2}Sc)(\phi)=\lambda \phi$ with $\lambda\geq 0$,  then 
$$Sc( h^\rtimes (y, t))= Sc(h_0, y) -\frac {2}{\phi} \Delta \phi(y)=\frac {2}{\phi}\left(- \Delta+\frac {1}{2}Sc(h_0,)\right)(\phi)=\lambda>0m$$ see sections \ref{bubbles5}.

The operation  
$$Y\leadsto  Y^\rtimes$$ 
  is   applied in the present case  to   stable minimal hypersurfaces $Y\subset X$,  
where the resulting passage 
 $X\leadsto  Y^\rtimes$ can be regarded as {\it symmetrisation} of    $X$  (or rather of  infinitesimal neighbourhood of $Y\subset X$), because

  {\it  the metric $ h^\rtimes$  is invariant under the natural action  of $\mathbb T$  on $ Y^\rtimes$ and    
 $$ Y^\rtimes/\mathbb R=Y  \subset X$$}.

This  $ h^\rtimes=h_0(y)+\phi^2(y)dt^2$   defined with the first eigenfunction $\phi$ of the   $-\Delta +\frac {1}{2}Sc$  on $Y$  was introduced  by Doris  \textbf Fischer-\textbf Colbrie and  Rick  \textbf Schoen\footnote {{\sl The structure of complete stable
minimal surfaces $Y$ in 3-manifolds of non-negative scalar curvature}.}  who used it for \vspace {1mm}

 {\it classification  of complete  stable minimal surfaces in $3$-manifolds $X$ with 
 
 $Sc(X)\geq 0$, including $X=\mathbb R^3$}. \vspace {1mm}

Then $ h^\rtimes$  was used in  [GL(complete) 1983], where, with an incorporation of 
Schoen-Yau's inductive descent, this  allowed  higher dimensional applications of the following kind.

\vspace {1mm}

{\sf Given a Riemannian metric $g$ on a  product manifold $X=X_0\times \mathbb T^k$,  a consecutive   
symmetrization 
$$X=  X_{0}\leadsto  X_{1}= Y_1^\rtimes/\mathbb Z \leadsto X_{2} =Y_2{^\rtimes/\mathbb Z} \leadsto...$$
delivers  a  {\it $  \mathbb T^k$-invariant}   metric $\bar g$ on  $\bar X_{k} =Y_{- k}\times \mathbb T^k$, where $Y_{- k}\subset X$ is a submanifold of codimension $k$ which is {\it homologous to} $X_0=X_0\times t_0\subset  X$ and such that the ($\mathbb T^k$-invariant) scalar curvature  $Sc(\bar g)$ on $\bar X_{k}$  is {\it bounded from below} by $Sc(g)$  on $Y_{-k}=\bar X_{k}/\mathbb T^k\subset X$.}\vspace {1mm}

Thus, for instance, one obtains   a somewhat different proof of the Schoen-Yau theorem
for $n\leq 7$:\vspace {1mm}

 {\sf  {\it no metric} $g$ on   $X=\mathbb T^n$ can have  $Sc(g)>0$, because all {\it  $\mathbb T^n$-invariant} metrics on   $\mathbb T^n$  are Riemannian flat.}\vspace {1mm}

{\it Non-Compact Case.} An  apparent   bonus of this argument is its   applicability to {\it non-compact complete manifolds}. \vspace {1mm}
%%%%%%%%%%%%%%%%%&&&&&&&&&&&&&

\textbf {Example:  Non-domination of $ \mathbb T^n$ by $Sc>0$.}   {\sf The $n$-torus  admits no domination by {\it complete} manifolds  $X$ with $Sc(X)>0$.}\footnote{Here, as at  other similar occasions,  singularities of minimal hypersurfaces and of $\mu$-bubbles create complications  for $n=dim(X)\geq 8$.

In the present case, if $X$ is spin, this non-domination property follows by a Dirac operator argument from section 
 6 in [GL(complete)  1983].
 
 If $n=8$ the  perturbation argument from  [Smale(generic regularity) 2003]  takes care of things.
 
 If $n=9$ one can still   apply     Dirac operators to  non-spin 
  manifolds, exploiting  the fact that  
  singularities of  hypersurfaces are at most 1-dimensional,    while the obstruction to spin (the second Stiefel-Whitney class) is 2-dimensional, see section 5.3 in [G(billiards) 2014].

If $n\geq 8$ the recent desingularization results presented in [Lohkamp(smoothing) 2018] and in [SY(singularities) 2017] 
apply to all $X$.}

 For instance, {\sf if a closed subset in the torus $Y\subset \mathbb T^n$ is contained in a topological ball $B\subset \mathbb T^n$,} then 
 \vspace {1mm}
 
 \hspace {6mm}{\sl the complement $T^n\setminus Y$   admits no complete metric with 
$Sc >0$. }

\vspace {1mm}\vspace {1mm}

 The main role of  the above $\mathbb T^k$-symmetrization, however,   is {\it not for the proof of topological  non-existence theorems} of metrics with  $Sc>0$
on  closed  or non-compact  complete manifolds,  but for the  {\it geometric  study of such metrics} on,  possibly  {\it non-compact and  non-complete,} manifolds $X$.

 In fact, this symmetrization  applies
  to stable minimal hypersurfaces $Y\subset X$  {\it  with prescribed as well  as   free boundaries}, say with  $\partial Y\subset \partial X$   and also to 
   {\it stable  $\mu$-bubbles}. \footnote {See section 12 in  [GL(complete) 1983],  [G(inequalities) 2018] and 
  sections \ref {separating3}, \ref {bubble5}).}       
  \vspace {1mm}

%%%%%%%%%%%%%%%%%%%%

\subsubsection{\color{blue}  Averaged Curvature of Levels of Harmonic Maps} \label {harmonic1}
 
 %%%%%%%%%%%%%%%%%%%

 Recently, Daniel  Stern [Stern(harmonic) 2019] found a version of the 3d Schoen-Yau argument for the levels of {\it non-constant harmonic maps} $f:X\to \mathbb T^1$, where,  instead of  the second variation  formula for $area(Y)$,  one uses

{\sf  the Bochner identity,  which expresses the  Laplace   of  the norm of the  

gradient of $f$ in terms of the Hessian of $f$ and the Ricci curvature,}  
 $$\frac {1}{2}\Delta|\nabla f|^2)= |Hess(f)|^2 +Ricci_X(\nabla f, \nabla f).$$

Thus,  Stern proved that the average Euler characteristics of these levels $Y_t=f^{-1}(t), t\in \mathbb T^1$ satisfies:

\textbf {Harmonic Map Inequality.}  
$$4\pi\int_{\mathbb T^1}\chi(Y_t)dt \geq \int_{\mathbb T^1} dt\int_{Y_t}( |df(y,t)|^{-2} |Hess f(y,t)|^2 +Sc(X, (y,t)))dy.$$
This  shows that  
$$4\pi\int_{\mathbb T^1}\chi(Y_t)dt\geq \int_{\mathbb T^1}dt \int_{Y_t}Sc(X, (y,t))dy. $$
and implies, among other things,  that\vspace {1mm}

{\sl if  the universal covering of a  compact   3-manifolds with positive scalar curvatures is connected at infinity, then 
the one-dimensional  cohomology $H^1(X;\mathbb Z)$ vanishes.}\footnote{It is known that compact 3-dimensional manifolds with $Sc>0$ are connected sums of space forms and $S^2\times S^1$, see [GL(complete) 1983] and 
 [Genoux(3d classification) 2013].}\vspace {1mm}

  Indeed,  if $H^1(X;\mathbb Z)\neq 0,$  then $X$ admits a non-constant harmonic map to the circle $\mathbb T^1$, where  non-singular levels $Y_t\subset \mathbb X $ {\it can't  contain    spherical  components}, because  
  lifts of such a component
   to the universal covering of $X$ would bound  balls on which  (the lift of) $f$ would be constant by the maximum principle for harmonic functions. \footnote {In this respect, the surfaces    $Y_t$ are  radically different from  minimal surfaces and $\mu$-bubbles which tend to localize 
 around narrow necks in $X$, e.g. in  "thin" connected sums  $\mathbb T^3\# S^3$ described in  section \ref {thin1}.}
  \vspace {1mm}

{\it Vague Questions.} {\sf Is there an algebraic link between S-L-W-(B)  and the above Bochner formula that would connected  Dirac  operators with  harmonic maps?

 Do {\it  Dirac harmonic} and/or similar maps bear a relevance to the scalar curvature problem?}

%%%%%%%%%%%%%%%%%%%%

\subsubsection {\color {blue}  Seiberg-Witten Equation}\label {Seiberg1}

%%%%%%%%%%%%%%%%%%%%
  The third kind of $\Phi$ are solutions to the  4-dimensional {\it Seiberg-Witten equation} of 1994, that is  the  Dirac  equation coupled with a certain non-linear equation and where the relevant formula is essentially the same as in {\color {blue}[I]}.

Using these, Claude  LeBrun\footnote {[LeBrun(Yamabe)  1999]:  {\sl Kodaira Dimension and the Yamabe Problem.}}     established a  non-trivial (as well as sharp)  \vspace{1mm}

 {\it   \textbf {Fundamental 4D lower bound} \color {blue} on $\int_XSc (X,x)^2 dx$ for Riemannian manifolds $X$ diffeomorphic
to  algebraic surfaces of general type.}

\subsubsection {\color {blue} Hamilton-Ricci Flow}{Hamilton1}

The  {\it Hamilton Ricci flow} $\Phi=g(t)$ of Riemannian metrics on a manifold $X$, that is defined by  a  {\it parabolic}  system of equations,  also  
delivers a geometric information on the scalar curvature, where the main algebraic identity for $Sc(t)=Sc(g(t))$
reads
$$\frac {dSc(t)}{dt}=\Delta_{g(t)}Sc(t) + 2Ricci(t)^2\geq \Delta_{g(t)}Sc(t) +\frac{2}{3} Sc(t)^2,$$
which implies by the maximum principle  that the minimum of the scalar curvature grows with time as follows:
$$Sc_{\min}(t)\geq  \frac {Sc_{\min}(0)}{1-\frac{2t Sc_{min}(0)}{3}}.$$

\vspace{1mm}

If $X=(X, \underline g)$ is a closed $3$-manifold of constant sectional curvature $-1$, then, using the Ricci flow, Grisha Pereleman proved 

    \textbf {Sharp 3D Hyperbolic Lower Volume Bound.} {\color {blue} All Riemannian metrics $g$ on $X$ with $Sc(g)\geq -6=Sc(\underline g)$  satisfy
  
  \hspace {35mm} $Vol(X,g)\geq Vol(X, \underline g)$}.
  
  (See Proposition 93.9 in 
[Kleiner-Lott(on Perelman's) 2008].)\vspace{1mm}

And, more recently,  Richard Balmer, Paula Burkhardt-Guim and Man-Chun Lee, Aaron Naber  and  Robin Neumayer   applied the Ricci flow for regularization of of  (limits of) metrics with $Sc\geq \sigma$.\footnote {See  [Bamler(Ricci flow proof) 2016],  [Burkhart-Guim(regularizing Ricci flow) 2019],   [Lee-Naber-Neumayer](convergence) 2019] and  section 
\ref {C0-limits3}.}

(The logic of the Ricci flow, at least on the surface of things, is quite different 
from how it goes in the above  three  cases that rely on {\it elliptic} equations:

{\sf the   quantities $\Phi$ in the former  result from    geometric or topological 
{\it complexities}  of underlying manifolds $X$, that is  necessary for the very existence  of these $\Phi$,  while the   Ricci flow,  as  a road roller, leaves a uniform terrain behind itself as   it  crawls along  erasing  complexity.})
\vspace{1mm}

{\it Question.}  Do 3D-results obtained with the  Ricci flow generalize to $n$-manifolds which have  $Sc\geq \sigma$  and which   come with {\it free isometric actions 
of the tori} $\mathbb T^{n-3}$?

  For instance, let $X^3$ be a 3-dimensional Riemannin manifold which admits a hyperbolic metric $\underline g$ with sectional
  curvature $-1$
  and let  $X=X^3\rtimes  \mathbb T^1$ be a warped  product (with  $\mathbb T^1$-invariant metric), such  that   $Sc(X)\geq -6$.

Is the volume of  $X^3=X/\mathbb T^1$ is bounded from below by that of $(X^3, \underline g)$?

(It is not even clear if  the inequality  $Sc(X^3\rtimes  \mathbb T^1)\geq -6$
imposes {\it any lower bound} on the Riemannin metric  $g$ of $X^3$. Namely, 

Can  
such a $g=g_\varepsilon$ satisfy  $g\leq \varepsilon \underline g$ for a given $\varepsilon>0$?\footnote {An elementary  proof of  such a bound on $g$ is suggested in  [G(foliated) 1991].})

%%%%%%%%%%%%%%%%%%%%

\subsubsection{ \color {blue}Modifications of Riemannian Metrics by a Single Function} \label{single1}

Riemannian metrics $g$ on an $n$-manifold $X$  are  given locally   by $\frac {n(n-1}{2}$  functions 
$g_{ij}(x)$, where the scalar curvature $Sc(g) $   is a (messy) non-linear function of  these $g_{ij}$ and their   first and second derivatives. 

  There are several constructions of Riemannian metrics  on $X$  and of  modifications of a given metric $ g_0$ on $X$ by means of a {\it single}  function $\phi(x)$,   where the the scalar curvature of the resulting metric $g(\phi)= g(\phi,   g_0)$
is  expressed by a "nice"  non-linear second order differential  applied to $\phi$.

The simplest and most studied case of this is  the  conformal transformation $g\mapsto \varphi^2g$, where 
for $n\geq 3$ the scalar curvature of this metric  is given by the (Yamabe?) equation
$$Sc(\varphi^2g_0)=-\frac{4(n-1)}{n-2}\varphi^\frac{n+2}{2}\Delta \varphi^\frac{n-2}{2} +\varphi^2 Sc(g_0),$$
where $\Delta=\Delta_{g_0}$ is the Laplace  on
 functions $\phi=\phi(x)$  on  the Riemannian manifold   $(X,g_0)$.

We  present some     properties  of  this equation,  due to Jerry Kazdan and Frank  Warner, in section \ref {conformal2},  which are used  
in   the proof of Schoen-Yau's {\it non-existence}  theorem  for metrics with $Sc>0$ on tori in  sections \ref  {SY1}, \ref{SY+symplectic2}.

Also we briefly discuss in  \ref {conformal2} similar transformations of metrics, where the scaling 
takes place only in some preferred directions,  e.g. 
 in a single direction,  where the scalar curvature satisfies a non-linear parabolic    (Bartnik-Shi-Tamm) equation,   special  {\it solutions} of which 
  used for  the proofs of {\it non-extension} theorems for metrics with $Sc>0$, see   section \ref {fill-ins3}. 

 Finally,  recall   K\"ahler metrics defined with single functions  via the $\partial\bar\partial$ , where, as we   mention in section \ref {SYS1},  Yau's solution of the  Calabi conjecture delivers "interestingly thick" metrics with $Sc>0$ on complex  algebraic manifolds.   
 \vspace{2mm}
%%%%%%%%%%%%%%%%%%%%%%%%%

\section{Curvature Formulas  for  Manifolds and  Submanifolds.} \label{formulas2}

%%%%%%%%%%%%%%%%%%%%%%%%%%

 We enlist in this section several classical formulas of Riemannian  geometry   and indicate their (more or less) immediate  applications.

\subsection {\color {blue} Variation of the Metrics and Volumes in  Families of Equidistant Hypersurfaces} \label{equidistant2}

{\color {blue}(2.1. A) {\large \sf Riemannian  Variation Formula.}  } Let $h_t$, $t\in [0, \varepsilon]$, be a family of Riemannian metric on an $(n-1)$-dimensional manifold $Y$ and let us incorporate $h_t$ to the metric  $g=h_t+ dt^2$ on $Y\times [0,\varepsilon].$

 Notice that an arbitrary  Riemannian metric on an $n$-manifold  $X$  admits such a representation in
normal geodesic coordinates in a small (normal) 
neighbourhood of any given compact hypersurface  $Y\subset X$.

 The $t$-derivative of $h_t$ is equal to {\sl twice the   second fundamental form}  of the hypersurface $Y_t=Y\times \{t\}\subset Y\times [0,\varepsilon]$, denoted and regarded as a quadratic differential  form on  $Y=Y_t$, denoted
$$A^\ast_t=A^\ast(Y_t)$$  
  and regarded as a quadratic differential  form on  $Y=Y_t$.

 In writing,%where
%the notation $A(Y_t)$ is reserved for the corresponding shape ; 
$$\partial_\nu h=\frac{dh_t}{dt}=2A_t^\ast,$$
 or, for brevity, 
 $$\partial _\nu h=2A^\ast,$$
where 

\hspace {15mm} {\it $\nu$ is the unit normal field to $Y$ defined as $\nu=\frac {d}{dt}$.}\vspace{1mm}

In fact, if you wish, you can take this formula for the definition of  the   second fundamental form of $Y^{n-1}\subset X^n$.

Recall, that the {\it principal values}  $\alpha^\ast_i(y)$, $i=1,...,n-1$, of the quadratic form  $A^\ast_t$ on the tangent space $T_y(Y)$,  that are the values of this form on the orthonormal vectors $\tau^\ast_i\in T_i(Y)$, which {\it diagonalize} $A^\ast$,  are called
{\it the principal curvatures} of $Y$, and that the sum of these is called {\it the mean curvature} of $Y$,
$$mean.curv(Y,y)=\sum_i\alpha^\ast_i(y),$$
where, in fact ,
$$\sum_i\alpha^\ast_i(y)=trace (A^\ast) =\sum_iA^\ast(\tau_i) $$
for {\it all } orthonormal tangent frames $\tau_i$ in $T_y(Y)$ by the Pythagorean theorem.

\vspace {1mm}

{\sc \color {magenta}  Sign Convention.}
 The first derivative of $h$  changes  sign under reversion of the $t$-direction. Accordingly the sign of  the quadratic form $A^\ast(Y) $ of a hypersurface $Y\subset X$  depends on the {\it coorientation} of $Y$ in $X$,  where our  convention is such that  \vspace {1mm}

{\sf the   boundaries of {\it convex}  domains have {\it positive (semi)definite} second fundamental 
forms $A^\ast$, also denoted $\mathrm {II}_Y$,
hence, {\it positive} mean curvatures}, with respect to 
{\color {magenta} \it the  outward} normal vector fields.\footnote{At some point, I found out  to my dismay, that this is opposite to the standard convention  in the  differential geometry. I apologise to the  readers who are used to the commonly accepted  sign.} 
 
\vspace {1mm}
  
{\color {blue}(2.1.B)    {\large \sf First Variation Formula.}} This concerns the $t$-derivatives of the $(n-1)$-volumes of  domains $U_t=U\times \{t\} \subset Y_t$, which are computed  by tracing the above {\color {blue}(I) } and which are related to the mean curvatures as  follows.

$$\partial_\nu vol_{n-1}(U)= \frac{dh_t}{dt}vol_{n-1}(U_t)=\int _{U_t} mean.curv(U_t)dy_t\footnote{This come with the {\it \color {magenta} minus} sign in most (all?) textbooks, see e.g. [White(minimal) 2016], [Cal(minimal( 2019].}\leqno \mbox { {\Large \color {blue}$[ \circ_U]$}}  $$
where $ dy_t$ is the volume element in $Y_t\supset U_t$.

This can be equivalently expressed with the fields $\psi\nu=\psi \cdot \nu$  for $C^1$-smooth functions $\psi=\psi(y)$ as follows
 $$\partial_{\psi\nu} vol_{n-1}(Y_t)=\int _{Y_t} \psi(y)mean.curv(Y_t)dy_t \footnote{This remains true  for Lipschitz functions but if $\psi$ is (badly) non-differentiable, e.g.  it is equal to  the characteristic function of a domain
 $U\subset Y$, then the derivative $\partial_{\psi\nu} vol_{n-1}(Y_t)$ may become (much) larger than this integral.  }\leqno \mbox { {\Large \color {blue}$[ \circ_\psi]$}}  $$

 \vspace {1mm}

Now comes the first formula with the Riemannian curvature in it. \vspace{1mm}

%%%%%%%%%%%%%%%%%%%%

\subsection {\color {blue}Gauss' Theorema Egregium}\label {Gauss2}
%%%%%%%%%%%%%%%%%%%%%%%%%

Let $Y\subset X$ be a smooth hypersurface in a Riemannian manifold $X$. Then the sectional curvatures of $Y$ and $X$ on a tangent  2-plane $\tau \subset T_y(Y)\subset T)y(X)$ $y\in Y$,  satisfy
$$ \kappa(Y,\tau)= \kappa(X,\tau)+ \wedge^2A^\ast(\tau^{}),$$ 
where $\wedge^2A^\ast(\tau)$ stands for the product of the two principal values of the second fundamental form  form $A^\ast=A^\ast(Y)\subset X$ restricted to the plane $\tau$,
$$ \wedge^2A^\ast(\tau)= \alpha_1^\ast(\tau)\cdot \alpha^\ast _2(\tau).$$

This,  with the definition the scalar curvature by the formula $Sc=\sum \kappa_{ij}$,  implies   that 
$$Sc(Y,y)= Sc(X,y)+ \sum_{i\neq j}\alpha^\ast_i(y)\alpha^\ast_j(y)-\sum_i \kappa_{\nu,i},$$
where:

$\bullet$   $\alpha^\ast_i(y)$, $i=1,...,n-1$ are the (principal)  values of the second fundamental form on the diagonalising  orthonormal frame of vectors $\tau_i$ in $T_y(Y)$;

$\bullet$   $\alpha^\ast$-sum is taken over all ordered pairs $(i,j)$ with $j\neq i$;

$\bullet$ $\kappa_{\nu,i}$ are the sectional curvatures of $X$ on the bivectors $(\nu, \tau_i)$ for $\nu$ being a unit (defined up to $\pm$-sign)  normal vector to $Y$;
 
$\bullet$  the sum of $\kappa_{\nu,i}$  is equal to the value of the Ricci curvature of $X$ at $\nu$,
$$\sum_i \kappa_{\nu,i}=Ricci_X(\nu,\nu).$$
(Actually, Ricci can be defined as this sum.)

Observe that both sums are independent of  coorientation of $Y$ and that  in the case of $Y=S^{n-1}\subset \mathbb R^n=X$ 
this gives the correct value $Sc(S^{n-1})=(n-1)(n-2)$.

Also observe that 
$$\sum_{i\neq j}\alpha_i\alpha_j=\left(\sum_i\alpha_i\right)^2-\sum_i\alpha_i^2, $$
which shows that
$$Sc(Y)= Sc(X) +(mean.curv(Y))^2- ||A^\ast(Y)||^2-Ricci(\nu,\nu).$$

In particular, if $Sc(X)\geq 0$ and $Y$ is {\it  minimal}, that is  $mean.curv(Y)=0$, then 
$$ Sc(Y)\geq -2Ricci(\nu, \nu).\leqno {\rm \color  {blue}({Sc\geq-2Ric})}$$

{\it Example.} The scalar curvature of a  hypersurface $Y\subset \mathbb R^n$ is expressed in terms of the mean curvature of $Y$, the (point-wise) $L_2$-norm  of the  second fundamental form  of $Y$  as follows. 
$$Sc(Y)= (mean.curv(Y))^2- ||A^\ast(Y)||^2$$
for  $||A^\ast(Y)||^2=\sum_i(\alpha_i^\ast)^2,$
while  $Y\subset S^n$ satisfy
$$Sc(Y)= (mean.curv(Y))^2- ||A^\ast(Y)||^2+(n-1)(n-2)\geq (n-1)(n-2)-n\max_i(c^\ast_i)^2.$$
It follows that  {\it minimal} hypersurfaces $Y$ in $\mathbb R^n$, i.e. these with $mean.curv(Y)=0$, 
have {\it negative scalar curvatures}, while hypersurfaces in the $n$-spheres with all principal values 
$\leq  \sqrt{n-2}$ have $Sc(Y)>0$.

 \vspace {1mm}

 Let $A =A(Y) $ denote  {\it the shape } that is the  symmetric    on $T(Y)$  associated with $A^\ast$ via  the Riemannian scalar product   $g$  restricted from $T(X)$ to $T(Y)$,
 $$A^\ast (\tau,\tau)=\langle A(\tau),\tau  \rangle_g  \mbox{ for all }\tau\in T(Y).$$
 
 %%%%%%%%%%%%%%%%%%
 \subsection {\color {blue} Variation of the   Curvature of Equidistant Hypersurfaces  and Weyl's  Tube  Formula}\label {Weyl2}
%%%%%%%%%%%%%%%%%%%%%%%%

{\color {blue}(2.3.A)    {\sf\large  Second Main Formula of Riemannian Geometry.}}\footnote{The first main formula is {\it Gauss' 
Theorema Egregium. }} Let $Y_t$ be a family of hypersurfaces $t$-equidistant to a given $Y=Y_0\subset X$.
Then the shape operators $A_t=A(Y_t)$  satisfy:

$$\partial_\nu A=\frac{dA_t}{dt}=-A^2(Y_t)-B_t,$$
where $B_t$ is the symmetric  associated with the  quadratic differential form $B^\ast$ on $Y_t$,  the values of which on the  tangent unit vectors 
$\tau\in T_{y,t}(Y_t)$  are equal to the values of the  {\it sectional curvature} of $g$ at (the 2-planes spanned by)  the bivectors 
$\left (\tau, \nu=\frac {d}{dt}\right)$.\vspace{1mm}

 {\it Remark.}  Taking this  formula  for the   {\it definition} of the sectional  curvature, or just systematically using it,  delivers fast  clean  proofs  of the basic {\it Riemannian comparison  theorems} along with their standard corollaries, by far more efficiently than what is allowed by the 
cumbersome   language of      Jacobi fields  lingering on the pages of most textbooks on Riemannian geometry. \footnote{ Thibault Damur    pointed out to me that  this formula, along with the rest displayed on the pages in this section, are systematically  used by physicists in books  and in articles  on  relativity.  For instance,  what we present  under heading of 
 "Hermann Weyl's  Tube  Formula",  appears  in   [Darmos(Gravitation einsteinienne) 1927] with the  reference to Darboux' textbook
 of 1897.} 
\vspace {1mm}

 Tracing this formula yields \vspace {1mm}

 {\color {blue}(2.3.B)  {\sf \large  Hermann Weyl's  Tube  Formula.}} 
 $$trace\left(\frac{dA_t}{dt}\right)=-||A^\ast||^2-Ricci_g\left (\frac{d}{dt},\frac{d}{dt}\right),$$
 or $$trace (\partial_\nu A)= \partial_\nu trace (A) =-||A^\ast||^2 -Ricci(\nu,\nu),$$
 where
 $$||A^\ast||^2=||A||^2=trace( A^2),$$
where, observe,  $$trace (A)=trace (A^\ast) = mean.curv=\sum_i \alpha^\ast_i$$
  and where $Ricci$ is the quadratic form on $T(X)$ the value of which on a unit vector  $\nu\in T_x(X)$ is equal to the trace of the above  $B^\ast$-form (or of the  $B$)  on the normal hyperplane $\nu^\perp \subset T_x(X)$ (where $\nu^\perp=T_x(Y)$ in the present case).

Also observe  -- this  follows from the definition  of the scalar curvature as $\sum \kappa_{ij}$ -- that 
 $$Sc(X)=trace (Ricci)$$
  and  that the above formula  $Sc(Y,y)= Sc(X,y)+ \sum_{i\neq j}\alpha^\ast_i\alpha^\ast_j-\sum_i \kappa_{\nu,i}$
  can be rewritten as 
  $$Ricci(\nu,\nu)=\frac{1}{2} \left(Sc(X)-Sc(Y) -\sum_{i\neq j}\alpha^\ast_{i}\cdot \alpha^\ast_j\right)=$$
$$=\frac{1}{2} \left(Sc(X)-Sc(Y) -(mean.curv(Y))^2+||A^\ast||^2\right)$$
where, recall,   $\alpha^\ast_i= \alpha^\ast_{i}(y)$, $y\in Y$,  $i=1,...,n-1$, are the  principal curvatures of $Y\subset X$, where $mean.curv(Y)=\sum_i \alpha^\ast_{i}$ and where
 $||A^\ast||^2=\sum_i  (\alpha^\ast_{i})^2$.

\vspace {1mm}
%%%%%%%%%%%%%%%%%%%%%%
 \subsection {\color {blue}  Umbilic Hypersurfaces  and Warped Product Metrics} \label {warped2}
%%%%%%%%%%%%%%%%%%%%%
A hypersurface $Y\subset X$ is called {\it umbilic }  if all principal curvatures of $Y$ are mutually equal at all points in $Y$. 

For instance, spheres in the {\it standard} (i.e. complete simply connected)  {\it spaces with constant curvatures}  (spheres $S_{\kappa>0}^n$, Euclidean spaces $\mathbb R^n$ and hyperbolic spaces $\mathbf H^n_{\kappa<0}$) are  umbilic.

In fact these are special case of the following class of spaces .\vspace{1mm}

{\it Warped Products.}  Let  $Y=(Y,h) $ be  a smooth  Riemannian (n-1)-manifold and $\varphi=\varphi(t)>0$, $t\in [0, \varepsilon]$  be   a smooth positive function.   Let 
$g=h_t +dt^2= \varphi^2h+dt^2$ be the corresponding metric on $X=Y\times [0,\varepsilon]$.

 Then  
the hypersurfaces $Y_t=Y\times \{t\} \subset X$  are umbilic  with  the principal  curvatures of $Y_t$  equal to  $\alpha^\ast_i(t)= \frac {\varphi'(t)}{\varphi(t)}$, $i=1,...,n-1$ for
 $$\mbox {$A_t^\ast = \frac {\varphi'(t)}{\varphi(t)}h_t$ for $\varphi'=\frac {d\varphi(t)}{dt}$ and $A_t$ being multiplication by  $\frac {\varphi'}{\varphi}$ }.$$

  The  Weyl formula   reads in this case    as follows.
   $$(n-1)\left (\frac {\varphi'}{\varphi}\right)'  =-  (n-1)^2\left (\frac {\varphi'}{\varphi}\right)^2-\frac {1}{2}\left (Sc(g) - Sc(h_t) -(n-1)(n-2)\left (\frac {\varphi'}{\varphi}\right)^2\right).  $$
   Therefore,

  $$Sc(g)=\frac {1}{\varphi^2}Sc(h)-2(n-1)\left (\frac {\varphi'}{\varphi}\right)' -n(n-1)\left(\frac {\varphi'}{\varphi}\right)^2=$$
   $$=\frac {1}{\varphi^2}Sc(h)-2(n-1)\frac {\varphi''}{\varphi}-  (n-1)(n-2)\left (\frac {\varphi'}{\varphi}\right)^2,\leqno {\mbox {(\Large ${\color {blue}\star}$})}$$
  where, recall, $n=dim(X)=dim(Y)+1$  and the mean curvature of $Y_t$ is 
  $$mean.curv(Y_t\subset X)=(n-1)\frac {\varphi'(t)}{\varphi(t)}.$$

{\it Examples.} (a) If  $Y=(Y,h)=S^{n-1}$ is the unit sphere, then

$$Sc_g= \frac {(n-1)(n-2)}{\varphi^2}  -2(n-1)\frac {\varphi''}{\varphi} -(n-1)(n-2)\left(\frac {\varphi'}{\varphi}\right)^2,
$$
which  for $\varphi = t^2$ makes the expected  $Sc(g)=0,$
since $g=dt^2+t^2 h$, $t\geq 0$,  is the Euclidean metric  in the  polar coordinates.

If  $g=dt^2+\sin t^2 h$, $-\pi/2 \leq t\leq  \pi/2$, then $Sc(g)= n(n-1)$   where this $g$  is the spherical metric on $S^n$.  \vspace {1mm}

(b) If  $h$ is the (flat) Euclidean metric on $\mathbb R^{n-1}$ and $\varphi=\exp t$, then 
 $$Sc(g) =-n(n-1)=Sc(\mathbf H^n_{-1}).$$

(c) What is slightly less obvious, is that if  
$$\varphi(t) =\exp \int_{-\pi/n}^t -\tan \frac {nt}{2} dt, \mbox { }  -\frac {\pi}{n}<t < \frac {\pi}{n},$$
  then the scalar curvature of the metric $\varphi^2h+dt^2$, where $h$ is flat, is {\it constant positive}, namely  $Sc(g)=n(n-1)=Sc(S^n)$,
 by elementary calculation\footnote{See \S12 in [GL(complete) 1983].}

\vspace {1mm}

{\it Cylindrical  Extension Exercise.} Let  $Y$ be a smooth manifold,  $X=Y\times\mathbb R_+$, let  $g_0$ be  a Riemannian metric  in a neighbourhood  of the boundary $Y=Y\times \{0\}=\partial X$,
 let $h$ denote the Riemannian metric in $Y$ induced from $g_0$  and let $Y$ has {\it constant mean curvature} in $X$ with respect to $g_0$. 

Let $X'$ be a (convex if you wish)  ball in the standard (i.e  complete simply connected) space   with  constant sectional curvature  and of the same dimension $n$ as $X$, let $Y'=\partial X'$ be its boundary sphere,  let, let $Sc(h)>0$ and let  the  mean and the scalar  curvatures of  $Y$ and $Y'$ are related by the following (comparison) inequality.
$$\frac {|mean.curv_{g_0}(Y)|^2}  {Sc(h, y)}< \frac {|mean.curv(Y')|^2}  {Sc(Y')} 
\mbox { for all } y\in Y.\leqno {[<]}$$

  Show that \vspace {1mm}

 {\sl  if $Y$ is compact, there exists 
  a    smooth positive function  $\varphi(t)$, $0\leq t<\infty$,   which is constant at infinity and  such that the 
 the warped  product metric $g=  \varphi^2h+dt^2$ has

 the same {\sf Bartnik data} as $g_0$, i.e.

 $$\mbox {$g|Y=h_0$ and
  $mean.curv_g(Y)= mean.curv_{g_0}(Y),$}$$}

Then show that

{\sf  one {\it can't make} $Sc(g)\geq Sc(X')$ in general,  if  [<] is relaxed to the corresponding  {\it non-strict} inequality, where an example is provided by 
 the Bartnik data of $Y'\in X'$ itself}.\footnote { It follows from [Brendle-Marques(balls in $S^n$)N 2011]  that the the cylinder $S^{n-1}\times \mathbb R_+$  admits a complete Riemannian metric $g$ cylindrical at infinity which has   $Sc(g)> n(n-1)$,  and which has the same Bartnik data as 
the boundary sphere $X'_0$ in the hemisphere $X'$ in the unit $n$-sphere. But the non-deformation result from [Brendle-Marques(balls in $S^n$) 2011], suggests that this might be impossible for the Bartnik data of  {\it small} balls in the round sphere.}  

  {\it Vague Question.} {\sf What are "simple natural" Riemannian metrics $g$ on $X=Y\times\mathbb R_+$
  with given Bartnik data $(Sc(Y), mean,curv(Y))$, where      $Y\subset X$ is allowed  {\it variable} mean  curvature, and what  are possibilities for  lower bound on the     scalar curvatures of such $g$ granted   
     $|mean.curv(Y, y)|^2/  Sc(Y, y)< C$, e..g. for $C=|mean.curv(Y')|^2/Sc(Y')$  for $Y'$ being a sphere in a space of constant curvature.}

\vspace {1mm}

%%%%%%%%%%%%%%%%%%

\subsubsection {\color {blue} Higher Warped Products} \label {warped+2}

%%%%%%%%%%%%%%%%%%%

Let  $Y$ and $S$ be Riemannian manifolds with the metrics denoted  $dy^2$ (which now play the role  of the above $dt^2$)  and   $ds^2$ (instead of $h$), 
let  $\varphi>0$ be a smooth function on $Y$, and let 
 $$g=  \varphi ^2 (y)ds^2+dy^2$$ 
be the corresponding  warped metric on $Y\times S$, 

 Then
 $$Sc(g)(y,s) =Sc(Y)(y)+\frac {1}{ \varphi(y)^2} Sc(S)(s)  -\frac  {m(m-1)}{ \varphi^2(y )}
 || \nabla  \varphi(y)||^2-\frac {2m}{\varphi(y)}\Delta  \varphi(y),\leqno {\mbox {(\Large ${\color {blue}\star\star}$})}$$
  where $m=dim (S)$ and $\Delta= \sum \nabla_{i,i}$ is the Laplace  on 
  $Y$.
  
  To prove this, apply   the above c ({\Large {\color {blue} $\star$}}) to $l\times S$ for   naturally parametrised geodesics   $l\subset Y$ passing trough $y$ and then average over the space of these $l$,  that is the unit tangent sphere of $Y$ at $y$. \vspace{1mm}

The most relevant example here  is 
where $S$  is the real line  $\mathbb R$ or the circle $S^1$ also denoted $\mathbb T^1$ and where  ({\Large {\color {blue} $\star$}}) reduces to 
 $$Sc(g)(y,s) =Sc(Y)(y)-\frac {2}{\varphi}\Delta  \varphi(y). \footnote{The roles of $Y$ and $S=\mathbb R$   and notationally  reversed here with respect   to those in    ({\Large {\color {blue}$\star$})}} \leqno {\mbox {(\Large ${\color {blue}\star\star}$})_1} $$

For instance, if the  $L=-\Delta +\frac {1}{2}Sc$ on $Y$ is  strictly positive, that is the lowest eigenvalue $\lambda$ is strictly positive  and if  $\varphi$ equals to   the 
corresponding  eigenfunction of $L$,  
then 
$$-\Delta  \varphi =\lambda\cdot\varphi-\frac {1}{2}Sc\cdot \varphi$$
and  
$$Sc(g)=2\lambda>0,$$

The  basic feature of the metrics     $\varphi ^2 (y)ds^2+dy^2$ on $Y\times \mathbb R$ is that they are  
 {\it $\mathbb R$-invariant}, where  the quotients   $(Y\times \mathbb R)/\mathbb Z= Y\times 
 \mathbb T^1 $ carry   the corresponding {\it $\mathbb T^1 $-invariant} metrics, while {\it the $\mathbb R$-quotients are isometric to $Y$.}
 
Besides $\mathbb R$-invariance,  a characteristic feature of warped product metrics  is {\it integrability} of  the tangent  hyperplane  field  normal  to the $\mathbb R$-orbits,  where 
$Y\times \{0\}\subset Y\times \mathbb R$, being normal to these orbits,  serves as an  integral variety  for this field. 

Also notice that $Y=Y\times \{0\}\subset Y\times \mathbb R$ is totally geodesic with respect to 
the metric  $\varphi ^2 (y)ds^2+dy^2$, while the ($\mathbb R$-invariant)   {\it curvature} (vector field) {\it of the  $\mathbb R$-orbits} is equal to the {\it gradient field $\nabla \varphi$}  extended from  $Y$ to  $Y\times \mathbb R$.
coordinates

In what follows, we emphasize  $\mathbb R$-invariance and   interchangeably speak of  $\mathbb R$-invariant metrics on $Y\times \mathbb R$ 
 and metrics warped  with factors $\varphi^2$ over $Y$.

\vspace {1mm}  
{\it Gauss-Bonnet  $g^\rtimes$-Exercise.}   Let the above $S$ be the Euclidean space $\mathbb R^N$ (make it $\mathbb T^n$ if you wish to keep compactness) with coordinates $t_1,...,t_N$,   let
  $$\Phi(y)=(\varphi_1(y),...,\varphi_i(y),..., \varphi_N(y))$$ 
be an $N$-tuple of smooth positive function on a Riemannian mnanifold  $Y=(Y,g)$ 
and define the  (iterated t warped product)  metric $g^\rtimes =g^\rtimes_\Phi$  on $Y\times S$  as follows:
$$g^\rtimes =g(y)+\varphi^2_1(y)dt_1^2+\varphi^2_2(y)dt_2^2+...+\varphi^2_N(y)dt_N^2$$

Show that the scalar curvature of this metric, which,  being   $\mathbb R^N$-invariant,  is regarded as a function on $Y$, satisfies:
$$Sc(g^\rtimes,y)=Sc(g) -2\sum_{i=1}^N\Delta_g \log \varphi_i  - \sum_{i=1}^N(\nabla_g\log \varphi_i)^2
-\left(
\sum_{i=1}^N \nabla_g\log \varphi_i\right)^2,
$$ 
 thus 
$$\int_Y Sc(g^\rtimes,y)dy\leq \int_Y Sc(g,y)dy,$$
and, following  [Zhu(rigidity)  2019], obtain the following\vspace {1mm}

 {\it "Warped"  Gauss-Bonnet Inequality  for Closed Surfaces $Y$}:
$$ \int_Y Sc(g^\rtimes,y)dy\leq 4\pi\chi(Y)$$
{\sf for the  (iterated) warped product metrics $g^\rtimes=g^\rtimes_{\phi}$ for all positive $N$-tuples of  $\Phi$  of positive functions on $Y$.
\footnote {See [Zhu()  2019] and sections \ref{Gauss-Bonnet5}, 
\ref{Gauss-Bonnet7}  for applications and generalizations.}}
%%%%%%%%%%%%%%%%%%%

\subsection {\color{blue}   Second Variation Formula}\label {2nd variation2}
%%%%%%%%%%%%%%%%%%%%%%%%%%%

The Weyl  formula also yields  the following formula for the {\it second derivative} of the $(n-1)$-volume  of a cooriented hypersurface $Y\subset X$ under a normal deformation of $Y$ in $X$, where the scalar curvature of $X$ plays an essential role.

The deformations we have in mind are by vector fields directed by geodesic normal to $Y$, where in the simplest case the norm of his field equals one.

In this case we have an equidistant motion  $Y\mapsto Y_t$  as earlier and the second derivative 
of $vol_{n-1}(Y_t)$, denoted here $Vol=Vol_t$,  is expressed in terms of of the shape   $A_t=A(Y_t)$ of $Y_t$ and the Ricci curvature of $X$, where, recall $trace (A_t)=mean.curv(Y_t)$ and 
$$\partial_\nu Vol=\int_{Y} mean.curv(Y)dy$$
 by the first variation formula.
 
Then, by Leibniz' rule, 
$$\partial^2_\nu Vol=\partial_\nu\int_{Y}trace (A(y))dy =\int_{Y}trace^2 (A(y))dy +\int_{Y}trace (\partial_\nu A(y))dy,$$
and where, by Weyl's formula, 
$$trace (\partial_\nu A)  = -trace  (A^2) -Ricci (\nu, \nu)$$ 
for  the  normal unit  field $\nu$.

Thus, 
$$\partial^2_\nu Vol=\int_Y (mean.curv)^2- trace (A^2)-Ricci(\nu,\nu),$$
which,  combining this with the above expression 
$$Ricci(\nu)=\frac{1}{2} \left(Sc(X)-Sc(Y) -(mean.curv(Y))^2+||A^\ast||^2\right),$$
shows that
$$\partial^2 _\nu Vol=\int \frac {1}{2}\left ( Sc(Y)-Sc(X)  + mean.curv^2 - ||A^\ast||^2 \right ). $$

In particular, if $Sc(X)\geq 0$ and $Y$ is minimal, then, 
$$\int _YSc(Y,y)dy\geq 2\partial^2_\nu Vol \leqno  {\rm\color  {blue}({\int \hspace {-0.7mm}Sc\geq 2\partial^2 Vol})}$$
(compare with  the {\color{blue} $(Sc\geq -2Ric)$} in 2.2).\vspace{1mm}

{\it \color {red!80!black}Warning.}   Unless $Y$ is minimal  and despite the notation $\partial^2_\nu$, this derivative depends on how the normal filed  on $Y\subset X$ is extended to a vector filed on (a neighbourhood of $Y$ in) $X$.\hspace {1mm}

{\it Illuminative Exercise.} Check up this formula for  concentric spheres of radii $t$ in the spaces with constant sectional curvatures that are   $S^n$, $ \mathbb R^n$ and $\mathbf H^n$.

\vspace{1mm}

Now, let us  allow a non-constant geodesic field normal to $Y$, call it  $\psi\nu$, where $\psi(y)$ is a smooth function on $Y$ and write down the full second variation formula as follows:
$$\partial^2_{\psi\nu} vol_{n-1}(Y)=\int_Y ||d\psi(y)||^2dy +R(y) \psi^2(y)dy$$
for 
$$R(y)=\frac {1}{2}\left ( Sc(Y,y)-Sc(X,y)  +M^2(y) - ||A^\ast(Y)||^2 \right ),\leqno  \mbox { {\Large \color {blue}$[ \circ\circ]$}}$$
where $M(y)$ stands for the mean curvature of $Y$ at $y\in Y$ and $||A^\ast(Y)||^2=\sum_i(\alpha^\ast)^2$, $i=1,...,n-1.$

Notice, that the "new" term $\int_Y ||d\psi(y)||^2dy$   depends only on the normal field itself, while the $R$-term  depends on the extension of $\psi\nu$ to $X$, unless 

\hspace{24mm} {\it $Y$ is minimal, where   {\Large \color {blue}$[ \circ\circ]$} reduces to }
$$\partial^2_{\psi\nu} vol_{n-1}(Y)=\int_Y ||d\psi||^2 +\frac {1}{2}\left ( Sc(Y)-Sc(X)  - ||A^\ast||^2 \right ) \psi^2.\leqno  \mbox { {\Large \color {blue}$[ \ast\ast]$}}$$

Furthermore, if  $Y$ is volume minimizing  in its neighbourhood, then $\partial^2_{\psi\nu} vol_{n-1}(Y)\geq 0$; therefore,
$$\int_Y( ||d\psi||^2 +\frac {1}{2} ( Sc(Y)) \psi^2\geq\frac {1}{2}\int_Y( Sc(X,y)  + ||A^\ast(Y)||^2)\psi^2dy \leqno  \mbox { {\Large \color {blue}$[ \star\star]$}}$$
for all non-zero functions $\psi=\psi(y)$.

Then, if we recall that 
$$\int_Y ||d\psi||^2dy =\int_Y\langle-\Delta\psi, \psi \rangle dy,$$
we will see that  {\Large \color {blue}$[ \star\star]$} says that \vspace{1mm}

{\it \hspace {-4mm} {\color {blue} the  }  $\psi\mapsto -\Delta \psi +\frac {1}{2} Sc(Y)\psi  $  {\color {blue} is greater than}\footnote{$A\geq B$  for selfadjoint operators signifies that $A-B$ is positive semidefinite.} $\psi\mapsto \frac {1}{2}(Sc(X,y)  + ||A^\ast(Y)||^2)\psi$.}\vspace{1mm}

Consequently, \vspace {1mm} 

{\it \hspace {9mm}{\it \color {blue} if $Sc(X)> 0$, then 
the    $ -\Delta +\frac {1}{2} Sc(Y) $ on $Y$ is positive.}

\vspace{1mm}

{\it Justification of the $||d\psi||^2$ Term.} }Let $X=Y\times \mathbb R$ with the product metric and 
let $Y=Y_0=Y\times \{0\}$ and $Y_{\varepsilon \psi}\subset X$ be the graph of the function $\varepsilon\psi$ on $Y$.
Then 
$$vol_{n-1}(Y_{\varepsilon \psi})=\int_Y\sqrt {1+\varepsilon^2||d \psi||^2}dy=
vol_{n-1}(Y)+\frac{1}{2}\int_Y\varepsilon^2||d \psi||^2 +o(\varepsilon^2)$$
by the Pythagorean theorem 

and 

$$\frac  {d^2 vol_{n-1}(Y_{\varepsilon \psi})}{d^2\varepsilon}=||d \psi||^2+o(1).$$
by the binomial formula.

This proves   {\Large \color {blue}$[ \circ\circ]$}  for product manifolds and the general case follows by  
{\it linearity/naturality/functoriality} of   the formula {\Large \color {blue}$[ \circ\circ]$}.
\vspace {1mm}

{  \large \color {magenta}\it Naturality Problem.} All "true formulas"  in the Riemannian geometry should be derived with minimal, if any,   amount of calculation --  only on the basis of their "naturality" and/or of their validity in  simple examples, where these formulas are obvious.

Unfortunately,  this "naturality principle" is absent from the textbooks on differential geometry, but,  I guess,
it may be found in some algebraic articles (books?). 

{\it Exercise.} Derive  the second main formula {\color{blue} 2.3.A} by pure thought  from its manifestations in the examples in  the above {\it illuminative exercise}.\footnote {I haven't myself solved this exercise.}

%%%%%%%%%%%%%%%%%

\subsection  {\color {blue}  Conformal Laplacian and the Scalar Curvature of Conformally and non-Conformally  Scaled Riemannian Metrics}\label {conformal2}   
%%%%%%%%%%%%%%%%%%%%
 Let $(X_0,g_0)$ be a compact Riemannian manifold of dimension $n\geq 3$ and let  $\varphi=\varphi(x)$ be a smooth positive function on $X$.
   
   Then, by a straightforward calculation,\footnote{There must be a better argument.}
$$Sc(\varphi^2g_0)= \gamma_n^{-1}\varphi ^{-\frac{n+2}{2}}  L(\varphi^\frac{n-2}{2}), \leqno {\LEFTcircle}$$
where $L$ is the {\it \color{blue} conformal  Laplace } on $(X_0,g_0)$
$$ L(f(x))= -\Delta f(x)+ \gamma_nSc(g_0,x)f(x) $$
for the ordinary Laplace (Beltrami) $\Delta f=\Delta_{g_0}f=\sum_i\partial_{ii}f$
and  $\gamma_n= \frac{n-2}{4(n-1)}$.\vspace{1mm}

Thus,  we conclude to the following.  \vspace{1mm}

 \textbf {Kazdan-Warner Conformal Change Theorem.} \footnote {[Kazdan-Warner(conformal) 1975]: {\sl  Scalar curvature and conformal deformation of Riemannian structure}.} {\sf Let  $X=(X, g_0)$ be a  closed Riemannian manifold,  such the the conformal  Laplace  $L$ is positive. }
 
 {\it Then $X$ admits a Riemannian metric $g$ (conformal to $g_0$) for which $Sc(g)>0$.}\vspace{1mm}

{\it Proof.}  Since $L$ is positive, its first eigenfunction, say $f(x)$ is positive\footnote{We explain
this in section \ref{eigenfunctions2}. }
and since $L(f)=\lambda f, $  $\lambda >0$,  

$$Sc\left ( f^\frac{4}{n-2}g_0\right)= \gamma_n^{-1}  L(f) f^{-\frac{n+2}{n-2}}=\gamma_n^{-1}f^{\frac {2n}{n-2}}>0.$$

\vspace {1mm}%%%%%%%%%%%%%&&&&&&&&&

{\it \textbf {Example}}: {\color {blue}\large  Schwarzschild metric.}  If $(X_0,g_0)$   is the Euclidean $3$-space,   and  $f=f(x)$  is positive  function, then 

\hspace{20mm}{\it  the sign of 
 $Sc(f^4g_0)$ is equal  to  that of $-\Delta f$.}\vspace {1mm}
 
   In particular, since the function $\frac{1}{r}= (x_1^2+x_2^2+x_3^2)^{-\frac {1}{2}},$ is harmonic,
{ \it  the  Schwarzschild metric {\color {black}  \hspace{0.5mm} $g_{Sw}=\left (1+\frac {m}{2r}\right)^4g_{0}$} \hspace{0.5mm} has {\it zero} scalar curvature.}\vspace {1mm}

{\sf If $m>0$, then this metric is defined for all $r>0$ and it  is invariant under the involution
$r\mapsto \frac {m^2}{r}$. 

If $m=0$, this the flat Euclidian metric. 

If $m<0$, then  this metric is defined only for $r>m$ with a singularity ar $r=m$.}
\vspace{1mm}

\vspace {1mm}

{\sf \large  Non-Conformal Scaling.}  Let $X=(X,g)$ be a smooth $n$-manifold, and let $\mathbb R^\times_x\subset GL_x(n)$, $x\in X$, be a smooth family of diagnosable (semisimple)  1-parameter subgroups in the linear groups $GL_x(n)=GL_n$ that act in the tangent spaces $T_x(X)$.

Then the the multiplicative  group of functions $\phi :X\to \mathbb R^\times $ acts on the tangent bundle  $T(X)$ by 
$$\tau\mapsto  = \phi(x)(\tau)\mbox {  for } \phi(x)\in  \mathbb R^\times = \mathbb R_x^\times \subset GL_x=GL(T_x(X))$$ 
 and, thus on the space of Riemannin metrics $g$ on $X$.

The main instance of  such an action is where  the tangent bundle is orthogonally split,  $T(X)=T_1\oplus T_2$, and $\phi$ acts by scaling on the subbundle  $T_2$.

It is an not hard to write down a formula for the scalar curvature of $g_1+\phi^2g_2$,
but it is unclear  what, in general, would be a workable criterion for solvability of  the inequality 
 $Sc(g_\varphi)>0$ in $\varphi$, e.g. in the case where $X=X_1\times X_2$ and the subbundles $T_1 $ and $T_2$ are equal to the  tangent  bundles of  submanifolds  $X_1 \times x_2\subset X$, $x_2\in X_2$, and $x_1\times X_2\subset X$, $x_1\in X_1$.

Yet, in the case of $rank(T_2)=1$, this equation introduced, I believe, by Robert Bartnik in  [Bartnik(prescribed scalar) 1993]  was successfully applied to extension  of metrics with $Sc>0$ (see section \ref {fill-ins3})\footnote{Other special cases of this are (implicitly) present in the geometry of  Riemannin  warped product, in the process  of {\it smoothing corners with} $Sc\geq \sigma$  and in the  {\it transversal  blow up} of foliations  with $Sc>0$. }

\subsection {\color {blue}  Schoen-Yau's    Non-Existence Results   for  $Sc>0$ on SYS Manifolds    via    Minimal (Hyper)Surfaces  and  Quasisymplectic $[Sc\ngtr 0]-$Theorem} \label {SY+symplectic2}
{\sf Let   $X$ be  a  three dimensional  Riemannian manifold with $Sc(X)>0$ and $Y\subset X$  be  an orientable   cooriented surface with minimal area in its integer homology class. 

Then the inequality {\color  {blue}$({\int \hspace {-0.7mm}Sc\geq 2\partial^2 V}$)} from section \ref {2nd variation2}, which says in the present case that 
$$\int _YSc(Y,y)dy> 2\partial^2_\nu area(Y),  $$%%%%??
implies that

 \hspace {20mm} {\it $Y$ must be a topological sphere}}. \vspace {1mm}

In fact, minimality of $Y$ makes $\partial^2_\nu area(Y)\geq 0 $, hence $\int _YSc(Y,y)dy>0$, and the sphericity of $Y$ follows by the Gauss-Bonnet theorem.

And since all integer  homology classes in closed orientable Riemannian  $3$-manifolds admit area minimizing 
representatives by  the geometric measure theory developed by Federer, Fleming and  Almgren, we arrive at the following conclusion. \vspace {1mm}

 {\color {blue} $\bigstar_3$}   \textbf{Schoen-Yau $3d$-Theorem}.  {\it All  integer  2D homology classes in   closed Riemannian $3$-manifolds  with $Sc>0$
 are  spherical.
  
 For instance, the   3-torus  admits no metric with $Sc>0$.}

\vspace {1mm}

The above argument appears in  Schoen-Yau's  15-page paper   [SY(incompressible)  1979], most of which is occupied by    an independent proof of the existence and regularity  of minimal $Y$. 

In fact, the existence  of minimal surfaces  and their regularity needed for the above argument has been known since late (early?) 60s\footnote{Regularity of volume minimizing hypersurfaces  in manifolds $X$ of dimension $n\leq 7$, as we mentioned earlier, was proved  by  Herbert Federer  in [Fed(singular) 1970], by reducing the general case of the problem to that    of minimal cones  resolved by   Jim Simons in  [Simons(minimal) 1968].} but, what was, probably, missing prior to the  Schoen-Yau  paper  was  the innocuously  looking corollary of Gauss' formula  in 2.2,  
$$Sc(Y)= Sc(X) +(mean.curv(Y))^2- ||A^\ast(Y)||^2-Ricci(\nu,\nu)$$
and the issuing inequality 
$$ Sc(Y)> -2Ricci(\nu, \nu)$$
for minimal $Y$ in manifolds $X$ with $Sc(X)> 0$.

For example, Burago and Toponogov,  come  close to the above argument, where, they bound from below the injectivity radius of Riemannian  $3$-manifolds $X$ with $sect.curv(X)\leq 1$ and $Ricci (X)\geq \rho>0$ by
$$ inj.rad(X)\geq  6e^{-\frac {6}{\rho}}, $$ 
where this is  done by carefully  analysing minimal surfaces $Y\subset X$  bounded by, a priori very short, closed geodesics in $X$,
and where an essential step in the proof is the lower bound on the first eigenvalue of the Laplace  on  $Y$ by  $\sqrt {Ricci(X)}.$\footnote{[BurTop(curvature bounded  above)1973],{\sl On 3-dimensional Riemannian spaces with curvature bounded
above}.}

 \vspace {1mm}
 
 {\it Area Exercises.} Let $X$ be homeomorphic to $Y\times S^1$, where $Y$ is a closed orientable surface  with the Euler number $\chi$. 
 
 (a) Let $\chi>0$, $Sc(X)\geq 2$  and show that there exists a   surface  $Y_{o}\subset X$  homologous to $Y\times \{s_0\}$, such that  
 $ area(Y_{o})\leq 4\pi$.\footnote{See [Zhu(rigidity)  2019] for a higher dimensional version of this inequality.}

(b) Let  $\chi<0$, $Sc(X)\geq -2$ and show that all      surfaces $Y_\ast\in  X$  homologous to $Y\times \{s_0\}$ have 
$area(Y_\ast) \geq -2\pi \chi.$

(c) Show that (a) remains valid  for complete manifolds $X$ homeomorphic to $Y\times \mathbb R$.\footnote{I haven't solved this exercise.}

\vspace {1mm}%%%%%%%%%%%%%%&&&&&&&&&&&
 
{\color {blue} $\bigstar^{codim1}$} \textbf {Schoen-Yau Codimension 1 Descent Theorem,} [SY(structure) 1979].
 {\sf Let $X$ be a compact orientable   $n$-manifold with 
$Sc>0$.}
 
 {\sl If $n\leq 7$, then all  integer homology classes $h\in H_{n-1}(X)$ are representable by   compact  oriented  $(n-1)$-submanifolds $Y$ in $X$, which  admit  metrics with $Sc>0$.}

{\it Proof.} Let $ Y$ be a volume minimizing hypersurface representing $h$, the existence and regularity of which is guaranteed by a Federer 1970-theorem\footnote{[Federer(singular) 1970]:  {\sl The singular sets of area minimizing rectifiable currents with
codimension one and of area minimizing flat chains modulo two with arbitrary
codimension}.}  and recall that by 
  {\Large \color {blue}$[ \star\star]$} in \ref{2nd variation2} the 
  $-\Delta  +\frac {1}{2} Sc(Y)  $ is positive.
Hence, the conformal Laplace    $-\Delta +\gamma_nSc(Y)$
is also positive for
$\gamma_n= \frac{n-2}{4n-1}\leq \frac {1}{2}$ and the proof follows by Kazdan-Warner conformal change theorem.
\vspace{1mm}

 %%%%%%%%%%%&&&&&&&&&&
 {\color {blue} $\bigstar_{\mathbb T^n}$} \textbf{  Mapping to the Torus Corollary. } {\sf If a closed orientable $n$-manifold $X$ admits a map to the torus $\mathbb T^n$  with {\it non-zero degree}, then $X$ admits {\it no metric with $Sc>0.$}}\vspace {1mm}
 
 Indeed, if a  closed submanifold   $Y^{n-1}$ is  {\it non-homologous to zero} in this  $X$ then it  (obviously)  admits a map to $\mathbb T^{n-1}$  with non-zero degree.
Thus,  the above  allows an   inductive reduction of  the problem to the  case of $n=2$, where the Gauss-Bonnet theorem applies.\vspace {1mm}
  
  {\it \textbf  {SYS-Manifolds}.} Schoen and Yau say  in  [SY(structure) 1979] that their codimension 1 descent theorem 
 delivers a  topological  obstruction  to $Sc>0$ on a class of manifolds,  which is, even in the spin case, \footnote{ A smooth connected   $n$-manifolds $X$ is {\it spin} if the frame bundle over $X$ admits a double cover extending the natural double cover of a fiber,  where such a fiber is equal to  
  the linear group,  (each of the two connected components of) which   admits a  a unique   non-trivial double cover  $\tilde GL(n)\to GL(n)$.} is  not covered  by the twisted Dirac operators methods.
  
  This claim  was confirmed by Thomas Schick, who  defined, in homotopy theoretic terms, integer  homology classes in aspherical spaces, say $h\in H_n (\underline X)$  and who proved using the codimension one descent 
   theorem that these $h$ for $n\leq 7$  can't be dominated by compact orientable $n$-manifolds with $Sc>0$.
  
In  more geometric  terms, the  $n$-manifolds $X$,  to which   Schick's argument applies, we call them {\it Schoen-Yau-Schick},
can be described d as follows.

{\sf A closed orientable $n$-manifold  is {\it Schoen-Yau-Schick}  if it admits a smooth map $f:X\to\mathbb T^{n-2}$, such that the homology class of the  pullback of a generic point,  
  $$h=[f^{-1}(t)]\in H_2(X)$$  
is {\it  non-spherical,} i.e.  it is not in the image  of the {\it Hurewicz homomorphism} $\pi_2(X)\to H_2(X)$.  }

Then Schick's  corollary to  Schoen-Yau's theorem reads. \vspace {1mm}
  
  %%%%%%%%%%%%%%%&&&&&&&&&&&&&
  
  {\color {blue} $\bigstar_{SYS}$} \textbf {Non-existence Theorem for SYS Manifolds.}  {\it  Schoen-Yau-Schick manifolds of dimensions $n\leq 7$ admit no metrics with $Sc>0$.} \vspace {1mm}

(b)  {\it Exercises.}  (b$_1$)   Construct  examples of SYS manifolds of dimension $n\geq 4$, where all maps $X\to \mathbb T^n$ have zero degrees.

{\it Hint}: apply surgery to  $\mathbb T^n$.

 \vspace {1mm}

(b$_
2$) Show that if the first  homology group $H_1(X)$ of a SYS-manifold has no torsion, then a finite covering of $X$ admits a map with  degree one to the torus  $\mathbb T^n$.

 \vspace {1mm}

 (c) The limitation $n\leq 7$ of  the above argument is due a  presence of  singularities of minimal subvarieties in $X$ for $dim(X) \geq 8$. 
 
 If $n=8$, these singularities were proven    to be unstable by Nathan Smale; this improves $n\leq 7$ to $n\leq 8$ in  {\color {blue} $\bigstar_{SYS}$}
  
  More recently,  as we mentioned earlier,  the  dimension restriction was  removed for all $n$  by Lohkamp and by Schoen-Yau;  the arguments in both papers are difficult  and I have not mastered them.\footnote {See [Smale(generic regularity) 2003], SY(singularities)  2017],  [Lohkamp(smoothing) 2018] and section  \ref {singularities3}.}
      \vspace {1mm}
  
Although  the Dirac operator arguments don't apply to  SYS-manifolds, they do deliver  topological  obstructions to  
$Sc>0$, which, according to the present state of knowledge,  lie beyond the range of the minimal surface techniques.
Here is an  instance  of this. \vspace {1mm}

%%%%%%%%%%%%%%%%%%&&&&&&&&&&&&

{\Large \color {blue}$\otimes_{\wedge^k\tilde\omega}$}  \textbf{Quasisymplectic Non-Existence Theorem.} 
 {\sf Let $X$ be a compact   {\color {blue}$\otimes_{\wedge^k\tilde\omega}$}-{\it manifold} of dimension $n=2k$,
 i.e. $X$  is orientable and it  carries a {\it closed 2-form} $\omega$  (e.g. a symplectic one), such that  $\int_X\omega^k\neq 0$, and such that  the lift $\tilde \omega$ of $\omega$ to   the universal covering $\tilde X$ is {\it exact}, e.g.  $\tilde X$ is contractible.\footnote{ It's enough to have $\tilde X$ spin.}

 Then {\it $X$ admits no metric with $Sc>0$.}} \vspace {1mm}

This applies, for instance, to  {\it even dimensional tori},  to {\it aspherical $4$-manifolds with} $H^2(X, \mathbb R)\neq 0$ and to {\it products} of such manifolds\footnote{Recently,  Chodosh  and  Li  proved that 

{\it compact aspherical manifolds of dimensions 4 and 5 admit no metrics with positive scalar curvatures.}
(See [Chodosh-Li(bubbles) 2020], [G(aspherical) 2020]  and  section\ref{5D.3})

But this remains problematic for products of pairs of aspherical 4-manifolds.}
but {\it not} to general SYS-manifolds.
\vspace {1mm}

{\it Idea of the Proof.} Assume without loss of generality that $\omega$ serves as the curvature form  of a complex line bundle$L\to X$and let $\tilde L\to \tilde X$ be the lift of $L$ to the universal covering $\tilde X\to X$. 

Since the curvature $\tilde \omega$ of $\tilde L$, is exact the   bundle $\tilde L$is  topologically trivial, hence it can be represented by $k$-th tensorial power of another line bundle, 
$$L=(L^{\frac{1}{k}})^{\otimes k},$$
 where the curvature of $L^{\frac{1}{k}}$ is $\frac{1}{k}\tilde\omega.$
By  Atiyah's {\it $L_2$-index theorem}, there are {\it non-zero harmonic $L_2$-spinors} on $\tilde X$
  {\it twisted with  $L^{\frac{1}{k}}$} for infinitely many $k$, but the  twisted Schroedinger-Lichnerowicz-Weitzenboeck-(Bochner) formula 
 applied to large $k$ doesn't allow such spinors for $Sc(\tilde X\geq \sigma>0$.\footnote {Atiyah's theorem from  [Atiyah(L2) 1976] needs a slight adjustment here, since the action of  the  fundamental group $\Gamma=\pi_1(X)$ on $\tilde X$ doesn't lift to $L^{\frac{1}{k}}$; yet the fundamental group of  the (total space) of the unit  circle bundle    of $L$  does naturally act on   $L^{\frac{1}{k}}$.
 Also, there is no difficulty in   extending Lichnerowicz' vanishing argument to the $L_2$ case, see 
\S$9\frac {1}{8}$ in [G(positive) 1996].}
\vspace {1mm}

{\it Exercise.}  Show that if $X$ is {\color {blue}$\otimes_{\wedge^k\tilde\omega}$},  then the    classifying map $X\to {\sf B} (\Pi)$, where ${\sf B} (\Pi)=K(\Pi,1)$ is  the classifying space for the group $\Pi=\pi_1(X)$, sends the fundamental
homology  class $[X]$ 
 to a {\it non-torsion} class in $H_n({\sf B}(\Pi)).$ 

  \vspace {1mm}
  
  {\sf \large Problem.} {\sl Is there a unified approach that would apply to  $SYS$-manifolds and 
  to the above {\color {blue}$\otimes_{\wedge^k\tilde\omega}$}-manifolds   $X$, e.g. symplectic ones with contractible universal coverings?}

For instance,  

{\sl do products of  $SYS$ and   {\color {blue}$\otimes_{\wedge^k\tilde\omega}$}-manifolds ever 
carry metrics with positive 

scalar curvatures?}

%%%%%%%%%%%%%%%%%%
\subsection{\color {blue} Warped $\mathbb T^\rtimes$-Stabilization and Sc-Normalization}\label{warped stabilization and Sc-normalization2} 
%%%%%%%%%%%%%%%%%%%%%%
Many geometric properties of  Riemannian manifolds $X=(X,g)$ implied by the inequality $Sc(g)\geq \sigma$  follow (possibly in a weaker form) from the same inequality for a larger manifold, say  $X^\ast$, that,  topologically,    is  the product of $X$ with the  a torus, $X^\ast =X\times \mathbb T^N$   for some $N=1,2,...$, where the Riemannian metric $g^\ast$ on $X^\ast$ is invariant under the action of  $ \mathbb T^N$ and where $X^\ast/\mathbb T^N$ is isometric to $X$.
\vspace{1mm}

{\it Surface Examples. }  Let  $X=(X,g)$ be  a closed
 surface and  $g^\ast$  be a  $\mathbb T^N$-invariant metric on $X\times \mathbb T^N$, such that
  $$(X\times \mathbb T^N, g^\ast)/\mathbb T^N=(X,g).$$
 
%%%%%%%%%%%%%%%%%&&&&&&&&&&&&

(a) {\it \textbf {Sharp Equivariant    Area  Inequality}.} {\sf If  $Sc(g^\ast)\geq \sigma>0$, then a special  case a theorem  by Jintian Zhu,\footnote {See [Zhu(rigidity) 2019] and     \ref {Gauss-Bonnet5}, \ref {Gauss-Bonnet7} for related inequalities.} says that} \vspace{1mm}

 {\it the area of $X$ is bounded  the same  way as it is for $Sc(g)\geq \sigma$,  
$$area(X)\leq \frac {8\pi}{\sigma}.$$}
Moreover,  

\hspace {10mm} {\it the equality holds only if  $X^\ast$ is the  {\sf isometric} product $X\times \mathbb T^N$.}

\vspace{1mm}
%%%%%%%%%%%%%%%%%&&&&&&&&&&&&
(b) ({\sf Weakened}) {\it \textbf{$\mathbb T^\ast$-Stable 2d Bonnet-Myers Diameter Inequality}}. {\sf If $Sc(g^\ast)\geq \sigma$, then
 $$diam (X)\leq  2\pi\sqrt{\frac { N+1}{ (N+2)\sigma}} \hspace {1mm}< 
 \frac {2\pi}{\sqrt\sigma}.  \leqno {[\color {blue} BMD}]$$} \vspace {1mm}

{\it Proof.}  Given two points $x_1,x_2 \in X$, take two small $\varepsilon$-circles $Y_{-1}$  and $Y_{+1}$  around them, let $X_\varepsilon \subset X$ be the band between them and 
 and apply (the relatively elementary $\mathbb T^N$-invariant  case of)  the {\it  $\frac {2\pi}{n}$-Inequality} from  section \ref {bands3}.\footnote{Also see \S2 in [G(inequalities) 2018]   and the proof of theorem 10.2 in   [GL(complete)  1983].}

\vspace {1mm}

{\it Non Trivial Torus Bundles.} The inequality  {[\color {blue} BMD}] is valid for (all) Riemannian $(N+2)$-manifolds   $X^\ast$ with free isometric  $\mathbb T^N$-actions: \vspace {1mm}

{\it if $Sc(X^\ast)\geq \sigma>0$, then $diam(X^\ast/\mathbb T^N)\leq  2\pi\sqrt{(N+1)/(N+2)\sigma}$.}\vspace {1mm}

In fact, the above proof applies, since, {\it topologically}, the part of $X^\ast$  that lies  over the  band $X_\varepsilon\subset X$ {\it is } the product, $X_\varepsilon\times \mathbb T^N$.\vspace {1mm}

It is {\it \color {blue}  unclear}, however, if the areas of  $X^\ast/\mathbb T^N$ are bounded in terms of $Sc(X^\ast)$    for all such $X^\ast$. 

And, as we shall see later,  possible non-triviality of  torus bundles  create complications   for other 
problems with  scalar curvature.\vspace {1mm}

{\it General Question.} The above  examples suggests that  quotients $X$  of manifolds $X^\ast$  with  $Sc(X^\ast)\geq \sigma$ under free isometric  actions of tori have  {\it  similar } geometric properties to those of manifolds which have $Sc\geq \sigma$ themselves.
But it is unclear how far this similarity goes.\vspace {1mm}

{\it Example.} let $X$ be a closed surface and  $X^\rtimes=X\rtimes \mathbb T^1$ be a warped product as described below.

{\sf Does the inequality $Sc(X^\rtimes)\geq 2$ yield an upper bound on {\it all of geometry} of 
$X$?} 

For instance,  

\hspace {10mm}{\sf is there  a  bound on the number of  unit discs needed to cover $X$?}

\vspace {1 mm}

(If   $Sc (X)\geq 2$, then $X$ admits a {\it distance decreasing} homeomorphism from the unit sphere $S^2$, that
  can be constructed  using  the  family of boundary curves of concentric discs 
with center at some point in $X$.)

\vspace {1mm}

{\it Warped Products.}  As far as geometric applications are concerned, the relevant $X^\ast$ are (iterated) {\it warped products}, we denote  them $X^\rtimes $ and  call {\it warped $\mathbb T^N$-extensions of $X$}, that are  characterized by the existence of  {\it isometric} sections
$X\to X^\rtimes $ for  $X^\rtimes \to X=X^\rtimes/\mathbb T^N$.

Clearly,   metrics $  g^\rtimes$  on these $X^\rtimes $ are 
$$g^\rtimes=g+\varphi^2_1(x)dt_1^2+\varphi^2_2(x)dt_2^2+...+\varphi^2_N(x)dt_N^2$$
for some positive functions $\varphi_i$ on $X$.
%$\tilde X\times = X \times \mathbb R^N$.   

Among these we distinguish {\it  $O(N)$-invariant} warped extensions, where the $\mathbb Z^N$ covering manifolds $\tilde X^\rtimes = X \times \mathbb R^N$, where  
$$\tilde X^\rtimes/\mathbb Z^N= X^\rtimes,$$
 are invariant under the action of the orthogonal group $O (N)$. Thus, $\tilde X^\rtimes$ are acted upon by the full isometry group of $\mathbb R^N$, that is $\mathbb R^N\rtimes O (N)$.

Equivalently, the  metric in such an  $X^\rtimes$ is a   "simple"  warped product:
$g^\rtimes=g+ \varphi^2d||\bar t||^2$ for  $\bar t=(t_1,t_2,...,t_N)$,  the scalar curvature of which, as we know, \ref {warped2} is 
$$Sc(g^\rtimes)(x,\bar t) =Sc(X)(x)-\frac {2N}{\varphi(x)}\Delta_g  \varphi(x)-\frac  {N(N-1)}{ \varphi^2(x )}$$
and which is most simple (and useful) for $N=1$, where 
 $$Sc(g^\rtimes)(x,\bar t) =Sc(X)(x)-\frac {2}{\varphi(x)}\Delta_g  \varphi(x). \leqno { [\color{blue} \rtimes_\varphi]}$$
for the Laplace (Beltrami)  $\Delta_g$ on $X=(X,g)$.\vspace {2mm}

%%%%%%%%%%%%%%%&&&&&&&&&&&&&
{\large \it \color {blue}  $[ \rtimes_\varphi]^N$}-\textbf{Symmetrization Theorem.}  {\sf Let $X=(X,g)$   be  a closed oriented Riemannian manifold of dimension $n=m+N$  and 
let
 $$X\supset X_{-1}\supset ...\supset X_{-i} \supset... \supset X_{-N},$$  
be a descending chain of closed oriented submanifolds,
where each $X_{-i}\subset X$ is equal to a transversal intersection of $X_{-(i-1)}$ with  a smooth closed oriented hypersurface    $H_i\subset X$, 
 $$H_i\cap X_{-(i-1)}=X_{-i}.$$}
    If $n\leq 7$, then 

{\it there exists a closed oriented $m$-dimensional  submanifold $Y\subset X$ homologous to   
$X_{-N}$ and  a warped  product
$\mathbb T^N$-extension  $Y^\rtimes$ of $Y=(Y,h)$ for the Riemannian metric $h$ on $Y$ induced from $g$ on $X$, such that 
the scalar curvature of  $Y^\rtimes$, that is, being $\mathbb T^N$-invariant, is represented by a function on $Y$,  is bounded  from below by the Scalar curvature of $X$ on $Y\subset X$,
$$Sc(Y^\rtimes,y)\geq Sc(X,y),\mbox { } y\in Y.$$}

 {\it Proof.} Proceed by induction on codimension $i=1,2,,....N$ and construct submanifolds 
 $$X\supset Y_1\supset... \supset Y_i\supset...\supset Y_N=Y\subset X$$
  as follows.
 
 At  the first step, let 
  $Y_{1}  \subset X$ be a volume  minimizing, hence stable,  hypersurface homologous to $X_{-1}$
  where, the positivity of the second variation implies the positivity of the 
  $$-\Delta +\frac {1}{2} (Sc(Y_1)-Sc(X)|_{Y_1},$$
  for the Laplace  $ \Delta=\Delta_{h_1} $   on $Y_1$  with  the metric $h_1$ induced from $X$   and 
  let $\psi_1>0$ be the first eigenfunction of this  with the positive eigenvalue $\lambda_1$, thus 
  $$-\Delta\psi= \left (\lambda-  \frac{1}{2} (Sc(Y, h_1)-Sc(X))\right)\cdot \psi_1.$$

 Here, let
  $h_1^\rtimes(y) =h_1(y)+\psi^2dt^2$ 
  be the warped product   metric on $Y_1\times \mathbb T^1$  and observe 
 $$Sc(h_1^\rtimes, y)= Sc (h_1, y)-\frac {2}{\psi}\Delta\psi_1=Sc(X,y) +2\lambda_1.$$

 Then,  at the second step,  let $Y_2\subset Y_1$ be a hypersurface, such that  $Y_2\times \mathbb T^1\subset  Y_1\times \mathbb T^1$ is volume minimizing for the metric $h_1^\rtimes$, which  is equivalent 
  for $Y_2$ 
  to
  be volume minimizing in   $Y_1$ with respect to the metric $\psi^{l_1}_1h_1$  for 
  $l_1=\frac {2}{n-1}$. 

  Thus we obtain 
  $Y'_2$, where the corresponding metric on $Y'_2\times \mathbb T^2$ is
   $$h'_2+\psi_1^2dt_1^2 + \psi_2^2dt_2^2.$$
  
  Repeating this $N-2$ more times, we arrive at $Y_N'$ and an (iterated) warped product metric 
   $$h'_N+\sum_{i=1}^N\psi_i^2dt_i^2\mbox   {  on $Y'_N\times \mathbb T^N,$} $$
 which    can be symmetrised further to the required $h^\rtimes$   by applying the above infinitely many times  to hypersurfaces $Y'_N\times T^{N-1}\subset Y'_N\times T^{N}$ for  all  subtori   $T^{N-1}\subset Y'_N\times T^{N}$.\footnote{  See  in, \S12[GL(complete)1983], [G(inequalities) 2018] and also  the sections  \ref{separating3}, \ref{separating5} for details of this argument and for  generalizations.} (The luxury of the extra $O(N)$-symmetry is unneeded for most purposes.)

  \vspace {1mm}
  
{\it Exercise.}   Apply {\large \it \color {blue}  $[ \rtimes_\varphi]^N$}-symmetrization to $n$-manifolds 
 with    isometric   $\mathbb T^{n-2}$-actions  and prove the above 
 equivariant area inequality by reducing it to the warped product case  that was  already  settled in section   \ref {warped+2}. \vspace {1mm}

{\it Symmetrization by Reflections and Convergence Problem.} Let $Y$ be a closed minimal co-orientable (i.e. two sided)  hypersurface in a Riemannian manifold. If $Y$ is   locally volume minimizing, then it admits arbitrarily 
small neighbourhoods $V_\varepsilon \supset Y$ in $X$ with smooth {\it strictly mean convex} boundaries. 
 Then by reflecting such a varepsilon in the two boundary components, one obtains manifolds $\hat V_\varepsilon $ with isometric actions of $\mathbb Z\rtimes \mathbb Z_2$.
 
If these $Y$ are non-singular, e.g. if $dim(X)\leq 7$, then one can take solutions of the isoperimetric problem  for these   $V_\varepsilon$, where one  minimize the volumes of both components of the boundaries of $V_\varepsilon$ per given (small) volume contained  between them and $Y$.
In this case,  $\hat V_\varepsilon $,  $\varepsilon\to 0$, converge to smooth Riemannian  manifolds $V^\rtimes$ with isometric actions of 
$\mathbb R$ and with their scalar curvatures bounded from below by $Sc(X)|_{Y}$.

 If $Y$ is singular, the boundaries of  these   $V_\varepsilon$, even if singular, \footnote{   If $n=8$, then, by adapting Nathan Smale's argument, one can show that  these  $V_\varepsilon$ are non-singular for an open dense set of values of $\varepsilon$; but this is problematic for $n\geq 9$.} can be smoothed with positive mean curvatures, but it is unclear if they converge to a reasonable object for $\varepsilon \to 0$:
what is {\it \color {red!60!black}  missing  for convergence}  is a {\it Harnack type inequality} for the boundary components 
of
  $\partial_1, \partial_2  \subset \partial V_\varepsilon$, that is  a uniform bound for the ratios of the distances 
  $$\frac{dist(y, \partial_i)}{dist(y', \partial_i)}, y,y'\in Y,$$ 
$i=1,2$, and /or of distances $dist(x,x', Y)$, $x,x'\in\partial_i$.
 
Notice, that   "symmetrization by reflections", albeit open to generalizations to singular $Y$, is not, apparently,  applicable, to stable $\mu$-bubbles $Y$, where the warped product construction does apply.
 \footnote {See \S8 in [G(billiards) 2014],  \S4.3 in [G(inequalities) 2019] and section \ref {variation5} for more about all this.}

\vspace{1mm}

{\it Symmetrization versus Normalization.}   $\mathbb T^\rtimes$-Symmetrization of metrics  $g$  typically)  makes  their scalar curvatures  {\it  constant} by  paying  the price of modification of the topology of the  underlying manifolds, $X\leadsto X \times \mathbb T^1$. 

As far as sets of "interesting"   maps between Riemannian  manifolds are concerned a  similar effect effect is  achieved by keeping the same manifold $X$ but modifying  the metric by  $g=g(x)\leadsto  g^\circ =g^\circ(x)=Sc(X,x)g(x)$.

In fact, we   shall see later  in many examples, that  

 {\sf there is a close (but not fully understood)  similarity between the sets of {\sf $\lambda^\circ$-{\it Lipschitz} maps  $(X, g^\circ)\to (Y, h^\circ )$
and of  
 {\it $\mathbb T^1$-equivariant 
  $\lambda^\rtimes $-Lipschitz} maps $(X\times \mathbb T^1,g^\rtimes)\to (Y\times \mathbb T^1,h^\rtimes)$ for     $\lambda^\circ$  and $\lambda^\rtimes$ related in a certain way.}}

%%%%%%%%%%%%%%%%%%%%%%%%
\subsection {\color{blue} Positive Eigenfunctions and the Maximum Principle}  \label{eigenfunctions2}
%%%%%%%%%%%%%%%%%%%%%%%%%%%%%%

Let $X$ be a {\it compact connected}  Riemannian manifold and  let 
$$\Delta f= \sum_i\nabla_{ii}f= {\sf  trace}\hspace {0.5mm}{\sf Hess} f= {\sf div} \hspace {0.5mm}{\sf grad} f$$
denote the Laplace (Beltrami)  on $X$, which, recall, is a {\it negative} , since 
$$\int_X \langle f,\Delta f \rangle dx =-\int_X ||{\sf grad}f||^2 dx\leq 0$$
by Green's formula.

\textbf { Non-Vanishing Theorem.} {\sf  Let  $s(x)$ be  a smooth function, such that
the  
 $$L=L_s:f(x)\mapsto-\Delta f(x)+s(x)f(x)$$
   is {\it non-negative}, that is  $\int_X  \langle f(x),Lf(x)\rangle dx \geq 0$ for all $f$ or, equivalently,  if 
$L$  the  lowest eigenvalue $\lambda=\lambda_{min}$ is     $\geq 0$}.\footnote {This is equivalent since our  $L$ has {\it discrete spectrum.}}

Then\vspace {1mm}

 \hspace {0mm} {\it the eigenfunction $f(x)$ associated with $\lambda$ doesn't vanish anywhere on $X$}.
\vspace {1mm}

Start with   two lemmas.\vspace {1mm}

1. {\it $C^1$-Lemma.}  {\sf If the minimal  eigenvalue of the  
$f(x)\mapsto Lf(x)=-\Delta f(x)+s(x)f(x)$ on a compact Riemannian manifold is {\it non-negative},
$\lambda=\lambda_{min}\geq 0$, then 
the {\it absolute value} $|f(x)|$ of the eigenfunction $f$ associated with   $\lambda$  is 
{\it $C^1$-smooth.}}

\vspace {1mm}

2. {\it $\Delta$-Lemma.} {\sf Let $f(x)$ be a   {\it non-negative}   continuous function on a  Riemannian manifold, such that 

(i) $f(x)$  {\it vanishes} at some point in $X$, 
$$\mbox{ $f(x_0)=0 $,   $x_0\in X$},$$

(ii) $f(x)$  is {\it not identically zero} in any neighbourhood of  the  point $x_0\in X$,
\vspace {0.6mm}

(iii)  $f(x)$ is everywhere   $C^1$-{\it smooth}  and   it is   $C^2$-{\it smooth} at the points $x$ where 

it doesn't vanish.}
 \vspace {0.5mm}
 
\hspace {-6mm} {\it Then there exists a sequence of points $x_1,x_2,...\in X$  {\sf  convergent to} $x_0$, 
        where $f(x_i)>0$ and such that 
$$ \frac { \Delta f(x_i)}{f(x_i)}\to \infty,\mbox  {   for } i\to \infty.$$}
\vspace {1mm}

{\it Derivation  of Non-vanishing Theorem from the Lemmas.} Since  $|f|$ is $C^1$ by the first lemma, the $\Delta$-lemma,  applied  to $|f(x)|$,   shows that there exists a  point $x$, where $f(x)\neq 0$
 and  $$\frac{\Delta f(x)}{f(x)}=\frac{\Delta |f(x)|}{|f(x)|} > |s(x)|,$$
that  is incompatible with $ -\Delta f(x) + s(x)f(x)=\lambda f(x)\geq 0$  for $\lambda\geq 0$.\vspace {1mm}

{\it Proof of  $C^1$-Lemma.}  Recall that the eigenvalues of the  $L=L_s=-\Delta+s$
are equal to the critical values of {\it the energy functional} 
 $$E(f)= \int_X( ||{\sf grad}f(x)||^2 +s(x))f^2(x) dx$$
on the  sphere 
$$||f||^2= \int_Xf^2(x)dx=1$$ 
in the Hilbert space $L_2(X)$ and   the  critical points of $E$ are represented by eigenfunctions

Indeed, 
 $$E(f)=\langle f, Lf \rangle=\int_X\langle f(x),Lf(x)\rangle dx$$
by   Green's formula and the differential of the quadratic function $f\mapsto\langle f, Lf \rangle$ on  the  sphere $||f||^2=1$  is
$$(dE)_f(\tau) = \langle \tau,Lf\rangle\mbox  {   for all for all $\tau$ {\it normal to} $f$}.$$ 
Thus, vanishing of $dE$ at $f$ on  the unit sphere says, in effect, that $Lf$ is 
a multiple of $f$, i.e. $Lf=\lambda f$.

All this makes sense in the present case, albeit the space $L_2(X)$ is {\it infinite dimensional}  and $L$ an {\it unbounded}, because $L$ is an  {\it elliptic} operator, which implies, for compact $X$, that  

{\it the spectrum of $L$  is discrete,  
 bounded from below and   all eigenfunctions  are smooth.}

In particular  -- this is all we need,

{\it all  minimizes  of $E(f)$  on the unit sphere, that are, a priori, only Lipschitz continuous, are  smooth.}\footnote{Recall that  our  "smooth" means $C^\infty$ and all our Riemannian manifolds are assumed smooth.} 

 \vspace{1mm}
 Now, observe that, 
  
  \vspace{1mm}
 
 {\it taking absolute values of smooth functions $f(x)\mapsto |f(x)|$ doesn't change  their    
  energies, as well as their $L_2$-norms,  
  $$\Vert |f|\Vert= \Vert f\Vert= \sqrt{\int_X|f|^2(x)dx},$$ 
 $$E(|f|)=E(f)=\int_X( ||{\sf grad}|f|(x)||^2 +s(x))|f|^2(x) dx,$$}

Indeed,   absolute values     $|f|(x)$  are Lipschitz for Lipschitz $f$, hence, they are almost everywhere differentiable functions, such that
${\sf grad}|f|(x)=\pm{\sf grad}f(x)$ at all differentiability points $x$ of $|f|$.

 It follows that  the absolute value of the eigenfunction $f$ with the smallest energy $E(f)=\lambda_{min}$ is also a minimizer; hence, {\it this $|f|$ is smooth.} QED.\vspace{1mm}

{\it Poof of $\Delta$-Lemma.} 
The common strategy for locating points $x\in X$ with "sufficiently positive"  second differential of a function $f(x)$ is by  using   simple auxiliary functions   $e(x)$  with this property and looking for  
minima  points for  $f(x)-e(x)$.

The basic example of such a  function  $e(x)$ in one variable is $e^{-Cx}$, $x>0$, for large $C$,
where $\frac {e''}{e}=C^2$, and where observe that the ratio $\frac {e''}{e'}=C$  also becomes large for large $C$.

It follows that that the Laplacians of the corresponding radial functions in small $R$-ball $B_y(R)$  in Riemannian manifolds $X$, 
 $$e(x)=e_C(x)=e_{y,C}(x)= e^{-C\cdot r_y(x)} \mbox { for } r_y(x)=dist(y,x)\leq R$$
satisfy
$$\Delta e(x)\geq C^2 e(x)-  C\cdot mean.curv ( \partial B_y(r),x)\mbox { for } r= r_y(x)=dist(y,x)$$

Now, in order to find a point $x$ close to a given $x_0\in X$ where $f(x)=0$,   take $y\in X$
very close to $x_0$, where $f(y)>0$,  let $B_y(R)\subset X$ be the {\it maximal} ball,
such that $f(x)>0$ in its interior, let
$$e(x)=e_C(x)=e^{-C\cdot r_y(x)}-e^{-C\cdot R}$$
 and observe that $ e(x)$ vanishes on the boundary of the ball  $B_y(R)$ and is strictly positive in the interior. Moreover 
 $$e(x)\geq \varepsilon\rho,$$
for all $x$ on the geodesic segment between $y$ and $x_0$ within distance $\geq \rho$ from $x_0$
 for all $\rho_0\leq R$.
 
Notice that   this $ \varepsilon=\varepsilon_C$ albeit  {\it strictly positive},  tends to zero for $C\to \infty$.

Assume without loss of generality that $x_0$ is the only point in $B_x(R)$ where $f(x)$ vanishes (if not, move $y$ closer to $x_0$ along the geodesic segment between the two points),  
let  $C$ be {\it very very large} and see what happens  to $f(x)$ and $e(x)$ in the vicinity     of $x_0\in \partial B_y(R)$, say  in the intersection  $$U_0=B_y(R)\cap B_{x_0}(R/3).$$

Observe the following.

$\bullet$  Since   $f(x)>0$ for  $x\in B_y(R)$,    $x\neq x_0$, and  since   $e_C(x)\to 0$ for $C\to \infty$ for 
$r_y(x)=dist(y,x)\geq r_0>0$,
  {\sl the function $e(x)=e_C(x)$, for large $C$,  is bounded by $f(x)$ on the boundary of $U_0$, 
$$e(x)\leq f(x), \mbox { }  x\in \partial U_0,$$
where $e(x)< f(x)$ unless $x=x_0$.}

$\bullet $  Since $f$ is differentiable  at $x_0$  and assumes minimum at this point, 
the differential $df$ vanishes at $x_0$, which makes  
$f(x)=o(\rho)$   for $\rho=dist (x, x_0)$, 
{\sl there is a part   of  (the interior of) $ U_0$,  where $e(x)>f(x).$}

Hence, the difference $f(x)-e(x)$ assumes minimum at an interior point $x=x_{y,C}\in  U_0$,
such that $x=x_{y,C}\to x_0$ for $C\to \infty$ and 
$$\frac{\Delta f(x)}{f(x)}\geq\frac { \Delta e(x)}{e(x)}\to\infty.$$
The proof of the $\Delta$-lemma  and  of the non-vanishing theorem are thus concluded.
\vspace {1mm}

{\it Discussion.} The non-vanishing theorem, which, probably,  goes back to Rayleigh,   is often used 
without being  even explicitly stated as, for instance, by Kazdan and Warner in their  "conformal change" paper. But I couldn't find an explicit reference on the web, except for  the paper 
   by Doris Fischer-Colbrie  and Rick Schoen, where they  prove  such a  non-vanishing for non-compact manifolds needed for their

\hspace {5mm} {\it non-existence theorem for non-planar    stable minimal surfaces in $\mathbb R^3$.} 

 Their argument  relies  on the {\it "strong maximum principle"} for the   $L$, for which they refer to pp. 33-34 of the canonical Gilbarg-Trudinger  textbook,   where the relevant case of this principle  is   stated (on p. 35 in the 1998 edition which is available on line)    after  the proof of theorem 3.5 as follows. 
\vspace {1mm}

{\sl "Also, if $u = 0$ at an interior maximum (minimum), then it follows from the proof of the theorem that u = 0, irrespective of the sign of $c$."}\vspace {1mm}

(The  assumptions of the  theorem specifically   rule out   $c$ with variable signs, where this $c=c(x)$ is the coefficient at the lowest term  in the equation    
 $Lu=a^{ij}(x)D_{ij}u +b^iD_iu+c(x)u=0$ introduced  on p. 30.)
\vspace {1mm}

What  is actually  proven in this book on  about  twenty lines on p. 34, is  a version of "$\Delta$-lemma"  for  $L$.

In our  proof, we reproduce what  is written on these lines, except for "{\sl direct 
calculation gives}"   that is replaced  by  an explicit evaluation of  $\Delta e(x)$ \footnote{In truth,   the only non-evident aspect  of the argument  resides with the identities 

\hspace{-4mm}$(e^{-Cx})'=  -Ce^{-Cx}$ and $(e^{-Cx})'' =  (-Ce^{-Cx})'=C^2e^{-Cx}$ with the issuing inequalities  $(e^{-Cx})''>>e^{-Cx} $ and $(e^{-Cx})''>> |(e^{-Cx})'|$, which {\it can't be done by just staring} at the exponential  function.
(The  appearance of  $e^x$, that is an isomorphism between the {\it additive}  $\mathbb R$ and {\it multiplicative} $\mathbb R_+^\times$
with all its  counterintuitive  properties, is  amazing here -- there is  nothing visibly multiplicative in $\Delta$; besides, the geometric proof of the existence of $e^x$ via 
 the  {\it conformal} infinite cyclic covering map  $\mathbb C\to \mathbb C\setminus \{0\}$ and analytic continuation 
is    non-trivial.)

The rest of the proof is geometrically effortless: you {\it  just  look} at the graph $\Gamma_e$ of the function  $e(x)=\exp-C\cdot dist(y,x)$ in a small $R$-ball $B\subset X$  outside zero set of $f$  with the center of your choice, such that $B$  touches this set at $x_0$, and let $C=C_i\to \infty$.
Then you {\it see}  a tiny region  
 in this ball close to $x_0$, where $\Gamma_e$ mounts  above $\Gamma_f$, and  you take  the point in  $X$  just  under the {\it top} of this mountain, i.e. where the      distance measured vertically between   the two graphs is maximal,  for you $x=x_i$.}\vspace {1mm}
 
 The  following (obvious) corollary to the non-vanishing theorem  will be used for construction of stable  symmetric  $\mu$-bubbles  in sections \ref{existence5}, \ref{separating5}.\vspace {1mm}

\textbf {Uniqueness/Symmetry Corollary.} {\it If $X$ is compact connected, then the lowest eigenfunction $f$ of the  $L$ is unique up to scaling. Consequently, if $L$ is invariant under an action of an isometry group on $X$, then, even if $X$ is disconnected, there exists a positive  $f$  invariant under this action.}
\vspace {1mm}

{\it Exercises.} (a)  \textbf{Multi-Dimensional Morse   Lemma.}  Show that two {\it non-coinciding}  volume minimizing hypersurfaces in the same indivisible homology integer homology class of an orientable manifold $X$  have {\it empty intersection and that, consequently,  volume minimizing  hypersurfaces must be {\it  invariant under symmetries} of $X$}.\footnote {This was used by Marston Morse to show that

{\sf  if 
 the $(n-1)$-dimensional homology group of some covering of a  compact Riemannian $n$-manifold, doesn't vanish
 then  the universal covering $\tilde X$ of $X$ contains an infinite minimal hypersurface  
  the image of which under the covering map  $\tilde X \to X$  is compact.}
   
   Morse was concerned in his paper 
"{\sl Recurrent Geodesics on a Surface of Negative Curvature}"  with the case of $n=2$ but his argument, transplanted to the environment of the geometric measure theory, applies to  manifolds of all dimensions $n$.}

(b) Generalize this to   $\mu$-bubbles, that are boundaries    of  domains $V$ in a Riemannian  manifold  $X$ that {\it minimize} the functional 
$$V\to vol_{n-1}(\partial V)-\int_V \mu(x)dx$$
 for a smooth function   $\mu(x)$. 
(Unit spheres $S^{n-1}\mathbb  R^n$ are {\it not minimizing}   $\mu$-bubbles  for $\mu=(n-1)dx$.)

(b)  \textbf {Courant's Nodal Theorem.} Show that the    that is the number of connected components of  the complement to the "$k$-th nodal set",  i.e. the  zero set of the $k$-th eigenfunction of $L=L_s=\Delta+s$ on a compact 
connected manifold,  can't have more than $k$ connected components. 

{\it Question.} Is there a counterpart to this for  non-quadratic functionals in  spaces of functions, or, even better, 
 spaces of hypersurfaces?

%%%%%%%%%%%%%%%%%%c
 \section {Topics, Results, Problems}\label{topics3}
%%%%%%%%%%%%%%%%%%%%%%%%%

%%%%%%%%%%%%%%%%%%%%%%%%%%%
\subsection {\color {blue}  Scale Persistent Criteria for $Sc\geq \sigma$ for Smooth and non-Smooth 
  Metrics}\label {Sc-criteria3}

%%%%%%%%%%%%%%%%%%%%%%%%%%%

  {\it Scale persistence } of a  geometric property {$\sf P$} applicable to compact   $n$-dimensional  Riemannian manifolds $V$ with boundaries,  means the following:   
 
 {\sf if  such a {$\sf P$} is satisfied by  small neighbourhoods of  all points in a Riemannian $n$-manifold $X$
 then it is satisfied for all domains $V$ in $X$.}

 Two classical examples of these are the following characteristic properties of surfaces with 
 {\it non-negative 
sectional curvatures $\kappa$ } and of  $n$-dimensional manifolds with {\it non-negative Ricci curvatures}.\vspace {1mm}

In the case of  the sectional curvature we formulate such a property  as a    {\it comparison inequality }   for {\it geodesic quadrilaterals}  as follows.
 \vspace {1mm}

$\Square_{\kappa\geq 0}$  {\it All convex, e.g. geodesic,  quadrilaterals in surfaces with non-negative curvatures, call them  $\Square \subset X$,   have the greatest   of their four angles  at least $90^\circ$}:
$$\max_{i=1,...,4}\angle_i(\Square)\geq \frac {\pi}{2}.$$

In fact, the sum of the four angles of such a $\square $ must be $\geq 2\pi$ by  the {\it Gauss-Bonnet inequality} for  compact 
 surfaces $V$ with (quadrilateral in the present case) boundaries $\Theta$:
 
{\it  if  $\kappa(V)\geq 0$, then
$$\int_\Theta curv(\Theta,\theta)d\theta\leq 2\pi.$$}
(The curvature of $\Theta=\partial V$  at a vertex of $V$  with the angle $\alpha$ is the  point-measure with the weight $\pi-\alpha$.)

It is also clear that $\Square_{\kappa\geq 0}$ is {\it sufficient}, as well as necessary, for $\kappa\geq 0$: \vspace {1mm}

{\it if $\kappa(X,x)<0$, then there exist (small) geodesic  quadrilaterals in $X$ around x with all angles 
$<\frac {\pi}{2}.$}\vspace {1mm}

Thus,{\it 

\hspace {20mm} local validity of $\Square_{\kappa\geq 0}$   implies the global one.} \vspace {1mm}

(Also notice that if  all four angles  of a convex  $\square$ with $\kappa(\square)\geq 0$ are $\leq \frac {\pi}{2}$, then this  $\square$ is isometric to a plane Euclidean rectangular.)

\vspace {1mm}

Next, turning  to Ricci, observe  that the   inequality $\Circle_{Ricci\geq 0}$ stated below  says, in effect, that the mean curvatures of the boundaries  of compact 
manifolds with $Ricci\geq 0$ can't be greater than these  of  Euclidean balls of comparable size.

\vspace {1mm} 

{\large \color {blue} ${\sf P}_{Ricci\geq 0}$.}  {\sf If   $Ricci(V)\geq 0$, then 
the minimum of the {\it mean curvature} of the  boundary of $V$  is related  to
 the {\it inradius} of $V$  by the inequality
$$\inf_{v\in \partial V}mean.curv(\partial V, v)\leq\frac {n-1}{inrad(V)}, \mbox { } n=dim(V),$$}
where 
{\sf $$inrad(V)=\sup_{v\in V}dist (v, \partial V),$$ 
 where our sign convention for the  mean curvature is such that convex domains have 
$mean.curv\geq 0$
and where  

{\sf the (opposite)   inequality} {\it  
$$mean.curv(\partial V, v)\geq\frac {1}{inrad(V)},$$
implies that $V$ is isometric to the Euclidean ball of  radius $R=inrad(V)$.}}
\vspace {1mm}

 All this  follows from   Hermann Weyl's tube formula  applied to concentric spheres in $V$ around the point 
 $v\in V$ farthest from the boundary.

\vspace {1mm}

Let us now  state   two inequalities   that   characterize $n$-manifolds $X$  with non-negative scalar curvatures, where 

{\sf the first one says that } {\it cubical domains  $V\subset X$  can't be   "more mean convex",  than  rectangular solids  in the Euclidean space},

and  the  second one, that    applies to domains $V\subset X$ with smooth boundaries, claims  
that  

{\sf these boundaries } 
{\it can't be  simultaneously greater in size and  "more mean convex"  than  convex hypersurfaces  in the Euclidean space. }

%which say, similarly to $\Square_{\kappa\geq 0}$  and $\Circle_{Ricci\geq 0}$,  but in a different way,   that boundaries of smooth domains $V\subset X$ can't be   "more mean convex"  than  convex hypersurfaces  in the Euclidean space. 

 \vspace {1mm}
 %%%%%%%%%%%%%%%%%%%&&&&&&&&&&&&&

{ \textbf I. } \hspace{0mm} {\color {blue}$ \blacksquare$}-{\textbf {Inequality.}} {\sf Let $V$ be a Riemannian $n$-manifold  diffeomorphic to the cube $[0,1]^n$.} 
 
{\it  If  $Sc(V)\geq 0$ and if all  $(n-1)$ faces  $\partial_i \subset \partial V$, $i=1,...,2n$, have
 $mean.curv (\partial_i) \geq 0$, then the  supremum of the dihedral angles between the  (tangent spaces of) $(n-1)$-faces at the  points in the $(n-2)$ faces  satisfy
$$\sup_{i,j}\angle_{ij}(V)\geq \frac {\pi}{2}.\leqno  {\color {blue}  \blacksquare}$$} 
  
 This may serve as a criterion for $Sc(X)\geq 0$,  since \vspace {1mm}

{(\large \color {magenta}$\square$}) {\it the inequality $Sc(X,x)<0$ implies the existence of a small (topologically) cubical mean convex neighbourhood $V\subset $X of $x$ which violates {\color {blue}{ $ \blacksquare:$}}

\hspace {19mm}all dihedral angles of $V$  are everywhere $>\frac {\pi}{2}$.}
 \vspace {2mm}
 %%%%%%%%%%%%%%%%%
 
\hspace{0mm}  { \color {blue}\EllipseSolid}
\vspace{-5mm}

 { \hspace {-6mm}{ \textbf {II}. } { \hspace {4.5mm}-{\textbf {Inequality}. } {\sf Let $V$ be a Riemannian manifold diffeomorphic to the $n$-ball
 the boundary of $V$ of which has  positive mean curvature,
$$mean.curv(\partial V)>0,$$
 let  $\underline V\subset \mathbb R^n$ be  a  convex domain with smooth boundary  and 
let  
$$f:\partial V\to \partial \underline V$$ be a diffeomorphism.\footnote { It is enough to assume that  $f$ is  a smooth map with {\it positive degree}  as it will become clear later on.}}
  
{\it  If  $Sc(V)\geq 0$, then  the differential of $f$  can't be everywhere  strictly  smaller than  the ratio of the mean curvatures of the  two boundaries:
there exists a point $v\in\partial V$, such that 
$$||df(v)|| \geq \frac {mean.curv(\partial V, f(v))}{mean.curv(\partial \underline V, v)}.\leqno
 \mbox { {\color {blue}\EllipseSolid$_\geq$}}$$}

  \vspace {2mm}

({\color {magenta}\Ellipse$<$})
\vspace{-5mm}

{\it \hspace {7mm} {\sf Conversely,}
 the inequality $Sc(X,x)<0$ implies  the existence of $V\subset  X$,  $\underline V\subset \mathbb R^n$ and of a diffeomorphism
  $f:\partial V\to \partial \underline V$,
such that 
$$||df(x)|| < \frac {mean.curv(\partial  V, f(v))}{mean.curv(\partial \underline V, v)} \mbox { for all } v\in \partial V.$$}
 
 We indicate the proofs  of {\color {blue}{ $ \blacksquare$}}  and  { \color {blue}\EllipseSolid} in the next section, and refer to section \ref{trace-norms3} for a generalization of  { \color {blue}\EllipseSolid} to topologically non-trivial   manifolds  $V$;
   below, we   turn  to manifolds with $Sc\geq \sigma\neq 0$.

 \vspace {2mm}

\textbf {Corollaries } of { \textbf I }  and { \textbf {II} } for {\it manifolds $X$ with  $Sc(X)\geq \sigma$  for $\sigma \gtrless 0$ }. {\sf The  inequalities {\color {blue}  $\blacksquare$} and { \color {blue}\EllipseSolid},
 when   applied to manifolds  $X$ multiplied by    surfaces  $S$ with scalar curvatures $- \sigma$, 
yield  } \vspace {1mm}

\hspace {24mm}{\it geometric criteria for  $Sc(X)\geq \sigma$  for all $\sigma$.}\vspace {2mm}

The geometric meaning  of this, if any, is obscure;  possibly, it   can be  expressed  in terms of 2-parametric families of domains  
 $V_s$, $s\in S$. 
 But  the  following  generalizations of
{ \color{blue}\EllipseSolid} to $\sigma>0$ and of { {\color {blue}  $\blacksquare$} to $\sigma<0$ are geometrically   transparent.

  % is available  for domains $V$ in manifolds  $X$ with  $Sc(X)>\sigma>0$, that is achieved by 
  %comparing   the mean curvatures of the boundaries $\partial V$  with  these of convex domains $\underline V$ in 
 %the $n$-sphere of radius $\sqrt {\frac {n(n-1)}{\sigma}}$, where the $\sigma$-counterpart of  {\color {blue}\EllipseSolid}  reads  as follows.
 
 \vspace {1mm}

   \textbf{{ \color{blue}\EllipseSolid}-Comparison Theorem for $\mathbf {Sc > 0}$}.  {\sf 
   Let $V$  and $\underline V$  be  compact Riemannian $n$-manifolds with smooth boundaries,
   where  $\underline V$ has constant  sectional curvature $+1$  and the  boundary $\partial \underline V$ is convex and
     and let $f:V\to \underline V$ be a diffeomorphism.}

 Then either  

{\it there  exists a point $v\in V$, where 
 the norm of the exterior square  of the differentials of $f$   is bounded from below by 
 $$||\wedge^2 df(v)|| \geq \frac{ Sc(V,v)}{n(n-1)} $$}
or, as earlier, 

{\it there exists a point $v'\in \partial V$, such that 
$$||df(v')|| \geq \frac {mean.curv(\partial  V, f(v'))}{mean.curv(\partial  \underline V, v')}.\leqno
 \mbox  {\color {blue}\EllipseSolid$'_{>0}$}$$}}

\vspace {1mm}
 
 {{ {\color {blue}  $\blacksquare$}}-Comparison Theorem for $\mathbf {Sc(V)> \sigma<0}.$} Let 
  $$(\mathbf H^n, g_{hyp})= \mathbb R^1\rtimes \mathbb R^{n-1}= 
  (\mathbb R^1\times \mathbb R^{n-1}, dt^2 +e^{2t}dx^2)$$ 
 be 
 the hyperbolic space with sectional curvature -1 represented as  the warped product in the normal  horospherical coordinates, let 
 $$\underline V=[0,1] \times  [0,1]^{n-1}\subset  \mathbb R^1\rtimes \mathbb R^{n-1}=\mathbf H^n$$
  and 
  observe that all dihedral angles  in $\underline V$ are $\frac {\pi}{2}$, all "side faces" are geodesic flat,   
  while the "bottom" $\{0\}\times    [0,1]^{n-1}\subset \underline V$ and the "top" $\{1\}\times  [0,1]^{n-1}\subset \underline V$,
  have mean curvatures $-(n-1)$ and $n-1$ respectively.
  
{\sf The corresponding   comparison inequality for cubical Riemannian manifolds
 $V$  diffeomorphic  to  $[0,1]\times [0,1]^{n-1}$ reads.}

 {\it Let all dihedral angles of $V$ be $\leq \frac {\pi}{2}$,  let all ("side") faces $\partial_i\subset V $, except for 
 $\partial_0=\{0\}\times    [0,1]^{n-1}$ and $\partial_1=\{1\}\times  [0,1]^{n-1}$,  have non-negative mean curvatures and let 
 $$ mean.curv(\partial_0)\geq -(n-1)\mbox { and  } mean.curv(\partial_1)\geq n-1.$$
 Then the scalar curvature of $V$ can't be everywhere greater than  that of $\mathbf H^n$,
  $$\inf _{v\in V}  Sc(V,v)\leq -n(n-1).\leqno {\color {blue}  \blacksquare'_{<0}}$$}
 
 {\it Remarks.} 
 (a)The {\it proofs} of these are indicated in the sections \ref{reflection3} below. 
 
 (b) {\color {red!40!black}Probbaly} -- figuring this out this way or another can't be too difficult -- these 
 {\color {blue} {\small \EllipseSolid}$'_{>0}$} and   {\color {blue}$\blacksquare'_{<0}$} {\it characterizes} $Sc\geq \pm 1$. 

 (c) Granted (b),  either of {\color {blue}  $\blacksquare'_{<0}$} or  {\color {blue}\EllipseSolid$'_{<0}$ } \hspace {1mm}can 
 be  used for characterization of $Sc\geq \sigma$  for all $\sigma$ by passing to  products of $X$ with $S^2$ or $\mathbb H^2$ as we did earlier. \hspace {1mm}
 
  (c) The proof of  {\color {blue} $\blacksquare'_{<0}$} for $n\geq 9$, which relies on stable $\mu$-bubbles,  needs (a slight generalization of)  the desingularization theorem  from [SY(singularities) 2017] or of such a result from  [Lohkamp(smoothing) 2018].

 (d) {\color {red!40!black} Probably}, a combination of ideas from  [Min-Oo(hyperbolic) 1989]  and from recent papers by Cecchini, Zeidler, Lott and 
 Guo-Xie-Yu on index theorems for manifolds with boundaries\footnote{See 
  [Cecchini-Zeidler(Scalar\&mean) 2021], [Lott(boundary) 2020], [Guo-Xie-Yu(quantitative K-theory) 2020.} may provide  an  alternative proof of  
 $\blacksquare'_{<0}$  for all $n$.

  {\it {\color {blue}  $\blacksquare$} And 
  {\color {blue}\EllipseSolid}\hspace {1mm}
 for Continuous Riemannian Metrics.}
One can define mean convexity and, more generally,  lower bound of the mean curvatures from below 
 for  boundaries $\partial V$ of domains $V$ in a metric space $X$, whenever one has a notion of the volume/measure  for $\partial V$ as follows.\vspace {1mm}
  
{\sl $mean.curv(\partial V, v)> m$, $v\in \partial V$, if there exists a sequence of subdomains $V_i \subset V$ with the following properties.}

 {\sf(i) The difference between  $V$ and $ V_i$ contains a neighbourhood  $v$  in $V$ for all $i$  and it   converges to  $v$ for $i\to \infty$, i.e.  $V\setminus V_i$ is contained in the $\delta_i$-ball  around $v$ for $\delta_i\to 0$.

(ii) the volume of $\partial V_i$ is bounded in terms of the volume of the part of $\partial V$ outside $V_i$ and 
the  Hausdorff distance between the boundaries of $V$ and $\partial V_i$ as follows:
$$vol(\partial V_i) < vol(\partial V)- m \cdot vol(\partial V\setminus V_i)\cdot dist_{Hau}(\partial V, \partial V_i).$$} 

With   this "mean curvature", the definitions of  {\color {blue}  $\blacksquare$}  and {\color {blue}\EllipseSolid}\hspace{1mm} as well as    {\color {blue}  $\blacksquare_{<0}$}  and {\color {blue}\EllipseSolid}$_{>0}$}\hspace {1mm}
automatically extend  to continuous, even only Borel,  Riemannian metrics. \vspace {1mm}

{\it \color {magenta!80!black} Question}. Do    {\sf {\color {blue}  $\blacksquare$}  and {\color {blue}\EllipseSolid}
define the same concept of $Sc(g)\geq 0$   for {\it continuous} Riemannian metrics g?}\vspace {1mm}
%%%%%%%%%%%
\subsubsection {\color {blue} Reflection  Orbihedra and  Trapped Minimal Hypersurfaces} 
\label {reflection3}
%%%%%%%%%%
%%%

(1$_{\geq 0}$) {\sf \large Idea of the  Proof} of   {\color {blue}  $\blacksquare$.}  Reflect $V=(V,g)$  as a cube  in
 $ \mathbb R^n$  in the $(n-1)$-faces,  let   $ \hat V$ be the resulting {\it universal orbi-covering} manifold  with  an action of the relection group 
  $\Gamma$ that is  isomorphic  of  a finite extension the group $\mathbb Z^n$ for cubical $V$, and let $\hat g$ be the  (singular path) metric on $\hat V$ that $\Gamma$-invariantly extends $g$ from  $V$ naturally  embedded to $ \hat V$.

If  the mean curvatures of all codimension one faces are $\geq 0$, 
all  dihedral angles of $V$ are  $\leq \frac {\pi}{2}$ and the   dihedral angle at a point $v$ on some codimension two face of $V$   satisfies the strict inequality    $\angle_{ij}(v),\frac{\pi}{2}$, 
then one can {\it approximate   $\hat g$ by smooth $\Gamma$-invariant metrics} $ \tilde g $ for, such that 
$Sc( \tilde g)>\sigma$ for $\sigma=\inf_{v\in V}(Sc(V,v)$.\footnote 
{Working this out in detail  requires some patience, see  [G(billiards) 2014]) and [Nuchi(cube) 2018].}

Thus, if we assumed that  $Sc(V)\geq 0$,  this  $ \tilde g  $  would descend to a metric on the 
 torus  $\hat V/\mathbb Z^n$ with $Sc>0$ and the proof of  {\color {blue}  $\blacksquare$}  follows by contradiction: {\it there is no  metrics with positive scalar curvatures on the torus}.
\vspace {1mm}

(2$_{\leq 0}$.) {\it About   {\color {blue}  $\blacksquare_{<0}$}.}  Here we reflect $V$ around the $2(n-1)$  "side faces" and, thus,  after smoothing, reduce    {\color {blue}  $\blacksquare_{<0}$} to the comparison inequality for hyperbolic cusps  $\mathbb R^1\rtimes \mathbb T^{n-1}=\mathbf H^n/
\mathbb Z^{n-1}$. \footnote {See \S$5\frac{5}{6}$ in [G(positive) 1996], \S 9 in  [G(inequalities) 2018].}

The available proofs  of this  inequality,   apply only for $n=dim(X)\geq 3$, while the case $n=2$ follows by a simple  argument that we suggest as an elementary exercise to the reader.

{\it \color {red!39!black}  Remark.} The proof  of the comparison inequality for hyperbolic cusps 
 relies  on stable $\mu$-bubbles $\hat Y$ between pairs of $(n-1)$ tori in the cusps  $\mathbb R^1\rtimes \mathbb T^{n-1}$ where  $mean.curv(\hat Y)\geq n-1$.
 
 This  suggests a similar  proof directly in $V$  with   relevant $\mu$-bubbles  $Y\subset V$   having  free boundaries in   the {\it side boundary} of $V$ that is in the union of the  "side"  faces  i.e. all,  except for the two corresponding to the bottom $\{0\}\times [0,1]^{n-1}$  and the top of the cube $\{1\}\times [0,1]^{n-1}$.

But since  the mean curvatures of the side faces are only  assumed {\it non-negative}, such a bubble  with mean curvature  $\geq {n-1}$ may   meet a side face in the interior  and render  the argument invalid.

It is unclear how to  formulate the existence of  needed $Y\subset V$  without direct appeal to reflections 
in the "side mirrors".\vspace{1mm}

(3)  {\sf \large Other   Reflection Groups.} The above construction applies to all  {\it Riemannian   reflection orbifolds}, that  are   manifolds $(V, g_0)$ with corners that serve as fundamental domains $\Delta$  of reflection groups $\Gamma$, that act on the corresponding orbi-covering  manifolds  
 $\hat V$ as follows.
 
 Let  $g$ be a Riemannian metric on such a $V$, which   satisfies the following $2\frac {1}{2}$ conditions.

{\sf  $\bullet_1$ the   codimension 1 faces $\partial_ i$  of $V$ are $g$-mean convex:  
$$mean.curv_g(\partial_i V)\geq 0;$$
  
  $\bullet_4$  the   $g$-dihedral angles $\angle _{ij}$ of $V$ at the codimension 2 faces of $V$
 are bounded by the corresponding $g_0$-angles of $V$, 
 $$\angle_{ij}(V,g)\leq \angle_{ij}(V, g_0);\footnote {All  dihedral $ g_0$-angles  are    $\frac {\pi} {k}$, $k= 3,4,...$,
 where  $k$  are   half-orders of  the    stabilizer subgroup of the corresponding  faces  $\partial_{ij}$.
 Thus,  all dihedral angles of $(V,g)$ must be $\leq \pi/2$. }$$
   
 $\bullet^<_2$ there is a point $v$ on some codimension 2 face of $V$, where the above inequality is strict,
 
 $$\angle_{ij}(V,g)(v)<\angle_{ij}(V).$$}

 Then, as earlier for the cubical $V$, one can show, that 

 {\it  the  natural  singular metric  $\hat g$ on $\hat V$ can be approximated by smooth $ \tilde g$ 
 with $Sc ( \tilde g)>\inf_v Sc(V,g)(v)$.}\vspace {1mm}

{\it About Examples.} There are few $V\subset\mathbb R^n$ and  Euclidean  reflection groups, to which the above applies. In fact, all such     $V$ are  the products of segments and  triangles  with the angles
 (60$^\circ$, 60$^\circ$, 60$^\circ$),  (60$^\circ$, 30$^\circ$, 90$^\circ$) and   (45$^\circ$, 45$^\circ$, 90$^\circ$). 
 
But there are lots of non-Euclidean  orbifold  $V$, e.g.  with right-angled corners, (see  [Davis(orbifolds) 2008]),   the universal orbi-coverings $\hat V$ of which are {\it hyper-Euclidean}  and, hence, admit no $\Gamma$-invariant metrics with $Sc>0$ (see sections \ref {twisted1}, \ref {twisted3}). Therefore, \vspace {1mm}

\hspace {2mm} {\it the conditions $\bullet_1$,$\bullet_2$ and $\bullet^<_2$ imply that
$\inf_{v\in V}Sc(g, v)\leq 0$ for these $V$.}
 
\vspace {1mm} 

But, in general, the following problem, solutions of  special cases of which  are spread throughout this paper,
remains widely open, \vspace {1mm}

 {\large{\color {red!80!black}$\mathbf\pentagon$}}-{\sc \color{red!50!black}\textbf  {Problem.}}    {\sf Let $V$ be a compact manifold with  corners,
e.g. a closed  manifold, or,  at the other end of the spectrum,  diffeomorphic to a convex polyhedron in the Euclidean space. Find necessary (ideally, necessary and sufficient)  conditions on $V$ for the existence of  a Riemannian metric $g$ on $V$, such that:

(i)  the scalar curvature is bounded from below by a given real number, $Sc\geq \sigma $,

 (ii)  the mean curvatures of the  codimension 1 faces, call them $V_i$,  are are similarly bounded from below, $mean.curv_g(V_i)\geq \mu_i$   
 
 (iii)   the dihedral angles of all codimension 2 faces are bounded by given numbers, say,  $\angle_{ij}\leq \alpha_{ij}$.}
 \vspace {1mm}
 
The above   {\color {blue}  $\blacksquare$}-comparison theorem provides an instance of such a condition with $\sigma\leq 0$     (this, moreover,  characterizes metrics  with $Sc\geq \sigma$), but  the  
  inequality,  { {\color {blue}\EllipseSolid$'_{>0}$} for $\sigma>0$, unlike  {\color {blue}  $\blacksquare$}
 involve the  distance  geometry of $V$.
 
 It $n=2$, then it is not hard to show (an exercise to the reader)   that
 
 {\sl  if $\sigma\geq 2$ then 
 {\it \color {red!49!black}no} $k$-gonal Riemannian  (surface)  $V$ may have the faces (edges) with   curvatures $\geq \mu\geq- \varepsilon $  and the angles $\leq \alpha\leq \varepsilon $ for  a sufficiently small $\varepsilon =\varepsilon_k>0$.}
 
 But it is unclear if this condition is $\mathbb T^\rtimes$-stable, i.e. extends to $\mathbb T^{n-2} $-invariant (warped product) metrics  $g^\rtimes$   on $V\times \mathbb T^{n-2}$, and thus,  would allow the   reduction  of higher dimensions $n$   to $n=2$ by the (warped FCS)  $T^\rtimes$-symmetrization. 
 
 (The full solution of the   {\large \color{red!80!black}$ \pentagon$}-problem  remains  unsettled even for $n=2$.) 

 (6) {\sf \large Minimal Hypersurfaces in Cubical $V$.}  At the  core of the proof of   {\color {blue}  $\blacksquare$} lies non-existence of metrics with $Sc>0$ on the torus $X=\hat V/\mathbb Z^n$,  which in turn, can be proved in two different ways: by the Schoen-Yau's descent with minimal hypersurfaces or with a use of  twisted Dirac operators on $X$
 
To pursue the latter, one has to   describe/construct/analyse  twisted harmonic spinors on $\hat V$ in terms of  $V$ with no use  to the orbi-covering  condition $V\leadsto \hat V$, such that it would be applicable to (more) general manifolds $V$ with corners.\footnote {Relevant harmonic spinors on 
  $\hat V$ may be not $\Gamma$-invariant but interesting classes of such spinors are.}
 
The picture of    
    minimal hypersurfaces in $X$ is more transparent in this respect, where those homologous to the coordinate subtori in $X$  may originate from  $V$, namely from
 
 {\it minimal hypersurfaces  $Y\subset V$, which   separate pairs of opposite $(n-1)$-faces in $V$.}\footnote{Complete minimal subvarieties in $\hat Y\subset \hat  V$ correspond  to non-compact  singular $Y\in V$ that reflect in the codimension 1 faces alike   to   billiard trajectories in the case $dim (\hat Y)=1$.  }

In general, such  $Y$ do not exist, since they may escape the interior of $Y$ in the course of volume minimization,  but if $mean. curv(\partial_i V)> 0$ and $\angle_{ij}(V)< \frac { \pi}{2}$, then  the "boundary walls" $\partial_i$ "trap" $Y$  inside $V$.

Indeed, the first  inequality shows that, in the course of minimisation, the interior of   $Y$ can't touch  $\partial V$  by the maximum principle  and  the second one keeps the boundary of   $Y$  away from   faces  $Y$ is suppose to  separate.\footnote{In the case of non-strict inequalities, the minimal $Y$ may touch these two faces,  even coincide with one of them but the interior of $Y$ can't touch the boundary of $V$ by the  maximum principle.}

What is non-obvious here is the nature of  
 singularities at the boundary of $Y$ which may  create complications in   consecutive inductive  steps of descent  
method, even for $n\leq 7$, where there is  no singularities away from the $n-2$-faces of $V$.

Recently, however, Chao Li  [Li(rigidity) 2019] established  necessary regularity property of minimal $Y\subset V$ at the corners
of $V$ and thus gave a direct proof of {\color {blue}  $\blacksquare$} for $n\leq 7$ by Schoen-Yau's inductive decent method with minimal hypersurfaces separating pairs of opposite $(n-1)$-faces in $V$.

An  advantage of the  direct approach is  applicability   to  a class  of non-cubical manifolds $V$ with corners, which are {\it not amenable to reflections}, namely to products $V=[0,1]^{n-2}\times \hexagon\subset \mathbb R^n$, where $\hexagon\subset \mathbb R^2$  is a convex polygon.
  
  {\it no metric $g$ on such a $V$, for which the codimension 1 faces are mean convex and all dihedral angles are bounded by the Euclidean ones of $V$, can have $Sc(g)>0$.}(See  section 3.16 for more about it.)

However,  reflections  reveal a fuller picture of the geometry of $V$, not limited to minimal hypersurfaces  between opposite faces,  
 but  also including those reflected in various $(n-1)$-faces which correspond to the  minimal $\Gamma$-invariant 
 hypersurfaces in the universal orbi-covering $\hat V$ of $V$.  
\vspace{1mm}
  
 {\it  Plateau Billiard Problem.} {\sf Given a Riemannin  manifold  $V$  with (smooth or cornered)  boundary, study  minimal subvarieties in $V$  with the {\it reflection boundary condition}.}

  \subsubsection{\color {blue} MC-Normalization of Hypersurfaces with Positive 
{\color {black}M}ean {\color {black}C}urvatures and $Sc$-Normalized  Convex  Area Extremality Theorem}\label {MC-normalization+area extremality}
 
  %%%%%%%%%%%%%%%
  (5) {\sf \large Reduction of    {\color {blue}  $\blacksquare$} to { \color {blue}\EllipseSolid} and the Proof of  \color {blue}\EllipseSolid.}
   Such a  reduction, which provides an alternative proof of {\color {blue}  $\blacksquare$},  is achieved  by smoothing the corners of $V$.    Then  { \color {blue}\EllipseSolid}  is  proved by doubling  $V$ and   applying  the following.
   
   \vspace {1mm}
   
{\color {blue}\textbf {[\Ellipse}]} \hspace {1mm}\textbf{ Convex Area  Extremality Theorem.}  {\sf  Let $\underline X\subset \mathbb R^{n+1}$  be a compact smooth convex hypersurface,  let $\underline g$  be the Riemannian metric on $\underline X$ induced from the ambient Euclidean space $\mathbb R^{n+1}\supset \underline X$ and let $g$  be another Riemannian metric on  $\underline X$  with non-negative scalar curvature, $Sc(g)\geq 0$. 

Denote by $  \underline g^\circ$ and $g^\circ$ {\it normalizations} of these metrics by their respective scalar curvatures, 
$$  \underline g^\circ(\underline x)= Sc(\underline g,\underline x) \cdot \underline g(\underline x) \mbox { and } 
    g^\circ(\underline x)= Sc(g,\underline x) \cdot g(\underline x).$$
   
  (These metrics vanish, where the scalar curvatures vanish.)
  
   {\it If $n$ is even,\footnote{{\color{red!59!black}Conjecturally}, this  parity assumption is unneeded.} then there exists smooth surface $S\subset \underline X$, on which  both functions $Sc(\underline g,\underline x)$ and 
   $Sc(g, \underline x)$ are strictly positive and such that
   $$area_{g^\circ}(S) \leq area_{\underline g^\circ}(S).$$} 
   In words,\vspace {1mm}

  {\sf \color {blue!40!black} The Sc-normalization of {\it {\color {blue!40!magenta} no} Riemannian  metric} with non-negative scalar curvature on a {\it convex}  Euclidean hypersurface  can't  be {\it area wise greater} than the   $Sc$-normalized {\it original metric}   on  this hypersurface  that is  induced from the Euclidean space.}}\vspace {1mm} 
   
  This is a special case of   {\it Spin-Area Convex  Extremality Theorem}  (see  {\color {blue} [$X_{spin}{^\to}$\Ellipse]} in  sections \ref{trace-norms3}, \ref{area extremality3} that is   derived  from curvature estimates for the twisted  Schroedinger-Lichnerowicz-Weitzenboeck-(Bochner) formula
   due to  
Uwe  Goette  and Sebastian Semmelmann. Earlier, these  estimates and the issuing extremality  were established by Marcelo Llarull for the spheres $S^n$ for all  $n$, while the idea of $Sc$-normalisation, which is crucial for geometric  applications,  was suggested by Mario Listings.\footnote {See [Goette-Semmelmann(symmetric)2002],  [Llarull(sharp estimates)1998], [Listing(symmetric  spaces) 2010]  and compare with 
  \S $5\frac {4}{9}(D)$ in [G(positive) 1996]. 
Also note that     
   recently, 
    John Lott [Lott(boundary 2020] suggested a direct proof of a non-Sc  an non-MC-normalized  version of { \color {blue}\EllipseSolid}   by establishing     {\it index and vanishing theorems for Dirac operators on manifolds with boundaries}. {\color {red!40!black} Probably}, a minor adjustment of  his  argument  will deliver the full  normalized { \color {blue}\EllipseSolid}.}
\vspace {1mm}

  (6) {\sf \large Problem.} {\sf Find an elementary (whatever this means) proof of   {\color {blue}  $\blacksquare$}    in the case where $(V,g)$ admits an isometric embedding to $\mathbb R^n$. }

{\it Exercise.} Give a direct proof of   {\color {blue}  $\blacksquare$}  for {\it convex} $V\subset \mathbb R^n$.  

{\it Hint.} Show that if all $\frac {n(n-1)}{2}$  dihedral  angles at a vertex $v \in V$ are $<\frac{\pi}{2}$, then the spherical measure of the set of the  supporting hyperplanes of $V$ at $v$ is $>\frac{1}{2^n}vol(S^{n-1})$. 
\vspace {1mm}

It is known (private communication by Karim Adiprasito) that 

{\sf no convex polyhedron admits an infinitesimal deformation, which  decreases its dihedral angles but it is unknown if  a  polyhedron  $P'$  combinatorially equivalent to  $P$  may have all its   dihedral angles $\angle'_{ij} < \angle_{ij}$}. 

 {\it \large \color{red!50!black} Conjecturally}}, 
{\sf  there is no such $P'$  even among   {\it curve-linear} polyhedra with}    {\it mean convex  faces.  } 
 
 At the present time, this is confirmed for    for  special     polyhedra $P$, e.g. those  with all dihedral 
 angles $\leq \frac {\pi}{2}. $  (see section \ref{corners3}).
  
  \vspace {1mm}

(7) {\sf \large On the Proof of {\large \color {magenta}$\square$}  and {\color {magenta}\Ellipse}.}
Construction of   a  small  strictly mean convex  $V\subset X$ with rectangular corners needed for 
 {\color {magenta}  $\square$} proceeds by induction on $n$, where the resulting $V$ looks like a solid 
  $[0,l_1]\times )\times [0,l_2\times ...\times [0,l_n]$ with $ l_1>>l_2>>...'>>l_n$.
  (See [G(billiards) 2014] or do it yourself.)
  
Then the  proof of  {\color {magenta}\Ellipse}  follows
 by smoothing these corners (another exercise for the reader). 
 
 Small domains   $V\subset X$, especially for {\color {magenta}\Ellipse}, obtained this way are fairly artificial. It would be nicer to have 
 exp-images of ellipsoids from $T_x(X)$ at a point where $Sc(X,x)<0$,   or small perturbation of these  in $X$.
  
  But, probably, {\sf such a $B$ can't be a ball,  unless $X$ has  constant sectional curvature.}
  
  \vspace {1mm}

(8)  {\sf \large Normalization of  Metrics by Mean Curvatures }.  The relations
{ \color {blue}\EllipseSolid$_\leq$} and 
 { \color {blue}\EllipseSolid$_=$} for $f:\partial V\to\partial W$ becomes more transparent if the Riemannian metrics in 
the hypersurfaces $\partial V\subset V$ and in $\partial W\subset \mathbb R^n$ induced from the ambient spaces, 
call them $g$ on $\partial V$ and $h$  on $\partial W$, are  rescaled by    (the squares of) their  mean curvatures, 
denoted here  $m(v)= mean.curv(\partial V,v)$,  $v\in \partial V$  and  $m(w) =mean.curv(\partial W,w)$, $w\in \partial W,$  
$$ g=g(v)\mapsto g^\natural = m(v)^2\cdot g(v), v\in \partial V, \mbox { and  } h= h(w)\mapsto h^\natural =m(w)^2\cdot h(w),  w\in \partial W.$$
 Now, the inequality { \color {blue}\EllipseSolid$_\leq$}  says that the map $f$ is {\it distance non-increasing}
 with respect to the  {\it MC-scaled }  metrics  $g^\natural$ and  { and  }  $ h^\natural $, while { \color {blue}\EllipseSolid$_=$} expresses the isometry between  these metrics  established by $f$.

\vspace{1mm}

{\it Exercise.} Let $V\subset  \mathbb R^n$ be a convex polyhedron and $V_i\supset V$ be a  decreasing  sequence    of  larger  convex  subsets in   $ \mathbb R^n$  with smooth boundaries, which converge to $V$, i.e. 
 $$V_1 \supset V_{2}\supset ... \supset V_{i}\supset...\supset V \mbox { 
  and }
 \bigcap_{i}  V_{i}=V. $$

Describe the limit of the metrics spaces $(\partial V_i, m^2_i \cdot g_i)$,  where $m_i=m_i(v),$ 
$v\in \partial V_i$,  denote the mean curvatures of the boundaries  $\partial V_i$ and $g_i$ 
are the Riemannian metrics on these $\partial _i V$ induced from $\mathbb R^n$.

(To make it more specific, let $\partial V_i$ be (very closely) pinched between the boundaries of $ \varepsilon_i$- and 
$(\varepsilon_i+\varepsilon_i^i)$-neighbourhoods of $V$, i.e.
$$U_{\varepsilon_i}(V)\subset V\subset U_{\varepsilon_i +\varepsilon_i^i}(V),$$
where $ \varepsilon_i\to 0$  for $i\to \infty$.)

Do the same for the (induced Riemannian)  metrics $g_i$ on  $\partial V_i$ {\it normalized by their scalar curvatures},  $g_i\leadsto Sc(g_i)\cdot g_i$, and then for the other {\it symmetric functions 
$s_k$, $k=1,2,...n-1$, of the principal  curvatures}  $\alpha_1, \alpha_2,...\alpha_{n-1}$ of 
 $\partial V_i$, i.e. $g_i\leadsto s_k^{\frac{2}{k}}\cdot g_i$.

\vspace {1mm}

(10)   {\large \it \color {magenta!70!black} Problem.}  {\sf Is  there a   theory of  singular  spaces with $Sq\geq \sigma$, that is   built  on the basis of 
  { \color {blue}$\blacksquare$}, { \color {blue}\EllipseSolid} or more powerful  inequalities?}
  \vspace {1mm}
  
  {\it Example.} Let $U\subset \mathbb R^{n+1}$ be a convex subset, the principal curvatures $\alpha_i$ of the boundary $X=\partial U$  of which satisfy
$$\sum_{i>j} \alpha_i\alpha_j\geq \sigma\geq 0.$$
 at all regular points of $X$.

Albeit in general singular, $X$  can be  neatly approximated by   the $C^{1,1}$-regular boundary $X_\varepsilon $ of the $\varepsilon$ neighbourhood of $U$, where  
 the induced  Lipschitz continuous
Riemannian  metric $g_\varepsilon$  on $X_\varepsilon$ $g_\varepsilon$  be $C^0$  can be  approximated--     -- this is  obvious in the present case -- by $C^\infty$-metrics with $Sc  \geq \sigma$. 

Thus, 
(the intrinsic) path metric on $X$ {\it must share all its geometric properties with smooth manifolds 
which have $Sc\geq \sigma$}.  \vspace {1mm}

{\it Exercise.}  Show that   $X$ satisfies the inequalities   { \color {blue}$\blacksquare$}, { \color {blue}\EllipseSolid}.\footnote{I haven't done this exercise.}

  \vspace{1mm}
(11) {\it\color {blue}  Category Theoretic Perspective on  normalized Riemannian metrics.}  The above suggests that  the geometry of Riemannian manifolds 
$X=(X,g)$, where $Sc(g)>0$  without bounary is well  depicted by the {\it {\sf Sc}-normalised  metric   
$Sc(X)\cdot g$} and that  maps, which are {\it 1-Lipschitz  with respect to the   {\sf Sc}-normalised  metrics} can be taken for morphisms in the 
  category of manifolds with $Sc>0$. 

Now, If $X$ does have a  boundary and this boundary is mean convex, the normalization of 
$X$ by $Sc(X) $  and of $\partial X$ by $mean.curv (\partial X)$  do not agree  on $\partial X$.

Alternatively, one may use positivity of the mean curvature of $\partial X$ for blowing up the metric 
of $X$ near $\partial X$  keeping $Sc>0$, as it is done  in {fill-ins3}.

\vspace {1mm}
  (11) {\sf {\color {red!50!black}  Polyhedral Localization Conjecture}}. Let $\underline V\subset \mathbb R^n$
  be a convex polyhedron an let 
  
 $V$ be a compact oriented  Riemannian $n$-manifold with corners,  of "combinatorial class" of   $\underline V$, which means there exists a {\it proper corner continous map} $f:V\to \underline V$
 of degree one, where "proper corner" means "face respecting": the  $ (n-1)$ faces $V_i\subset \partial V$ ae equal the pullbacks of the corresponding faces $\underline V_i\subset \partial \underline V$.
 
 Let the boundary of $V$  be  { \it more  mean convex}  than that of $\underline V$, i.e.
  all $ (n-1)$-faces $V_i\subset \partial V$   have non-negative mean curvatures. 
  and the dihedral angles of $V$ along the $ (n-2)$-faces  are strictly bounded by the corresponding angles in 
 $\underline V$,
 $$\angle(V_i,V_j)\leq \angle(\underline V_i, \underline V_j).\footnote{This   strict "<", rather than more natural  "$\leq $", is used here to avoid possible technical complications  with the rigidity problem (see sections \ref{corners3}  and \ref{dihedral4.})}$$
 
 Then, {\sf {\color {red!50!black} conjecturally}, 
 
{\color {blue} [{\large $\boxbox$}]} {\sf $V$ contains domains $U_\varepsilon \subset V$  with corners, which have   {\it arbitrarily small diameters}, 
 $$diam(V_\varepsilon)\leq \varepsilon>0,$$ 
 which  are  also in the combinatorial class of  $\underline V$  and  such that  their boundaries are {\it more mean convex} than that of $\underline V$.}\vspace {1mm}
 
{\sf \textbf {"Cubical"  Remark/Theorem}}}. The  {\color {blue}$ \blacksquare$}-{inequality} 
together with its  contraposition {(\large \color {magenta}$\square$}), that is the existence 
of cubical mean convex domains with acute dihedral angles in manifolds with $Sc<0$,
(see section\ref{Sc-criteria3}) imply the validity of {\color {blue} [{\large $\boxbox$}]} for  cubical, 
$\underline V=[0,1]^n$ and $V$ {\it homeomorphic} to $\underline V$.

Now, a close look at the proofs of the  {\color {blue}$ \blacksquare$} in  section \ref{reflection3} the  {\color {blue}$ \blacksquare$}  (also see {\ref{corners3}  and  \ref{dihedral4})  apply  to more general $V$:

The proofs with the Dirac operator works for {\it all spin} manifolds $V$, while  the calculus of variation methods needs no spin, but it
become cumbersome for $n\geq 8$ and especially so for $n\geq 9$, due to possible singularities of minimal hypersurfaces of dimension$\geq 7$ (see section \ref {singularities3}  for more about it.)

However, these proofs, especially the Dirac-operator theoretic one, do not directly pinpoint the 
 small cube with acute angles in $V$, as they also need the local property  {(\large \color {magenta}$\square$}) of $Sc<0$.\footnote{Our  conjecture  doesn't  mention  any curvature and we want the proof to be also like that.}

In fact, one can construct  $U_\varepsilon$,  at least for $n\leq 8$, to get such a cube $U_\varepsilon$ arguing with minimal hypersurfaces, roughly as follows.

Given an {\it admissible} $U\subset V $, i.e. a  mean convex cubical domain with acute dihedral angles, let us push  the  faces  $U_i$ in one after another     little by little   keeping  $mean. curv. 0$ and $\angle_{ij}\leq\frac{\pi}{2}$.
 
This process stops when we arrive at some $U=U_{min}$, where all faces a (locally) volume minimizing with free 
boundaries on the unions of th remaining faces. Now one can slightly move each  face, say $U_1$ in, 
$U_1\leadsto U_{1,\delta}$ 
keeping the dihedral angles equal $\frac {\pi}{2}, $ but now such that  $U_{1,\delta}$ is 
everywhere  {\it mean concave} rather  than mean convex. Thus, we obtain a smaller  admissible 
domain, say $U'\subset U=U_{min}$ namely the band between $U_1\leadsto U_{1,\delta}$ in  $U$.

If $n\leq 8$ one can indeed arrange  this process  to arrive at an $\varepsilon$-small "cube"  $U=U_\varepsilon$, 
but, in general, this "cube" may  be only as small as the singularities of the minimal hypersurfaces are:  the  "cube minimization process" can, a priori,   converge to an $(n-8)$-dimensional closed subset $U_\bullet \subset  V$. \vspace {1mm}

{\sf \large\color {green!40!black}  About Dimension  $n=9$.}
{\sl If $n=9$, the domains $U_\varepsilon$,
are {\it spin}. }

Indeed, $U_\varepsilon$  are localized near  a  subset  $U_\bullet$ of dimension $\leq 1$, while the obstruction to spin, that is the second Stiefel-Whitney number $w_2\in H^2(V; \mathbb Z_2)$, resides in dimension 2.

At this point,  the Dirac theoretic argument applies  to $U_\varepsilon$  and, together with   {\large \color {magenta}($\square$)}, yields {\color {blue} [{\large $\boxbox$}]} for cubical $\underline V$ with   $dim(\underline V)=9$.  (see [G(billiards) 2014]).

  \vspace {1mm}
 %%%%%%%%%%%%%%%%%%%%%%
  
 \subsubsection {\color {blue}   $C^0$-Limits  of Metrics $g$ with $Sc(g) \gtrless \sigma$ }\label {C0-limits3}

 %%%%%%%%%%%%%%%%%%%%
 Let $X$ be a smooth Riemannian manifold,
 let  $G=G(X)$ the space of $C^2$-smooth Riemannian  metrics $g$  on $X$ and let 
 $G_{Sc\geq \sigma}\subset G$ and  $G_{Sc\leq \sigma}\subset G$, $-\infty < \sigma<\infty$,
be the subsets of metrics $g$ with $Sc(g)\geq \sigma$ and with $Sc(g)\leq \sigma$ respectively.

Then:\vspace{1mm}

{\color {blue}  \large A}:   \textbf {$C^0$-Closure Theorem.}  {\sf The subset   $G_{Sc\geq \sigma}\subset G$ is closed in $G$ with respect to the $C^0$-topology:} 

\vspace {1mm}

{\it uniform limits  $g=\lim g_i$ of metric $g_i$ with $Sc(g_i)\geq \sigma$ have $Sc\geq \sigma$,
provided these $g$ are $C^2$-smooth in order to have their scalar curvature defined.}

\vspace{1mm}

{\color {blue}  \large B}:  \textbf {$C^0$-Density Theorem.} {\sf The subset   $G_{Sc\leq \sigma}\subset G$ is dense in $G$ with respect to the $C^0$-topology}.

Moreover, 
{\it all $g\in G$ admit  fine (which is stronger than uniform for non-compact $X$) approximations by metrics with scalar curvatures $\leq \sigma$.}\vspace{1mm}

{\it Short  Proof of} {\color {blue}  \large A.} Let us show that violation of  { \color {blue}\EllipseSolid} for a smooth metric $g$ on a manifold $X$, that is ({\color {magenta}\Ellipse}) from the previous section,  implies this  for $g_\varepsilon$  for sufficiently small
$\varepsilon=||g-g_\varepsilon||_{C^0}$. 

Indeed let the boundary $\partial V$ of a compact strictly mean convex  domain $V\subset X=(X,g)$ admits a smooth map $f$ of degree one to the boundary of a convex $W\subset \mathbb R^n$,  the norm of the differential of which satisfies:
$$||df(v)||<  \frac {mean.curv(\partial W,f(v))}{mean.curv(\partial V,v)}.$$
If $g_\varepsilon$ is close to $g$, then there exists a smooth $V_\varepsilon\subset X$, the boundary of which is {\it $\delta$-close} to 
$\partial V$ and its 
$g_\varepsilon$-mean curvature is {\it $\delta$-close} to the $g$-mean curvature of $\partial V$,  where 
$\delta \to 0$ for $\varepsilon\to 0$, and where  "$\delta$-close"  means the following.
diffeomorphisms exists a smooth  $(1+\delta)$-Lipshitz map\footnote{A map  between metric spaces, $f:X\to Y$,  is
 {\it $\lambda$-Lipschitz} if   
$dist_Y(f(x_1),f(x_2))\leq \lambda dist_X(x_1,x_2)$);   a   $\lambda$-Lipschitz map  is  {\it $\lambda$-bi-Lipschitz} if it is one-to-one and the inverse map is also  $\lambda$-Lipschitz.}
$ \nu: \partial V_\varepsilon\to \partial V$, i.e. $||d\nu||\leq 1+\delta$, such that 
$$dist_g(x, \nu(x))\leq \delta  \mbox { for all } x\in \partial V_ \varepsilon$$
as well as 
 $$|mean.curv_{g_\varepsilon} (\partial V_ \varepsilon, x)-  mean.curv_{g} (\partial V, \nu(x))|\leq \delta .$$

In fact, one can take the $g$-normal projection of the $\delta$-neighbourhood of $\partial V$ to $\partial V$
restricted to $\partial V_ \varepsilon$ for this  $\nu$, where, observe, this projection $\partial V_ \varepsilon\to \partial V$, albeit {\it not necessarily  
a diffeomorphism }  for small $\varepsilon\to 0$,  can be {\it $C^0$-approximated by diffeomorphisms.} 
\footnote{The existence of such a $V_\varepsilon$ and its properties must be a common knowledge among experts on the geometric measure theory  but I couldn't find a reference and written down a proof in section 10.2 of  [G(Hilbert) 2011].}\vspace{1mm}

{\it About  Alternative Proofs.} Instead of   { \color {blue}\EllipseSolid}  one can   use     {\color {blue}{ $ \blacksquare$} } but the available   argument in   [G(billiards) 2014]) is unpleasantly  convoluted. 

A streamlined  proof based on {\it Hamilton-Ricci flow}  was suggested by 
Richard Bamler and 
 further developed by Paula Burkhardt-Guim  who has shown, in particular, that \vspace {1mm}

{\color {blue} (\Large $\star$}) {\it if a continuous  metric $g$ on a smooth manifold $X$ admits a $C^0$-approximation by metrics $g_\varepsilon $ 
with $Sc(g_\varepsilon)\geq 
\varepsilon\to 0$, then $X$ admits a smooth metric with $Sc\geq 0$.}

 Moreover,

{\it $g$ can be $C^0$ approximated by metrics with $Sc\geq 0.$}

 \vspace{1mm}
 
 Thus,    
 
 {\sl continuous  metrics which are $C^0$-limits of of smooth metrics metrics $g_i$ with  
$\lim Sc(g_i)\geq -\varepsilon\to 0$  have the same kind of geometries as metrics 
with $Sc\geq 0$.}

 \vspace{2mm}

 {\it \color {blue!60!black} Question.} { \sf Do {\it Lipschitz metrics}\footnote{ A measurable Riemannian metric $g$ 
 on a smooth $n$-manifold $X$ is {\it Lipschitz} if it is locally   {\it bi-Lipschitz}  equivalent 
 to  the Euclidean metric on (a domain in) $\mathbb R^n$, see  [Ivanov(Lipschitz)  2008].  Notice that the natural domains  $X$ for such $g$ are {\it  Lipschitz}, rather than smooth, manifolds that are defined by   {\it bi-Lipschitz  atlases} on $X$,  see [NS{Lipschitz) 2007]. }}
 are similar to continuous ones in this respect  for suitable limits $g_i\to g$?}

\vspace{1mm}\vspace{1mm}

{\it About}{ \color {blue}  \large B}. This is a special case of   {\it the curvature $h$-principle}  discovered  by Joachim Lohkamp,\footnote {[Lohkamp(negative  Ricci curvature) 1994].}  whose proof  in   depends on a
(more or less) direct, yet,  elaborate, geometric construction, which  also    shows that \vspace{1mm}

{\color {magenta!60!black} (\Large $\star$}) {\it the  metics with $Ricci<0$ are 
$C^0$-dense in he space of all Riemannian metrics on $X$.}\vspace{1mm}

(If, contrary  to   {\color {blue}  \large A}, the space of metrics with $Sc\geq 0$ were dense, 
there would be no hope for a non-trivial  geometry of such metrics.)

\vspace{1mm}

 \vspace {1mm}

(The  $C^0$-closure theorem  for the scalar curvature looks  similar  to \vspace {1mm}

  {\it \color {blue!59!black} \textbf{Eliashberg's  $C^0$-Closure Theorem}}, which claims that   \vspace {1mm}
  
  {\it  $C^0$-limits  of   symplectomorphisms,
  are again symplectomorphisms, provided  they 
  
  are  $C^1$-smooth and $C^1$-invertible. }\vspace {1mm}

 But, unlike  $Sc\geq 0$,  non-smooth  such limits are significantly  {\it more flexible} and geometrically {\it  less constrained}  than  smooth symplectomorphisms\footnote  {[Buhovsky-Opshtein($C^0$-symplectic)2014], [Bu-Hu-Sey($C^0$ counterexample)  2016]; yet,  some    symplectic  geometry, if  properly understood, passes the $C^0$-barrier [Bu-Hu-Sey($C^0$-symplectic)  2020]. }

   \vspace {1mm}

 {\sf \large Weak convergence of metrics and convergence of manifolds.} Besides uniform convergence, 
 there are other metric conditions on    sequences  of metrics that preserve positivity of the scalar curvature in the limit, where the simplest unkown case is the following.   
 
Let smooth  Riemannian metrics $g_i$ {\it converge in measure} to an also  smooth  $g$,  i.e. the measure of the subset,  where 
the $|g(x)-g_i(x)|\geq \varepsilon $  tends to $0$ for $i\to \infty$. \vspace {1mm}

 {\it Do the inequalities  $Sc(g_i)\geq 0$, imply that  also $Sc(g)\geq 0$}?\vspace {1mm}
  
   This is most likely to hold true  
  if the {\it Lipschitz distance}\footnote{This is  the maximum of the Lipschitz constants of the identity map
  $V\to V$ with respect the pairs of metrics, $(V, g)\to(V,g_i)$ and  $(V, g_i)\to(V,g)$.}
 between $g$ and   $g_i$ remains  bounded by a constant independent  of $i\to\infty$. \footnote {Beware of examples implied by theorem 1.4 in [Brun-Han(large and small) 2009]).} \vspace {1mm}

  {\it In  geometry}, however, natural limits are not these of metrics but those of     {\sf   Riemannian  manifolds,    
with   no fixed topological   background,}
 $$X_i=(X_i,g_i)\to X=(X,g),$$
where,  relevantly for us, such     limits, even when  drastically changing topology, may  preserve positivity of the scalar curvature;   yet,   some natural geometric  limits, e.g.  the {\it Hausdorff} and 
{\it intrinsic flat} ones may uncontrollably change  scalar curvature.\footnote {
See sections \ref {C0-limits3}, \ref 
{weak limits6}  for  examples (and counter examples),  of various kind of behaviour of the scalar curvature under  such convergence.}

     \vspace {1mm}

{\it  \textbf {{\color {red!50!black} Conjecture}}: Quantification of $C^0$-Convergence.} {\sf Let $g_0$  and $g_\varepsilon $ be smooth Riemannian metrics on the ball $ B^n$, such that the $C^0$-norm of the difference $g_\varepsilon -g_0$  is bounded by $\varepsilon$,  or  (almost) equivalently  the identity map 
 $(B^n, g_0)\to(B^n,g_\varepsilon)$ 
is $\left (1+\frac {\varepsilon}{2}\right)$-bi-Lipschitz.}

{\sl  Then there exist positive constants (large) $c_0>0$ and (small)  $\varepsilon_0> 0$, which depend only on $g_0$,  such that if $\varepsilon\leq \varepsilon_0$, then  the scalar curvature of $g_\varepsilon$ at the center of the ball satisfies,
$$ Sc(g_0(0))\geq    \inf_{x\in B^n}(Sc(g_\varepsilon(x)))  -c_0\varepsilon.$$}

{\it Motivating Example.} If  $g_\varepsilon=(1\pm\varepsilon)^2g_0$, then  
$||g_\varepsilon-g_0||=2\varepsilon +o(\varepsilon)$ and $|Sc(g_\varepsilon)-Sc(g_0)|=O(\varepsilon). $ 

\vspace{1mm}

{\it Exercise.}  Prove this conjecture for $n=2$.\vspace{1mm}

{\it Remark.}  {\color {red!40!black} Probably}, a close look at  the   proof of { \color {blue}  A}  will  yield the conjecture for 
radial (i.e. $O(n)$-invariant)  metrics $g_0$ (compare with approximation corollary in \S$5\frac {5}{6}$ from [G(positive)1996]  as well as the inequality $Sc(g_0, 0)\geq    \inf_{x\in B^n}(g_\varepsilon)  -c_0\varepsilon^{\frac {1}{2}}$.

%%%%%%%%%%%%%%%%

\subsection {\color {blue} Spin Structure, Dirac Operator, Index Theorem,   $\hat A$-Genus, $\hat \alpha$-Invariant and   Simply Connected Manifolds with and without $Sc>0$}\label {spin index3}
%%%%%%%%%%%%%%%%%%%%%%%%%%%%%%%%%%%%%%

Let $L\to X$ be a  real orientable   vector bundle of rank $r$  and  $F\to X$ be the oriented frame bundle  of $L$. If $r\geq 2$ the fundamental group of the fiber $F_x=SL(k)$ is infinite cyclic and if $k\geq 3$ this group is cyclic of order 2. In both  cases, $F$ comes with  a canonical double cover $\tilde F_x\to F_x$.\vspace {1mm}

{\sf The bundle $L$ is called {\it spin}, if $\tilde F_x\to F_x$  extends to a double cover $\tilde F\to F$,

and  smooth  orientable manifold  $X$ is {\it spin} if its tangent T(X) bundle is spin.}\vspace {1mm}
  
 Extension of the covering $\tilde F_x\to F_x$, if it exists, is, in general, non-unique. In the case of 
 of $L=T(X)$ such an extension  is called a  {\it spin structure} on $X$.

 When you speak of spin, it is common in geometry and for a good reason,  to reduce the structure group of $L$  from $SL(r)$ to $SO(r)\subset  SL(r)$  and to deal with the orthonormal frame bundle $OF\to X$  instead of $F$, where the  double cover group $\tilde SO(r)=OF_x$ is called  
  {\it spin group} $Spin(r).$ \vspace {1mm}

 {\it Example.} The tangent bundle of the 2-sphere is spin, but the Hopf bundle over $S^2$ is not, since 
$OF$,   that is  $S^3$  for the Hopf bundle, is simply connected.

Similarly -- this an exercise in elementary topology, \vspace {1mm}

 {\sf   an  oriented bundle $L$ of rank two over an oriented surface $X$ is spin if and only if  its {\it Euler class}, that is  the {\it  self-intersection number of $X\subset L$} is {\it even}; if  $X$ is non-orientable, then $L$ is spin if the {\it  second  Stiefel-Whitney class}, that is the  self-intersection number mod 2 of  $X\subset L$ 
{\it vanishes}. In either case $L$ is spin if and only if the Whitney sum of $L$  with the trivial line bundle $l\simeq X\times \mathbb R^1$ is trivial,  $L\oplus l\simeq X\times\mathbb R^3.$ }\vspace {1mm}
In general,\vspace {1mm}

{\sf  a bundle $L$ over a manifold $X$ of dimension $n\geq 3$ is spin, if and only if its restriction   to all surfaces in $X$ is spin, which is again equivalent to the vanishing of  the 
second  Stiefel-Whitney class $w_2(L)$.}\footnote{ The value of  $w_2(L) \in H^2(X;\mathbb Z_2)$ on a homology class $h\in H^2(X;\mathbb Z_2)$ is, almost by definition, equal to zero if and only if the restriction of $L$ to surfaces in $X$ that represent $h$ is trivial. 

 Geometrically, the double
cover $\tilde F_x\to F_x$ extends to $F$ over the complement to a subvariety  $\Sigma\subset X$ of codimension two, the homology class of which is Poincare dual to $w_2(X)$. This $\Sigma\subset X$
is waht stands on the way of applying Dirac theoretic methods to non-spin manifolds.}\vspace {1mm}

{\it Half-spin Bundles.} There exit two (remarkable) irreducible unitary representations of the group $Spin(r)$ for  $r=2k$ of complex dimensions $2^{k -1}$, say  $S^{\pm}(r)$. 
 Accordingly, Riemannian spin manifolds, (i.e. with spin structures on them)   $X$ support  two   $Spin(n)$   bundles $\mathbb S^\pm$ with the fibers $S^{\pm}(r)$ that are associated with principal spin bundle $\tilde SO\to X$ for the double covering representing the spin structure on $X$.
 We   let 
$\mathbb S=\mathbb S^+\oplus \mathbb S^-$ and call this $\mathbb S$ the {\it spin bundle.}\footnote {
In reality,   $\mathbb S$ comes first and then  splitting    $\mathbb S^-\oplus\mathbb S^+$ follows,
see section \ref{Clifford3}.}
 \vspace {1mm}

The Dirac  operator
$$\mathcal D:  C^\infty(\mathbb S)\to  C^\infty(\mathbb S)$$
is  a first order differential operator constructed in  a canonical geometrically invariant way
universally applicable to all $X$   (see section ).

This is an {\it elliptic selfadjoint} operator, which interchanges  $C^\infty(\mathbb S^+)$ and  $C^\infty(\mathbb S^-)$
where the  operators 
  $$\mathcal D^+:C^\infty(\mathbb S^+) \to C^-(\mathbb S^-)\mbox  { and } \mathcal D^-:C^\infty(\mathbb S^-) \to C^-(\mathbb S^+)$$
   are mutually adjoint.

We  explained  already in section \ref{Dirac1}  how, following Lichnerowicz, that the   {\it \textbf {Atiyah-Singer index theorem}} for the Dirac's  $\cal D$ and    the S-L-W-(B) identity 
$${\cal D}^2=\nabla^2+\frac {1}{4}Sc,$$ 
imply that \vspace {1mm}

{\sf there are smooth {\it closed simply connected} manifolds $X$ of all dimensions $n=4k$, $k=1,2,...$, that admit no metrics with $Sc>0$.}\vspace {1mm}

The simplest example of these for  $  n=4$ is 
 {\it the Kummer surface} $X_{\sf Ku}$  given by the equation 
 $$z_1^4+z_2^4+z_3^4+z_4^4=0$$
 in the complex projective space $\mathbb CP^3$.

 In fact,    all complex surfaces of {\it even degrees} $d\geq 4$  as well as their Cartesian products,  e.g      $X_{\sf Ku} \times ... \times X_{\sf Ku}$ {\it admit no metrics with $Sc>0$.}

Also we know that the  Atiyah-Singer $\mathbb Z_2$-index theorem   of 1971  allowed an extension of  Lichnerowicz'  argument  to manifolds of dimensions $8k + 1$ and  $8k + 2$, e.g. to exotic spheres  in    

{\it \textbf {Hitchin's theorem}}:
{\sf there exist manifolds $\Sigma$ homeomorphic (but no diffeomorphic!) to the spheres $S^n$,  for all 
$ n = 8k + 1, 8k + 2$,  $ k=1,2,3...$, which admit  no metrics with $Sc>0$.}
 \vspace{1mm}

 (What makes the {\it differential}  structures  of Hitchin's  {\it topological} spheres  $\Sigma$ incompatible with   $Sc>0$ is that to these   $\Sigma$  are {\it not boundaries of spin manifolds.})
 
  \vspace{1mm}

 %Notice  that, by  {\it Yau's  solution of the  Calabi conjecture,}  the Kummer surface admits a metric with $Sc=0$,  even with $Ricci=0$, but, probably, (I guess this must be known)  there is no   metrics with $Sc=0$ on  Hitchin's exotic spheres    $\Sigma$.

 \vspace {1mm} 
The actual  Lichnerowicz-Hitchin theorem says that if a certain  {\color {blue!80!red}topological invariant $\hat \alpha  (X)$} {\it doesn't vanish},  then  $X$ {\it admits no metric with $Sc>0$},  since, by the Atiyah  and Singer index formulae,\footnote{ The Dirac  operator is defined only on {\it spin} manifolds;    we postulate  at the  present moment  that $\hat \alpha  (X)=0$ for non-spin manifolds $X$. 

(In fact, if $n=dim(X)=4k$, this   $\hat \alpha  (X)$ is a certain linear combination of the {\it Pontryagin numbers} of $X$,  called {\it $\hat A$-genus} and denoted   $\hat A[X]$.

Accidentally, since all {\it compact  homogeneous} spaces $X=G/H$, except for tori, support  metrics with $Sc>0$,  Lichnerowicz'  theorem says that they {\it either non-spin or 
 $\hat A[X]=0$.})} 

$$\hat \alpha  (X)\neq 0\Rightarrow Ind({\cal D}_{|X})\neq 0 \Rightarrow\mbox {$\exists$ {\sf   harmonic spinor   $\neq 0$ on}      $   X$}, $$
which is incompatible with the identity ${\cal D}^2=\nabla^2+\frac {1}{4}Sc$ for $Sc(X)>0$

Conversely, \vspace{1mm}

 {\color {blue} \HandRightUp} {\sf  if $X$ is a {\it simply connected} manifold of {\it dimension $n\geq 5$},
and if $\hat \alpha  (X)=0$, then, 
 an  application of "thin surgery"} (see section \ref{thin1})  {\sf to suitably chosen   generators   $O(n)$- and $Sp(n)$-  cobordism groups   in dimensions $n\geq 5$, where these generators carry  metrics with  $Sc>0$,  yields\footnote {[GL(classification) 1980], [Stolz 1992].} that  $X$ admits a metric with positive scalar curvature.} 

Thus, for instance 

{\it all simply connected manifolds of dimension $n\neq 0,1,2,4 \mod 8$     admit metrics with $Sc>0$,}\footnote{If $dim(X) =3,$  this follows from Perelman's  solution of the Poincar\'e' conjecture.}
since  $\hat \alpha  (X)=0$ is known to vanish for these $n$.\footnote {As far as the exotic spheres $\Sigma$ are concerned, these $\Sigma$
admit metrics with $Sc>0$ if and only if $\hat \alpha(\Sigma)=0$, i.e. if $\Sigma$ bound spin manifolds, which directly follows by the codimension 3 surgery of manifolds with $Sc>0$ described in [SY(structure) 1979] and in [GL(classification) 1980]. Moreover, many of these $\Sigma$, e.g. all 7-dimensional ones, admits metrics with {\it non-negative  sectional curvatures} but  the full extent of "curvature  positivity"  for exotic spheres  remains problematic (see [JW(exotic) 2008]  and references therein.}

 \vspace {1mm}

 \vspace {1mm}

{\it Topology of   Scalar Flat.}  By {\it Yau's  solution of the  Calabi conjecture,}  the Kummer surface admits a metric with $Sc=0$,  even with $Ricci=0$, but there is  \vspace {1mm}
 
 \hspace{14mm}{\it no   metrics with $Sc=0$ on  Hitchin's exotic spheres    $\Sigma$.}  \vspace {1mm}

In fact,    \vspace {1mm}

{\sl if  a compact simply-connected scalar-flat  manifold $X$ of dimension
$\geq 5$ admits no metric with $Sc>0$,{\footnote  {These $X$  are Ricci flat, [Bourguignon (these)  1974], [Kazdan[complete 1982].} then there are cohomology classes $\alpha\in H^2(X)$ 
and $\beta\in H^4(X)$,  
 such that 
$$\langle \exp \alpha \smile \exp \beta \smile p_1(X)\rangle\neq 0,$$
where $p_1(X)$ is the first Pontryagin number, }}
 [Futaki(scalar-flat)1993], [Dessai(scalar flat) 2000]. \vspace {1mm}

And if $X$ is non-simply connected then \vspace {1mm}

{\it a finite covering of $X$  isometrically splits  into the product of a flat torus and 

the above kind simply connected manifold,} \vspace {1mm}

\hspace {-6mm}as it follows from  Cheeger-Grommol  splitting
 theorem + Bourguignon-Kazdan-Warner perturbation theorem.
 \vspace {1mm}

{\it \color {blue}  A Few Words on  $n=4$ and on $\pi_1\neq 0$.} If $n= 4$ then, besides   vanishing of  the $\hat \alpha $-invariant (which is equal to a non-zero multiple of the   first Pontryagin number for $n=4$),  positivity of the scalar curvature also  implies  the vanishing of the {\it Seiberg-Witten invariants} (See   lecture notes by Dietmar Salamon,  [Salamon(lectures) 1999]; also we say more about it in section  \ref {Seiberg3}).

\vspace {1mm}

If $X$ is a closed spin manifold of dimension $n\geq 5$ with the fundamental group $\pi_1(X)=\Pi$ ,
then, again by an application of the  thin surgery,  \vspace {1mm}

{\sl the existence/non-existence of a metric $g$ on $X$ with $Sc(g)>0$ is an invariant of the spin bordism class $[X]_{sp}\in bord_{sp}$(\sf B{$\Pi)$ in the classifying  space {\sf B}$\Pi$,}}\vspace {1mm}

 \hspace {-6mm} where, recall,  that (by definition of "classifying")
the universal covering of  {\sf B}$\Pi$ is contractible and  $\pi_1$({\sf B}$\Pi)=\Pi$. \footnote {See lecture notes   [Stolz(survey) 2001].} \vspace {1mm}

There is an avalanche of papers, most of them  coming under the heading of "{\sf Novikov  Conjecture}", with various criteria  for the class  $[X]_{sp}$, and/or for the corresponding homology class $[X]\in H_n$ ({\sf B}$\Pi)$ (not) to admit  $g$ with $Sc(g)>0$ 
on manifolds in this class, where these criteria usually (always?)  linked to  generalized  index theorems for  twisted Dirac operators on $X$  with several levels of sophistication in arranging this "twisting".

Yet, despite the recent progress in this direction for dimensions 4 and 5\footnote {See   [Chodosh-Li(bubbles) 2020],  [G(aspherical)  2020] and section \ref {5D.3}}  proving/disproving the following for $n\geq 4$ remains  beyond the present day means. 
%this  remains unknown even for spin manifolds   of 
%dimension $\geq 4$.
\footnote{The case $n=3$ follows  from the topological classification of compact 3-manifolds $X$ with positive scalar curvature {\it these are connected sums of quotients of spheres $S^3$ and products $S^2\times S^1$ by finite isometry groups} 
[GL(complete) 1983],[Ginoux(3d classification) 2013].)} %proving/disproving the following remains  beyond the present day means.

 \vspace {1mm}

  ({\color {green!70!black} Naive?}) \textbf {{\color {red!50!black} Conjecture}}.\footnote {This, 
as many other our conjectures,  is  based on a  limited class of examples with no 
 idea of  where to look
 for counter examples.}
{\sf \color {green!30!black}If a closed  oriented $n$-manifold $X$ admits a continuous map to an {\it aspherical  space {$\sf B$,}}\footnote { That is the universal covering of {$\sf  B$} is contractible, hence,   ${\sf B}$ is ${\sf B}(\Pi)$  for $\Pi=\pi_1({\sf B})$.}
such that the image of
the rational fundamental homology class of $[X]_{\mathbb Q}$ in  the rational homology\footnote {Bernhard Hanke pointed out to me that  non-vanishing of this image  in homology with {\it finite coefficients}, e.g. for finite groups $\Pi$,  may also prohibit  $Sc>0$,  but this remains obscure even on the  level of conjectures.}  $ H_n({\sf B};\mathbb Q)$ doesn't vanish, then $X$ admits no metic $g$ with $Sc(g)>0$.}

\vspace {1mm}

(We shall   describe the status of this problem   together with the {\it Novikov conjecture} in section \ref 
{Novikov3}.)

%%%%%%%%%%%%%%%%%%%%

\subsection {\color {blue}Unitary Connections, Twisted Dirac  Operators  and Almost Flat    Bundles Induced by   $\varepsilon$-Lipschitz  Maps }\label {twisted3}
%%%%%%%%%%%%%%%%%%

We turn now to twisted Dirac operators $\mathcal D_{\otimes L}$ that act on  tensor products $\mathbb S\otimes L$ for vector bundles $L\to X$ with linear (most of the time, unitary)  connections $\nabla$.

One can think of such  a  $\mathcal D_{\otimes L}$  as an {\it infinitesimal family} of 
$\cal D$-s parametrized by $L$, where the action takes place along $\mathbb S$ with no differentiation in the $L$-directions.

For instance if  $L=(L, \nabla)$  is a trivial flat bundle, $L=X\times L_0$, where $L_0$  is a vector space (fiber),  then  
$C^\infty(\mathbb S\otimes L)=C^\infty(\mathbb S)\otimes L_0$
 and the   $\mathcal D_{\otimes L}$
doesn't act on $L$ at all: 
$$\mathcal D_{\otimes L}(f\otimes l)=\mathcal D(f)\otimes l\mbox {   for all vectors } l\in L_0.$$

In general, the   $\mathcal D_{\otimes L}$ differs from that in  the flat case by a zero order term, which is, bounded by the curvature of $L$ and, strictly speaking, is  defined only locally, where the bundle $L$ is topologically trivial. But exactly  this impossibility of global comparison of 
  $\mathcal D_{\otimes L}$ on   $C^\infty(\mathbb S\otimes L)$  with  $\mathcal D$  on $C^\infty(\mathbb S)\otimes L_0$ creates a correction term in the index   formula.

  This correction, unlike the background operator $\cal D$,  carries no  subtle topological information 
  about $X$, such as $\hat A(X)$ for $n=4k$, which  is not a homotopy invariant for $n>4$
  and even less so about $\hat \alpha(X)$  for $n=8k+1, 8k+2$, which is not even invariant under  p.l. homeomorphisms  and which is far removed from anything even  remotely, geometric about   $X$, while the 
  topology (Chern classes) of $L$   reflects the  area-wise  size of the metric $g$ on $X$, which, in turn, influences 
   {\it homotopy theoretic properties} of $X$ linked to the fundamental  group. 
   
   The following  definition gives  you  a fair idea of what  kind of properties these are. 
  
  \vspace {1mm}

{\it Profinite  Hypersphericity.}
 A  Riemannian $n$-manifold $X$ is {\it profinitely  hyperspherical} if

{\sf given an   $\varepsilon>0$, there exists an orientable   finite covering $\tilde X=\tilde X_\varepsilon$,  
which admits

 an {\it $\varepsilon $-Lipschitz} map between \footnote {A map between metric
  spaces, $f:X\to Y$,  is 
$\varepsilon $-Lipschitz if $dist_Y (f(x_1), f(x_2))\leq \varepsilon dist_X (x_1,x_2)$ for all $x_1,x_2\in X$. For instance, "1-Lipschitz" means "distance non-increasing". $\varepsilon $-Lipschitz  for smooth maps $f$ between Riemannian manifolds is equivalent to $||d(f(x)||\varepsilon$,  $x\in X$.} $\tilde X\to S^n$ of {\it non-zero} degree.}

 This property of  {\it compact} manifolds  (the definition of this  hypersphericity extends too open manifolds)
  doesn't depend on the Riemannian metric on $X$.  
 Moreover

{\sf If $X_1$ is  profinitely  hyperspherical and $X_2$ admits a map of non-zero degree to $X_1$ 
then, obviously,  $X_2$ is also profinitely  hyperspherical;  in particular, this property is  a {\it homotopy invariant}   of $X$.}\vspace {1mm}
 
 {\it Example.}   {\sf Manifolds   $X$, which admit  {\it locally expanding self-maps} $E:X\to X$,
 e.g.   the $n$-torus 
 $\mathbb T^n$, where
  the endomorphism   $t\mapsto Nt$ locally  expands the metric by  $N$, 
are   profinitely  hyperspherical. }\vspace {1mm}

Indeed,  such an $E$ defines a {\it globally} expanding homeomorphism,  call it $\hat E$,     from  $X$ onto a finite covering $\tilde X=\tilde X(E)$,  where the inverse map 
$\hat E^{-1}:\tilde X\to X$  contracts as much as $E$ expands.

 Therefore, the covering corresponding to the  $i$-th iterate  of $E$ 
comes with an  $\varepsilon_i$-Lipschitz map to $X$,  where   $\varepsilon_i\to 0$ for $i\to\infty$ and compositions of these with a map $X\to S^n$ of non-zero degree also have  $deg\neq 0$, while  their Lipschitz constants go to zero.\footnote{Further   examples of this phenomenon  and issuing topological obstruction to $Sc > 0$     for manifolds with {\it residually finite} fundamental  groups  are given  in  [GL(spin) 1980] under the heading of  "enlargeability".     Since the residual finiteness condition  was eventually  lifted,
this terminology now applies to a  broader  class of  manifolds, including spaces $X$  the {\it universal covers} of which admit contracting self-maps of positive degrees,  see section\ref{negative sectional4}}
\vspace {1mm}

Now, if you recall Atiyah-Singer   index theorem for the twisted Dirac  operator and  $\mathcal D_{\otimes L}$ and the  (untwisted) 
  S-L-W-(B)
formula $\mathcal D^2=\nabla \nabla^\ast +\frac {1}{4}Sc$\footnote{This $\nabla$ stands for the Levi-Civita connection in the spin bundle.}  you arrive  at the following. 

\vspace {1mm}

{\color {blue}$\mathbf {[Sc\ngtr 0]}$}: \textbf  {Provisional 
Proposition}.\footnote {This will be significantly generalized later on.}{\it Compact orientable\footnote{If $X$ is non-orientable, take  an oriented double cover of it.} profinitely  hyperspherical  spin manifolds $X$ of all dimensions  $n$ 
  support 
 no metrics with $Sc>0$.} \vspace {1mm}

{\it Proof.}  This is obvious once said.
Indeed, a {\sl simple special case of the Atiyah-Singer   index theorem} says that,
\vspace {1mm} 

{\sf  if a
complex vector bundle $L$ of rank $k$ over a compact  orientable {\it spin}  Riemannian manifold $X$ of dimension $n=2k$,  has  {\it non-zero   Euler} ({\it Chern}) {\it number}, that is the self-intersection index of the zero section $X\to \hookrightarrow \underline L$,} then 

{\it the twisted Dirac  
$D_{\otimes L}: C^\infty(\mathbb S\otimes L)\to C^\infty(\mathbb S\otimes L)$ has {\it non-zero kernel}, for all linear connections in $L$,}
provided,
  
{\sf the number $k$ is {\it odd}, and 
the restriction of $L$  to the  complement to a point in $X$ is a {\it trivial bundle.}\footnote{These are minor technical conditions, the role of which is  to avoid undesirable consequences of
  possible cancellation in  the index formula (see section \ref{Dirac4}).
For instance if  $X$ can be embedded or immersed into  $ \mathbb R^{2k+1}$, or if it admits a metric with positive scalar curvature  then {\sf even} $k$ is allowed.
 (Observe in passing that  these $X$ are spin.)}}\vspace {1mm}

Then, by elementary algebraic topology, 

{\sf the $2k$-sphere supports a complex vector bundle of rank $k$, say $\underline L\to S^{2k}$,   which 

has  {\it non-zero } Euler (Chern) number}, 

\hspace {-5mm}and  

{\sf  bundles $L=f^\ast(\underline L)\to X$  induced from   $\underline L$ by  continuous maps 
$f:X\to S^{2k}$ have 

 their Euler numbers  $e(L)=deg(f) e(\underline L)$.}\vspace {1mm}

It follows that finite  coverings $\tilde X_\varepsilon$ of $X$  admit smooth $\varepsilon$-Lipschitz-maps $f_\varepsilon :\tilde X_\varepsilon\to S^n$ with arbitrary small $ \varepsilon$  and such that the  twisted Dirac  operators
$\mathcal D_{\otimes L_\varepsilon}$ on   $\tilde X_\varepsilon$ for  $L_\varepsilon=f_\varepsilon^\ast(\underline L)$, have non-zero kernels  for all connections in 
$L_\varepsilon$.
 
Apply this to connections  $\nabla_\varepsilon$ in 
 $ L_\varepsilon$ induced by $f_\varepsilon$  from a fixed  smooth linear (unitary if you wish)  connection $\underline \nabla$ in  $\underline L\to S^{2k}$, let  $\varepsilon\to 0$  and observe that, since the maps  $f_\varepsilon$ converge to constant ones on all  unit balls in $\tilde X_\varepsilon$,  the bundles  $(L_\varepsilon,\nabla_\varepsilon)$  converge to trivial ones with trivial flat connections on all balls.  
 Therefore the difference between the Dirac  operator $\mathcal D_{\otimes L_\varepsilon}$ and 
 $\mathcal D$ twisted with the trivial flat bundle $L_{flat}$ of rank $k$ becomes arbitrary small
for $\varepsilon \to 0$,  and the S-L-W-(B) formula applied to $\mathcal D_{L_{flat}}$ shows that 
$ \inf_X {Sc(X)  = \inf_{\tilde X} Sc(\tilde X_\varepsilon})\leq 0$. 
 
 This completes  the proof for $n=4l+2$ and the general case follows by (shamelessly) taking the product $X\times \mathbb T^{3n+2}$.
 
Well..., this is convincing but it is not quite a proof. We still have to define $\mathcal D_{\otimes L}$
and to make sense of the  "difference" between the operators $\mathcal D_{\otimes L_\varepsilon}$ and $\mathcal D_{\otimes L_{flat}}$ that are defined in {\it different spaces}.  We do all this below  closer to the end  of this section.\vspace {1mm}

{\it Why Spin?}   The essential new  information delivered by $\mathcal D_{\otimes L}$ does not  visibly  depend on the spin structure (unlike to how it is with the Dirac  operator $\mathcal D$ itself).\footnote {Sometimes, e.g. for  lower bounds  on the (area) norms of  differentials of  maps $ X\times X_{\sf kum}\to S^n$, $n=dim(X)$, for metrics $g$ on  $X\times X_{\sf kum}$ with large 
scalar curvatures, the spin is irreplaceable.}
 
However, one {\sf \color {magenta}doesn't know} how to get rid of the spin condition, in the cases where it appears irrelevant.  For instance, there is no single known  area-wise bound on the  size of  a non-spin manifold with a large scalar 
curvature.\footnote {In truth,  this applies only   to    {\it non-spin$^\mathbb C$} manifolds, where
spin$^\mathbb C$ means that the second Stiefel-Whitney class is equal to  the  mod 2 reduction of the Chern class of a complex line bundle $L\to X$.

 Such bounds are available  for spin$^{\mathbb C}$ manifolds.
For instance (a special case of) {\it Min-Oo extremality/rigidity theorem} says that

{\sf if   the scalar curvature  a Riemannian metric $g$ on  on $\mathbb CP^m$ is (non-strictly) greater  than that of the Fubini-Study metric, $Sc(g) \geq Sc(g_{FuSt})$, and 
$area_g(S) \geq area _{g_{FuSt}}(S)$
for all  smooth surfaces $S\subset \mathbb CP^m,$  than $g=g_{FuSt}$.} 

(The complex projective  spaces   $\mathbb CP^m$ are non-spin for even $m$, yet they are all 
spin$^{\mathbb C}$).}

All in all, although  "twisted Dirac"  proofs are short and simple, their  nature  remains  obscure.

 Partly, this is  why  we explain below with such a care  standard "trivial" properties of the "twist"
$\mathcal D\leadsto \mathcal D_{\otimes L}$, hoping this may help  us  to visualize 
{\it something}  behind this  "trivial" that  makes the  Dirac's $\cal D$ work in geometry,  "something", which  is only tangentially related to the Dirac  operator  itself and, 
if untangled from $\cal D$ with  its bondage to spin, 
would open  up new possibilities. 
\vspace {2mm}
%%%%%%%%%%%%%%%%%%
\subsubsection {\color{blue}Recollection on Linear Connections and Twisted  Differential Operators} \label{connections3} 
%%%%%%%%%%%%%%%%%%%%%%%

 A connection   in a smooth fibration  $L\to X$ is a {\it retractive homomorphism}  from the tangent bundle $T(L)$ to the subbundle  $T_{vert}=T_{ver}(L)\subset T(L)$ of the vectors tangent to the fibers of $L$.\footnote{Here, "retractive" means being  the identity on $T_{vert}.$}

 Denote this by
$$\hat \nabla: T(L)\to T_{vert}\subset T(L),$$
and observe that $\hat\nabla$ is uniquely defined by its kernel,  that is what is called a horizontal 
 subbundle,  $T_{hor}=T_{hor}(L)\subset T(L)$ that is complementary to $T_{vert}$ such that 
  $T(L)=T_{vert}\oplus T_{hor}$.

If $L$ is a trivia?l (split) fibration $L=X\times L_0$,  then it comes with the trivial or split flat connection, where $T_{hor}$  is the bundle of vectors tangent to the graphs of constant maps  $X  \to L_0$, $l\in L$.

A connection is called {\it flat}  at $x_0\in X$ if, over a neighbourhood $U\subset X$ of $x$,  it is {\it isomorphic} to the trivial flat connection  on $X\times L_{x_0}$, for the fiber $L_{x_0}$ of $L$ over $x_0$.

If  the fibration $L$ carries a  fiber-wise geometric structure $\mathscr S$, say, {\sf linear, affine, unitary, etc},
then "flat" signifies that the implied isomorphism, that is a fiber preserving diffeomorphism
$L_{|U}\to U\times L_x$, preserves $\mathscr S$, i.e. it is  fiber-wise {\sf linear, affine, unitary, etc.}

A connection $\hat\nabla$ in $L$ is called  $\mathscr S$: {\sf linear, affine, unitary, etc} if, for each $x\in X$,   there exist   
a  flat $\mathscr S$-connection $\hat\nabla_{x, flat}$  {\it adapted to} $\hat\nabla$ at $x$, i.e.such that 
 the 
 restriction of $\hat\nabla_{x, flat}$ to the fiber $L_x\subset L$, denoted  $(\hat\nabla_{x,flat})_{|L_x}$ is equal to $\hat\nabla_{|L_x}$.

{\it Twisting Differential  s}. A first order   differential  between (sections of) vector bundles 
(linear fibrations) $K_1$ and $K_2$ over a manifold $X$, is a linear map
$$ D: C^\infty(K_1) \to   C^\infty(K_2),$$
such that the value  $D f (x)\in K_2$ depends only on the differential 
 $df(x): T_x(X)\to T_{f_x}(K_1)$ for all $x\in X$.

For instance, a linear connection  in $L$ defines a differential , denoted just $\nabla$, from $L$ to the bundle 
$Hom(T(X),L)=T^\ast(X)\otimes L$, that is the  composition of  the differential $df: T(X)\to T(L)$
with  $\hat\nabla : T(L)\to T_{vert}$ combined  with the  canonical identifications of   all (vertical) tangent 
spaces  of the fiber $L_x$ with $L_x$ itself.

Such a $\nabla$ uniquely determines (linear) $\hat\nabla$,  it is also called "connection".
where 
the values $\nabla f(\tau)$ at  tangent vectors $\tau$ are written as (covariant) derivatives $\nabla_\tau f$.

{\it Basic Example.} If $\nabla$  is the flat split connection in $X\times L_0$, then this  is applies to sections $X\to X\times L_0$,  that are the graph of  maps $f:X\to L_0$, as 
the ordinary differential $df:T(X)\to L_0$.

If a section $f:X\to L$ vanishes at a point $x\in X$, then, clearly,   $\nabla f(x)=\nabla_{flat}f(x)$ for all {\it nabla}.

It follows that  the difference between two connections in $L$,
$\nabla_1-\nabla_2$, is a {it zero order}  defined by a homomorphism 
$\Delta=    \Delta_{1,2}: L\to Hom(T(X),L)$, that can be thought of as a $Hom(L,L)$-valued 1-form on $X$.

 Thus any $\nabla$ in a flat, e.g. split, bundle is $df+\Delta$. \vspace {1mm}

If $\nabla$ is a {\it flat split connection}, in $L=X\times L_0$, then the twisted  $D_{\otimes L}: C^\infty(K_1\otimes L)\to C^\infty(K_1\otimes L)$
 is defined via the
identity $C^\infty(K\otimes L_{split})= C^\infty(K)\otimes L_0$, 
as it was  explained above for the Dirac operator.

If  $\nabla$ is {\it flat}, then    $D_{\otimes \nabla}=D_{\otimes (L, \nabla)}$ is defined on all neighbourhoods where this connection splits and   local $\Rightarrow$ global by locality of differential s. 

Finally, for {\it general}   $(L, \nabla)$,  the twisted 
   $D_{\otimes \nabla}(\psi) $  for sections $\psi:X\to K_1\otimes L$ is determined by 
   its values at all points $x\in X$ that are defined as follows
   $$  D_{\otimes \nabla}(\psi) (x)= D_{\otimes \nabla_{x, flat}}(\psi) (x)$$
   for  flat connections  $\nabla_{x, flat}$ adapted to $\nabla$ at $x$.

  Since the difference  $\nabla- \nabla_{flat}$ is a {\it zero order } for all connections $\nabla$ in  flat split bundles $L=(X\times L_0 \nabla_{flat})$,  the same is true for $D$ twisted with $\nabla$:
   the difference 
  $$\Delta_\otimes =D_{\otimes \nabla}-D_{\otimes \nabla_{flat}}$$ 
  is a zero  order   --  "{\sf vector potential}"  in  the physicists' parlance.

 A  similar representation  $D_{\otimes \nabla}=D_{\otimes  \nabla_{flat}}+\Delta_\otimes$ for topologically non-trivial bundles $L$ is achieved 
 as follows.

Let $L^\perp\to X$ be  a bundle complementary to $L$ such that  the Whitney sum of the two bundles  topologically splits, 
$$L\oplus L^\perp=L^\oplus \simeq  X\times( L_0\oplus L^\perp_0)$$
 and let $\nabla^\perp$ be an arbitrary connection in $L\oplus L^\perp$
 and
 
 Define  the connection  $ \nabla^\oplus=\nabla\plus\nabla^\perp$ in $L^\oplus$ by the rule
$$ \nabla_\tau ^\oplus (l+l^\perp)=\nabla_\tau l+ \nabla_\tau^\perp l^\perp$$
and observe that the $\nabla^\oplus$-twisted operator $D_{\otimes \nabla^\oplus}$,  that maps the space  
 $$C^\infty (K_1\otimes L^\oplus)=C^\infty (K_1\otimes L)\oplus C^\infty (K_1\otimes L^\perp)$$
 to  
 $$C^\infty (K_2\otimes L^\oplus)=C^\infty (K_2\otimes L)\oplus C^\infty (K_2\otimes L^\perp) $$
 respects this splitting:
 $$D_{\otimes \nabla^\oplus}=D_{\otimes \nabla}\oplus D_{\otimes \nabla^\perp}.$$
Thus, if not the  $D_{\otimes \nabla}$  itself, then its  $\oplus$-sum with another  is\vspace {2mm}

\hspace{30mm}   $D_{\otimes \nabla_{flat}}$ +  {\sf zero order }.
\vspace {2mm}
%%%%%%%%%%%%%%%%%%%
\subsubsection{\color {blue}  {\color {blue}$\mathbf {[Sc\ngtr 0]}$} for Profinitely  Hyperspherical   Manifolds,  Area Decreasing Maps and Upper  Spectral Bounds for Dirac Operators}\label{area+Vafa3}
%%%%%%%%%%%%%%%%%%%%%%%%%%

{\it Conclusion of the proof of provisional proposition {\color {blue}$\mathbf {[Sc\ngtr 0]}$}  from {\ref{twisted3}.}}  Return to the   bundles $L_\varepsilon=f^\ast(\underline L)\to X$   induced  by smooth  $\varepsilon$-Lipschitz maps $f:X\to S^n$, $n=dim(X)=4l+2$,
with non-zero degrees and $\varepsilon \to 0$
 from a  complex vector bundle   $\underline L\to S^n$, 
with 
the Euler number $e(\underline L) \neq 0$. 

Let  $\underline L^\perp\to S^{2k}$ be a bundle {\it complementary} to  $\underline L\to S^{2k}$,
i.e. the sum $\underline   L\oplus \underline L^\perp  $ is a trivial bundle, 
endow  $\underline   L$ and  $ \underline L^\perp  $ with a  connections  $\underline \nabla$  and  $\underline \nabla^\perp$ and let 
 $\nabla_\varepsilon^\oplus$ be the  connection on the (topologically trivial!) bundle 
 $$L^\oplus_\varepsilon=f^\ast(\underline L\oplus \underline L^\perp)$$
 induced from  $\underline \nabla^\oplus= \underline \nabla \oplus \underline \nabla^\perp,$
where the latter is defined  by the  component-wise differentiation rule:
 $$\underline \nabla^\oplus(\phi, \psi) =(\underline \nabla \phi, \underline \nabla^\perp \psi)\mbox { for   sections }  (\phi,\psi) = \phi+\psi  :S^n\to \underline   L\oplus \underline L^\perp. $$

Then (see the   proof of {\color {blue}$\mathbf {[Sc\ngtr 0]}$}) the twisted Dirac operator  decomposes into the sum of (essentially) untwisted $\cal D$ and a zero order  (vector potential)
$$\mathcal D_{\otimes \nabla_\varepsilon^\oplus}=\mathcal D_{\otimes \nabla_{flat}}+\Delta_{\varepsilon}$$
where $\nabla_{flat}$ is the flat split connection in the bundle $L^\otimes_\varepsilon  $ with the  splitting induced by  $f_\varepsilon$ from a splitting of  $\underline L\oplus \underline L^\perp$,
obviously (but most  significantly), $\Delta_{\varepsilon}\to 0$ for $\varepsilon\to 0$.

Now,  the (untwisted) S-L-W-(B) formula, 
   applied to $ \mathcal D_{\otimes \nabla_{flat}}$ 
    says that 
$$\mathcal D_{\otimes \nabla_\varepsilon^\oplus}^2 = \nabla_{flat, \mathbb S} \nabla_{ flat, \mathbb S}^\ast + \frac{1}{4} Sc+\Delta^\smallsquare_\varepsilon,$$
where $\nabla_{flat, \mathbb S}$  denotes the flat connection $\nabla_{flat, \mathbb S}$   in the twisted spin bundle associated with $\nabla_{flat}$.

The correction term $\Delta^\smallsquare_\varepsilon$ in this formula  is a first order differential  (it depends on how  you trivialise the   bundle  $\underline   L\oplus \underline  L^\perp$)
 which  tends  to $0$  for $\varepsilon \to 0,$
$$\Delta^\smallsquare_\varepsilon\to 0 \mbox {  for }\varepsilon \to 0.$$

A priori, the $\varepsilon$-bound on the differential of $f_\varepsilon$ doesn't make the  coefficients of the   $\Delta^\smallsquare_\varepsilon$ small, but an obvious smoothing allows an approximation of $f_\varepsilon$ by maps  their derivatives of  which of  {\it all 
 orders  converging to $0$.}

Because of this, we may assume  $\Delta^\smallsquare_\varepsilon\to 0$ in the {\it strongest conceivable
 sense}, while al  is needed  is  that $\Delta^\smallsquare_\varepsilon\to 0$ becomes negligibly small compare to $\nabla_{flat, \mathbb S} \nabla_{ flat, \mathbb S}^\ast + \frac{1}{4} Sc$, which implies 
{\it strict positivity} of the  
 $\mathcal D_{\otimes \nabla_\varepsilon^\oplus}^2 = \nabla_{flat, \mathbb S} \nabla_{ flat, \mathbb S}^\ast + \frac{1}{4} Sc+\Delta^\smallsquare_\varepsilon,$ for $\varepsilon$ much smaller than the lower bound $\sigma=\inf_{x\in X} Sc(X,x)>0$.

Thus, the condition $Sc(X)>0$  leads to a contradiction with the index formula, which in this case,  as we already know  from  the   proof of {\color {blue}$\mathbf {[Sc\ngtr 0]}$}    yields non-zero 
 harmonic $\nabla_\varepsilon$-twisted, hence  $\nabla_\varepsilon^\oplus$ twisted, spinors, because  
the subbundle $\underline L\subset  \underline L\subset$ is invariant under the parallel transport by the connection  $\underline \nabla^\oplus= \underline \nabla \oplus \underline \nabla^\perp$, by the very definition of the sum of connections  and this  property is inherited  by   the induced connection 
$ \nabla_\varepsilon^\oplus$.
 
This concludes the proof of   {\color {blue}$\mathbf {[Sc\ngtr 0]}$} for $n=4l+2$ and, as we have already explained,    the general case follows by stabilization $X\leadsto X\times \mathbb T^{3n+2}$.

\vspace {1mm}

{\it Area Contraction instead of Length Contraction.} Say that $X$ is 
$\wedge^2$-{\it profinitely hyperspherical} if, instead of $\varepsilon$-Lipschitz property of maps 
$f_\varepsilon\tilde X_\varepsilon\to S^n$ of no-zero degree, we require that the second 
exterior power of the differential of $f_\varepsilon$ is bounded by $\varepsilon^2$,
$$||\wedge^2df_\varepsilon(x)||\leq \varepsilon^2.$$ 
 Geometrically, this means that 
 $ f_\varepsilon $ decreases the areas  of the smooth surfaces in $X$ by factor $\varepsilon^2$, (This, obviously, is satisfied by    $\varepsilon$-Lipschitz maps.)

It is clear, heuristically, that  the   Dirac   operator {\it twisted with} $\nabla_\varepsilon$ in this case,
 similarly how it is for   $\varepsilon$-Lipschitz maps, is  {\it close to the untwisted} $\mathcal D$;     this \vspace {1mm}
 
 {\sf rules out positive scalar curvature for  $\wedge^2$-profinitely hyperspherical manifolds.}
\vspace {1mm}

  However, the above proof with the  complementary  bundle $L^\perp$
doesn't apply here; to justify  heuristics, one has to  pursue  algebraic similarity between 
$\nabla$ and the ordinary differential $d$ a step further.   

This can be done by pure thought,  on the basis of general principles only, (no tricks like $L^\perp$)  but 
writing down this "thought"  turned out     more space and time consuming than what is needed  for  
(a few lines of)  the twisted version of the  S-L-W-(B) formula, as we shall see   in section \ref{K-area3}.

So,  we conclude here  with two remarks.

(i) It is { \color  {magenta} unknown} if  "length contractive" is {\it topologically} more   restrictive than "area contractive".

For instance one has no idea if there exist $\wedge^2$-profinitely hyperspherical manifolds which are {\it not } profinitely hyperspherical. 
  
(ii) Representation of  $\nabla$-twisted   differential operators by vector-potentials $\Delta$  in larger bundles has further uses, such  as  {\it Vafa-Witten's lower bounds on the  spectra of Dirac operators}.  For instance,\vspace {1mm}

{\sf \color {teal!50!black}if a  compact  Riemannian spin  $n$-manifold $X$   admits a distance decreasing map to $S^n$ of degree $d$, then the number $N$ of eigenvalues $\lambda$ of the Dirac  on $X$ in the interval
$   -C_n\leq \lambda \leq  C_n$ satisfies $N\geq d$, where 
 $C_n>0$ is a universal constant.\footnote{See  \S 6 
in [G(positive) 1996]
for related spectral geometric  inequalities.}}

%%%%%%%%%%%%%%%%%%%%%

\subsubsection {\color {blue} Clifford Algebras, Spinors, Atiyah-Singer Dirac Operator and  Lichnerowicz Identity
 }\label{Clifford3}
 %%%%%%%%%%%%%%%%%%%%%%%%%%

 The Dirac  on $\mathbb R^n$ is a particular first order differential , 
 which acts on the space of smooth  $\mathbb C^N$-valued functions, 
  $$\mathcal D: C^\infty(\mathbb R^n, \mathbb C^N)\to  C^\infty(\mathbb R^n, \mathbb C^N),$$ where $N=2^{\frac {1}{2}n}$  for even  $n$ and   $N=2^{\frac {1}{2}(n-1)}$  for odd $n$ and  where    
  essential  properties  of this $\cal D$ are as follows.

 \textbf{I.} {\it Ellipticity.} The   $\cal D$ is an  {\it elliptic}, which  means that the initial value (Cauchy) problem 
for the equation $Df=0$
 is {\it formally  uniquely solvable} for {\it all} initial data on {\it all}   smooth hypersurfaces in $\mathbb R^n$, where "formally" can be replaced by "locally"  in the real analytic case.
 
\vspace{1mm}
 
 {\it Basic Example.} {\sf The Cauchy-Riemann} (system of two) equation(s) $D_{CR}f=0$ for maps $f:\mathbb R^2\to \mathbb C^1$,  defines conformal orientation preserving maps $  \mathbb R^2\to \mathbb C$. These  are called {\it holomorphic} if $\mathbb R^2$ is "identified" with $\mathbb C$, where the ambiguity inherent in this identification is  responsible for  spin. 
 
 {\it $D_{CR}$ is elliptic}:    {\sf real analytic functions 
locally   uniquely extend from real analytic curves in   $\mathbb C^1$ to holomorphic functions.}

 Let us describe ellipticity in linear algebraic terms applicable to all  (systems of) partial differential equations of first order for maps between smooth manifolds,  $f:X\to Y$.
 Such a system, call it $S$, is characterised by  subsets  in the spaces of linear maps 
 between the tangent spaces of $X$ and $Y$  at all  $x\in X$ and $y\in  Y$, denoted
 $\Sigma_{x,y}\subset Hom(T_x\to T_y)$,  where $T_x=\mathbb R^n$, $n=dim(X)$, where  $T_y=\mathbb R^m$, $m=dim(Y)$ and where  $f$ {\it satisfies} $S$ if $df(x)\in \Sigma_{x,f(x)}$  for all $x\in X$.
   
 Let $R_L : Hom(\mathbb R^n, \mathbb R^m) \to Hom(L, \mathbb R^m)$ denote the restriction 
  of linear maps to $\mathbb R^m$ from  $\mathbb R^n$ to  linear  subspaces $L\subset \mathbb R^n$, that is $R_L:h\mapsto h_{|L}:L\to \mathbb R^m$.\vspace {1mm}

 Call  a smooth submanifold  $\Sigma\subset Hom(\mathbb R^n, \mathbb R^m)$   {\it elliptic} if 
 the map $R_L$ {\it diffeomorphically} sends  $\Sigma$ {\it onto $Hom(L, \mathbb R^m)$}
 for all {\it hyperplanes} $L\subset \mathbb R^n$.
 
 \vspace {1mm}
 
 Now,  a PDE system $S$ is called {\it elliptic} if the  subsets 
 $$ \Sigma_{x,y}\subset Hom(T_x\to T_y)=H_{n,m}=Hom(\mathbb R^n, \mathbb R^m)$$   
 are elliptic for $x\in X$ and $y\in Y$.

 Put it  another way, let $K_p\in H_{n,m}$,  $p\in \mathbb RP^{n-1}$, be the family of  $m$-dimension linear subspaces that
are the kernels of the linear maps $R_{L_p}: H_{n,m}\to H_{n-1,m,p}=Hom(L_p, \mathbb R^m)$ parametrized by   the projective space $\mathbb RP^{n-1}$ of hyperplanes $L=L_p\subset \mathbb R^n$. Then ellipticity says that 

{\sf $\Sigma$ {\it transversally} intersect $K_p$ at  a single point for 
all $p\in \mathbb RP^n$.}\vspace {1mm}

Finally, back to the  linear  case,  observe that      systems  $Df=0$ for maps 
$$f: \mathbb R^n\to\mathbb R^m,\mbox { $x\in \mathbb R^n$}$$ 
are 
depicted by    {\it  linear} subspaces  
$$\Sigma =\Sigma_x\subset Hom(T_x( \mathbb R^n), \mathbb R^m= T_{0}( \mathbb R^m)), \mbox  {$x\in\mathbb R^n$}$$
  and   ellipticity  says in these terms  that 
  
  {\sf firstly, $dim(\Sigma)=n$}
 
\hspace {-6mm}  {\sf secondly}  \vspace {1mm}

\hspace {-2mm}{\color {teal} \textbf \EllipseShadow}  {\it the  linear maps   $h: T_x(\mathbb R^n) \to \mathbb R^m$ have  $rank(h)=n$ for 
all non-zero $h\in \Sigma$}.
  
  {\sf and finally} 

{\it  differential operators between sections of vector bundles over a smooth manifold  $X$ are elliptic if these properties are verified locally over small neighbourhoods of all points in $X$.} 
 
 \vspace {1mm}

{\it Exercises}. (a) {\sf Twisting with $\nabla$}. Show that   

\hspace {30mm}  {\sf $D$ is elliptic $\Rightarrow $  $D_{\otimes\nabla}$ is elliptic:}
  
 twisting with connections doesn't hurt ellipticity. 
 \vspace {1mm}

(b) {\sf Symmetric but non-Elliptic.} Figure out what makes  the exterior differential 
$$d: C^\infty (\bigwedge^k(T(X)) \to C^\infty (\bigwedge^{k+1}(T(X))$$
 on  $(2k+1)$-dimensional manifolds {\it non-elliptic.}

\vspace {2mm}

\textbf {II.} {\it Symmetry and Spinors.} The  Dirac  operator $\cal
D$  on $\mathbb C^\infty(\mathbb R^n, \mathbb C^N)$,  is  {\it $Spin(n)$-invariant},   where 
 \vspace {1mm}
 
$\bullet_1$ {\sf $Spin(n)$ denotes  the double cover of the special orthogonal group  $SO(n)$,
 
$\bullet_2$ the group $Spin(n)$ acts on $\mathbb R^n$ via   the (2-sheeted covering map) homomorphism 
$Spin(n) \to SO(n)$,

 $\bullet_3$ the  action of $Spin(n)$ on 
$\mathbb C^k$, called {\it spin representation},   is {\it faithful}:  it {\it doesn't factor} through an action of $SO(n)$,\footnote 
{The spin representation, as we shall explain below, is {\it irreducible} for odd $n$ and it splits into  
{\it two irreducible half-spin  representations} for even $n$. There are no {\it faithful} representations of  $Spin(n)$ in lower dimensions (except for $n=1, 2$), where,
 apparently, this faithfulness is  necessitated by  ellipticity of $\mathcal D$.}

$\bullet_4$  "invariant" here means {\it equivariant} under the (diagonal) action of $Spin (n)$ on 
the space of maps $\psi:\mathbb R^n\to \mathbb C^N$, that is 
$$g(\psi)(x)=g(\psi(g(x)), \mbox { } g\in Spin(n),  \footnote {To visualize this, think of the graphs $\Gamma_\psi\subset \mathbb R^n\times \mathbb C^N$ moved by the diagonal actions of $g\in Spin(n)$ on this product.}  $$
and "{\it equivariant}" says that 
$$ \mathcal D( g(\psi))= g( \mathcal D( \psi)).\footnote {This  Dirac  operator has "constant coefficients", which means  is {\it invariant under parallel translations $t_y$} of $\mathbb R^n$ that act on our maps: 
 $D (t_y(\psi))=  t_y(D( \psi))$   for $(t_y(\psi))(x)= \psi(x+y)$, $x,y\in \mathbb R^n$.}$$

 {\it Cauchy-Riemann Example.}
The  group   $Spin(2)$  diagonally  acts on 
on maps $f: \mathbb R^2\to \mathbb C^1$, where all actions   (representations) of  $Spin(2)=\mathbb T =U(1)\subset \mathbb C^\times$  
on $\mathbb C^1$ are possible: these are $t(z)=t^mz$, $m=...-1,0,1,2,...$. (There are no such 
possibilities for   for $n>1$.)

The corresponding operators $\bar\partial=\bar\partial_m$ are all  {\it locally non-canonically} isomorphic (this  makes them often confused in the literature), but this $m$ 
 (spin quantum number),  becomes the major feature of the  
  $\bar\partial_i$ globally, where it controls its  very existence  and 
  its  index.}\vspace {2mm}

 \textbf {III.} {\it Spin Representations and   Clifford Algebras $Cl_n=Cl(V)$.}\footnote {The basic reading   on this subject matter is the book   [Lawson\&Michelsohn(spin geometry) 1989]  and a (very) brief outline of the main points is contained  in [Min-Oo(K-Area) 2002],
[Min-Oo(scalar) 2020].}The lowest dimension complex vector 
  space, 
  where $Spin(n)$,  can  linearly  {\it faithfully} act is 
  $ \mathbb C^{2^k}$ for $k=\frac {1}{2}n$ for $n$ even and   $k=\frac {1}{2}(n-1)$ for odd $n$, where  such an action (representation) is  obtained   by realizing   
   $Spin(n)$ as a  subgroups in the multiplicative semigroups of the  {\it Clifford algebra}, denoted  
  $Cl_n =Cl(\mathbb R^n)=Cl(\mathbb R^n,-\sum_1^nx_i^2)$.
   \vspace {1mm}
   
   Recall that $Cl_n$ is an
     {\it unital\footnote {This means possessing  a unit in it} associative algebra  $A$} over the field of real numbers with a  distinguished {\it Clifford basis} that is  linear subspace  $V=V_{Cl}\subset Cl_n$  {\it endowed with a Euclidean structure}, that is represented by a {\it negative  definite}    quadratic form.\footnote {It is negative  to agree with  the  Laplacian $\sum_i\partial_i^2$, which  is a negative operator.} 
     
      We denote  the  {\it Clifford product} by  $a_1\cdot a_2$ an let  "1" stand for the unit in $A$.
    
\vspace{1mm}

{\sf (There is nothing especially exciting   about $Cl_n$ understood as "just an algebra", especially if you tensor it  with $\mathbb C$,  which  we  do at the end of the day anyway.  For instance,   we shall see it presently,  $Cl\otimes \mathbb  C$ is isomorphic  to a matrix algebra  for even $n$ and to the sum of two matrix algebras for odd $n$. 

What  gives to a particular favour to $Cl_n$ is the distinguished linear subspace $V\subset Cl_n$, which,  on the  one hand, {\it generates all} of $Cl_n$, on the other hand,
the  matrices corresponding to all $v\neq 0$ in $V$,  have {maximal possible} ranks, since 
all non-zero $v\in V$  are {\it invertible} in the multiplicative semigroup $CL_n^\times$.
This "maximal rank property" is  exactly what   makes the Dirac  operator {\it elliptic} and, because of this,  so powerful in the  Riemannian geometry.)}

 \vspace{1mm}

 The fundamental feature of  the  pair $(A,V)$ is that $A=Cl(V)$ is   {\it  functorially} determined by $V$: \vspace {1mm}
 
 \hspace {-4mm}  {\sf isometric embeddings $V_1\to V_2$ canonically extend to monomorphisms $A_1\to A_2$.}
     \vspace {1mm}

   \hspace {-6mm}where this Clifford functor is uniquely characterised by the following two properties.

 A. {\it  $V=V_{Cl}$ is a Basis in $A$.} The subspace $V$ generates  $A$ as an 
 $\mathbb R$-algebra.
   
  B. {\it Specification.} The algebra  $Cl_1=Cl(\mathbb R^1)$ is  isomorphic to $(\mathbb C, i\mathbb 
  R)$, for  $i=\pm\sqrt {-1}$.

 (It is {\color {red!40!black} \sf  impossible} to {\it mathematically}, distinguish $i$ and  $-i$;   this unresolvable  $\pm$i-ambiguity is grossly amplified, at least   psychologically,,  when it comes to spinors.
 \footnote{To be blameless, write $\pm i$  (even better, $\{ \pm i,\mp i\}$)
and  never dare utter  "left   ring ideal"  and    "right group action", even in absence of  
  left-handed (left-minded?) persons.  (Defending such an action by    {\it biological molecular homochirality}  and   {\it  parity violation by weak interactions}  
   is not recommended for being   politically incorrect.)  
  
  Jokes apart, arbitrary   terminological   conventions  presented as  mathematical definitions
  sow confusion  and 
   undermine "rigor" in mathematics.
    
    Who  are the  lucky ones who are able to tell if $f\circ g$ means 
   $f(g(x))$ rather than   $g(f(x))$  or  vice versa?

 Can   {\sl encoding   formulas by Peano's integers},  e.g. in   the proof of    {\it G\"odel's incompleteness  theorem},  be accepted as "logically rigorous", unless you   face  the issue  of   "directionality" inherent    in the decimal representation of integers?})

 In simple words, the Clifford squares of all unit vectors $v\in V$ are equal to $-1$, or, equivalently,
 $$v\cdot v=-||v||^2=\langle v, v\rangle  \mbox { for all $v\in V$}.$$
 
  \vspace {1mm}

 A\&B.  {\it Anti-commutativity.} The  Clifford product is {\it anti-commutative on orthogonal vectors.}
    $$ v_1\cdot v_2=-v_2,v_1,\mbox {  for } \langle v_1,v_2\rangle=0. $$
    
    Indeed, since $||v_1-v_2||^2=||v_1+v_2||^2$ for orthogonal vectors, 
    bilinearity of the the Clifford product implies that
    $$0=(v_1-v_2)^2-(v_1+v_2)^2=-v_1\cdot v_2-v_2\cdot v_1+v_1\cdot v_2+v_2\cdot v_1=2(v_1\cdot v_2+v_2\cdot v_1).$$

  {\it Exercise.} Show that 
  $$\mbox {$v_1\cdot v_2+v_2\cdot v_1=-2\langle v_1,v_2\rangle$ for all
    $v_1,v_v\in V$.}$$
  
  \vspace{1mm}
   
  \textbf {IV.} {\it Groups  $Pin(n)$  and $Spin(n)$  and $G_n$.}  The group $Pin(n)$ is   
   defined in  $Cl_n$-terms   as the subgroup of the multiplicative semigroup of $Cl_n^\times \subset Cl_n$ {\sf multiplicatively generated by 
   the unit vectors $v\in V\subset Cl_n$.}  \vspace{1mm}
   
   {\sf  The subgroup  $Spin(n)\subset Pin(n)$  consists of the products of {\it even numbers} of 
 
   unit  vectors from $V$.}  
 \footnote {This parallels   the definition of 
$SO(n) \subset O(n) $ as the subgroup consisting of products of {\it even numbers}  of reflections of $\mathbb R^n$. In fact, $Spin(n)$ equals the connected component of the identity in $Pin(n)$ and 
  $Pin(n)/Spin(n)=O(n)/SO(n)=\mathbb Z_2=\{-1,1\}$.}
\vspace {1mm}

 {\it Existence} \&{\it Uniqueness.} 
Let us explain why the algebra  $Cl(\mathbb R^n)$,  if exists at all, is  large enough to (multiplicatively) contain the   group $Spin(n)$ that double covers the special  orthogonal  group  $SO(n)$.
\footnote {To appreciate non-triviality of the problem,
try to construct  geometrically  more than two, say three,   anti-commuting linear isometric involutions represented by  
  reflections around linear subspaces  in some   Euclidean space.} 

 Observe that the {\it Clifford relations} \footnote {This must be  written in Clifford's   unpublished note {\sl On The Classification of Geometric Algebras} see [Diek-Kantowski (Clifford History)1995]  for further references.}
$$\mbox {$e_i\cdot e_j=-e_j\cdot e_i$ and $e_i^2=-1$} \leqno {\mathbf {[Cl]}}$$ 
for  an orthonormal frame $\{e_i\}\subset V$, $i=1,...,n,$\vspace {1mm}

on the one hand, imply  $v_1\cdot v_2+v_2\cdot v_1=-2\langle v_1,v_2\rangle$ for all
    $v_1,v_v\in V$, hence, 
    
   \hspace {28mm} {\it fully characterize Clifford's algebras, }\vspace{0.5mm}
    
    on the other hand, 
define \vspace{0.5mm}

\hspace {10mm}{\sf a finite group $G_n$ of order $2^{n+1}$ that is  a  central extension of $\mathbb Z_2^n$},
\vspace {1mm}

\hspace {-6mm}  with an additional generator   (central element)  $c$ of order 2 and  the following relations, 
$$\mbox {$ ce_i=e_ic$,  $c^2=1$,    $e_ie_j=ce_je_i$ and $e_i^2=c.$} \leqno {\mathbf {[Cl_c]}}$$
where  the central element $c$ in $G_n$ corresponds to $-1\in Cl_n$.

Non-triviality  of this $G_n$ is apparent, since  letting $c=1$  defines  a {\it surjective} homomorphism $G_n\to\mathbb Z_2^n$ with kernel   $\mathbb Z_2$.  

(What is not immediately apparent, is a pretty combinatorics of shuffling  indices in  $e_{i_1}e_{i_2}... e_{i_m}\in G_n$, $i_1<i_2<...<i_m$, under multiplication   by $e_k$, which is 
 rightly appreciated by people working on  
 quantum computers.)

One  look at  $G_n$    is sufficient to make it  obvious that there is    a {\it homomorphism from $G_n$ to the multiplicative (semi) group $Cl_n^\times $}  of the Clifford algebra  (with the image in  $Pin(n)\subset Cl^\times $),  such that
\vspace {1mm}

{\sf  the algebra  homomorphism from the 
{\it real group algebra} $\mathbb R(G_n)$\footnote{$\mathbb R(G)$ is the space of formal linear combinations  $\sum_{g\in G}  c_gg$ with the obvious product rule,  where  the identity element $id\in G$ serves as  the unit of this algebra.

 Alternatively, $\mathbb R(G)$ is defined  as the algebra of  linear operators $\mathbb R(G)$ on functions $\psi(g)$  that is generated by translations on the space of functions on $G$, for 
  $\psi(g)\mapsto \psi(g'g)$,  $ g'\in G$. 

The same space  $\mathbb R(G)=G^\mathbb R$ of functions on $G$ with the action of $G$ by $ \psi(g)\mapsto \psi(g'g)$  is called (not very inventively) the {\it regular $\mathbb R$-representation}
of $G$,   where just   "regular representation" stands for 
   regular $\mathbb C$-representation.}
to $Cl_n$ associated with this group homomorphism $G_n\to Cl_n^\times$  is {\it surjective}   and the kernel of this homomorphism is defined  
  by the relation $c=-1$, that is $$Cl_n= \mathbb R(G_n)/(c+1),\footnote {Recall that $c\in G_n\subset \mathbb R (G_n)$ is the central involution in $G_n$ and "1" is the unit in the algebra   $\mathbb R (G_n)$ that is  represented by the unit function,
   where 
 $(c+1)\subset \mathbb R(G_n)$  denotes the   ideal generated by $c+1\in  \mathbb R(G_n)$.
 (The quotient algebra $ \mathbb R(G_n)/(c+1)$   has the same underlying linear space as  the group algebra 
 $\mathbb R(G_n/(c))$,
 for the normal subgroup  $(c)\subset G_n$ generated by $c$, but  multiplicatively
 $ \mathbb R(G_n)/(c+1)$  is  much different from the  (commutative) group  algebra of $G_n/(c)=\mathbb Z_2^n$.) }$$}\vspace {1mm}
 
 Amazingly,  nowhere, except for a few papers on quantum computers,    $G_n$ is called  "finite Clifford group",\footnote{ The terms "Clifford group", sometimes   "naive  Clifford group",   are reserved for the  subgroup $G$ of the multiplicative semigroup of $Cl$, the action of which on $Cl$ by conjugation for $a\mapsto g\cdot a\cdot g^{-1}$ preserves $V$.}     while the  authors of the  only mathematical papers found on the web (unless I missed some) call 
$G_n$ a "{\it Salingaros vee group.}"\footnote {See [AbVaWa(Clifford  Salingaros Vee)2018] for  more general definitions and references to the the original 1981-82 papers by  Nikos   Salingaros. (I don't know what is written  in these papers, since these are not  openly accessible on the web.)

Also,  amazingly, no survey or tutorial on Clifford algebras I located   on the web makes any use or even mentions $G_n$. Possibly, there is something about it in textbooks, but none seems to be 
openly  accessible.  }\vspace {1mm} 
 
The structure this  "vee  group" $G_,$,  which  tells you everything about $Cl_n$, is  transparently seen in the combinatorics of its multiplication table, where  $g\in G$  are  written as  {\it lexicographically ordered products} of $e_i$ and (if it is there) $c$. Here are a few relevant properties of $G_n$. 

 \vspace {1mm}

{\sf All  elements in 
$G_n$  have orders 2 and/or 4. 

The commutator subgroup $[G_n,G_n]=\{g_1g_2g_1^{-1}g_2^{-1}\}$ equals to the 2-element (central)  subgroup $\{1,c\}$.

If $n$ is even,  it coincides with the center of $G_n$;
$$center(G_n)=[G_n,G_n]=\{1,c\}$$

If $n$ is odd, the center  has order 4. For instance $G_1=\mathbb Z_4$; in 
general, the extra central element for $n=2k+1$ is the product  $e_1e_2,...,e_n$.

If $n$ is even, then the number  $N_{conj}(G_n)$ of the conjugacy classes of $G_n$ is 
$2^n+1$  where $2^n$  of them comes from  $\mathbb Z_2^n$ and the extra one is that of $c$.  If $n$ is odd, there are $2^n+2$ classes, where centrality of   $e_1e_2,...,e_n$  is responsible for the additional one. }

\vspace{1mm}

  \textbf {V.} {\it Representations of the Group  $G_n$.}  The space $\Psi_n=\mathbb C(G_n)=\mathbb C^{G_n}$ of {\it complex}  functions on $G$ splits into the sum 
$\Psi_n=\Psi_n^+\oplus \Psi_n^-$, 
where $\Psi_n^+$ consists of {\it $c$-symmetric functions $\psi(g)$} that are  {\it invariant} under the action of the central $c\in G_n$, i.e. $ \psi(g)=\psi(cg)$ and where the functions  $\psi\in \Psi_n^-$ are {\it antisymmetric},
  $\psi(cg)=-\psi(g)$.
 
  The  space   $\Psi_n^+$ obviously identifies with the space $\mathbb C(\mathbb Z_2^n)$ of functions on  
  the Abelian group $\mathbb Z_2^n$, where action of  $G_n$ factors through the homomorphism 
  $G_n\to \mathbb Z_2^n$.
  
  Since the commutator  subgroup of $G_n$ is equal to $\{1,c\}$, all 1-dimensional  representations of $G_n$ are contained in $\Psi_n^+$.

Frobenius

Now, the {\it number one   theorem} in the   representation theory of finite groups reads:\footnote{This  must be attributed  to 
Frobenius (1896), since it follows by  his character theory, see \url {file:///Users/misha/Downloads/Curtis2001_Chapter_RepresentationTheoryOfFiniteGr.pdf}
 
 Unfortunately, this theorem has no name an can't be {\it instantaneously} found on Google.  }\vspace {1mm}
\vspace {1mm}

{\it the  regular representation  of $G$  uniquely decomposes into the sum of sub-representations 
 $G^\mathbb C\bigoplus_i  R^2_i$, $i=1,2,...,N=N_{irrd}(G)$, where each  $R^2_i$ is (non-canonically) isomorphic to the sum of $k_i$-copies of an 
irreducible representation $R_i$ of dimension $k_i$ and where every irreducible representations of $G$ is isomorphic to one and only one of $R_i$.

Accordingly,  the group algebra of $G$  {\sf  (the same linear space $G^\mathbb C$, but now with the  group algebra structure)}
  decomposes into the sum of  matrix algebras
$$\mathbb C(G)=\bigoplus_i End(\mathbb C^{k_i}).$$}
This is an exercise in linear algebra.
What is less obvious is that \vspace {1mm}

{\it The  number $N_{irrd}(G)$ of mutually non-isomorphic irreducible complex  representations of $G$ is equal to the number of the conjugacy classes in $G$.
 $$N_{irrd}(G)=N_{conj}(G)\mbox {   for all finite group $G$.}$$}
Consequently,

{\it the sum of the squares of the dimensions of the irreducible representations of $G$ is equal to the order of the group $G$,
$$\sum_ik_i^2=card(G).\footnote {See {\url {https://projecteuclid.org/download/pdf_1/euclid.lnms/1215467411} and  the character sections in \url {https://web.stanford.edu/~aaronlan/assets/representation-theory.pdf} and \url{https://arxiv.org/pdf/1001.0462.pdf}.}}$$}

%\vspace{1mm}

If we apply this to $G_n$ for $n=2k$, we shall see that,  besides the one dimensional representations, this group has a {\it single      irreducible} one  of dimension $2^k$, call it $S_n$, which enters the regular representation with multiplicity $2^k$.

Now, clearly, \vspace{1mm}

{\it the  $2^k$-multiple $S_n$-summand of the regular representation is exactly the space $\Psi_n^-$ of antisymmetric functions $\psi$ on $G_n$.}

\vspace{1mm}

Equally clearly, \vspace{1mm}

{\it the space   of antisymmetric functions $ \psi(g)=-\psi(cg)$ on $G_n$ (here we speak of $\mathbb R$-valued functions $\psi$)   is $G_n$-equivariantly isomorphic to $Cl_n$.}

\vspace {1mm}

\textbf {VI.} {\it Clifford Conclusion.} Since   the Clifford algebra $Cl_n$ is, as an algebra,  generated by 
$G_n\subset Cl_n$, the  representation $S_n$ of $G_n$ in $\mathbb C^{2^k}$,  that is a multiplicative homomorphism  $G_n\to  End(C^{2^k})$,  extends to an algebra  
homomorphism $ Cl_n\to  End(C^{2^k})$; hence, to

\hspace {20mm} {\it an irreducible  representation of $Pin(N)$ in $\mathbb C^{2^k}$},

\hspace {-6mm}  which extends (irreducible!) representation $S_n$ of $G_n\subset Pin(n)$.

This is  called  {\it the spin representation} and  still denoted  $S_n$.

\vspace {1mm}

\hspace {10mm}{\it Why Clifford  Algebra? Why algebras are needed here at all?} \vspace {1mm}

What we used for
 the construction of the spin representation   $S_n$  of $Pin(n)$ in $\mathbb C^{2^k}$ for even  $n=2k$ are  the  two following  simple, not to say "trivial",  but indispensable (are they?) algebra theoretic facts.

{\sf (i)  The linear  actions of $Pin(n)$ and  $G_n$ on the space $\Psi_n^-$ (and also on  $Cl_n$)    generate the same subalgebras of operators on this space.

(ii) If an  algebra $A$ of  operators on a linear  space $ M$, e.g.   $M=\mathbb C^{N^2}$, is isomorphic to the (full matrix) algebra of endomorphisms of another space,
$$A \simeq End(L),$$
then  $M$ is $A$-equivariantly isomorphic to  $End(L)$ for, say "left", action of the algebra  $End(L)$ on itself. }

(Also we were jumping back and forth between $\mathbb R$-linear  and $\mathbb C$-linear spaces and actions, but with nothing non-trivial happening on the way.)

The correspondence $\Phi:L\leadsto A=End(L)$ {\it is a functor} from the category of vector spaces over
 $\mathbb R$ to the category of unital $\mathbb R$-algebras, but $L$ can be reconstructed from  
$End(L)$ only {\it up to a homothety} $l\mapsto rl$, $r\in \mathbb R^\times$,where the projective space 
$P=L/\mathbb R^\times$ can be identified with the space of   {\it maximal left ideals} in $End(L)$\footnote{
Left ideals  $I\subset End(L)$ corresponds to   linear subspaces $L_I\subset L$,  such that 
$a\in I \Leftrightarrow  a_{|L_I}=0$.} 

{\sf (Because of this  ambiguity, one can't globally define  the Dirac  operator on a non-spin manifold  
$X$, because there is no vector bundle that would support $\cal D$.

And  although the   the projectivized  spin bundle $\mathcal {PS}\to X$ with a real projective space as the fiber is still there,  this fibration {\it admits no continuous section}  $X\to \mathcal {PS}$  --  non-zero 
 second Stiefel-Whitney class is an obstruction to that.)}
\vspace{1mm}

\textbf {VII.} {\it Subgroup $G_n^0\subset G$ and half-Spin Representations.}  Let $\mathbb Z^n\to \mathbb Z_2$ be the (only) non-zero homomorphism,  which is invariant under   permutations of $e_i$,  denote by  
$deg:G_n \to \mathbb Z_2=\{-1,1\}$ be the composition of this with  the  homomorphism $G_n\to \mathbb Z^n_2$ which sends $c\to 1$ and 
let  $G_n^0$ be the {\it  kernel of this "degree" homomorphism.}

In terms of $Cl_n$, this is the intersection of  the subgroups $G_n$ and $Spin(n)$ in $Pin(n)$,
 $$G_n^0 = G_n\cap Spin(n)\subset Pin(N)\subset Cl_n.$$

{\it Exercise.} Show that $G_{n+1}^0$ is isomorphic to $G_{n}$.

{\it Hint.} Send $e_i\in G_n$, $i=1,...,n,$ to the products $e'_{n+1} e'_i$ for  $e'_1,....,e'_{n+1}\in G_{n+1}$.\vspace {1mm}

 Let $\hat e=e_1e_2...e_n$ and let us split the representation space $L=\mathbb C^{2^k}$ of $S_n$ for even $n=2k$ into $\pm 1$-eigenspaces of  $\hat e$,
 $L= L^+\oplus L^-$

If $n$ is even then  this $\hat e$ anti-commute with all $e_i$, that is $\hat ee_i=ce_i\hat e$.

It follows that, for  $n=2k$, 

\hspace {30mm}{\sf all   $ e_i$ that  act via  $S_n$ on $L$ send $\mathbb L^+\leftrightarrow L^-$}

and 
 
{\it the  restriction of  the representation  $S_n$  on $L=\mathbb C^{2^k}$ from the group $G_n$  to the 

subgroup   
$G_n^0\subset G_n$ 
 sends $L^+\to  L^+$ and $L^-\to L^-$.}\vspace {1mm} 
 
 Furthermore, since the representation $S_n$ is  irreducible for $G_n$, \vspace {1mm} 
 
 \hspace {19mm}{\sf the representations  $S_n^\pm$ on $L^\pm$ are irreducible for $G_n^0$.}\vspace {1mm} 
 
Extend these  representations by linearity to the subalgebra   $Cl_n^0\subset Cl_n$ generated by 
$G^0_n\subset Cl_n$, observe that $Cl_n^0$ contains  the  group $Spin(n)$ and restrict from 
  $Cl_n^0$ to $Spin(n)$. Thus, for $n=2k$,  we obtain \vspace {1mm}

  {\it two faithful  irreducible representations, called half-spin representations
$S^\pm$

 of the group $Spin(n)$ of dimensions $2^{k-1}$.}

\vspace {1mm}

{\it Remark/Question.} The above  shows that a linear space of dimension $ <2^k$ can't have $2k$ anti-commuting anti-involutions.
{\sl Is there a direct geometric  proof of this?} 

(The answer must be  known to some people.)
\vspace {1mm}

\textbf {VIII.} {\it Clifford's $Spin(n)$ Covers  $SO(n)$,}  What remains (for $n=2k$) to show is that this $Spin(n)$, which is  defined as the subgroup of the multiplicative group of the Clifford algebra $Cl_n$ generated by products of  even numbers of unit vectors $V\in V\subset Cl_n$,  double covers  the special orthogonal group  $SO(n)$. 

To do this we define an orthogonal (i.e. linear isometric) action of  all  of $Pin(n)\supset Spin(n)$ on the 
($n$-dimensional Euclidean) subspace  $V=V_{Cl}\subset Cl_n$ as follows.

Let $\alpha: Cl_n\to Cl_n$ be the automorphism that linearly  extends  $v\mapsto -v$ on $V\subset Cl_n $
and let
$$p(v)=\alpha(p)\cdot  v\cdot p^{-1}\mbox { for $v\in V$ and $p\in Pin(n)$}.$$

It is clear that if $p$ is a unit vector  in $V$, then the transformation $v\mapsto p(v)$
sends $V$  to itself by reflection in the hyperplane $p^\perp\subset V$ normal to $p$.
(You can think  of this $p\in Pin(n)$ as the square root of such a reflection.\footnote{If you omit
 $\alpha$, the resulting transformation square$v\mapsto pvp^{-1}$ becomes   {\it minus reflection} in
 $p^\perp$. Thus, if $n$ is odd, all  of $P(n) $ ends up in $SO(n)$. 
 
 Since one  wants $Pin(n)$ to cover 
 the full orthogonal group $O(n)$ one brings in this $\alpha$.})

Since $\alpha$ is an {\it automorphism} of the Clifford algebra, the map from $Pin(n)$ to the group 
$O(n)$, regarded as  the group of linear Euclidean isometries of $V=(V,\sum_ix^2_i)$,
is a {\it homomorphism} of groups, which sends $Spin(n)$ onto this $SO(n)$.

To conclude,  we need to  show  that the kernel of the homomorphism
$Pin(n)\to O(n)\subset End(V)$ is equal to $\{1,-1\}\subset Cl(n)$, which is done 
  by  induction on $n$  starting from  $Pin(1)=\mathbb Z_4= \{1, i, -1, i\} $ and $\alpha (i)=-1$,  
  and using  the following. \vspace{1mm}.
 
 {\it  Lemma}. {\sf If $\alpha(p)\cdot  v\cdot p^{-1}=v$ for a unit  vector 
 $v\in V$, then $p$ is contained in the  subalgebra $Cl(v^\perp) \simeq Cl_{n-1}$  generated by the 
 hyperplane $v^\perp\subset V$}. \vspace{1mm}
 
 {\it Proof.} Decompose  the Clifford algebra into  sum  of four linear subspaces,  
 $$Cl_n=A_0\oplus v\cdot A_1\oplus A_1\oplus v\cdot A_0,$$
 where $A_0\subset Cl(v^\perp)$ is equal to  the $+1$-eigenspace of $\alpha$,  i.e. where   $\alpha(a)=a$, and  $A_1\subset Cl(v^\perp)$ is the $-1$-eigenspace. 
 
  Observe that all $a_0$ in  $A_0$  are linear combinations of products of  of {\it even numbers} of vectors from $V$, while all   $a_1\in A_1$  are  combinations of {\it odd}  products.
 
 Now, by keeping track of parity of products we see that  the relation 
 
\hspace {-6mm}  $\alpha(p)\cdot  v\cdot p^{-1}=v$ divides into  two equalities,
 $$\mbox {$(a_0+v\cdot a_1')\cdot v=v \cdot (a_0+v\cdot a_1')$ and   $(a_1+v\cdot a_0')\cdot v=-v \cdot (a_1+v\cdot a_0')$},$$
 which  imply  that $a'_1=0$  and  $a_0'=0$.
 
Indeed, since  $v$ {\it commutes} with $a_0$ and anti-commute with $a_1$, 
 $$(a_0+v\cdot a_1')\cdot v=v \cdot (a_0+v\cdot a_1')\Rightarrow v\cdot a_1'\cdot v=v \cdot v\cdot a_1'\Rightarrow -v \cdot v\cdot a_1'=v \cdot v\cdot a_1', $$
and  $v \cdot v=-1\Rightarrow a_1'=0$. 

Similarly, one shows that also $ a_0' =0$ and {\sf lemma follows.}

 \vspace{1mm}
 
Finally,  we are through with {\it even} $n$:  \vspace{1mm}

{\sf \color {blue!50!black} the double cover group $Pin(n)\to O(n)$ for $n=2k$  comes with  a   faithful 

irreducible  complex representation $S_n=S_{2k}$ in $\mathbb C^{2^{2k}}$, called {\it spin representation}.\footnote{There in no faithful representation of $Pin(n)$ in a lower dimensional space, since even the  subgroup $G_n\subset Pin(n)$ admits no such representation.}} \vspace{0.6mm}

{ \sf \color {blue!30!red!30!black}The restriction of $S_n$ to $Spin(n)\subset Pin(n)$, that is  the double cover of 

$SO(n)\subset O(n)$,  splits into the sum $S_{n}=S_{\frac {1}{2}n}^+\oplus S_{{\frac {1}{2}n}}^-$ of two complex 

conjugate\footnote {We didn't prove these are complex conjugate but this follows from their construction} representations,   called {\it half spin representations}.}\footnote { Arriving at this point took   unexpectedly long -- 
not a page  or two as I had expected.} \vspace {1mm}

\textbf {IX.} {\it About Odd $n$.}    A quick way to arrive at the spin representation $S_{2k}$ of the group $Spin(n)$  in $\mathbb C^{2^k}$  for $n=2k+1$  is by   imbedding $Spin(n)\hookrightarrow Cl_{n-1}^\times$ and then restricting the  spin representation $S_{n-1 =2k}$   the Clifford algebra $Cl_{n-1}$ to the so embedded   $Spin(n)\subset Cl_{n-1}^\times$.

To achieve this,  we start, somewhat paradoxically, with  a (somewhat artificial) embedding $Cl_{n-1}\to 
Cl_n$ that sends $Cl_{n-1}$ onto the {\it even part }
$Cl^0_{n}\subset Cl_{n}$, that is the  $+1$-eigen space of  the  automorphism    $\alpha:Cl_{n}\to Cl_{n}$ of the Clifford algebra  induced by the central symmetry  $v\mapsto  -v$ of  the Clifford base  subspace $V=V_{Cl} \subset  Cl_n.$

It is (nearly) obvious that $Cl^0_n$ is a {\it subalgebra} in $Cl_{n}$  and that (this is slightly less obvious) this subalgebra is isomorphic to $Cl_{n-1}^0$.

To  prove the latter,  imbed $Cl_{n-1}$ to $Cl_{n}$ with the image $Cl^0_{n}$ as follows.

Map the orthogonal complement  $v^\perp \subset V\subset Cl_{n}$ of  a unit vector 
$v\in V$ back to $Cl^0_{n}$    by $e\mapsto v\cdot e$ for all  $e\in  v^\perp$ and show that this map  extends to an {\it injective  algebra homomorphism} $Cl_{n-1}=CL(v^\perp)\to Cl^0_{n}$.

All you need for  this is an (easy)   check up of the identities 
$$\mbox {$(v \cdot e)^2=-1$  and $v \cdot e \cdot v \cdot e'=-v \cdot e'\cdot v \cdot e$}$$
for all $v,v'\in v^\perp$ (implicit in   the above exercise about the  homomorphism  $G_n\to G^0_{n+1}$).
 
Finally,  since  that the group $Spin(n)$, by its very definition, is contained in $(Cl_n^0)^\times $ it goes to 
 $Cl_{n-1}^\times $ by inverting the  isomorphism $Cl_{n-1}\to  Cl_n^0$. QED.

\textbf {IX.} {\it Spin Representation of $Pin(n)$ fo Odd $n$.} Just for completeness sake, let us explain why 

{\sf the complexified  Clifford algebra $Cl_{2k+1}$, which has dimension $2^{2k+1}$,  is isomorphic to the  sum of of two matrix algebras 
$End(\mathbb C^{2^{k-1}})$.}

 Recall that the group $G_{2k+1}$ has exactly  two irreducible non-one-dimensional representations, where   the sum  of their  dimensions is
 $2^k$.
 
 In fact both representation must have the same  dimensions, because of 
 another     fundamental (also nameless?) theorem: 

{\it the dimensions of all irreducible representations of a finite group $G$  divide the order 
order of  $G$.}\footnote {See \url {https://math.stackexchange.com/questions/243221/proofs-that-the-degree-of-an-irrep-divides-the-order-of-a-group}  for several proofs.})
\vspace {1mm}

Therefore the non-Abelian part of the group algebra of $G_{2k+1}$, hence the Clifford algebra $Cl_{2k+1}$
is the sum of two matrix algebras of the same dimension. QED. 

As a consequence, we get\vspace {1mm}

{\it  two irreducible representations of the group  $Pin(2k+1) $ of dimensions $2^{k-1}.$}

\vspace {2mm}

\textbf {X.} {\it Example: Pauli "Matrices".} The first interesting case   of  $S_n$ is an irreducible 2-dimensional  complex  representation $S_2$ of the group $G_2$, hence of  $Pin(2)$, where he latter  is the non-trivial central $\mathbb Z_2$-extension of the circle group $\mathbb T^1=U(1)$.

To obtain such a  representation   all you need is to   find  two {\it anti-commuting anti-involutions}  
$\sigma_1, \sigma_2$ of $\mathbb C^2$  corresponding to the generators of $e_1,e_2$ of
 the (sub)group  $G_2\subset Cl_2\supset Pin(2)$.

This is kindergarten math:   let $\underline \sigma_1, \underline \sigma_2$ be  anti-commuting {\it involutions} of the real plane $\mathbb R^2$, namely reflections in two lines with the 45$^\circ$  between them. 
Their compositions, $\underline \sigma_1\underline \sigma_2$ and $\underline \sigma_2\underline \sigma_1$ are rotations by  90$^\circ$
 in the opposite directions, thus $\underline \sigma_1$ and $ \underline \sigma_2$ anti-commute:  
$$\underline \sigma_1\underline \sigma_2=- \underline \sigma_2\underline \sigma_1.$$
The {\it anti-involutions}  $\sigma_1=i\underline \sigma_1$  and $\sigma_2=i\underline \sigma_2$, $i=\sqrt{-1}$, of $\mathbb C^2$  with $\sigma_3=\sigma_1\sigma^2$ coming along are  your  Pauli guys.

\vspace {1mm}

{\it $\hat \otimes$-Remark.} This example can be amplified by taking tensor products, for 
$$Cl_{m+n}= Cl_m\hat \otimes Cl_n,$$
where $\hat \otimes$ stands for  {\it $\mathbb Z_2$-graded} tensor product,  for which
$$(a\otimes b)\cdot (a'\otimes b')=(-1)^{deg(a')\deg(b)}(a\cdot a')\otimes ( b\cdot b').$$

This allows a sleek construction of the spin representations but it  doesn't make it  more geometrical than 
the one via $G_n$.
\vspace {2mm}

\textbf {X.} {\it Clifford Moduli and Dirac  operators.} It is convenient at this point to call a linear space  $L$ with an action $S$ 
of  the Clifford algebra  $Cl(V)$  {\it "Clifford $V$-module}" and to  write just $S$ instead of  
$L=(L,S)$.

Also observe at this point that the actual action of  $V\subset Cl(V)$  on  such an $S$   reduces to a single  linear  map  $cl:V\otimes S\to S$,  where  the Clifford action is denoted by "$\cdot$",
 $$cl(v\otimes s) = v\cdot s.$$
 
Now, recall, that such a map  defines (and is defined by)   a first order  differential  on the space of smooth maps $\psi:V\to S$, denoted
$D: C^\infty(V,S)\to C^\infty(V,S)$, that is  the composition of this $cl$  with the differential
$d: C^\infty(V)  \to C^\infty (H)$  for $H=Hom (V, S)$ as we explained in the previous section.

Since  $v^2=-||v||^2$, $v\in V$, all $v\neq 0$   are {\it invertible} in the multiplicative semigroup 
 $End^\times(S)$; thus, \vspace{1mm}
 
\hspace {15mm}  {\sl the linear operators $D$   are elliptic for all $Cl_n$-moduli $S$.}\vspace{1mm}

These  $D$ can be also defined  with   orthonormal frames $\{e_i\}\subset V$ by
$$D(\psi)=\sum_{i=1}^ne_i\cdot \partial_i\psi,$$ 
which shows that  $D^2=-\Delta^2=-\sum_i \partial_i\partial_i$, since 
$$D^2=\sum_{i,j} e_i\partial_i e_j \partial_j=\sum_{i,j} e_i\cdot e_j\partial_i \partial_j=
\sum_i e_i\cdot e_i\partial_i\partial_j  +\sum_{i\neq j}(e_i\cdot e_j\partial_i \partial_j+ 
e_j\cdot e_i\partial_j \partial_i)= -\sum_i \partial_i\partial_i. $$

or, where   the symmetry is apparent,  by integration over the unit sphere 
$\{||v||=1\}\subset V$,
$$D(\psi(v))=const_n\int_{||v||=1} v\cdot \partial_v\psi(v)dv,$$
and if $V=\mathbb R^n$.
\vspace {1mm}

It follows by a simple symmetry  consideration or by a  one  line  computation that
$$D^2=-\Delta =-\sum\partial_{e_i}^2.$$

{\it Exercise.} Prove, directly that   
$$\int_{||w||=1} w\cdot \partial_wdw \int_{||v||=1} v\cdot \partial_vdv=
 const_n\int_{||v||=1} -\sum\partial_{v}^2dv.\footnote{I myself got lost in this calculation.}$$

\vspace {1mm}

{\it Dirac operator $\cal D$ on   Spinors.}   This $\cal D$  is defined with   the spinor representation $S_{2k}$ in $\mathbb C^{2^k}$, 
$$\mathcal D: \mathbb S_{2k} \to \mathbb S_{2k},$$
where the  "spinors" are understood here as smooth maps $\psi: \mathbb R^n\to S_{2k}$
for $n=2k$ or $n=2k+1$.

If $n$ is even, the spin representation splits into two adjoint representation, accordingly  $\mathbb S_{2k}= \mathbb S^+_{2k} \otimes \mathbb S^-_{2k}$, where the action of the Clifford algebra interchanges  
$\mathbb S^+_{2k} \leftrightarrow \mathbb S^-_{2k}$. 
It follows that $\mathcal D=\mathcal D^+\otimes \mathcal D^-$ for
$$\mathcal D^+ :\mathbb S^+_{k} \to \mathbb S_{k}^-\mbox { and }  
  \mathcal D^-:\mathbb S_{k}^- \to \mathbb S_{k}^+,$$
the operators  $\mathcal D^+$ and $  \mathcal D^-$ are {\it formally adjoint.}

\vspace {1mm}

 \textbf {XI.}  {\it $\cal D$ on Manifolds and Schr\"odinger-Lichnerowicz-Weitzenb\"ock-Bochner Formula.} Let $X$ be a Riemannian spin manifold of dimension $n$  and let $\mathcal S_{2k}$ be the spin bundle 
associated   with the principal spin bundle over $X$ that is the double cover of the orthonormal frame bundle, where this cover is what defines the spin structure on $X$.

Let $\nabla$ be the  Riemannian  Levi-Civita connection, which is, observe,   simultaneously  and coherently 
defined on all bundles associated with the tangent bundle. (It is irrelevant whether his is done via the principal $O(n)$-bundle or $Spin(n)$-bundle.)

We know  (this applies to all bundles with connections, see section\ref{connections3}) that  this $\nabla$ decomposes at each point $x\in X$ into the sum 
$\nabla=\nabla_{flat}+E_\nabla$, where $E_\nabla=E_{\nabla, x} $ a  smooth endomorphism of $\mathcal S_{2k}$ over a (small)  neighbourhood of $x\in X$, which {\it vanishes} at $x$. 

 This allows   a  "functorial transplantation"  of the above $\mathcal D=\mathcal D_{flat}$ to  an   $\mathcal D_\nabla$ on the space $\mathbb S$ of   sections of the bundle $\mathcal S_{2k}$,   where $\mathcal D_\nabla$  infinitesimally agree with $\cal D_{flat}$ at each point  $x\in X$,  
 $$\mathcal D_\nabla=\mathcal D_{flat}= E_\mathcal D,$$
for  a smooth (locally defined)  endomorphism  $E_\mathcal D=E_{\mathcal D,x }$ of $S_{2k}$, which {\it vanishes at $x$.}

If $n$  is even,  then, clearly,   $\mathcal S_{2k}=\mathcal S^+_{2k}\oplus \mathcal S^-_{2k}$ and the operator 
$\mathcal D_\nabla$, denoted just $\mathcal D$ from now on, splits accordingly: 
 $\mathcal D =\mathcal D^+\otimes \mathcal D^-$ for
$$\mathcal D^+ :\mathbb S^+_{k} \to \mathbb S_{k}^-\mbox { and }  
  \mathcal D^-:\mathbb S_{k}^- \to \mathbb S_{k}^+,$$
where the operators   $\mathcal D^+$ and $  \mathcal D^-$ are {\it formally adjoint.}

 Since the  $\mathcal D^2_{flat}=\mathcal D^2_{flat,x}$, which is defined locally, is equal  to $-\Delta=\nabla_{flat}\nabla^\ast_{flat} $ at each $x$,
the square  of $\mathcal D=\mathcal D_{\nabla}$,  now {\it globally}, satisfies what is called   "{\it Weitzenboeck   identity}" (this applies to all "geometric operators") 
$$\mathcal D^2  = \nabla\nabla^\ast+R_{ \mathcal D},$$
where $\nabla^\ast$ stands  for  the differential  {\it formally adjoint} to $\nabla$  (this spinor $\nabla$ acts from (sections of) $\mathcal S_{2k}$ to (sections of) the bundle $Hom(T(X), \mathcal S_{2k}$),   where
$R_{ \mathcal D}=R_{\nabla, \mathcal S,  \mathcal D}$ is a selfadjoint endomorphism of the bundle $\mathcal S_2k$.

It would be nice to continue this  line of this reasoning and see without calculation that,  why this $R_{ \mathcal D}$,  is a multiplication by a scalar.  Regretfully, I couldn't do this and have resort to the (standard) symbolic manipulations. \footnote {It  doesn't   feel right when you can't do   mathematics solely in your mind:     a piece of paper for this purpose  is no more satisfactory than a  digital  computer.}   

To perform these we, observe that the bundle of the Clifford algebras $Cl(T_x(X))$ acts on $\mathcal S_{2k}$, where this action agrees with the covariant differentiation $\nabla$ in $\mathcal S_{2k}$. 
Then we see that, for all orthonormal framed of tangent vectors $e_i$, $i=1,...,n$, the Dirac operator is
$$\mathcal D= \sum_i e_i\cdot \nabla_{i}  \mbox {  for  }   \nabla_{i}=  \nabla_{e_i}$$
and 

$$\mathcal D^2 =\sum_{i,j}e_i\cdot \nabla_{i}e_j\cdot \nabla_{j}=\sum_{i,j}e_i\cdot e_j\cdot \nabla_{i} \nabla_{j}=
\sum_{i =j} e_i\cdot e_j\nabla_{i} \nabla_{j}+
\sum_{i\neq j}e_i\cdot e_j\cdot \nabla_{i} \nabla_{j}=$$

$$=-\sum_{i} \nabla_{i} \nabla_{i}+ \sum_{i<j}e_i\cdot e_j\cdot( \nabla_i \nabla_{j}-\nabla_{j}\nabla_{i}) =
\nabla\nabla^\ast+ \sum_{i<j}e_i\cdot e_j\cdot R_\mathcal S(e_i\wedge e_j),$$
where $R_\mathcal S(e_i\wedge e_j)$ is the curvature of the bundle $\mathcal S_{2k}$ written as  a  2-form on $X$ with values in $End(\mathcal S_{2k})$.

The first term in this formula, $\nabla\nabla^\ast$ is the {\it Bochner Laplacian} in the bundle  $\mathcal S_{2k}$
which a selfadjoint non-strictly positive .

This $\nabla\nabla^\ast$, regarded as  a real operator,  is characterized by the integral  identity 
$$\int_X \langle \nabla\nabla^\ast\phi(x), \psi(x)\rangle dx-\langle \nabla\phi(x),\nabla\psi(x)\rangle dx=0$$
which is satisfied, whenever  one of the two functions has a compact support.

The proof of this formula, which makes sense and is valid for all  vector bundles with orthogonal connections, 
contains two  ingredients, where the first   {\it algebraic} one  consists in finding

 {\sf an invariant representation of the integrant as the  differential of  an $(n-1)$ form}

and the second ingredient is, of course, {\sf  Green's formula.}

In fact,  all  algebra  needed in our   is the following Leibniz formula for the Laplace Beltrami 
 $$\Delta  \langle \phi(x), \phi(x)\rangle = \langle \nabla\nabla^\ast\phi(x), \phi(x)\rangle+  \langle \phi(x), \nabla\nabla^\ast \phi(x)\rangle +2 \langle \nabla\phi(x), \nabla\phi(x)\rangle.$$
This  implies  all positivity of $\nabla\nabla^\ast$  we need.

Next,  turn to the
 curvature term $\mathcal R=  \sum_{i<j}e_i\cdot e_j\cdot R_\mathcal S(e_i\wedge e_j)$
 in the above  Bochner -- Weitzenb\"ock  formula for $\mathcal D^2$, that is 
 an  endomorphism  $\mathcal R :\mathcal S_{2k}\to \mathcal S_{2k} $, which, being self adjoint as a real , is  represented by a family of   {\it symmetric} linear operators $\mathcal R_x :(\mathcal S_{2k})_x\to (\mathcal S_{2k})_x$, $x\in X$, in the fibers  $(\mathcal S_{2k})_x\simeq S_{2k}=\mathbb C^{2^k}=\mathbb R^{2^{k+1}}$, while the curvature operators  $R_\mathcal S(v_1 \wedge v_2)$ themselves are {\it antisymmetric},  for all bivectors  
 $v_1\wedge v_2  \in   \bigwedge^2 T_x(x)= \bigwedge^2\mathbb R^n$,
 since they represent the action of the Lie algebra of the group $Spin(n) \subset SO(2 ^{k+1})$ on $\mathbb R^{2^{k+1}}$.
   
   In fact, a closer look shows\footnote {See formula  4.37 on p. p110 in  [Lawson\&Michelsohn(spin geometry) 1989].} that
 $$  R_\mathcal S (v_1\wedge v_2)=\frac{1}{2}\sum_{i<j}\langle R(v_1\wedge v_2)(e_i),e_j\rangle e_i\cdot e_j$$
where   $R(e_i\wedge e_j): T(X)\to T(X)$ is  the  curvature of our connection $\nabla$ as it acts on the tangent bundle of $X$.

(The  presence of  "$\frac {1}{2}$"  agrees with the idea of the bundle $\mathcal S_{2k}$   being a "the square
 root" of the tangent bundle $T(X)$, hence having one half of the curvature of  $X$, which is clearly  seen for  the Hopf complex line  bundle   $L\to S^2$, where  $L\otimes_\mathbb C L$  is isomorphic to the tangent bundle
$T(S^2)$ and, accordingly, $curv(L)=\frac {1}{2} curv(S^2)$.)

\vspace {1mm}
 
 Everything up to  this point was applicable to an arbitrary Euclidean  vector bundle $T\to X$ of rank $m$ with  a spin structure, i.e. a double cover of the associate principal $SO(m)$-bundle and 
 the action of bundle of the Clifford algebras $Cl(T)$ on the corresponding spin bundle with the fibers 
 $\simeq \mathbb C^{2^l}$, for $m=2l$ or $m=2l+1$, where  the Dirac  operator defined via an orthogonal connection in $T$ enjoys all  formulas we have presented so far.
 
 But in the case of $T =T(X)$  the  symmetries of the curvature tensor encoded  by Bianchi identities allow  the following   
 simplification of  $\mathcal R$.
\vspace {1mm}

{\it Lichnerowitz' Identity.} 
$$ \mathcal R=\sum_{i<j}e_i\cdot e_j\cdot R_\mathcal S(e_i\wedge e_j)=
 \frac{1}{2}\sum_{i<j, k<l}\langle R(e_k\wedge e_l)(e_i),e_j\rangle e_i\cdot e_j  =\frac {Sc}{4} ;$$  
Thus, 
$$\mathcal D^2\phi(x)=\nabla\nabla^\ast\phi(x) + \frac{1}{4}Sc(X,s)\phi(x)  \mbox{ for all sections }
\phi: X \to \mathcal S_{2k}. $$

{\it Why it is so.} The action of the linear  group $GL(n)$ on the space  $RCT \simeq \mathbb R^{\frac {n^2(n^2-1)}{12}}$
of (potential) Riemannian  curvature tensors splits into three irreducible representations
$RCT=Sc\oplus Ri\oplus W$, where $Sc$ is the trivial one dimensional representation, $Ri$ the space of traceless  symmetric 2-forms and  $W$ the space of Weyl tensors.
Accordingly,  every smooth $n$-manifolds $X$  supports 
 three (curvature)   differential operators    of the second order  from the space $G_+$ of positive definite quadratic differential forms $g$ on $X$ to the space of sections of vector bundles over $X$ associated with the tangent bundle  $T(X)$ via one of these representations, such that\vspace{1mm}
 
{\sf $ \bullet_{lin} $  these  operators are linear in the second derivatives of $g$;

$ \bullet_{inv} $ these operators are equivariant  under the  action of  the diffeomorphism group $diff(X)$}
operator
and where 

{\sf these operators and their scalar multiples are the only ones  with such quasi-linearity and invariance properties}   

On the other hand the  $\cal R$ is also constructed in a $diff(X)$-equivariant manner 
but it operators on the spinor  bundle $S_{2k} $, where the double cover of $GL(n)$ can't act.\footnote {Lemma 5.23. p 132 in [Lawson\&Michelsohn(spin geometry) 1989].}
This suggests that  there   is  {\it  no non-scalar  intertwining  } from   the space of curvature tensors on $X$ to the space of symmetric operators on $S_{2k} $, but  since  I didn't figure out  how to  prove  this without a few lines of  manipulations with Bianchi identities,  let us accept this for a fact.\footnote {Or see "Proof of Theorem 8.8" on  page  161 in [Lawson\&Michelsohn(spin geometry) 1989].} 

\vspace {1mm}

%%%%%%%%%%%%%%%%%%%%%

\subsubsection {\color {blue} Dirac Operators with Coefficients in  Vector Bundles,  Twisted 
S-L-W-B Formula and  $K$-Area }\label {K-area3}

%%%%%%%%%%%%%%%%%%%%%%%%

Let $\mathcal D_{\otimes L}$ be the  Dirac  twisted with a complex vector   bundle $L\to X$ with a 
unitary connection $\nabla^L$ on it.
Then,  as earlier, we have the general Bochner-Weitzenb\"ock formula 
  $$\mathcal D^2_{\otimes L}= \nabla\nabla^\ast + \sum_{i<j}e_i\cdot e_j\cdot R_\otimes(e_i\wedge e_j),$$
where this $\nabla=\nabla^\otimes$ is the connection in the tensor product of the spinor bundle with $L$ that is defined by the Leibniz rule, 
$$\nabla^\otimes(s\otimes l)= \nabla^\mathcal S  \otimes l+ s\otimes \nabla^L (l);$$
hence,  the curvature  $\cdot R_\otimes$ of this connection, that is the commutator of  
the $\nabla^\otimes$-differentiations,  also   behaves by this  rule:
$$R_\otimes (e_i\wedge e_j)(s\otimes l)=R_{\mathcal S} (e_i\wedge e_j)(s)\otimes l+l\otimes R_L(e_i\wedge e_j)( l)$$
which brings us to the following. \vspace {1mm}

{\large \color{blue}$\mathbf {[\mathcal D_\otimes]}$} \textbf {Twisted S-L-W-B Formula:}
 $$\mathcal D^2_{\otimes L}(\sigma\otimes l)= \nabla\nabla^\ast(\sigma\otimes l) +  \frac{Sc(X)}{4}(\sigma\otimes l) +\sum_{i<j}e_i \cdot e_j
 \cdot \sigma\otimes R_L(e_i\wedge e_j)(l). $$ 
 \vspace {1mm}

   {\it A basic  application} of this formula is the bound on the  area-size of manifolds with $Sc\geq \sigma>0$ expressed in terms of vector bundles  over $X$. \vspace {1mm}
  
  { \color {blue} [$K_\largestar$]} \textbf {Bound on  K-Area by Scalar Curvature.} {\sf Let $X$ be compact orientable  Riemannian Manifold 
 with positive scalar curvature and let $L\to X$ be a complex vector bundle with the unitary connection.}
 
 {\it If   the  norms of the curvature operators $R_x(e_1\wedge e_2):T_x(X_x\to T_x(X)$ of this connection are bounded by
 $$ ||R_x(e_1\wedge e_2):T_x(X_x\to T_x(X)||\leq \kappa_n\cdot Sc(X,x)$$
  for all $x\in X$, all unit bivectors $e_1\wedge e_2$  in  the tangent spaces $T_x(X)$ and a universal strictly positive constant $\kappa_n>0$,
 then, {\sf provided $X$ is spin}, all Chern numbers of the bundle $L$ {\color{red!30!black} vanish.}}

  \vspace {1mm}
  
  {\it Proof.} If some Chern number of $L$ doesn't vanish, then  an easy computation with Chern classes and    the index formula  shows\footnote {For details and further applications see  [GL(spin) 1980], \S 4-5 in chapter IV in   [Lawson\&Michelsohn(spin geometry) 1989], \S 4-5 in  [G(positive) 1996], [Min-Oo(K-Area) 2002] and sections  \ref{K-area3},   \ref {twisted4}, \ref{cowaist4}.}
   that  there exists  an associated bundle $L'$, such that the curvature $R'$ of the connectin in $L'$ satisfies 
  $$ ||R'_x(e_1\wedge e_2):T_x(X_x\to T_x(X)||\leq  const_n\cdot ||R_x(e_1\wedge e_2):T_x(X_x\to T_x(X)||$$
  and such that the {\it index of the twisted  Dirac  operator} on the spinor bundle tensored with $L'$,
   $$ \mathcal D^+_{\otimes L'}:\mathbb S^+\otimes L' \to  \mathbb S^+\otimes L',  $$ 
   {\it 
  doesn't vanish}. 
  
  But if 
   $$||R'_x(e_1\wedge e_2):T_x(X_x\to T_x(X)||< \frac {1}{4}\cdot  \frac {2}{n(n-1)} \cdot Sc(X,x).$$ 
   then, according to $ {\color{blue}[\otimes]}$
   the  $\mathcal D^2_{\otimes L'}$ is positive and the poof follows by contradiction.

\vspace {1mm}
 
 At  first sight this    { \color {blue} [$\largestar$]}  looks as an artifact of  symbolic manipulations with curvatures of   vector bundles, an insignificant   generalization of 
 the Lichnerowicz theorem, as devoid of an actual geometric information about $X$ as this theorem is.
 
  %\vspace {1mm}
  
  But,  surpassingly, although the proof of    { \color {blue} [$\largestar$]} is 90\% the same\footnote{The  proof of   { \color {blue} [$\largestar$]} , unlike that of  Lichnerowicz' theorem,  needs only 10\% of   the power of the Atiyah-Singer  theorem -- the  easy part of it:   non-trivial  variability of the index of $\mathcal D_{\otimes L}$ with variations of 
 (the  Chern classes of)  $L$, rather than a  more subtle aspect of the formula which involves 
 $\hat A$-genus of $X$.}
  as that by 
  Lichnerowicz, the information contents of the  two  statements are vastly different --  almost nothing in common   between  them: 
  
  Lichnerowicz is 99\% about {\it delicate smooth topological invariants} of manifolds with $Sc>0$, while  { \color {blue} [$\largestar$]}
   reveals raw geometric essence of   $Sc(X)\geq \sigma>0$, which, 
    as it becomes   a {\it positive curvature} condition,  { \it limits the size} of $X$.\footnote{Positivity of the sectional (and  Ricci) curvature, imposes bounds the {\it first and the second derivatives} of the  growths   of balls in  respective manifolds.}

  Below   is a specific instance of this.

 \vspace {1mm}
 
  \textbf{ Rough Area (non)-Contraction Corollary.} Given a {\it compact} Riemannian manifold $\underline X$, {\sf there  exists  a positive constant $\kappa=\kappa_{\underline X}>0$}, which restricts     how much  manifolds   
  $X$ with $Sc\geq \frac {1}{\kappa}$ {\it   can  be area-wise greater  than $\underline X$}, {\sf which is  expressed by}  
 {\it a bound on a  possible  decrease of areas} of surfaces in $X$ under {\sf  "topologically significant"  maps $X\to \underline X$.} \vspace {1mm}

 In precise language,

{[\Large $\star$}] {\sl let $X$ be  an oriented Riemannian manifold  with $Sc(X)>0$ and $f:X\to \underline X$ a smooth map, such that 
  the  norm of {\it the second exterior power} of the differential} of $f$,
  $$ \wedge^2df: \bigwedge^2 T(X)\to  \bigwedge^2 T(\underline X),$$
 is {\it bounded by the reciprocal of the scalar curvature of $X$ times $\kappa_{\underline X}$,}
 $$ ||\wedge^2df(x)||< \frac { \kappa_{\underline X}}{ Sc(X,x)},\mbox { for all    } x\in X.$$

{\it Then, {\sf provided $X$ is  \color{red!30!black}spin}, the image $h$ of of the fundamental homology class of $X$ in the homology of $\underline X$,
 that is $$\mbox {$h= f_\ast [X]\in H_n (\underline X)$, $n=dim(X)$,}$$ 
 is {\color{red!30!black} torsion}.}
  \vspace {1mm}
  
  {\it Proof.} By basic topology (a corollary to a theorem by Serre),     an {\it even dimensional 
  non-torsion}  homology class $h$ in $\underline X$ is "detected" by  a complex vector bundle: that is a  
  $\underline  L\to  \underline X$, such that  some characteristic cohomology  class  $\underline c$ of  $\underline  L$, doesn't vanish on $h$ $$\underline  c(h)\neq 0.$$
  
If  $h=f_\ast [X]$,
 then
  $f^\ast ( \underline  c)[X]$, which  serves as a {\it characteristic number} of the  induced  bundle $L=f^\ast(\underline  L)\to X$,  is equal to  $\underline  c(h)$; hence it  {\it doesn't vanish} either.
  
Now,  arguing as in the proof  of {\color {blue}$\mathbf {[Sc\ngtr 0]}$} for profinitely  hyperspherical manifolds (see  section \ref{twisted3}), 
 let   $\underline  \nabla$ be a unitary connection in $\underline  L$  and observe that  the
norm of the  curvature $R$ of the induced connection in $L$, which  is, after all, is a 2-{\it form},  is bounded by the curvature $\underline R$
of $\underline \nabla$, 
$$||R_x||\leq ||\wedge^2df(x)||\cdot ||\underline R_x||.$$
Thus, if  $n=dim(X)$  is even,  the proof follows from   { \color {blue} [$\largestar$]}  and the odd case reduces to the even one by taking the products of both manifolds with the circle.

  \vspace{1mm}
  
 {\it Remarks and Exercises.} \textbf {(a)} We use  the word  "{\it K-area}" to express the idea   that 
 
 {\sf  if $X$ contains   "homologically significant" families of surfaces with  {\it small areas}, then 
$K$-cohomology classes of $X$  can't be represented by  bundles  with connections, which have  small curvatures} 

\hspace {-3mm}and   
where

{\sf  the norm of 
 $\wedge^2 df$  measures by how much $f$ {\it contracts/expands these areas}.}\footnote {See  [G(positive) 1996], [Min-Oo(K-Area) 2002]  and  sections \ref {K-area3}, \ref {cowaist4} for more about this $K$-area.}
    
   Yet, we shall eventually switch to an uglier but more appropriate word "K-cowaist$_2$".

    \vspace{1mm}
 
 \textbf {(b)}  Let $X$ and $Y$ be closed oriented surfaces with Riemannian metrics on them and let $f_0: X\to Y$ be a continuous map of degree $d$.  Show that $f_0$ is homotopic to a smooth  {\sl strictly area decreasing} map $f$, i.e.  where $||\wedge^2 df(x)||<1$ for all $x\in X$, if and only if $area(X)>d\cdot area(Y)$.

 \textbf {(c)} The principal case in the above corollary, which yields    most topological  applications,  
 \footnote {See [GL(spin) 1980],   [GL(complete) 1983], [Lawson\&Michelsohn(spin geometry) 1989].}
 is  where $\underline X$ is the $n$-sphere $S^n$  and where the non-torsion condition amounts to {\it non-vanishing of the degree} of $f:X\to S^n$.

 In fact,  as one knows by a theorem of Serre,  the multiple of every    cohomology class $h$  in $\underline  X$ with $h\smile h=0$  can be induced from the the fundamental class of $S^n$ by a smooth map  $\underline  X \to S^n$, the general case of this corollary, for all dimensions, can be (with a minor effort)  reduced  $\underline  X = S^n$.

  \textbf {(d)} We call this corollary "rough", since the (lower) bound on $\kappa_{\underline X}$ its proof delivers   is far from optimal;
  
  Optimal bounds, however,  are available, albeit only  in a few cases, including $\underline X=S^n$ as we shall see in the following  sections.\vspace {1mm}
  
 {\it Questions. }\textbf {  (A) } {\sl  Is the spin condition in} {[\Large $\star$}]  {\sl redundant?} 
   
  {\sl  Or  the opposite is true}:  {\sf  if an  orientable  non-spin $n$-manifold $X$  admits a metric $g_0$ with $Sc(g_0) >0$, then it  carries  metrics $g_\varepsilon$, for all $\varepsilon>0$, 
    with $Sc(g_\varepsilon)\geq 1$, for which  allow   smooth maps $f_\varepsilon : (X, g_\varepsilon) \to  S^n$  with $deg( f_\varepsilon)\neq 0$, and 
    $||\wedge^2df|| \leq \varepsilon$?}
   
   \textbf {(B)} {\sl Can  the torsion conclusion  in} {[\Large $\star$}]  {\sl  be replaced by  {\it "$p$-torsion} for some particular  $p$, preferably for  $p=2$ 
 and,  in lucky cases, even  by just $f_\ast[X]\neq 0$?}

(It is not even clear   if  this can be done with a bound on $||df||$  rather than  on 
$\wedge^2||df||$, where  there is a chance for a successful use of  minimal hypersurfaces.)

%%%%%%%%%%%%%%%%%%%%%%%%%%%%%%%
\subsection{\color {blue} Sharp Lower Bounds on $sup$- and  $trace$-Norms of  Differentials  of  Maps  from Spin manifolds with $Sc>0$ to Spheres.} \label {trace-norms3}
%%%%%%%%%%%%%%%%%%%%%%%%%%

There is no single numerical invariant faithfully representing the size of $X$, but there are several  ways of  comparison the sizes of different  manifolds.

In the case, where  two Riemannian  metrics  are defined on the same background manifold, say $g$ and $\underline g$ on $\underline X$,  one compares these  at a point $ \underline x$ by simultaneously diagonalizing them and recording  
 the ratios of their values on the  vectors $e_i$ from the  common orthonormal frame $\{e_1,e_2,....,e_n\}\subset T_{\underline x}(\underline X)$, that are the numbers
 $$ \lambda_i(\underline x)= \lambda_i( \underline g/g,\underline x)= \frac{||e_i||_ {\underline g}}{ ||e_i||_g} .$$
 
In terms of these numbers,  the inequalities  $\lambda_i(\underline x)\leq 1$, $ \underline x\in \underline X$, say that $g\geq \underline g$, 
while the inequalities  $\lambda_i\lambda_j(\underline x)\leq 1$ convey that $g$ is (only)  {\it area wise} ({\sf non-strictly})  {\it greater} than $\underline g$,  where, of course, the former implies the latter.

Another way to compare the metrics  is by using  the {\it trace of $\underline g$ relative to $g$}, denoted  
$$trace (\underline g/g)=\sum_1^n\lambda_i, \mbox { } n=dim(  \underline X),$$ 
where the inequality 
 $$\frac {1}{n} trace (\underline g/g)\leq 1$$ 
expresses the idea of $g$ being greater than  $ \underline g$.

This    "trace-wise greater"  is 
less restrictive, yet,  moderately so, than  the "ordinary greater"  $g\geq \underline g$, for
 $$ g\geq \underline g \Rightarrow  g\underset {tr} \geq \underline g\Rightarrow  g\geq \frac {1}{n^2} \underline g. $$
 (Notice that  $\lambda_i(\underline g/c^2g)=\frac {1}{c}\lambda_i(\underline g/g)$.)
 
 A more relevant for us is the "area trace"
$$trace_{\wedge^2} (\underline g/g)=\sum_{i\neq j} \lambda_i\lambda_j$$
where "trace area-wise  greater" inequality reads 
 $$\frac {1}{n(n-1)} trace_{\wedge^2} (\underline g/g)\leq 1,$$
which is related to the "untraced area-wise greater" ratio  by the relations
$$[ g\underset {\wedge^2}\geq\underline g] \Rightarrow  [g\underset {tr_{\wedge^2}} \geq \underline g]\Rightarrow \left [g\geq \frac {1}{n(n-1)} \underline g\right]. $$

%%%%%%%%%%%%%%%%%%%%%
\subsubsection {\color {blue}Area Inequalities  for Equidimensional Maps:Extremality and Rigidity}\label{area extremality3}
%%%%%%%%%%%%%%%%

In order to apply the above  to Riemannian  metrics $g$ and $\underline g$ on  {\it different} manifolds $X$  and $\underline X$ we relate them  by  a smooth map, say $f:X\to \underline X$, where the  principal case is of $dim(X)=dim(\underline X)=n$ and where, to make sense of what follows, the map $f$ must be "{\it homotopically onto}",  that is {\it not homotopic to a map into a proper subset in   $\underline X$}.

 If both manifolds are    {\it orientable}  -- they are assumed  compact without boundaries at this point  --  this is equivalent to {\it non-vanishing of the degree $deg(f)$} of the map,
\footnote{The implication  [$deg(f)\neq 0]  \Rightarrow [f$  is   homotopically onto],  which  is obvious by the modern standards,  is by no means trivial.  For instance, "homotopically onto" for the identity  map of the $n$-sphere is equivalent (one line kindergarten argument) to the Brauer fixed point theorem for the $(n+1)$-ball.} 

If non-orientability is easily taken care of  by just passing to orientable  double covers, what does  cause a problem is the {\it spin condition}, the relevance of which  the following two geometric theorems  remains problematic. 

\vspace {1mm}
 {\color {blue} \large  [$X_{spin}{^\to}$\Ellipse]}  {\textbf{Spin-Area Convex  Extremality Theorem.}}  {\sf Let $\underline X\subset \mathbb R^{n+1}$ be 
  a smooth  compact convex hypersurface and let $\underline g$
   be the Riemannian metric on $\underline X$ induced from $\mathbb R^{n+1}$.
   Let $X=(X,g)$ be a  compact  orientable Riemannian $n$-manifold with $Sc\geq 0$ and let
  $f:X\to  \underline X$ be a smooth map of {\it non-zero degree.}

    Let $g^\circ=Sc(g)\cdot  g$ and  $\underline g^\circ= Sc(\underline g)\cdot \underline g$ be the corresponding $Sc$-normalized  metrics}\vspace {0.6mm}

{ \it If $X$ is {\sf spin} and $n$ is {\sf even}, then the map $f$  can't be   strictly area decreasing, that  
 is the metric $g^\circ$ is  {\color {blue!49!red}not} area-wise greater, than the induced metric $f^\ast (\underline g^\circ)$
 on $X$.}
 
Put it another way, 

{\it there necessarily exists a point $x\in X$, where the norm of the second exterior  power  of the  differential of $f$ is bounded from below by the scalar curvature of $X$ as follows
 $$Sc(\underline X, f(x))\cdot  ||\wedge^2df(x)||\geq Sc(X,x), $$
which, in terms of $\lambda_i^\circ=\lambda_i(f^\ast (\underline g^\circ))/g^\circ)$,
 reads
$$\max_{x\in X,i\neq j}  \lambda_i^\circ(x)  \lambda_j^\circ(x)\geq 1.$$}

\vspace{1mm}

In the simplest case, where  $\underline X$ is  the unit sphere  $S^n\subset \mathbb R^{n+1}$,   this theorem can be refined as follows.
\vspace{1mm}

 {\color {blue} \large  [$X_{spin} \to\Circle]$} {\textbf{Spherical  Trace Area Extremality Theorem}.} {\sf Let $ X$ be a compact orientable Riemannian spin manifold of dimension $n$ and $f:X\to S^n=\underline X$ be  a map  with  $ deg(f)\neq 0$.}

 Then   $f$  {\it can't be {\color {blue!50!red} trace} area-wise  strictly decreasing  with respect to the $Sc$-normalized  metrics  $g^\circ = Sc(g)\cdot g$ on $X$ and 
$\underline g^\circ= Sc(\underline g)\cdot \underline g^\circ= n(n-1)ds^2$, which,  in terms of the exterior power of $f$, says that 
there is a point $x\in X$, where the trace-norm of the second exterior  power  of the  differential of $f$ is bounded from below by the scalar curvature of $X$ as follows
 $$2  ||\wedge^2df(x)||_{trace}\geq Sc(X,x), $$
that is  
$$\frac{1}{2n(n-1)}\sum_{i\neq j}  \lambda_i^\circ(x)  \lambda_j^\circ(x)\geq 1\mbox {
for $\lambda_i^\circ=\lambda_i(f^\ast (\underline g^\circ))/g^\circ)$.}$$}

\vspace {1mm}

 {\it Remarks} (a) Neither  {\color {blue} \large  [$X_{spin}{^\to}$\Ellipse]} nor   {\color {blue} \large  [$X_{spin} \to\Circle]$}  seem obvious even, where $X$ is also a  convex hypersurface in $\mathbb R^{n+1}$.

{\it Question.}  Are  there  counterparts  of  {\color {blue} \large  [$X_{spin}{^\to}$\Ellipse]} and/or of    {\color {blue} \large  [$X_{spin} \to\Circle]$} for symmetric function  $s_k$  of the principal  curvatures $\alpha_1,\alpha_2, ..., \alpha_n$  of convex hypersurfaces $X$ and $\underline X$?
(We shall return to this question in (b) of \ref{mean convex3}.)

(b) The condition  $n=2k$, which is unneeded for  {\color {blue} \large  [$X_{spin} \to\Circle]$}, {\color {red!40!black}probably} is also redundant for  {\color {blue} \large  [$X_{spin}{^\to}$\Ellipse]}.

(c) These two theorem  will  be  later generalized in several directions.\vspace {1mm}

(d$_{1}$) One may allow non-compact, and sometimes even non-complete manifolds $X$ with  suitable conditions 
on maps  $f$, in order to have their degrees being  properly defined.

(e$_{2}$) In the case, where $dim(X)= dim(\underline X)+4l$,  
the condition  $deg(f)\neq 0$ can be replaced by $\hat A[f^{-1}(x)]\neq 0$  for a generic point  $x\in X$  of a smooth map $f: X\to   \underline X$.   
(c$_{3}$) Instead  of a convex hypersurface in  $ \mathbb R^{n+1}$,  one may take  a more general  Riemannian manifold for $\underline X$, namely   one with   {\it a non-negative curvature operator} and -- this is, {\color {red!40!black} probably}, unnecessary --  {\it with non-zero Euler characteristic.}    
 
\vspace {1mm}

(f) {\sf Who  is extremal?} These two  extremality theorems   can be thought of as {\it properties of $X$}, saying that "{\it large scalar curvature makes $X$ small}".  

From another perspective,  {\it these theorems  are  about $\underline X$}, saying that {\it $\underline X$ can't be enlarged without making its scalar curvature smaller at some point.} 

This suggest  two avenues of generalizations  that we shall explore in the following sections. 

1. Widen  the  class of manifolds $X$ and maps $f:X\to \underline X$,  which satisfy the above or similar theorems  
and, regardless of the scalar curvature,  study   invariants of  manifolds $X$ responsible for  existence/non-existence of metrically contracting,  yet  topologically significant, maps from $X$ to "standard" manifolds  $\underline X$
such as the spheres,   for instance.

2. Find further  instances   of {\it extremal} manifolds $\underline X=(\underline X,\underline g)$ with  
$Sc(\underline g)>0$, i.e.  where   no Sc-normalized  metric $g$ can be greater the so normalized $\underline g$,
$$Sc( g)\cdot g  \ngtr  Sc(\underline g)\cdot \underline g$$
and study  properties of such metrics.\footnote {See   [Sun-Dai(bi-invariant)2020]  for the proof of the extremality of
 {\it bi-invariant metrics} on compact Lie groups in the class of {\it left invariant} metrics.}%{\color {red}(new !!!!!!!!)}}

\vspace{1mm}\vspace{1mm}

{\it A few  Words about the Proofs.} \footnote {For detailed  poofs the  above mentioned  results  see [Llarull(sharp estimates)  1998],  [Min-Oo(Hermitian) 1998],   [Min-Oo(K-Area) 2002],    [Goette-Semmelmann(symmetric) 2002],    [Goette(alternating torsion) 2007], [Listing(symmetric  spaces) 2010];  also  we   say a bit  more  about  this   in sections\ref{Dirac4}, \ref
{twisted4}.}.  The logic here is the same as in the proof of  the rough area (non)-contraction corollary  from the previous section, where the sharpness of the bound on $\wedge^2df$ 
is achieved by a  choice of the bundle $\underline L\to \underline X$ with a non-zero top Chern class 
 with a connection $\underline \nabla$ with  minimal possible curvature, that  allows the necessary strong  bound on the "twisted curvature" term 
 $\sum_{i<j}e_i \cdot e_j
 \cdot \sigma\otimes R_L(e_i\wedge e_j)(l)$ in the Schr\"odinger-Lichnerowicz-Weitzenb\"ock-Bochner formula for the Dirac  operator on $X$ tensored with induced connection $\nabla=f^\ast(\underline \nabla)$
 in the bundle $L=f^\ast(\underline L)\to X$,
  $$\mathcal D^2_{\otimes L}(\sigma\otimes l)= \nabla\nabla^\ast(\sigma\otimes l) +  \frac{Sc(X)}{4}(\sigma\otimes l) +\sum_{i<j}e_i \cdot e_j
 \cdot \sigma\otimes R_L(e_i\wedge e_j)(l).$$

The natural choice of $\underline L$  -- this was suggested by Blaine Lawson  40 years 
ago --\footnote {I recall this well, since I was taken by surprise by the properties of this bundle,  which has the minimal curvature (one half of that of the tangent bundle of the sphere) among all unitary bundles with non-trivial Euler class.})  is one of the 
{\it Bott generator bundles}, that are  the  $\frac{1}{2}$-spinor bundles $\underline L^\pm=\mathcal S^\pm(\underline X)$ (with
$rank_\mathbb C (\underline L)=2^{k-1}$  for $n=2k$),   which, being the "moral square roots" of the tangent bundle $T(\underline X) $, have  their curvatures equal 
to the one half of that of  $T(\underline X) $. (This   is  clearly   seen for $n=2$ where 
$\underline L^+$ is the Hopf complex line bundle over $S^2$.

What makes   $\underline L^\pm$  promising candidates for  S-L-W-B-extremality, is the fact that
$\underline  L^\pm$-twisted Dirac  operator on the manifold $\underline X$ itself {\it does have harmonic spinors} but only  {\it barely so}: these spinors are {\it parallel} as they correspond to {\it constant functions} and/or to  {\it constant multiples of the Riemannian volume $n$-form} on $\underline X$.

The extremality property of  $\underline L^\pm$  was confirmed by Llarull in the case of $\underline  X=S^n$ and  -- this  was by no means  expected  --
by Goette and Semmelmann for manifolds $\underline  X$ with positive curvature operators, 
while the possibilities of $Sc$-normalization and of tracing  $\wedge^2df$,  were suggested by Listings.
(Although there is {\it no technical novelties in the proofs} of the Sc-normalised and traced modifications  of 
 {\color {blue} \large  [$X_{spin}{^\to}$\Ellipse]}  and  {\color {blue} \large  [$X_{spin} \to\Circle]$} 
these  significantly {\it widen the range of applications} of these extremality theorems.)

Besides facing algebraic complexity of the "twisted curvature" one has to ensure  the existence of  non-zero $L^\pm$-twisted harmonic spinors on $X$ for  $ L^\pm=f^\ast \underline L^\pm$.
 
 The index  formula
 guarantees this for $n=2k$ and, under  an  additional condition on $f$, also for $n=4l+1$, 
 but in general the  existence   of such  spinors  for all metrics on $X$  and all $n$ remains {\color{magenta}
\sf problematic.}

\vspace {1mm}

{\color {blue} $\Circle_\blacksquare .$} {\it The Proof of} {\color {blue}  $[X_{spin} \to\Circle]$} {\it for odd $n=dim(X)$.}  Given a map $X\to S^n\subset S^{n+1}$, radially (and obviously) extend it to the map $X\times \left[-\frac {\pi}{2},\frac {\pi}{2}\right]\to S^{n+1}$ with the bottom and the top  of the cylinder $X\times \left[-\frac {\pi}{2},\frac {\pi}{2}\right]$ sent to the poles of $S^{n+1}$,
$X\times \left\{\mp\frac {\pi}{2}\right\}\to \mp 1.$
 
One can proceed three ways from this point. 
 
1.   Endow $X\times  \left[-\frac {\pi}{2},\frac {\pi}{2}\right]$ with the (spherical suspension)  warped product metric $\hat g$
with the same warping factor as that for the spherical cylinder
 $S^{n+1}\setminus \{-1,+1\}$  and observe  that, say   in the case of 
$Sc(X)\geq n(n-1)=Sc(S^n)$, this  metric has greater scalar curvature than that of $S^{n+1}$. 

Then,  by an easy argument,  an  $\varepsilon$-small  $C^0$-perturbation of this metric $\varepsilon$-near the boundary extends,  for all $\varepsilon>0$,  to a complete metric $\hat g_\varepsilon$  on  the infinite  cylinder $X\times (-\infty,+\infty)$, such that $Sc(\hat g_\varepsilon)\geq n(n+1)-\varepsilon$ and such that the geometry of $X\times (-\infty,+\infty), g_\varepsilon)\geq n(n+1)-\varepsilon$ is cylindrical for $|t|\geq \frac {\pi}{2}+\varepsilon$ infinity 
with the scalar curvature $\geq n(n+1)+1.$

Thus the 
 untraced inequality {\color {blue}  $[X_{spin} \to\Circle]$}  applies to the  product $S^n\times S^1(R)$  $R\geq 2$ obtained by closing this cylinder at infinity and letting $ \varepsilon\to 0$.
 
2. Apply the traced inequality   {\color {blue}  $[X_{spin} \to\Circle]$} to maps $X^n\times S^1(R)\to S^{n+1}$, where $S^n\times S^1(R)$ comes with  with the product metric,  and let the radius  of the circle  $R\to\infty$. (This is, essentially, how it was done in [Llarull(sharp estimates) 1998].)

3. Regard a map  $X^n\times S^1\to S^{n+1}$ of non-zero degree as a family of maps $f_s: X\to 
S^{n+1}$ and 
 use the spectral flow index theorem for the family of operators on $ X=X\times s$  parametrized by $S^1$.\footnote{Such  argument was used in [Vafa-Witten(fermions) 1984]  for lower bounds on  spectral gaps  for the Dirac  operator, succinctly exposed   in   [Atiyah(eigenvalues) 1984]  and applied in \S 6 
in [G(positive)  1996]
to spectral  bounds for the   Laplace operators on odd dimensional Riemannian manifolds.

Also spectral flow for Dirac operators combined with a {\it refined Kato inequality} is used in   [Davaux(spectrum) 2003]  for the proof of  sharp upper bounds  on  the scalar curvatures of Riemannian metrics on compact  manifolds which admit  hyperbolic metrics.}

 . 
 \vspace{1mm}
 
 {\it Exercise.} Fill in the  details in (1) and (2).

 \vspace{1mm}

{\it Question} {\sf Is there a more direct  ($K^1$-theoretic?) proof of the inequality   {\color {blue}  $[X_{spin} \to\Circle]$}  for odd $n$ with no direct  reference to $S^{n+1}$ and desirably of  {\color {blue}  [$X_{spin}{^\to}$\Ellipse]} as well for, odd $n$, e.g. by   a spectral  flow argument?} \vspace{1mm}

 {\it Infinite Dimensional Remark.} Both,  spherical suspension in 1 and the cylindrical one in 2, when  repeated 
N-times  times and can be interpreted in the limits for $N\to \infty$  as properties of 

$1^\infty$  infinite dimensional manifolds 
$X^\infty $ with $Sc(X^\infty)\geq Sc (S^\infty)$;

$2^\infty$  $Sc(X^\infty)\geq Sc(S^n\times \mathbb R^{\infty-n}$;

 inequalities  are implemented in both cases by certain special Fredholm-type maps $X^\infty \to S^\infty$.

Conversely, one can prove an infinite dimensional version of  {\color {blue}  $[X_{spin} \to\Circle]$} for
limits of the above maps, say for 

{\sf Fredholm maps 
from a Hilbertian manifold $X$ to the Hilbertian sphere, $f:X\to S^\infty$,such that $deg(f\neq 0$ and such that 
there exists a sequence of equatorial spheres 
$$S^{N_1} \supset S^{N_2}\supset ... \supset S^{N_i}\supset ...\supset S^\infty,$$ 
 where the union $\bigcup_i S^{N_i} $ is {\it dense} in $S^\infty$ and such that 
 the pullbacks $X_i=f^{-1}(S^{N_i}) \subset X$ are smooth submanifolds of dimensions $N_i$, the scalar     curvatures of which with the induced metrics 
satisfy $Sc(X_i)-N_i(N_i-1)-0$  for $i\to\infty.$}

{\it Infinite Dimensional  Questions.} What is the most  general/natural infinite dimensional inequality   {\color {blue}  $[X_{spin} \to\Circle]$}?

Is there a direct proof of such an inequality with no use of finite dimensional  approximation?

Are there natural  Hilbertian and/or non-Hilbertian spaces $X$ to which such an inequality may apply?
\vspace{1mm}

{\it Stability Remark.}  Probably, (I haven't thought trough this)  the reduction argument $even\leadsto odd$ implies certain stability of harmonic spinors  on 
$(2m-1)$-manifolds $X$  twisted with spherical   spinors, that are section of the induced bundle $ f^\ast(\mathbb S(S^{2m-1}))$  by  maps $f:X\to S^{2m-1}$
with $deg(f)\neq 0$.

Another    (seemingly unrelated)  instance of   {\it stability  of harmonic spinors (seemingly)  independent  of the index theorem}  is present in 
 
 {\sf Witten's   argument in his proof of the Euclidean positive mass theorem as well in Min-Oo's proof of the hyperbolic one.}
 
  {\it \large Probably,}  {\sf there are many examples of  stable (twisted) harmonic spinors on {\it compact manifolds}, where this  stability  is not  not predicted,  at least not directly, by the index theorem.}\footnote {To make sense of this  one has to properly specify  the meaning of  "stability"   not to run into (counter) a example, see  {\sl Harmonic Spinors and Topology} by Christian B\'ar,\url{https://link.springer.com/chapter/10.1007/978-94-011-5276-1_3}}
\vspace {1mm}

{\it \textbf {Area Rigidity Problem: Examples and Counter Examples.}} Given a smooth convex hypersurface $\underline X\subset \mathbb R^{n+1}$ and let $\underline g$ be the induced Riemannian metric on $\underline X$.

 {\sf Describe (all) Riemannian$n$-manifolds $X=(X,g)$   along with smooth maps  $f: X\to \underline X$, such that  
$$Sc(g,f(x)) \leq Sc(\underline g,x)\cdot ||\wedge^2df(x)||$$
at all $x\in X$ and also $X$ and $f$ where 
$$Sc(g,f(x)) =Sc(\underline g,x)\cdot ||\wedge^2df(x)||.$$}

 In the "ideal rigid" case, at least for $Sc(\underline X)>0$,   one wants all such maps to be   {\it locally isometric with respect to  the $Sc$-normalised metrics} $g^\circ=Sc(g)\cdot g$  and 
  $\underline g^\circ =Sc(\underline g) \cdot\underline g$.
(This,  if I am not mistaken,  is the same  as {\it local homothety} with respect to the original metrics: 
  the induced Riemannin metrics $f^\ast( \underline g)$ on $X$  are constant multiples of  $ g$, i.e.
$g=\lambda \cdot f^\ast( \underline g)$)

But the true  picture is  more interesting than  this "ideal". 
 Here is what one can say in this regard.

\vspace{1mm}

(A) If $n$=2 then the equality  $Sc(\underline g,f(x)) =Sc(g,x)\cdot ||\wedge^2df(x)||$ says that
$f$ is {\it locally area preserving} with respect to $ g^\circ$ and $\underline g^\circ$; hence, the space of 
such maps is (at least) as large as the group of area preserving  diffeomorphisms of the disc.

(B) If $n\geq 3$, then locally area preserving maps are locally isometric and, in fact,

{\it \textbf  "{{\color {blue!60!black} Ideal rigidity}}", i.e. the implication
$$Sc(\underline g,f(x)) \leq Sc(g,x)\cdot ||\wedge^2df(x)||\Rightarrow g=\lambda \cdot f^\ast( \underline g),$$
 was proven by Mario Listing} 
 {\it under the following assumptions:}\footnote {See theorem 1 in  [Listing(symmetric spaces) 2010], and compare with Theorem 4.11 in [Llarull(sharp estimates) 1998].}

  {\sf  $\bullet$ $\underline  X$  is a closed {\it strictly} convex hypersurface  of dimension
$n\geq 3$,  where  this  "strictly" signifies that all principal curvatures are $>0$ (rather than  non-existence of straight segments  in   $\underline  X$);

$\bullet$  $X$ is a closed connected orientable {\it spin} manifold and    $deg(f)\neq 0$.} \vspace{1mm}

Now let us look at {\it non-strictly convex} hypersurfaces of dimensions $n\geq 3$.\vspace{1mm}

(C) Let a hypersurface  $\underline  X\subset \mathbb R^{n_0+m}$ be the product
$$\underline  X=\underline  X_0\times \mathbb R^m$$ 
where $\underline  X_0\subset \mathbb R^{n_0}$ is a smooth hypersurface.
 Then all (self) maps 
 $$f =(f_0, f_1): \underline  X \to \underline  X_0\times \mathbb R^{m}=\underline  X,$$
such that $||df_1||\leq 1$,  satisfy
 $Sc(\underline g,f(x)) =Sc(g,x)\cdot ||\wedge^2df(x)||.$\vspace{1mm}

If $m\geq 2$,  there are no {\it closed}  Euclidean  hypersurfaces displaying such non-rigidity
(unless I am missing obvious Euclidean examples)\footnote{There are these in $\mathbb R^{n_0}\times \mathbb T^{m}$.} but 
this   non-rigidity,  of cylinders, i.e. for $m=1$,   can be cast  into
   a  compact form; also   this  
can be done to conical hypersurfaces as follows.\vspace{1mm}

(D) Let $C\subset \mathbb R^{n+1}$   a smooth convex cone and let $\underline  X \subset C$
be a smooth  closed convex  hypersurface, such that the intersection  $\underline  X \cap \partial C$
contains a  conical annuls $A$ in the boundary of $C$  pinched  between two spheres,
$$A= \{a\in \partial C\}_{R_1\leq ||a|| \leq R_2}.$$
Thus, the boundary of  $\underline  X \subset C$ consists of three parts:

{\sf the {\it side boundary} that is the intersection $\underline  X\cap \partial C$;

{\it bottom} $\underline  X_1\subset   \underline  X$ that lies  on the $R_1$-side  in  the interior of   $C$,  i.e.
 $||\underline  x_1||< R_1$,  for $\underline  x_1 \in \underline  X_1$,
  
{\it top} of $\underline  X_2\subset   \underline  X$ that lies  on the $R_2$-side  in  the interior of     $C$, i.e.
$||\underline  x||> R_2$   for $\underline  x_2 \in \underline  X_2$.}

Scale up the top of $\underline X$ and set:
$$X=(\underline X\setminus \underline  X_2)\cup\lambda X_2, \mbox { } \lambda>1.$$

This $X$ admits an obvious  (infinite dimensional) family of   diffeomorphisms  $f: X\to  \underline  X$, that 

fix the bottom,   

return
back the top by $x\to \lambda^{-1}x$, 

send all straight  radial segments in the side boundary of $X$
to  themselves, 

 satisfy  the equality  $Sc(\underline g,f(x)) =Sc(g,x)\cdot ||\wedge^2df(x)||.$
 \vspace {1mm}
 
{\color {red!40!black}Probbaly,} (C) and (D) give a fair picture of possible kinds of   {\it not-quite-rigid} $\underline X$ with $Sc(\underline X) >0$   {\it  in  the class of convex} $X$, but it is not so clear for the class of all $X$ with $Sc(X)>0$

 %%%%%%%%%%%%%%%%%%%%%%%%%%%%

\subsubsection {\color {blue} Area Contracting Maps with Decrease of Dimension} \label{decrease of dimension3}

%%%%%%%%%%%%%%%%%%%%%%%%%%%%%

The  lower bounds on the norms  $||\wedge^2 df||$ for  equividimensional maps $f:X\to \underline X$  with {\it non-zero degree} generalize to maps,  where $dim(X)>dim( \underline X)$ with an appropriate generalization of the concept of
 degree. 
 
 For example, the proofs  of the rough Area (non)-contraction property (section \ref{K-area3}) and of   both its above 
  refinements {\color {blue}  [$X_{spin}{^\to}$\Ellipse]} and    {\color {blue} [$X_{spin} \to\Circle]$},
which say that  such norms can't be too small at all points in $X$, \vspace{1mm}

 \hspace {5mm}$||\wedge^2df(x)||\nless \frac {Sc(X, x)}{Sc(\underline X, f(x)||}$ and $
 ||\wedge^2df(x)||_{trace}\nless 2\frac{Sc(X,x)} {n(n-1)} $ correspondingly,\vspace{1mm}
 
 \hspace {-6mm} extend with (almost) no change to maps $f: X^{n+4l}\to  \underline X^n$ with {\it non-zero
 $\hat A$-degrees}, which means non vanishing of the $\hat A$-genera of the pullbacks 
 $f^{-1}(\underline x)\subset X^{n+4l}$ of generic points $\underline x \in \underline X^n$.
\footnote {This is done in   
 [GL(spin)1980],  [Llarull(sharp estimates)  1998], [Goette-Semmelmann(symmetric) 2002]  and in  [Goette(alternating torsion)2007] for bounds on $||\wedge^2df||$, but the corresponding lower bound on    $||\wedge^2df||_{trace}$   is missing from  [Listing(symmetric  spaces) 2010]; however, as I see it, there in no problem with this either.  }
 
 For instance: \vspace{1mm}

{\color {blue} \large  [$X_{spin}\overset{\hat A} \to\Circle]]$} {\textbf{  $\hat A$-Extremality Theorem}.} {\sf Let $ X$ be a compact orientable Riemannian spin manifold of dimension $n+4l$ and 
$f:X\to\underline X= S^n$ be  a  smooth map, such that the $\hat A$-genus of the  $f$-pullback of a regular point from $S^n$ doesn't vanish, 
$$\hat A[f^{-1}(\underline x_0)]\neq 0, \mbox { }\underline x_0\in S^n.$$.}

{\it Then there exists a point $x\in X$, where the  the trace-norm of the second exterior  power  of the  differential of $f$ is bounded from below by the scalar curvature of $X$ as follows,
 $$2  ||\wedge^2df(x)||_{trace}\geq Sc(X,x). $$}

Since $$2  ||\wedge^2df(x)||_{trace}=\sum_{i\neq j}\lambda_i(x) \lambda_j(x)\leq n(n-1)\max_{i\neq j}\lambda_i(x) \lambda_j(x) =n(n-1) ||\wedge^2df(x)||, $$
this implies that if $Sc(X)\geq n(n-1)=Sc(S^n)$, then  the map  $f$ {\it can't be strictly area decreasing}.

 \vspace {1mm}

\textbf {Generalization to  $\mathbf { \hat \alpha}$}. The above remains true with $ \hat \alpha$ instead of  $\hat A$, e.g. 
where the  pullback of a regular point  $ f^{-1}(\underline x_0)\subset X$ is diffeomorphic to Hitchin's exotic 
sphere $\Sigma^n$ for $n=8k+1, 8k+3$.\footnote{Such a  $\Sigma^n$ is homeomorphic to the ordinary sphere  $S^n$, but  \it  doesn't bound a  spin manifold.}

 \vspace {1mm}

{\it \textbf{ Question}.} Does the conclusion of the above theorem remain true if the nonvanishing of $\hat A[f^{-1}(\underline x_0)]$ is replaced by the following

  \vspace{1mm}

{\sf the   pullbacks  $(f')^{-1}(\underline x)$,   for all smooth maps $f': X^{n+m}\to  \underline X^n$ homotopic to $f$ and

 all $f'$-non-critical $x\in\underline X^n$, {\it admit no metrics with $Sc>0$.}} \vspace{1mm}

This is  beyond the present day techniques,  already for manifolds 
$X^{n+m}$  homeomorphic to $S^n\times Y^m$, where $Y^m$ is SYS-manifold.

 But if $Y^m$ is {\it the torus} or, more generally an  {\it enlargeable manifold}, e.g. if it admits    {\it a metric with non-positive sectional curvature}, then Dirac theoretic techniques on {\it complete manifolds} (see sections  \ref{relative3},
 \ref{twisted4}) delivers the proof of the following.
 
\vspace {1mm}

{\color {blue}$\bigtimes\mathbb R^m$}-\textbf { Stabilized Mapping Theorem.} {\sf Let $X^{n+m}$ be a  complete orientable Riemannian spin manifold with 
$Sc( X^{n+m})\geq \sigma >0$ \footnote {In view of [Zhang(Area Decreasing) 2020], one can, probably, relax this to $Sc( X^{n+m})\geq 0$.}
 and let  $\underline X^n$ be a smooth convex hypersurface in $\mathbb R^{n+1}$.
   Let 
 $f_1: X^{n+m}\to  \underline X^n$ and $f_2:  X^{n+m}\to \mathbb R^m$
be smooth maps, where  $f_2$ is a proper\footnote{This   is the  usual "proper": pullbacks of compact subsets are compact.} {\it distance  decreasing} map} and where
the "product map",
$$(f_1,f_2):X^{n+m}\to  \underline X^n\times \mathbb R^m$$
has {\it non-zero degree.}

{\it Then, if $n$ is even, there exists a point $x\in X$, where  
 $$||\wedge^2df(x)||\geq \frac {Sc(X, x)}{Sc(\underline X, f(x)).}$$

Furthermore, if $ \underline X^n=S^n$, one can allow {\sf odd $n$} and replace  the above inequality by the stronger
 one: 
 $$2||\wedge^2df(x)||_{trace}\nless Sc(X,x). \footnote  {See[Cecchini(long neck) 2020], [Cecchini-Zeidler(generalized Callias) 2021],  
[Cecchini-Zeidler(Scalar\&mean) 2021]   for 
 more general results applicable to manifolds $X^{n+m}$  with boundaries and to all closed manifolds $Y^m$, the non-existence of metrics with $Sc>0$  on which follows from non-vanishing of {\it Rosenberg index}.} $$}

\vspace{1mm}

There is a particularly useful     corollary  of this theorem, where  $X^{n+m}=Y^n\rtimes \mathbb T^m$   is   {\it a  $ \mathbb T^\rtimes$-extension} of  a manifold  $Y^n$, that is the product 
$Y^n\times \mathbb T^m$ with a warped product metric $dy^2+\phi(x)^2 dt^2$  and where the map $f:X\times \mathbb T^m$
 factors  as  $Y^n\times \mathbb T^m \to Y^n \to \underline X$  for the coordinate projection $Y^n\times \mathbb T^m \to Y^n $
 
 For instance,  such a   $\mathbb T^\rtimes $-stabilized mapping theorem  for $m=1$ together    with  the
{ \it $\mu$-bubble  separation theorem} (sections  \ref {separating3},  \ref {separating5}), yield {\it a sharp  area mapping inequality} for a class of  {\it manifolds  $X$ with boundaries}, e.g. for $X=Y\times [-1,1]$.

%%%%%%%%%%%%%%%%%%%%%
\subsubsection{\color{blue}Parametric  Area Inequalities for Families of Maps }\label{parametric3}
 
%%%%%%%%%%%%%%%%%%%%%%%%%%
{\large \color {red!30!black}  Introduce parameters wherever possible} is a motto  of modern mathematics;    
Grothendieck concept of  {\it topos} -- a category of sets parametrized by a   "topological site" --  is the
 most general manifestation of this.

The first instance of this in the present context is an application of 
 the index theorem
 to the
 
   {\sf family of flat complex line bundles $L_p$ over the torus  parametrized by the dual} ({\it Picard})  {\sf torus} 
   
   and thus showing 
 that the torus $\mathbb T^{2m}$  with an {\it arbitrary Riemannian  metric $g$} supports a {\it non-zero harmonic spinor} twisted with a {\it flat} unitary bundle; hence, {\it  no  metric $g$ on the torus may have $Sc(g)>0$} by the (untwisted) S-L-W-B formula.  
 \footnote {This idea goes back to George  Lusztig's  paper 
 {\sl Novikov's higher signature and families of elliptic s} where it is used for a  proof of the homotopy invariance of "torical" Pontryagin  classes.}

Today, this idea  is expressed in terms of  elliptic operators $\mathcal E_{\otimes A}$ with {\it coefficients in 
$C^\ast$-algebras $A$},   which, for commutative $A$,  are algebras  of continuous functions on toplogical spaces  $P$ parametrizing families  of operators $\mathcal E_p$,  $p\in P$.\vspace{1mm}

Closer  home,  we want to determine a homotopy bound on a {\it space of maps} $f:X\to \underline X$ in terms of $\inf Sc(X)$ and the 
the norms $||\wedge^2d f||$ of these maps. 

Here is an  instance of what we are looking for.\vspace {1mm}

 {\color {blue} \large  [$X\times P \to\Circle$]} \textbf {Sharp Parametric  Area Contraction Theorem}. {\sf Let $X$  be an orientable    spin manifold of  dimension $n$,
 let   $P$ be an $m$-dimensional orientable pseudomanifold, let $g_p$, $p\in P$, be a  $C^2$-continuous  family  of smooth Riemannian metrics   and let $f: X\times P\to S^{n+m}$ be a continuous map, where 
 all maps $f_p=f_{|X_p} :X=X_p=X\times \{p\} \to S^{n+m}$ are $C^1$-smooth.}
 
 {\it Then there exists a point $(x,p)\in X$, where  the $g_p$-trace norm of the exterior square of the  differential of $f_p(x)$ is bounded from below by

 $$2||trace (\wedge^2df_p)||=\sum_{i\neq j}^n \lambda_i(x)\lambda_j(x) \geq Sc(g_p,x)$$
 for some $(x,p)\in X\times P$.}
 
Consequently,  

{ \it the  inclusion $\mathcal I_{\{g\}}$ of the space $\mathcal F_{\{g\}}$ of pairs $(g, f)$, where $g$ is a Riemannian metric on $X$ and   $f: X\to S^{n+m}$   $f$  is a  smooth  map, such that 
$$2||trace (\wedge^2df_p)||=\sum_{i\neq j}^n \lambda_i(x)\lambda_j(x) < Sc(g_p,x) \mbox { for all } (x,p) \in X\times P,$$
 to the space  of all continuous maps $X\to S^{n +m}$,
$$ \mathcal  I_{\{g\}}\mathcal F_ {\{g\}}\hookrightarrow \mathcal F_{cont}(X, S^{n+m}),$$
is contractible.}
\vspace {1mm}

{\it Outlines of  two Proofs.}  1. Apply   the parametric  index theorem to the Dirac  operators on $X_p$
twisted with bundles $L_p\to X_p$ induced from  the same bundle $\underline L=\mathcal S^\pm(S^{n+m})\to S^{n+m}$ that was used in the proof of the area contraction theorems in section \ref {area extremality3} and confirm    curvature estimates needed for the
 twisted  S-L-W-B  formula.
 
 (If $n$ is odd, one has to argue as in {\color {blue} $\Circle_\blacksquare $}   in the proof of {\color {blue}  $[X_{spin} \to\Circle]$} for odd $n$ in  section \ref {area extremality3}.)

2.  Reduce the parametric problem to  the non-parametric  trace  extremality theorem
 {\color {blue} \large  [$X_{spin} \to\Circle]$} from section  \ref {area extremality3}  applied to maps $X^{n+m}\to S^{n+m} $.

To do this,  assume  $P$ is a manifold\footnote{ In the general case, by using a Thom's theorem, replace $P$ by a manifold $P'$ mapped to $P$ 
with non-zero degree} and  let $h_\lambda$, $\lambda\geq 0$, be a family of Riemannian metric on $P$ such that 
$g_\lambda\geq \lambda \cdot h_0$ and $Sc( h_\lambda)\leq  \lambda^{-1}$ and send  $\lambda\to\infty$.
Then, due to {\it additivity of trace},  application of  {\color {blue}  [$X_{spin} \to\Circle]$} yields  {\color {blue}   [$X\times P \to\Circle$].}

\vspace {1mm}

{\it Remarks.}(a) If instead of the trace norm of $df$ we had used the sup-norm, this argument would give you a non-sharp inequality, namely  with the  extra  constant $\frac {(n+m)(n+m-1)}{n(n-1)}$.

  (b) {\it Non-product families.}   Let $\{X_p\}$ be a continuous family of  compact connected orientable  Riemannian $n$-manifolds parametrized by an orientable $N$-psedomanifold $P\ni p$, that is $\{X_p\}$ is represented by a fibration
   $\mathcal X=\{X_p\}\to P$ with the fibers $X_p$.

 Let  $f: \mathcal X\to S^{n+N}$, where $n=dim(X_p)$ and  $N=dim(P)$, be continuous   map  the {\it restrictions} of which to all $X_p$ are,  smooth  {\it area non-decreasing}, e.g. 1-Lipschitz maps,  the differentials  of which are continuous in $p\in P$, and let the degree of $f$   be {\it non zero.}

 \vspace{1mm}

  {\sf If the fiberwise tangent bundle $\{T(X_p)\}$ of $\mathcal X$ is {\sf spin}, then the above mentioned  parametric index theorem to the Dirac operators on $X_p$ implies that}
  
   {\it the  infimum of the scalar curvatures of all $X_p$ satisfies 
   $$\inf_{x\in X_p,p\in P} Sc(X_p,x)\leq n(n-1).$$}

Moreover, in the extremal  case of    $\inf_{x\in X_p,p\in P} Sc(X_p,x)= n(n-1), $ 
 one can 
show  that  {\it some of $X_p$ is isometric to $S^n$.}

(If $P$ is a  smooth manifold, such that    $\mathcal X$ is  {\sf spin}, then all this can be proved with   the  index theorem
on   $\mathcal X$.)

 \vspace{1mm}
(c) {\it Maps to Fibrations.}  Let $\underline {\mathcal X}\to  P$ be a sphere bundle with the fibers 
$S^{n+N}_p=S^{n+N}$ and $f:\mathcal X \to \underline {\mathcal X}$ a fiberwise   map,
$$f=\{f_p :\mathcal X_p \to S^{n+N}_p\}.$$
Then, with a suitable defined condition "$deg(f)\neq 0$",the above inequality on the scalar curvatures of the fibers $X_p$ remains valid.

To  see this,  reduce (c) to (b) as follows.

Let $ \underline {\mathcal X^\perp}\to P$ be the complementary $S^m$-bundle, that is the 
{\it join bundle} $\underline {\mathcal X}\ast \underline {\mathcal X}^\perp$ with the fibers 
$S_p^{n+N+m+1}=  S_p^{n+N}\ast S^m$ is {\it trivial}, and observe that 
the map $f$ canonically suspends to a fiberwise  map 
$$X\ast  \underline {\mathcal X}^\perp\to \underline {\mathcal X}\ast \underline {\mathcal X}^\perp,$$
which,   due to the triviality of the fibration $ \underline {\mathcal X}\ast \underline {\mathcal X}^\perp$,  defines a map 
$$ ^\ast\hspace {-1mm}f:\mathcal X\ast \underline {\mathcal X}^\perp\to S^{n+N+m+1}.$$

Since the scalar curvatures of the fibers $\mathcal X_p\ast \underline {\mathcal X}_p^\perp$ are bounded from below by the curvature of $S^{n+N+m+1}$ (see exercise [{\Large $\ast$}] in section \ref {deceptive1}) one can  use (b), 
where, 
as in the  reduction   of the odd dimensional case of maps $X\to S^n$ to $n$ even in 
 {\color {blue} $\Circle_\blacksquare $}  in section \ref{area extremality3}, the 
 fibers  $\mathcal X_p\ast \underline {\mathcal X}_p^\perp$ and thus the space  $\mathcal X\ast \underline {\mathcal X}^\perp$ must be completed  by slightly perturbing the metric and then  extending it cylindrically at 
infinity with (arbitrarily) large scalar curvature. 
\vspace {1mm}

{\it Exercises.} (c$_1$) Use the trace norm on $\wedge^2df$ and  reduce (c) to  (b) with the    fiberwise version of  
{\color {blue} $\Circle_\blacksquare $}.

(c$_2$) Directly define "$deg(f)$" and prove (c) with  the parametric index theorem.

\vspace {1mm}

(d)  {\it Families of Non-Compact Manifolds.} The above generalizes to families of complete manifolds $X_p$ and   maps  $f: \mathcal X\to S^{n+N}$,  which are  (locally) {\it  constant at infinities} of all $X_p$ (degrees   are  well defined for such maps $f$), where, the parametric  relative index theorem,   according to 
 [Zhang(area decreasing) 2020], applies whenever all $X_p$ have (not necessarily uniformly)
 positive scalar curvatures and where     the conclusion concerns  the scalar curvatures of $X_p$ 
on the support of the differential $df$ on the  manifolds $X_p$
$$\inf_{x\in supp (df_{|X_p}),p\in P} \frac{ Sc(X_p,x)} {2  ||\wedge^2df_{|X_p}(x)||_{trace}}\leq1,  
$$

(e) {\it Foliations.} There is a  further  generalizations of  (b) to smooth foliations  $n$-dimensional leaves 
 on compact orientable  $(n+N)$-dimensional  manifolds   $\cal X$,  with  smooth Riemannian  metrics on them $X$.i
 
 Namely, 
 {\sf let    $\mathcal X\to S^{n+N}$ be a smooth map of non-zero degree.}
 
  {\it If either the manifold $\mathcal X$ is spin or the tangent bundle to the leaves is spin, then 
  there exists a point $x\in \mathcal X$,such that the scalar  curvature of the leaf $X=X_x\subset \mathcal X$ passing trough $x$ at $x$ is related to the differential of $f$ restricted to $X$ by the inequality
  $$Sc(X,x)\leq 2   ||\wedge^2df_{|X}(x)||_{trace}.$$}
   
 This is proven with $n(n-1)||df||^2$ instead of   $||\wedge^2df_{|X}(x)||_{trace}$ by Guangxiang Su [Su(foliations) 2018] and extended to complete manifolds in  [Su-Wang-Zhang(area decreasing foliations) 2021] 
 by  sharpening  the arguments by  Alain Connes  and Weiping Zhang. (The  proofs in
  these papers, if I red them correctly, allows   a use of $||\wedge^2df_{|X}(x)||_{trace} $  rather than   $||df||^2$. 
\vspace{1mm} 

{\it Examples.}  Most natural (homogeneous) foliations  with non-compact leaves support no metrics with 
$Sc>0$ by Alain Connes' theorem,  but their products with spheres $S^i$, $i\geq 2$ carry lots of such metrics, to which Su's theorem  applies.

\vspace{1mm} 

{\it Questions.}  
Does this theorem remain  valid for foliations with smooth fibres but only $C^k$-continuous in the transversal direction, such for instance,  as stable/unstable foliations of Anosov systems?

(Notice in this regard that another Connes' theorem, which  generalizes   Atiyah $L_2$-index theorem  and applies to foliations with {\it transversal measures,} 
needs  these foliations to be only  $C^3$-continuous in the transversal direction,  compare with discussion in sections  $9\frac {2}{3}$, $9\frac {3}{4}$ in  [G(positive)1996].) 

What is the comprehensive inequality  that would  include all of the above from  (b) to (e)?

\vspace {1mm}

{\it Families  with Singularities.} Is there a meaningful version of the above for families $X_p$, where some $X_p$   are singular,  as it happens, for instance, for Morse functions $\mathcal X\to \mathbb R$?
 
Notice in this regard that Morse  singularities, are, essentially, conical, where  positivity of $Sc(X_p)$
for singular $X_p$ in the sense of section \ref {max-scalar5} can be enforced by a choice of   a Riemannian metric in $\mathcal X$. \footnote{These are cones over $S^k\times S^{n-k-1}$, $n=dim X_p$, where the scalar curvature of such a cone  can be made positive, unless $k\leq 1$ and   $n-k-1\leq 1$.}

Conversely,  positivity of $Sc(X_p)$, for all $X_p$ including  the  singular ones,  probably,  yields a smooth metric with $Sc>0$ on  $\cal X$.

And it must be more difficult (and more interesting)   to decide if/when a  manifolds with $Sc>0$ admits a    Morse function, where all, including singular,  fibers have 
positive scalar curvatures or, at least, positive operators $-\Delta +\frac {1}{2}Sc$.

%%%%%%%%%%%%%%%%%%%%%%%
\subsubsection {\color {blue}Area Multi-Contracting Maps to Product Manifolds and Maps to Symplectic Manifolds}\label {multi-contracting3} 
%%%%%%%%%%%%%%%%%%%%%%%%%
A guiding principle in the scalar curvature geometry reads:\vspace{1mm}

{\color {blue!60!black}{\sf  If certain  geometric and/or topological properties of Riemannian manifolds $X_i$, $i=1,2,....,k$  imply that 
$\inf Sc(X_i)\leq \sigma_i$, then such a property of  Riemannian manifolds $X$ {\it homeomorphic  the products} $\bigtimes_iX_i=X_1\times ...\times X_k $ implies that 
$\inf Sc( X)\leq \sum_i\sigma_i$.}}\vspace{1mm}

\textbf {1. Topological non-Existence Example.}  {\sf  If  $X_1$ and $X_2$ admit no complete  metrics with $Sc>0$, and if $X_2$ is compact, then  in many, probably, not  in all cases  the product $X_1\times X_2$ admits no such metric either, (this seems to fail for 
 SYS-manifolds).

A a prominent instance of this -- here and everywhere with scalar curvature --  is  $X_2$ equal to the $N$-torus $\mathbb T^N$.}\vspace{1mm}

 \textbf {2. Length Contraction Example.}{ Let $\underline X_i$ $ i = 1, ..., k,$ be orientable (spin) length extremal Riemannian manifolds with $Sc(\underline X_i) \geq 0$,
 which means that   all
smooth maps of {\it non-zero} degrees from orientable (spin) Riemannian $n_i$-manifolds $X_i$  with $Sc(X_i)>0$ to $\underline X_i$
$$f_i: X_i\to \underline X_i,$$
satisfy
$$\inf_{x_i\in X_i}\frac{Sc(X_i,x_i) }  {Sc(\underline X_i, f(x_i)) ||df_i(x_i)||^2}\leq 1.$$
Then -- this is expected in many cases --  the Riemannian  manifold 
$\underline X=\bigtimes_i\underline X_i$  is also (spin) length extremal.
 (This is, probably,  true for all {\it known} examples of {\it spin} length extremal} manifolds   
 $\underline X_i$.)

Moreover all smooth maps from    orientable  (spin) Riemannian manifolds $X$ to the product   
$\underline X= \bigtimes_i \underline X_i$ defined by a  $k$-tuple of maps $X\to X_i,$
 $$ \Phi=(\phi_1,...,\phi_k) : X\to \bigtimes^k_i\underline X_i,$$ 
which have  non-zero degree  should satisfy  the following stronger inequality,
$$\min_{i=1,...,k}\left(      \inf_{x\in X}\frac{Sc(X,x) }  {Sc(\underline X_i, \phi_i(x)) ||d\phi_i(x)||^2}     \right) \leq 1.  $$
And  in the ideal world  one expects even more: 
$$\left(      \inf_{x\in  X}\frac{Sc(X_i,x_i) }  {Sc(\underline X, \Phi(x))\left( \sum_{i=1}^k||d\phi_i(x)||\right)^2}     \right) \leq k^2.  $$\vspace {1mm}

One also expects this   product property for {\it area} rather than {\it length}, that is with the norm of the exterior power of the differentials, $||\wedge^2d\phi_i(x)||$ instead of $||d\phi_i(x)||^2$, which is  (partly) justified by what follows.

\vspace {1mm}

\textbf {Rough  Multi-Area  non-Contraction Inequality.} {\sf Let $\underline  X$ be a compact Riemannian manifold decomposed into  product  of Riemannian manifolds of positive dimensions, 
$$\underline  X=\underline    X_1\times...\times \underline  X_i\times ... \times \underline  X_k, \mbox  { } 
dim(\underline  X_i)\geq 1, $$
 let $X$ be a compact orientable  spin  manifold of dimension $n\leq dim(\underline  X)$}
and let $X\to \underline X$ be a  smooth map defined by a $k$-tuple of maps to $\underline X_i$,  
$$f=(f_1,..., f_i,...f_k):X\to   \underline X=\underline    X_1\times...\times \underline  X_i\times ... \times \underline  X_k.$$

{\sl If the image of the fundamental homology class under $f$,
$$f_\ast[X]\in H_n(\underline  X)$$ 
is  {\sf non-torsion}, then the scalar curvature of $X$ is bounded  by  the area contraction by $f$, as follows 
$$\min_i\inf _{x\in X} \frac {Sc(X,x)}{||\wedge^2df_i(x)||}\leq \sigma,$$
where the constant $\sigma$ depends on $\underline  X$ but not on $X$}.\footnote {If
 $\underline X$ is infinite dimensional, e.g. this is the Grassmann manifold of $m$-planes in the Hilbert space, then  $\sigma$ may depend on $n=dim(X)$.}
\vspace{1mm}

{\it Proof.} Since $f_\ast[X]$ is non-torsion,
there exist cohomology classes $h_i\in H^{n_i} ( \underline X_i;\mathbb Q)$, 
$\sum_in_i=n$, such that 
the cup product $h^\ast\in H^n(\underline X)$  of their lifts to $\underline X$ doesn't vanish on 
 $f_\ast[X]_\mathbb Q)$.
 
By  multiplying   $\underline X_i$, where $k_i$ are odd, by circles and multiplying $X$ by the product of these circles, we reduce   the situation to the case, where all $k_i$  as well $n=dim(X)=\sum_in_i$ are {\it even}.

Then, by the rational  isomorphism between the  $K$-theory and ordinary cohomology,   

{\sl there exist complex vector bundles $ \underline L_i\to  \underline X_i$, such that the  Chern 
character of  the tensor  product $\underline L\to  \underline X$  of the pull-backs  of $ \underline L_i$ 
to $\underline X$  doesn't vanish on $f_\ast[X]_\mathbb Q$ either.}

It follows,  that the {\it index of the Dirac  operator}  
 on $X$ with values in the $f$-induced bundle    $L^\ast=f^\ast( \underline L)$  -- we  assume that $X$ is spin and the this  is defined -- or in some associated bundle $L^\star\to X$  {\it doesn't vanish}.
(This is  elementary algebra as in  the definition f the $K$-area.)

Endow the bundles $L_i$ with unitary connections and observe,  as we did earlier,   that the norm of the curvature of the corresponding connection in $L^\star \to X$  is bounded by a 
constant  $C$ which  depends  only  on $\underline X$  and on the norms $||\wedge^2 df_i||$, but not in any other way on $X$ and on $f$.

Therefore,  by the twisted  Schr\"odinger-Lichnerowicz-Weitzenb\"ock-Bochner   formula  the index of 
 $\mathcal D_{\otimes L^\star}$ would vanish for $Sc(X)>> C$ and the proof follows.

  \vspace {1mm}
   \textbf { Rank 1 Corollary.} {\it If $Sc(X)>0$ and (the differentials of) all maps $f_i$  
  have ranks $\leq 1$ then $f_\ast[X]_\mathbb Q=0$}.
 
 This follows from the inequality 
 $\sigma(0,0,...,0)\leq 0$ and  the definition of $\underline \sigma$. 
 \vspace {1mm}

 For instance, this shows again that   

{\it  continuous maps from orientable Riemannian   spin  manifolds $X$ with $Sc(X)>0$ to $T^m$
 send the fundamental  homology classes $[X]\in H_n(X)$ to zero in $H_n(T^m)$,} \vspace {1mm}
since   tori are products of circles  and  maps to circles have ranks $\leq 1$.

 (Maps $f$ with all their  components $f_i$ of rank one, may be  themselves  smooth embeddings
$X\to \underline X$.)
 
 \vspace {1mm}

\textbf {Sharp Multi-Area Inequalities.} {\sf Let $\underline X_i$, $i=1,...,k$, be compact orientable  Riemannian manifolds, either  
with {\it non-negative  curvature s} or Hermitian ones with positive Ricci curvatures. 
Let $X$ be a compact orientable manifold and 
let 
$$f=(f_1,...,f_k): X\to \underline X= \bigtimes_{i=1}^k\underline X_i$$
be a map a positive degree.
Let $||\wedge^2df_i||$ stands either for the {\it norm} of the second exterior power of the differential 
of the map $f_i: X\to  \underline X_i$ or , in the case where $\underline X_i$ is the sphere $S^{n_i}$, it  for the {\it averaged trace} of $\wedge^2df_i$ defined as earlier:
  $$\frac{1}{n(n-1)} ||trace (\wedge^2df_i(x))||=\frac{1}{n(n-1)} \sum_{\mu\neq \nu}^n \lambda_\mu(x)\lambda_\nu(x).$$}   
(The latter is non-greater than the former.)\vspace{1mm}

{\it \textbf {{\color {red!50!black} Conjecture}.}} {\sf
There exists  a point $x\in X$, such that 
$$Sc(X,x)\leq Sc(\underline X,f(x))\cdot  \sum_i ||\wedge^2df_i(x)||.\leqno {\color {blue} (\bigstar)}  $$}

\textbf {1.} Start with  enumerating the cases,  where this conjecture   {\it was proved} for maps from {\it spin manifolds}  $X$ to  {\it unsplit into products}    manifolds
$\underline X$, i.e. for $k=1.$

\textbf {1.A}. {\it $\underline X$ is the $n$-sphere $S^n$.}

The main computation and reduction of the case $n=2m-1$ to $n=2m$ via the map $X\times \mathbb T^1\to S^{2m}$
was performed   in  [Llarull(sharp estimates)  1998]. Then the scale invariant  trace form of  Llarull's inequality was established in [Listing(symmetric  spaces) 2010] for even $n$,  and as we explained in section \ref{area extremality3}
the trace form of the area inequality  allows an automatic reduction  $n=2m-1\leadsto n=2m$.

\textbf {1.B.} {\it $\underline X$ is a Hermitian symmetric space with $Ricci(\underline X)>0$.}
This was proved for symmetric  $\underline X$ in [Min-Oo(Hermitian) 1998] and extended to all Hermitian  spaces with
$Ricci(\underline X) \geq 0$ and  $Ricci(\underline X,x_0)> 0$ at some point in [Goette-Semmelmann(Hermitian) 2002]. 

\textbf {1.C}. {\it $\underline X$  has  non-zero Euler characteristic.}
Proved in  [Goette-Semmelmann(symmetric) 2002] and brought to the scale invariant form in [Listing(symmetric  spaces) 2010].
\vspace{1mm}

\textbf {2.} {\it Stabilization by $\mathbb T^N.$} Whenever the  inequality  ${\color {blue} (\bigstar)}$ 
is established  for  manifolds $X_o$  of dimension $n_o$ and  maps $X_o\to\underline X_o$ by confronting the index theorem with  the twisted  Schr\"odinger-Lichnerowicz-Weitzenb\"ock-Bochner   formula (there is no  known alternative  for this) then  
this argument also applies to maps $f: X\to \underline X_o\times \mathbb T^N$, $dim(X)=n_o+N$.
\vspace{1mm}

To show this, recall that  the $N$-tori  $\mathbb T^N$ for $N$ even,  support  (almost flat)  unitary bundles $\underline L_\varepsilon$  for all  $\varepsilon>0$, (and  similar families of flat  bundles a la  Lusztig)
with 

(a) {\it non-zero Chern characters}
 man
 and, at the same time  with
 
 (b)  {\it curvature operators with norms $\leq \varepsilon$.} 

Now, suppose that  ${\color {blue} (\bigstar)}$ follows with  the  Dirac   $\cal D$   on $X_o$ twisted with the  bundle $L_o\to X_o$  induced from a bundle  $\underline L_o\to \underline X_o$  by a map $ f_o:X_o\to  \underline X_o$.

Then observe that  the same argument   applies
 to $\cal D$ on 
 on $X$ twisted with the bundle $ L\to X$
  induced by   a map $f: X\to \underline X_o\times \mathbb T^N$ of non-zero degree   from  the tensor product  
  $\underline  L_o\otimes \underline L_\varepsilon\to \underline X_o\times 
  \mathbb T^N$ by letting $\varepsilon \to 0$.
 
 Indeed, (a)\&$(deg(f)\neq 0)$ imply   non-vanishing of $ index(D_{\otimes L})$,while (b) guaranties the same bound on the $L$-curvature  term in the twisted  S-L-B-W   formula  for $\varepsilon \to 0$, as in the $L_o$-curvature 
for $D_{\otimes L_o}$.\vspace{1mm}

{\it Remark.} As  we mentioned above, one can use families of {\it flat} bundles over $\mathbb T^N$, (or more generally, suitable Hilbert  moduli over the  $C^\ast$-algebra of $\pi_1(\mathbb T^N)$) which have a  advantage of giving  (slightly)  sleeker proofs of rigidity theorems.

\vspace{1mm}

\textbf {3.} The above argument, probably,  applies to general manifolds $\underline X_1$  with  bundles $\underline L_1\to \underline X_1$ instead of $ \underline L_\varepsilon \to \mathbb T^N$, where an essential point is checking that
the curvature contribution to the   S-L-B-W   formula from the induced bundle $L=f^\ast(\underline L_0\otimes \underline L_1)\to X$ for maps $f:X\to  \underline X_o\times  \underline X_1$
is bounded by the sum of the  corresponding  contributions from $f^\ast_o( \underline L_o)$  and $f^\ast_1( \underline L_1)$
for maps $f:X_o\to  \underline X_o$ and $f_1:X_1\to   \underline X_1$.

We suggest the reader  will verify this, while  we  turn ourselves  to a special case,  where the necessary linear algebraic computation has been already done.\vspace{1mm}

\vspace{1mm}

\textbf {4.} {\it Maps to Products of  2-Spheres and to Symplectic Manifolds.} {\sf Let 
$$\underline X  =\bigtimes_{i=1}^k S^2_i,$$ 
$\underline S^2_i$ are spheres  with  smooth Riemannian metrics,    let $X$ 
be a compact orientable  Riemannian manifold of dimension $2k$  and let
$$f=(f_1,...,f_k):X\to \underline X$$
be a  smooth map.

Let $\underline \omega_i$  be the  area forms of  $S^2_i$,  thus,
 $\int_{ \underline S^2_i}\underline \omega_i  =area (\underline S^2_i)$, and let $\omega_i$ be the 2-forms   on 
$X$ induced   from $\underline \omega_i$ by $f_i\to S^2_i$.

Observe that $||\wedge^2df_i(x)||=||\omega_i(x)||$ equals the maximal area dilation by $f_i$   at $x$  of surfaces $S\ni x$ in $X$.}

{\it $f$ has non-zero degree, then there  exist  a point $x\in X$, where the scalar curvature of $X$
is bounded  un terms of $||\wedge^2df_i(x)||$ as follows,
$$Sc(X,x)\leq8\pi  \sum_i \frac {||\wedge^2df_i(x)||}{area(\underline S^2_i)},\leqno {\color {blue} (\bigstar_2)} $$
where the equality holds if and only if $X$ is the product of Euclidean spheres
$X=\bigtimes_{i=1}^k S^2(r_i)$ with no restrictions on  their radii  $r_i$   and on the Riemannian metrics in $\underline S^2_i$.}\vspace{1mm}

{\it Proof.} Start by observing that the right   hand  side of ${\color {blue} (\bigstar_2)}$  doesn't depend on 
 the choice of Riemannian metrics on $\underline S^2_i$ and we  may assume all $\underline S^2_i$ isometric to the unit sphere $S^2=S^2(1)$.
 
Let $\underline L\to \underline  X=(S^2)^k$ be the tensor product of the  pullbacks of the Hopf bundle over $S^2$ under the $k$  projections  $\underline  X\to S^2$ and observe that the curvature form of this 
(complex unitary line) bundle $\underline L \to X$ is: 
$$curv(\underline L)=\frac {1}{2} \sum_i\omega_i.$$
Therefore, for all $x\in X$, the diagonal decomposition of form $\omega _x$  in   an orthonormal 
   basis in the tangent space $T_x(X)$,  orthonormal basis $(\tau_i,\theta_i)$, $i=1,...,k$, 
$$\omega =\sum_i\lambda_i \tau_i,\wedge \theta_i, \mbox { }\lambda_i\geq 0 $$
satisfies 
$$\sum_i\lambda_i\leq \sum_i  ||\wedge^2df_i(x)||.$$

It follows (theorem1.1 in [Hitchin(spinors) 1974])  that if
 $$Sc(X,x)>8\pi  \sum_i \frac {||\wedge^2df_i(x)||}{area(\underline S^2_i)}, $$
 then 
$X$ supports {\it  no non-zero  harmonic spinors} twisted with $L$.

On the other  hand the  top term in the  Chern character of $L$ is  {\it non-zero} and  the index theorem 
says that $X$ {\it does support} such a spinor, and, as everywhere   in this kind of argument, the proof follows by contradiction.
\vspace{1mm}

{\it Symplectic  Manifolds  and $\underline \omega$-Extremality.} The above argument equally applies to maps of non-zero degree  between $2k$-dimensional orientable manifolds,   $f:X \to\underline X$, where $ \underline X$ is endowed with a closed 2-form  $\underline \omega$, such that 

$\bullet$ the cohomology  class  $\underline c=\frac {1}{2\pi}[\underline \omega]\in H^2(\underline X;\mathbb R$ is {\it integral}:
$\int_S[\underline \omega]\in 2\pi \mathbb Z$ for all closed oriented surfaces in $\underline X$
(the basic example is one  half  of the area form on $S^2$);

$\bullet$ the product 
of $\exp c=1 +c+\frac {c^2}{2}+ ...+\frac {c^k}{k!}$
where $c=f^\ast(\underline c)\in H^2(X)$ for the cohomology homomorphism $f^\ast : H^2(\underline X) \to  H^2(X)$,
with the {\it Todd  class} $\hat A(X)$
 (a polynomial  in Pontryagin classes of $X$, see section  \ref{Dirac4})  doesn't vanish on the fundamental  homology class of $X$ 
 $$(\exp c)\smile \hat A[X]\neq 0.$$
 (For instance  $c^k\neq 0$ and  $\underline X$ is {\it stably parallelizable}, which, by {\it Hirsch immersion  theorem}, is equivalent to the existence of a smooth immersion  $\underline X\to \mathbb R^{2k+1}$, while $c^k\neq 0$.)

 {\it $\kappa_\star$-Invariant.}  Let   $\underline X=(\underline X,\underline \omega,\underline h)$  be a smooth   manifold,   where:
 
 $\underline \omega$
is a {\it differential   2-form} 
 on  $\underline X$, e.g. a symplectic one, i.e. where $\underline \omega$ is closed, the dimension of $\underline X$ is even and $\underline\omega^m$,
 $m=\frac {dim(\underline X)}{2}$, nowhere vanishes  on $\underline X$
 
  and where $\underline h\in H_n(\underline X)$
is a {\it distinguished  homology class}.

Define $\kappa_\star(\underline X)$n  as the infimum of the numbers $\kappa>0$, such that  all smooth  maps of from  all  closed orientable Riemannian  {\it spin}\footnote{This can be relaxed to 
properly formulate $spin^c$.} manifolds of dimension $n$  to  $\underline X$, 
$$f:X\to\underline X,$$ 
which send  the fundamental homology class of $X to \underline h$,
$$f_\ast[X]=\underline h,$$
satisfy
$$ \inf_{x \in X} Sc(X,x)\leq 4\cdot \kappa \cdot  trace(\omega(x)), $$
where $ \omega=f^\ast(\underline \omega))$ is the $f$-pullback of the form $\underline \omega$
and  
$$ trace(\omega(x))= \sum\lambda_i$$ 
for the above $g$-diagonalization of  $ \omega$. 

(See \S$5\frac {4}{5}$ in [G(positive)1996] and 
section3.4 in [Min-Oo(scalar) 2020]  for integral versions  of   this invariant.)

A Riemannian manifold  $X$ is called {\it $\underline \omega$-extremal} if it admits a smooth map
$f:X\to\underline X$, such that $f_\ast[X]=\underline h$ and 
$$Sc(X,x)=4\cdot \kappa_\star \cdot  trace(\omega(x)),\mbox { for all }  x\in X.$$

The above  proof of ${\color {blue} (\bigstar_2)} $ actually shows that the product of 
spheres $\underline X=(S^2)^k$ 
is {\it $\underline \omega$-extremal}  for the sum of the area forms $\underline \omega_i$  of   the $S^2$-factors of $\underline X$,  $$\underline \omega =\sum_i\underline \omega_i,$$ 
where  $\underline h\in H_{2k}(\underline X)$ is the {\it fundamental} class $ [\underline X]$, 
  where $\kappa=\frac {1}{2}$  and where any   {\it symplectomorphism}  $X=\underline X\to\underline X$ can be taken    for $f$.

\vspace{1mm}

{\it Remarks.} (a)  The above  is a reformulation of  a special case of {\it area  extremality}\footnote{ {\it Area  extremality} of a Riemannian manifold $X=X(g)$ (essentially)  means that all metrics 
$g'$ with $Sc(g') >Sc(g)$ on $X$ must have $area_{g'}(S)<area (S_g)$ for some surface $S\subset X$.
If $X$ is a K\"ahler manifold then $\omega$-extremality (obviously) implies area extremality for the K\"ahler form $\omega$ of $X$.}
 theorems from [Min-Oo(Hermitian) 1998],  [B\"ar-Bleecker(deformed algebraic) 1999] and  [Goette-Semmelmann(Hermitian) 1999], where  the authors  establish the $\underline \omega$-extremality of several classes of {\it K\"ahler manifolds} including compact Hermitian  symmetric spaces,  K\"ahler manifolds $X$ with $Ricci(X)>0$ and also of certain  complex  algebraic   submanifolds  
 $X\hookrightarrow \underline X=\mathbb CP^N$,   with the Fubini-Study form $\underline \omega$ on $\mathbb CP^N$.

\vspace{1mm}

(b)   Besides  multi-area  contraction  inequalities there are similar multi-length inequalities,  such as  the  multi-width   $\square^n$-inequality  from section \ref {multi-spreds7}, where the  (stronger) multi-area  contraction  inequality  doesn't apply.

\vspace {1mm}

{\it {\color {red!50!black} Conjecture}} {\sf All (most?) $\underline  \omega$-extremal  manifolds are {\it K\"ahlerian},  or closely associated with with  {\it K\"ahlerian}  or similar manifolds, such, e.g. as 
 K\"alerian$\times \mathbb T^m$.}
\vspace{1mm}

{\it Admission.}  I don't even see, why  the forms $\underline \omega $ in all  extremal cases  must be closed but not, say,  "maximally non-closed", such as generic ones.

\vspace{1mm}

{\it Question.} Are there further sharp inequalities between  (norms of) differentials $df_i$ for maps 
$$f=(f_1,...f_i,...f_k): X\to \underline X= \bigtimes _{i=1}^k S^{n_i}, \mbox { }\sum_in_i=n,$$ 
with $deg(f)\neq 0$  and  (the lower bound on) $Sc(X)$ besides
$$\inf_{x\in X}\frac {Sc(X,x)}{Sc(\underline X,f(x))\cdot  \sum_i ||\wedge^2df_i(x)||}\leq 1  $$
from the  above  {\color {blue} $(\bigstar)$}  and/or  its $||\wedge^2df_i(x)||_{trace}$ counterpart?

Namely, what are conditions on   numbers $\sigma$ and $b_1, ...b_i,...b_k$, such that there exists a compact orientable (spin) manifold $X$ of dimension $n=\sum_in_i$  with $Sc(X)\geq \sigma$ and a smooth map $f=(f_1,...f_i,...f_k): X\to \underline X= \bigtimes _{i=1}^k S^{n_i}$  with $deg(f) \neq 0$,  such that  $||\wedge^2df_i(x)||\leq b_i$ for all $x\in X$?

%%%%%%%%%%%%%%%%%%%%%%%%%%%%%

\subsection {\color {blue} Sharp  Bounds on   Length  Contractions of  Maps from    Mean Convex  Hypersurfaces} \label{mean convex3}

%%%%%%%%%%%%%%%%%%%%%%%%%%%

 The    Atiyah-Singer  theorem, when applied   to the {\it double  \DD$(X)$ of a compact manifold $X$}  with  boundary,    delivers a non-trivial geometric information on $X$ as well as  on the boundary  $Y=\partial X$.

For instance, if $mean.curv(Y)> 0$,
 then,  as we explained in section \ref {mean1},  the  natural,  continuous, metric $g$  on   \DD$(X)$ can be   approximated by $C^2$-metrics $g'$ by    smoothing $g$  along the "$Y$-edge" {\it without a  decrease of the scalar curvature} in a rather canonical manner.
Here is an instance of what comes this way. 

\vspace {1mm}

{\color {blue} [$Y_{spin} \to$ \EllipseShadow]}  { \textbf {Mean Curvature  Spin-Extremality Theorem.}}  {\sf Let 
 $X$ be a compact  Riemannian  manifold of dimension $n$ with orientable  {\it mean convex  boundary}\footnote{This mean the mean curvature of    the boundary is non-negative, where the sign convention is such that boundaries of convex domains in $\mathbb R^n$ are mean convex.}  $Y$ and   let  $\underline Y \subset \mathbb R^n$ be a smooth compact convex hypersurface. 
 
  Let $h$ and $\underline h$ denote the Riemannian metrics in $Y$ and $\underline Y$ induced from their ambient manifolds and let $h^\natural$ and $\underline h^\natural$ be their MC-normalizations  (see
  section \ref {MC-normalization+area extremality3}),
  $$h^\natural(y)= mean.curv(Y,y)^2 \cdot h(y)   \mbox {  and } \underline h^\natural(y)= mean.curv(Y,y)^2 \cdot \underline h(y).$$}

{\it Then, provided the manifold $X$ is {\sf spin},  all $\lambda$-Lipschitz maps  $f:Y\to \underline Y$ with $\lambda<1$ are contractible. }
  
In other words,  

{\it if a smooth (Lipschitz is OK) map  $f:Y\to \underline Y$ has a non-zero degree,
then there exists a point $y\in Y$, where the norm of the  differential of $f$ is bounded from below as follows:
$$ ||df(y||\geq \frac{mean.curv(Y,y)} {mean.curv( \underline Y, f(y))}.\footnote {Here we agree that $\frac {0}{0}=\infty$.}$$}

 If $\underline Y=S^{n-1}$,  this extremality, as in the case of the scalar curvature,  can can be sharpened with a use of the trace norm of the differential $df$ ...,  except that I have not verified the computation and leave  the following "theorem" with a question sign.    

\vspace {1mm}

 {\color {blue}  [$Y_{spin} \to${\small\CircleShadow}]} {\textbf{Mean Curvature Trace  Extremality Theorem(?)}\footnote {Probably, the quickest way to remove "?", at least for even $n$, is by  adapting/refining the argument from  [Lott(boundary) 2020]and/or from [B\"ar-Hanke(boundary) 2021].}.} {\sf Let $ X$ be a compact orientable Riemannian spin manifold of dimension $n$ with orientable boundary $Y$ and $f:Y\to S^{n-1}=\underline y$ be  a map  with  $ deg(f)\neq 0$.}

 Then   $f$  {\it can't be {\color {blue!50!red} trace}-wise  strictly decreasing  with respect to the $MC$-normalized  metrics  $h^\natural=mean.curv(Y)^2h$  for the Riemannian metric $h$ on  $Y$ induced from $X\supset Y$  and 
$\underline h^\natural= (n-1)^2ds^2$ on $S^{n-1}$,} that is 
there is a point $x\in X$, where {\it  the trace-norm of the  differential of $f$ is bounded from below by the mean curvature of $X$ as follows:
$$\frac {1}{(n-1)}\sum_{i=1}^{n-1}  \lambda_i^\natural(y) \geq 1\mbox {
for $\lambda_i^\natural=\lambda_i(f^\ast (\underline h^{\natural}))/h^{\natural})$,}$$}
which means that the trace-norm of $df$ with respect to the original (non-normalized) metrics satisfies:
$$ \frac {1}{(n-1)}  ||df(y)||_{trace}\geq \frac {mean.curv(\underline Y, f(y))}{mean.curv(Y, y)}. $$

The simplest and most interesting  common corollary of these  two  theorems  is the following.
\vspace {1mm}

 {\color {blue} \SunshineOpenCircled} ({\it Seemingly Elementary}) \textbf {Example.} {\sf If the mean curvature of a   smooth 
  hypersurface $Y\subset \mathbb R^n$ is {\it bounded from below by $n-1$}, that is the mean curvature of the unit sphere $S^{n-1}\subset \mathbb R^n$, then {\it all $\lambda$-Lipschitz map $f: Y\to \mathbb R^n$, where  $\lambda <1 $, are contractible.}}\footnote{It is impossible not to ask oneself what happens for $\lambda=1$,  i.e. where $f$ is  distance non-increasing.  You bet, such an $f$ is  either contractible, or it is an {\it isometry}. Indeed (almost)  all our  extremality theorems  are accompanied by rigidity results  in the equality  cases, as we shall see later on. 
  
  But it  is non-trivial to formulate and hard to solve    the {\it stability problem}:  {\sf what happens to geometries of  hypersurfaces  $X_\varepsilon\subset \mathbb R^n$ with $mean.curv(X_\varepsilon )\geq n-1$  and to  
  $(1+\varepsilon)$-Lipschitz maps to $S^{n-1}$ with  non-zero degrees, when $\varepsilon\to 0$.}}

(If "Lipschitz" is understood with respect to the {\it Euclidean distance function} on $X$,   rather than the larger one which is    associated  with the {\it induced Riemannian metric}, the proof easily follows from {\it Kirszbaum theorem}.)
\vspace {1mm}

 {\it About the Proof of {\color {blue} \EllipseShadow}}.   Let    $X$ lie in a (slightly larger) Riemannian $n$-manifold $X_+\supset X$ without boundary, let  $Y_\varepsilon^{n+l-1}\subset X_+\times \mathbb R^l$   be the boundary of the $\varepsilon$-neighbourhood of $X\subset X_+\times \mathbb R^l$  and let us 
 similarly, define  $\underline Y_\varepsilon^{n+l-1}\subset \mathbb R^{n+l}=  \mathbb R^n\times \mathbb R^l$ as the boundary of  the $\varepsilon$-neighbourhood of $\underline X\subset \mathbb R^{n+l}=  $ for $\underline X\subset \mathbb R^n$ with  boundary $Y$.
 
 Observe -- this needs a little computation as in  section 1.4  --  that the lower bounds on the  scalar curvatures of  the "interesting parts"  
 $$\mbox {$Y_{\varepsilon\varepsilon}^{n+l-1}\subset Y_{\varepsilon}^{n+l-1} $ and 
 $\underline Y_{\varepsilon\varepsilon}^{n+l-1}\subset \underline Y_\varepsilon^{n+l-1}$}$$
  which are $\varepsilon$-close to the original  $Y\subset Y_\varepsilon^{n+l-1}$ and  $\underline Y\subset   \underline Y_\varepsilon^{n+l-1}$,  are perfectly controlled by their mean curvatures, while their complements, being flat in the ambient manifolds,  have the same scalar  curvatures  as $X$ and  $\underline X$, where the latter is equal to zero.

 Then extend    $f: Y\to  \underline Y$  a map 
 $$f_\varepsilon :Y_\varepsilon^{n+l-1} \to \underline Y_\varepsilon^{n+l-1},$$
 such that  the "interesting part"of $Y_\varepsilon^{n+l-1}$ goes  to that of    $Y_\varepsilon^{n+l-1}$  and the complement of one to the complement of the other and such that 
 the "interesting part"  of this extensions is done in a most economical  manner along  normal geodesics to $Y\subset Y_{\varepsilon\varepsilon}^{n+l-1}$ and to $  \underline Y  \subset  \underline Y_{\varepsilon\varepsilon}^{n+l-1}.$
  
  If we  do it with a proper care then, for a small enough $\varepsilon$ and  $l$ with the same parity as $n$,  we shall be able to apply the
   spin-area convex  extremality theorem  {\color {blue}  [$X_{spin}{^\to}$\Ellipse]} from the section \ref 
   {area extremality3} to the map $f_\varepsilon$, which that would need a preliminary smoothing  of
   the manifolds  $Y_\varepsilon^{n+l-1}$ and $\underline Y_\varepsilon^{n+l-1}$  by tiny $C^1$-perturbations (these manifolds themselves are only $C^1$-smooth), where, while   
while smoothing the hypersurface   $\underline Y_\varepsilon^{n+l-1}$  convex,  smoothing of  $\underline Y_\varepsilon^{n+l-1}$  must {\it keep the flat part flat.} 

Because  of the latter, the point $y_\varepsilon\in  \underline Y_\varepsilon^{n+l-1}$, where
$$Sc(\underline Y_\varepsilon^{n+l-1}, f(y\varepsilon)\cdot  ||\wedge^2df(\varepsilon)||\geq Sc( Y_\varepsilon^{n+l-1},y\varepsilon), $$
provided by   {\color {blue}[$X_{spin}{^\to}$\Ellipse]} must be necessary located in the "interesting region" $Y_{\varepsilon\varepsilon}^{n+l-1}$;   then the needed inequality for the mean curvature of $Y$ 
  will be satisfied by the point  $y\in Y$ nearest  to $y_\varepsilon.$ %(We shall return(?) to this in section????) 

{\it Remark about}  {\color {blue}  [$Y_{spin} \to${\small\CircleShadow}]}. To carry out the above argument one needs a  generalization of of the spherical trace inequality  {\color {blue}   [$X_{spin} \to\Circle]$}
from the previous section to manifold $\underline X$ that don't have full  $O(n+1)$-symmetry  of $S^n$.

In the present case the relevant metric $\underline g$ is $O(n)$ invariant   and one needs a separate 
bounds  on the two parts of the  trace norm of   $\wedge^2 df$:

 {\sf the first part     comes  from $\frac {n-1(n-2)}{2}$   bivectors $e_i\wedge e_j$ with $e_i$ and $e_j$, $i,j=1,...,n-1$,  tangent the $S^{n-1}$-spherical  $O(n)$-orbits and the second one  from the $n-1$ remaining  $e_i\wedge e_n$ with the vector  $e_n$ normal  to these orbits.  } 
  
This is an instance of   a more general principle:

{\sf \color {blue!50!black} to achieve the sharpest inequality, one should choose  the norm for  measuring  $df$ in accordance with the the symmetries of  the manifold $\underline X$.}

 We shall see later on other instances  of this 
"principle", e.g. for maps to  products of spheres in section   \ref {multi-contracting3}. \vspace{1mm}
\vspace {1mm} 

{\it \textbf {On Non-spin Manifolds and on $\sigma<0$.}}  {\sf \color {red!50!black} Conjecturally}, 
{\sf if the boundary $Y=\partial X$ of a  compact orientable Riemannin $n$-manifold $X$ with 
$Sc\geq -n(n-1)$ admits a smooth map $f$ with {non-zero} degree to the boundary of the $R$-ball in
 the hyperbolic $n$-space with sectional curvature $-1$,
$$f:Y\to \partial B(R)\subset \mathbf H^n(-1),$$
and if 
$$ mean.curv (Y)\geq n-1 \mbox  {  and } ||df||\leq 1,$$ 
then the map $f$ is an   {\it  isometry}. Moreover,  $f$ extends  to an isometry $X\to B(R)$.}\footnote
{Granted $f$ {\it is}  an isometry (with respect to the induced Riemannin metrics in $\partial X\subset X $ 
and $\partial B(R)\subset B(R)$), an    isometry $X\to B(R)$ follows from {\it Min-Oo's hyperbolic rigidity theorem} from section \ref{negative3}.}

\vspace{1mm}

We shall   prove a partial result   in this direction  with a use of 
{\it stable  capillary $\mu$-bubbles}, which may also apply to maps to more general 
hypersurfaces in $\mathbf H^n(-1)$  (see  section  \ref {capillary warped5}),  but it remains unclear how 
to approach the trace-norm version of this conjecture.\vspace{1mm}

{\it \color {blue} Questions and  Exercises}. (a) {\sf Is there  an elementary   proof of this inequality  for $n\geq 4$?}\vspace{1mm}\vspace{1mm}

(b) Besides the lower bound on the  mean curvature, that is the sum of the principal curvatures, $\sum_i \alpha_i$, the "size" of a hypersurface $Y$ is 
bounded by the scalar curvature $\sum_{i\neq j} \alpha_i\alpha_j$ and also - this is obvious 
by the product of the principal curvatures $\prod_i \alpha_i$.\vspace{1mm}

\hspace {-4mm}{\it Are there similar inequalities for  other elementary  symmetric functions of $\alpha_i$.} \vspace{1mm}

(If $Y\subset \mathbb R^n$ is  {\it convex}, i.e. all $\alpha_i\geq 0$,   then $\prod_i \alpha_i$ minorizes the rest of  elementary symmetric functions, which gives a   
 trivial proof of  {\color {blue} \SunshineOpenCircled} and similar inequalities for other  symmetric functions for distance decreasing maps from convex hypersurfaces to $S^n$.)
 
  the above theorems   for {\it convex}  hypersurfaces  $Y\mathbb R^n$.)

But it is unclear if, for instance, there is a bound on this radius in terms of 
$\sum_{i>j>k} \alpha_i\alpha_j\alpha_k$ for $n\geq 5$    when  this sum is positive.)
\vspace{1mm}

(d) Let $Y_0\subset \mathbb  R^n$ be a smooth  compact cooriented submanifold with boundary $Z=\partial Y_0$, such that 

{\sf  the mean curvature of $Y_0$ with respect to its coorientation satisfies  
$$mean.curv(Y)\geq n-1=mean.curv (S^{n-1}).$$}

Show that 

{\sl every distance decreasing map
$$f: Z\to  S^{n-2}\subset \mathbb R^{n-1}$$
is contractible,}

 where "distance decreasing"refers to the  distance functions on $Z\subset \mathbb R^n$ and on  $S^{n-2}\subset \mathbb R^{n-1}$  coming from the ambient Euclidean spaces $ \mathbb R^n$ and $\mathbb R^{n-1}$.\vspace{1mm}

{\it Hint.} Observe that the maximum   of the  principal curvatures of $Y_0$ is $\geq 1$ and show that 
{\it the filling radius} of $Z\subset \mathbb R^n$  is $\leq  1$.\footnote{ This means that $Z$ is homologous to zero in
 its 1-neighbourhood.}   \vspace{1mm}

(e) {\it \color {blue} Question.} {\sf Does  contractibility of $f$      remains valid if the distance decreasing property   of $f$ is defined  with the (intrinsic) spherical distance in  $S^{n-2}$ and with the   distance  in $Z\subset Y_0$    associated  
  with the {\it intrinsic metric}  in $Y_0\supset Z$, where $dist_{Y_0}(y_1,y_2)$ is defined as  the infimum of   length of curves in $Y_0$ between $y_1$ and $y_2$?}

\textbf{(f)} Formulate and prove the mean curvature counterparts of the theorems            {\color {blue}   [$X_{spin}\overset{\hat A} \to\Circle]]$},  {\color {blue}$\bigtimes\mathbb R^m$}  and {\color {blue}  [$X\times P \to\Circle$]}  for maps $X^{n+m}\to\underline X^n$ 
and $X^{n}\to\underline X^{n+m}$  from  sections  \ref{area extremality3}   and \ref {mean convex3}, either by the above $Y_{\varepsilon\varepsilon}^{n+l-1}$-construction or by generalizing Lott's argument for manifolds with boundaries.

\textbf{(h)} {\it Question.}  Is there a version (or versions) of the mean curvature extremality theorems for   maps to products of convex hypersurfaces in the spirit of area multi-contracting maps in section  \ref {multi-contracting3}

%%%%%%%%%%%%%%%%%

 \subsection {\color {blue} Riemannian Bands  with $Sc>0$ and {\color {blue} \large$\frac  {2\pi}{\mathbf  n}$}-\textbf {Inequality.} }\label {bands3}

%%%%%%%%%%%%%%%%%%%%%%%

We saw in the previous sections how a use of twisted Dirac operators leads to geometric bounds, including certain  sharp ones, on the size of compact Riemannian spin manifolds. Such  bounds usually (always) extend to non-compact complete manifolds, but until  recently    no such result was available for non-complete manifolds and/or for manifolds with boundaries. 
\footnote{Several such results  have appeared in the papers 
[Zeidler(bands) 2019], 
[Zeidler(width) 2020]  and  [Cecchini(long neck) 2020],[Cecchini-Zeidler(generalized Callias) 2021],  
[Cecchini-Zeidler(Scalar\&mean) 2021],  [Guo-Xie-Yu(quantitative K-theory) 2020],  which we briefly discuss letter on .}

On the other hand  minimal hypersurfaces were used in [GL(complete)1983] for obtaining rough bounds for non-complete manifolds;    below, we shall see how  such hypersurfaces (and $\mu$-bubbles in general) serve for getting  {\it sharp} geometric inequalities of this kind. \vspace {1mm}

{\it {\color {blue!80!red}  Bands}}, sometime we call them {\it \color {blue} capacitors}, are  manifolds $X$  with two distinguished disjoint non-empty  subsets in the boundary $\partial(X)$, denoted 
 $$\partial _-= \partial_-X\subset \partial X \mbox { and } \partial _+= \partial_+X\subset \partial X.$$

A band is called {\it \color {blue!80!red}  proper}  if $\partial_\pm$ are unions of connected components of $\partial  X$  and 
$$\partial _-\cup \partial _+=\partial X.$$

The basic instance of such a band is the segment  $[-1,1] $, where $\pm \partial =\{\pm1\}$. 

Furthermore, {\it cylinders }
$X=X_0\times  [-1,1] $ are also bands with $\pm \partial =X_0\times\{\pm1\}$, where such a band is proper if $X_0$ has no boundary.\vspace{1mm}

  {\it \color {red!20 !blue} Riemannian bands} are  those  endowed with  Riemannian metrics and \vspace {1mm}
  
  {\it the width of a Riemannian band $X=(X,\partial_\pm)$} is defined as 
  $$width (X)=dist (\partial_ -,\partial_+),$$
  where this distance is understood as the infimum of length of curves in $V$  between $\partial_-$  and $\partial_+$. \vspace {1mm}

We are mainly  concerned at this point with    {\sl compact}  Riemannian bands $X$  of dimension $n$, such that \vspace {1mm}

 {\color {blue}\SquareShadowBottomRight\textSFliii$_{Sc \ngtr 0}$} {\sf {\large {\it \color {red!50!black}no }}  closed  embedded  hypersurface
$Y\subset X$, which {\it separates} $\partial_ -$ from $\partial_+$, admits a   $\mathbb 
T^1$-stabilization $Y^\rtimes$ with positive scalar curvature, i.e. no complete (warped product)  metric on the product  $Y\times \mathbb T^1$ of the form $dy^2+\phi(y)^2dt^2$ has   $Sc(h^\rtimes)>0$.}

\vspace {1mm}

(Since   $Y$ is  compact,   the existence of this (warped product) metic $h^\rtimes$ with   $Sc(h^\rtimes)>0$ is equivalent to the existence of a metric $h$ 
with $Sc(h)>0$ on $Y$ itself, 
since  the  conformal Laplacian   $-\Delta+ {area extremality3n-2}{4(n-1)} Sc$ is more positive that the  
 $-\Delta+ \frac{1}{2}Sc $  implied by   positivity of $Sc(h^\rtimes)$.)
\vspace{1mm}

{\it Representative Examples} of  compact  bands  with this property are:

{$\bullet_{\mathbb T^{n-1}}$  {\sl toric bands} which are homeomorphic to  $ X=\mathbb T^{n-1}\times  
[-1,1]$;\vspace {1mm}

{$\bullet_{SYS}$ manifolds $X$  homeomorphic to a Schoen-Yau-Schick manifolds times $[-1,1]$;\vspace {1mm}

$\bullet_{\hat \alpha}$ these, called  {\it $\hat \alpha$ bands, }  are diffeomorphic to $Y\times  
[-1,1]$, where the $Y$ is a closed spin $(n-1)$-manifold with non-vanishing $\hat \alpha$-invariant (see  \ref  {spin index3} the  {\color {blue}IV} above);\vspace {1mm}

$\bullet_{\mathbb T^{n-1}\times \hat \alpha}$ these are bands {\it diffeomorphic to  products} $X_{n-k}\times \mathbb T^k$, where  $\hat \alpha (X_{n-k})\neq 0$.}

\vspace {1mm}

(A  characteristic  non-compact example with a similar property  is

$\bullet_{\times \mathbb R\setminus\{\mathbb Z\}}$: $X$ is homeomorphic to the product $\mathbb T^{n-2}\times\mathbb R\times [-1,1]$ 
minus a discrete subset.)\footnote { The property {\small \color {blue}\SquareShadowBottomRight\textSFliii$_{Sc \ngtr 0}$}
for  {\sl toric bands} and for SYS-bands  follows from the Schoen-Yau codimension 1 descent theorem (see section \ref {SY+symplectic2}), in  the case 
$\bullet_{\hat \alpha}$ this is the  Lichnerowicz-Hitchin theorem (section  \ref{spin index3})  and $\bullet_{\mathbb T^{n-1}\times \hat \alpha}$ is a corollary to  theorem 2.1 in [GL(spin) 1980], while a "complete"  version of this  property for  the non-compact 
$(\mathbb T^{n-2}\times\mathbb R\times [-1,1])\setminus \{\mathbb Z\}$ is  an example, where theorem 6.12. from
[GL(complete) 1983] applies. (See sections \ref{obstructions4},  \ref {obstructions5} for more about these   and more general examples.}

\vspace {2mm}

{\color {blue} \large$\frac  {2\pi}{\mathbf  n}$}-\textbf {Inequality.} {\sf Let  $X$ be a  {\sl proper compact Riemannian bands $X$  of dimension $n$ with $Sc(X)\geq \sigma>0$.}}

{\it  If no closed hypersurface in $X$
 which separates $\partial_ -$ from $\partial_+$ admits a metric with positive scalar curvature, then
$$width (X)=dist(\partial_-,\partial_+)\leq 2\pi \sqrt \frac {(n-1)}{n\sigma} = \frac { 2\pi}{ n}\cdot \sqrt \frac {n(n-1)}{\sigma}
\leqno{  \mbox {
{\Large $[$}{$\varocircle$}{\small$ _\pm \leq\frac { 2\pi}{ n}$}{\Large$]$}}} $$

In particular if  $Sc(X)\geq Sc(S^n)= n(n-1)$,
then $$width (X)\leq  2\pi\sqrt \frac {(n-1)}{n\sigma} = \frac { 2\pi}{ n}.$$
Moreover, the equality holds in this case  only  for warped products $X=Y\times \left ( - \frac {\pi}{n}, \frac {\pi}{n}\right)$\footnote{Here,  since   $X$ is non-compact, the width  is  understood as the distance between the two ends of $X$.}    with metrics $\varphi^2h+dt^2$, where the metric  $h$ on $Y$ has $Sc(h)=0$ and where 
  $$\varphi(t) =\exp \int_{-\pi/n}^t -\tan \frac {nt}{2} dt, \mbox { }  -\frac {\pi}{n}<t < \frac {\pi}{n},$$}
   as in  section  \ref {warped2}.
 \vspace {1mm}

 {\it About the  Proof.} If  a hypersurface $Y\subset X$, which separates $\partial_-$ from  $\partial_+$ contains a descending chain  (flag) of  closed oriented  hypersurfaces},     
 $$Y\supset Y_{-1}\supset ...\supset Y_{-i} \supset...,$$
where 
 where each $Y_{-i}\subset X$ is equal to a transversal intersection of $Y_{-(i-1)}$ with  a smooth closed oriented sub-band     $H_i\subset X$, of codimension one,
 $$H_i\cap X_{-(i-1)}=X_{-i}$$
and where 
 $ Y_{-i}$ represent {\it non-zero} classes in  the homology $H_{n-1-i}(X)$,  then one can  proceed  by the inductive Schoen-Yau's kind of  descent method (see sections \ref{SY+symplectic2})     with  minimal hypersurfaces 
  $$...X_{-i}\subset  X_{-(i-1)}\subset ... \subset X_{-1}\subset X,$$ 
 where these  $X_{-i}$ are  $\mathbb T^\rtimes$-symmetrised as in 
the  $[ \rtimes_\varphi]^N$-symmetrization theorem in section \ref{warped stabilization and Sc-normalization2}
where  $X_{-i}$ in our band $X$ 
  have  "free"  (pairs of)  boundaries  contained  in $\partial_\mp(X_{-(i-1)})$, and such that  the intersections $ X_{-i}\cup Y$ are homologous to  $ Y_{-i}$.

This argument  delivers the sharp version of  {\color {blue} $\frac  {2\pi}{\mathbf  n}$} for  {\it over-toric bands}, i.e. those which admit 
maps $X\to\mathbb T^{n-1}$, $n=dim(X)$, with non-zero degrees of their restriction to $\partial_\mp$, but when it comes to SYS-bands,  one gets only a weaker lower bound on  $width(X)$, that is  by $\frac  {4\pi}{n}$, instead of   $\frac  {2\pi}{  n}$.  

The same weakening of   {\color {blue}$\frac  {2\pi}{\mathbf  n}$} takes place if 
separating hypersurfaces  $Y\subset X$, are {\it enlargeable},  e.g. if the interior of   $X$, assumed compact,  admits a  complete metric with non-positive sectional curvature.  And if separating $Y$  are {\sl \color {blue} SYS times enlargeable}, one has to be content with $\frac  {8\pi}{  n}$.\footnote {This is  worked out in \S2-6 of  [G(inequalities) 2018].}

\vspace {1mm}

In  section \ref{separating3},  we present   a  more efficient   argument,  where, instead of working with chains of  minimal hypersurfaces, we     show  in one step that if $width(X)\geq 2\pi \sqrt \frac {(n-1)}{n\sigma}$, then a  certain   {\it stable $\mu$-bubble $Y_{st}\subset X$, which separates $Y_-$ from $Y_+$,
supports a metric with $Sc>0$.}

Besides, the sharp  {\color {blue}$\frac  {2\pi}{\mathbf  n}$} for  wide class  of {\it spin} bands was recently proven by Zeidler, 
 Cecchini and   Guo-Xie-Yu with new index/vanishing  theorems on {\sf \color {blue!50!black}Dirac operators with potentials on manifolds with boundaries.}
\footnote {See 
[Zeidler(bands)  2019], [Zeidler(width)  2020],
 [Cecchini(long neck) 2020] and the most recent  [Guo-Xie-Yu(quantitative K-theory) 2020],[Cecchini-Zeidler(generalized Callias) 2021],  
[Cecchini-Zeidler(Scalar\&mean) 2021].}

   \vspace {1mm} 
 
 {\it Remarks.}(a)
If  hypersurfaces separating $\partial_-$  from $\partial_+$ in $X$ are enlargeable, 
e.g. if $X$ is homeomorphic to $\mathbb T^{n-1}\times [0,1]$, then 
 a non-sharp version of {\color {blue} $\frac  {2\pi}{ n}$}-inequality,
$$dist(\partial_-,\partial_+)\leq 2^n\pi \sqrt \frac {(n-1)}{n\sigma}  $$
follows from theorem 12.1 in[GL (complete)1983].

  (b) One might  think  that the sharp $\frac  {2\pi}{n}$-inequality,   must be obvious for domains in the unit sphere $S^n$
 homeomorphic to $\mathbb T^{n-1}\times [-1,1]$ and   for   bands 
with {\it constant sectional curvatures} in general;   to my surprise, I couldn't find a  direct  proof of it even for $X$ is homeomorphic to $\mathbb T^{n-1}\times [0,1]$.
 
%%%%%%%%%%%%%%%%%%

 \subsubsection {\color {blue} Quadratic Decay of Scalar Curvature on Complete Manifolds with $Sc>0$.}
\label{quadratic3}

%%%%%%%%%%%%%%%%%%%%%%
{\it\textbf  {QD-Exercise}.} {\sl   \color {blue!60!black}  Quadratic   Decay Property.} {\sf Let $X$ be a complete non-compact  Riemannian $n$-manifold  and $X_0\subset X$ a compact subset, such that {\it there is no  domain}   $X_1\subset X$, which contains $X_0$ and   the boundary $\partial X_1$ of which (assumed smooth)  {\it  admits a  metric} with $Sc>0$, e.g.
$X$ is homeomorphic to $\mathbb T^{n-2}\times\mathbb R^2$.}
 
 Show that  there exists a constant  $R_0=R_0(X, x_0)$, such that 
 
 {\it the minima of the scalar curvature  of $X$ on  concentric  balls  $B(R)=B_{x_0}(R)\subset X$  
 around a point $x_0\in X$, satisfy
$$  \min_{x\in B(R)} Sc(X,x)\leq \frac{4\pi^2}  {(R-R_0)^2}  \mbox  {  for all } R\geq R_0. $$}
 
{\it Hint}. Apply $\frac  {2\pi}{n}$-inequality  to the annuli  between the spheres or radii $R$ and 
$R$ for a suitable constant $c$.

(Compare this with the quadratic decay theorem  in section 1 of 
[G(inequalities) 2018] and   see [Wang-Xie-Yu(decay) 2021]  for estimates of the scalar curvature decay rates   by contractibility radius  and the diameter control of the asymptotic
dimension and observe that, if $X$ is homeomorphic to $\mathbb T^{n-2}\times\mathbb R^2$,  than the quadratic decay with the constant  $2^{n+1}\pi^2$ follows from [GL(complete 1983].)\vspace{1mm}

{\it\textbf  {Critical Rate of Decay  {\color {red!50!black} Conjecture}}.}
There exists  a  universal  {\it critical constant $c_n$,}  {\sf conceivably}, $c_n=\frac {4\pi^2(n-1)}{n}$,  such that:

 [a]   {\sf if a smooth manifold  $X$ admits a complete  metric $g_0$ with $Sc(g_0)>0$, then, for all $c<c_n$,  it admits a complete metric $g_\varepsilon$, with $Sc(g_\varepsilon)>0$ and  at most    $c$-{\it sub-quadratic} scalar curvature decay,
$$Sc(g_\varepsilon, x)\geq \frac {c}{dist(x,x_0)^2}  \mbox { for a fixed $x_0\in X$ and  all $x\in X$ with $ dist(x, x_0) \geq 1$};$$  
 and  
  
 [b]    if $X$ admits a complete  metric $g_0$ with $Sc(g_0)>0$ and  $c$-sub-quadratic for $c>c_n$
 scalar curvature decay,
 $$Sc(g_\varepsilon, x)\geq  \frac {c_n}{dist(x,x_0)^2}  \mbox { for }  dist(x, x_0) \geq 1,$$
 then it admits a complete metric with $Sc\geq \sigma>0$.}
 
Moreover, 

{\sf for all continuous functions $\omega = \omega(d)$, there exists a complete metric $g_\omega$  on $X$, such that 
$$Sc(g_\omega, x)\geq \omega(dist(x, x_0))\mbox { for a fixed point $x_0$ and all $x\in X$}.$$ }

Here is a related {\it \textbf {compactness conjecture}}, which expresses the following  idea:\vspace {1mm}

{\sl \color {blue!40!black} The existence of a complete metric with $Sc\geq \sigma>0$ on an $X$ is detectible by topologies of   compacts parts $V$ of $X$}: \vspace {1mm}

{\sf if, for all compact subsets  $V\subset X$ and all  constants $\rho>0$, there exists a  (non-complete)  metric on $X$ with $Sc\geq 1$,  such that the closed  $\rho$-neighbourhood  $U_\rho(V)\subset X$ is {\it compact}, then 
$X$ {\it admits a complete Riemannian metric with $Sc\geq 1$.}}

%%%%%%%%%%%%%%
 \subsection {\color {blue} Separating Hypersurfaces  and the Second  Proof of the $\frac {2\pi}{n}$-Inequality}
\label{separating3}
  %%%%%%%%%%%%%%%%%%%%%%%%

The main ingredient in the proof of  the general $\frac {2\pi}{n}$-Inequality is the following. 
\vspace {1mm}
 
 {\color {blue} {\large \sf I{\large \sf \color {red!80!black}\sf I}\sf I}} \textbf { $\mu$-Bubble Separation Theorem.} \vspace {1mm} {\sf Let $X$ be an $n$-dimensional,  Riemannian band, possibly non-compact and non-complete.

   Let  $$Sc(X,x)\geq \sigma(x)+\sigma_1, \mbox {  },$$
   for a continuous function $\sigma=\sigma(x)\geq 0$ on $X$ and a constant $\sigma_1> 0,$
   where $\sigma_1$  is related to  $d=width(X)=dist_X(\partial_-,\partial_+)$ by the inequality
 $$\sigma_1 d^2>   \frac {4(n-1)\pi^2}{n}.$$    } 
(If scaled  to $\sigma_1=n(n-1) $, this   becomes $d> \frac {2\pi}{n}$.)

{\it Then   there exists a smooth  hypersurface $Y  \subset X$,   which separates
  $\partial_-$ from  $\partial_+ $, and a smooth positive  function $\phi$ on   $Y$, 
such that the scalar curvature of the  metric $g_\phi=g^\rtimes_\phi = g_{Y_{-1}}+\phi^2dt^2$ on $Y\times \mathbb R$ is bounded from below by 
$$Sc(g_\phi,x)\geq \sigma(x).$$}

{\it Derivation of  $\frac {2\pi}{n}$-Inequality from  {\color {blue} {\large \sf I{\large \sf \color {red!80!black}\sf I}\sf I}}.}
 If a band $X$ with $Sc\geq \sigma>0$ has  $width(X)=dist(\partial_-,\partial_+)> 2\pi\sqrt \frac {(n-1)}{n\sigma}$, then  {\color {blue} {\large \sf I{\large \sf \color {red!80!black}\sf I}\sf I}} implies the existence 
of a separating hypersurface $Y$  and a function $\phi(y)$, such that  $Sc(g^\rtimes_\phi)\geq \varepsilon $  for a {\it small} $\varepsilon >0$.

\vspace {1mm}

{\it About the Proof of {\color {blue} {\large \sf I{\large \sf \color {red!80!black}\sf I}\sf I}}}.  If $X$ is compact and $n\leq 7$, we take a {\it $\mu$-bubble  $Y_{min}$} for $Y$,   that  is the minimum of the functional
$$Y\mapsto vol_{n-1}(Y)-\mu[Y, \partial_-]$$
defined in the space of separating hypersurfaces $Y\subset X$,  where
$[Y, \partial_-]\subset X$ denotes the region in $X$ between $Y$ and $\partial _-\subset \partial X$  and where the key point is to choose $\mu$  suitable for this purpose.

What  is required of $\mu$ is that 

$\bullet$  the boundaries $\partial_\pm$ must  serve as  barriers for our variational problem and thus ensure  the existence of $Y_{min}$; 

$\bullet$  positivity of the second variation should imply the positivity of the  $\Delta +Sc(Y_{min}) -\sigma$ on $Y$.

 This is achieved with  $\mu$, that   is modeled after 
the  measure $\underline \mu$  on $\mathbb T^{n-1}\times [-1,1]$, (the density of) which is equal  the mean curvatures of 
 the hypersurfaces 
$\mathbb T^{n-1}\times \{t\} $  with respect to  the warped  product metric  $\varphi^2h+dt^2$ for 
$$\varphi(t) =\exp \int_{-\pi/n}^t -\tan \frac {nt}{2} dt, \mbox
 { }  -\frac {\pi}{n}<t < \frac {\pi}{n}.\footnote {This is the same $\varphi(t)$ that was used in section 12 in  [GL(complete) 1983] for proving a rough lower bound on  the norms of the  differentials of smooth maps of non-zero degrees  from {\it non-complete} Riemannian   manifolds $X$ with $Sc(X)\geq 1$ to $S^n$  for $n=dim(X)\leq 7$.}$$
\vspace{1mm}

 {\color {blue} {\large \sf I{\large \sf \color {red!80!black}\sf I}\sf I}$_\circlearrowleft$}
   \textbf    {Separation with Symmetry.}  {\it If the Riemannian band is isometrically acted upon by a compact group $G$,
then the separating hypersurface $Y\subset X$ and the function $\phi$ on $Y$ can be chosen invariant under this action.}

\vspace{1mm}

{\it Proof.}   Use  the  multi-dimensional   Morse lemma (see section \ref{eigenfunctions2});  alternatively,   
 apply   more elementary  uniqueness/symmetry property of the lowest eigenfunction of the (linear elliptic) second order variational  (linear  elliptic)   $\Delta_Y+s$ on a hypersurface $Y$,  which minimizes the functional $vol_{n-1}(Y)-\mu[Y, \partial_-]$ among $G$-invariant separating hypersurfaces  $Y\subset X$.

\vspace {1mm}

{\it Remark.}  In our case, the group $G$ is the torus $\mathbb T^k$, which  freely acts on $X$, and the equivariant $\mu$-bubble problem (trivially) reduces to the ordinary one on the quotient space $X/\mathbb T^k$. 

To make use of this  for the next step of $\mathbb T^\rtimes$-symmetrization,  one only  needs to check -- this is an exercise to the reader  --  that the corresponding warped product  with $\mathbb T^{k+1}$ will have the same scalar curvature as one gets by doing this in $X$ itself.

\vspace {1mm}

{\it Compact/Non-compact.} If $X$ is non-compact, then,  as usual, we exhaust $X$ by compact submanifolds with boundaries, 
proceed as in the compact case (these compact bands an not proper, part of their boundary is not contained in 
$\partial X=\partial _-\cup \partial_+$, but this causes no problem) and then pass to the limit. This is routine.

\vspace {1mm}

{\it \textbf {Example of  Corollary.}} {\sf Let $X=(X,g)$ be an  $n$-dimensional manifold with uniformly positive  scalar curvature, $Sc(X)\geq \sigma>0$, and let 
$f:X\to \underline X= \mathbb R^{n-m}$   be a smooth proper  (infinity to infinity) 1-Lipschitz (i.e.  distance non-increasing) map. }

{\sl Then the homology class of  the pullback of the generic point, $f^{-1}(\underline x)\subset X$,
is representable by a compact  submanifold $Y\subset X$, such that the product $Y\times \mathbb T^m$ 
admits  a $\mathbb T^m$-invariant  { \sf   (warped product)} metric  $h^\rtimes$  ($h=g_{|Y}$) with $Sc(h^\rtimes)>0$.

Consequently, $Y$ itself  admits a metric with $Sc>0$.}\vspace{2mm}

{\it Singularity Problem for $dim(X)>7$ and the Second  Proof of the  $\mu$-Bubble Separation Theorem.}  By the standard theorems of the geometric measure theory,  the minimizing $\mu$-bubble $Y\subset X$ exists for all $n$ but it may have singularities of codimension $7$.\footnote {The most general existence theorem of this type applicable to all codimensions 
is in the technically difficult  Almgren's  1986 paper "{\sl Optimal Isoperimetric Inequalities}".

The  existence and regularity  theorem we need in codimension  one are  easier, they follow  by 
   the usual technique of integer  currents  and regularity theorems,  see  [Ros(isoperimetric) 2001]), [Morgan(isoperimetric)(2003); the arguments  from this papers  which are  applied there  to the more  traditional formulation of the isoperimetric 
    problem, can be carried over  to our $\mu$-bubble  setting  with no  problem;  alternatively, one can use the language of
Caccippoli sets.}

 (The first instance of this is 
the vertex of the famous cone  from the origin over $S^3\times S^3\subset S^7\subset \mathbb R^8$.) 

If  $n=8$, then   (a minor generalization of) Natan Smale's
{\it generic regularity theorem}  takes care of things, but if $n\geq 9$ one needs to adapt  {\it Lohkamp's minimal smoothing} 
results and/or  techniques to our case. My, rather superficial, understanding of Lohkamp works  suggests that this is possible, but it  can't  be safely applied unless   everything is   written out in full detail.

  I feel more comfortable  at this point with  generalizing  
  theorem 4.6 from  the    Schoen-Yau paper [SY(singularities) 2017],  where it  used in    the
    inductive descent method with {\it singular minimal}  hypersurfaces, to our minimizing $\mu$-bubbles. 
    
 Such a generalization feels plausible and, if it's true, this must be obvious to Schoen and Yau. (I guess, the same can be said
about  what Lohkamp thinks about generalization  of his theorem to $\mu$-bubbles.)

Granted this,   one gets 
the  sharp {\color {blue}$\frac  {2\pi}{\mathbf  n}$}-inequality for SYS-bands, and in fact  for  all bands $X$
which satisfy {\color {blue}\SquareShadowBottomRight}, where the poof of non-existence of metrics of positive scalar curvatures on  separating hypersurfaces $Y\subset X$ is obtained by exclusively by inductive decent  with no 
appeal to Dirac operators and related invariants, such as the $\hat A$-genus and the $\hat \alpha$-invariant.\footnote{
This is complementary  to what can obtained by Dirac operators methods of Zeidler and Cecchini.}

\vspace {1mm}

{\it Minimal Hypersurfaces in Non-compact Bands.} An essential advantage of $\mu$-bubbles over minimal hypersurfaces is that the former are easier to "trap" them and  prevent from fully  sliding away  to infinity than the former. 

For instance if $X$ is a complete non-flat  manifold with positive sectional curvature which is conical at infinity
then it contains no  complete (even locally)  volume minimizing hypersurfaces, but it contains lots of stable complete (and  compact) $\mu$-bubbles. 

However, a version of the   {\color {blue}$\frac  {2\pi}{\mathbf  n}$} can be proven for non-compact complete bands by reduction to "large" compact  {\it non-proper} bands $X$,  where the boundary is divided into three parts
$$\partial X=\partial_+\cup \partial_-\cup\partial_{side}$$
where $\partial_+=\partial_+(X)$ and $ \partial_-=\partial_-(X)$ are disjoint with controlled lower bound on the 
distance between them, while  $\partial_{side}=\partial_{side}(X)$, which intersects both $\partial_+$ 
and $\partial_-$ is supposed to be far away from the bulk of the intended minimal 
hypersurfaces in $X$.\vspace {1mm}

 {\it Example.}  Let $\underline X$ be the cylinder $B^{n-1}(R)\times [-1,+1]$, where $B^{n-1}(R)$ the Euclidean $R$-ball of dimension $n-1\geq 2$  and where 
 $\partial_{side}(\underline X)=S^{n-3}(R)^{n-1}(R)\times [-1,+1]$  for  the equatorial sphere 
 $S^{n-3}(r)\subset S^{n-2}(r) =\partial B^{n-1}(r).$

Let $\underline \partial_r\subset \underline X$, $r>R$  be the $r$-cylinder concentric to $\partial_{side}(X)$, that is 
 $$\underline \partial_{side(r)}=S^{n-2}(R)^{n-1}(R)\times [-1,+1].$$

Minimal hypersurfaces   $Y=Y_r\subset X$    we shall meet in $X$ will be similar to those   $\underline Y\subset  \underline X$, which have  their boundaries contained in 
$\underline Z_r=\partial_+( \underline X)\cup \partial_-( \underline X)\cup \underline \partial_{side(r)}$  and 
which represent non-zero homology classes in  $H_{n-1}( \underline  X,  \underline  Z_r.$

\vspace {1mm}

Namely,  let $X$ be a  compact  orientable non-proper  $n$-dimensional band.
   Let 
$f =(f_1,f_2): X\to \underline X= B^{n-1}(R)\times [-1, +1$ be a smooth map which sends 
$\partial_\mp(X)\to  \partial_\mp(\underline X)$  and  $\partial_{side} (X) \to \partial_{side} (\underline X)$ and 
such that 

$\bullet_1$ the map $f_1:X\to B^{n-1}$ is   $1$-Lipschitz;

$\bullet_2$ the map $f_2:X\to[-1,+1]$  is $\lambda$-Lipshitz for $\lambda>0$

$\bullet_3$ the map $f$ has {\it non-zero} degree.

Observe  
 that 
 
 $\bullet_1$ implies that  
$$dist(\partial_{side} (X),\partial_ {side(r)}(X)\geq R-r;$$

 $\bullet_2$  makes   
 $$width (X)=dist(\partial_-(X),\partial_+(X))\geq d=d_\lambda= \frac {2}{\lambda};$$

 $\bullet_3$  shows that if an oriented   hypersurface $Y\subset X$ with $\partial Y\subset Z_r$, 
 represents a {\it non-zero} homology  class in $H_{n-1}(X,Z_r)$, then it  necessarily intersects  $\partial _{side(r)}(X)$.

In fact, $Y$ intersects every $(n-3)$-dimensional submanifold $Z'\subset Z_r$ (observe that $dim(Z_r)=n-2$ for generic maps $f$) which separates $\partial_-(Z_r)= Z_r\cap \partial_-(X)$ 
from $\partial_+(Z_r)= Z_r\cap \partial_+(X)$.

 %%%%%%%%%%%%%%%%
\subsubsection {\color {blue}Paradox with Singularities} \label{singularities3}

%%%%%%%%%%%%%%%%%%%%%%%

 Singularities must enhance the power of 
minimal hypersurfaces and stable $\mu$-bubbles rather  than to reduce it, since the large  curvatures  of hypersurfaces $Y\subset X$ (these curvatures  are infinite at singularities) {\it add} to the positivity of the second variation  .

Thus, for instance, if a $Y\subset X$, where $dim(X)=8$  and $Sc(X)\geq -1$,  is a stable minimal hypersurface with a singularity at $y_0\in Y$  and if no smooth  submanifold  in the homology class of $Y$ admits a metric with $Sc>0$, e.g. $X$ is homeomorphic to the torus and $Y$ is non-homologous to zero, then 
scalar curvature of $X$ can't be non-negative outside a small neighbourhood of $y_0\in X$.

Yet, there is no known   argument for $dim(X)\geq 9$  fully implementing this idea. 

\vspace {1mm}

{\it  On  $n=dim(X)=8$.}   If $dim(X)=8$ then  stable minimal hypersurfaces and $\mu$-bubbles  $Y\subset X$ have isolated singularities which can be removed by small generic perturbation  
  as in [Smale(generic regularity) 2003]  as  follows.    \vspace {1mm}

{\it Theorem.} {\sf Let $Y_0\subset X$ be a cooriented  compact  isolated volume minimizing hypersurface 
   and let et $X_t =[X_0,Y_t]\subset X$  be  the bands between $Y_0$  and hypersurfaces $Y_t$,   which are  positioned    close to $Y_0$ on their "right sides"  in $X$,  and which  minimize the function 
 $Y\mapsto vol_{n-1}(Y)-t\cdot vol[Y_0, Y]$ 
for $0\leq t\leq \delta$ for a small $ \delta>0$.}

{\it If $n=8$, then submanifolds $Y_t$ are non-singular for an  open dense set of $t\in [0,\delta]$.}\vspace {1mm}

{\it Outline of the Proof.} The key (standard) facts one needs here are as  follows.

1. {\it Monotonicity.} If the sectional curvature of  $X$ is bounded by $\bar \kappa^2$, then 
the volume of intersections  of $m$-dimensional minimal subvarieties $Y\subset X$ with $r$-balls  
$B_{y_0}(r)\subset X$ centered at $y_0\in Y$  satisfy
$$ \frac {d r^{-m} vol_m (Y\cap B_x(r)}{dr} \leq  const_n r \bar \kappa.\mbox { for  all }  r\leq r_0=r_0(X, y_0)>0.$$

2. {\it Corollary.} The densities of (singularities of) minimal $Y\subset X$  are {\it semicontinuous:}

{\sf if  be a sequence of pointed manifolds with uniformly bounded geometries,  $(X_i, x_i),$  Haussdorf converges to $(X,x)$ and if  minimal subvarieties $Y_i\subset X_i$,   which contain the points $x_i$,  current-converge to $Y\subset X$, then 
$$\limsup dens(Y_i, x_i)\leq dens(X,x),$$  
where, recall,
$$dens(Y,x)= \lim_{r\to 0}  r^{-m} vol_m (Y \cap B_x(r)), mbox { } m=dim(Y).$$}

3. {\it Weak Compactness}: {\sf The set $\mathcal Y_A$ of minimal subvaraities $Y\subset X$ with volumes bound by a constant $A$ is 
compact in the current topology  for all  $ A<\infty$.}

4. {\it Codimension one Intersection Property.}  {\sf Minimal codimension one cones  
$C_1,C_2\subset \mathbb R^n$  {\it necessarily  intersect} by  the {\it maximum principle}.}

5.  {\it Split Cone Property.} {\sf Let $C\subset \mathbb   R^n$  be a minimal  cone. Then  either the density of this cone at the apex $0\in C$ is maximal $+\varepsilon$, 
$$dens  (C,0)\geq dens(C,c)+\varepsilon\mbox {  for  all $0\neq c\in C$ and some }\varepsilon=\varepsilon(C)>0,$$ 
or the cone split, i.e. $C=C_{ -1}\times \mathbb R^1$  for a minimal cone $C_{-1}\subset \mathbb R^{n-1}.$}\vspace {1mm}

Now, turning  to the proof,   let all $Y_t$  have singularities at some points $y_t\to y_0\in Y_0$, $t\to 0$,  and assume without loss of generality, this is possible  due to 2, that  
$$dens(Y_t, y_t)= dens (Y_0,y_0).\footnote { Our $Y_t$ are $\mu$-bubble rather  than minimal, but 
 this makes no difference at this point.}$$

Let  $\lambda$-scale these 
$Y_t$ at  $y_0$, thus making   $\lambda Y_0$, converge to a minimal cone, call it  $Y'_0\subset T_{y_0}(X)\mathbb R^n$,  and let $Y'_t$ be what remains of the limits of other $Y_t$.

Since these  $Y'_t$ don't intersect $Y'_0$, none of $Y'_t$ is conical, which is only possible if the singularities of $Y_t$ slide tangentially along $Y_0$     for $t\sim \lambda^{-1}$ by the distance $c(t)$, such that 
$c(t)/\lambda\to \infty $ for $\lambda\to\infty$. It follows that if all $Y_t$ were singular, these singularities would accumulate in the limit  to a ( one dimensional or larger) singularity of $Y'_0$ of constant density equal to that of $dens(Y_0,y_0)$. Therefore, the cone $Y'_0$  splits, since 
$n=8$, it is non-singular 
and the proof follows by contradiction.

{\it  On  $n=dim(X)\geq 8$. } (a) Schoen-Yau in their desingularization argument apply descent  by warped  $\mathbb T^\rtimes $-symmetrised/stabilized  minimal hypersurfaces  
$$X=X^n\supset Y^{n-1} \supset ... \supset  Y^{n-i}\supset....\supset Y^2,$$
where minimization and $ \mathbb T^\rtimes$-stabilization (essentially)  apply to  non-singular parts of these $Y$  and where 
 the main difficulty,  as far as I can see, is to show that $Y^{n-i}$  can't be eventually  sucked  in the singularity of $Y^{n-1}$,\footnote{If $n\leq 9$, this problem for overtorical $X$ can be handled with Dirac  operators, as in   section 5.3 in [G(billiards) 2014].} and where the outcome  of  this process - the surface $Y^2$  --  is non-singular.\footnote{
 Schoen and Yau   articulate  their main results  (theorems 4.5 and 4.6 in    [SY(singularities) 2017])  for compact  SYS-manifolds,  although the basic arguments  of their paper    are essentially local   and and  apply to  a wider class of manifolds.}

    (b) The main desingularization  result by Lohkamp in 
[Lohkamp(smoothing) 2018], is
    
     {\it approximation theorem of  volume  minimizing codimension one  cones $C^{n-1}\subset \mathbb R^n$ by smooth minimal hypersurfaces} (generalizing Smale's  result in the case of 
     cones)
   with the following

 {\it Splitting Corollary.}  {sf Let $X$ be a compact orientable  Riemannian manifold with $Sc(X)>0$.    Then all 
 homology classes in $H_{n-1}(X)$   are representable by hypersurfaces $Y\subset X$, which support metrics with $Sc>0$.}\vspace{1mm}

{\it Remarks} (a) As far as the topology of compact manifolds with $Sc>0$   this result is more general 
   than that by Schoen and Yau.
   
   For instance it implies that 
   
   {\it products of Hitchin's spheres  and connected sums of tori with non-spin manifolds admit no metrics with $Sc>0$}.   
   
    Nor alternative proof of this kind of results is available.

(b) As far as I understand,\footnote {My understanding of the results by Lohkamp as well as those by Schoen and Yau is limited, since I haven't mastered the proofs from  [SY(singularities) 2017]) and  from [Lohkamp(smoothing) 2018].} Lohkamp's smoothing allows applications of  our  $\mu$-bubble  arguments  to manifolds of  all dimensions $n$, with  possible exceptions for {\it rigidity  theorems } 
for  {\it non-compact} manifolds.

(c) The above 1-5  seems to suffice for smoothing  conical singularities (am I missing hidden subtleties?) but it is unclear to me how  Lohkamp's splitting corollary for 
$n\geq 9$  follows from it.

 %%%%%%%%%%%%%%%%%%%%%%
\subsubsection {\color {blue}  $\mathbb T^\rtimes$-Stabilized  Scalar Curvature and Geometry of  Submanifolds of Codimensions One, Two and Three} \label {codimensions 1, 2, 3.3} 
%%%%%%%%%%%%%%%%%%

Besides   {\color {blue} \large$\frac  {2\pi}{\mathbf  n}$}, there are other immediate applications of the separation theorem {\color {blue} {\large \sf I{\large \sf \color {red!80!black}\sf I}\sf I}}. \vspace {1mm}

 { \color {blue}\textbf {[1]}} \textbf {Compact Exhaustion Corollary.}  {\sf Let $X$ be a complete Riemannian manifold with 
 $Sc(X)\geq \sigma >0$.}
 
 {\it Then $X$ can be exhausted by compact domains $U_i$ with smooth boundaries $Y_i=\partial U_i$
$$U_1\subset U_2\subset... \subset  U_i\subset  ...\subset X, \mbox {  }  \bigcup_iU_i=X,$$
such that $U_{i+1}$  is  contained in the $\rho$-neighbourhood of $U_i$ for all $i=1,2,...$ and 
and where all $Y_i$ admit $\mathbb T^{\rtimes }$-extension $Y_i\rtimes \mathbb T^1$ 
with 
$$Sc(Y_i\rtimes \mathbb T^1)\geq \frac {\sigma}{2}.$$}

{\it Poof.}  Let $S(10), S(20),... S(10i),...\subset X$  be concentric spheres 
 around a point $x_0\in X$, let  $Y_i$ be hypersurfaces in the annuli  $[S (10i), S(10(i+1))] $ between these spheres, which     separate  $S (10i)$ from  $S(10(i+1))$  and which  enjoy  the properties  supplied by  {\color {blue} {\large \sf I{\large \sf \color {red!80!black}\sf I}\sf I}}. Then take the domains in $X$ bounded by $Y_i$ for $U_i$. 

\vspace {1mm}

 { \color {blue}\textbf {[2]}} \textbf {Codimension 2  Corollary}.  {\sf Let $X$ be a (possibly non-compact)   connected   orientable  $n$-dimensional Riemannian  manifold with boundary,  let $\underline X$ be a   compact  connected  orientable surface  with boundary and with an arbitrary metric compatible with topology  and let $\Psi:X\to \underline X$ be a smooth {\it distance decreasing} map which sends the boundary $ \partial X$ to $\partial \underline X$.}
 
{\it  If $Sc(X)\geq \sigma+\sigma _1$,  $\sigma,\sigma _1>0$ and the inradius of $\underline X$  is bounded from below by 
$$inrad (\underline X)=\sup_{\underline x\in \underline X} dist(\underline x, \partial  \underline X)>  \frac {2\pi}{\sqrt \sigma}, $$
 then $X$ contains an oriented  codimension two (possibly disconnected) submanifold $Y\subset X$,} {\sf  which, if $X$  is non-compact, is  {\it properly} embedded  to $X$} {\it   and which is homologous} 
 {\sf for the    homology group $H^{ncpt}_{n-2}(X) $}  with {\it infinite} {\sf supports in the  case of non-compact $X$)}
  {\it  to the pullback  $\Psi^{-1}(\underline x)\subset X$  of a generic point  $\underline x\in \underline X$, and such that $Y$ with the induced  Riemannian metric from  $X$ admits a $\mathbb T^2$-extension,  that is the product $Y\times \mathbb T^2$ with the metric $g_\phi=dy^2 +\phi^2 (dt_1^2+dt_2^2)$,  such that 
 $$Sc(g_\phi)\geq \sigma_1.$$}

 \vspace {1mm}
 
  {\it Proof. }   Let $X_1\subset X$ be the   {\large \sf \color {red!80!black}\sf I}-hypersurface that, according to  {\color {blue} {\large \sf I{\large \sf \color {red!80!black}\sf I}\sf I}}, separates the boundary of $X$   from the 
   $f$-pullback  of the (small disc around)  the point $\underline x\in  \underline X$ furthest from   the boundary
    (as in the proof of $\mathbb T^\rtimes $-stabilized Bonnet-Myers diameter inequality   [{\color {blue} BMD}]  in section \ref{warped stabilization and Sc-normalization2} and  apply
      {\color {blue}$\frac  {2\pi}{\mathbf  n}$} to the infinite cyclic covering of $X_1\rtimes \mathbb T^1$ induced by the natural  cyclic covering of $ \underline  X$ minus this point.

 { \color {blue}\textbf {[2$'$]}} \textbf {Codimension 2  Sub-Corollary}. 
{\sf Let $X$ be a closed  orientable  $n$-dimensional Riemannian  manifold with $Sc(X)\geq \sigma>0$,
 let $\underline X$ be a   closed surface with an arbitrary metric compatible with topology and 
 let $\Psi:X\to \underline X$ be a smooth {\it distance decreasing} map.}
 
 {\it If {\sf no} closed oriented  codimension two  submanifold $Y\subset X$ {\sf homologous to the pullback}  $\Psi^{-1}(\underline x)\subset X$  of a generic point  $\underline x\in \underline X$
  admits a metric with $Sc>0$, then 
  the diameter of the surface  $\underline X$ is bounded in terms of $\sigma$ as follows.
$$ diam( \underline X)< \frac {2\pi}{\sqrt \sigma}.$$}

{\it Proof. } Let
$\underline  x_0, \underline x_1 \in \underline  X$ be   mutually furthest  points  and apply the above to the pullback $X_-$ of the  complement $\underline  X_-$ to
a small disc  in $\underline  X$  around $\underline  x_0$.   

 \vspace{1mm}

 { \color {blue}\textbf {[3]}} \textbf {Area non-Contraction Corollary}. {\sf Let $X$ be a proper compact orientable   Riemannian band of dimension $n+1$, let  
 $\underline X\subset \mathbb R^{n+1}$ be a smooth convex hypersurface and let $f:X\to \underline X$ be a smooth map the restriction of which to $\partial _-\subset \partial X$ (hence, to $\partial_+$ as well) has {\it non-zero} degree.}
 
{\it  If $X$ is spin and if $n$ is even,\footnote {As we have  said already several times, these conditions must be redundant.} then  there exists a point $x\in X$, where the  exterior square of the differential  of $f$  is bounded from below in terms of 
    $d=width(X)=dist(\partial_-,\partial_+)$  and  the scalar curvature $Sc(X,x)$ as follows.  
$$Sc(\underline X, f(x))\cdot  ||\wedge^2df(x)||\geq Sc(X,x)- \frac {4(n-1)\pi^2}{nd^2}.$$}

Furthermore, {\it if $\underline X= S^n$, then, now for odd as well as for even $n$,
the trace norm of $\wedge^2df$ satisfies:
$$2 ||\wedge^2df(x)||_{trace}\geq Sc(X,x)- \frac {4(n-1)\pi^2}{nd^2}.$$}

{\it Proof.} Apply the  { \sl $\mathbb T^m$-stabilized area/mapping  extremality theorem} 
(\ref  {area extremality3},  \ref{multi-contracting3})   for $m=1$   to $Y\rtimes \mathbb T^1$
where $Y\subset X$ is the separating hypersurface from {\color {blue} {\large \sf I{\large \sf \color {red!80!black}\sf I}\sf I}}.

\vspace {1mm}

\vspace {1mm}

 \vspace{1mm} 
 {\it Exercises}.  (\textbf a)  {\sf \color {blue!60!black} {\sl Codimension 3 Linking Inequality.}} {\sf Let $X$ be a closed  orientable  $n$-dimensional Riemannian  manifold with $Sc(X)\geq \sigma>0$,
 let $\underline X$ be the $3$-sphere
with an arbitrary metric compatible with topology and 
 let $f:X\to \underline X$ be a smooth {\it distance decreasing} map.}
  Show that

 {\it  if {\sf no} closed oriented  codimension three  submanifold $Y\subset X$ {\sf homologous to the pullback}  $f^{-1}(\underline x)\subset X$  of a generic point  $\underline x\in \underline X$
  admits a metric with $Sc>0$, then 
  the distances between all pairs of embedded circles  $S_1, S_2\subset \underline X$ with 
  {\sf non zero linking numbers} between them satisfy: 
  $$ dist( S_1, S_2)< \frac {2\pi}{\sqrt \sigma}.$$}

{\it Hint}. Use  the argument from  the   proof of the  {codimension 2  corollary} { \color {blue}\textbf {[2]}} 
and consult   [Richard(2-systoles) 2020]\footnote {Our codimension 2 area bounds, including this  exercise, are motivated by  Richard's  bound on  systoles  
of 4-manifolds with $Sc>\sigma$ proved in this paper.}  
 \vspace{1mm}
 
 (\textbf b)  {\sf \color {blue!60!black} {\sl Area non-Contraction in Codimension  3.}} {\sf Let $X$, $\underline X$ and $f:X\to \underline X$ be as in (\textbf a),  let  $\underline X_1\subset \mathbb R^{n-2}$ be a smooth closed convex hypersurface  and let $f_1: X\to  \underline X_1$ be  a smooth map, such that the "product" of the two 
 maps,
 $$(f,f_1): X\to \underline X\times \underline X_1,$$
 has {\it non-zero degree}.}
 Show that

 {\it  if $X$  is spin and $n$ is odd (thus, $dim(\underline X_1)$ even) 
  then  there exists a point $x\in X$, where the  exterior square of the differential  of $f$  is bounded from below in terms of 
    $d=width(X)=dist(\partial_-,\partial_+)$  and  the scalar curvature $Sc(X,x)$ as follows.  
$$Sc(\underline X, f(x))\cdot  ||\wedge^2df(x)||\geq Sc(X,x)- \frac {4(n-1)\pi^2}{nd^2},$$
 for $d$ equal the supremum of the distances between pairs of linked circles in $\underline X$.}

%%%%%%%%%%%%%%%%%
\subsubsection {\color {blue}  On  Curvatures of Submanifolds in the unit Ball $B^N\subset \mathbb R^N$}
 \label{focal3}
%%%%%%%%%%%%%
{\color {red} (The earlier versions of this  section contained  errors.)}  \vspace{1mm}

Here is our \vspace {1mm}

{\it \textbf {Problem.}} Given a closed  smooth $n$-manifold $X$ and a number $N>n$,

 \vspace{1mm}

\hspace  {-0mm}{\sf \color {blue!50!black}    {\sf evaluate   the minimum of the  curvatures of     smooth immersion  of $X$ to the 

unit $N$-ball,
$$f:X\hookrightarrow B^N=B^N(1) \subset \mathbb R^N.$$}} 
 \vspace{1mm}

We shall briefly describe in this section   what  is known and and what is unknown about this problem and refer to section 3 and 7 in  [G(inequalities) 2018] and to [G(growth of curvature) 2021] for  more general discussion and for the proofs.

   \vspace {2mm}
   
   \hspace {5mm} {\sc Six Examples of Immersed and Embedded  Manifolds 
   
   \hspace {32mm}  with Small Curvatures}   \vspace {1mm}

Just to clear the  terminology, we agree that a  smooth map $f:X\to Y$ is an {\it immersion}  if  the differential $df: T(X)\to T(Y)$ is {\it  injective}  on all tangent spaces 
$T_x(X)\subset T(X)$.

An immersion $f$ of a compact manifold is an {\it embedding} if it has no double
 points, $f(x)\neq f(y) $ for $x\neq y$. 

If $Y$ is a {\it Riemannian} manifold, e.g.  $Y=\mathbb R^N$, 
then  
the curvature of this $f$, denoted 
$$curv_f(X)=curv_f(X\hookrightarrow Y)=curv(X\hookrightarrow Y)=curv(X),$$
is   the    {\it supremum of the  $Y$-curvatures of all  geodesics} in $X$, where  "geodesic"  is understood with respect to the Riemannian metric in  $X$ induced from $Y$.

\vspace{1mm}

1. {\it Clifford    Embeddings.} Here, $X=X^n$ is  the product  of $m$   spheres of  dimensions $n_i$,  $\sum _{i=1}^mn_i=n$,
 all  of    the  radius $r=\frac{1}{ \sqrt m}$, 
 $$X= S^{n_1}\left (\frac{1}{\sqrt m}\right)\times ...\times S^{n_i}\left (\frac{1}{\sqrt m}\right)\times ... \times S^{n_m}\left (\frac{1}{\sqrt m}\right)$$
and $$f_{Cl}:X\hookrightarrow S^{N-1}\subset B^N(1)\subset \mathbb R^N,\mbox { } N=m+\sum_in_i,$$
is  the obvious embedding, that is the $\frac {1}{\sqrt m}$-scaled Cartesian  product of the imbeddings
$S^{n_i}(1)\subset \mathbb R^{n_i+1}.$

Clearly, 
$$curv_{f_{Cl}}(X\hookrightarrow B^N)=\sqrt m$$ 
and  the curvature of $X$ in the unit  sphere is
$$curv_{f_{Cl}}(X\hookrightarrow S^{N-1})=\sqrt {m-1}.$$

Two natural questions arise:\vspace {1mm}

{\sf  Can the products of spheres be immersed to the unit ball 
with smaller curvatures?

Are there non-spherical, immersed or embedded,  submanifolds $X\hookrightarrow  B^N(1)$ with 
 $curv(X )<\sqrt 2$? }

\vspace {1mm}

 A definite answer is available only  for  immersions $X^n\to S^{n+1}$ by a   {\it theorem of Jian Ge}. 
\footnote {See [Ge(linking) 2021] and {\color {blue}$\largelozenge$} in this section.}
   \vspace {1mm}

 $[\Circle \times \Circle]$   {\sl Clifford's   are the  only codimension one   immersed non-spherical  submanifolds $X$ in the spheres    
with curvatures $curv(X\hookrightarrow S^{n+1})\leq 1$.}\vspace {1mm}

%{\sl If an  immersion $f$  from  a closed  $n$-manifold $X$ to  $S^{n+1}$, 
%$n\geq 2$ 
%satisfies 
%$curv_f (X) \leq 1$,
%then either $f$ is a Clifford embedding or it is a diffeomorphism of $X$ to a convex hypersurface in the sphere $S^{n+1}$.}\vspace {1mm}

But if $m\geq 3$  then there are immersions of {\it non-spherical} $n$-manifolds   to 
$S^{n+m-1}$  with smaller curvature.

2.  {\it Veronese embeddings} of projective spaces to spheres,
$$f_{Ver}: \mathbb RP^n\to S^{\frac {(n+1)(n+2)}{2}-2}\subset  B^{\frac {(n+1)(n+2)}{2}-1}=B^{\frac {(n+1)(n+2)}{2}-1}(1)$$
satisfy 
$$curv_{f_{Ver}}\left ( \mathbb RP^n \hookrightarrow S^{\frac {(n+1)(n+2)}{2}-2}\right)=\sqrt\frac  {n-1}{n+1}<1.\leqno 
{\color {blue}\left [\sqrt\frac  {n-1}{n+1}\right]}$$
and 
$$curv_{f_{Ver}}\left ( \mathbb RP^n \hookrightarrow  B^{\frac {(n+1)(n+2)}{2}-1}\right) =
\sqrt{\frac  {n-1}{n+1}+1}<\sqrt 2.$$

{\color {red!50!black}Conjecturally}, these have the  minimal curvatures among 
{\it all non-spherical} $n$-submanifolds   in the unit spheres and unit balls, where the minimum for all $n$ is achieved (only conjecturally)  by  
 Veronese's  projective plane in unit 4-sphere, where 
$$curv_{f_{Ver}}(\mathbb RP^2\hookrightarrow S^4)=\frac {1}{\sqrt 3}=0.577350... \mbox {  }\leqno 
{\color {blue}\left [\frac {1}{\sqrt 3}\right]},$$
and 
$$curv_{f_Ver_2}(\mathbb RP^2\hookrightarrow B^5)=\frac {2}{\sqrt 3}=1.15470... \mbox { }.$$

3. The $\frac {1}{\sqrt l}$-scaled  Cartesian power of the Veronese  map
$$F=\frac {1}{\sqrt l}\cdot f_{Ver}^{\times l}:  X^{2l}=\underset {l}{\underbrace { \mathbb RP^2\times ... \times  \mathbb RP^2}}\to S^{4l-1}\subset B^{4l}(1)$$
competes with the Clifford embedding, for
$$curv_F(X^{2l} \hookrightarrow B^{4l})=\sqrt l \cdot \sqrt { \frac{l}{3}+1} 
<\sqrt {2l}.$$

$4.$ {\sf If $N  \geq (1+\Delta)^n$, say for $\Delta>10$ then {\it all $n$-manifolds $X$}
admits immersions  
$$f:X\hookrightarrow S^N$$ with
$$curv_f(X) \leq C_\Delta,$$
 where  $C_\Delta<\sqrt 2$ for all $n$ and where
 $$C_\Delta\to \sqrt {2-\frac {6}{n+2}}   \mbox { for } \Delta\to \infty$$ }
 with the rate of convergence which,   {\color{red!50!black} a priori,} may depend on $n$. 
 
It  is   {\sf {\color{red!50!black} unclear} if the "true" $C_\infty$ is, actually,   {\it smaller} than   $\sqrt {2-\frac {6}{n+2}}$ and   it is also 
  {\sf {\color{red!50!black} unclear} what   happens to $C_\Delta$ for  $\Delta$ close to zero.}

\vspace {1mm}

5. It easily follows from the above that 

 {\sl  if the dimension $n_m$ of the last factor in a product of spheres 
$$X^n=\bigtimes_{i=1}^mS^{n_i}, \mbox { } \sum_{i=1}^m n_i =n,$$ 
is much greater then the remaining ones, say, roughly,  
$$n_m\geq \exp \exp \sum_{i-1}^{m-1}n_i,$$
then $X^n$ admits an immersion 
$$f: X^n\hookrightarrow B^{n+1}(1)$$
such that
$$curv_f(X^n) <2\sqrt 3.$$}}
 This is {\it smaller} than Clifford's $\sqrt m$ starting from $m=12$.\vspace {2mm}

It is  {\sf {\color{red!50!black} unclear}, however, if  these $X^n$ admit {\it embeddings} to the unit  ball   with $curv(X^n\hookrightarrow B^{n+1})\leq 100,$  for example.}\vspace {1mm}

6. There are {\it no topological bounds on curvatures} of {\it immersed} submanifolds 
 of  a {\it given  dimension} $n$: 

 {\sf if an $X^n$ admits a smooth immersion to $\mathbb R^N$, then it also admits an immersion 
to the unit ball with $curv(X^n\hookrightarrow B^N) < const_n$.}

But  all we can say about this constant is,  roughly, that
 $$ 0.1n<const _n< 10n^\frac{3}{2}.$$

{\it  Imbeddings}, at least these with codimension one,  are different  from immersions in this 
regard.

 For instance, if $X=X^n$ is  disconnected and
 contains $m$ {\it mutually non-diffeomorphic} components, then, clearly, 
 $$curv_f(X\hookrightarrow B^{n+1})\geq const_n m,\mbox  { } const_n\geq \frac {1}{(10n)^n}, $$
 for all embeddings $f:X \hookrightarrow B^{n+1}(1)$.
 
 It is also not hard to  construct similar {\it connected} $X$  for $n\geq 6$ and,  probably, 
  for all $n\geq 3$.  
  
 {\color {red!50!black} Conceivably} the same is possible for imbeddings with  higher
 codimensions   $k$,  at least for $k<<n$, where
 one expects that,   say for $k<\frac {n}{3}$ and a given, {\it arbitrarily large}, constant $C>0$, there 
 exists 
 
 {\sf a  connected $n$-dimensional submanifold  $X\subset \mathbb R^{n+k}$, such that   
  all imbeddings $X\hookrightarrow  B^{n+k}(1)$ satisfy
  $$curv(X\hookrightarrow  B^{n+k})\geq C. $$ }

 But  it should be noted that   
 
 {\it all  connected orientable surfaces}   embed to the  unit ball  $B^3$  with curvatures $\leq 100$
  
  and

{\sf the  {\it connected sums} $X$  of copies of   products of spheres with any number of summands admit   {\it embeddings}
 $$f:X\hookrightarrow B^{n+1}(1),\mbox {  $n=dim(X),$}$$
 with 
$$curv_f(X)\leq 100 n^\frac{3}{2}.$$}
 \vspace {1mm}

{\it \color {red!50!black}Questions.} {\sf Do all smooth $n$-manifold  admit embeddings to the unit $2n$-ball  with  
$$ curv(X^n\hookrightarrow B^{2n})\leq 100?$$
 
 Do the products of spheres 
 $$X=\bigtimes^m_{i=1} S^{n_i}, \mbox { where all   $n_i\geq 2$, e.g.   $X=(S^2)^m$},$$  
  embed to $B^N(1)$,
$N=1+\sum_i n_i$ with $curv(X)\leq 100$?}

\vspace {2mm}

\hspace {25mm} {\sc Lower bounds on $curv(X)$.}\vspace {1mm}

A.  It is obvious  that 

{\it immersed $n$-manifolds  $X\hookrightarrow B^N(1)$ 
with $curv(X)\leq1+\delta$ for a small $\delta>0$ keep close to an equatorial $N$-sphere in  $S^{n}\subset S^{N-1}= \partial B^N$; thus,  they are  diffeomorphic to  $S^{m}.$}\vspace {1mm}

In fact, it 
is is not hard to show, that

      \hspace{19mm}       {\it$\delta=0.01$, is small enough for this purpose,}

    \hspace{-6mm} while, {\color {red!50!black} conjecturally,}  this must hold for 
$$\delta<\frac {2}{\sqrt 3}=  1.15470...\mbox { }  $$
with the Veronese surface being the extremal one.\vspace {1mm}

B. Also {\color {red!59!black} conjecturally,}

$[\Circle\times 
\Circle]_?$  {\sf the inequality
$curv_f(X)<\sqrt 2$  
for codimension one immersions
$f:X \to B^{n+1}$
must imply that $X$  is diffeomorphic  to $S^n$} (with the equality for non-spherical  $X$ achieved by   the Clifford embeddings).

This is apparently unknown even for $n=2$..\vspace {1mm}

\vspace {1mm}

C. Let $X$ be an $n$-dimensional   $\nexists$-PSC manifold, i.e. {\it admitting no metric with} $Sc>0$, e.g. Hitchin's sphere or a connected sum of $n$-tori.

Then a simple application of Gauss's Theorema Egregium,\footnote {Compare with [Guijarro-Wilhelm(focal radius) 2017].} shows that \vspace {1mm}

{\it immersions    $f:X\to S^N$ satisfy
$$curv_f(X)\geq \sqrt\frac {n-1}{N-n} $$
and 
$$curv_f(X)\geq\sqrt {1-\frac {1}{n}}.$$
 for all $n$  and $N$.}
 
 Here, observe, it is  as it should be: no contradiction with the  above 4, for
$$ 1-\frac {1}{n} \leq 2-\frac {6}{n+2}$$
for all $n\geq 2$ with the the equality for $n=2$.\vspace {1mm}

D. If  $X=X^n$ is   $\nexists$-PSC, then all  immersions    $f:X\to B^N= B^N(1)$ satisfy
$$curv_f(X)\geq \frac {1}{C_\circ}\sqrt{\frac {n-1}{N-n}+1} $$
where $C_\circ>0$ is a universal constant that is  defined as

 \hspace {10mm} {\it the minimal possible  increase  of 
curvatures of curves} 

 \hspace {-6mm}under smooth immersions $ B^N\to S^n=S^N(1)$. More precisely, $C_\circ$ is the infimum of the numbers  $C>0$,  for which

{\it  there exits an immersion  $g: B^N\subset S^N$, such that all curves $S\subset B^N$ with curvatures 
 $$curv_{B^N}(S)\leq\sqrt { 1+\kappa^2}$$ 
 are sent to curves with curvatures 
 $$curv_{{S^{N}}}(g(S))\leq C\kappa.$$}

This $C_\circ$, most probably, is  assumed by {\it  a radial} (i.e. $O(n)$-equivariant) map $g$ and then it must be easily computable; without computation, one can  get
$$C_\circ<4.\footnote{A natural candidate for $g$   is a {\it  projective map}, where
$curv_{S^n}(g(S)\leq const_g curv_{B^n}(S)$ for {\it all} curves $S\subset B^n$.
 But since we are essentially  concerned  only with what happens to curves  with    $curv>1$, the best  $g$ doesn't have to be projective -- it might be conformal, for example.}$$

 \vspace {1mm}
 
 E. \textbf {Conjecture + Theorem.} {\sf  If   If  $X=X^n$ is   $\nexists$-PSC, then {\color {red!50!black}   conjecturally} all  immersions    $f:X\to B^N= B^N(1)$ satisfy 
$$curv_f(X)\geq const\frac {n}{N-n}. \leqno {\color {blue}\left[ \frac {n}{N-n}    \right]}$$} 

E$_1.$ {\it It is esay to see } in this regard that the $\frac {2\pi}{n}$-inequality yields this conjecture 
for $N=n+1, n+2$:\vspace{1mm}

{\it if $N=n+1$, then 
$$curv_f(X)\geq \frac {N}{2\pi}=\frac {n+1}{2\pi}.$$}
and if $N=n+2$,
then 
$$curv_f(X)\geq \frac {N}{4\pi}=\frac {n+2}{4\pi}.$$

Here  we must recall   that our proof of the $\frac {2\pi}{n}$-inequality in section \ref {bands3}
is unconditional only for $N\leq 8$, where these inequalities are not especially informative. 
And if $N\geq 9$,  our proof   relies on not formally published "desingularization"  results by Lohkamp  and by Schoen-Yau. 
 
Fortunately,  there are now  two  Dirac theoretic proofs for a large class of  $\nexists$-PSC manifolds of all dimensions, including  $n$-tori  $\mathbb T^n$ and  connected sums of these for,  example.\footnote {See [Cecchini-Zeidler(generalized Callias) 2021] and   
  [Guo-Xie-Yu(quantitative K-theory) 2020].}

E$_2.$ If  $X$ is {\it enlargeable} e.g.  the  connected sum of the $n$-torus with another closed manifold, then    a minor generalization of the Schoen-Yau "desingularization" theorem allows a  proof of  
the following version of {\color {blue}$\left[ \frac {n}{N-n}    \right]$}  for $N=n+3$:
$$curv(X\hookrightarrow B^N)\geq const_3N,$$
where, roughly, $const_3>\frac {1}{16\pi}.$

Also, granted a more serious (but realistic) generalization of the Schoen-Yau result or a version of Lohkamp's theorem, one can prove a similar inequality for $N=n+4$.
$$curv(X\hookrightarrow B^N)\geq const_4N$$
with $const_4>\frac {1}{400\pi}.$

Finally, assuming that  one can ''go around  singularities'' of stable  $\mu$-bubbles, and that (this is more serious)

 {\it \color {red!20!black} the filling radii of $n$-manifolds $Y$ with $Sc(Y) \geq \sigma>0$  satisfy
$$filrad (X) \leq  100\frac {n}{\sqrt\sigma},$$ }
one can show for all $n$ and $k=N-n$ that  
$$curv(X\hookrightarrow B^N)\geq const_kN$$
where one needs  $const_k$
 about $\frac{1} {500^{500}k}.$\vspace {1mm} 

F. All of the above  equally applies  to immersions of {\it products of enlargeable  manifolds $X_0$  with spheres}, say to 
$$f: X=X_0^{n_0} \times S^{n_1}\to B^{n_0+n_1 +k},$$
where we {\color {red!50!black}conjecture}  that 
$$curv_f(X\subset B^{n_0+n_1 +k})\geq const \frac{n_0}{n_1+k}\leqno {\color {blue}\left[ \frac{n_0}{n_1+k}   \right]} $$
and where the case $n_1+k\leq 4$ is within reach.}
 (Notice that {\color {blue}$\left[ \frac{n_0}{n_1+k}   \right]$} implies {\color {blue}$\left[ \frac {n}{N-n}    \right]$}.) \vspace {2mm}

    \hspace {38mm}                         {\sc  Four Questions} \vspace {1mm}
    
     {\sf  I. Are there lower bound on $curv_f(X)$ unrelated to the scalar curvature?\vspace {1mm}
    
    II. What is the minimal dimension  $N=N(n)$ such that all $n$-manifold can be immersed to the unit $N$-ball  with curvatures $\leq$ 1 000?\vspace {1mm}

   III What is the minimal $C=C(n)$ such that the $n$-torus can be immersed to the unit $(n+1)$-ball   with 
   $$curv(\mathbb T^n\hookrightarrow B^{n+1})\leq C?$$  \vspace {1mm}

IV Can  the Cartesian $n$-th power of the 2-sphere  be immersed  to  the unit $(2n+1)$-ball
$$X= \underset {n}{\underbrace {S^2\times ...\times S^2}}\hookrightarrow B^{2n+1}$$
with 
$$curv(X\hookrightarrow B^{2n+1})\leq 100?$$}\vspace {1mm}

 Looking back on  the above   examples, questions and  conjectures, one may be  disconcerted  by their  chaotic   irregularity. But   this only   highlights  the patchiness of  our  present-day  knowledge of 
the basic  geometry of submanifolds in  Euclidean spaces.\vspace {1mm}

{\color {blue}$\largelozenge$} {\it Wide bands with sectional curvatures $\geq 1$.} {\sf Let a proper  compact  Riemannian  band $Y$ (see \ref{bands3})  of dimension $n+1$  admit an immersion to  a  complete  $(n+1)$-dimensional  Riemannian manifold  $Y_+$ with sectional curvature
$$sect.curv(Y_+) \geq  1,$$
and let the width of $Y$ with respect to the induced 
Riemannian metric  satisfy  
$$width(Y)=dist (\partial_-Y,\partial_+Y)> \frac {\pi}{2}.$$}
Then  
\vspace {1mm}
 
{\it $Y$ contains a subband $Y_-\subset Y$  of width $d=width (Y)> \frac {\pi}{2}$,  
  which is homeomorphic  to the spherical cylinder  $S^{n}\times [0,1]$.}\vspace {1mm}

{\it Acknowledgement. }  A similar result for $n=3$ is proved  in  [Zhu(width) 2020], while 
our argument below follows  that of   Jian Ge  from [Ge(linking) 2021], who sent me his preprint prior to  publication.\vspace{1mm}
\vspace {1mm}

{\it Proof.}  Let  $Y_-$
be the  intersection of   the   $d$-neighbourhoods of the  $\partial_\mp$-boundaries of $Y$,
$$Y_-=U_d(\partial_-)\cap U_d(\partial_+), $$  and observe that the 
 $\partial_\mp$-boundaries of this $Y_-$ are {\it  concave}  for $\kappa\geq 1$ and $d> \frac {\pi}{2}$.
Therefore, $\partial_\mp$ are diffeomorphic to $S^{n-1}$ and 
 the immersions
$$ \partial_\mp \to Y_+$$ 
 extend to immersions of $n$-balls, such that the  {\it locally convex} boundaries of these are equal to  $\partial_\mp$ (with their   coorientations opposite to those in $Y_-$).
\footnote{Recall that {\it a closed  immersed  locally convex} hypersurface in a complete Riemannian manifold of dimension $n\geq 3$  with
 sectional curvatures $ >0$ bounds an immersed ball.} 
 
 It follows, that if $Y_+$  is simply connected, then then  the immersion $Y_-\to Y_+$  is  one-to-one and the complement $Y_+\setminus Y_-$  consists  of two convex balls with 
 distance $>\frac {\pi}{2}$ between them.
 
Hence, $diam(Y_+)>\frac {\pi}{2}$  and   $Y_+$ is homeomorphic to $S^{n+1}$ by   the {\it Grove-Shiohama diameter theorem};  consequently,  $Y_-$ is homeomorphic  to $S^{n} \times [0,1]$. QED.
 \vspace{1mm}

 {\it Remark.} (a) The  conclusion of the  theorem, probbaly,  holds if $sect.curv(Y_-)\geq 1$ and 
 $sect.curv(Y_-)\geq 0$, since the proof of the diameter theorem seems to work in this case.
 
 (b)  It also   {\color {red!50!black} doesn't seem  difficult}   to prove the  rigidity theorem a la {\it Berger-Gromoll-Grove}   in  case of an
  {\it open band  with} 
 $width(Y)=\frac {\pi}{2}$, where 
  the only alternatives to the homeomorphism of $Y$ to $S^{n}\times (0,1)$  should be  as follows:
  
  $\bullet$ $Y$ is isometric 
   the open $\frac {\pi}{4}$-neighbourhood of a Clifford submanifold   
  $$S^{n_1}\times S^{n_2}\subset S^{n+1}\mbox {   }  n_1+n_2=n;$$ 
 
 $\bullet\bullet$ $Y_+$   is isometric to the  {\it projective 
space} over complex numbers,  quaternion numbers or Cayley numbers  and $Y$ is isometric to the open $\frac {\pi}{2}$-ball minus the center
 in such an $Y_+$.

   In fact, the poof of this  rigidity seems quite  easy in the case of the interest (the above $[\Circle \times \Circle]$),   where $Y$ is equal to the 
   (normal) $\frac {\pi}{4}$-neighbourhood of a hypersurface  $X\subset S^{n+1}$ with $curv(X)\leq 1$.
    \vspace {1mm}

 {\it Questions.}   (i)   Is the manifold $Y_+\supset Y$  indispensable?   
 Do there exist "non-obvious" bands with $sect.curv\geq 1$ and with $width \geq \frac {\pi}{2}?$

 (ii) Given a closed $n$-manifold $X$, e.g. a  product of spheres, $X=\bigtimes_iS^{n+1}$,    what is the supremum  of the widths of the Riemannian   bands $Y$ 
 homeomorphic to $X\times [0,1]$  with $sect.curv(Y)\geq 1$?

%%%%%%%%%%%%%%%%%%%%%%%%%%%%%%%

\subsection {\color {blue}  Multi-Width   of Riemannian Cubes }\label {multi-width3} 

%%%%%%%%%%%%%%%%%%%%%%%%%
{\sf Let $g$  be a Riemannian metric on the cube $X=[-1,1]^n$ and let $d_i$, $i=1,2,...,n$, denote the $g$-distances  between the pairs of the opposite faces denoted $\partial_{i\pm}=\partial_{i\pm}(X)$
   in this cube $X$, that are the length of the shortest  curves between $\partial_{i-}$ and 
   $\partial_{i+}$ in $X$.  \vspace {1mm}

 {  \color{black}  $\square^n$-\textbf {Inequality}.}  {\it  If $Sc(g)\geq n(n-1)=Sc(S^n)$, then  

 $$\sum_{i=1}^n \frac{1 }{d_i^2}\geq  \frac  {n^2}{4\pi^2} \leqno {\color {blue}\square_{\sum}}$$
In particular, 
$$\min_i dist (\partial_{i -},\partial_{i +})\leq \frac {2\pi}{ \sqrt n}.\leqno {\color {blue}\square_{\min}}$$}}
 (On the surface of things, this inequality is {\it purely geometric}  with no topological string attached.   But in truth, the combinatorics of the cube fully reflects   toric topology in it.)
\vspace {1mm}
 
{\it \textbf {$\frac {2\pi}{n}$-Corollary.}} {\sl If $X$ is   a proper  orientable  non-compact band  with 
$Sc(X)\geq n(n-1)$, which  admits a proper
 1-Lipschitz map  $f: X\to\mathbb R^{n-1}$, such that the  restriction of $f$ to  the $\partial_\pm$-components of  the boundary, of $X$,
  $$\partial_-(X),\partial_+(X) \to  \mathbb R^{n-1},$$ 
  have  non-zero degrees,(these two degrees are mutually equal) then  
 $$width(X)=dist(\partial_-,\partial_+)\geq \frac {2\pi}{n}.\footnote{A   proof of this     
 for bands with locally bounded geometries  can be performed using minimal  hypersurfaces rather than $\mu$-bubbles  as it is  (briefly and sloppily) indicated in section 11.7 in [G(inequalities 2018].}$$}

 {\it  The proof} of {\color {blue}$\square_{\sum}$ }  proceeds by  inductive dimension descent with $\mathbb T^\rtimes$-symmetrization with the use of the    "separation with symmetry" theorem {\color {blue} {\large \sf I{\large \sf \color {red!80!black}\sf I}\sf I}$_\circlearrowleft$}  from section \ref{separating5}. %{\color {red} (added in June!!!!!)} 
  \footnote {A Dirac theoretic proof of this inequality is given in the recent paper   [Wang-Xie-Yu(cube inequality) 2021].}
 \vspace {1mm} 
 
  \vspace {1mm}

 {\it \textbf {Generalization.}}  We shall apply this argument  in \ref {separating5} to more general "cube-like" manifolds   $X$, such as products of surfaces with square-like decompositions of their boundaries and also  to products 
  $Y_{-m}\times [-1,1]^{n-m}$, where this  yields  inequalities mediating between     $\square_{\sum}$ and the $\frac { 2\pi}{ n}$-inequality.

  \vspace{1mm}
  
  { \it \color{blue}  $\square^2$-Example.} {\sf Let $Z$ be a compact
   connected  orientable surface with non-empty connected  boundary   where this (circular) boundary $S=\partial Z$  is decomposed into four segments meeting at their ends,
  $$S=S_{1+}\cup S_{2+}\cup S_{1-}\cup S_{2-}.$$
Let $g$  be a Riemannian metric on $Z\times \mathbb T^{n-2}$  with
 $Sc(g)\geq \sigma>0$.}
 
 {\it Then the $g$-distances between the products of the pairs  of the opposite (i.e. non-intersecting) segments in $S$ by the torus $\mathbb T^{n-2}$,
 denoted $\partial_{i\pm}= S_{i\pm}\times \mathbb T^{n-2}\subset Z\times \mathbb T^{n-2}$, $i=1,2$,  
 satisfy:
 $$\min_{i=1,2} \left (dist_g( \partial_{1-}, \partial_{1+}),    
 dist_g(\partial_{2-}, \partial_{2+})\right)  \leq
  2\sqrt 2 \pi\cdot \sqrt\frac  {n-1}{n}\cdot \sqrt\frac  {1}{\sigma}.  \leqno {\color {blue}[2\sqrt 2]} $$}
 
 {\it Proof.} Pass to $Z\times \mathbb R^{n-2}$   for the universal covering 
   $\mathbb R^{n-2}\to  \mathbb T^{n-2}$  and apply the   $\square^n$-{inequality} to 
   $Z\times [-d, d]^{n-2}$  for $d\to \infty$.

 \vspace{1mm}
 
  {\it\large \color {blue!50!black} Hemi-spherical Corollary.} {\sf Let $X$ be a Riemannian manifold with $Sc(X)\geq n(n-1)=Sc(S^n)$,  which admits a  $\lambda_n$-Lipschitz,  (i.e. {\color {blue!40!black}$dist(f(x) f(y))\leq \lambda_ndist(x,y)$}) homeomorphism onto the hemisphere  $S^n_+$,
  $$f:X\to S^n_+.$$
Then $$\lambda_n\geq \frac {\arcsin \beta_n}{ \pi \beta_n}> \frac {1}{ \pi} \mbox 
 { for } \beta_n=\frac {1}{\sqrt n}.$$}

{\it Proof.} The hemisphere   $S^n_+$ admits an obvious  cubic decomposition with the (geodesic) edge length $2\arcsin \frac {1}{ \sqrt n}$ and ${\square_{\min}}$ applies to 
the pairs of the  $f$-pullbacks of the faces of this decomposition.
 
 \vspace {1mm}

{\it Remarks and Exercises.} (a) This lower  bound on $\lambda_n$ improves those in \S 12 of [GL(complete) 1983] and in \S 3 of [G(inequalities) 2018].

Moreover the {\it sharp} inequality for Lipschitz maps to  the punctured sphere stated in  the next section   implies that  
$\lambda_n\geq \frac{1}{2}$  for all $n$.\vspace {1mm} 

But it remains  {\it \large problematic}  if, in fact, $\lambda_n\geq {1}$  for all $n\geq 2$.\vspace {1mm} 

(b)   Show that $\lambda_2\geq 1$. \vspace {1mm}

(c) The proof of the  inequality $\square_{\sum}$ in section \ref {separating} applies to {\it proper} ((boundary$\to$ boundary) {\it   $\lambda$-Lipschitz} maps with {\it non-zero degrees} from  all compact   connected orientable manifolds $X$ to $S^n_+$, while  the proof via punctured spheres needs $X$ to be spin.

\vspace {1mm}

 (d)  Show that  {\sf the   Riemannian metrics  with  sectional curvatures $\geq 1$ on the square   $[-1,1]^2$    satisfy 
$$\min_{i=1,2} dist (\partial_{i -},\partial_{i +})\leq\pi.\leqno {\square^2_{\min}}.$$}

  \vspace {1mm}

 (e)  Construct iterated warped product metrics $g_n$ on the $n$-cubes $[-1,1]^n$  with $Sc(g_n)=n(n-1)$,  where, for $n=2$,  both   $d_i$, $i=1,2$,  are equal to  $\pi$  and such that 
  $$d_i>2\arcsin \frac {1}{ \sqrt n},\mbox { } i=1,...,n, \mbox  { for all  }n=3,4,... ,\mbox { } .$$

(f) Show, that  {\color {blue}$\square_{\min}$} is equivalent to the {\it over-torical }  case of {\color {blue} \large \it $\frac { 2\pi}{ n}$-Inequality.}
modulo constants. Namely,\vspace {1mm}

(i). {\sf If   a Riemannian $n$-cube  $X$ has  $\min_i dist (\partial_{i -},\partial_{i +})\geq d$,
then it contains an $n$-dimensional  Riemannian band  $X_\circ\subset X$, where
$dist (\partial _-X_\circ,\partial _+X_\circ)\geq \varepsilon _n\cdot d$, $\varepsilon _n>0$, and where 
$X_\circ$ admits a continuous map to the $(n-1)$ torus,  $f_\circ : X_\circ \to \mathbb T^{n-1}$, such that all closed  hypersurfaces $Y_\circ \subset X_\circ$ which separate $\partial _-X_\circ$ from $\partial _+X_\circ$ are sent by $f_\circ$ to  $ \mathbb T^{n-1}$ with {\sl non-zero degrees}.\vspace {1mm}}

(ii).  Conversely, { \sf let   $X_o$ be a  band,  where $dist (\partial _-X_o,\partial _+X_o)\geq d)$  and which 
 admits a continuous map to the $(n-1)$ torus,  such that the hypersurfaces $Y_o \subset X_o$, which separate $\partial _-X_o$ from $\partial _-X_o$, are sent  to  this torus  with {\sl non-zero degrees.  
 
 Then there is a (finite if you wish)  covering  $\tilde X_o$ of  $X_o$, which contains a domain $X_\smallsquare\subset \tilde X_o$,  
where this domain admits a continuous proper    map  of degree one onto the $d$-cube $f_\smallsquare: X_\smallsquare\to(0, d)^n$,
 such that the $n$  coordinate projections of this map, $(f_\smallsquare)_i :X_\smallsquare\to (0, d)$,  are distance decreasing. }}
 
%%%%%%%%%%%%%%%%%%%%
 \subsection {\color {blue}  Extremality and Rigidity of Punctured Spheres } \label {punctured3}

%%%%%%%%%%%%%%

{\sf  Let $\underline X$ be the unit  sphere $S^n$ minus two opposite points $\pm x_0\in S^n$
 and let 
 $ \underline g=g_{sphe}$ denote  the spherical metric (of constant curvature +1) restricted to this 
  $\underline X= S^n\setminus \{\pm x_0\} \subset S^n$. }\vspace {1mm}
 
\textbf {Double Puncture Extremality/Rigidity Theorem.}  {\it If a smooth metric $g$ on $\underline X$  satisfies 
$$g\geq \underline g  \mbox  { and } Sc(g)\geq n(n-1) =Sc(\underline g),$$  \vspace {1mm}
then $g=\underline g$.}

This is shown by applying the  spin-area   extremality theorem {\color {blue}[$X_{spin}{^\to}$\Ellipse]} from section 
 \ref{area extremality3}  (one needs here only the spherical case of it but sharpened by  rigidity in the case of equality) 
 to the $\mathbb T^1$-symmetrization of   a certain  stable  $\mu$-bubble, $Y\subset  S^n\setminus \{\pm x_0\}$, which separates the punctures $\pm x_0\in S^n$. 
 
 (See section \ref {log-concave5}  for the proof of this 
for general  {\it spin} manifolds  with $Sc\geq n(n-1)$ {\it properly  mapped} to  $S^n\setminus \{\pm x_0\}$ with $deg\neq 0$, where, recall,  
the details of this proof for  $n\leq 8$  are yet to be worked out.)

\vspace{1mm}

 {\it Remark.} If the above   metric $g$ on $\underline X=S^n\setminus \{\pm x_0\}$  is {\it complete}, one can prove that the inequalities $g\geq \underline g$   and  $Sc(g)\geq n(n-1)$ imply that $g\geq \underline g$ for the complements 
 $\underline X=S^n\setminus \Sigma$ for certain 
 subsets $\Sigma$ larger than $\{\pm x_0\}$.

 For instance, Llarull's inequality  implies  this for all  {\it finite subsets} $\Sigma\subset S^n$ 
and  a similar (purely index theoretic)  argument yields this for

{\it piecewise smooth 1-dimensional subsets} (graphs)
 $\Sigma\subset S^n$, such that the {\it monodromy transformations}  of the principal tangent $Spin(n)$-bundle (that is double cover of the orthonormal tangent frame-bundle over all closed
  curves in $\Sigma$ are {\it  trivial} (e.g. $\Sigma$ is contractible).

 \vspace{1mm}

 But  if one makes {\it no   completeness assumption}, our proof  is limited to 
  $\Sigma$ being either empty, or consisting of  a single point or of  a pair of opposite points.

 \vspace{1mm}

{\it Exercise.} Prove with the above that no metric $g$ on the hemisphere $(S^n_+, \underline g)$  
can satisfy the inequalities $g\geq 4\underline g$ and $Sc(g)>n(n-1)$. Then  directly show 
that if $n=2$ then the inequality  $g\geq \underline g$ and $Sc(g)\geq 2$ imply that 
$g=\underline g$.\vspace{1mm}

{\it \large \it Questions.} (a)  {\sf Does the implication
$$[g\geq \underline g]\&  [Sc(g)\geq n(n-1)] \Rightarrow g=\underline g $$  \vspace {1mm}
ever hold for $S^n\setminus \Sigma$ apart from the above cases?

(b) Can the sphere  $S^n$ with $k$-punctures carry a metric $g$, such that   $[Sc(g)\geq n(n-1)] $ and such that the $g$-distances between these  punctures are all $\geq 10^{nk}$?}\footnote {The negative answer   was recently  delivered by Cecchini's   {\it long neck theorem}, see section \ref{Roe3}.}

%%%%%%%%%%%%%%
\subsection {\color {blue}   Slicing and Sweeping  $3$-Manifolds and Bounds on their
Widths and Waists .} \label {slicing3D.3}

%%%%%%%%%%%%%%%%%%%%

If  $n\geq 4$, then then all  known bounds  on the size of $n$-manifolds $X$ with $Sc(X)\geq \sigma>0$ are  expressed by 
{\it non-existence}  of   "topologically complicated but geometrically simple" maps from  these $X$ to  "standard 
manifolds" $\underline X$. \vspace{1mm}

But if $n=3$ then \vspace {1mm}

{\sf  complete 3-Manifolds $X$ with scalar curvature $Sc(X)\geq \sigma>0$ are

 known to  satisfy the following properties \textbf {\color {blue} {A, B, C}}  }

 \vspace {1mm}

{\sf \textbf {\color {blue}A.}} \textbf {Uryson's 1-Width Estimate.}  {\sf Let $X$ be a complete Riemannian
3-manifold   with $Sc(X)\geq \sigma>0$.}

{\it Then there exists a continuous  map $f:X \to P^1$, where $P$  is   a 1-dimensional polyhedral space (topological graph), 
such that the diameters of the pullback of all points are bounded by 
$$diam(f^{-1} (p))\leq \frac {24\pi}{\sqrt \sigma}. \leqno {\color {blue}[width_{3/1}]} $$
\vspace {1mm}

{\sf {\color {blue} \sf \textbf A$'$}.} Moreover, if the rational homology group  $ H_1(X,\mathbb Q)$ {\sf vanishes}, then the diameters of the  connected components of the levels of the  distance function 

\hspace {41mm}$x\mapsto dist( x_0, x)$ 

\hspace {-6mm}are bounded by $ \frac {8\pi}{\sqrt \sigma}.$ for all $x_0\in X$.}\vspace {1mm}

We prove a $ \mathbb T^\rtimes$-stabilized version of  {\color {blue} $\mathbf A'$} for manifolds with mean convex boundaries in the next section and then derive 
 {\color {blue} $\mathbf A$}, also in the   $ \mathbb T^\rtimes$-stabilized form needed for applications,  for all 3-manifolds $X$
 with $Sc(X)\geq 6$.\footnote{{\color {blue} $\mathbf A'$}  is proven in [GL(complete)1983], where the condition $H_1(X,\mathbb Q)=0$ was erroneously omitted.  Also see sections (E)-(E$'_2$) in Appendix 1 in [G(filling). 1983].}\vspace {1mm}

{\sf \color {blue}  \textbf A$''$.}  \textbf {Corollary.} {\it The filling radius of  a complete $3$-manifold $X$ with $Sc(X)\geq \sigma$ is bounded by
$$ fil.rad(X)\leq  C\cdot  {\frac {1}{\sigma}} \mbox { for } C\leq   24\pi.$$}
(We shall show in the next section that  $C\leq   8\pi.$)

\vspace {1mm}

{\it Exercises.} (a)  Map the unit sphere  $S^n\subset \mathbb R^{n+1}$  onto the cone  $P^1=P^1_n$ over the vertex set of the regular simplex  inscribed into  $S^n$,  such that the diameters of the pullbacks of all points are
 $\leq \pi-\delta_n$ for $\delta_n>0$.
 
 ({\sf Probably, The Uryson   1-width of $S^n$ is realized by such a map}.)\vspace {1mm}

 (b) Let $X$ be a compact Riemannin  3-manifold with $Sc(X)\geq \sigma>0$ and let $f:X\to X$ be a continuous map. Show that if the 1-dimensional homology of $X$ is torsion, e.g. zero, then there exists a point $x\in X$, such that $dist(x, f(x))\leq 6\pi\sqrt { \frac{2}{\sigma}}.$
 
 {\it Hint.}  See (E$'_5$) in Appendix 1 in [G(filling) 1983].\vspace {1mm}
 
 {\color {blue} \textbf B.} \textbf {Topological $\bigsqcup S^2$-Sweeping Theorem.}  {\sf A compact $3$ manifold admits a Riemannian metric with   with $Sc>0$}
  {\it if and only if there exists a finite covering $\tilde X\to X$ and a Morse function $\tilde f:\tilde X\to\mathbb R$,  such that  the pullbacks
 $\tilde f^{-1}(t)\subset X$ of all  non-critical   $t\in \mathbb R$ are disjoint unions or spheres $S^2$.}
 
\vspace {1mm}

{\it About the Proof.} This is a reformulation of the classification theorem for compact manifolds $X$ with $Sc>0$, which says, in effect,    that 
these $X$ 

\hspace {2mm}{\it 
admit finite coverings $\tilde X$  diffeomorphic to connected sums of  $S^2\times S^1$}, 

\hspace {-6mm}and which follows from 
{\it non-existence of aspherical components} in the prime decompositions of manifolds with $Sc>0$, and  {\it Perelman's solution of Thurston conjecture}. \footnote {See [Gl(complete) 1983] and  [Ginoux(3d classification) 2013].}
\vspace {1mm}

 {\color {blue} \textbf B$'$.} {\it \textbf {$S^2$-Sweeping Complete Manifolds}.} {\sf If $X$ is a complete oriented 3-manifold  with $Sc(X)\geq \sigma>0$,} then, 
 instead of a finite covering $\tilde X$,  one  constructs 
 
 {\it a 3-polyhedron $\hat X$, a proper  piecewise smooth locally 
 finite-to-one  map 
 $\hat \Phi: \hat X\to X$ and a proper  piecewise linear positive  function $\hat f: \hat X\to \mathbb R_+$,} such that 

(i) {\sf the map $\hat \Phi$ sends a (non-compact) homology class from the {\it rational} 3-dimensional homology group of  $\hat X$ with infinite supports to the fundamental class of $X$};

(ii) {\sf the connected components of the  pullbacks  $\hat f^{-1}(t)\subset \hat X$ for all $t\in \mathbb R$ are either single points or   joints of  2-spheres. }

\vspace {1mm}

(This is suggestive of what can be expected for $n>3$.)\vspace {1mm}

 {\color  {blue} \textbf C.} \textbf  {Sharp Area Slicing  Inequality.} {\sf Let $X$ be a Riemannian  $3$-manifold diffeomorphic to $S^3$ or to a connected sum of several  $S^2\times S^1$. }
 
 {\it If  $Sc(X)\geq 6$, then $X$ admits a Morse function $f$, the non-singular levels of which are disjoint union of spheres, where  the areas of all these spheres are bounded by $4\pi$.}
 
\vspace {1mm}

{\it About the Proof.}  We already  know the  all stable minimal surfaces $Y$ in $X$ have areas bounded by 
$\frac {4\pi}{3}$ by  Schoen-Yau's rendition of the second variation  inequality (sections \ref{2nd variation2} and \ref{SY+symplectic2}   + the Gauss-Bonnet theorem. Furthermore, this  inequality combined with   {\it Hersch's upper bound  on  the first non-zero eigenvalue} of the Laplace  on surfaces $Y$ diffeomorphic to $S^2$  with $area(Y)\geq 4\pi= area (S^2)$, that is
       $$\lambda_1(Y)\geq 2=\lambda_1(S^2)$$
implies that minimal surfaces $Y\subset X$  with Morse index 1  have their areas bounded by $4\pi$. \footnote {See J. Hersch, {\sl Quatre propri\'et\'es isop\'erim\'etriques de membranes sph\'eriques homog\`enes}, C.R.
Acad. Sci.Paris S\'er. A-B 270 (1970), A1645-A1648  and  [Marques-Neves(min-max spheres in  3d)  2011].}

Then "the almost extremal   $\bigsqcup S^2$-Morse slicing" $f:X\to \mathbb R$,  that almost minimizes the area of the maximal pullback sphere is the required one. \footnote {See [Lio-Max (waist inequality) 2020], where this is proved   using the mean  curvature flow.

Probably, this can be also proved by the  {\it Sacks-Uhlunbeck  direct  minimization}, where  bubbling  creates disconnectedness of the levels of $f$. 
(Apparently, if $X$  diffeomorphic to $S^3$  contains no stable minimal surfaces, then it admits a  Morse function with two critical points and areas of all levels bounded by $4\pi$. But, in general 
high   disconnectedness  of the levels of $f$   is   inevitable, even for $X$ diffeomorphic to $S^3$.)}

\vspace {1mm}

 {\color  {blue} \textbf C.} \textbf  {Liokumovich-Maximo  Area+Diameter Slicing  Inequality.}   {\sf Let $X$ be a compact Riemannian 3-manifold with $Sc\geq 6=Sc(S^3)$.} 
{\it Then $X$  admits a Morse function,  the connected components $\Sigma$ of  all nonsingular levels $f^{-1}(t)\subset X$, $t\in \mathbb R$ satisfy}: 

(i) $area(\Sigma)\leq 64\pi)$,

(ii) $diam(\Sigma) \leq \frac {40\pi}{\sqrt 6}$,

(iii) $genus(\Sigma)\leq 13$.

\vspace {1mm}

\textbf {Corollary}.  {\it $X$ admits a map $F:X\to \mathbb R^2$, such that the lengths of the pullbacks of 
all points are bounded by a universal constant $C\leq 100$.}

   {\sf Consequently},{ $X$ {\it contains  a
stationary geodesic net of $length \leq  C$.}\vspace {1mm}

For the proof we refer to [Lio-Max (waist inequality) 2020].

\vspace{1mm}

All {\it known} manifolds with $Sc\geq \sigma>0$ satisfy counterparts of these   \textbf {\color {blue} {A, B, C}}  
for all  dimensions $n$, which suggests the following conjectures.

\vspace {1mm}

 {\color  {blue} \it \textbf {Topological $S^2$-Sweeping   {\color {red!50!black} Conjecture}}.}  {\sf Let $X$ be a complete, e.g. compact, orientable  $n$-manifold, $n\geq 3$   with 
 $Sc(X)\geq \sigma>0$. }
 
 {\it Then there exists an $n$-polyhedron $\hat X$, a proper  piecewise smooth locally 
 finite-to-one  map 
 $\hat \Phi: \hat X\to X$ and a proper  piecewise linear map $\hat f: \hat X\to  P^{n-2}$,}  where  
 $P^{n-2}$   is an $(n-2)$-dimensional  polyhedral space} ({\sf pseudomanifold maybe?}),  such that 

(i) {\sf the map $\hat \Phi$ sends a (non-compact if $X$ is non-compact) homology class from the {\it rational} n-dimensional homology group of  $\hat X$ (with infinite supports if $X$ and $\hat X$ are non-compact)  to the fundamental class of $X$};\footnote{This condition ensures  "homotopical surjectivity" of  this map, that is non-existence of its (proper in the non-compact case)   homotopy to a map into a subset   $ Y  \subsetneqq X$. I am not certain if another  such condition is  relevant here.}

(ii) {\sf the connected  of the  pullbacks  $\hat f^{-1}(t)\subset \hat X$ for all $t\in \mathbb R$ are either single points or   joints of  2-spheres. }

{\it \color {blue}Corollary.}  {\sf If a compact orientable manifold $X$ admits such an $\hat \Phi: \hat X\to X$, then
all continuous maps from $X$ to aspherical spaces  induce {\it zero homomorphisms} on $H^n(X)$.}

(This remains unknown for  manifolds $X$ with $Sc(X)>0$  of dimensions $n\geq 4$, but   a weaker property -- 
{\it non-contractibility of the universal covering of $X$} -- is confirmed by the Chodosh-Li theorem we prove in the next section.)
\vspace{1mm}

{\it \color {blue} \textbf {Reversed $S^2$-Slicing  {\color {red!50!black} Conjecture}.}} {\sf Let a complete (e.g.compact) smooth  manifold $X$ admits a (proper in the non-compact case)    piecewise linear map with respects a smooth triangulation of $X$  to a pseudomanifold $P^{n-2}$, such that the connected components of the pullbacks of the points are   either single points or   joints of  2-spheres.} 

{\it Then $X$ admits a metric with $Sc>0$.} 

(Possibly, one should  add  an extra condition on  singularities of such a map.)

\vspace {1mm}

\vspace{1mm}

 {\color  {blue} \it \textbf  {Width/Waist  {\color {red!50!black} Conjecture}.}} {\sf All complete $n$-manifolds $X$ with $ Sc(X)\geq n(n-1)$ admit continuous maps to 
polyhedral spaces of dimension $n-2$, say,
$F:X\to P^{n-2}$, such that  
$$ diam (F^{-1}(p))\leq const_n  \mbox    { and } vol_{n-2}(F^{-1}(p))\leq const'_n \mbox { for all }p\in P^{n-2}.$$}

Observe in this regard the following.

$\bullet $ The existence of a proper  map $F:X\to P^{n-1}$ with  
$diam (F^{-1}(p))\leq C$ would imply that $fill.rad(X)\leq C$, that remains unknown if   $n\geq 4$  even for universal coverings of compact $n$-manifolds with $Sc>0$.

 $\bullet $ The bound  $fill.rad(X)\leq C<\infty$ implies that the balls in $X$ {\it can't  be all contractible}; moreover,
 
 {\it  given a continuous  function $R(r)\geq r$, $r\geq 0$,  and a number $r_0,$  there exits  a ball of  radius $r\geq r_0$ in $X$, which is non-contractible in
 the concentric $R(r)$-ball.}
 
This {\it uniform non-contractibility property},  remains conjectural for $n\geq 6$ but we prove it for  $n=4,5$ in the next section.  
 
 The only known result of this kind, which implies a (sharp) bound on the injectivity radius of manifolds with $Sc\geq \sigma$ is\vspace {1mm}

{ \color {blue} Green-Berger Integral Scalar Curvature Inequality.}
{\it Among all compact manifolds $X$ with given $vol (X)$ and the integral $\int _XSc(X,x)dx$, the   round
spheres maximize the average distance between conjugate points
on geodesics.}\footnote { M. Berger, {\sl Lectures on geodesics in Riemannian geometry}, Tata Institute
of Fundamental Research, 1965.

Berger's proof of this  applies to  complete non compact {\it amenable} manifolds $X$ with $Sc(X)\geq n(n-1)$, e.g. with {\it subexponential volume growth},  thus providing a bound 
$inj.rad\leq C_n<\infty$.  But I don't see offhand   how to prove  such a bound for non-amenable  $X$.}

%%%%%%%%%%%%%%%%
\subsubsection {\color {blue} Filling Radii of  $3$-Manifolds, Hyperspherical Radii,  Enlargeability and Uniform  Asphericity}\label{filling+hyperspherical+asphericity3}

%%%%%%%%%%%%%%%%%%%%%%

Recall that the Uryson k-width  $width_m(X)$ of a metric space $X$ is the infimum of the numbers $d\geq 0$, such that 
 $X$ admits a continuous map  a   $m$-dimensional polyhedral space  $P^m$,  such that that the diameters of 
the pullbacks of all points $p\in P^m$ are bounded  by $d$.  

 {\it \textbf {Lemma {\sf (A)}}.}  {\sf  Let  a {\it proper}\footnote {Closed  bounded subsets are compact.} locally contractible  metric space  $X$  be covered by closed locally contractible  subsets  $X_i$, $i\in I$, such that

$\bullet_1$ there is {\it no triple intersections} between $X_i$;

$\bullet_2$ the {\it connected components}  $Y=Y_{ijk}\subset X_i\cap X_j$  of all  double intersections are    locally contractible,\footnote  {Probably, "locally contractible" is an unnecessary precaution,
but in our case $X_i$ are  manifolds with boundaries,  where the intersections   $X_i\cap X_j$ are unions of connected components of common boundaries of $X_i$ and $X_j$.}  and  their {\it rational homology $H_1(Y_{ijk};\mathbb Q)$ vanish};   
 
$\bullet_3$  the diameters  of all these $Y$ are  {\it bounded by} 
   $$diam(Y_{ijk})\leq \delta_{ijk}<\infty.$$

Then 
 $$width_1(X) \leq  \sup_{ijk} (2\delta_{ijk}+ width_1(X_i)).$$}

 {\it Proof.} Let   $\chi_i:X_i\to P^1_{i}$ be a continuous map, such that the pullbacks  $\chi_i^{-1}(p)$, $p\in P^1_i$, are  bounded by 
 $$diam(\chi_i^{-1}(p))\leq d_i \mbox { for all } p\in P^1_i $$
 and observe that  the union $U_p$ of the $\chi_i$-pullback  $\chi_i^{-1}(p)$ of  $p\in P^1_i$  with the components $Y_{ijk}\subset X_i$ for which $\chi(Y_{ijk})\ni p$ satisfies:
 $$diam(U_p) \leq diam (\chi_i^{-1}(p))  + 2\sup_{j,k} \delta_{ijk} \leq d_i+ 2\sup_{j,k} \delta_{ijk}\mbox  { for all } p\in P^1_i.  $$
 
 Since  $H_1(Y_{ijk};\mathbb Q)=0$ the $\chi_i$-maps  from $Y_{ijk}$ onto their images in $P^1_i$ are {\it contractible}. Therefore,

 {\sf given an arbitrary small neighbourhood $U'_{i}\subset X_i$ of the union $\bigcup_{j,k}Y_{ijk}\subset X_i$,
  there exists a map   $\chi'_i: X_i\to P_i^1$ homotopic
  to $\chi_i$ , such that }
 
 $\bullet$ $\chi'_i$ is {\it constant} an all $Y_{ijk}\subset X_i$;

$\bullet$ $\chi'_i$ is {\it equal to $\chi_i$ outside}  $U'_{i}$, where we may assume that 
  $U'_{i}=\bigcup_{j,k} U'_{ijk}$ for  the  connected components   $U'_{ijk}\supset Y_{ijk}$  of $U'_i$;
 
 $\bullet$ the image $\chi'_i(U_{ijk}\subset P^1_i$ is contained in the image  $\chi_i(U_{ijk})$.
  
It follows that 
$$diam(\chi'_i)^{-1}(p)\leq diam(U_p)+\varepsilon_i \leq  d_i+ 2\sup_{j,k} \delta_{ijk}+\varepsilon_i,$$
 where  $\varepsilon_i$ can be made   arbitrarily small for small $U_i$.

Finally, we glue  the graphs $P^1_i$ and $P_j$ at the points $p_{ijk}\chi'_i(Y_{ijk}\in P^1_i$
and $p_{jik}\chi'_j(Y_{jik}\in P^1_j$ and let $P^1=\bigcup_iP^1_i$  be the resulting graph. 
 
Then the obvious map 
$$\chi': X=\bigcup_iX_i\to P^1=\bigcup_iP^1_i$$
satisfies
$$diam((\chi')^{-1}(p))\leq \sup_i diam((\chi_i')^{-1}(p)) \leq\sup_{ijk}
 (2\delta_{ijk}+ width_1(X_i)  +\varepsilon_i)$$
which  concludes the proof    for $\sup_i\varepsilon_i\to 0$.\vspace {1mm}

 {\it \textbf {Remark {\sf (A$'$)}}.}  The condition $H_1(Y_{ijk};\mathbb Q)=0$ looks  strange here  but I don't know if it 
can be omitted.

On the other hand,  the 2-width of $X$ can be bounded by
$$width_2(X)  \leq \max (  \sup_i width_1(X_i), \sup _{i,j,k} diam(Y_{ijk}))\leqno {\color {blue}[width_2]}  $$
with no such condition,
where the relevant 2-polyhedron where $X=\bigcup_iX_i$  goes  is 
obtained  from the cones $C_{jk}(X_i)$ by  gluing the apexes  of $C_{jk}(X_i)$  with these of 
$C_{ik}(X_j)$.
\vspace{1mm}

Next  -- this is a provisional definition adapted to our present purpose  --    define {\it  maximal  $\square$-width} 
of a metric space  $S$ {\it homeomorphic to the circle},  denoted $max.width_\square (S)$, as the maximum of the numbers $D$, such that $S$ admits a decomposition into four segments,  with the same combinatorial arrangement as that of the four faces of the square $[-1,1]^2$ and such that the distances between  both pairs of  opposite (i.e. non-intersecting)  segments are  {\it bounded from below by  $ D$. }
\vspace {1mm}

{\it \textbf {Lemma}} {\sf (B)}. {\sf Let  a circle   in a proper  metric space, say $S\subset X$,  with the induced metric satisfies
$$ max.width_\square (S)\geq D.$$
Then there exists a pair of  non-negative 1-Lipschitz functions on $X$ that define a  {\it proper}  map 
$$\Psi_\square =(\psi_1,\psi_2):X\to \mathbb R^2_+, $$
such that}

{\it  $S$ is sent by $\Psi_\square$ to the complement of the interior   of the square $[0,D]^2\subset  \mathbb R^2_+$,
where   the induced homology homomorphism
$$\mathbb Z=H_1(S)\to H_1( \mathbb R^2_+\setminus   (0,D)^2)=\mathbb Z$$
is an isomorphism.}

{\it Proof.} Let 
$$S=S_{1+}\cup S_{2+} \cup  S_{1-}\cup S_{2-},$$ 
where    
$$dist (S_{1+},S_{1-})\geq D\mbox {   and } dist (S_{2+},S_{2-})\geq D$$
and observe that the map defined by
$$\psi_i(x)= dist(x, S_{1+}),  \mbox {  } i=1,2,$$
is the required one. 
\vspace {1mm}

{\it \textbf {Corollary}} {\sf (B$'$)}. {\sf The rational filling radius of $S$ in $X$ is  bounded from below by 
$$fil.rad(S, X; \mathbb Q)\geq \frac {1}{2}width_\square (S).$$} 

This means that {\it no multiple of the curve  $S$ bounds in the $\rho$-neighbourhood} $U_{\rho}(S)\subset X$
for $\rho<width_\square (S),$ or, more formally,
 the homology boundary homomorphism $H^2(U_{\rho}(S), S)\to H_1(S)=\mathbb Z$ {\it vanishes} for these 
(small) $\rho$.\vspace {1mm}

These (A) and (B) and (B$'$), albeit useful, are   boringly trivial, but the following one is mildly amusing.\vspace {1mm}

{\it \textbf {Lemma}} {\sf (C)}. {\sf Let $X$ be a  locally contractible  {\it path metric} space  \footnote {The distances between pairs of  points are  equal to the infima of  lengths of curves  between them.} 
and let  $\gamma(x)=dist (x,x_0)$ be the distance  function on $X$  to  some point  $x_0\in X$.}

{\it If  all embedded circles $S\subset X$ have 
$$\max.width_\square\leq D,$$
then the diameters of all connected components of the levels $\gamma^{-1}(t)\subset X$, $t\geq 0$ are bounded by $3D$. }

{\it Proof.} Let $\gamma(x)=dist (x,x_0)$ for some point  $x_0\in X$ and let us show that the diameters of all connected components of the levels $\gamma^{-1}(t)\subset X$, $t\geq 0$ are bounded by $3D$.
Indeed, assume without loss of generality that $t\geq \frac {3}{2}D$ and 
 let $x_1,x_2$ be two points in a {\it connected component} of $\gamma^{-1}(t)$.
   
    Let $\widetilde {x_1,x_2}$ be a segment joining $x_1$ with $x_2$ in a small neighbourhood of  this component  and let $[x_1,x_0]$   and $ [x_2,x_0]$ be almost shortest  segments between $x_0$ and $x_i$, $i=1,2$, where we assume without loss of generality that 
    the union $S$ of these three segments 
    $$S=\widetilde {x_1,x_2}\cup [x_1,x_0]\cup  [x_2,x_0]\subset X$$
     makes  a topological circle.
    
  Let   $[x_i,x'_i]\subset [x_i,x_0]$  be subsegments of length $D+\varepsilon$
  and let 
  $$\widetilde{x'_1,x'_2}= [x'_1,x_0] \cup [x'_2,x_0] \subset [x_1,x_0]\cup [x_2,x_0]$$ 
   be the union of the complementary segments.  
  Now, 
  $$dist (\widetilde {x_1,x_2}, \widetilde{x'_1,x'_2})>D$$
  while   
  $$dist([x_1,x'_1], [x_2,x'_2])\geq dist(x_1,x_2)-2D-\varepsilon',$$
  which, due to our assumption $max.width_\square(S)\leq D$, 
   implies for $\varepsilon, \varepsilon'\to 0$ that $dist(x_1,x_2)\leq 3D.$
  QED.
  \vspace {1mm}
  
   {\it \textbf {Corollary }} {\sf (C$'$)}. {\it The inequality 
$$\max.width_\square\leq D$$
 for   all embedded circles $S\subset X$  implies that 
   the Uryson  1-width of $X$ is bounded by:
$$width _1(X)\leq 3D.$$} 
 \vspace {1mm}
 
 {\it Proof.} Factor the map $\gamma:X\to \mathbb R$ as 
$$ X\overset {\alpha} \to P_0^1 \overset {\beta}\to \mathbb R,$$
such that  the levels of $\alpha$  are equal to the connected components of $\gamma$,
  and approximate (this is trivial)  the (1-dimensional!) space $P^1_0$ by a polyhedral $P^1$ as it is done for this purpose in  the proof of corollary  10.11 in [GL(complete) 1983].
 
   \vspace {1mm}

  {\it \textbf {Corollary }} {\sf (C$''$)}. {\sf If  the first cohomology of $X$ vanishes then the above graph $P^1$ is a {\it tree.}}
  
  Consequently, {\it $X$ is quasiisometric to a tree} {\sf in this case.}
   \vspace {1mm}
   
  { \it  Proof.}  Slightly modify $\alpha$ to  make the levels of $\alpha$ path connected  and  observe that   continuous onto  maps   with path connected fibers are surjective on the fundamental  groups. 
   
   \vspace {1mm}

   {\it \textbf {Corollary}} {\sf (D)} {\sf  If all closed curves immersed to $X$ bounds in their $\rho$-neighbourhoods,}  then the {\it  the Uryson  1-width of $X$ is bounded by:
$$width _1(X)\leq 6\rho.$$}

  {\it \textbf {Corollary}} {\sf (E)}. {\it If $X$ has infinite 1-width, then it contains closed curves $S$ with arbitrarily large maximal $\square$-widths.}\vspace {1mm}
  
    {\it \textbf {Example}} {\sf (E$'$)}. {\sf Universal coverings $\tilde X$ of compact spaces $X$ with {\it non-virtually free} fundamental groups  $\pi_1(X)$ contains circles with arbitrarily large maximal $\square$-widths.}\footnote 
  {Compare with Corollary   in the section 1.2 C in [G(foliated) 1991].}

  {\it \textbf {Lemma}} {\sf (F)}. {\sf Let $X=(X,g)$ be  a complete Riemannian $3$-manifold  with a boundary  and let   $X^\rtimes =   X\rtimes \mathbb T^m= (X\times \mathbb T^m, g^\rtimes=g+\phi^2dt^2$
  be a  $\mathbb T^\rtimes$ extension  of $X$  with {\it mean convex} boundary and  such that  $Sc( g^\rtimes)\geq \sigma>0$.}  
  
 {\it Then all   immersed  circles  $S\subset X$,  {\sf homologous to  zero (e.g. contractible ones)} have their 
 rational  filling radii bounded by 
 $$fil.rad(S, X'\mathbb Q)<  \frac {2\pi}{\sqrt \sigma}. $$}

  {\it Proof.}  Since $S$ is homologous to zero, it bounds an orientable surface in $X$;  let  $Z\subset X$ be such a surface with boundary $\partial =S$, which  minimizes the {\it $g^\rtimes$-area} of $Z$
  that is the $(m+2)$-volume of the hypersurface  $Y\times \mathbb T^m\subset X^\rtimes$, that is  
 equal to the   area of $Y$  with respect to the {\it  conformal metric}  $\psi(x)\cdot g(x)$ on $X$,
where $\psi(x)=vol (\{x\} \times \mathbb T^m)=(2\pi)^m \phi^{m}(x)$. 
  
  Since the minimal  hypersurface   $Z\times \mathbb T^m\subset X^\rtimes$ is stable it admits a  $\mathbb  T^1$ extension  
  with the scalar curvature bounded from below by that of $ X^\rtimes$ and the proof follows from  codimension 2 corollary  to $\frac {2\pi}{n}$-inequality.  (See \textbf {\color {blue}  [2] }  in section \ref{codimensions 1, 2, 3.3},   where the surface is denoted $\underline X$ rather than $Z$.) 
 \vspace {1mm}
   
   This,   together with the above  corollary {\sf (D)} yields the following. \vspace {1mm}
   
    {\it \textbf {Proposition}} {\sf  (F$'$).}  {\sf Let $X=(X,g)$ be  a complete Riemannian $3$-manifold  with a boundary, such that the  homology group of $X$ is torsion.
  
    If $X$ admits a $\mathbb T^\rtimes$-extension    $X^\rtimes =   X\rtimes \mathbb T^m= (X\times \mathbb T^m, g^\rtimes=g+\phi^2dt^2)$
    with {\it mean convex} boundary and with  $Sc( g^\rtimes)\geq \sigma>0$,  then the {\it first Uryson width of 
    $X$ is bounded by  
    $$width_1(X)\leq \frac {12\pi}{\sqrt \sigma}.\footnote{Our present proof of this inequality  follows  that  of  corollary 10.11
    in [GL(complete) 1983] where it is stated in the non $\mathbb T^\rtimes$-stable form and where the 
    $H_1$-torsion condition, albeit implicitly used in the argument,  was erroneously omitted.}$$}}

   \vspace {1mm}
     
    {\it  \textbf {3D Classification Corollary}} {\sf (F$''$).}  {\sf  A compact  $3$-manifold  admits a metric with $Sc(X)>0$ {\it only if }   it contains {\it no aspherical  connected summand} in its Kneser-Milnor  prime decomposition.}
     ("If" is also true by Perelman's  theorem.)
      \vspace {1mm}
   
    {\it Proof}. In view of    {\sf (C$''$)},  the universal covering of $X$ is quasiisometric to a tree,  and then an application of 
     Stallings' theorem  shows that the fundamental group is virtually free. QED,.     
      \vspace {1mm} 
   
   {\it  On Domination.}  The   proof  of this classification theorem in [Gl(complete), 1983], which used
  the index theorem, albeit less  straightforward,  yields  more:     \vspace {1mm}

  {\sf If a closed  orientable 3-manifold  $X$ contains an {\it aspherical manifold  in its prime decomposition}, then it 
    {\it can't be dominated} by  a complete manifold $\hat X$ {\it with $Sc>0$}: 
   
   {\sf all  maps between such orientable  manifolds, $\hat X \to X$, have  {\it degrees zero}.}}

       \vspace {1mm}

     {\it \textbf {Theorem }} {\sf  (G).}   {\it If a complete orientable Riemannian $3$-manifold  $X=(X,g)$  with a boundary  admits a $\mathbb T^\rtimes$-extension    $X^\rtimes =   X\rtimes \mathbb T^m= (X\times \mathbb T^m, g^\rtimes=g+\phi^2dt^2)$
    with {\it mean convex} boundary and with  $Sc( g^\rtimes)\geq \sigma>0$,  then the {\it first Uryson width of 
    $X$ is bounded by  
    $$width_1(X)< \frac {36\pi}{\sqrt \sigma}.$$}}
  
  {\it Proof.}  Decompose $X$ into a union of  submanifolds with common boundaries.
  $$X=\bigcup_i X_i,$$
   where two such submanifolds intersect by a stable
  {\it $g^\rtimes$-minimal} surface (as in lemma {\sf (F)}) and such that all  connected stable  $g^\rtimes$-minimal surfaces in all $X_i$ are homologous to the connected components of the boundaries of these $X_i$.

We  know that all these surfaces are  spherical and there diameters are bounded by  
$\frac {12\pi}{\sqrt \sigma}$.

It is also clear that the rational homology groups $H_1(X_i,\mathbb Q)$ vanish and the proof follows from 
lemma {\sf (A)} and proposition {\sf  (F$'$)}.

 {\it \textbf {Corollary }} {\sf  (G$'$).}  {\sf If $X$ is compact  with  no boundary, then its 
  {\it absolute  filling radius} is bounded by 
 $$fillrad(X)< \frac {18\pi}{\sqrt \sigma}.$$}
This means that 

{\it there exists an orientable $4$-dimensional  pseudomanifolds $V$ with boundary and a metric $dist_V$ on $V$,  such that the boundary  $\partial V$  is {\sf isometric to} $V$  and 
$$dist_V(v, \partial V)<\frac {18\pi}{\sqrt \sigma}.$$  }  \vspace {1mm}
 
 {\it Proof.} If a piecewise linear map\footnote {This is understood with respect to some smooth triangulation of $X$ map. Notice  that our  $\chi$ in the definition of width was assumed {\it continuous} rather  than piecewise , 
 linear, but it can be approximated by piecewise linear ones.}    $\chi: X\to P^1$ has  $diam(\chi^{-1}(p)\leq d$, $p\in P^1$,
 then { it the cone $F=C_\chi$ of $\chi$} is a pseudomanifold,  which carries an obvious metric  which has the required properties.

  {\it \textbf { Remark  }} {\sf  (G$''$).} It follows from {\color {blue}$[width_2]$} in  Remark {\sf (A$'$)}
  that  $$fillrad(X)< \frac {6\pi}{\sqrt \sigma},$$
 
 {\it \large Probably}, 
  $$fillrad(X)< \frac {\pi}{2}\sqrt\frac {6} {\sigma}.$$
 
\vspace {1mm}
 
\textbf {Enlargeability in Dimension  3}. Let us conclude this section by proving that \vspace {1mm}
 
 {\it the universal coverings of compact  aspherical 3-manifolds $X$  are 
 enlargeable.}\vspace {1mm}
   
Notice  that   "enlargeability" (defined below) is, a priori, stronger, than non-existence of  a metric with $Sc>0$;
this remains conjectural for higher dimensional aspherical manifolds.
 
Also notice  that in dimension 3,  one  can easily prove this property  for manifolds of each of the  7 non-elliptic Thurston's  geometries separately, and then (this is also easy) show that   enlargeability is stable under JSJ-decomposition. 
 
  Our point   here   is to furnish a  direct elementary proof. \footnote  {The concepts of enlargeability and related "bad  ends"   are discussed   in   [GL(complete) 1983]  around theorem  8.1, in  [Lawson\&Michelsohn(spin geometry) 1989]  around theorem IV.6.18] and also in    [G(positive) 1996] in  \S\S  $9\frac{1}{4}$,      $9\frac {3}{11}$, where similar properties are proved for "multiple largeness";  later this   appears in    [G(inequalities)  2018], section 4,  under the name of {\it iso-enlargeability}.}
  \vspace {1mm}    
    
\textbf {Definitions of the Hyperspherical Radius, Hypersphericity, Enlargeability and Uniform  Lipschitz Asphericity.}  The {\it {\color {red!20 !blue} hyperspherical radius $Rad_{S^n}(X)$}}, of a 
 {\it closed} orientable  Riemannian $n$-manifold $X$ as 
{\sl the supremum of  the radii  $R>0$ of $n$-spheres, such that $X$ admits a {\it non-contractible 1-Lipschitz}, i.e. distance non-increasing, map $f: X\to S^n(R)$.}\vspace {1mm}

More generally, if $X$ is an {\it open} manifold,  this definition still make sense for  maps $f:X\to S^n$,  which are {\it locally constant  at infinity},\footnote {On can drop 'locally" if $X$ is connected at infinity.} i.e. outside  compact subsets in $X$.   Similarly, if $X$ allowed a boundary, then $f$ should be constant on all components  of this boundary. 

Notice that  a (locally constant at infinity if $X$ is open)  map $f$ from an {\it orientable} $n$-manifold   $X$ to the sphere  $S^n$  is {\it contractible} (in the space of  locally constant at infinity maps in the open case)  if and only if $f$ has {\it zero degree}.

In view of that, we define $Rad^{deg1}_{S^n}(X)\leq Rad_{S^n}(X) $ as  the supremum of  $R$ for  1-Lipschitz, maps $f: X\to S^n(R)$  of degrees 1.

A manifold $X$ is called {\it  hyperspherical} if $Rad_{S^n}(X)=\infty$ and $X$ is {\it enlargeable} if it admits coverings with arbitrary large hyperspherical radii.

A manifold $X$ is called   {\it $deg1$-hyperspherical} if, $Rad^{deg1}_{S^n}(X)\leq Rad_{S^n}(X) =\infty,$ or,
in different terms if $X$ 
$\lambda$ -Lipschitz dominates 
 (the fundamental homology class of)  the unit sphere $ S^n$ for all $\lambda>0$.

 \vspace {1mm}

A metric space $S$ is called {\it  uniformly Lipschitz $k$-aspherical} if  $\lambda$-Lipschitz maps from the unit  sphere $S^k$ to $X$,  are extendable to $\Lambda (\lambda)$-Lipschitz  maps from the unit ball $B^{k+1} $  that bounds $S^K$, i.e. $\partial B^{k+1}=S^k$,  for all $d>0$, where   
  $\Lambda(d)=\Lambda_X(d)$ is a continuous  (control)  function. 
  
  \vspace {1mm}

 \textbf  {Example}. If  $X$  has bounded geometry,\footnote {Probably, a lower bound on the sectional curvature suffices.}   e.g. $X$ is   a  covering of a compact  locally contractible space,   and if $X$ is $k$-aspherical, i.e.  $ \pi_k(X)=0$, then it is    uniformly 
  Lipschitz $k$-aspherical. 
  
 {\it Exercise .}  Construct complete uniformly contractible surfaces,  which are not Lipschitz  uniformly $1$-aspherical.
 
 Also construct  complete uniformly contractible Riemannian 3-manifolds where the sectional is asymptotically
 non-positive:
 $\kappa(X,x)\leq \varepsilon(dist(x,x_0))$,  where $\varepsilon(d)$ is  a positive function   which goes to $0$ for $d\to\infty$.

 \vspace{1mm}

  {\it \textbf {Lemma }} {\sf  (I).} {\sf Let $X$ be a complete orientable
  Riemannian $n$-manifold, and  let     $  Y_i\subset  X$  be  smooth  closed connected  orientable codimension 2  submanifolds with trivial normal bundles, (e.g. $H_{n-2}(X)=0$)  and let $U_i=U_{\rho_i} \supset Y$ be the $\rho_i$-neighbourhoods of $Y_i$ in  $X$  where 
  $\rho_i\to\infty$ for   $i\to\infty$.
  
 Let $X$ be  $2$-aspherical  and Lipschitz  uniformly $1$-aspherical,   let  $Y_i$  admit   
 $deg1$-hyperspherical coverings $\tilde Y_i$.}

 {\it If the  inclusion homomorphisms $\pi_1(Y_i)
 \to\pi_1(U_{\rho_i})$ are injective and if   the manifolds    $Y_i$  are not rationally  homologous to zero in $U_{\rho_i}$,  then $X$ is   
  $deg1$-hyperspherical.} \vspace{1mm}
  
  {\it Proof.} Let  $ U'_i\subset U_i$  be a small tubular neighbourhood of $Y$, where, observe the boundary $\partial U'_i$ topologically
 splits as $Y\times S^1$. 
  
  Since $X$ is 2-connected and the class $[Y]\in  H_{n-2}(Y)$ goes to a non-zero class in 
  $H_{n-2}(U_i;\mathbb Q)$, the linking number between closed curves in the complement   $U_1\setminus Y$  with $Y$ defines a homomorphism from $\pi_1(U_1\setminus Y)$ to  $\mathbb Z$, which {\it doesn't vanish on the class
  of the circles} $\{y\}\times S^1\subset \partial U'_i\subset U_i\setminus Y$, and  {\it vanish on 
  $\pi_1(Y)=\pi_1(Y\times\{s\}$.}
  
  It follows  that the inclusion homomorphism $\pi_1(\partial U'_i)\to \pi_1(U_1\setminus Y$ is {\it  injective}. 
  
  Now, let $\widetilde {U_i\setminus Y}$ be the covering of $U_1\setminus Y$, the restriction of which to 
  $\partial U'_i=Y\times S^1$ equals to  $ \tilde S^1 \times \tilde Y_i$ for $\tilde S^1=\mathbb R^1$ and the above 
  $deg1$-hyperspherical covering  $\tilde Y_i$ of $Y_1$ and 
  show, this is easy, \footnote{See \S\S 5, 6 in [GL(complete) 1983] and   IV. 6 in [Lawson\&Michelsohn(spin geometry) 1989].}  that 
  
  {\sf since $Rad^{deg1}_{S^{n-2}}(\tilde Y_i)=\infty$, the $deg1$ hyperspherical  radius of  $\widetilde {U_i\setminus Y}$  bounded from below only by $\rho_i$, say as follows,
$$Rad^{deg1}_{S^n} (\widetilde {U_i\setminus Y_i}) \geq \frac{1}{10} \rho_i.$$}

  Finally, since $X$ is Lipschitz  uniformly $1$-aspherical the covering map
  $$\widetilde {U_i\setminus Y_i} \to  \partial (U_i\setminus Y_i\subset X$$
   is one-to-one  "deeply
    inside" $\widetilde {U_i\setminus Y_i}$, i.e.   sufficiently far from $Y_i$  and $\partial (U_i\setminus Y_i)$, ("how deeply"or har "far"  depends on the (control) function $\Lambda(\lambda)$) and the proof follows.

   \vspace {1mm}

   {\it \textbf {Corollary}} {\sf  (I$'$).} {\it Compact orientable aspherical $3$-manifolds $X$ are enlargeable.}
   \vspace {1mm}
  
  Indeed, as we know,  the universal coverings $\tilde X$  of these $X$ contain closed curves  with arbitrarily large 
  maximal $\square$-widths;  these  can be taken for the above $Y_i$.

   \textbf {Theorem} {\sf {(J)}.} {\it If a compact orientable 3-manifold  $X$ contains an aspherical  summand in its prime decomposition then no non-zero multiple of (the fundamental class of) $X$ can be dominated by a 
  complete orientable manifold $\hat X$ with $Sc(\hat X).0$.}
  
  In simple words, 
  
  {\it all continuous maps $\hat X\to X$ constant at infinity have zero degrees.}   \vspace{1mm}

If $X$ is enlargeable and    $\hat X$ dominates $X$,   then  $\hat X$ is also enlargeable, and  one knows  that enlargeable {\it spin}  $n$-manifolds
support no metrics with $Sc>0$  for all $n$ (theorem 6.12 in [GL(complete) 1983]).  Since all 3-manifolds are spin, the proof follows.

%%%%%%%%%%%%%%%%%%%%%%%%%

\subsubsection {\color {blue}Geometry and Topology of Complete  3-Manifolds  with $Sc>0$}  \label{complete  3D.3}
 
%%%%%%%%%%%%%%%%%%%%%%%%%%

 Start with a simple proof of the following  result  by Laurent Bess\`eeres, G\'erard Besson, and Sylvain Maillot    [Be-Be-Ma(Ricci flow) 2011]. \vspace{1mm}

 \textbf  {  (A) Theorem.} {\it  Complete 3-manifolds $X$  with $Sc(X) \geq \sigma>0$ are infinite connected  sums of 
spacial space forms $S^3/\Gamma$ and copies of $S^2\times S^1$.} \vspace {1mm}

{\it Proof.}  By the  compact exhaustion corollary  from section  \ref{codimensions 1, 2, 3.3},  $X$ decomposes  into infinite  connected sum  of {\it  compact} manifolds $X_i$, and   theorem  {\sf (J)}   from the previous section implies that there is no 
aspherical  summands in these $X_i$.  Then the conclusion follows by Perelman's  theorem.

\vspace{1mm}

{\it Remarks. } (i) The argument in  [Be-Be-Ma(Ricci flow) 2011] depends  on a generalization of Perelman's arguments to non compact manifolds.  Also,  G\'erard Besson recently told me  that     Jian Wang found a proof    of \textbf  {  (A) }  with minimal surfaces.  \vspace {1mm}

(ii) A close look at  proof of \textbf  {  (A)}  shows that  $X$ decomposes into a union of compact submanifolds $X_i\subset X$,   such that

$\bullet$  $X_i$ intersect with $X_j$, for all $i\neq j$,  over the common  components of their   boundaries;

$\bullet$ the boundaries of $X_i$ are union of spheres the areas and the  intrinsic diameters of which are  bounded by a constant depending only on $\sigma$;

$\bullet$ the diameters of all $X_i$ are also  bounded by a constant depending only on $\sigma$.

 \vspace {1mm}

\textbf  {(A$'$) Corollary.} {\it No non-torsion  homology  class $\underline h\in H_3(\underline X)$   in an aspherical space can be dominated by a complete 
3-manifold $X$  with $Sc(X)\geq \sigma>0$.} \vspace {1mm}

{ \it Proof.} Given a map $X\to \underline X$, homotop it to $f'$  constant on representatives of all 
non-contractible 2-spheres in $X$ and thus reduce the problem to the, case where $X$ is a single spherical space form

Alternatively, argue algebraically and use the fact that all   finitely generated subgroups in  $\pi_1(X)$   are virtually free. 
   \vspace {1mm}

\textbf {(B)} {\it  \textbf  {Generalization to}} $\mathbf {Sc>0}.$ {\it If a $3$-manifold $X$ admits a complete metric with 
$Sc>0$ then    all   finitely generated subgroups in  $\pi_1(X)$   are virtually free.}
     \vspace {1mm}

{\it Proof.}  Let $\tilde X\to X $  be a covering with non-virtually free fundamental group and let 
  $\bar X \subset \tilde X $ by the compact {\it Scott core}  of $\tilde X$. Then, by the loop theorem,
  the boundary of $\bar X \subset \tilde X $ is incompressible and the proof follows from 
  theorem  6.12  in [GL(complete) 1983].

Alternatively, one can prove that  $X$ contains a (compact or complete)  {\it stable minimal surface}, which is {\it non-simply connected}, while  one knows   (see the proof of Wang's  theorem below) that   such surfaces don't  exist  in  3-manifolds with $Sc>0$.

\vspace{1mm}

 {\it \textbf  { Remarks,Examples and Open Problems for $\mathbf{Sc>0}$}}.  
 \textbf{(a)} 
 The
     apparent {\it irreducible}, i.e. non-trivially indecomposable into connected sum, 
 example of an 
 open  manifold, which admits a  complete metric  with $Sc>0$ is
 $\mathbb R^2\times S^1$  with the (radial  warped product) metric 
$$\mbox {$g =dr^2 +\varphi(r)^2 d\theta^2$, $t\in [0,\infty) , \theta\in [0, 2\pi] $,  
with $\varphi(r) =  r^{2\alpha}$,}$$ 
 where the scalar curvature  
 $$Sc(g)(r)=-\frac {2\varphi''(r)}{\varphi(r)}=\alpha(\alpha-1)\frac {1}{r^2} $$
 is positive for $1<\alpha<2$ with quadratic decay for $r\to\infty$. 
  
 By the above, $\mathbb R^2\times S^1$  admits no  metric with $Sc\geq \sigma>0$; moreover, the 
 the curvature must decay at least as  $\frac{4\pi^2}{r^2}$  
according to  QD-exercise in section \ref {quadratic3}.

  \textbf{(a$'$)} If $n\geq 4$,  there are  similar  complete warped metrics with $Sc>0$ on $\mathbb R^2\times X^{n-2}$  for all (compact and open) manifolds $X^{n-2}$.

  \textbf{(b)} 
     It is {\it {\color {red!40!black} unknown}} (unless I am missing something obvious)
  if {\sf open  handle bodies of all genera, hence,  the interiors  all  compact 3-manifold $\bar X$ with boundaries, 
 which have a virtually  free  fundamental groups  $\pi_1(\bar X)$, admit  complete metrics with $Sc>0$.}

 \textbf{(b$'$)} If $X$ is an $n$-manifold for $n\geq 4$,  (maybe one should assume $n\geq 5$), which {\sf contracts to its 
 codimension 2 skeleton}, e.g. a  contractible one, then, {\it {\color {red!40!black} conjecturally,}} it  
 {\it admits a complete  metric   with $Sc>0$.}
  
   However, no such metrics are  known, for instance,  in the interiors of compact manifolds the boundaries of which admit metrics with   negative sectional curvatures $<0$.

\vspace {1mm}

The following result by Jian Wang shows  that obstructions to the  complete metrics with $Sc>0$ on $X$ may resides in  the complexity of  the  {\it proper homotopy type} of $X$.

 \vspace {1mm}

\textbf {(C)}  \textbf {Theorem.} {\it Complete {\sf contractible} $3$-manifolds with $Sc>0$  are simply connected at infinity} (see   [Wang(Contractible)  2019]  and  [Wang(topological characterization)  2021] in this volume).\vspace {1mm}

{\it Idea of the Proof.} Recall that  the first   contractible 3-manifold $X=X_{Wh}$ {\sl not simply connected at infinity},  which was discovered by Whitehead in 1935, is  equal to  the  union of an infinite  increasing sequence  of solid tori,  
$$X_{Wh}=\bigcup_k T_i,\mbox { } T_1\subset T_2\subset ... \subset T_k\subset.... \subset  X_{Wh},$$
where the {\it boundary} of  $T_k$, $k\geq 2$, is  {\it  not contractible in} $T_2$  for all $k\geq 2$.

Wang shows in this case that, 
given an arbitrary  complete Riemannian  metric on $X_{Wh}$,    there exist connected  stable minimal surfaces $\Sigma_k\subset T_k$  of genus zero   with boundaries $\partial \Sigma_k\subset T_k$, such that the number of connected  of the intersections
$\Sigma_k\cap T_1$ goes to infinity for $k\to \infty$.

Then, in the limit, he obtains a {\it connected stable} minimal surface  $\Sigma=\Sigma_\infty\subset X_{Wh}$ of genus zero
with a {\it complete} induced metric,  such that the intersection $\Sigma_k\cap T_1$ has {\it  infinite} area; 
this, in the case of $Sc(X>0)$,  contradicts to the {\it  Fischer-Colbrie{\sf\&}Schoen}
(Gauss-Bonnet-Cohn-Vossen) {\it inequality}
$$\int_\Sigma Sc(X,\sigma) d\sigma\leq 2\pi\chi(\Sigma).$$

 {\sf Hence \vspace{1mm}
 
\hspace {5mm}  {\it the   Whitehead manifold} admits no complete metrics with $Sc>0$.}\vspace {1mm}

 \textbf {(C$'$) \it  \textbf {Possible  Generalizations.}}  Wang's  argument  applies  (as far as  I understand it) 
 to connected sums $X=X_{Wh}\#X_1\#X_2...$ with other 3-manifolds and shows that these$X$ admit no complete metric with $Sc>0$.
 
{\it \color {red!40!black}   Conceivably,} Wang's argument can be  also applied to manifolds $X$, which {\it dominate}
  {\sf   the  fundamental homology class  $[X_{Wh}]$ (with infinite support). } 
 
 \vspace{1mm}
 
 If so, then by the (non-compact $\mathbb T^\rtimes$-stabilized  version of the) Schoen-Yau  inductive descent argument , the products $X_{Wh}\times \mathbb T^m$ (and probably, the products of $X_{Wh}$ with enlargeable manifolds in general) admit no complete  metrics with $Sc>0$ either.
(If  $m> 5$, one has to appeal to  Lohkamp's desingularization theorem.)  \vspace{1mm}

In fact, contractibility of  manifolds in Wang's theorem doesn't seem  that essential.

 {\it \color {red!30!black} Conjecturally},  {\sf if an orientable $3$-manifold can be exhausted by  compact submanifold  $V_1\subset V_2\subset... \subset  V_i \subset...   \subset X$, such that all components of the  complements $X\setminus V_i$ are aspherical with infinitely generated  fundamental groups and the inclusion homomorphisms 
 $\pi_1(X\setminus V_{i+1})\to  \pi_1(X\setminus V_{i})$ for all $i=1,2,....$, are injective (maybe, its enough to assume that the images of the inclusion  homomorphisms  $\pi_1(X\setminus V_{i+1}) \to\pi_1(X\setminus V_{1})$ are infinitely generated)}, then 
 $X$ {\sl admits no  complete metric with $Sc(X)>0$}.

Moreover, 
 
 {\sl  no non-zero multiple of the fundamental  homology class  $[X_+]$
   can be dominated, by a complete manifold with $Sc>0$,
  that is,  no  complete orientable 3-manifold $\hat X$ with $Sc(\hat X)>0$  admits a proper map to a  $\underline X_+$  with non-zero  degree.}

 \vspace{1mm}

 {\it Example.}  Let a  connected  orientable manifold  $X$ decompose  
 into a countable  union of compact aspherical submanifolds with  aspherical boundaries, 
 $X=\cup_iX_i$, such that
 
$\bullet$ {\sf every two $X_i$  intersect (if at all) over  several connected components of their boundaries,  where these intersections are denoted $Y_{ij}=X_i\cap X_j=\partial X_i \cap \partial X_j$.

 $\bullet$   the inclusion homomorphisms $\pi_1(Y_{ij})\to\pi_1(X_i)$ are  injective and their   images have  infinite indices in the fundamental groups  $\pi_1(X_i)$.}
 e.g.  as  in  the Whitehead manifold,  where $X_i$ are (the closures of  $T_{i+1}\setminus T_i$.

Then  the above conjecture implies that for $n=3$  {\it no   manifold $X_+$, which contains $X$ as a submanifold, admits  a  complete metric with $Sc>0$.}  \vspace {1mm}     \vspace{1mm}

{\it Remark about $n>3$}.   For all we know, the   $n$-manifolds  $X$  (minus the boundaries)  and $X_+\supset X$   in this 
example   don't admit     complete metric with $Sc>0$ enlargeable  for all $n$,
 but   this can be proved at the present moment  only in special cases, for instance, if    some   manifold 
 $Y_{ij}  \subset  \partial X_i$  is {\it enlargeable}, (  e.g. if $dim(X)=4$,  since compact aspherical 3-manifolds are enlargeable, see the previous section)   and  if  the inclusion homomorphism $\pi_1(Y_{ij})\to\pi_1(X_+)$
 is {\it  injective} (see  section \ref {obstructions4}).
   
 \vspace{1mm}

 \textbf {Attaching Cylinders to Stable  Hypersurfaces.} Let $X=(X,g)$ be a complete, e.g.  compact,  Riemnnian manifold with a boundary,

 Notice that completeness of $X$, i.e.  compactness of  closed  bounded 
subsets, implies completeness of the boundary with respect to the Riemannian distance function $dist_g$ in $X\supset Y.$ 

Let $Y\subset \partial X \subset X$ be a connected component of the boundary and let $Y\rtimes_\phi \mathbb R_+$
be the warped product with the metric $h_\phi=h+\phi^2dt^2$ for the Riemannian metric $h$ on $Y$ induced from $g$ on $X$ and a smooth positive  function $\phi=\phi(y)$.

 Observe that the boundary $Y\times 0\subset   Y\rtimes_\phi \mathbb R_+$ is isometric to $Y\subset X$ and let
 $$X_O=X\sqcup_YY\rtimes_\phi \mathbb R_+$$ 
 be obtained by  attaching  $Y\subset X$ to  
 $Y\times 0
 \subset   Y\rtimes_\phi \mathbb R_+$ by this isometry.

This $X_O$,which homeomorphic to the complement $X\setminus Y$  carries a  natural continuous Riemannian  metric which is complete   of $X$ and hence, $Y$, are complete.

 Now, if the  $Y\subset X$ is mean convex,  and if the scalar curvatures of both manifold are positive, 
   $Sc(g)>0$ and $Sc(h_\phi)>0$, then the metric on $X_O$ can be approximated  by smooth metrics with  $Sc>0$,  since  $Y\times 0\subset 
 Y\rtimes_\phi \mathbb R_+$ is {\it totally geodesic} in  $Y\rtimes_\phi \mathbb R_+$ (see section \ref {mean1}).
 This yields  the following.
  \vspace{1mm}
 
 \textbf {(D)}  \textbf {Proposition.} {\sf  Let  $X=(X, g)$ be a complete Riemannian manifold with $Sc>0$  and let 
 $Y\subset X$ be a cooriented  stable minimal hypersurface.
  The the complement $X\setminus Y$  admits a complete metric $G$ with $Sc(G)>0$, which is equal to $g$ outside a given neighbourhood $U\supset Y$ intersected with $X\setminus Y$.}
   \vspace{1mm}
 
  Let us apply this to 3-dimensional  manifolds $X$, where $Y$ is a topological 2-sphere  and where we benefit from the following homotopy theorem of  
 due to Laurent Bessi\`eres, G\'erard Besson, Sylvain Maillot, and Fernando Coda Marques.
  \vspace{1mm}
 
 \textbf {(D)}  \textbf {Theorem.}  {\it The  space of complete Riemannian metrics of bounded geometry and uniformly positive scalar curvature on an orientable 3-manifold is path-connected.}
 
 \vspace{1mm}
 
 It follows the the metric $G$ on $X\setminus Y$ can be homotoped outside a given compact subset to the standard cylindrical  metric $ds^2+dt^2$ on $S^2\times \mathbb R_+$, and then extended to the ball $B^3$  keeping the curvature positive all along. 
 
 The  one can attach the  unit  3-ball   to  the sphere $S^2\times \{t_0\}$ and smooth the resulting $C^1$ metric
 with $Sc>0$ .
 
  \vspace{1mm}

 {\it Remark.}  One  may use     the compact case of  \textbf {(D)},  namely  for $S^2 \times S^1$, where 
this was earlier proven in [Marques(deforming $Sc>0$)2012].

Besides, one doesn't   need  here the full power of the Ricci flow, since the relevant  deformation proceeds in  the space
of $S^1$-invariant metics, which  are moreover,   of  the forms $g+\phi^2dt^2$, and where the 3-D equations reduce to 2-dimensional ones for pairs $(g, \phi)$  of Riemannian metrics $g$ and functions  $\phi$ on  
$S^2$. 
 
 \vspace{1mm}
 
 \textbf {(E)}  \textbf  {Corollary}.  {\sf Let $X$ be a complete 3-manifold with $Sc(X) >0$.
 Then there  exists a  complete (disconnected) manifold $X^\sim$, such that} 
 
 $\bullet$ {\sl all connected  $X^\sim_i$ of $X^\sim$  are "simple":
 
   the complement  to a    embedded 2-sphere $S^2$  or to  a  properly (infinity $\to$ infinity)   embedded plane $\mathbb R^2$    in $X_i^\sim$, for all $i$,   is disconnected and at least one of the two components is homeomorphic to $S^3$ with finitely or countably  many punctures;

  $\bullet$ The complement  to a finite or countable set of disjoint complete stable  connected minimal surfaces  $\Sigma_{min}$ in in $X$  is 
  isometric to an open subset in $X^\sim$, where all these $\Sigma_{min}$ are simply connected, and if $Sc(X)\geq \sigma>0$ they are all compact, hence spherical. }

  \vspace {1mm}

{\it Proof.}  Let    $\Sigma$  be an   {\it "essential"} embedded 2-sphere or a  properly  embedded plane    in a complete 3-manifold $ X$, i.e. such that   the complement $X\setminus\Sigma$  is either connected or none of the two  components is homeomorphic to 
$S^3$ with punctures. Then either $X$ contains an essential stable minimal sphere  or an essential    stable  minimal plane.  \footnote{It is, certainly well known.  I apologize to the author for not being able to   find  his/her article. 

Notice, however, that this fact  is easy  in  our case, where $Sc(X)>0$, since all complete minimal surfaces necessarily are either spherical or planar in these $X$.}

 If  the scalar curvature of $X$ is uniformly positive, i.e. $Sc(X)\geq \sigma>0$,  then  all these minimal surfaces are spherical  and we  attach 3- balls to them as above. 
 
 In general, where $Sc(X)>0$, we attach 3-balls to the spherical  cutting surfaces and   cylinders to the planar ones.
 QED.

  \vspace {1mm}
  
  {\it Question}   {\it  Is  there  a version of the above  for,  say compact, $4$-manifolds with $Sc>0$?}

  \vspace {1mm}

{\it Remark.} The first things one needs is a "natural filling" of the spherical space forms $S^3/\Gamma$ by 4-manifolds (may be singular ones?)  with $Sc>0$,  something in the spirit  of discs bundles over $S^2$ that  fill in the diagonal lens
spaces $S^3/\mathbb Z_k$.
 %%%%%%%%%%%%%%%%%%%

\subsubsection{\color {blue}Non-Existence of  Uniformly Contractible  and  Aspherical  4- and  
 and 5-manifolds with $Sc>0$} \label  {5D.3}

 %%%%%%%%%%%%%%%%%%

Recall that a metric space $X$ is {\it uniformly contractible} if 
 there exists a function $R(r)\geq r$, Recall that a metric space $X$ is {\it uniformly contractible} if 
 there exists a function $R(r)\geq r$, called {\it contractibility control    function} 
 such that the  $r$-balls 
$B_{x}(r)\subset X$  around all $x\in X$,  are  contractible in the concentric balls $B_{x}(R(r))$.

For instance, if $X$ is bounded  then  "uniformly contractible"="contractible". 

Also obviously, but more interestingly, 
the same applies to  spaces $X$ that  with {\it cobounded}, (e.g.  compact)  isometry groups:
there is a constant $d$, such that, 

{\sf for every two points $x_1,x_2\in X$, there exists  isometry $I:X\to X$
such that

 $dist(x_1, I(x_2))\leq d$.}

In particular, 

{\it universal coverings of  compact aspherical manifolds are uniformly contractible}.

\vspace {1mm}

An essential property of these $X$ is a \vspace {1mm}

{\it bound on the filling radii of cycles  $Y\subset X$ in terms of the absolute filling radii

 of these cycles.}\vspace {1mm}

In fact, a  standard {\it  induction by skeletons extension} argument shows the following.

{\color {blue} $\bigstar$}  {\sf Let $W$ be a polyhedral space, $Y \subset W $  a polyhedral subspace and let $\varphi: W\to X$ be a  
continuous map.}

{\it If $X$ is uniformly contractible, then  $\varphi$ extends to a continuous map $\Phi:W\to X$, such that 
the distances from the points  $\Phi(w)$ to the image of $\varphi$ are bounded by
$$dist(w, \varphi(Y)\leq D(dist(v, Y)),$$
 where $D(d)$ is a continuous function that depends only on the  contractibility control    function 
 $R(r)$ of $X$.}  \footnote {Consult [G(filling) 1983],  [G(aspherical) 2020] for basics on  filling and uniform contractibility and see [Katz(systolic geometry)  2017],  [Guth (waist) 2014], [Wenger(filling) 2007],
  [DFW(flexible) 2003],
     [Dranishnikov(asymptotic) 2000]. 
       [Dranishnikov(macroscopic) 2010],
       [Dranishnikov (large scale) 1999].  
[Dra-Kee-Usp(Higson corona) 1998] 
 and section \ref{metric7} for related   topics.}

 \vspace {1mm}
The following immediate corollary to   {\color {blue} $\bigstar$}  will be used below for manifolds $X$ of dimension $n=4$.

\vspace {1mm}

 {\color {blue} $\bigstar_1$} {\it \textbf {Codimension 1 Filling Lemma}.  {\sf Let  $X$ be an $n$-dimensional orientable  pseudomanifold    a proper  metric
and let $U_1\subset U_2 \subset ...\subset U_i\subset ...\subset X$  be an exhaustion  of $X$ by compact 
sub-pseudomanifolds with boundaries.}} 

For instance, $X$ can be a {\it complete Riemannian manifold}  exhausted  by  compact domains with smooth boundaries.

{\it $X$ is uniformly contractible, then the absolute  filling radii of the boundaries of $U_i$ tend to infinity:
$$filrad(\partial U_i)\to\infty\mbox { for } i\to \infty. $$}

{\it Proof.} Let $S\subset X$ be an infinite path, i.e. a  curve, issuing at a  point $x_0\in U_1$ and tending to infinity.
Since   $S$ has  {\it  non-zero}  intersection indices with the boundaries of all $X_i$,  the   boundary 
$\partial U_i$ can bound in its $D$-neighbourhood only if  $D< dist (x_0, \partial U_i)$. Since 
$ dist (x_0, \partial U_i)\to \infty$  so does $D$ and the proof follows.

Now we recall that, according to  compact exhaustion corollary   ({ \color {blue}\textbf {[1]}}   section \ref{codimensions 1, 2, 3.3}),  that  complete Riemannian manifolds $X$ with $Sc(X)>\sigma>0$  it can be exhausted 
compact smooth domains $U_i$,  the boundaries  $Y_i=\partial U_i$ of which  admit $\mathbb T^{\rtimes 1}$-extension $Y_i\rtimes \mathbb T^1$ 
with  
$$Sc(Y_i\rtimes \mathbb T^1)\geq \frac {\sigma}{2}.$$

If $dim(X)=4$ and $dim(Y_i)=3$. then  all these $Y_i$ have their filling radii bounded by $$fillrad(X)< \frac {18\pi}{\sqrt {\sigma/2}}$$
by  corollary {\sf  (G$'$)} {from}  \ref {filling+hyperspherical+asphericity3}.{\color {red}(corrected!!!!!!!!!!)}%????%%%%%%%%%%%%%%%%%
 Hence, \vspace {1mm}

\textbf{[4D]A.}   {\it\color {blue!60!black}  complete uniformly contractible 4-manifolds $X$ can't have $Sc(X)\geq \sigma>0$.}
 \vspace {1mm}

This, applied to the  universal coverings  of compact manifolds yields \vspace {1mm}

\textbf{[4D]B.} \textbf {Chodosh-Li 4D Theorem}. {\it No compact aspherical 4-manifold admits a metric with  
$Sc>0$.}
\vspace {1mm}

\textbf{[4D]C.} {\it\color {blue} \textbf{ Generalization-Exercise}.} {\sf Let $\underline X$ be an orientable {\it uniformly rationally acyclic}   four dimensional pseudomanifold}, which means that 
  the  rational homology inclusion homomorphisms between the balls around all points $x\in X$,
$$H_i(B_x(r); \mathbb Q)\to H_i(B_x(R);\mathbb Q)$$
{\it vanish} for all $i=1,2,... $  and $R\geq R_\mathbb Q(r)$ for some (acyclicty control) function $R_\mathbb Q(r)$. 

{\sf Let $X$ be a complete orientable  Riemannian manifold and $f:X\to \underline X$ be a proper $1$-Lipschitz map  with {\it non-zero} degree.}

 Show that  there exist  constants  $C, R_0>0$, such that 
 
 {\it the minima of the scalar curvature  of $X$ on  concentric  balls  $B(R)=B_{x_0}(R)\subset X$  
 around a point $x_0\in X$, satisfy
$$  \min_{x\in B(R)} Sc(X,x)\leq  \frac {C}{R^2}  \mbox  {  for all } R\geq R_0. $$}
 
{\it Hint.} Adapt  the above proof to maps $X\to \underline X$ similarly to how this is done in    [G(aspherical) 2020]. 
 \vspace {1mm}

{\it Remark.} This implies that if $X$ is a compact orientable  4-manifold with $Sc(X>0)$, then  continuous 
maps from $X$ to {\it aspherical 4-dimensional pseudomanifolds} send the rational fundamental homology class of $X$ to zero.  But this remains unknown for maps from these $X$ to aspherical spaces in general.\vspace {1mm}

Let us  prove another corollary to   {\color {blue} $\bigstar$} needed which will be   used in dimension $n=5$.\vspace {1mm}

 {\color {blue} $\bigstar_2$} {\it \textbf {Codimension 2 Filling Lemma}. {\sf Let  $X$ be an $n$-dimensional orientable  pseudomanifold where the {\it singular locus of} $X$, (where it is {\it not} a manifold) has {\it codimension} 3, i.e.  
 the links of  the codimension  2 faces are connected  and let $X$ be endowed with a proper  path metric with respect to which $X$ is
 {\it uniformly contractible.}
 
  {\sl Then, for all $R>0$, there exits a  proper piecewise linear  1-Lipschitz map 
 $\Psi=\Psi_R: X\to \mathbb R^2$, such that}

{(\Large$\star$}) {\it all orientable codimension 2-sub-pseudomanifolds  $Y$, which are  contained in the 
 the pullback $\Psi^{-1}(B(R))$  of the $R$-ball  $B(R)=B_{\mathbf 0}(R)\subset  \mathbb R^2$  and  which are }
 
  \vspace {-0.1mm}
 
 \vspace {0.6mm}
 
 \hspace {16mm} {\sf \color {black!50!blue} homologous in $\Psi^{-1}(B(R))$ to  
  the pullbacks  
  
  \hspace {16mm} $ \Psi^{-1}(t)\subset X$,    of regular points 
  $t\in B(R)$
  of  $\Psi$,} \footnote{These points are dense in $\mathbb R^2$ and their pullbcks $ \Psi^{-1}(t)\subset X$  are  compact  sub-pseudomanifolds in $X$, which, for $t\in B(R)$, are all mutually homologous in  $\Psi^{-1}(B(R))$, i.e.  represent the same class in  the group  $H_{n-2}(\Psi^{-1}(B(R)) )$.}}   
   \vspace {0.6mm}

 \hspace {6mm} {\it  have their {\color {black!50!blue} absolute  filling radii bounded from below} by
 $$fillrad(Y)\geq r(R),$$
 where  $r(R)$ is a continuous function} (\sf which depends on the contractibility control function of $X$}),
 {\it such that 
 $$r(R)\to \infty \mbox { for } R\to \infty.$$}
 
 {\it Proof.}  The uniform contractibility of pseudomanifolds $X$ implies that their  Uryson 1-width  are infinite 
 $$ width_1(X)=\infty;$$
 otherwise, the hypersurfaces in $X$  would have bounded width, hence their filling radii as well in contradiction with {\color {blue} $\bigstar_1$.}

Next, by lemma {\sf (C)} from the previous section, there exists closed  curves $S\subset X$ with arbitrary large maximal $\square$-widths $D$, and let $\Psi_\square : X\to \mathbb R^2$ be the corresponding maps
delivered by lemma {\sf (C)},which, recall,  is $\sqrt 2$-Lipschitz and which sends $S$ onto the boundary
of the square $[0,D]^2\subset \mathbb R^2$ with degree 1. 

 Scale this map by $\frac {1}{\sqrt 2}$, shift  it to move   the center of the square to the 
 origin $\mathbf 0 \in \mathbb R^2$  and  take the resulting map for $\Psi$.
 
Since $X$ is contractible, the curves $S=S_D$ bound  orientable  surfaces $\Sigma=\Sigma_D$ 
  which have non-zero intersection indices with $Y$. Therefore, if $r$ is much smaller then $D$, yet 
  goes to infinity along with $D$, then the filling radii of $Y\subset \Psi^{-1}(B)$  also tend to infty.
Q.E.D.

Now,  if  $X$ is  a Riemannian $n$-manifold with $Sc(X)\geq \sigma>0$, then,   
 by codimension 2  corollary   { \color {blue}\textbf {[2]}} from section \ref{codimensions 1, 2, 3.3}, 
 it contains submanifolds $Y\subset \Psi^{-1} $ which admit $\mathbb  T^\rtimes $ extensions  $Y^\rtimes =Y\rtimes \mathbb T^2$ with  $Sc(Y^\rtimes)\geq \frac {\sigma}{2}$, which for 
  $n=5$ and $dim(Y)=3$, have    filling radii uniformly bounded  by  {\sf  (G$'$)} from the previous section; hence,\vspace {1mm}

\textbf{[5D]A.}  {\it\color {blue!60!black}  complete uniformly contractible 5-manifolds $X$ can't have $Sc(X)\geq \sigma>0$.}

 \vspace {1mm}
 
 Accordingly, one has the following
 
\textbf{[5D]B.}  \textbf{$5D$-Non-asphericity Theorem.} {\it No compact aspherical 5-manifold admits a metric with  
$Sc>0$.}

  \vspace {1mm}
  
  {\it Remarks, Generalizations, Problems.}  (a)  Albeit    these [5D]A\&B  {\it imply}  {[4D]A\&B (for $X^4\leadsto X^5=X^4\times \mathbb R^1$)
  the mapping versions   $X\to \underline X$  of them, {[4D]C} and  {[5D]C},\footnote{The statement and the proof of this is left to the reader.}   are {\it formally independent}, due to the  codimension 3 condition on singularities of   $\underline X$  for $dim(\underline X)=dim(X)=5$  that is needed 
   for  a homological   definition
  of the linking numbers between  curves  $S\subset X$ and codimension two   sub-pseudomanifolds $Y\subset X$.\vspace {1mm}
 
(b) This kind of "dual linking"  appears  in section 1.2 of [G(foliated) 1991] for  the purpose of "trapping" minimal foliations and also in 
\S$9\frac{3}{11}$ [G(positive) 1996], where $Y$ is  a circle,  for the proof  of enlargeability  and the Novikov conjecture for 3-manifolds.

  In the present context,    Chodosh and Li use it  for their proof of non-asphericity  of 4-manifolds with $Sc>0$ (enlargeability and the Novikov  conjecture remain  problematic for $n\geq 4$), where  $Y$ is a surface and  the bound 
on $filrad(Y)$ (in  terms which I don't quite understand) was derived in the first version of [Chodosh-Li(bubbles) 2020]  from the area bound due to  Zhu. 

Then Chodosh-Li's  linking idea,  combined with the $\mathbb T^\rtimes$-stabilized   bound on widths  of
 3-manifolds,  was applied to $n=5$ in  [Chodosh-Li(bubbles) 2020] (I didn't quite follow  how this is  done in their paper)  and in [G(aspherical) 2020],  where it is proved that   \vspace {1mm}
 
 {\it complete uniformly rationally acyclic (e.g. the universal coverings of compact aspherical) Riemannian manifolds $\underline X$ of dimension $5$ can't be  1-Lipschitz dominated\footnote {See section \ref {domination1} for the definition.}   with degrees $\neq 0$ by Riemannian manifolds $X$ with 
 $Sc(X)\geq \sigma>0$.}

  \vspace {1mm}
 
(c)  To extend the linking argument to $n=6$  one needs, as it is explained in [G(aspherical) 2020], either a  proof of  a {\it universal  bound on the  filling radii of 4-manifolds with $Sc\geq \sigma>0$} (this remains conjectural), or the existence of
 closed {\it surfaces $\Sigma$ (instead of curves $S$) in uniformly contractible manifold $X$, with  filling radii $fillrad(\Sigma,X)\geq \rho$, for all $\rho>0$},  i.e. non-homologous to zero in their
  $\rho$-neighbourhoods in $X$. (This remains problematic even for the universal covers of compact manifolds.) 
  
  At the present moment, one  has only limited results for $n\geq 6$  available  along these lines, e.g. 
  
(d)   {\sf non existence of metrics with $Sc>0$  on closed aspherical manifolds $X$ of dimension $n\geq 5$, the} {\it fundamental groups of which contain  subgroups isomorphic to $\mathbb Z^{n-4}$}, see section \ref{m-radii of uniformly contractible7}.
 \vspace {1mm}
 
 (e)  A closer look at the above argument  shows the following. 
 
  Let  $\underline X$  be an orientable  $n$-pseudomanifold   with a proper (bounded subsets are compact) metric.  
 {\sf If $\underline X$ is uniformly contractible (uniformly rationally acyclic will do) and if either $n=4$ or if $n=5$ 
 and the singularity of  $\underline X$ has codimension 3 (or more)}.  
 
 {\it  Then, if  complete Riemannian $n$-manifold $X$ with $Sc(X)\geq \sigma>0$  admits a  proper $1$-Lipschitz map $f: X\to \underline X$ then
 $\inf _{x\in X}Sc(X,x)\leq 0.$}
 
 Moreover,   
 
  {\it the scalar curvature of $X$, assuming it is positive, can't  decay subquadratically, or even slow quadratically:
  one can't have
  $$Sc(X,x)\geq \frac {C}{dist (x, x_0)^2}$$
  $$\mbox {  for a fixed point $x_0\in X$, all $x\in X$ with $dist (x,x_0)\geq 1$   and a positive constant $C=C(\underline X)$.}$$ 
  
\vspace {1mm}
 
 \textbf {Corollary.}  {\it Compact aspherical manifolds of dimensions 4 and 5 with punctures  admit no complete metrics with 
 $Sc\geq \sigma>0$.}

 {\it \textbf {Question}} {\sl} Are there  complete metrics with $Sc>0$ on these punctured  manifolds}?

\vspace {2mm}

%%%{\color {red} (\textbf{[5E]} below is added on June 27 2021!)!!!!!!}\vspace {1mm}
%%%%%%%%%%%%%%%%%%%%%%%%\%%%%
%%% added June 27

\textbf{[5E].} { \color{blue!50!black} \textbf {Classification of Non-aspherical 4- and 5-Dimensional Manifolds} 
\textbf{with $Sc>0$.}}
The classification theorem for 
3-manifolds with positive scalar curvatures was generalized in [Chodosh-Li-Liokumovich
(classification) 2021] as follows.\vspace {1mm}

  {\Large {\color {blue} $\bullet$} }  {\sf Let $X$ be a closed  connected  Riemannian   $n$-manifold with  {\it infinite} fundamental group $\pi_1(X)$
  and such that the higher homotopy groups $\pi_2(X), $..., $\pi_{n-2}(X)$ {\it vanish}.}

{\it If $X$ admits a metric with $Sc>0$}, {\sf then, assuming $n=4$ or  $n=5$,
{\sf a finite covering   of $X$ is homotopy equivalent to the}} {\it connected sum of several copies of $S^{n-1}\times S^1$.}

\vspace {1mm}

Let us   extend the above    non-asphericity arguments   to  classification and prove  {\Large {\color {blue} $\bullet$} }  by observing   the 
following. \vspace {1mm}

1.    The (obvious induction by skeleta)   proof of the above {\color {blue} $\bigstar$}  actually shows 
that 

 {\sf  if  a  compact polyhedral space, e.g. a compact manifold  $X$ has {\it trivial} homotopy groups $\pi_2(X),...,\pi_k(X)$,  then the filling radii $R$ of  all  $m$-dimensional 
 submanifolds  (and subpseudomanifolds, if you wish)  of dimensions $m\leq k$ in the universal covering of  $X$, say $Y\subset \tilde X$, are bounded in terms of their absolute filling radii $r$, 
 $$ R\leq D(r)\mbox{  for   $R= filrad(Y\subset X)$ and  $r= filrad (Y)$}$$   
and where $D=D_X=D_{X,k} $
 is the iterated  contractibility control function for  $l$-dimensional subpolyhedra    in $\tilde X$, for $l=2,3,  k$.}

2. Notice  that,  

\hspace {-1mm}  {\it { \color {red!40!black} unless} the fundamental group of a compact manifold $X$ is virtually free,}\vspace {0.6mm}

\hspace {-6mm} the conclusion of  the  codimension 2 filling lemma {\color {blue} $\bigstar_2$} holds for  $\tilde X$:
 that is, 
\vspace {1mm}

  {\sf {\color {black!50!blue}  there   exists of  a  proper piecewise linear  1-Lipschitz map 
 $\Psi=\Psi_R: \tilde X\to \mathbb R^2$, such that 
   all orientable codimension 2-sub-pseudomanifolds  $Y$, which are  contained in the 
 the pullback $\Psi^{-1}(B(R))$  of the $R$-ball  $B(R)=B_{\mathbf 0}(R)\subset  \mathbb R^2$  and  
 which are 
\color {black!50!blue} homologous in $\Psi^{-1}(B(R))$ to  
  the pullbacks  
   $ \Psi^{-1}(t)\subset X$,    of regular points 
  $t\in B(R)$
  of  $\Psi$,   
have their {\color {black!50!blue} absolute  filling radii bounded from below} by
 $$fillrad(Y)\geq r(R), \mbox { } r(R)\to \infty \mbox { for } R\to \infty.$$}}
   
In fact, according to example  {\sf (E$'$)} in  \ref {filling+hyperspherical+asphericity3}.\vspace {1mm}

 {\sf  universal coverings  of compact spaces $X$ with {\it non-virtually free} fundamental 
 
  groups  $\pi_1(X)$ do  contain circles with {\it arbitrarily large maximal $\square$-widths} }

\hspace {-6mm} and the presence of   such "large circles"  was all we used in the proof of {\color {blue} $\bigstar_2$}. (In truth, we also needed these circles to be  {\it homologous to zero}, which is automatic for $\tilde X$ is simply connected.)
\vspace{1mm}

Now, if $Sc(X)\geq \sigma>0$, then, as earlier,  the  "large" codimension 2 cycle $Y\subset \tilde X$ delivered by the codimension 2 filling lemma can be represented by a submanifold $Y'$  which is  positioned close to $Y$ and such that  a  $\mathbb T^\rtimes
 $-stabilization of it  has $Sc\geq \frac {\sigma}{2}$.
 
 Since, as we know,\footnote {If $n=4$ this follows from the  $\mathbb T^\ast$-stable Bonnet-Myers diameter inequality proved  in section   \ref {warped stabilization and Sc-normalization2} and if $n=5$ this is stated in corollary (G$'$) in section  \ref {filling+hyperspherical+asphericity3}.}  the inequality $Sc(Y'\rtimes \mathbb T^2)\geq \sigma/2$ implies that 
the filling radius of $Y'$  is bounded by  $filrad (Y')\leq const/\sqrt \sigma$, 
 it follows by contradiction  that  \vspace {1mm}
 
{\color {teal}({\Large $\star$})} \hspace {15mm}{\it the fundamental group $\pi_1(X)$ {\color {red!50!black}is} virtually free. } \vspace {1mm}

 We   conclude the proof with the
  following elementary topological lemma.

 \vspace {1mm}

{\color {teal}({\Large $\ast$})} {\sf If a closed orientable  manifold $\hat X$ of dimension $n$ has  {\it free fundamental group}  and  {\it zero
$\pi_2(X),...,\pi_k(X)$,} $k=n-2$, then it is homotopy equivalent to a connected sum of copies of $S^{n-1} \times S^1$.} QED.

 \vspace {1mm}

{\it Exercises.} (a$_1$) Show that the conclusion of {\color {teal}({\Large $\ast$})}   remains valid for $n<2k$,
e.g. for $k=3$ and $n=6$.

(If $n=5$, then these manifolds are actually  {\it diffeomorphic} to connected sums of $S^4\times S^1$, see  [Gadgil-Seshadri(isotropic)2008], [Kreck-Su](5-manifolds) 2017]  and references therein.)

(a$_2$)  Show that  {\color {teal}({\Large $\ast$})} also holds for orientable {\it pseudomanifolds} $\hat X$,
 the singular loci  of which have codimensions 3.
 
 (a$_3$) Formulate and prove a counterpart of  (a$_1$) for pseudomanifolds with  singular loci  of  codimensions $l$.
 
 %[Kreck-Su](5-manifolds) 2017] Matthias Kreck, Yang Su, {\sl  On 5-manifolds with free fundamental group and simple boundary links in $S^5$,}  arXiv:1602.01943.  vspace {2mm}

%[Gadgil-Seshadri(isotropic)2008] S. Gadgil, H. Seshadri, {\sl On the topology of manifolds with positive isotropic curvature}, arXiv:0801.2221 \vspace {2mm}

(b) Generalize {\Large \color{blue}$\bullet$} by replacing  "{\it $X$ admits a metric with $Sc>0$}"
    with {\it "$X$ can be dominated by a complete Riemannin manifold with $Sc\geq  \sigma>0$"}.
    
    Generalize further  to manifolds  $X$,  {\it the universal coverings} $\tilde X$  of which admit such  {\it Lipschitz dominations}, i.e.  such that 
    
     {\sf there exist  complete orientable  Riemannin manifolds $X'$  with $Sc(X')\geq \sigma>0$ and  quasi-proper (e.g. proper) 1-Lipschitz maps 
$X'\to \tilde X$ with non-zero degrees.}

Then do the same for {\it pseudomanifolds} $X$, the   singular loci  of  which have codimension $3$.

%%%%%%%%%%%%%%%%

\subsection {\color {blue}Asymptotic Geometry with  $Sc>0$,  Positive Mass Theorem and  Penrose Inequality}
\label {asymptotic3}
%%%%%%%%%%%%%%%%%%%%%
Let us      show that  complete Riemannian $n$-manifold $X$ with $Sc(X)\geq 0$   can't grow faster  the Euclidean space $\mathbb R^n$,  which is  regarded as the cone over the unit sphere $S^{n-1}$ with  
 the Euclidean metric represented (in polar coordinates)  is   $g_{\mathbb R^n}=dr^2+r^2ds^2$.

\textbf {1.}  {\it \textbf  {Conical Example}. } {\sf Let $Y$ be a  Riemannian  manifold of dimension $(n-1)\geq 2$ with 
$Sc(Y)\leq (n-1)(n-2)=Sc(S^{n-1})$ and let $g$ be the Riemannian metric on $X=Y\times [0, \infty)$ asymptotic  to a conical,
one, namely
$$g=g(y,r) = dr^2+ \lambda r^2 dy^2 + g_o(y,r),$$
where $g_o(y,r)=o(1)$, or, in words, $ g_o(y,r)$  converges to $0$  for $r\to \infty$.  

 This,  in  the scale invariant terms, means  that   the differential  $dI$ of the   identity map 
$$ I: (X, g)\to (X, dr^2+ \lambda r^2 dy^2)$$
{\it converges to isometry} for $t\to \infty$,
$$||dI(y,r)I| \to \mbox {  and also }   ||(dI)^{-1}(y,r)I| \to 1,$$
  that is $I$ is    {\it $\lambda(r)$-bi-Lipschitz for  $\lambda(r)\to1$.}}

\hspace{15mm} {\it  Granted the above, if  $Sc(X)\geq 0$,  then $\lambda\leq 1$.}  \vspace{1mm} 

{\it Proof.}  The condition $g_o(y,r)=0(r^2)$ is equivalent to the $C^0$-convergence of the $\varepsilon$ scaled metric $g$ to the background conical  (scale invariant)   metric  
$$\varepsilon^2  g(y,r)\to dr^2+ \lambda r^2 dy^2\mbox {  for }\varepsilon\to 0.$$
Hence, $Sc(dr^2+ \lambda r^2 dy^2)\geq 0$ by the $C^0$-closure theorem (section \ref{C0-limits3}).
But if 
$\lambda  >1$, then, for $dim(X)\geq 3$,  the  conical metric  $dr^2+ \lambda r^2 dy^2$ on $X$ 
has $Sc(dr^2+ \lambda r^2 dy^2)<0.$ QED.

 \vspace{1mm}

\textbf {2.}  {\it \textbf  {Asymptotically Schwarzschild}. } Recall (see section\ref {conformal2}) that  { \it the scalar curvature  (of the space slice of the) Schwarzschild metric with mass $m$,} 
$$ g _{Sw_m}=g_{Sw}=\left (1+\frac {2m}{r}\right)^4g_{Eucl}.$$
is zero and that:

{\sf if $m>0$, the metric $g_{Sw_m}$  is defined on $\mathbb R^3$  minus zero, and  it is complete,

if $m=0$, this is the flat Euclidean metric; 

if $m<0$, then  this metric is defined only for $r>m$ with a singularity ar $r=m$.}

For all $m$, the  metric $g_{Sw_m}$ is asymptotically  Euclidean  (conical), where,  $g_{Sw_m}$ grows slightly slower then the Euclidean metric  one and if $m<0$ it growth slightly faster.

This seen by rewriting this metric as 
 $$g_{Sw_m}=\left(1-\frac {2m}{r}\right)^{-1}dr^2+r^2 ds^2$$
and  computing the  mean curvature  of the $r$-spheres  with respect to the Schwarzschild metric (see [Brewin(ADM) 2006]),
$$mean.curv(S_{Sw_m})=\frac {2}{r} \left(1-\frac {2m}{r}\right)^{\frac {1}{2}}=\frac {2}{r}-\frac {2m}{r^2}+O\left(\frac {1}{r^3}\right),$$
where $\frac {2}{r}$ is equal to the Euclidean  mean curvature of the sphere $S^2(r)$.

Observe that the difference between  the Euclidean and  Schwarzschild metrics  and  their  first  derivatives  in the Euclidean coordinates satisfy
$$ g_{Sw_m }-g_{Eucl}=
\frac{2m}{r}dr^2-O\left(\frac {1}{r^2}\right)$$
and 
$$\partial^1(g_{Sw_m }-g_{Eucl})=\partial^1g_{Sw_m }=\frac{\partial g_{Sw_m }}{dr} =\frac{2m}{r^2}dr^2+O\left(\frac{1}{r^3}\right)$$

Now   we are going  to estimate the scalar curvature of a metric $g= g(s,r)$ that  is asymptotic to  $g_{Sw_m}$,  
where we start with the following.  \vspace{1mm} 

{\it Observation.} Let a  $g_0$  be a Riemannian metric in a neighbourhood $U_0\subset \mathbb R^n$ of the origin $0\in \mathbb R^n$,  let $S\subset U_0$ 
be a smooth hypersurface passing through the origin and let  $g$ be another smooth Riemannian metric 
in  $U_0$, which is $\varepsilon$-close to $g_0$ with its first Euclidean derivatives  at the origin,
$$ ||g_{0 }(0)-g(0)|| +||\partial^1(g_{0}(0) -g(0))||  \leq  \varepsilon\leq 1.$$
 
  Then 
  
  {\it the mean curvatures of $S$ at the origin with respect to  the two  metrics satisfy:
$$|mean.curv_g(S,0)-mean.curv_{g_0}(S,0)|\leq C\varepsilon,$$
where the constant $C>0$ depends only  on the $g_0$ and its fist derivatives at the origin,
i.e. 
$$C\leq  C_n (1+||g_{0 }(0)-g_{Eucl}|| +||\partial^1(g_{0}(0)||)$$
for a universal constant $C_n$.}

{\it Proof.} Check it for metrics $g_0=a_0^2dx_i^2+ b_0^2dy^2$ and $g=a^2dx^2+ b^2dy^2$
for smooth positive functions $a_0, b_0,a, b$ in the $(x,y)$-plane  and for the parabola $S=\{y=cx^2\}$; then reduce the general case to this  special  one.

 \vspace{1mm} 

\textbf {3.} \textbf {Positive Mass Corollary.} {\sf Let $g_0=g_0(s,r)$  be a warped product metric on the cylinder 
$S^{n-1}\times [2m,\infty)$
written as 
$$g_0= (1-\alpha(r)) dr^2 +r^2ds^2,$$
 where $0<\alpha(r) <m$ is a smooth  function  with 
 $\sup_r\frac {d\alpha(r)}{dr}<\infty$. 
 Let $g_i$ be a sequence of  smooth  Riemannian  metrics defined in neighbourhoods 
 $$V_i\subset S^{n-1}\times [2m,\infty)$$ of  $r_i$-spheres 
 $$S_i=S^{n-1}(r_i)=S^{n-1}\times \{r_i\}\subset 
 S^{n-1}\times [2m,\infty),$$ 
 where $r_i\to\infty$  with $i\to\infty$, such that the differences between 
 $g_i$  and $g_0$ and  their first derivatives measured with respect to the Euclidean metric 
 $dr^2+r^2ds^2$ 
  are asymptotically bounded as follows,
 $$\frac {||g_i(s, r_i)-g_0(r_i,s)||}{\alpha(r_i)}\underset{i\to\infty} \to 0\mbox {  and } 
 \frac {||\partial^1g_i(s, r_i)-\partial^1g_0(r_i,s)||}{\alpha (r_i)} \to 0.$$}

{\it If a complete orientable  (possibly disconnected) Riemannian   spin $n$-manifold $X$  
 contains closed  smooth embedded hypersurfaces $ \Sigma_i$ which admits isometries 
 $f_i\to \Sigma_i \to S_i$,  which  preserves their mean curvatures, 
  $$mean.curv(\Sigma_i,x)=mean.curv(S_i,f_i(x),$$
then $X$ can't have non-negative scalar curvature
$$\inf_{x\in X} Sc(X,x)<0.$$}

{\it Proof.} Observe that the   mean curvatures  of $S_i$ with respect to $g_0$ are
$$mean. curv_{g_0}(S_i)= (1-\alpha(r))^{-\frac {1}{2}}.$$
Then, according to the  above observation, 
mean curvatures of $\Sigma_i$ are strictly bounded away  from below by those of the Euclidean $r_i$-spheres,
$$(mean. curv_{g_i}(\Sigma_i))> \frac {n-2}{r_i^2}$$
and the proof follows from the mean curvature  spin-extremality theorem ({\color {blue} \EllipseShadow})
in section \ref{mean convex3}) applied to spheres.  \vspace{1mm}

{\it Example.}   If $\alpha(r)=\left(1-\frac {2m}{r}\right)^{-1}$ the above applies to  to  complete manifolds 
with asymptotically   Schwarzschild's metrics  $g$ and shows that positivity of the scalar curvature of $g$ makes $m\geq 0$. \vspace{1mm}

{\it \color {red!50!black}  Worrisome Remark.} Our assumptions on  the asymptotics of $g$ are  suspiciously
weak  compared to  the commonly used in the literature, e.g. by Schoen and Yau  in
[SY(positive mass) 1979].  This  makes  me  wonder if I haven't make a silly  mistake in my interpretation 
of  "derivatives of metrics". \footnote {Only while preparing these notes, I attempted  to penetrate  the meaning of these   mysterious "derivatives",  the geometry of which still remains above my understanding. Probably, if these have any meaning,   it should reside with physics (which I don't know) rather  than with geometry.} 

\vspace{1mm}

\vspace{1mm}

\textbf {4.} {\it $C^0$- Variation.} Let the  above neighbourhoods $V_i$   be the   annuli (bands) between the spheres of radii $r_i$ and
 $c_ir_i$,
$$V_i=S^{n-1}\times [r_i,c_ir_i]\subset S^{n-1}\times [2m,\infty), \mbox { } c_i>1. $$
If 
 the difference 
 $||g_i(s, r)-g_0(r)||$ divided by the width of such a band becomes   {\it sufficiently  small} on 
$V_i$  
$$\frac {||g_i(s, r_i)-g_0(r_i,s)||}{(c_i-1)r_i}\leq \varepsilon _i,$$
then,  regardless  of any bound on the derivatives of $g_i$, such a  band $(V_i ,g_i)$   contains a $\mu$-bubble $S_{i,min} \subset V_i$  for a (density)  function $\mu_i(s,r)$, which is close to the $g_0$-mean curvature of the spheres 
$S^{n-1}(r)$
$$|\mu_i(s,r)-mean.curv_{g_0}(S^{n-1}(r)|\leq \delta_i$$

For instance if $g_0$ is the Schwarzschild metric with mass $m<0$,  and if
 $\varepsilon_i$ is very small, e.g. $o(r_i^{-4})$, $r_i\to \infty$, then  the  argument, as in  the proof of the  approximation corollary
in \S$5\frac {5}{6}$  from[G(positive1996],\footnote {This  argument 
was  motivated by trying to geometrically understand  Min-Oo's hyperbolic positive mass theorem.} based on  
 Llalrull's inequality,  shows that  $\inf Sc(g)\leq 0$, and it is not hard to show that  the mean curvature extremality theorem allows the 
  same conclusion with $\varepsilon_i=o(r_i^{-1})$.

  \vspace{1mm}

{\it About  Rigidity.} Our argument, unlike these   by Schoen-Yau and by  Witten,  is poorly adapted to the case, where $g$ is asymptotically close, even when it very close, to  the Euclidean  metric, 
that is Schwarzschild  with the mass $m=0$. Apparently, rigidity of this kind is hard to  derived from the geometry of finite object without passing to the  infinite limit at an earlier stage of the argument.

\vspace{1mm}

{\it Remark/ Questions.} {\sf What happens  to  twisted harmonic spinors (best seen in  Lott's rendition of the mean curvature  spin-extremality theorem) that lie at the bottom of our argument in the limit for $r\to \infty$?}

They don't seem to converge in an  obvious way to Witten's spinors, but do they?

\vspace{1mm}

{\sf Does the positive mass rigidity holds for $C^0$-perturbations of the Euclidean metric?}
\vspace{1mm}

Although no  available techniques is capable to prove this even for very fast decay of $||g-g_0||$, we formulate  the  following. \vspace{1mm}

{\it  Euclidean  $C^0$-Rigidity {\color {red!50!black} Conjecture}.}  {\sf If  a smooth Riemannian metric  $g$ on $\mathbb R^n$

(a)  satisfies 
$$||g(x)-g_{Eucl}(x)||= o\left (\frac {1}{||x||}\right),  \mbox { }   x\to \infty,$$
 or

 (b) if the identity map 
 $$ (\mathbb R^n, g)\to (\mathbb R^n, g_{Eucl})$$
is $\lambda$-bi-Lipschitz for  some $\lambda<\infty$ 
and the difference of the two distance functions  
$$ dist_g(x_1,x_2)-dist_{g_{Eucl}}(x_1,x_2)$$
is  bounded, on $\mathbb R^n\times \mathbb R^n$,  then either 
 $$\inf_{x\in \mathbb R^n} Sc(g(x))<0,$$
or $g$  is Riemannian flat.}

(If $g$  satisfies (a) and  is everywhere $C^0$-close to the Euclidean $g_{Eucl}$,  one may try the 
Hamilton-Ricci flow.)

{\it Admission.} It is unclear,  not  even conjecturally,  how close these sufficient rigidity  conditions (a)  and (b)  are close to necessary ones..

\vspace {1mm}

\textbf {5.} {\it On History and Recent Developments.} The following special case of   the positive mass conjecture (unsolved by that time) was emphasized  by Robert Geroch  
 in  his expository article [Geroch(relativity)  1975]  for geometers.

{\sf  The   Euclidean metric on $\mathbb R^n$ admits no compactly supported perturbations  with increase of the scalar curvature.} 

Moreover,

 {\it If a metric $g$ on $\mathbb R^n$ with $Sc(g)\geq 0$ is equal to  $g_{Eucl}$ outside a compact subset in $\mathbb R^n$, 
 then  $(\mathbb R^n, g)$ is isometric to  $(\mathbb R^n,g_{Eucl})$.}

\vspace {1mm} 

This, of course,   "trivially" follows from non-existence  of non-flat metrics  with $Sc\geq 0$ on tori,  since compactly supported  perturbations of the flat metric on $\mathbb R^n$ yield similar perturbations of flat tori.
\footnote{The reduction to tori is  {\it amazingly}  simple,   where this  "amazing" brings it far  from    "trivial".}

Originally  Schoen and Yau directly  proved a stronger {\it   positive mass/energy theorem},  that   claims positivity of the  {\it ADM-mass},\footnote{ In their paper [SY(positive mass) 1979] the authors refer to some  earlier results, e.g. to   Jang, P.S.: J. Math. Phys. 1, 141 (1976), but its hard to say what's in there since it is not   openlc available on line.} which means that the 
 
 {\it total (i.e. integral)  mean curvature of the Euclidean spheres $S^2(R)$  with respect to $g$, 
  is bounded,  for large $R\to \infty$, by $8\pi R$.}
 \footnote {This interpretation of  the ADM-mass is explained in [Brewin(ADM) 2006], where the autor referres to Brown and York  for the origin of this idea.}

Two  years later,  Schoen and Yau extended their argument, based on  
non-compact minimal surfaces,  to  manifolds of dimensios $n\leq 7$,
while  Witten  suggested a proof applicable to  spin manifolds of all dimensions.
 
 Witten's    argument, that uses    perturbations   of invariant (non-twisted) harmonic spinors on $\mathbb R^n$, was worked out in details by Bartnik and  it  was   adapted by Min-Oo  to hyperbolic spaces.

\vspace{1mm}

Later,  Lohkamp   found  a (relatively)  simple reduction  of the general,  and technically more challenging, case of the positive mass theorem to that of 
 compactly supported  perturbations, thus reducing the problem  to   $Sc \ngtr 0$ on tori.

Most recently, the positive mass theorem was extended to a class of incomplete manifolds.   (See [Lesourd-Unger-Yau(arbitrary ends) 2021], where there are references to the earlier work by these authors.)
 \footnote {We don't even  attempt 
to    convey the basics of  physics and mathematics behind   the positive mass/energy idea, with dozens(hundreds?)  papers  dedicated to it,  besides  the early ones we mentioned: [SY(positive mass) 1979],  [Witten(Positive Energy) 1981],  [Bartnik(asymptotically flat) 1986],  [Min-Oo(hyperbolic) 1989],  
 [Lohkamp(hammocks) 1999];     we refer    to the  survey   [Herzlich(mass) 2021] and to   {\sl Positive energy theorem}  in Wikipedia for an    overview of this subject matter.}

\vspace{1mm}

{\it \color {blue} Problems.} {\sf  What are other (homogeneous?) Riemannian spaces that admit no (somehow) localised  deformations with increase of the scalar curvatures?

What are most general asymptotic conditions on such deformations that would allow their localization?}\vspace {1mm}
 
\textbf {6. Penrose Inequality.}
 Recall that \vspace {1mm}
 the  Schwarzschild metric with mass $m>0$, \vspace {1mm}
$$ g _{Sw_m}=\left (1+\frac {2m}{r}\right)^4g_{Eucl},$$ \vspace {1mm}
defined in the 3-space minus the origin,
is invariant under the the (conformal) reflection of $\mathbb R^3$ around the sphere $S^2(\rho)\subset \mathbb R^3$
 of radius $\rho=\frac {m}{2}$, that is 
  $$(s,r)\mapsto\left  (s,\frac  {\rho^2}{r}\right).$$
 This show that  he  Schwarzschild metric is complete and that 
 the sphere $S^2(\rho)$ is totally geodesic in geometry of  $g_{Sw}$, 
 with area 
 $$area_{g_{Sw_m}}(S^2(\rho)= \pi \rho^2 \left (1+\frac {\rho}{\rho}\right)^4= 16\pi m^2.$$

In 1973 Penrose formulated in [Penrose(naked singularities) 1973] a conjecture concerning black holes in general relativity with an  evidence in its favour, that  would, in particular  imply the following.  \vspace {1mm}

{\it \color {blue} Special  case of the  Riemannian   Penrose Inequality.} {\sf Let $X$ be complete Riemannian $3$-manifolds with compact boundary $Y=\partial X$, such that 

$\bullet$ $X$ is isometric at infinity to  the Schwarzschild space of mass $m$ at one of its two ends   at infinity;

$\bullet$  the scalar curvature of $X$ is everywhere non-negative: $Sc(X)\geq 0$;

$\bullet$ the boundary $Y$ of $X$ has zero mean curvature;\footnote {It suffices to assume  that the the boundary  is  {\it mean convex},  i.e.  its   mean curvature relative to the normal field pointing outward  is  positive.}

$\bullet$  no minimal surface in $X$ separates  a connected  component of $Y$ from infinity.}

{\it Then the area of   $Y=\partial X$ is bounded by the mass of the Schwarzschild space as follows.\footnote {This version of the  Penrose conjecture  is taken from the modern literature.
It is unclear, at least to the present author, when, where and by whom an influence of positivity of  scalar curvature in 3D on  geometry of surfaces, which was, probably, known to physicists  since the early 1970s  (1960s?)  was explicitly   formulated in mathematical terms for the first time.}
$$area (Y)\leq 16\pi m^2.$$}.

This, in a greater generality  was proven by Hubert Bray in   [Bray(Penrose inequality) 2009].

\vspace{1mm}

%%%%%%%%%%%%%%%%%%%%

\subsection { \color {blue}  Extensions and Fill-ins with $Sc>0$}\label{fill-ins3}
%%%%%%%%%%%%%%%%

The positive mass/energy  results from  in the previous two sections concerning   asymptotically flat and asymptotically hyperbolic  spaces, as well as sharp bounds on  the  size of mean convex hypersurfaces
from section \ref{mean convex3} are solutions of special cases of the following two general problems.

\textbf {A. {\color {blue!60!black} \it \textbf {Extension Problem for $Sc\geq \sigma$.}}} Let $X$ be a smooth manifold
with a boundary $Y=\partial X$, , let $h$ be a Riemannin metric on $Y$ and let $\sigma(x)$ and $\mu(y)$
be smooth functions on $X$ and on $Y$.

{\sl  What are necessary and what are sufficient  conditions for the existence of a complete (if $X$ is non-compact)  Riemannin metric $g$ on $X$, which extends $h$, 
 $$g_{|Y}=h,$$
 with respect to which the mean curvature of  $Y\subset X$ is equal to $\mu$,
$$mean.curv_g(Y)=\mu,$$
and such that 
$$Sc(X,x)\geq \sigma(x)?$$}

\textbf {B. {\color {blue!60!black} \it \textbf {Fill-in Problem for $Sc\geq \sigma$.}}}
Let $Y=(Y,h) $ be a Riemannian manifold and $\mu(y)$ be  a smooth function on $Y$.

{\sl Under what condition(s) does there exist, for a given number $\sigma $,  a  complete   Riemannian manifold  $X=(X,g)$  with $Sc(g)\geq \sigma$  with boundary $\partial X=Y$, such that 
$$g_{|Y}=h\mbox{  and } mean.curv_g(Y)=\mu,$$ 
and where, if $Y$ is compact, one may (or may not) require that $X$ is also compact?}
%%%%%%%%%%%%%%%%%
\subsubsection{\color {blue}  Construction of Extensions of Metrics  with $Sc>0$} \label {extensions3}
%%%%%%%%%%%%%%%%
Prior to enlisting known obstruction to  extensions and fill-ins  with $Sc>0$  in the next section, let is  describe    known  instances of existence of such extensions  and formulate several  questions. 

\vspace {1mm}

 {\it Remarks.}  (a) One could,   instead of the 
{\it Bartnik data} $(h, \mu)$ on $Y$,  prescribe a germ $g_\circ$ of a Riemannin metric on 
 an infinitesimal  neighbourhood of $Y$ in $X$, since, by the proof Miao's  gluing  lemma in section \ref {mean1},   
 a metric $g_0$  on $X$ with the same Bartnik data on $Y$ several as $g_\circ$ can  be deformed to $g$ with the same germ at $Y$ as $g_\circ$  without decease  of  the  scalar curvature.

(b) Following the general logic of the scalar curvature problems, one is concerned   not only with the  shear  existence of metrics   $g$, (manifolds $ X=(X,g)$ in the case  \textbf B),  but  with the space of all $g$ which have  given Bartnik data on $Y$ and $Sc(g)\geq \sigma$.
 
 Also, one may ask for a metric  $g$  with some its metric  invariant(s) (e.g. the hyperspherical radius)  bounded from below.

(c) If one drops $\mu$ from Bartnik Data $(h(y) \mu(y)$  then one expect no constraint for $Sc(g) $ on 
$X$ at all, where a recent definite result   in this regard, due to  Yuguang Shi, Wenlong Wang, Guodong Wei in    [SWW(total mean)   2020]  (responding to "embarrassing question"   from an earlier version of this manuscript) is as follows. \vspace {1mm}

\textbf { SWW Extension Theorem.} {\sl \color {blue!50!black}  All smooth Riemannin metrics $h$ on the boundary 
$Y=\partial X$  of a compact $n$-manifold $X$, extend to metrics $g$ on $X$ with $Sc(g)>0$.} \vspace {1mm}

The main technical ingredient of  the proof  is the following, \vspace {1mm}

{\it \textbf { SWW  Lemma.}} {\sf Let $h_0$ and $h_1$ be smooth Riemannin metrics on a compact Riemannin manifold $Y$,  and let $M_1$ be a  constant.
 
 If $h_1>h_0$, then then there exists a  smooth metric $g_{\circ}$ on  the cylinder $Y\times [0,1]$  with $Sc>0$, which} {\it extends}
   $h_{0}$ {\it and $\lambda \cdot h_{1}$, 
 $$\mbox {$g_\circ |_{Y\times \{0\}}=h_0$ and $g_\circ |_{Y\times \{1\}}=dot h_1$.}$$
 and such that } 
 
 {\it the mean curvature of the 1-end of the cylinder is bounded from below by $M_1$,
  $$mean.curv_{g_\circ}(Y\times\{1\}\subset Y\times [0,1])\geq M_1.$$}
 
 {\it Derivation of the theorem from the lemma.} By the  {\it $h$-principle for  open manifolds}, there exists a Riemannin metric $g_1$ on $X$ with $Sc(g_{1})>0$. 
 
 Let $g_{\circ}$ be the metric on the cylinder $Y\times [0,1]$ delivered  by the lemma for $h_{0}=h$ and $h_1=g_{|Y\times \{1\}}$, such that the $g_\circ$-mean curvature of the  1-end $Y\times \{1\}=Y$  of the the cylinder is greater than the minus $g_1$-curvature of the boundary $\partial X=Y$
  $$mean.curv_{g_\circ}(Y,y)= mean.curv_{g_1}(Y,y)$$
 Multiply the metric $g_1$ by $\lambda\geq 1$  from the lemma and isometrically attach the cylinder to $(X,\lambda\circ g_1)$.
By Miao's gluing lemma \ref{mean1}, the (continuous Riemannian) metric on 
$$X\sqcup_{Y\times \{1\}} Y\times [0,1]$$
can be approximated by a smooth metric $g$ on $X\sqcup_{Y\times \{1\}} Y\times [0,1]$ with $Sc(g)>0$  and, 
by an obvious identification $X=X\sqcup_{Y\times \{1\}} Y\times [0,1]$,
 the proof  of the theorem follows.  
 
 \vspace{1mm}  
 {\it On the Proof of the Lemma.}\footnote { I want to thank Yuguang Shi who explained to me several points in the proof of this lemma  and pointed out an error in my first  version of the proof of the "weak lemma".}   The  metric $g_{\circ}$   
 is constructed in  [SWW(total mean)   2020] in   the form
 $$g_\circ=g_u = (1-t)h_0+th_1 +u^2dt^2,$$
where the needed function  $u=u(y,t)$  is obtained as  a solution of a (non-linear parabolic) equation  expressing $Sc(g_u)$ in terms of the function  $u$ and its first an second derivatives.)
 \vspace {1mm}

  {\it $\sigma$-Remark.}   {\sf The metric $g_\circ$ delivered by the argument in [SWW(total mean)   2020]  can be  chosen  
 with {\it arbitrarily  large scalar curvature}}
 $$\mbox { $Sc(g_\circ)\geq \sigma$ for a given $\sigma>0.$}$$  
 
  \vspace {1mm} 
 
 {\it $\sigma$-\textbf{Corollary. }} {\it All  smooth Riemannin metrics $h$ on the boundary 
$Y=\partial X$  of a compact $n$-manifold $X$, extend to metrics $g$ on $X$ with $Sc(g)>\sigma $  for all $\sigma>0$.} \vspace {1mm}
 
 {\it Proof.}   By the   $h$-principle,  one gets    $g_1$ on $X$ with $Sc(g_{1})>\sigma$, where then the induced metric 
$h_1$ on $Y= \partial X$ can be made greater  than a given $h$ by the {\it (Nash)-Kuiper stretching construction.}
 \vspace {1mm}
 
 The following   elementary  proposition  also yields SWW theorem (albeit only with small  $\sigma$) via the     gluing 
 argument from [SWW(total mean)   2020]  . \vspace{1mm}

 \textbf {Weak SWW  Lemma.} {\sf There exists  positive constants $\delta_\nu>0$, for  $\nu>0$ and 
 a  family  of smooth positive monotone increasing functions $\lambda_\nu$  on the segments $[0,\delta_\nu]$
 
 $$\mbox { $\lambda_{\nu}(t)$,  $0\leq t\leq \delta_\nu$, $\nu\geq 0$, $\lambda_{\nu}(0)=1$}$$
  with the following property.}
  
  {\it Let $\underline h_{t}$ , $0\leq t\leq 1$, be a smooth family of Riemannian metrics on a compact manifold $Y$. 
 Then  the scalar curvature of the metric 
 $$g_\nu=\lambda^2 _{\nu} (t)\cdot \underline h_t+dt^2\mbox {  on } X_\nu=Y\times [0,\delta_\nu]$$ 
becomes arbitrarily large for large $\nu$,
 $$ Sc(g_\nu)\geq \sigma=\sigma(h_t, \nu)\to\infty \mbox {  for }  \nu\to \infty$$ 
 and also the mean curvature of the $\delta_\nu$-boundary of $X$ becomes  large
 $$ mean.curv_{g_\nu}(Y\times \{\delta_nu\})\geq  M=M(h_t, \nu)\to\infty \mbox {  for }  \nu\to \infty.$$}
 
 {\it Proof.} Recall the function   
  $$\varphi_\nu(t) =\exp \int_{-\pi/\nu}^t -\tan \frac {\nu t}{2} dt, \mbox { }  -\frac {\pi}{\nu}<t < \frac {\pi}{\nu},$$
from section \ref {warped2},  let 
$$\lambda^\circ_\nu(t) = \frac {\varphi_\nu(t)}{\varphi_\nu(t_0)}, \mbox  { } t\in [t_0, t_1],$$ 
 and 
$$ \lambda_\nu(t)= \lambda^\circ_\nu(t +t_0),
\mbox  {  } t\in [0, t_1-t_0].$$

Now an elementary computation\footnote {The formulas one needs, collected in sections  \ref{equidistant2},
\ref {Weyl2},\ref{Gauss2}, \ref{warped2} are: Riemannian  variation formula: $\frac{dh_t}{dt}=2A_t^\ast$,
Second Main Formula: $\frac{dA_t}{dt}=-A^2(Y_t)-B_t,$ and Gauss’ theorema egregium, while the 
relevant computation is sufficient to perform for the case of $h_t=h$, where $h$ is a flat metric (as in 
example (c) in \ref {warped}) and then argue by continuity. In fact,  a similar computation free argument
can be applied to the   metric with constant curvature 1 on $S^2$.}
 shows that if 
$$\mbox { $t_0=  -\frac {\pi}{\nu} +\frac {1}{\nu^3}$ 
and $t_1=-\frac {\pi}{\nu} +\frac {1}{\nu^2}$},$$ 
then the family $ \lambda_\nu(t) $, $t\in [0, \delta_\nu=t_1-t_0] $  is the required one.

\vspace {1mm} 
%%%%%%%%%%%%%%%
{\it  Remark.} The  lower bounds on the scalar curvature of $g$    and on  the  mean curvature of $Y\times \{1\}$  in   the sublemma 
depend only on the lower bounds on the  scalar curvatures of the metrics $h_t$ on $Y$ 
and on the mean curvatures of the submanifolds $Y_t=Y\times [0,t] \subset [0,t] $ with respect to the metric 
$\underline g=h_t=dt^2$ on $Y\times [0,1]$.
Thus the  sublemma remains valid for non-compact manifolds, where   the  scalar curvatures of the metrics $h_t$ 
and on the mean curvatures of the submanifolds $Y_t=Y\times [0,t] \subset [0,t] $  are bounded from below.\vspace{1mm}

{\it \textbf {Corner Corollary to SWW Theorem.}}\footnote {A version of this   was suggested in in section 6 in [G(boundary) 2019] as an approach to   "Unproven {\sf (non-extendability with $Sc>0$)} Corollary", which we will prove by a  different   argument in section \ref {capillary warped5}.}
{\sf Let $X=(X,g_0)$ be a  smooth manifold with corners. Then $X$ admits a metric $g$ with $Sc(g)>0$ and  such that all codimension 1 faces are mean convex and all dihedral angles are bounded from above by given positive numbers.}

{\it Proof.} It is obvious that there  exists a  Riemannin metric  with $Sc>0$ in a small  neighbourhood $U\subset X$  of the boundary 
$\partial X\subset X$ with respect to  which $\partial X$ is mean  convex with arbitrarily small dihedral angles. Then 
the theorem  applies to a  domain $X_0\subset X$ with smooth(!) boundary $\partial X_0\subset U$ and the proof is 
concluded  with 
Miao's  gluing lemma. \vspace {1mm}

{\it Exercise.} Let $R: G_+(X)\to H(Y)$  be the restriction map, $g\mapsto g_{|Y}$,  
  from the space $G_+(X)$ of metrics $g$ on $X$ with $Sc(g)>0$ to the space $H(Y)$ of (all) Riemannian metrics $h$ on $Y$.
   shows that 
{\sl   $R$ is a Serre fibration,} 
 
(It is not so clear 
{\sf if   the Serre  fibration property remains satisfied if   $G_+(X)$ is replaced by a subspace $G_+(X, U_0, g_0)$ of metrics that are equal to a given $g_0$ away from a small neighbourhood $U_0\subset X$  of $Y\subset X$.} \footnote{It is easy to see   if you 
 replace $H(Y)$ by the quotient space   $H(Y)/\mathbb R$ for the action of the multiplicative group  $\mathbb R$ on metrics by $r:H\mapsto r\cdot h$.})

 %(The positive answer   would allow  a generalization of  non-extendability results for $Y=\partial X$ diffeomorphic  to$S^{n-1}$  to general  $Y$) 

\vspace {1mm}

{\it \color {green!40!black}     (Naive?)} {\it Questions.}    Let $X$ be a compact manifold with a boundary. 

%(1) %{\sf {\color {red}   check with Shi} {\sf Do all Riemannian  metrics $h$ on the boundary extend to $g$ on  $X$ with 
%$Sc(g)\geq   \sigma_+$ for a {\it given  $\sigma_+>0$}}?
(1)  {\sf Does, assuming   $n=dim(X)\geq 3$, (it may be safer to assume $n\geq 5$) the manifold $X$ admits   a Riemannin metric $g$ such that
$$Sc(g)\geq\sigma\mbox {  and }  mean.curv_g(\partial X)\geq M_-$$
 for  {\it given  $\sigma_+>0$   and $M_-<0$?}}
 
 Observe the following  in this regard.
 
(i)  If $n$=2, such a $g$   seems to exist for {\it all}   $\sigma_+>0$   and $M_-<0$ {\it only} if $X$ is homeomorphic  to the 
  disc, cylinder or the M\"obius band.
  
(ii) It   is obvious that $g$ exists for all  $\sigma_+>0$   and $M_-<0$ if     $X$ contracts to the $(n-2)$-dimensional polyhedral subset $P\subset X$.

(iii) It is unclear if such   metrics exist, for all $\sigma_+>0$   and $M_-<0$,  on the $n$-torus  minus an open ball and/or     on an  $X$  homeomorphic to a compact hyperbolic manifold with a  totally geodesic boundary.
  
(iv)  If such $g$ don't always  exist,  then the  supremum of $\frac{\sigma_+}{|M_-|}$, for which  such a $g$  does  exist on an $X$, makes a  non-trivial  topological  invariant of $X$, which, one can only dream of this , would assume   several different values at certain $X$.   

(v) This may be too good to be true, but this  invariant  does make sense for {\it Riemannian} manifolds $Y=(Y,h)$,  where the above metric $g$ must extend $h$ and where the maximum of the ratios   $\frac {\sigma_+}{|M_-|}$, where such a $g$ exists is an interesting (is it?) invariant of $(Y,h)$,  evaluation of which  may be possible for  specific  manifolds $Y$, such as compact  symmetric spaces, for instance.
\vspace {1mm}

(2) Let $X$ be a compact orientable manifold with two boundary components, say $\partial X=Y_0\sqcup Y_2$ 
and let $h_0$ and $h_1$ Riemannian metrics on $Y_0$ and on $Y_1$ and let $f:Y_1\to Y_0$ be a smooth  strictly distance decreasing  map ($||df||<1$) of degree 1 (e. g. a diffeomorphism)  and let $M_0  <0$ and $M_1>0$ be  two numbers 
such that $M_0+M_1<0$.

Does the pair of     metrics  $(h_0, h_1)$ extend to a metric $g$ on $X$ with $Sc(g)>0$  and such that  the $g$-mean  curvatures of $Y_0$ and $Y_1$ 
are bounded from below by $M_0$ and $M_1$  respectively?

 %%%%%%%%%%%%%%%%%%%%%

\subsubsection{\color {blue} Obstructions to Fill-ins with  $mean.curv\geq M$ and   $Sc\geq \sigma$}\label {nonfill3}

%%%%%%%%%%%%%%%%%%%%

I. {\sf \large BMN-Counter  Example.} Motivated by Min-Oo's conjecture to the contrary,  Simon Brendle,  Fernando C. Marques  and  Andre Neves constructed  in  [Bre-Mar-Nev(hemisphere) 2011] a $C^2$-small perturbation  of the standard Riemannian  metric on the hemisphere  $S^n_+$, $n\geq 3$, that enlarges its scalar curvature while keeping unchanged  the metric and the (zero) second fundamental form  on the boundary sphere $S^{n-1}=\partial S^n_+$.
\vspace {1mm}

II.  {\sf \large BM-Non-Perturbation Theorem}  Brendle and  Marques proved in [Brendle-Marques(balls in $S^n$)N 2011] that small balls in $S^n$ admit no such perturbations and conjectured that

there is a critical radius $r_n>0$,  such that\vspace {1mm}

{\sl if a  compact Riemannin manifold $X$ with a boundary has 
 $Sc\geq n(n-1)$, and if the  mean curvature $mean.curv\partial X$  is  bounded from below by that of the $r$-ball $B^n(r)\subset S^n$, $r\leq r_n$, then $X$ is isometric to this ball.}\vspace {1mm}

III. { \sf \large STEMW Total Mean Curvature Rigidity  Theorem.}  Michael  Eichmair, Pengzi Miao and Xiadong Wang  generalized an earlier result by Yuguang Shi and Luen-Fai Tam\footnote {See [EMW(boundary) 2009] and  [Shi-Tam(positive mass) 2002]}  and proved the following.
 \vspace {1mm}

 {\sf Let $\underline X \subset \mathbb R^n$  be a star convex domain, e.g. a convex  one, such as the unit ball, for example,  
and let $X$ be a compact Riemannian manifold, the boundary $Y=\partial X$ of which is isometric to the boundary $ \underline Y=\partial  \underline X$.}

{\it If $Sc(X)\geq 0$ and if  the total  mean curvature  of $Y$ is bounded from below by that of  $\underline Y$,
$$\int_Y mean.curv(Y,y)dy \geq \int_{\underline Y}mean.curv(\underline Y, \underline y)d\underline y,$$
then $X$ is isometric to $\underline X$.} \vspace {1mm}

This is  proven in the above cited papers  by extending $g$ (from a small neighbourhood of $Y$ in $X$) to a complete asymptotically flat  metric $g_+$ on $X_+\supset X$  with $Sc(X_+)\geq 0$, where $Y$ 
serves as the     boundary of   the closure of $X_+\setminus X\subset X_+$, and  such   that \vspace {1mm}

\hspace {-0mm}{\it ADM-$mass(g_+)<0$   for  $\int_Y mean.curv(Y,y)dy > \int_{\underline Y}mean.curv(\underline Y, \underline y)d\underline y$}  \vspace {1mm}

\hspace {-6mm}and then applying  the positive mass theorem, where,  originally this was for $n\leq 7$. 
But this restriction, due to possible  singularities on minimal hypersurfaces,   may be now removed in view of the  recent results by Lohkamp and Schoen-Yau.

 \vspace {1mm}

{\it \color {blue} {\color {red!50!black} Conjecture}.} {\sf Let $X$ be a compact  Riemannian manifold with $Sc\geq \sigma$. Then 
the integral mean curvature of the boundary  $Y=\partial X$ is bounded by
 $$\int_Ymean.curv (Y,y)dy\leq const, $$
where this $const$ depends on $\sigma$ and on the (intrinsic) Riemannian  metric on  $Y$
induced from that of $X\supset Y$.}

\vspace {2mm}

IV. {\sf\large Non-Fill  for Euclidean Hypersurfaces.}  It is shown  in  [SWWZ(fill-in) 2019],   [SWW(total mean)   2020]  among other things  that a  pointwise version of    STEMW
holds for {\it non-spin}  Riemannian $n$-manifolds $X=(X,g)$  with   boundaries $Y=\partial X$ which admit 
  {\it smooth topological embeddings} to
$\mathbb R^n$:

\vspace {1mm}

{\color {blue}(A)}  {\sf if $Sc(X)\geq 0$, then the lower bound on the mean curvature of $Y$ is bounded in terms  of topology of $Y$ and  (geometry of)  $g$,
$$\inf_{y\in Y} mean.curv(Y,y)\leq  const(Top(Y), g).$$}
Furthermore,

{\color {blue}(B)} {\sf  if  $Y$ is diffeomorphic to $S^{n-1}$ and  the induced Riemannian metric $ g_{|Y}$ on $Y$  is {\it  homotopic in the set of metrics on $Y$  with $Sc>0$  to one with constant sectional curvature}, then 
{\sf $$  \int_Ymean.curv(Y,y)dy\leq  const'_n(g_{|Y}),$$
 that confirms the above conjecture in a special  case.}}
 
 This is proven  by extending $g$ from a small neighbourhood of $Y$ in $X$ to a complete asymptotically flat metric $g^+$ with $Sc\geq 0$, where $Y$ 
serves as the     boundary of the closure of $X_+\setminus X \subset X_+$  and  such   that the ADM mass of  $g^+$ is negative  provided {\it the mean curvature of (or its integral over) $Y$ is sufficiently large}. Then the  the positive mass theorem applies.

 \vspace {1mm} 
 
 {\it Remark.} Probably, by incorporating Lohkamp reduction of the positive mass theorem to the flat at infinity case (see section \ref{asymptotic3})
 one can make   $g_+$   {\it flat}, rather than only  {\it  asymptotically  flat  with mass$\leq 0$} at infinity, where this may generalize  to manifolds $Y$ that are  not necessarily diffeomorphic to $S^{n-1}$.

   \vspace {1mm} 
  
  V. {\sf \large Pointwise non-Fill-in for Compact $Y$.} Pengzi Miao found a simple derivation  of the following version of (A)   from  the SWW extension theorem   for  {\it all} $Y$
   [Miao(nonexistence of fill-ins) 2020]:
 $$\inf_{y\in Y} mean.curv(Y,y)\leq  const(Top(X), g).$$

 {\it Proof.} Given $X=(X,g)$ with $Sc(g)\geq 0 $   and $mean.curv_g(Y)\geq \mu_+$, 
 $Y=\partial X$ let $X'$ be a connected sum of $X$ with the $n$-torus  and let $g'$ be a metric with $Sc>0$ 
 such the restriction of $g'$ to $Y=\partial X'=\partial X$ is equal to the restriction of $g^{|Y}$, where the existence of such a $g'$ is guaranteed by the SWW  extension theorem. 
  
  Observe that the supremum  $\mu_\ast'= \sup_{g'} \inf_{x'\in Y}mean.curv_{g'}(Y, x')$   depends on  the topology $X$
 and on the restriction of  $g$ to $Y$.
 
 Also observe   that  the manifold $X \sqcup_Y X'$ obtained by gluing $X$ and $X'$ along $Y$ admits  no metric with 
 $Sc>0$ by the Schoen-Yau theorem.
 
 But if $\mu_++\mu'_\ast>0$, the natural continuous metric $g\&g'$ on   $X \sqcup_Y X'$ can be smoothed with
 $Sc>0$ by Miao gluing theorem; hence, $\mu_+\leq -\mu'_\ast$. QED.
  
\vspace {1mm}

VI. Another derivation of pointwise non-fill-in theorem from SWW  extension theorem is with a use of the Corner Corollary
from the previous section. Indeed, if the mean curvature of $\partial Y$ is sufficiently large, one can modify  the Riemannin metric on
  $X$  (by attaching an external color to $X$ along   $Y\partial X$) keeping $Sc\geq 0$ and creating cubical corner 
structure on the boundary with dihedral angles $<\frac{\pi}{2}$,  as in    "Unproven Corollary" from section 6 in [G(boundary) 2019]. 

Then, by the reflection argument from section 
\ref  {reflection3}, the problem reduced to Schoen-Yau theorem on non-existence of metrics with $Sc>0$   on manifolds which admits maps with  non-zero degrees to tori.

VII.  In the case of spin manifolds $X$ a more  precise non-fill inequality  follows from the  mean curvature spin-Extremality theorem  in section \ref{mean convex3}, and a  $\mu$-bubble  approach to the non-spin  case is indicated in section \ref{capillary warped5}.

\vspace{1mm}

{\it Questions.}
Let $X$ be compact $n$-manifold with boundary, let $Y_i\subset \partial X$ be  be the connected components of the boundary. (For instance, $X$ is the $n$-torus minus two open balls and $\sigma=1$.)
 
(a) Given numbers $\sigma$ and $M_i$, when does there exist a Riemannin metric $g$ on $X$, such that
$Sc(X)\geq \sigma$ and the mean curvatures of $Y_i$ are bounded from below my $M_i$, 
$$mean.curv_g(Y_i)\geq M_i?$$
 
 (b) Let all $Y_i$ be diffeomorphic to the sphere $S^{n-1}$ and let, besides $\sigma$ and $M_i$, we are  given 
 positive numbers $\kappa_i$. 
 
 When does there exist a Riemannin metric $g$ on $X$, such that
$Sc(X)\geq \sigma$, the induced metrics  $g_{|Y_i}$ have constant  sectional  curvatures $\kappa_i$ and
$$mean.curv_g(Y_i)\geq M_i?$$
 
 (c) Let now  Riemannin   metrics $g_i$ on $Y_i$ be given.  When does there exist a Riemannin metric $g$ on $X$, such that
$Sc(X)\geq $,  $g_{|Y_i}-g_i$  and
$$mean.curv_g(Y_i)\geq M_i?$$

 %%%%%%%%%%%%%%%%%%%%%

\subsection {\color {blue} Manifolds with Negative Scalar Curvature Bounded from Below }\label {negative3}
%%%%%%%%%%%%%%%%%%%%%%%
 If a    "topologically complicated" closed Riemannian  manifold $X$, e.g. an  aspherical one with a {\it hyperbolic fundamental group}, has $Sc(X)\geq \sigma $ for $\sigma <0$, then a certain  "growth" of the universal covering $\tilde X$ of $X$ is expected to be {\it bounded from above}  by 
$const\sqrt{-\sigma}$ and accordingly, the "geometric size" -- ideally $\sqrt[n]{vol(X)}$-- must  be {\it bounded from below} by $const'/\sqrt{-\sigma}$.

If $n=3$ this kind of lower bound are easily available  for  areas of stable  minimal surfaces of large genera via Gauss Bonnet theorem 
by the Schoen-Yau argument  from [SY(incompressible) 1979].

Also  Perelman's proof of the geometrization conjecture  delivers  a sharp  bound of this  kind for manifolds $X$ with hyperbolic $\pi_1(X)$ and   similar results for $n=4$  are possible with the Seiberg-Witten theory for $n=4$ (see section \ref {Seiberg3}).

No such estimate has been established yet for $n\geq 5$ but the following results are available.
\vspace {1mm}

\textbf{Ono-Davaux Hyperbolic Spectral Inequality}.\footnote {See  [Ono(spectrum) 1988], [Davaux(spectrum) 2003].} {\sf Let $X$ be a closed Riemannian manifold 
and let $\tilde X\to X$ be some Galois covering  of $X$,  e.g the universal covering, such that   all smooth functions $f(\tilde x)$ with compact supports on $\tilde X$ of $X$ satisfy
$$ \int_{\tilde X} f(\tilde x)^2d\tilde x  \leq\frac {1}{\tilde \lambda_0^2 } \int_{\tilde X} ||df(\tilde x)||^2d\tilde x.$$}
 (The maximal such $\tilde \lambda_0\geq 0 $ serves as the lower bound on  the spectrum of the Laplace  on the universal covering $\tilde X$ of $X$).  

{\it If $\tilde X$ is spin and if  one of the following two conditions (A) or (B) is satisfied, then
$$\inf_{x\in X} Sc(X,x)\leq \frac {-4n\tilde \lambda_0}{n-1}.\leqno {\color {blue}[Sc/\tilde \lambda_0]}$$}

{\sl \large Condition (A)}. {\sf  The dimension of $X$ is $n=4k$ and the $\hat \alpha$-invariant from section \ref {spin index3} (that is a certain linear combination of Pontryagin numbers called  {\it $\hat A$-genus}) doesn't vanish.}\vspace {1mm}

{\sl \large Condition (B)}. {\sf The  manifold $\tilde X$ is {\it hypereuclidean}: it properly Lipschitz dominate 
the Euclidean space, i.e.  $\tilde X$ is orientable and it admits a proper  distance  decreasing map to $\mathbb R^n$
with  {\it non-zero} degree.}\vspace{1mm}

\vspace{1mm}

{\it Idea of the Proof.} By {\it Kato's inequality} (and/or by {\it the Feynman-Kac formula}, see   \ref {Kato6}),  the
 lower bound on $\tilde \lambda_0$ implies a similar bound on the Bochner Laplacian $\nabla^2$  on $\tilde X$, hence, a corresponding  
 bound on  the (untwisted) 
  Dirac operator  expressed  by the SLW(B)-formula   $\mathcal  D=\nabla^2+\frac {1}{4}Sc.$

This, 
confronted with the $L_2$-index theorem, yields {\sl \large Condition (A)} and  {\sl \large Condition (B)}   for $n$ even follows by 
 similar argument for $mathcal D$ on $\tilde X$ twisted with suitable almost flat bundles, while the
   sharp inequality  for odd $n$ needs a an odd dimensional version of the $L_2$-index theorem and a  
    delicate analysis of the spectral flow for  a  family of Dirac operators (see  [Davaux(spectrum) 2003]).\vspace{1mm}

{\it Remarks.} (a)  The inequality {$\color {blue}[Sc/\tilde \lambda_0]$} is sharp: if $X$ has constant negative curvature $-1$,  then
$$\mbox { $-n(n-1)=Sc(X)= \frac {-4n\tilde \lambda_0}{n-1}$}$$
 for $\tilde \lambda_0=\frac {(n-1)^2}{4}$, that is 
the bottom of the spectrum of $\mathbf H^n_{-1}=\tilde X$.\vspace {1mm}

(b) The rigidity sharpening of {$\color {blue}[Sc/\tilde \lambda_0]$} is proved in   [Davaux(spectrum) 2003] in the case $A$ and it seems that  a  minor   readjustment of the   argument from this paper  would work  in the case $B$ as well.

(c)  Since the spectrum of the Laplacian is Lipschitz continuous under $C^0$-deformations of Riemannian metrics,  the 
 Ono-Davaux hyperbolic spectral inequality implies, for instance,  that 
 
 {\sl if a metric $g$  on a compact $n$-manifold $X$ is $\lambda$-bi-Lipschitz, $\lambda\geq 1$, to a metric $g_0$ with sectional curvatures $\kappa\leq -1$,
then $$\inf_{x\in X}Sc(g, x)\leq -\frac {const_n} {\lambda^{2}}, \mbox { }   const_n>0.$$}

 On the other hand, the spectrum of $\Delta$ drastically  drops down, if for instance one takes connected sums of $X$ with spheres $S^n$   attached to $X$ by long  narrow "necks"  by means of the thin surgery  with only  a minor perturbation of the infimum of the scalar curvature.\vspace{1mm}

 \textbf {Non-Amenable  Hypereuclidean Manifolds with $Sc\geq \sigma<0$.} Probably, 
 the above bound on the scalar curvature in case {\large \sf  (B)} remains true for all complete Riemannin manifold $\tilde X$ ,with no spin assumptions and with no action of any (deck transformation) group on it.
 
Below is a result in this direction, which we formulate in geometric rather than analytic terms.  

 A Riemannian $n$-manifold $X$ is  
called (uniformly) {\it $\alpha$-non-amenable} if all compact smooth domains $U\subset X$ satisfy the {\it linear isoperimetric inequality 
with constant $\alpha$,}
$$ vol(U)\leq \alpha\cdot  vol_{n-1}(\partial U).$$

It is easy  to see that  

{\sf if a {\it complete  $\alpha$-non-amenable} Riemannin  $n$-manifold $X$ has {\it bi-Lipschitz bounded local geometry}, 
i.e. all $\delta$-balls in $X$ are $\lambda$-bi-Lipschitz homeomorphic to the Euclidean $\delta$-ball for some
positive numbers $delta$ and $\lambda$ depending on $X$, then $X$ can be {it exhausted} by  compact  smooth domains $V_i$,
$$V_1\subset V_2\subset ... \subset V_i \subset... \subset X$$
such that boundaries $Y_i\partial V_i$ satisfy 
$$mean.curv(Y_i)\geq \alpha-\varepsilon _i, \mbox  {  where  $\varepsilon _i\to 0$ for $i\to \infty$}.$$

Indeed,  let $\mu_i(x)$ be a sequence of smooth functions on $X$, such that 

{\Large $ \ast$} {\it all $\mu_i$ is very large} at a point $x_0\in X$, 

 {\Large $\ast$}  all $\mu_i(x)<\alpha-\epsilon_i $  for $x$ very far from $x_0$,  

{\Large$\ast$ } the gradients of all  $\mu_i(x)$ is very small  at all $x\in X$,

{\Large$\ast$ } $\mu_i(x)\to \alpha$ for $i\to \infty$ and $x\to \infty$.}

Then our conditions on $X$ imply the existence of  $\mu_i$-bubbles $Y'_i=\partial V'_i\subset X$, 
where $V_i$ {\it exhaust} $X$   and where $Y'_i$  can be smoothed to the required $Y_i=\partial V_i$  (compare with 
1.5(C) in [G(Plateau-Stein) 2014 ]).

This, combined with  {\it multi-width mean curvature inequality} from section \ref{capillary warped5}, 
yields the following. \vspace {1mm}

{\it \textbf {Rough Negative   bound on Sc(X)}}. {\sf Let $X$ be an $\alpha$-non-amenable  {\it hypereuclidean} Riemannian  $n$-manifold.} 
{\it If $n\leq 7$,  then the  infimum of the scalar curvature of $X$ is bounded by $\alpha$ as follows
$$\inf_{x\in X}Sc(X,x)\leq -const_n \alpha^\frac {2(n-1)}{n}$$}
for some $const_n>0$.\vspace {1mm}
 
Let us indicate an application of this to  manifolds $X$ discretely a cocompactly acted upon by , 
countable groups $\Gamma$, e.g. to universal coverings of compact manifolds, where
$\Gamma$ {\it uniformly no-amenable}, where

{\sl  a finitely generated group  $\Gamma$ is 
{\it uniformly non-amenable} if there exists an $\underline \alpha=\underline\alpha(\Gamma)>0$, such that, for all 
symmetric finite  generating subset 
$\Delta  \subset \Gamma$ the cardinalities of  all finite subsets $V\subset \Gamma$  are bounded 
by the cardinalities of their $\Delta$-boundaries, 
$$card(V)\leq \underline\alpha \cdot card (\Delta \cdot S\setminus S.$$}

{\it Example.} Non-virtually solvable subgroups of the linear group $GL(N,\mathbb C)$ are 
uniformly non-amenable (see [Breuillard-Gelander(non-amenable) 2005] and references therein)

To use this, we observe that if  a Riemannin $n$-manifold is discretely and isometrically  acted upon by such a $\Gamma$ 
with compact quotient $X/\Gamma$, then the isoperimetric constant $\alpha=\alpha(X)$ is  bounded from below in terms of  
$\underline\alpha(\Gamma$  and a bound on the local Lipschitz geometry of $X$.

Namely, given numbers    $\lambda>0$, $d$  and $\underline\alpha>0$, there exists $\alpha=\alpha_n(\lambda, , d\underline\alpha)>0$, 
 such that 
if all unit balls in $X$ are $\lambda$  bi-Lipschitz homeomorphic to the unit  Euclidean ball,
and the diameter of the quotient space is bounded by $diam( X/\Gamma)\leq $, 

then,  by an easy argument,  \vspace {1mm}

\hspace {7mm}{\it the inequality $\underline\alpha(\Gamma) \geq  \underline\alpha$     implies that     $X$ is $\alpha$-non-amenable.}     

\vspace {1mm}

{\it Remark/{\color {red!50!black} Conjecture}.} The  conditions on the Lipschitz geometry and the diameter are unpleasantly restrictive.
{\color {red!50!black} Conjecturally} all one needs is a {\it bound on the volume of $X/\Gamma$.}

\vspace {1mm}
\vspace {1mm}

The last theorem in tis section  formulated below, was, historically, the first result on  the geometry of 
$Sc\geq \sigma$ for {\it negative} $\sigma$. \vspace {1mm}

Thus, \vspace {1mm}

 { \it  if $X$ is hypereuclidean, then the infimum of the  scalar curvature  $ Sc(X)$ is bounded by a {\it strictly negative} constant 
which depends only on $\Gamma$ and the bound on the local Lipschitz geometry of $X$.}

 \vspace {1mm}

{\it Exercise.}  Show that the universal covering  $\tilde X$ of the $n$-torus $X$ with an arbitrary Riemannian metric 
can be exhausted by {\it over-cubical}  domains $V_i\subset \tilde X$  {\it with corners}, i.e.  such that 
 they admit face preserving  ({\it corner proper} in terms of section\ref{corners3}) maps  
$$f_i :V_i\to [0,1]^n$$ 
with   {\it degree} 1 and such that all $(n-1)$-faces  of all $V_i$ have {\it positive}  mean curvatures 
and the {\it dihedral angles} of all  $V_i$ along  the  $(n-2)$-faces {\it are  $\leq \frac {\pi}{2}$.}
 
 {\it Hint.}  Cut the manifold $X$  by a minimal hypersurface in the homology class of $\mathbb T^{n-1}\subset \mathbb T^n\underset{homeo} \simeq X$,
 then cut the resulting band by a sub-band homologous to $\mathbb T^{n-2}\times  [0,1] $ etc. 
 If $n\leq 7$ this terminates  in  a cubical $V_1\subset \tilde X$, and by applying the same procedure to  finite coverings 
of $X$  with fundamental subgroups $i\cdot\mathbb Z^n\subset \mathbb Z^n=\pi_1(X_)$ we obtain an exhaustion of $\tilde X$ by over-cubical $V_i$ with {\it minimal} $(n-1)$-faces and all {\it dihedral} angles $\pi/2$.

These $V_i$ may have, however, not very smooth faces and an extra work is needed to smooth them.

 And if $n\geq  8$,  such  $V_i$ , come, in general,  with more serious singularities, but one can smooth them keeping the 
 $(n-1)$-faces  mean convex and the dihedral angles  $\leq \pi/2$, as it is  done in  [G(Plateau-Stein) 2014].

\vspace {1mm}

{\it\color {red!50!black} Remark/Conjecture.} It is not impossible (but unlikely)  that {\it all contractible} manifolds $ \tilde X$ 
which admit {\it cocompact isometric group actions} 
also admit similar over-cubical exhaustions, where   this seem quite   realistic  for {\it enlargeable} $X$. 

 Also  other "large" manifolds $\tilde X$ without   any group actions, e.g.  complete simply connected manifolds with non-positive sectional curvatures  admit such exhaustions or, at least, contain arbitrarily  large man convex  overtorical domains with dihedral angles $\leq \pi/2$.\vspace {1mm}

{\it Question.} {\sf What are possible  values of dihedral angles of large {\it non-over-cubical} domains 
with corners in various  manifolds?}

For instance, it seem not hard to show in this regard that  the 2-plane with a metric {\it bi-Lipschitz homeomorphic to the hyperbolic plane} can be exhausted by convex  $k$-gons,  for all $k=2, 3, 4,...$ with all {\it angles $\leq \varepsilon$ for
all $\varepsilon>0$.}

Also it seems {\color {red!50!black}not impossible} that, for {\it all  convex  polyhedral domains}  $P\subset \mathbb R^n$,  the above  universal covering $\tilde X\to X\underset{homeo} \simeq \mathbb T^n$ can be exhausted  by mean convex  "{\it over $P$-domains}" $V_i$  (admitting face respecting maps  $V_i\to P $ with  degrees 1), such that   the  {\it dihedral angles of all $V_i$ are bounded  by the corresponding angles in $P$,}
$$\angle_{kl} (V_i)\leq \angle{kl} (P),$$
and where, moreover, unless $X$ is Riemannin  flat,  one can find/construct such  $V_i$ with $\angle_{kl} (V_i)<\angle{kl} (P)$.

\vspace {1mm}\vspace {1mm}

 The last theorem   in this section we state below was, historically,  the first result on {\i geometry} of $Sc\geq \sigma$ for $\sigma<0$.\footnote {Strictly speaking, the first, for alI know,  {\it topological-geometric}  constraint on $Sc\geq \sigma <0$  appears  in [Ono(spectrum) 1988], but his argument resides within the realm of $Sc\geq 0$.} \vspace {1mm}

\textbf{ Min-Oo Hyperbolic Rigidity Theorem.} {\sf Let $X$ be a complete Riemannian manifold,  which  is isometric at infinity (i.e. outside a compact subset in $X$) to the hyperbolic space $\mathbf H^n_{-1}$.}

\hspace {10mm}{\it If $Sc(X)\geq -n(n-1)=Sc(\mathbf H^n_{-1})$,  then $X$ is isometric to $\mathbf H^n_{-1}$.}

\vspace {2mm}

{\it About the Proof.} The original argument  in  [Min-Oo(hyperbolic) 1989],  which generalizes  Dirac-theoretic  Witten's proof of the positive mass/energy theorem for asymptotically Euclidean (rather than hyperbolic) spaces, (see section \ref{asymptotic3})
needs $X$ to be {\it spin}. 

But  granted spin,  Min-Oo's  proof allows more general   asymptotic (in some sense) agreement between $X$ and $\mathbf H^n_{-1}$ at infinity.

If one wants to   get rid of spin, one can use minimal hypersurfaces or $\mu$-bubbles, where it  is convenient,   to pass to a quotient space  $\mathbf H^n_{-1}/\Gamma$,  where $\Gamma$ is
a {\it parabolic} isometry group isomorphic to $\mathbb Z^{n-1}$, and where  the quotient    $\mathbf H^n_{-1}/\Gamma$ is the  {\it hyperbolic cusp-space}, that is $\mathbb T^{n-1}\times \mathbb R $ with the metric $e^{2r}dt^2 +dr^2$.\footnote {The logic of what we do here  is similar to   the proof  of   rigidity of $\mathbb R^n$ by passing to
 $\mathbb T^n=\mathbb R^n/\mathbb Z^n$  and thus, reducing the problem to   the scalar curvature $\geq 0$ rigidity of the flat tori.}

Then one applies  the {\sf rigidity theorem  for the flat  metrics on tori with $Sc\geq 0$}   to $\mathbb T^1$-{\sl symmetrised  stable $\mu$-bubbles} in    manifolds $X$ {\it isometric to $\mathbf H^n_{-1}/\Gamma$  at infinity}, where these bubbles {\it separate   the two ends}  of $X$ and where $\mu=(n-1)dx$. Thus one   shows  that

 \vspace {1mm}

{\it n-manifolds $X$ with $Sc(X)\geq -n(n-1)$,   which  are isometric to $\mathbf H^n_{-1}/\Gamma$ {\it at infinity}, are,  isometric to  $\mathbf H^n_{-1}/\Gamma$ {\it everywhere}.}\footnote {Instead of  using $\mu$-bubbles as in  
 \S$5\frac {5}{6}$ of    [G(positive) 1996], one can proceed  here by  inductive descent  with 
 $\mathbb T^\rtimes$-symmetrised  minimal hypersurfaces with free boundaries, as in the proof of the  $\frac {2\pi}{n}$-inequality indicated in section 3.6;
    see  section \ref {separating}  for this and  for  more general results of this kind.}

 \vspace {1mm}

Finally, a derivation of a  Min-Oo's  kind  {\it hyperbolic positive mass theorem} without the spin condition  from the rigidity theorem follows by an extension  of the Euclidean  Lohkamp's  argument from  to the hyperbolic spaces, due to  {\sl Andersson, Cai, and Galloway}.\footnote{See  [Lohkamp(hammocks) 1999]  and    [AndMinGal(asymptotically hyperbolic) 2007].}

 %Thus, the  Ono-Davaux inequality can't be used for proving  Min-Oo  kind theorems, but, conceivably, this can be amended by an argument in the spirit of  [Cecchini(long neck) 2020].

\vspace {1mm}

{\it Questions.}{\sf Can one put the index theoretic and  associated  Dirac-spectral considerations  on equal footing  with Witten's and Min-Oo's kind of  arguments on stability of harmonic spinors with a given asymptotic behavior under deformation/modifications of manifolds away from infinity?

Can Cecchini kind    long neck argument(s)   be extended to $\sigma<0$?}\footnote { Notice that the long neck proofs in [Cecchini(long neck) 2020] and in [Cecchini-Zeidler(generalized Callias) 2021], similarly to these in [Min-Oo(hyperbolic) 1989], depend  on Dirac operators with potentials.}

\vspace {1mm}

 \vspace {1mm}

\hspace {38mm} {\sc Three  conjectures } \vspace {1mm} \vspace {1mm}

{ \color {blue}\textbf {[\#$_{-n(n-1)}$]} }  {\sl     \color {blue!40!black}  Let $X$ be a closed orientable Riemannian  manifold of dimension $n$ with $Sc(X)\geq -n(n-1).$}\vspace {1mm}

Then  the following topological invariants of $X$ must be  bounded by the volume of $X$, and, even more optimistically,   (and less realistically), where  the constants  are  such that   the equalities are  achieved  for compact  hyperbolic  manifolds with sectional curvatures $-1$.

Namely,  granted { \color {blue}\textbf {[\#$_{-n(n-1)}$]} } one expects the following.
\vspace {1mm}

1.   {\it  \textbf {Simplicial Volume {\color {red!50!black} Conjecture}}}: 
{\sf \color {blue!10!black}  There exist orientable $n$-dimensional {\it  \color {blue!40!black}  pseudomanifolds} 
$ {\sf X}^\smalltriangleup_i$ and continuous  maps $f^\smalltriangleup_i: {\sf X}_i^\smalltriangleup\to X$ 
with degrees $$deg(f^\smalltriangleup_i)\underset {i\to \infty}\to \infty,$$ 
such that  
the numbers  $N_i$  of {\it  \color {blue!40!black}  simplices} in the triangulations of  ${\sf X}^\smalltriangleup_i$ and the degrees $deg(f^\smalltriangleup_i)$  are related to the volume of $X$ by the following inequality:
$$ N_i\leq C^\smalltriangleup_n \cdot deg(f^\smalltriangleup_i)\cdot vol(X).$$}

{\rm 2}.  {\it  \textbf {The L-Rank Norm {\color {red!50!black} Conjecture}:} {\sf \color {blue!10!black} 
 There exist,  for all sufficiently large  $i\geq i_0=i_0(X)$, smooth  orientable $n$-dimensional   {\it  \color {blue!40!black}  manifolds }
$ {\sf X}^\circ_i$   and continuous maps $f_i: {\sf X}_i^\circ\to X$, 
with degrees 
$$deg(f^\smalltriangleup_i)\underset {i\to \infty}\to \infty,$$ 
such that the minimal possible numbers  $N_i$ of 
{\it  \color {blue!40!black}  the cells} in the  cellular decompositions of $ {\sf X}^\circ_i$ and the degrees of the maps $f^\smalltriangleup_i$  are related to the volume of $X$ by the following inequality:}
$$ N_i\leq C^\circ_n \cdot deg(f^\circ_i)\cdot vol(X).$$}

3.   \textbf{ Characteristic Numbers {\color {red!50!black} Conjecture}. }  {\sf if, additionally to  {  \color {blue}\textbf {[\#$_{-n(n-1)}$]}},  the manifold  $X$ is  {\it  \color {blue!40!black} aspherical},  then {\it  \color {blue!10!black} the Euler 
characteristic} $ \chi(X)$  and {\it  \color {blue!10!black} the Pontryagin numbers $p_I$ of $X$ are bounded by 
$$ |\chi(X)|,  |p_I(X)| \leq C^\smallsquare_n \cdot vol(X).$$}}

\vspace {1mm} 

{\it Remarks.}  (i)  {\color {red!50!black} Conjecture} 1 makes sense for an $X$,  in so far as  $X$ has  {\it non-vanishing} simplicial volume  
$||X||_\triangle$,  e.g.  if  $X$   admits a  metric  with {\it negative sectional curvature} or a
 locally symmetric metric  with {\it negative Ricci curvature}. \footnote {See [Lafont-Schimidt(simplicial volume) 2017]  and   the monograph  [Frigero(Bounded Cohomology) 2016] for the definition and basic properties of the simplicial volume.}	

\vspace {1mm}

(ii) The $L$-rank norm  $||[X]_L||$ is   defined in \S$8\frac{1}{2}$ of [G(positive) 1996] via the  Witt-Wall $L$-groups of the fundamental group of $X$.

 This  $||[X]_L||$  is known to  be {\it non-zero} for compact  locally symmetric spaces with non-zero Euler characteristic as it follows from [Lusztig(cohomology) 1996].\footnote{ In the simplest case, where  $X$ is the product of  $k$  closed surfaces $S_1,S_2,...,S_k$   with negative Euler characteristics, non-vanishing of  $||[X]_L||$ is proven in [G(positive) 1996]: 
 
 {\sl If a manifold $X^\circ$ admits a map of degree $d$ to such an $X$, then  $X^\circ$ can't be decomposed into less than
  $$N=const_k\cdot  d\cdot |\chi(S_1)|\cdot | \chi (S_2)\cdot ... \cdot |\chi ( S_k)|, \mbox { }  const_k>0,$$ 
  cells.}}

In fact, all {\it known} manifolds $X$  with  $||[X]_L||\neq 0$ admit maps of non-zero degrees
  to locally symmetric spaces with non-zero Euler characteristics.
  
 And  nothing is  known about zero/non-zero possibility for the values of the L-rank norm for manifolds
  with negative sectional  curvatures   of  odd dimensions $> 3$. 
  
  (Vanishing of  $||[X]_L||$ for all 3-manifolds $X$ trivially follows from  the Agol-Wise theorem on virtual fibration of 
  hyperbolic 3-manifolds over $S^1$.)
  \vspace {1mm}

{\it Question.}   What are realtions between the   $||X||_\triangle$ and $||[X]_L||$?  Are there natural  invariants  mediating between the two? \vspace {1mm}

(It is tempting to suggest that    $||X||_\triangle\geq ||[X]_L||$, since being a triangulation is (by far) more restrictive than being just a  cell decompositions, but since  $||[X]_L||$,  unlike $||X||_\triangle$ defined with {\it manifolds}, rather than with {\it pseudomanifolds} mapped to $X$,  this is unlikely to be true in general.)\vspace {1mm}

 {\it \textbf {Integral Strengthening of the Three {\color {red!50!black} Conjecture}s}}.  The above conjectural  inequalities  1,2,3,  for the three  topological invariants, call them here $inv_i $, $i=1,2,3$,  may, for all we know,   hold (with no a priori assumption $Sc(X)\geq - n(n-1)$)
 in the following integral form,
{ \color {blue!60!black} $$inv_i\leq const_i\cdot \int _X|Sc_-(X,x)|^{\frac{n}{2}}dx,$$} where $Sc_-(x)=\min (Sc(x), 0)$,   but no lower bound on this integral
is anywhere  in sight for $n\geq 5$. \footnote{One doesn't even know if there are such bounds for  $||X||_\triangle$ and/or $||[X]_L||$ in terms of the {\it full  Riemannian curvature tensor} $R(X,x)$, namely the bounds  
    $$||X||_\triangle, ||[X]_L|| \leq const_n \int_X||R(X,x)||^\frac {n}{2}dx.$$} (See section \ref {Seiberg3}  for what is known for $n=4$.)

 %%%%%%%%%%%%%%%%%
 \subsection { \color {blue} Positive Scalar Curvature, Index Theorems  and the Novikov Conjecture} \label {Novikov3}

%%%%%%%%%%%%%%%%%

Given a proper (infinity goes to infinity)  smooth map between smooth oriented manifolds,
$f:X\mapsto \underline X$ of dimensions $n=dim(X)= 4k+ \underline n$ for  $\underline n=dim(\underline X)$,
let $sign (f)$ denote the signature of the  pullback $ Y_{\underline x}^{4k}=f^{-1}( \underline x)$  of a generic point $ \underline x\in \underline X$, that is the signature of the (quadratic)  intersection form on the homology $H_2(  Y_{\underline x}^{4k};\mathbb R),$ where, observe, orientations of $X$ and $ \underline X$ define an orientation of  $ Y_{\underline x}^{4k}$ which is needed for the definition of the intersection index.

Since the $f$-pullbacks of generic (curved) segments $[\underline x_1, \underline x_2]\subset \underline X$  are manifolds with boundaries  $Y_{\underline x_1}^{4k}-Y_{\underline x_2}^{4k}$ , (the minus sign means the reversed orientation), 
  $$sign(Y_{\underline x_1}^{4k}) = sign(Y_{\underline x_2}^{4k}),$$
  as it follows from  the Poincar\'e duality for manifolds with boundary by a 
   two-line argument.
   Similarly, one sees that $sign(f)$ depends only on the proper homotopy class $[f]_{hom}$  of $f$. 
 
  Thus, granted $\underline X$ and a  proper homotopy class of maps $f$, the signature  $sign[f]_{hom}$ serves as a {\it smooth invariant} denoted  $sign_{[f]}(X)$, (which is actually  equal to  the value of some polynomial in Pontryagin classes of $X$ at the homology class
  of $Y_{\underline x_2}^{4k}$ in the  group $H_{4k}(X)$).
  
   If $X$ and   $\underline X$ are closed manifolds, where  $dim(X)>dim (\underline X)>0$,  and if  $\underline X$,  is {\it simply connected}, then, by the Browder-Novikov theory, as one varies the smooth structure of $X$ in a given homotopy class $[X]_{hom}$  of $X$, the values of $sign_{[f]}(X)$ run through {\it all integers}     $i=sign_{[f]}(X) \mod 100n!$ (we exaggerate for  safety's sake), provided $dim(\underline X)>0$ and $Y_{\underline x}^{4k}\subset X$ is non-homologous to zero.

 However, according to  the (illuminating special case of the) {\color {blue} \it Novikov conjecture}c \vspace {1mm}
 
  if $\underline X$ is a { \it closed aspherical} manifold\footnote {{\it Aspherical} means that  the universal cover of $\underline X$ is contractible} then 
 this $sign_{[f]}(X)$ {\it depends only on the homotopy type of $X$.} \footnote {Our topological formulation, which is  motivated by the history  of the Novikov conjecture, is  deceptive: in truth, Novikov conjecture  is 90\% about infinite groups, 9\%  about geometry and  only  1\% about manifolds.}  \vspace {1mm}

 Originally, in 1966, Novikov proved this, by an an elaborated  surgery argument,  for the torus    $\underline X=\mathbb T^{\underline n}$, where $X=Y\times \mathbb T^{\underline n}$ and $f$ is the projection $Y\times \mathbb T^{\underline n}\to  \mathbb T^{\underline n}.$
 
In 1972,  Gheorghe  Lusztig found a  proof for  general $X$ and  maps $f:X\to T^{\underline n}$ 
 based on {\it the Atiyah-Singer index theorem for families  of differential operators $D_p$ parametrised by topological spaces $P$,} where the index takes  values not in $\mathbb Z$  anymore but  in the $K$-theory of $P$, namely, this index is defined as the $K$-class of the (virtual) vector bundle over $P$ with the fibers $ker(D_p)-coker(D_p)$, $p\in P$, (Since the operators $D_p$ are Fredholm, this makes sense despite possible   non-constancy of the ranks  of  $ker(D_p)$  and $coker(D_p)$.)

 \vspace {1mm}

 The family $P$  in  Lusztig's proof in  [Lusztig(Novikov) 1972]  is composed of the {\it signature s} on $X$ twisted with complex line bundles $L_p$, $p= P$, over $X$,  where these $L$ are induced by a map $f:X\to\mathbb   T^{\underline n}$  from {\it flat}   complex unitary  line bundles $\underline L_p$ over $\mathbb  T^{\underline n}$  parametrised by  $P$ (which is the 
 $\underline n$-torus   of homomorphism $\pi_1( \mathbb  T^{\underline n})=\mathbb Z^{\underline n}\to \mathbb T$). 

Using the the  Atiyah-Singer index formula, Lusztig computes the index of this , shows that it is equal to $sign(f)$ 
and deduce  from this    the homotopy invariance of  $sign_{[f]}(X)$.\vspace{1mm} 

What is relevant for our purpose is  that Lusztig's computation equally applies to the Dirac   operator  twisted with $L_p$ and shows the following. \vspace{1mm} 

{\it Let  $X$ be a closed  orientable  spin manifolds of even  dimension $ \underline n$ and  $f:X \to \mathbb T^{\underline n}$ be continuous map  of non-zero degree. Then
$$ind({\cal D}_{\otimes \{L_p\}})\neq 0.$$}

Therefore, there exits a point $p\in P$, such that  $X$ carries a harmonic  $L_p$-twisted spinor

But if $Sc(X)>0$,  this
is incompatible with the  the Schroedinger-Lichnerowicz-Weitzenboeck-(Bochner) formula \ which says for {\it flat} $L_p$ that $${\cal D}_{\otimes L_p}= \nabla_{\otimes L_p}^2+\frac {1}{4} Sc(X).$$

Thus, \vspace {1mm} 

{\it  the existence of a map $f:X\to \mathbb  T^{\underline n}$ with $deg(f)\neq 0$ implies that $X$ carries 

no metric with $Sc>0$.} \vspace {1mm}

Moreover, Lusztig's computation applies to manifolds $X$ of all dimensions   $n= \underline n+4k$, shows that if  a generic pullback manifold $Y_p^4=f^{-1}(p)\subset X$ (here $f$ is smooth) has {\it non-vanishing  $\hat \alpha$-invariant } defined in section \ref {spin index3} (that is the  $\hat A$-genus for 
$4k$-dimensional manifolds), then the index  $ind({\cal D}_{\otimes \{L_p\}})$ doesn't  vanish either and, assuming $X$ is  spin,   it  {\it can't carry metrics with $Sc>0$.}
\vspace{1mm}

{\it Remark on $X=(X, g_0)=\mathbb T^{\underline n}$.} If $(X,g_0)$ is isometric to the torus,  then  the only $g_0$-harmonic $L_p$-twisted spinors on $X$ are parallel ones, which allows a direct computation of the index of  ${\cal D}_{\otimes \{L_p\}}$. Then the result of this  computation extends  to all Riemannian metrics $g$ on $\mathbb  T^{\underline n} $ by the invariance of the index of 
 ${\cal D}_{\otimes \{L_p\}}$ under deformations of $\cal D$, where the essential point is that, albeit the harmonic spinors of 
the (untwisted) $\cal D$ may (and typically do) disappear under a  deformation   ${\cal D}_{g_0}\leadsto{\cal D}_g$, they re-emerge as    harmonic spinors 
of ${\cal D}_g$ twisted with a non-trivial flat bundle $L_p$.\vspace{1mm}

The index theorem for families can be reformulated with $P$ being replaced by  the algebra $cont(P)$ of all continuous functions on $P$, where in Lusztig's case the algebra $cont(\mathbb  T^{\underline n})$ is Fourier isomorphic to the algebra $C^\ast(\mathbb Z^{\underline n})$
of bounded linear operators on the Hilbert space space $l_2(\mathbb Z^{\underline n})$ of square-summarable functions on the group  $\mathbb Z^{\underline n}$, which commute with the action of  
  $\mathbb Z^{\underline n}$ on this space.

A remarkable fact  is that a significant portion  of   Lusztig's argument  generalizes  to all discrete groups $\Pi$ instead of $\mathbb Z^{\underline n}$, where the algebra $C^\ast(\Pi)$ of bounded operators on $l_2(\Pi)$ is regarded as  the algebra of  continuous  functions on a
 "non-commutative space" dual to $\Pi$ (that is the actual  space, namely that of   of homomorphisms $\Pi\to \mathbb T$ for commutative $\Pi$.) 

This  allows a formulation of  what is called in    [Rosenberg($C^\ast$-algebras - positive  scalar) 1984] the  {\it strong Novikov {\color {red!50!black} Conjecture}}, the relevant  for us special case of which reads  as follows. 
\vspace{1mm}

{\it \color {blue}${\cal D}_{ \otimes C^\ast}$}-\textbf{{\color {red!50!black} Conjecture}.}  {\sf  If  a smooth closed orientable Riemannian spin  $n$-manifold $X$ for $n$ even  admits a continuous map $F$  to the classifying space {\sf B}$\Pi$ of a group $\Pi$, such that the   homology homomorphism $F_\ast$  sends the fundamental  homology class $[X]\in H_n(X;\mathbb R)$  to  {\it non-zero} element 
  $h\in H_n( ${\sf B}$\Pi;\mathbb R),$}
 then\vspace{1mm}

 {\it the Dirac operator  on $X$ twisted with some flat unitary Hilbert bundle over $X$ has 
 non-zero kernel.}

(Here "unitary" means that the monodromy action of $\pi_1(X) $ on the Hilbert fiber $\cal H$ of this bundle  is unitary and where an essential structure in this $\cal H$  is the action of the algebra $C^\ast(\Pi)$,  which commute with the action of $\pi_1(X) $.)\vspace{1mm}

This, if true, would  imply, according to  the Schroedinger-Lichnerowicz-Weitzenboeck formula, the spin case of  the conjecture stated in section \ref {spin index3}.
saying that

\hspace {25mm} {\it $X$ admits no metric with $Sc>0$.} \vspace{1mm}
 
Also  "Strong Novikov"  would imply, as it was  proved  by Rosenberg,  the validity of the \vspace {1mm}

{ \it \textbf {Zero in the Dirac Spectrum {\color {red!50!black} Conjecture}.}} {\sf Let $\tilde X$ be a complete {\it contractible}  Riemannian manifold  the  quotient of  which under the action of the  isometry group $iso(\tilde X)$ is {\it compact}.}

  {\it Then the  spectrum of the Dirac   operator $\tilde {\cal D}$  on $\tilde X$
contains zero, that is, for all $\varepsilon >0$,  there exist $L_2$-spinors $\tilde s$ on $\tilde X$,   such that 
$$||{\tilde {\cal D}(\tilde s}||\leq \varepsilon ||\tilde s ||. $$}

This,    confronted with the Schroedinger-Lichnerowicz-Weitzenboeck formula, would show that $\tilde X$ can't have $Sc>0$.
\vspace{1mm}

 {\it \color {magenta!19!black} Are we to  Believe in these {\color {red!50!black} Conjecture}s}?  A version of the    {\sl Strong Novikov conjecture} for  a rather  general class of groups, namely those which {\it admit discrete  isometric actions on spaces with non-positive sectional curvatures}, was proven    by Alexander Mishchenko in 1974.
 
Albeit this has been  generalized since 1974 to  many other  classes of   groups $\Pi$ and/or  representatives   $h\in   H_n( ${\sf B}$\Pi;\mathbb R)$, (most recent results   and references  can be found in  [GWY (Novikov) 2019]) the sad truth is that 
 one has a poor understanding  of  what these  classes actually are, how much they overlap and  what part of the world of   groups they fairly represent.

At the moment, there is no   basis for believing in this conjecture  and there is no idea where to look for a counterexample either.\footnote{ Geometrically most complicated groups are those which represent  one way or another  universal Turing machines;  a  group,  the  $k$-dimensional   homology (L-theory?)  of which, say for $k=3$, models such a "random" machine, would   be a good candidate for a counterexample. }

 \vspace {1mm}
 
  The following is a more  geometric version of the above conjecture.
 
  \vspace {1mm}
 {\it \color {blue} Coarse $\cal D$-Spectrum {\color {red!50!black} Conjecture}.} {\sf Let $\hat X$ be a  complete {\it uniformly contractible} Riemannian manifold, i.e. there exists a function $R(r)\geq r$, such that the  ball 
$B_{\hat x}(r)\subset \hat X$, $x\in X$, of radius $r$ is contractible in the concentric ball $B_{\hat x}(R(r))$  for all $\hat x\in \hat X$ and all radii $r>0$.}

 {\it Then the  spectrum of the Dirac     operator on $\hat X$
contains zero.}\vspace {1mm}

This conjecture, as it stands,  must be,  in view of [DRW(flexible) 2003],  {\it false},  but finding a counterexample becomes harder if  we require the bounds $vol (B_{\hat x}(r)) \leq \exp r$ for all  $\hat x\in \hat X$ and  $r>0$.\footnote {See [F-W(zero-in-the-spectrum) 1999] for what is known about the similar conjecture by John Lott for the DeRham-Hodge .}

%D-spec

And although  this conjecture remains unsettled  for $n=dim(X)\geq 4$, its  significant  corollary -- 

{\it non-existence of complete uniformly contractible Riemannian $n$-manifolds with positive scalar curvatures}  

\hspace{-6mm} was recently proved  for n=4 and 5  by means of torical symmetrization of stable $\mu$-bubbles,\footnote  {See  [Chodosh-Li(bubbles) 2020]  and [G(aspherical) 2020].}

 %%%%%%%%%%%%%%%%%%%
 
 \subsubsection{\color {blue}  Almost Flat Bundles and  $\bigotimes_\varepsilon$-Twist Principle}\label {twist principle3}

 %%%%%%%%%%%%%%%%%%

Let us  recall Dirac operators  twisted with  {\it  almost flat  unitary bundles} and construction of such bundles over     {\it profinitely hyperspherical} manifolds  such as  $n$-tori, for example.
\vspace {1mm}
 
 Let $X$ be a Riemannian manifold and $L=(L, \nabla)$ be   a complex  vector bundle $L$ with unitary connection.
  If the curvature  of $L$  is $\varepsilon$-close to zero, 
 $$||\mathcal R_L||\leq \varepsilon,$$
 then, locally, $L$ looks, approximately as the flat bundle $X\times \mathbb C^r$, $r=rank_\mathbb C( L)$,
 and the Dirac   twisted with $L$, denoted $\mathcal D_{\otimes L}$,{ that acts on the spinors with values in $L$, is locally approximately 
 equal  to the direct sum  $\underset {r}{\underbrace { \mathcal D \oplus ...\oplus \mathcal D}}.$
 
 It follows that if $Sc(X)\geq \sigma>0$ and if $\varepsilon$ is much smaller than $\sigma$, then by the (obvious) continuity of the Schroedinger-Lichnerowicz-Weitzenboeck formula,
this twisted Dirac  operator has trivial kernel, $ker (\mathcal D_{\otimes L})=0$  and, accordingly,
$$ind  (\mathcal D^+_{\otimes L})=0, \footnote{Here we assume that  $n=dim(X)$ is {\it even}, which makes  $\mathcal D$ {\it split} as 
$ \mathcal D= \mathcal D^+ \oplus \mathcal D^-$,  such that  $ind  (\mathcal D^+)=-ind( \mathcal D^-)$, see section \ref {Dirac4}.}$$
 where,  by the Atiyah-Singer  index theorem, this index is equal to a certain topological invariant
 $$ind  (\mathcal D^+_{\otimes L})=\hat \alpha(X,L).$$
 
 For instance, if $X$ is  an even dimensional topological torus, and if  the top Chern class of $L$ doesn't vanish, $c_m(L)\neq 0 $ for $m=\frac {dim(X)}{2}$,  then $\alpha(X,L)\neq 0$ as well.
 
 On the other hand,  given a Riemannian metric $g$ on the torus  $\mathbb T^n$, $n=2m$, 
 and $\varepsilon>0$,  \vspace{1mm}
 
 {\color {blue!40!black}{\sl there exists
  a finite covering $ \tilde {\mathbb T^n}$ of  the torus,  which admits an $\varepsilon$-flat  vector  bundle $\tilde L\to \tilde  {\mathbb T^n}$ of $\mathbb C$-rank $r=m=\frac {n}{2}$ with  $c_m(L)\neq 0 $,}  \vspace{1mm}}
  
\hspace {-6mm}   where the "flatness"  of $\tilde L$, that is the norm of  the curvature $\mathcal R_{\tilde L}$ regarded as a 2-form with the values in  the Lie algebra
   of the unitary group $U(r)$, $r=rank_\mathbb C(\tilde L)$, is measured  with  the lift $\tilde g$ of  the metric $g$ to $\tilde  {\mathbb T^n}$. \vspace {2mm}
  
{\color {red!20!black} {\sf Indeed, let $\hat L\to \mathbb R^n$, $n=2m$, be a vector  bundle with a unitary connection, such that $\hat L$  is isomorphic (together with it connection)  at infinity to the trivial bundle  and such that $c_m(\hat L)\neq 0$, where such an $\hat L$ may be induced by a  map $\mathbb R^n\to S^n$,  which is    constant at infinity and  has degree one,  from a bundle $\underline L\to S^n$ with $c_m(\underline L)\neq 0$. 
  
  Let $\hat L_\varepsilon$ be the bundle induced from $\hat L$ by the scaling map 
 $x\mapsto \varepsilon x$, $x\in \mathbb R^n$. 
 Clearly, the curvature of  $\hat L_\varepsilon$ tends to 0 as  $\varepsilon\to 0$.
 
Since the  finite coverings $\tilde { \mathbb T^n}$ of the torus  converge to the universal covering $\mathbb R^n\to  \mathbb T^n$  this   $\hat L_\varepsilon$ can be transplanted to a bundle $\tilde L_\varepsilon \to \tilde  {\mathbb T^n}$  over a sufficiently large finite covering $\tilde { \mathbb T^n}$ of the torus, where  the top Chern number  remains unchanged and where the curvature of $\tilde L$ with respect to the flat metric on 
$\tilde { \mathbb T^n}$  can be assumed as  small as you wish, say $\leq \epsilon$. 

But then this very curvature with respect to the lift $\tilde g$ of a given Riemannian metric $g$ on  $ \mathbb T^n$ also will be small, namely $\leq const_g\epsilon$ and our claim follows.}}\footnote{Why do we need  {\sf {\large \color {magenta!40!black} twelve lines} to express, not even fully at that,   so an obvious idea?} Is it due to an imperfection of our mathematical language or it is something about our mind that makes  instantaneous images of  structurally protracted  objects?  Probably both, where the latter  depends on  the {\it parallel processing} in the human {\it subliminal} mind, which can't be well represented by any  sequentially structured  language  that follows our  {\it conscious} mind and where besides "{\it parallel}" there are many other properties of  "{\it subliminal}" hidden from our    conscious mind eye.}
 \vspace{1mm}
  
 With this,  we obtain  \vspace{1mm}
 
 {\it one of the (many)  proofs of nonexistence of metrics $g$ with $Sc(g)>0$ on tori.}
 \vspace{3mm}
  
{\it \color {blue!20!black} \textbf{Seemingly Technical Conceptual Remark}.}  The above rough qualitative argument admits a finer quantitative version, which depends on the twisted Schroedinger-Lichnerowicz-Weitzenboeck formula 
$$ \mbox { ${\cal D}_{\otimes L}^2=\nabla^2_{\otimes L} + \frac{1}{4}Sc(X) + {\cal R}_{\otimes L}$},$$
where ${\cal R}_{\otimes L}$   is an operator  on twisted spinors, i.e. on the bundle $\mathbb S\otimes L$,  associated with the curvature of $L$
and where an essential feature of ${\cal R}_{\otimes L}$ is a  bound on its norm  by the {it  norm}
$||\mathcal R_L||$  of the curvature  $\mathcal R_L$ of $L$, with a constant {\it \color {blue!60!black} independent} of the rank of $L$.

Thus, for instance, the above proof of nonexistence of metrics $g$ with $Sc(g)>0$ on tori, that  was  performed with
the twisted  Dirac  $\mathcal D_{\otimes \tilde L}$ over a {\it finite covering} $\tilde X$   of our torical $X$, can be brought { \it back to }  $X$ by pushing forward $\tilde L$ from the    $\tilde X$  to $X$, where 
this  push forward bundle  $(\tilde L)_\ast \to X$  has 
$$rank(\tilde L)_\ast = N\cdot rank (\tilde L)$$
for  $N$ being the number of sheets of the covering.

(The lift of   $(\tilde L)_\ast $ to $\tilde X$  is the Whitney's sum of $N$-bundles obtained from $\tilde L$ by the   deck  transformations of $\tilde L$.) \vspace {1mm}

This property of $\mathcal R_{\otimes L}$, in conjunction with the shape of the Atiyah-Singer  index formula, for the  Dirac    operator twisted with Whitney's  $N$-multiples

=
$$ L\oplus ... \oplus L=   \underset {N}{\underbrace {L\oplus ... \oplus L}},$$ 
which implies that in the relevant cases 
$$ind  (\mathcal D^+_{\otimes ({{L\oplus ... \oplus L}})})  = \alpha(X,   L\oplus ... \oplus L )=N\cdot \hat \alpha(X, L)+O(1),$$
allows $N\to \infty$ and even $N=\infty$ in a suitable sense, e.g.  in the context of infinite coverings  and/or of $C^\ast$-algebras as was mentioned  in  the previous section.

\vspace {1mm}

What is also  crucial, is that  
{\sf twisting with almost flat bundles } is a {\it \color {blue!60!black} functorial}
 operation, where this functoriality  yields  the following.\vspace {1mm}

{\color {blue}  $\bigotimes_\varepsilon$}-\textbf {Twist Principle}.} {\sl All} (known) {\it  arguments with Dirac operators for  {\color  {gray!60!black} non-existence } of metrics with $Sc\geq \sigma >0$   under  specific  topological conditions on $X$  can
 be} (more or less) {\it  automatically  transformed to {\color  {blue!60!black} inequalities} between $\sigma$  and certain {\color {blue!60!black}  geometric invariants} of $X$  defined via  {\color {blue!60!black}$\varepsilon$-flat bundles} over $X$.} \vspace {1mm}

 {\color {blue}  $\bigotimes_\varepsilon$}-{\it \textbf {Problem.}} {\sf Can one turn {\color {blue}  $\bigotimes_\varepsilon$}-Twist Principle} to a  {\color {blue}  $\bigotimes_\varepsilon$}-{\it theorem?}  \vspace {1mm}
 
 At the present moment, an application of  the {\color {blue}$\bigotimes_\varepsilon$-principle} necessitates  tracking {\it step by step}, let it be in a purely mechanical/algorithmic fashion,  a particular Dirac theoretic argument, rather than a direct  application of this principle to {\it the conclusion}  of such an argument. 
 
 What, apparently, happens here is that the true outcomes of  Dirac  operator proofs are {\it not}
 the  geometric theorems they assert, but certain linearized/hilbertized generalization(s) of these, possibly,  in the spirit of Connes' non-commutative geometry. 
  
To understand what goes on, one needs, for example,  to reformulate (reprove?)  Llarull's, Min-Oo's  and Goette-Semmelmann's inequalities in such a  "linearized" manner.\footnote {A promising approach is suggested by the convept of 
 {\it quantitative K-theory}, which was  successfully used in [Guo-Xie-Yu(quantitative K-theory) 2020] for a new proof of the $\frac {\2pi}{n}$-bounds in the width of  Riemannian bands
with $Sc\geq n(n-1)$. 

This theory encodes the geometric information on the underlying Riemannian manifold $X$ in term of the {\it propagation radius} $r$
 of operators in the {\it Roe translation algebra} that correspond  to  (linear combinations)  of {\it $r$-translations} of $X$  that are   self mappings $a:X\to X$ with $dist(a(x),x)\leq r$.
   
   This faithfully reflects the  {\it distance  geometry} of $X$, but the quantitative  $K$-theory, as it stands now,    can't  adequately capture   the area geometry; conceivably this can be achieved  by incorporation   ideas from 
   Cecchini's long neck paper into this theory.}
 
  \vspace {1mm} 
  
  {\it  Twists with non-Unitary Bundles.}  Available (rather   limited) results concerning  scalar curvature  geometry of manifolds $X$, which  support  almost flat non-unitary bundles and of (global spaces of possibly) non-linear   fibrations with almost flat connections  over $X$,  are discussed in section \ref {{unitarization6}}.

 \vspace {1mm}

{\it \large Flat or Almost Flat?}  Lusztig's approach to the Novikov conjecture via the signature operators twisted   with (families of)  {\it finite dimensional non-unitary flat} bundles was  superseded, starting with  the work by Mishchenko and Kasparov, by more general index theorems,  for  {\it infinite dimensional  flat unitary} bundles.

Then it was observed in  [GL(spin) 1980] and  proven in a general form  in  [Rosenberg($C^\ast$-algebras - positive  scalar) 1984])       that all these results can be  transformed to the corresponding statements about  Dirac operators on spin manifolds,  thus providing   obstructions
to $Sc>0$ essentially  for the same kind  of  manifolds $X$, where  the generalized signature  theorems were established.

Besides   following topology,  the geometry of the scalar curvature suggested a quantitive version of these topological theorems by allowing twisted Dirac and signature operators with {\it non-flat vector bundles}  with {\it controllably small curvatures}, thus providing geometric information on $X$ with $Sc\geq \sigma >0$,   which complements the  information on pure topology of $X$.

At the present moment, there are two groups  of papers on twisted (sometimes untwisted)  Dirac operators on  manifolds with $Sc>\sigma $. 

The first and a most abundant one goes along with the 
work on the Novikov conjecture, where  it is framed into the $KK$-theoretic formalism. 

A notable  achievement of   this is

 {\it  Alain Connes' topological  obstruction  for leaf-wise metrics with $Sc>0$ 
 
 on  foliations}, 
 
 \hspace {-3mm}  where 
 
 {\sf a geometric shortcut through  the $KK$-formalism of Connes' proof is unavailable 
 
 at the present moment. }\vspace {1mm}

Another direction is a  geometrically oriented  one, where we are not so much concerned with the 
$K$-theory  of the $C^\ast$-algebras of  fundamental groups $\pi_1(X)$, but with  geometric constraints on $X$ implied by the inequality $Sc(X)\geq \sigma$.

This goes close to what happens in the papers inspired by the general relativity,  where one is concerned  
with  specified (and rather special, e.g. asymptotically flat)   geometries at infinity of complete Riemannian manifolds and where 
 one  plays, following Witten  and Min-Oo,  with  Dirac  operators, which are asymptotically adapted at infinity  to such   geometries.   (In this context, the Schoen-Yau and the  related 
methods relying  of the {\it mean curvature flows}  are also used.)\vspace {1mm} 
 
In the present paper, we are primarily concerned with  {\sf \color    {black!60!green} geometry} of  manifolds, while   {\color {blue} topology} is  confined to {\it  an  auxiliary,} let it be  irreplaceable, role.

%%%%%%%%%%%%%%%%%%%%%%

\subsubsection  {\color {blue}Relative Index  of Dirac Operators on Complete 
 Manifolds} \label {relative3}

%%%%%%%%%%%%%%%%%%%%%%

 Most (probably, not all)  bounds  on  the scalar curvature of {\it closed} Riemannian  manifolds derived  with twisted Dirac operators  $\mathcal D_{\otimes L}$ have their counterparts for {\it complete} manifolds $X$, where  one  uses a relative version  of the  Atiyah-Singer theorem  for  {\it pairs of Dirac operators which agree at infinity} \footnote {See [GL(complete) 1983], [Bunke(relative index) 1992],  [Roe(coarse geometry)  1996]), and more recent papers [Zhang(Area Decreasing) 2020],  [Cecchini(long neck) 2020],  [Cecchini-Zeidler(generalized Callias) 2021],  
[Cecchini-Zeidler(Scalar\&mean) 2021]. } 
  the simplest and 
 the most relevant case of  this  theorem applies to vector bundles 
  $L\to X$  with  unitary connections which are {\it flat trivial at infinity}.
  
  In this case the pair in question is $(\cal D_{\otimes L},   \cal D_{\otimes |L|}$),  where $|L|$ denotes the trivial flat bundle $X\times \mathbb C^k\to X$ for $k=rank_{\mathbb C}(L)$, which comes along with an isometric connection preserving isomorphism between  $L$  and $|L|$  outside  a compact subset in $X$.
   \vspace {1mm}
  
{\it  \color {blue} $f^\ast$-Example}. Let  $f:X\to S^n$ be a smooth map which is locally constant at infinity (i.e. outside a compact subset) and let $\underline L\to S^n$ be a bundle with a unitary connection on $S^n$. 

Then the pullback bundle 
$f^\ast(\underline L)\to X$ is an instance of such an $L$.   \vspace {1mm}
  
 The relative index theorem, similarly to its absolute counterpart implies  that  if  the scalar curvature of $X$ is {\it uniformly positive (i.e. $Sc\geq \sigma>0$) at infinity} and if  \vspace {1mm}

 {\sf a certain  topological invariant, call it  
  {\it $\hat \alpha (X, L)$,\footnote{See section \ref{spin index3} for the definition of this invariant. } {\it doesn't vanish}, then either $X$ admits a non-zero  (untwisted) harmonic $L_2$-spinor $s$ on $X$, that is a solution of $\mathcal D(s)=0$, or there is a   non-zero  $L$-twisted harmonic $L_2$-spinor on $X$.}\footnote{ If we don't assume that $Sc(X)$ is uniformly positive at infinity, then one can only claim the existence of   either non-zero untwisted or non-zero twisted  {\it almost harmonic} $L_2$-spinors, i.e. satisfying 
  $\int_X\mathcal D^2(s)dx   \leq \varepsilon \int_X||s(x)||^2$  or $\int_X\mathcal D_{\otimes L}^2(s)dx \leq \varepsilon  \int_X||s(x)||^2$,
  for arbitrarily  small $\varepsilon>0$.   }} \vspace {1mm}

{\it \color {blue}$f^\ast$-Sub-Example.} {\sl Let  $L=f^\ast(\underline L)$ be as  in the{ \color {blue} $f^\ast$-example}, let  where   $n=dim(X)$ is even, and let the  bundle  $\underline L\to S^n$ has  non-zero top Chern class (e.g. $\underline L$ is the bundle of spinors on the sphere, $\underline L=\mathbb S_+(S^n)$).
   If the map $f:X\to S^n$ has {\it non-zero degree},   then  $\hat \alpha (X, L)\neq 0.$}
\vspace {1mm}

Finally,  since  the twisted Schroedinger-Lichnerowicz-Weitzenboeck formula (obviously)  applies to $L_2$-spinors, one obtains, for example,  as  an application of the {\color {blue}  $\bigotimes_\varepsilon$-Twist Principle}
the following relative version of 
the Lichnerowicz'  theorem for $k$-dimensional manifolds  from section \ref {spin index3},  that, let us  remind it, says  that 

\hspace {20mm} {\sl $\hat A[X]\neq 0\Rightarrow Sc(X)\ngtr 0$  for closed  spin manifolds $X$.}
\vspace {1mm}

{\it If a complete Riemannian   orientable spin manifolds $X$ (of dimension $n+4k$) admits a proper $\lambda$ -Lipschitz map $f:X\to \mathbb R^n$  for some $\lambda<\infty$, then the pullbacks of generic points $y\in \mathbb R^n$ satisfy
$\hat A[f^{-1}(y)]=0$.}\vspace{1mm}

This, in the case $dim(X)=n$, shows that \vspace{1mm}

{\it the  existence of  proper Lipschitz map $X\to \mathbb R^n$ implies that $\inf_xSc(X,x)\leq 0$.}\footnote 
 {This has a variety  of   generalizations and  applications, (see  e.g.  [GL(spin) 1980], [GL(complete) 1983], 
 [Roe(partial vanishing) 2012]  and references therein), such as  {\sf non-existence of metrics with $Sc>0$ on tori}. }
 
 Moreover, 
 
 {\it it follows from {\color {blue} Zhang's theorem} stated below,  that, in fact, $\inf_xSc(X,x)< 0$.}

\vspace {1mm}

\vspace{1mm}

\vspace {1mm}

 The relative index theorem combined with the  linear-algebraic analysis of
the $L$-curvature term $\mathcal R_{\otimes L}$ in the   twisted Schroedinger-Lichnerowicz-Weitzenboeck formula due to  Llarull , Min-Oo,  Goette-Semmelmann  and Listing    allows an extension of their inequalities  from  compact manifolds) to  non-compact complete manifolds.

For instance, 

   \vspace{1mm}

{\it \color {blue} {\huge  $\star$}}{\it   If a complete Riemannian   orientable spin manifolds $X$ (of dimension $n+4k$) with $Sc(X)>n(n-1)$  admis a  locally constant at infinity   $1$ -Lipschitz map $f:X\to S^n$,  then the pullbacks of generic points $y\in \mathbb S^n$ satisfy
$\hat A[f^{-1}(y)]=0$.}\footnote {
See  [Llarull(sharp estimates)  1998]  and also  sections  \ref {area extremality3}},  \ref {sharp algebraic4}, \ref {Llarull4}.}
\vspace {1mm}
 
%%%%%%%%%%

 {\it \color {blue}  Zhang's Extension of the Relative Index Theorem with Applications to maps 
 $X\to S^n$ .}   The above stated  relative index theorem      needs    {\it uniform} positivity of the  scalar curvature of $X$  at infinity,   i.e. the bound $Sc(X)\geq \sigma>0$.
 
This uniformity  condition was removed in  [Zhang(Area Decreasing)  2020] by using a small zero order  perturbation of the relevant  twisted  Dirac  at infinity making the resulting  positive at infinity and thus,  proving the following theorem.\vspace {1mm}

 {\sf  Let a complete  orientable  spin $n$-manifold $X$ of non-negative scalar curvature,
 $Sc(X)\geq 0$ and let $X$ 
 admit a smooth  {\it area decreasing} locally constant at infinity  (i.e. outside a compact subset)  map $f:X\to S^n$ of}      {\it  non-zero degree. } 
 
 Then \vspace {1mm}
  
{\huge \color {blue} $\star \star$}  {\it the scalar curvature of $X$ on the  support of the differential  of $f$}  (where $df\neq 0$) {\it satisfies:
  $$ \inf_{x\in supp(df)} Sc(X, x)\leq n(n-1),\footnote{This  also follows from Cecchini's long neck theorem stated in the next section.}$$
and if $n$ is even, then  
 $$\inf_{x\in supp(df)} Sc(X, x)< n(n-1),$$  
   unless  $X$ is compact and $f$  is an  isometry. }

 \vspace {1mm}

{\it Remark }  It remains unclear, even for compact $X$,  if the spin condition is essential, 
but the completeness condition  can be significantly relaxed as we shall explain in the next section.

%%%%%%%%%%%%%%%%%%%%%%

 \subsubsection {\color {blue} Roe's Translation  Algebra, Dirac Operators on Complete Manifolds with Boundaries and Cecchini's  Long Neck Theorem  for Non-Complete manifolds }\label {Roe3}

%%%%%%%%%%%%%%%%%%

 $C^\ast$-algebras bring forth   the following   interesting perspective on  {\it  coarse geometry}  of non-compact spaces  proposed by John Roe  following  Alain Connes' idea of non-commutative geometry of foliations.

 Given a metric space $\Xi$, e.g. a discrete group   with a word metric,  let ${ \cal T}=Tra(\Xi)$   be the semigroup of translations of $M$ that are  maps $\tau:\Xi\to \Xi$, such that 
 $$\sup_{\xi\in \Xi} dist(\xi, \tau(\xi)) <\infty.$$
 
The (reduced) Roe $C^\ast$-algebra $R^\ast(\Xi)$ is a certain completion of the semigroup  algebra $\mathbb C[{ \cal T}].$
 For instance if  $\Xi$ is a    group with a word metric for which, say  the left action of  $\Xi$ on itself is isometric, then the right actions lie  in  ${ \cal T}$ and $R^\ast(\Xi)$ is equal to the (reduced)  algebra $C^\ast(\Xi)$.\footnote{"Reduced"  refers to a minor technicality not relevant at the moment.  A more serious problem -- this is not joke --  is impossibility of   definition of "right" and "left"   
without an appeal to  violation  of mirror symmetry  by weak interactions.}

 Using this algebra, Roe proves in  [Roe(coarse geometry) 1996], (also see  [Higson(cobordism invariance) 1991], [Roe(partial vanishing) 2012])  a {\it partitioned index theorem},  which implies, for example, that\vspace{1mm}.
 
\PencilRight\hspace{1mm} {\it the toric  half cylinder manifold $X=\mathbb T^{n-1}\times \mathbb R_+$ admits no complete Riemannian metric with $Sc\geq \sigma>0$}.\footnote{I must admit I haven't  fully understood Roe's argument.}\vspace{1mm}

  \vspace{1mm} 
 
 Nowadays \PencilRight  \hspace {1mm}  can be  proved  with the  techniques of minimal hypersurfaces and of stable $\mu$-bubbles, (sections  \ref {bands3},  \ref {quadratic3}) 
   as well as withDirac theoretic techniques with potentials developed by  Zeidler and by  Cecchini and by the
Guo-Xie-Yu in the  framework of the quantitative K-theory, (see below) 
 where these techniques yield not only the bound $\inf_xSc(X,x)\leq 0$ but a {\it quadratic decay}  
of  the scalar curvature
on $\mathbb T^{\-1}\times \mathbb R_+$.

 Also notice in this regard that if $X$ is sufficiently "thick at infinity", then  \PencilRight  \hspace {1mm}  follows by a simple  argument with twisted Dirac operators and the standard  bound on the number of small eigenvalues  in  the spectrum of  the Laplace (or directly of the  Dirac)  operator in vicinity of $\partial X$,   which applies to all manifolds with  boundaries and which yields, in particular, (see section \ref {amenable4})  the following.

 \vspace{1mm}

\NibRight \hspace {1mm}   Let $X$ be a complete  oriented  Riemannian spin   $n$-manifold {\it with compact boundary}, such that
 
 {\sf there exists  a sequence of  smooth area decreasing maps $f_i:X\to S^n$, which are constant in a  (fixed) neighbourhood   $V\subset X$ of the boundary $\partial X$ as well as  away from   compact  subsets $W_i\subset V$,   and such that   $$deg(f_i)\underset {i\to \infty}\to \infty.$$}

{\it Then the scalar curvature of $X$ satisfies
$$\inf_{x\in X} Sc(X, x) \leq n(n-1).$$ }

  \vspace{1mm}

{\it \color {blue} Quantitative K-theory and Long Neck Principle.}  It seems that most (all?) results  for  {\it complete} Riemannian manifolds with  $Sc\geq \sigma$  have their  counterparts for manifolds $X$ with boundaries insofar as  this  concerns the part of $X$   that lies far from the boundary $\partial X$. 

Definite results  in  this regard  were  recently obtained  by  Hao Guo, Zhizhang Xie  and  Guoliang Yu
who, if I  understand this correctly, developed  a quantitative version of Roe's theory, and also by Rudolf Zeidler and Simone Cecchini who obtained   index theorems  for Dirac  operators with potentials on manifolds with boundaries.\footnote {See [Cecchini(long neck) 2020], [Guo-Xie-Yu(quantitative K-theory) 2020],
[Zeidler(bands)  2019],  
[Zeidler(width)  2020], [Cecchini-Zeidler(generalized Callias) 2021], [Cecchini-Zeidler(scalar\&mean) 2021].} 

Here is an instance of some of new results.\vspace{1mm}

{\color {blue}\textbf  {Cecchini's Bound on Hyperspherical Radii  of Long Neck manifolds}.}\footnote{This was a response to   a  question from  an earlier version of this manuscript.} {\sf Let $X$ be a compact $n$-dimensional orientable {\it spin}  Riemannian  manifolds with a boundary, let  $Sc(X)\geq \sigma_0>0$ and let $f:X\to S^n(R)$ be a smooth area decreasing  map, which  is   locally constant in a neighbourhood of the boundary $\partial X\subset X$  and which have {\it nonzero degree.}

Let the scalar curvature  of the support of the differential of $f$ be bounded from  below by
$\sigma$ (where typically   but not necessarily $\sigma  \geq \sigma_0$),
$$Sc(X,x)\geq \sigma, \mbox { } x\in supp(df).$$}

{\it If $f$ satisfies the following "long neck condition", 
$$dist(supp (df),\partial X)\geq  \pi\sqrt{\frac {n-1}{n\sigma_0}},$$
then the radius of the sphere $S^n(R$) is bounded, similarly to the case  of complete $X$, by  
$$R\leq \sqrt \frac {n(n-1)}{\sigma},$$}
{\sf where in the case of {\it odd} $n$ one  additionally assumes (this is, probably, redundant) that
$f$ is {\it constant} (not just locally constant)  in the neighbourhood of $\partial X\subset X$.}

   \vspace{1mm} \vspace{1mm}
  
  {\it Question.} What are (preferably sharp)  long neck counterparts of the  {\it  spin-area convex}  and {\it spherical   trace area} {\sf extremality theorems} from section \ref {area extremality3}

%%%%%%%%%%%%%%%%%%%%%%%%%
 \subsection{\color {blue} Foliations With Positive Scalar Curvature}  \label {foliations3}
  %%%%%%%%%%%%%%%%%%%%%%%%%%
  
  According to the philosophy (supported by a score of   theorems) of Alain  Connes much of the geometry and topology of manifolds with discrete group actions, notably, those concerned with index theorems for    Galois actions of fundamental groups on universal  coverings of compact manifolds, can be extended to foliations.   
  
  In particular, Connes shows in [Connes(cyclic cohomology-foliation) 1986] that  compact manifolds $X$  which carry   foliations $\mathscr L$  with  leaf-wise Riemannian metrics  with positive scalar curvatures behave in many respects as manifold which themselves admit such metrics.  
  
 For instance,  \vspace{1mm}

{\color {blue} \huge $\star$}   {\sl if  $\mathscr L$ is spin, i.e.  the tangent (sub)bundle  $T(\mathscr L)\subset T(X)$  of such an  $\mathscr L$ is {\it spin}  then,  
  
  by Connes'  theorem,  $\hat A[X]=0$.}  \vspace{1mm}
  
This generalises Lichnerowicz' theorem  from section \ref {spin index3(} for oriented  spin  manifolds of dimensions $n=4k$ , where, recall,   $\hat A[X]$  is the value of a certain rational  polynomial  $\hat A(p_i)$ in the Pontryagin classes $p_i\in H^{4i}(X:\mathbb Z)$ (see section \ref {Dirac4})  on the fundamental homology class $[X]\in H_n(X)$.\footnote {By definition, the values of the $p_i$-monomials $P_d= \smile_jp_{i_j}\in H^d(X)$,  $4\sum_ii_j=d$, on $[X]$ equals zero for all $d\neq n.$ }
  
  \vspace{1mm}

 In fact, the full  Connes'  theorem implies among other things   \vspace{1mm}

 {\sl  vanishing of the  $\smile$-products of the $\hat A$-genus $\hat A(p_i)$,  $j=0,1,...,k=\frac {n}{4}$,   with  all  polynomials in the Pontryagin classes  
 of the  "normal" bundle 
 $T^\perp(\mathscr L)=T(X)/T(\mathscr L)$, in the case where $\mathscr L$ is  spin.}

 \vspace{1mm}
 
 Connes' argument, which  relies on  Connes-Scandalis  {\it longitudinal index theorem for foliations}), 
 delivers a  non-zero almost harmonic spinor on some leaf of $\mathscr L$  and 
  an alternative and simpler proof of the existence of such spinors under suitable conditions     was  given  in   [Bern-Heit(enlargeability-foliations) 2018], 
  where $\mathscr L$, besides being spin, is required to have  {\it Hausdorff homotopy groupoid}.{\footnote {One 
      finds a   helpful explanation  of the meaning this condition in [Con(foliation) 1983] and in  the lectures   [Meinrenken(lectures) 2017]. }

 Another simplified proof of (a part of)   Connes' theorem was also   suggested in
[Zhang(foliations) 2016], where 
  the {\sl manifold $X$  itself}, rather than  the tangent bundle $T(\mathscr L)$ is assumed  {\it spin} \footnote {In the ambience of Connes' arguments  [Connes(cyclic cohomology-foliation) 1986], these two spin conditions reduce one to another.}  
  and where the existence of  almost harmonic spinor is proven on some auxiliary manifolds associated with $X$. 
\vspace {1mm}
  
One can get more mileage from  the index theoretic arguments in these papers by applying the  {\color {blue} $\bigotimes_\varepsilon$-Twisting Principle} from section \ref {twist principle3}, but this needs  honest checking all   steps 
 in  the proofs in there. This was (partly)  done in  [Bern-Heit(enlargeability-foliations) 2018],  in 
[Zhang(foliations:enlargeability) 2018],  [Su(foliations) 2018] and [Su-Wang-Zhang(area decreasing foliations) 2021] in the context of the index theorems used by the authors in  their papers. \footnote{
  I recall going through  Connes' paper long time ago and observing in  \S $9\frac {2}{3}$ in [G(positive) 1996]).  that    Connes' argument  yields the following. \vspace {0.5 mm} 
 
 {\it Complete manifolds $X$ with infinite K-cowaist$_2$}  (called "K-area" in  [G(positive) 1996]) , {\it e.g. $\mathbb R^n$,    carry no spin foliations, where the  induced     Riemannian metrics in the leaves satisfy    $Sc\geq \sigma>0$,}  
 
 but my memory is uncertain at this point.} \vspace {1mm} ,

Here is  a geometric   conjecture in this regard.\footnote {This may  follow  by what is done in The techniques (results?) from [Zhang(foliations:enlargeability) 2018], [Cecchini(long neck) 2020], [GWY (Novikov) 2019] may be useful for settling this question.   } \vspace {1mm}

 \vspace {1mm} 
  
  {\it \color {blue} \textbf {Long Neck  Foliated   {\color {red!50!black} Conjecture}}}.
   {\sf Let $X$ be a compact oriented  $n$-dimensional Riemannian manifold with a boundary, let   $\mathscr L$  be a smooth  $m$-dimensional, $2\leq m\leq n$, foliation on $X$, such that 
   the induced Riemannian metrics on the leaves of $\mathscr L$ have positive scalar curvatures, 
   $$Sc(\mathscr L)>\sigma_0, $$
    let 
 $f: X\to S^n$ 
  be a smooth  map, which  is   locally constant in a neighbourhood of the boundary $\partial X\subset X$  and which have {\it nonzero degree}, let the scalar curvature of $L$ on the support of the differential of $f$ be bounded from  below by the trace norm of the second exterior power of the differential of $f$ on the tangent bundle of $\mathcal L$ as follows
  $$Sc(X, x\in supp(f))\geq 2 trace( \wedge^ 2 df_\mathcal L)(x)).$$}
   
   {\it Then the leaf-wise distance $ D=  dist_{\mathcal L} (supp(df),\partial X)$\footnote {This $D$, that is the infimum of  the length of curves in the leaves between the intersections of these leaves with  $supp(df)\subset X$ and with $\partial X\subset X$, isIn general, greater than  the   distance $d=dist_X( (supp(df),\partial X)$ For instance, if  no leaf   intersects both subsets,  $supp(df)\subset X$ and $\partial X\subset X$, than $D=\infty.$}
   is bounded by some universal  function of $\sigma_0$,
$$d\leq \theta(\sigma_0).\footnote {Combined arguments from [Su-Wang-Zhang(area decreasing foliations) 2021]and 
 [Cecchini(long neck) 2020]  may lead to the proof if either $X$ or $\mathcal L$ is spin and
  $\theta(\sigma_0)=  \pi\sqrt{\frac {n-1}{n\sigma_0}}.$}$$}
   
   \vspace {1mm}
   
   {\sl\color {blue!50!black} Non-Integrable Question.} {\sf Is there a "good"  bound on the scalar curvature of  non-integrable subbundles $\mathcal T\subset T(X)$ of rank $m$ (instead of the tangent subbundles $T(\mathcal L)$ of foliations $\mathcal L$)?}

   Here, $Sc(\mathcal T, x)$ is defined as the sum of the sectional curvatures of $X$ in an orthonormal frame of bi-vectors in 
 the space  $\mathcal T_x$,  and where, besides the scalar curvature, such an inequality must contain a non-integrability correction term.
 
 If this correction term is sufficiently small in the $C^1$-topology, then the above conjecture   could apply to   families of approximate  integral manifolds of  $\mathcal T$;  however, the   resulting bound  on $Sc(\mathcal T, x)$  seems very rough.
 
But  what  we look for  is a  sharp or a nearly sharp inequality  approaching      model examples, such  as the standard  codimension one (contact) subbundles on the odd dimensional  spheres and codimension  three  subbundles on the $(4k-1)$-spheres.  
 
Next,  we want to work   out  a concept of scalar curvature of  {\it sub-Riemannian}  ({it Carnot-Caratheodory})  manifolds and show, for instance, that (self-similar)  nilpotent Lie groups admit no such metrics quasi-isometric to the standard (self-similar)  ones.

  \vspace {1mm}

{\it Stable Complementation Question} [{\Large \color {blue}$\star_?$}].  {\sf Let $(X,g)$ be a (possibly non-complete) Riemannian $n$-manifold   with a smooth foliation, such that scalar curvature  of the induced metric on the leaves satisfies  $Sc\geq \sigma>0$.

Does the product of $X$ by a Euclidean space, 
$X\times \mathbb R^N$, admit an  {\it $\mathbb R^N$-invariant} Riemannian metric $\tilde g$, such that $Sc(\tilde g)\geq \sigma$ and the quotient  map  $(X\times \mathbb R^N, \tilde g)/ \mathbb R^N\to (X,g)$ is 1-Lipschitz, or, at least, $const_n$-Lipschitz}?

\vspace {1mm} 

(See   \S$ 1\frac {7}{8}$ in [G(positive) 1996] and  section \ref {foliations6},  \ref{Connes6},  \ref {Hermitian Connes6}  for partial results in this direction based on the  geometry of    {\it Connes'  fibrations}.)

   \vspace{1mm}

 Notice that even the complete  (positive) resolution of [{\Large \color {blue}$\star_?$}] wouldn't 
 yield the entire  
 Connes' vanishing  theorem  from  [Connes(cyclic cohomology-foliation) 1986], nor 
 would this   fully reveal   the geometry of  foliated Riemannian  manifolds $X$ with   scalar curvatures of the leaves bounded from below,  e.g.  an answer to the  following questions.

   \vspace{1mm}

1.  {\sf Do compact Riemannian   $n$-manifolds with constant curvature $-1$  admit  $k$-dimensional foliations, 
  $2\leq k\leq n-1$, such that the scalar curvatures of the induced Riemannian  metrics in the leaves are bounded from below by 
  $-\varepsilon$ for a given  $\varepsilon>0$}?\vspace{1mm}

 2. {\sf What would be a foliated version of the Ono-Davaux Spectral Inequality?}\vspace{1mm}
%%%%%%%%%%%%%%%%

\subsection {\color {blue} Scalar Curvature in Dimension 4 via the  Seiberg Witten Equation} \label {Seiberg3}

%%%%%%%%%%%%%%%%%%

The simplest examples of  4-manifolds  where non-existence of metrics with $Sc>$ follows from non-vanishing of Seiberg-Witten invariants are complex algebraic surfaces $X$ in $\mathbb CP^3$ of degrees $d\geq 3$. (If $d$ is even and these  $X$ are spin,  this also follows  from  Lichnerowicz' theorem from section \ref {spin index3(}.) 
 
In fact, it was shown by LeBrun (see  [Salamon(lectures) 1999] and references therein) that \vspace{1mm}
 
{\it no minimal {\sf (no lines with self-intersections one)} K\"ahler  surface $X$ admits a Riemannian  metric with $Sc>0$, unless $X$ is diffeomorphic to $\mathbb CP^2$ 
   or to a ruled surface }.\vspace{1mm}

Furthermore, LeBrun following Witten  shows  in   [LeBrun(Yamabe)  1999] that \vspace{1mm}

{\sf if such an  $X$ has {\it Kodaira dimension} 2, which is the case, for instance,  for the  algebraic   surfaces $X \subset \mathbb CP^3$ of degree $d\geq 5$, }
then  \vspace{1mm}

{\it the total squared  scalar curvature is bounded by the first  Chern number of $X$,
$$\int _XSc(X,x)^2dx\geq 32\pi^2c_2(X),$$
where, moreover this inequality is sharp. }\vspace{1mm}

\vspace {1mm}

Although  one doesn't  expect anything comparable to the Seiberg-Witten equations for $n=dim(x)>4$,  one wonders if some coupling between the  twisted Dirac  $\mathcal D_{\otimes L}$ and an energy like functional in the space of 
connections in  $L$ may be instrumental  in the study of  the scalar curvature of $X$  and lead to   bounds on  $\int _XSc(X,x)^{\frac{n}{2}}dx$ for a  manifolds $X$ of dimension  $n>4$ and, even better
 on   $\int _X|Sc_-(X,x)|^{\frac{n}{2}}dx$ for 
 $Sc_-(X,x)=\min  (Sc(X,x),0)$.\vspace{1mm}

For instance, 

{\sf Let  a closed orientable Riemannian $n$-manifold  $X$ admits a map of non-zero degree to a closed locally symmetric manifold   $\underline X$ with negative Ricci curvature, e.g. with constant negative curvature.} 

{\color {blue!50!black} \sl Does then  the scale invariant  integral of the negative part of the scalar curvature is bounded from below  as follows: 
 $$\int _X|Sc_-(g,x)|^{\frac{n}{2}}dx \geq \int _{\underline X}|Sc(\underline X,\underline x)|^{\frac{n}{2}}d\underline x?$$ }

(Three conjectures  related to this one are formulated in    section \ref {negative3}.)\vspace {1mm}

 {\it Question.}  {\sf What is the   Seiberg-Witten 4D-version   of geometric inequalities
 on manifolds with boundaries and manifolds with corners?}
 
 %%%%%%%%%%%%%%%
 
 \subsection {\color {blue} Topology and Geometry of Spaces of Metics with $Sc\geq \sigma$.}\label {spaces of metics3}

%%%%%%%%%%%%%%%

{\it Non-connectedness}  of the space of metrics with $Sc>0$ 
 starts with the following observation.

Let   a closed $n$-manifold $X$ be decomposed as $X_-\cup X_+$ where  $X_-$ and $ X_+$ are  smooth domains  ($n$-submanifolds) in $X$   with a common boundary $Y=\partial X_-=\partial  X_+$ and where $X_\mp$ are equal to  regular neighbourhoods of disjoint  polyhedral subsets  
$P_\mp\subset X$ of dimensions $n_\mp$ such that $n_-+n_+=n-1$.

If $n_\mp\leq n-2$, then, by an easy elementary argument,  both manifolds $X_-$ and $X_+$ admit Riemannian metrics, say $g_\mp$, such that  

\vspace {1mm}

{\sf the  restrictions of these $g_\mp$ to $Y$, call them $h_\mp$, both  have   {\it positive} scalar curvatures.}  

\vspace {1mm}

{\sl And if $X$ admits no metric with positive scalar curvature, e.g. if $X$ is homeomorphic to the $n$-torus or to product of two Kummer surfaces, then $h_-$ and $h_+$ can't be joined by a homotopy of metrics with positive scalar curvatures.} 
\vspace {1mm}

Indeed, such a homotopy, $h_t$, $t\in[-1,+1]$  could be easily transformed to a metric on the cylinder $Y\times  [-1,+1]$ with positive scalar curvature  and with relatively flat boundaries isometric to $(Y,h_-)$ and $(Y,h_+)$, which would then lead in obvious way to a metric on $X=X_-\cup Y\times  [-1,+1] \cup X_+$ with $Sc>0$ as well.

\vspace {1mm}

The first case of disconnectedness of spaces of metrics with $Sc>0$  goes back to    Hitchin's paper   [Hitchin(spinors)1974] , where it is shown, among many  other things,   that

{\sl the sphere 
 $S^n$,
  $n =8k, 8k+1$, admits a diffeomorphism  $\phi:S^n\to S^n$, such that the pullback $g_1=\phi^\ast (g_0)$ 
   of the standard metric $g_0$  can't be joined with $g_0$ by a homotopy $g_t$ with $Sc(g_t)>0$, }
   
   \hspace {-6mm}where
     appropriate   $\phi$  are those  for which   
     
     {\sf  the exotic spheres obtained by gluing  pairs of $(n+1)$-balls  across their boundaries according to $\phi$ have  non-vanishing $\hat\alpha$-invariants} (see section \ref {spin index3(})  and where the proof relies on the index theorem for families of Dirac  operators.
 Similarly, Hitchin finds non-contractible loops in the spaces of metrics $g$ on $S^n$ with $Sc(g)>0$ for $n=8k-1, 8k$.

This kind of argument  combined with  thin surgery with $Sc>0$   and  empowered by  "higher" index theoretic  invariants of families of  diffeomorphisms,   leads   to the following results. \vspace{1mm}

 {\large\color {blue}\textbf  [}HaSchSt 2014{\color {blue}\textbf ]}.   {\sl If $m$ is much greater than $k$ then the {\it $k$th  homotopy group} of the space $\mathcal G_{Sc>0}(S^{4m-k -1})$ of Riemannian metrics with $Sc>0$ on the sphere $S^{4m-k -1}$ is {\it infinite.}}\vspace{1mm}

{\large  {\color {blue}\textbf [}EbR-W  2017{\color {blue}\textbf ]}}.  {\sl There exists a compact  Spin  6-manifold $X$ such that the  space $ \mathcal G_{Sc>0}(X)$  has each rational homotopy group infinite dimensional.}
 \footnote{It seems, judging by the references in  [Ebert-Williams(infinite loop spaces) 2017] and  [Ebert-Williams(cobordism category) 2019] ,  that all published results in this direction depend on the Dirac  operator techniques which  do not cover the above example, if we take  a {\it Schoen-Yau-Schick manifold}   for $X$.}

However, there is {\large \it no  closed  manifold} of dimension $n\geq 4$, which admits a metric with $Sc>0$ and  where the (rational) homotopy type, or even the set of connected components,    of {\large \it the space of such metrics is fully determined}.  {\footnote {Connectedness 
$ \mathcal G_{Sc>0}(X^3 )$ is proven in [Marques(deforming Sc > 0)2012] by means of the Ricci flow. Conceivably, a similar argument may reduce the study of the homotopy structure 
 of  $ \mathcal G_{Sc>0}(X^3)$ to the space of "standard metrics" with $Sc>0$ on $X^3$. }
\vspace{1mm}

 Let us formulate two specific questions motivated by the following vague one:\vspace{1mm}

\hspace {2mm} {\color {blue} \sf What is the "{\it topology of the geometric shape}" of the (sub)space of metrics with  $Sc\geq \sigma$?}\vspace{1mm}

 {\large \color {blue} {\it Question} 1. }  Given  a  Riemannian   manifold $\underline X$,    numbers $\lambda, \sigma >0$ and an integer $d\neq 0$, let $G(X; \underline X, \lambda, \sigma,d)$ be the space of  pairs $(g, f)$  where $g$ is a  Riemannian  metrics on a  $X$ with $Sc(g)\geq \sigma$  and $f:X\to  \underline X$ is a  $\lambda$-Lipschitz map of degree $d$.

What is the topology and geometry of this space and of the natural embeddings  
$$G(X; \underline X, \lambda_1, \sigma_1)  \hookrightarrow G(X; \underline X ,\lambda_2, \sigma_2)$$ 
for $\lambda_2 \geq  \lambda_1$ and  $\sigma_2\leq  \sigma_1$.

More specifically,

{\sf what is the the supremum $\sigma_+= \sigma_+(\lambda)=\sigma_+(\underline X, \lambda,d)$   of $\sigma$,
such that  the manifold $\underline X$ receives a $\lambda$-Lipschitz map 
$f:   X\to \underline X$ of degree $d(\neq 0)$, where $Sc(X)\geq \sigma$?}

Notice, for instance, that if $\underline X$ is the unit sphere, $\underline X=S^n$, then
$\sigma_+(\lambda)=n(n-1)/\lambda^2$ by Llarull's inequality. (Here as ever we need to assume $X$ is spin.)

\vspace {1mm}

{\large \color {blue} {\it Question} 2. } Let $\sf D$ be some natural distance function on the space $G$ of smooth  Riemannian metrics $g$ on a closed manifold $X$. For instance ${\sf D}(g_1,g_2)$ may be defined as $\log$ of the infimum of  $\lambda>0$, 
such that $$\lambda^{-1}g_1\leq g_2 \leq \lambda g_1.$$

 Let $\sf D_\sigma(g)$ denote the $\sf D$-distance from $g\in G$ to the subspace of metrics with $Sc\geq \sigma$
 and $\tilde {\sf D}_\sigma(g)$  be the  $\sf D$-distance  from the $diff(X)$-orbit of $g$ to    this subspace.  

What are  topologies, e.g.  homologies, of the {\it $a$-sublevels}, $a\geq 0$, of the functions $D_\sigma: G\to [0, \infty) $ and $\tilde {\sf D}_\sigma:  G\to \mathbb [0, \infty)$  and  of  the inclusions 
$$D_\sigma^{-1}(0, a] \hookrightarrow  D_\sigma^{-1}(0, b], \mbox { and  } \tilde {\sf D}^{-1}_\sigma(0, a]  \hookrightarrow \tilde {\sf D}_\sigma^{-1}(0, b] \mbox { for  } b>a? $$

Observe that the counterpart of the above $\sigma_+$ call it 
 $\sigma_+^+(\lambda)=\sigma^+_+(X,\lambda) $, satisfies 
$$\lim_{\lambda\to 1} \sigma_+^+(X, \lambda)=\inf_{x\in X}Sc(X,x)$$
by the $C^0$-closure theorem from  section \ref {C0-limits3}, and it is plausible that the function
$\sigma_+^+(\lambda)$ is  {\it H\"older continuous in $\lambda$.}

%%%%%%%%%%%%%%%%%%

%%%%%%%%%%%%%%%%%%%%%

%%%%%%%%%%%%%%%%%%

\subsection  {\color {blue}  Domination, Extremality and Rigidity of Manifolds with Corners} \label {corners3}

%%%%%%%%%%%%%%%%%%%%%%%%%
Recall  that   a {\it corner structure} on an $n$-manifolds $X$ is   defined by (a coherent sets of)
   diffeomorphisms  of small neighbourhoods of all points in $X$ to neighbourhoods of points  in convex polyhedra $P$  in 
   $\mathbb R^n$.

A  corner structure is called {\it simple} and/or {\it cosimpicial} if these $P$ are intersections of $m\leq n $ half-spaces in 
$\mathbb R^n$ in general position, i.e. such that the dimension of the intersection of their boundaries is equal to $n-m$.

Most (all?) theorems concerning  closed manifolds $X$ with $Sc\geq \sigma$ and, more visibly, manifolds with smooth   boundaries $Y=\partial X$,  have (some proven, some conjectural)  
 counterparts for Riemannian manifolds $X$ with {\it corners}  on the boundary,
where the mean curvature $mean.curv (\partial X)$ for the  smooth part of  $ \partial X$  plays the role of {\it singular/distributional} scalar curvature supported on $ \partial X$ and where
 the {\color {blue} \it dihedral angles  $\angle$} along  the corners, or  rather the {\color {blue} complementary angles}  {\color {blue} $\pi-\angle$},  can be regarded as {\it singular/distributional} mean curvature
supported on the corners.
 
 We bring several examples in tis section      illustrating this idea  starting with the following  definitions.

\vspace{1mm}

{\it \textbf {Domination with Corners}.}  A  proper continuous   map   between manifolds  with corners, 
$f:X\to \underline X$ called {\it corner proper} if the codimension $1$ faces $F_i \subset \partial X$ 
are equal to  the pullbacks $f^{-1}( \underline F_i)$ of the codimension 1 faces of $\underline F_i\subset  \partial\underline X$.

Such a map $f$ between {\it equidimensional} manifolds is called {\it proper  domination} if both 
manifolds are orientable and $f$  has non-zero degree. \footnote {This definition can be generalized by allowing proper maps  that sends some of the ends of $X$ to points and also maps  from spin manifolds  with with non-zero  $\hat A$-degrees, i.e.  with non-zero  $\hat A$-genera of  pullbacks of generic points,  but we don't  do it this here,  since we want to emphasize the corner aspect of the story.  }

\vspace{1mm}

{\it \textbf {Extremality}.} Given a class $\cal F$   of manifolds $X$ along with  dominating maps $f$  from $X$ to a Riemannin manifold $\underline X$ with corners,    
call  $\underline X$  {\it extremal} with respect to $\cal F$ if no map $f\in \mathcal F$ can be simultaneously {\it 

(a) "geometrically contracting"} 
and 

(b)  {\it "scalar and mean curvatures decreasing"}  at all points, 

where  an appropriate (but not the only one) specific meaning of these (a) and (b) is expressed by the following three  pointwise 
 inequalities,  call   them {\color {blue} \textbf {three} {\large  "$\leq$"}},  concerning 
 
{\sf the scalar curvatures versus area contraction in the interiors of the
  manifolds, 
  
  the mean curvatures of their boundaries in (the interiors of ) the  
   $(n-1)$-faces 
   
   $F_i\subset \partial X$,  
  
   the dihedral  
  angles along the $(n-2)$-faces 
$F_{ij}\subset F_i\cap F_j$.}

{\color {blue!50!black}$$ Sc(f(x))\geq   ||\wedge^2df(x)|| \cdot Sc(x), \mbox { $x\in X$,  }  \leqno{\color {blue}[codim =0]} $$
 $$ mean.curv (\underline F_i , f(y))\leq 
||df||\cdot mean.curv (F_i, y),  \mbox { $y\in F_i \subset \partial X$}, \leqno{\color {blue}[codim =1]}  $$   
$$\pi-\angle (\underline F_{ij}, f(z))\leq \pi=\angle ( F_{ij}, z), \mbox { $z\in F_{i,j} \subset F_i\cap F_j\subset \partial X$}.\leqno{\color {blue}[codim =2]} $$ } 

Observe that the fist  inequality {\color {blue} [codim=0]} is satisfied by all maps $f$, 
whenever $\underline X$ is scalar flat ($Sc(\underline X=0$) and $Sc(X)\geq 0$.    No condition on the  norms of the differentials $df$ or on
the  exterior powers $\wedge^2df$ is needed here.  However, we (usually) require in this case that $Sc(X,x)\geq 0$ even at the points $x\in X$ where $df(x)=0$.

Similarly the second inequality is automatic, if the faces   $\underline F_i$ are minimal  (mean.curv=0)   and the faces 
the boundary is mean convex,  $mean.curv F_i\geq 0$, 
e.g. where $\underline X$    is a convex polyhedron in $\mathbb R^n$. In  this case, however  we (usually)  require that $Sc(X)\geq 0$ and the boundary of $X$ is mean convex. 

% and $X$ is manifold with $Sc>0$ and and with smooth (no corners)  mean convex boundary.

Now,  an orientable  $n$-manifold $\underline X$ with corners is called {\it extremal} 
with respect to  a class $\cal F$ of domination maps $f$ from Riemannian manifolds\&maps  from  $\cal F$ 
if none of these  inequalities  {\color {blue} \textbf {three} {\large  "$\leq$"}} for $(X,f)\in \cal F$ can be strict at any point,  i.e.  {\color {blue} \textbf {three} {\large  "$\leq$"}} imply that 
$$ Sc(f(x))=   ||\wedge^2df(x)|| \leq  Sc(x), \mbox { $x\in X$,  }   $$
 $$ mean.curv (\underline F_i ) (X) f(y))= 
mean.curv (F_i, y),  \mbox { $y\in F_i \subset \partial X$,  }$$   
$$\angle (\underline F_{ij}, f(z))= \angle ( F_{ij}, z), \mbox { $z\in F_{i,j} \subset F_i\cap F_j\subset \partial X$}.$$

\vspace{1mm}

{\it Exercises.} (a) Show that the set of  {\it extremal} Riemannian metrics $\underline g$ on a smooth manifold with corners  $\underline X$ (extremality of $\underline g$  means that  for the manifold  $(\underline X, \underline g)$) is closed in the $C^2$-topology 
in the space of Riemannin metrics on $\underline X$, provided this  extremality is understood for a class $\cal F$ in which 
 manifolds $X$    have   $Sc(X)\geq 0$ and $mean.curv(F_i)\geq 0$. 

{\it Hint.}  Adapt  the redistribution of curvature arguments  from  section 11.2   in [G(inequalities) 2018].

(b)  Let  $g_0$   be a smooth Riemannin metric on a manifold $X$ with corners and let $x_0\in X$ be a point in $X$.

Show that there exists a smooth deformation $g_t$, $t\geq 0$, of $g_0$ supported in a given arbitrarily small neighbourhood $U_0\subset X$ 
of  $x_0$ and such that

$\bullet_0$ if $x_0$  lies in the interior of $\subset X$ then  the Scalar curvature $Sc(X,x_0)$ is strictly decreasing; 

$\bullet_1$  if $x_0$  lies in the interior of a codimension 1 face $F_i\subset X$ then  the mean curvature of 
$mean.curv_{g_t}(F_i,x_0)$ is strictly decreasing, while the scalar curvature of $X$ is nowhere decreasing;

$\bullet_2$ if $x_0$ is in the interior of a codimension 2 face $F_{ij}\cap F_j\subset F_j$, then the dihedral angle at this point 
$\angle_{g_t}(x_0)$ is strictly increasing while  the scalar curvature of $X$ and the mean curvatures of the faces are nowhere decreasing.

(i) the curvatures of $g$ are constant, $Sc(X)= \sigma$;

(ii) the faces of      all  edges $F_i$  are also constant, $mean.curv_g (F_i)=M_i$, 

and such that $g$ these $g$ are {\it locally extremal }: 

if a deformation of $g$ doesn't decrease the scalar curvature of $X$, of the mean curvatures of $F_i$, and of the complementary  angles between the edges $\pi-\angle_{i,j}$, 

\vspace{1mm}

{\it \textbf {Rigidity}.} A Riemannian manifold  $X$ with corners  is called  {\it rigid} in      $\cal F$ if {\color {blue} \textbf {three} {\large  "$\leq$"}}  imply that  small neighbourhoods $U_x\subset X$ of all points $x$ are {\it isometric} to some neighbourhoods  $\underline  U_{\underline x}$ (depending on of the the image points $\underline x=f(x)\
\in\underline X$.
\vspace{1mm}

Recall that in the scalar flat  case of complete manifolds rigidity often (but not always)  follows from extremality via 
 the Bourguignon-Kazdan -Warner perturbation theorem. Below is a possible generalization of this theorem  to manifolds with corners, that however has limited applications.  \vspace{1mm}

{\color {red!50!black}   Perturbation Conjecture.} {\sf Let $\underline X=(\underline X,\underline g_0)$ be a  complete  Riemannin  manifold  with corners, such that 
 $Sc(\underline g_0)=0$  and such that all codimension  1 faces are minimal, $mean.curv(F_i)=0$. } 
 
 {\sl Then either $Ricci(\underline X)=0$ and all faces $F_i$ are totally  geodesic, or the Riemannian metric  $g$ admits a   bounded deformation $g_t$, which increases the scalar curvature  and   the mean curvatures of the faces 
 $$\mbox {$Sc(\underline g_t)>0$   and $mean.curv_{\underline g_t}(F_i)>0$, for $t>0$},$$
and also  decreasing the dihedral angles,  $\angle_{ij}(\underline g_t)=\angle_{g_t}F_i,F_j)<\angle_{ij}(\underline g_0) $.}

\vspace{1mm}

{\it Remarks/Questions.} (a) It is unclear   what is  a similar perturbation property (if any)  for {\it non-scalar flat} (potentially extremal) manifolds with corners. 

(b)  It is easy  to see that extremal   surfaces $(\underline X)$ with corners are rigid:  these  have  constant curvatures and in the case of  $Sc\geq 0$, they have geodesic edges.

(c)   Quadrilaterals $\underline X$  in the hyperbolic  plane, such that 

{\sf (i) all angles $\frac {\pi}{2}$;

(ii)   two opposite   geodesic edges ($mean.curv=0$) of equal length, and   the 
two 

other segments  are concentric horospherical
(with $mean.curv=\pm1$)},

\hspace {-6mm}  {\it are rigid.}
  
  {\sf  Let an $X$ dominate $\underline X$, Then, if  and $Sc(X)\geq -2$, if
$mean.curv F_i\geq mean.curv  (\underline F_i$  and if  the angles between adjacent  edges in $X$ are all $\leq \frac{\pi}{2}$,
then $X$ is {\it isometric to a hyperbolic } (i)\&(ii)-{\it quadrilateral. }} 

It seems, there are no similarly  rigid hyperbolic  $k$-gons besides these  quadrilaterals.

\vspace{1mm}

{\color {red!50!black}\textbf {Extremality/Rigidity Problem.}} {\sf Identify/classify extremal and rigid Riemannian  manifolds $\underline X$ with corners for various classes $\cal F$ of manifolds $X$   and dominating maps $f:X\to \underline X$.}\vspace{1mm}

Two motivating examples, where this problems was solved, is the {\it rigidity  of flat metrics on closed manifolds}\footnote {This,  recall, in the case of non-spin manifolds $X$ of dimensions   $n\geq 10$,  
needs   Lohkamp's  or Schoen-Yau's desingularizations theorems.} and

the Goette-Semmelmann theorem, extended by Lott to compact Riemannian  manifolds $\underline X$  with  {\it smooth}
 that    claims that the following three conditions are {\it sufficient for   extremality of  an orientable   $\underline X$ in the  class of {\it\color {red!50!black} spin}  manifolds $X$ 
that  dominate $\underline X$  and have   $Sc(X)\geq 0$ and  mean convex boundaries.}

{\sf (1) The  curvature operator of  $\underline X$  is  {\it  non-negative}. 

(2) The boundary of $ \underline X$ is {\it convex}.

(3) The dimension $n$ of $\underline X$ is  {\it even} and  the  Euler characteristic of $\underline X$  is non-zero.}

Moreover, such an $\underline X=$\underline X=(\underline X, \underline g)$ $  is rigid in certain cases, e.g. if 
$$0< Ricci( \underline g)<\frac{1}{2} Sc( \underline g)\cdot  \underline g. $$

 {\it\color {red!50!black} Conjecturally,}  this holds for {\it all compact Riemannian  manifolds with corners}, which satisfy (1) and (2) and 
with no extra topological assumptions, i.e. possibly {\it non-spin and with $\chi(\underline X)=0$.} 

This may be too strong to be true even for Riemannin flat manifolds, where this reads as follows.

{\color {red!50!black}\textbf {Flat  Corner Domination Conjecture.}} {\sf Let $\underline X$ be a compact orientable  Riemannin flat 
$n$-manifold with corners, such that  all codimension 1 faces $\underline  F_i$ are flat, e.g. $X$ is a convex polyhedron in the Euclidean space $\mathbb R^n$.}

Then $\underline X$ is {\it rigid}:

{\sf if a  proper corner map $f$ of {\it non-zero degree}  from a compact Riemannian manifold $X$ with $Sc(X)\geq 0$,  with mean convex faces $F_i$ and with the dihedral angles between these faces at all points bounded by the corresponding angles $\angle (\underline F_i, \underline  F_j)$,} then 

{\it $X$ is also Riemannin flat, 
the faces $F_i$ are flat, the  dihedral angles between $F_i$ and $F_j$ are  equal  to $\angle (\underline F_i, \underline  F_j)$;  moreover, at all points  $x\in X$ the manifold $X$ is locally isometric to  $\underline X$ at $f(x)\in  \underline X$.}

\vspace {1mm}

Although  this   remains {\color {red!50!black} problematic} even in the category of  {\it convex polyhedra}, where rigidity is known only for  infinitesimal deformations,  see section  \ref {reflection3} and IV below,  the following results are  available.

\vspace {1mm}

I. {\Large \color{blue} $\times\hspace {-0.6mm} \blacktriangle^i$}-{\large \it \color{blue} Inequality.} {\sf  Let $X_0\subset \mathbb R^n$. {\sf Let $\underline  X$ be a compact orientable Riemannian flat $n$-manifold with corners, where all $(n-1)$-faces
 $\underline F_i$ are flat.}
 
 {\it If all dihedral angles $\angle _{i,j}=\angle (\underline F_i, \underline F_j)$ in $\underline X$   are $\leq \frac {\pi}{2}$
 then $\underline X$ is spin extremal}:
 
  {\sf if an  orientable   {\it spin} manifold $X$, which {\it dominates} $X$, i.e. comes with a proper corner map
 $f:X\to \underline X$ with {\it non-zero degree} and such that  
 
 $\bullet_0$ \hspace {1mm} $Sc(X)\geq 0$

   $\bullet_1$ \hspace {1mm} $mean.curv (F_i)\geq 0$ 

$\bullet_2$ \hspace {1mm} $\angle (F_i,F_j)\leq  \angle (\underline F_i, \underline F_j)$,

then  
$$\mbox {$Sc(X)= 0$, 
$mean.curv (F_i)= 0$, 
 $\angle (F_i,F_j)=  \angle (\underline F_i, \underline F_j)$.}$$}

\vspace{1mm}

{\it  Remark/Example.} (a) 
If $\underline X$ simply connected, thus,  is  isometric to a convex polyhedron in $\mathbb R^n$ 
then the condition $\angle _{i,j}\leq \frac {\pi}{2}$ implies (by an elementary argument) that $\underline X$ is   the product of simplices with dihedral angles  $\leq \frac {\pi}{2}$, such as the $n$-cube, fo instance.

{\it About the Proof.} The condition   $\angle (\underline F_i, \underline F_j)\leq \frac {\pi}{2}$, shows  (see section \ref {dihedral4}) that a suitable  
smoothing of  the boundaries of $\underline X$  and $X$ reduces the problem to the rigidity  in the smooth case.
For instance if $\underline X$ is a convex polyhedron one may use  the  mean curvature  spin  extremality  theorem  [$Y_{spin} \to$ \EllipseShadow]} from section \ref {mean convex3}.(If  $n$ is even, is  follows from 
 the above 
Goette-Semmelmann-Lott theorem.)\vspace{1mm}

{\it Exercise}.  
Directly prove the  {\Large \color{blue} $\times\hspace {-0.6mm} \blacktriangle^i$}- Inequality in the case, 
where the faces  $F_i \subset  X$  are {it convex}, rather than only mean convex.
 \vspace{1mm}

%{\it $\triangle$-Remark.} (a$^{\small \blacktriangle}$)  The only polytopes with $\angle_{ij}\leq \frac {\pi}{2}$ are products of simplices, such as the $n$-cubes $[0,1]^n$, for example.

{\it ${ \blacktriangle}-Remark$}.  If both $\underline X $ and $X$ are affine  $n$-simplices, then the implication
$$\angle_{ij}(X)\leq \angle_{ij}(\underline X)  \Rightarrow\angle_{ij}(X)= \angle_{ij}(\underline X)$$
follows from the Kirszbraun theorem with no need for the condition $\angle_{ji}\leq \pi/2$.

But there is no  no direct    elementary  proof of this (unless I am missing something obvious)
if $X$ has {\it convex}, rather than flat, faces

{\it Question.} Are  there  "good" local  boundary conditions for Dirac operators  on manifolds with corners  suitable for proving this kind of theorems similar to what is  done   by John Lott in [Lott(boundary) 2020] and by Christian B\"ar  with Bernhard Hanke in  [B\"ar]-Hanke(boundary) 2021] for manifolds with smooth boundaries? 

(Such  conditions  seem  plausible for orbifold like corners, especially for good orbifolds\footnote {Compare with what is done in  [Bunke(orbifolds) 2007] and in related 
paper cited in there.} 
  but the general case is not so clear.)\vspace{1mm}

II. {\it \color {blue} Reflection Orbifolds.} Let $\hat X$ be a smooth manifold acted upon by a (reflection) group $\Gamma$ generated by reflections in  cooriented hypersurfaces $\hat F_i\subset X$ and let $X\subset \hat X$ be the fundamental domain 
for this action that is   the intersection of  the  "half-spaces" $\hat X_i\subset \hat X_i$ bounded by $\hat F_i\subset \hat X_i$
in $X$. 

This $X=\hat X/\Gamma$  comes with a natural corner  structure and if the action of $\Gamma$ is isometric for a Riemannian metric $\hat g$ on $\hat X$, then all codimension 2-faces $F_{ij}\in X$ are endowed with angles 
 of the form $ \alpha_{ij}(\Gamma)=\frac  {\pi}{2l}$,  $l=1, 2,...$.

We have already  explained in section \ref {reflection3}  that \vspace {1mm}

{\sf if $\hat X$ admits no $\Gamma$-invariant metric with $Sc>0$ then $X$ satisfies the following } \vspace {1mm}

{\it  \color {blue} \textbf {No $Sc>0$ Property}.}
{\sf Let $g$  be a Riemannian metric $g$ on $X$, such that}

{\sl $\bullet_{Sc}$   
 the scalar curvature of $g$ is non-negative:   $Sc(g)\geq 0$;

$\bullet_{mean}$   all $(n-1)$-faces $F_i$ of $X$ are {\it mean convex}: $mean.curv_g(F_i)\geq 0$;

$\bullet_{\angle}$  The dihedral angles $\angle_{ij}$ of $X$   at all points of all $(n-2)$-faces  $F_{ij}=F_i\cap F_j \subset X$ are bounded by the canonical ones, 
$\angle_{ij}\leq \alpha_{ij}(\Gamma)$.}

{\it Then 

 $$\mbox {$Sc(g)=0$,  { } $mean.curv_g(Y_{reg})=0$,  and  $\angle_{ij}\leq \alpha_{ij}(\Gamma) $}$$}  \vspace{1mm}

{\it About the Proof}. This is shown by reflecting $X$ around its $(n-1)$-faces, smoothing around the edges and applying the corresponding result 
for closed manifolds  as it was done  in  [G(billiard] 2014]\footnote{When writing  this paper I overlooked  the paper by Brendle, Marques and Neves [ [Bre-Mar-Nev(hemisphere) 2011], where an essential step of smoothing codimension one corners appears as theorem 5.}  for cubical $X$,  and  where the general case needs an intervention of arguments from [G(inequalities) 2018], where  the  (non-spin) case $n\geq 9$ relies on [SY(singularities) 2017]. 

\vspace {1mm}

An immediate application of  of this to manifolds $X$ which dominate Euclidean reflection domains $\underline X$  is the following \vspace {1mm}

{ \textbf {Extremality of Euclidean Reflection Domains.}  {\it  If  $\hat X=\mathbb R^n$  and the reflections in  $\Gamma$ are isometric then the orbifold  $\underline X= \mathbb R^n/\Gamma$
 is extremal.}

\vspace {1mm}

{\it Remark.} Since the reflection  domains have their dihedral angles  $\angle_{ij}\leq \frac{pi}{2}$ their {\it spin}  extremality  follows from the above {\Large \color{blue} $\times\hspace {-0.6mm} \blacktriangle^i$}-{Inequality.}

{\it\color {blue}  \textbf {$\blacktriangle$-Rigidity.}}  This says that, in fact, 
{\it $X$ is Riemannin flat and all faces $F_{ij}$  are also flat.} \vspace {1mm}

{\it Proof .} The quickest   proof  of the  rigidity is  technical, namely it relies on the {\it regularization theorem} proven in  [Burkhart-Guim(regularizing Ricci flow) 2019]:

{\it If a continuous metric $g_0$ on a Riemannin manifold can be $C^0$-approximated by smooth  metrics $g_\varepsilon$, $\varepsilon>0$,   with 
$Sc(g_\varepsilon\geq \sigma_0-\varepsilon$ for $\varepsilon \to 0$, then   it can be approximated  by smooth metrics  with 
$Sc\geq \sigma_0$.}

We apply this to the $\gamma$-invariant  metric $ \hat g_0$ on $\hat X$ that extend the, a priori {\it non-Riemannian},  metric $g_0$ on $X\subset \hat X$, but 
but because  of the equalities  $\angle_{ij}\leq \alpha_{ij}(\Gamma) $ guarantied by the weak rigidity, this $g_0$ Riemannian and the  {\it regularization theorem} does  apply  and  then an easy argument  shows that the metric $g_0$ itself is Riemannian flat   and the faces $F_i$ are flat as well. \vspace{1mm}

{\it Remarks.}   (a)  The  rigidity for cubical  $X$ of dimension $\leq 7$  was originally proven by Ciao Li in [Li(rigidity) 2019]  and then  extended   in the second version of his paper to manifolds with the corner  structures  combinatorially isomorphic  to that in the product  of the   cube $\square^{n-2}$ by an acute angled triangle $\triangle$, where an essential 
 novel point in this paper  is the proof of a  sufficient regularity of minimal surfaces at the corners  that allows one to argue  as in the proof of the rigidity of flat tori. \footnote {I haven't read Li's paper carefully and I am not certain on how actually  he does it, but,  granted regularity,  the
$\mu$-bubble perturbation argument as  in section \ref{rigidity5} applies in the case considered  by Li.}

(The  products of cubes by general  triangles considered by Li are not, in general reflection orbihedra. On the other hand,
the above argument with reflections+Ricci flow,  implies, for instance, {\it rigidity of products} of several {\it regular} triangles, where no   present day   minimal hypersurface argument applies.)  
\vspace{1mm}
  
 (b) (\it Recapturing   Rigidity while Smoothing the Corners .} 
 \vspace{1mm}
 
 III. {\it\color {blue} Pyramids and Quasi-Prisms.} A counterpart  {\Large \color{blue} $\times \hspace  {-0.6mm} \blacktriangle^i$}-{\large \it \color{blue} inequality}
is known to hold for certain polytopes $P$ with   dihedral angles   $>\frac {\pi}{2}$, which, much as  the above products of simplices,  are {\it extremal} in the sense that  

{\it can't make the dihedral angles smaller, while  keeping the faces mean convex and the scalar curvature $\geq 0$}

The simples such extremal  $P$ are (convex)  $k$-gonal prisms, where for $k\geq 4$  some dihedral angles  are always   $>\frac {\pi}{2}$. This is  shown  in  [G(billiards)  2014] by looking at minimal surfaces with free boundary on the side-part of the boundary of $P$.

More generally Chao Li  [Li(comparison) 2017] proved  a similar property for convex {\it pyramids}  and 
{\it quasi-prisms} $P$ where the latter are convex polyhedra in $\mathbb R^3$, where all vertices are contained in a pair of parallel planes and where the proof follows by a construction and  analysis of  suitable  $\mu$-bubble (capillary  surfaces)  pinched between these planes.

(A technical  limitation  on   deformations of the  flat geometry in  $P$, a mild
lower bound on the dihedral angles between side faces  allowed by Li's argument,   was removed in his later paper.) 
 
   \vspace{1mm}

IV. {\it \color {blue} Polyhedral  Extremality Problem}.  The above kind of  extremality, even the local one,   remains problematic in general  even for {\it simple} $n$-polytopes, where  
   at most $n$ faces 
of dimension $n-1$ may meet at the vertices:   
 
{\sf  it is   unknown which  pairs of combinatorially equivalent    polytopes  $P$ and $P'$ (convex polyhedra)  may have their corresponding  dihedral angles  satisfying $\angle_{ij}\geq \angle'_{ij}$ without all corresponding angles being mutually equal.\footnote{As we mentioned in \ref {reflection3}, Karim Adiprasito told me  that  {\it Schl\"afli formula} (see  [Souam (Schl\"afli) 2004])   implies that  {\it no convex polytope admits an infinitesimal deformability on simultaneously decreasing all its dihedral angles}.} \vspace{1mm}}
 
   laer 
 V. {\it  \color {blue} Extremality and Rigidity of Hyperbolic Manifolds with  Corners.} It is  unclear, in general,  what kind of 
 extremality/rigidity one can expect from manifolds which may have   have negative scalar curvature at some points,  some non-mean convex faces and/or some dihedral angles  $>\pi$.
 
 But the  extremality of flat manifolds  $\underline X=(\underline X, \underline g= \underline g(x))$ with corners passes to the   hyperbolic cylinders 
 $$\underline X_d^\rtimes (-1)  =(\underline X\times [0,d],  g_{\exp}^\rtimes= e^{2t}\underline g(x) +dt^2), 0\leq t\leq d,$$  with constant sectional curvatures $\kappa(g_{\exp}^\rtimes)=-1$. Namely, these  $\underline X_d^\rtimes (-1)$ can't be dominated with manifolds with strict  increase of the scalar curvature, increase of the mean curvatures of the faces  and decrease of the dihedral angles, in the case of extremal $X$.

In fact,  the above reflection, doublings and smoothing arguments apply to these $ \underline X_d^\rtimes (-1)$ in conjunction with the existence and basic  properties of  stable  $\mu$-bubles $Y$  in the cylinders  $\underline X_d^\rtimes (-1)$, which separate the "bottom"  $\underline X \times \{0\}\subset \underline X_d^\rtimes (-1)$ from
 "top"  $\underline X \times \{d\} \subset \underline X_d^\rtimes (-1)$, which 
 have constant  mean curvatures $n-1$ and  such that some  warped products $Y\mathbb T^1$ have non-negative scalar curvatures.
  see section 
\ref {warped boundary5}, 
 
 However,  there are two technical caveats to this  reasoning. 
 
{\color {red!50!black} (1$_{reg}$)}   If $n+1=dim (\underline X_d^\rtimes (-1) ) \geq 8$ the bubles $Y$  may, a priori,  have stable singularities 
 where the present day state of desingularization art of Lohkamp-Schoen-Yau is not, at least not immediately, applicable to all cases of interest.

  {\color {red!50!black} (2$_{reg}$)} Even for $n\leq 7$, the bubbles $Y\subset \underline X \times \{d\}$ are not fully smooth,  at the corners, where the dihedral angles  $\angle_ij(x)\neq \frac {pi}{2k}$, and where, the  unconditional implication 
 $$\mbox {X is extremal $\Rightarrow \underline X_d^\rtimes (-1)$ is extremal}$$
and even more so 
$$\mbox {X is rigid  $\Rightarrow \underline X_d^\rtimes (-1)$ is rigid}$$
needs a  bit of technical reasoning.

 \vspace{1mm}
 
 {\it Motivations  for  Corners.} Besides opening avenues for generalisations of what is known for smooth manifolds,  Riemannian manifolds with corners and $Sc\geq \sigma$  may do   good to  the following. 

1. Suggesting  new techniques,  (calculus of variations, Dirac  operator) for   the  study  of   Euclidean polyhedra.

  2. Organising the totality of manifolds with $Sc\geq 0$ (or, more generally   with $Sc\geq \sigma$) into a 
nice category  ($A_\infty$-category?) ${\cal P}^{\square}$, that would include, as objects  manifolds $Y$ with Riemannian  metrics $h$ and functions $M$ on them
 and where {\it morphisms} are   {\it (co)bordisms}  ({\it h-cobordisms}?) $(X,g)$, $ \partial X =Y_0 \cup Y_1$, where $g$ is  a Riemannian metric on $X$ with $Sc\geq 0$, which restricts to $h_0$ and to  $h_1$ on $Y_0$ and $Y_1$  and where the  the mean curvature of $Y_0$ with inward coorientation is equal to $-M_0$  while  the mean curvature of $Y_1$  with the outward coorientation is equal to $M_1$.

Conceivably, the  variational techniques for   "flags" of hypersurfaces from [SY(singularities) 2017] or its generalisation(s),  may have a meaningful interpretation  in  ${\cal P}^{\square}$, while  a suitably adapted   Dirac   operator method may serve as a quantisation of  ${\cal P}^{\square}$.

\vspace {2mm}

 \hspace {-3mm} {\sf \textbf {  Comprehensive $Sc\geq  \sigma$ Existence     Problem   for Manifolds with Corners.}}\vspace {1mm}

{\sf Let $X$ be a smooth compact  manifold with corners and let  $X_i$, $i\in I$, \footnote {We  switched   the notation from $F_i$ to $X_i$ to place all faces, includung $X$ itself, on equal footing.}    be tote the faces here $X_i$ he set of faces of $X$ of all (co)dimensions, where, we agree that $X_0=X$ and let 
$\sigma_i :X_i\to [-\infty, \infty)$, $ \mu_{i,k}  :X_i\to  [-\infty, \infty)$ and $\alpha_{ijk}: X_i\cap X_j\to (0,2\pi)$  be continuous functions, where

$\mu_{ik}$  are defined for all $i$ and those $k$ for which $X_i$ serve as  codimension one faces in $X_k$;

$\alpha_{ijk}$ are defined for the  pairs of codimension one  sub-faces $X_i,X_j\subset X_k$, such that
$dim( X_i\cap X_j) =dim(X_i)-1=dim(X_j)-1=dim(X_k)-2$.}\vspace {1mm} 

{\sl When does there exist a smooth Riemannian metric $g$  on $X$, such that}\vspace {1mm} 

{\sf $\bullet_{Sc}$  the scalar curvatures  of $g$  restricted to $X_i$ are bounded from below by $\sigma_i$, that is 
$$Sc(g_{|X_i},x)\geq \sigma_i(x),\mbox { } x\in X_i;$$\vspace {1mm} 

$\bullet_{mean}$  the mean curvatures of $X_i\subset X_k$  with respect to $g$, for all $X_i$  and $X_k$, where $dim(X_i)=dim(X_k)-1$, are bounded from below by
$\mu_{ik}$;

$\bullet_{\hspace {-0.7mm}\angle}$ the dihedral angles between $X_i$ and $X_j$ in $X_k$ satisfy
$$\angle_g (X_i,X_j)\leq \alpha_{ijk}.$$}

More generally, one wants to understand  the topology (e.g. the homotopy type) of  the spaces 
$G(\sigma_i, \mu_{ik}, \alpha_{ijk})$ of  metrics $g$ on $X$,  which satisfy 
$\bullet_{Sc}$, 
$\bullet_{mean}$ and 
$\bullet_{\hspace {-0.7mm}\angle}$, as well as of the inclusions  
$$G(\sigma_i, \mu_{ik}, \alpha_{ijk})\hookrightarrow G(\sigma'_i, \mu'_{ik}, \alpha'_{ijk})$$
for $\sigma'_i \leq \sigma_i$,  $\mu'_{ik}\leq  \mu_{ik}$ and $\alpha'_{ijk}\geq \alpha'_{ijk}$
and restriction maps  from these spaces to the corresponding ones on submanifolds $Y\subset X$ compatible with the corner structures.
.\vspace {1mm}

{\it Exercise.} Let $X$ be  a smooth $n$-manifold with cornered boundary  $Y=\partial X$, and prescribed mean curvatures of the top dimensional faces $X_i$  and the dihedral angles between them, that is, in the above notation:  \vspace {1mm}

$\sigma_i =-\infty$, unless  $dim(X_i)=n-1$, 

$\mu_{i,k}=-\infty$, unless $k=0$, i.e.  $dim(X_i)=dim(X_k)-1= n-1$ and 

 $\alpha_{ijk}=2\pi$, unless  $dim (X_i)=din(X_j) = dim (X_k)-1 = n-1.$\vspace {1mm}

Show that $X$ admits a smooth metric $g$, such that 

{\it the  scalar curvature of $g$ is positive  in  a (small) neighbourhood of $Y\subset X$},  

\hspace {-6mm} and such that $g$ satisfies the above  conditions 
$\bullet_{mean}$  and $\bullet_{\hspace {-0.7mm}\angle}$.\vspace {1mm}

{\it Hint. } Construct $g$  in a  small neighbourhood of the  union of the  $i$-dimensional faces by  induction on $i=1,2,...,n-1$.
\vspace {1mm}

{\it \color {blue} Measure Valued Curvature.}  The mean curvature of the boundary $\partial X$  and the 
complementary dihedral angles $\pi-\angle_{ij}$ can be regarded as measures  with continuous densities on the faces   which represent singular scalar curvature, where this becomes especially clear if you think in terms of the double \DD$X$.

With this in mind,  the above problem can reformulated as follows:\vspace {1mm}

{\sf  given a triangulated $n$-dimensional manifold $X$(pseudomanifold?) and numbers $\sigma=\sigma (\Delta)$  assigned to all simplices  
$\Delta$   of codimensions 0, 1 and 2.

When does there exist  a  continuous  piecewise smooth Riemannian metric on $X$, such 
that 
its scalar curvature, understood as a measure, is bounded   from below    by $\sigma(\Delta)$  on  all of the  above 
  $\Delta$}, \vspace {1mm}
  
  \hspace {-6mm}where  the inequality $Sc(X|\Delta)\geq \sigma(\Delta)$ is understood as  earlier, namely,

(i) if   $dim(\Delta)=n$  this is  the usual $Sc \geq \sigma$;
 
(ii) if $dim(\Delta)=n-1$ this is the $\sigma$-bound on the sum of the mean curvatures of this 
$\Delta$ in the two adjacent $n$-simplices;

(iii) if $codim(\Delta)=n-2$ this is $2\pi$ minus the sum of the dihedral angles of the $n$-simplices adjacent to $\Delta$.

\vspace {1mm}

{\it Remark on Higher Codimension Singularity.} Strictly speaking the above applies not to triangulated but to stratified manifolds $X$, where

$\bullet $ there are  only strata of codimensions 0, 1 and 2,

$\bullet $  the codimension 2 strata  are smooth submanifolds in $X$,

$\bullet $   the codimension 1 
strata $\Sigma_{-1}$  are submanifolds with boundaries with all components of these boundaries being codimension 2 
strata $\Sigma_{-2}$, where different $\Sigma_{-1}$ with common components
 $\Sigma_{-2}$  of their boundaries meet transversally at these $\Sigma_{-2}$, 
 
$\bullet $  the Riemannian metics in question are piecewise smooth with respect to this stratification.

One may also to allow singularities of codimensions  $\geq 3$, but this is a different matter (compare with  (c) in 
\ref {max-scalar5}).

%%%%%%%%%%%%%%%

\subsubsection{\color {blue} Corners, Plateauhedra and Bubble Spaces}\label {plateauhedra3}

%%%%%%%%%%%%%%%%%%%%%

 Central  geometric examples  of manifolds with corners are convex polyhedra   in the Euclidean spaces and, 
 more generally domains in spaces with constant sectional curvatures $\kappa$  that are intersections on of half 
 spaces bounded by by umbilic  hypersurfaces,  that are  spheres, hyperplanes and, for $\kappa<0$, horospheres and equidistances of hyperplanes.

 \hspace{13mm}{\sf \large \color {red!30!black} What are Riemannin Counterparts of these?}
 
Below are candidates for answers.
\vspace{1mm}

$\bullet _b$ {\it Bubblehedra.} A bubblehedron $Q$ in a Riemannian  $n$-manifold $X$ is the boundary of a  domain $Q_>\subset X$ 
with corners
$$Q=\partial Q_> =\bigcup_ {i=1,2,...} Q_i$$
 where all $(n-1)$-faces $Q_i\subset \partial X$  have {\it constant mean curvatures $M_i$}, where all 
{\it dihedral angles $\angle_{ij}$ are constant}   along the $(n-2)$-faces $Q_i\cap Q_j$ and  where one may require  
these angles to be $\leq \pi$.

A special case of these where all $M_i=0$ are called {\it plateuhedra.}\vspace{1mm}

{\it Remark.} The common description  of minimal varieties in terms  of currents doesn't seem appropriate  for such $Q$, and even less so for similar arrangements of minimal subvarieties  $Q_i \subset X$  of   codimensions >1.
\vspace{1mm}

{\it Example 1: Normal Plateauhedra and Bubblehedra.} An attractive  instance of these is where all  dihedral angles $\angle_{ij}=\frac {\pi}{2}$ and where, moreover,
each  face $Q_i$ is (n-1)-volume minimizing with  free boundary in the union of the remaining edges, 
$$\partial Q_i\subset \bigcup _{j\neq i}Q_j.$$

A more general similar case is  where $M_i>0$ and  each $Q_i$ with  free boundary in $\bigcup _{j\neq i}Q_j.$.
minimizes $vol_{n-1} (Q_i)-M_i \cdot vol_n(Q_>)$. \vspace {1mm}

{ \it Example 2: Normal  Plateau Webs.} Let $X$ be a compact Riemannin manifold and let $Y_1\subset X$ be a a closed locally minimizing minimal hypersurface in $X$.
Next, let  $Y_2$  be a locally minimizing minimal hypersurface in the complement of $Y_1$  in $X$ 
 with free boundary in
$\partial Y_1$, i.e. 
$$Y_2\setminus \partial Y_2\subset X\setminus Y_1 \mbox { and } \partial Y_2\subset Y_1.$$

Then continue with minimal $Y_3, ..., Y_i$...,
where $Y_{i}\subset X$ lies in the complement of  all $Y_j$,$j<i$  and has free boundary in the union of these $Y_j$,  
$$Y_i\setminus \partial Y_i\subset X\setminus  \bigcup_{j<i}Y_j, \mbox { } \partial Y_i\subset \bigcup_{j<i}Y_j.$$

Thus we divide $X$ into mean convex domains  with $90^\circ$ dihedral angles.
\vspace {1mm}

{\it Questions.}  What do  combinatorics of such webs  $\{Y_i\}$  tell you about the topology and geometry of $X$?

How much does positivity of the scalar curvature $X$ restrict combinatorics of such a  $\{Y_i\}$?

If $X$ is  complete non-compact, where  "plateau"    may be too restrictive,  and asks:

  what geometric/topological  condition(s) on $X$ would guarantee the existence of  normal  {\it bubblehedra} $Q\subset X$   of  given combinatorial types?\vspace {1mm}

{\it \color {red!40!black} Conjectural} { \it Example 3:  Dodecahedral and Similar Exhaustions of Large Manifolds.}
If $n=3$ and $\tilde X$ is a the universal covering of a compact  Riemannian manifold  $X$,  where this $X$ admits a hyperbolic Riemannin metric (probbaly,  non-zero simplicial volume will do), then {\color {red!40!black} it seems} not hard to show that it can be exhausted by (compact domains $Q_>\subset \tilde X $ bounded by)  such normal  bubblehedra $Q\subset \tilde X$ of dodecahedral combinatorial types,  i.e. admitting proper corner maps of degree 1 to  (the boundary of) the dodecahedron).

Also "hyperbolically looking"  manifolds $\tilde X$ of dimensions $\geq 4$, e.g. the   universal coverings of compact manifolds $X$ with  {\it non-zero simplicial volumes}, can {{\color{red!40!black} probably}  be   exhausted  by  similar $Q$. 

 For instance, if $X$ is the product of surfaces of genera $\geq 2$, then such an exhaustion  is expected by normal bubblehedra  $Q$ of combinatorial types of products of $k$-gons. \vspace {1mm}

{\it Example 4; Local Riemannian Realization    of Euclidean $P$.}  Let $P$ be a convex polyhedron in a tangent space $T_{x_0}(X)=\mathbb R^n$, let us scale $P$ by a small
$\varepsilon>0$ and let 
$P_\varepsilon '\subset  X$ be the image of this $\varepsilon P\subset T_x(X)$ under the exponential map $\exp: \varepsilon P\subset T_x(X)\to X$.

This $P_\varepsilon' \subset X$, which is not a true but only an  {\it $\varepsilon$-approximate plateauhedron},  already may carry some  information about the scalar curvature $Sc(X,x)$,  in terms of the mean curvatures of its $(n-1)$-faces and the dihedral angles along its $(n-2)$-faces similar to the  {\Large \color {blue} $\circledast$} representation  of   the inequality $Sc(X, x_0) < Sc(X ', x_0 ')$ in section  \ref{deceptive1} by comparison  the integral mean curvatures of the $\varepsilon$-spheres around the points   $x_0$ and  $x_0 $ in two manifolds.

Next,  to make  this $P_\varepsilon '\subset X$ look prettier, one can    slightly perturb it  and thus turn it   
 into a true  bubblehedron $Q\subset X$ by solving the Plateau soap bubble problems with free boundaries  for all $(n-1)$-faces one after another\footnote{One can't, a priori, guarantee the full (not even $C^2$) regularity  of the $(n-k)$ faces for $k\geq 2$, 
 but in view of high   non-uniqueness   of these $Q$ explained below,  such regularity seems  non-impossible in many  cases.} 
 where, depending of what one wants, one can either make   all its {\it faces with zero mean curvatures}, or all its {\it dihedral angles equal those of the original $P$}. (This  seems  easy but I 
 didn't try to carefully check it.\footnote{To fully include {\Large \color {blue} $\circledast$} in this picture, one had to start with a $P$, which has     {\it spherical} as well as  planar faces.  But then   perturbing  $P_\varepsilon' $ into a bubblehedron $Q$ becomes a more delicate matter. For instance, if, as it is in {\Large \color {blue} $\circledast$}, our    $P\subset T_x(X) $ is a ball bounded by a single spherical "face",  then the corresponding $Q\subset X$ (bounded by a single hypersurface of constant mean curvature)   may (and usually will)  drift away from  the point $x$.}) 
 
We know, however that  if $Sc(X,x_0)>0$, there are  some constraints on  possible  values of the mean curvatures 
 $M_i(Q)=mean.curv(Q_i)$ and the dihedral angles $\angle_{ij}(Q) $ of such a $Q$, e.g.  
 if $\angle_{ij}\leq\frac {\pi}{2}$, then,  for small $\varepsilon\to 0$,  one can't have 
 all $M_i(Q)\geq   M_i(P) $ and $\angle_{ij}(Q)\leq \angle_{ij}(Q) $,   where {\color {red!50!black} conjecturally} this is true for all $P$.
 
Despite this,  in general, if $n\geq 3$,
 the space $\mathcal Q$ of all bubblehedra (or plateuhedra) $Q$  in a small neighbourhood of $P_\varepsilon'$ is typically   {\it infinite dimensional}. \vspace {1mm}
 
{\it Example 5:  Too Many  $Q$. } Let  a plateuhedron $Q$ in a Riemannin manifold $X$ contains  only two $(n-1)$-faces $Q_1$ and $Q_2$, which are compact  smooth  hypersurfaces in $X$ the common boundary of which  makes the only $(n-2)$-face of $Q$,
$$Q_{12}=Q_1\cap Q_2=\partial Q_1=\partial Q_ 2$$

Imagine that    $Q_1$ extends in $X\supset Q_1$ beyond its boundary to a minimal  $ Q_{1+}\supset Q_1$, such that

 $\bullet_{1+}$  the extended face   $ Q_{1+}\subset X$ is 
  {\it a strictly locally minimizing hypersurface}\footnote {"{\it Strict} local minimum" usually means "{\it isolated} local minimum" -- this is sufficient for most our present geometric purposes. But  if following an analytic  vein of thinking, "strict" should be understood as {\it strict positivity} of the second variation operator.}   with  respect to  its boundary $ Z=\partial  Q_{1+}$;
  
  $\bullet_{2}$ the face  $Q_2$ is  {\it strictly locally minimising   with free  boundary} in $Q_{1+}$.
  
   Then  small deformations $Z'$ of the (smooth  closed)  $(n-2)$-submanifold $Z\subset X$
    are (by an elementary elliptic  perturbation  argument)   accompanied by  {\it unique minimal} deformations
     $Q'_{1+} $ of $Q_{1+}$   {i.e.  submanifolds $Q'_{1+} $ are minimal) followed by {\it minimal} $Q_2'$ that are small deformations of $Q_2$,  such  that 
 the  boundary of such a  of $Q'_2$  is   {\it contained in} $Q'_{1+} $  and where   $Q'_2$ is       {\it normal}   to  $Q'_{1+} $  along $\partial Q'_2 $.

Thus the local moduli space of these $Q$ contains, as subspace, the {\it full  space of small functions on} $Z$ corresponding to small deformations of $Z$ 
{\it normal} to $Q_{1+}$.\footnote {Deformations of $Z$ within   $Q_{1+}$  don't affect $ Q$.}
\vspace {1mm}

 {\it Example 6: Too few $Q$.} Let  {\it both faces} $Q_1$ and $Q_2$ in the above example extend to minimal hypersurfaces beyond  their boundaries, say to   $Q_{1+}\supset Q_1$ and $Q_{2+}\supset Q_2$, and let both  be  {\it a strictly locally minimizing hypersurface } with  respect to  their  boundaries $ Z_1=\partial  Q_{1+}$  and  $ Z_2=\partial  Q_{2+}$.

Let $Z_1'$ and  $Z_2'$ be small perturbations of these  $\partial Q'_{1+}$  and   $ \partial  Q'_{2+}$, let   $Q'_{1+}$, and $Q'_{2+}$ be  the corresponding minimal perturbations of $Q_{1+}$, and $Q_{2+}$,  let   
$$Q'_{12}=  Q'_{1+} \cap Q'_{2+}$$
and let $\angle'_{12}$ be the dihedral angle between $Q'_{1+}$,  $  Q'_{2+}$  regarded as a function on the perturbed intersection $Q'_{12}$ of the two  $(n-1)$-faces of $Q$,
$$\angle'_{12}=\angle'_{12}(q'),\mbox{ } q'\in Q'_{12}.$$

Here, the situation is opposite to that in the previous example:

{\sl the operator (map) 
$$(Z'_1, Z'_2)\mapsto \angle'_{12}$$
from the space of small deformations of the boundaries of $Q_{1+}$  and  $ Q_{2+}$ to the space of functions on 
$Q_{12}$\footnote {Normally project perturbed intersections  $Q'_{12}$ to $Q_{12}$ and thus identify the spaces of 
functions on all   $Q'_{12}$ with the space f functions on the unperturbed $Q_{12}$.} is (by elliptic regularity) {\it compact}.}

Hence, only 

{\sf a {\it minority of functions on} $Q_{12}$ is realizable by dihedral angles of   (not quite) plateuhedra with  {\it minimally extendable faces.}}
 
 \vspace{1mm}

 {\it Ouroboros Example 7: Biting its Own  Tail.}} Let us describe a class of hypersurfaces, where the two opposite phenomena from  the above examples strike a balance  and make  the Plateau  problem "well posed",  in particular, allowing its {\it a  Fredholm representation}.\footnote {The concept of "Fredholm"  strikes as   artificial in the  present  geometric picture  and begs for something more adequate.  Perhaps,  I am  missing something in  the literature.}

 Let $Q\subset X$ be the image of  a compact $(n-1)$-manifold,\footnote {Much of what follows makes sense for $Q$ of codimension >1.} $n=dim(X),$ with boundary, say $\hat Q$, immersed to $X$,
 $$h:\hat Q\to X,\mbox { } h(\hat Q)=Q, $$ 
 such that 
 
$\bullet_{int}$ \hspace {1mm} the immersion $h$  is one-to one on the interior of $Q$, $$h: \hat Q\setminus \partial \hat Q \hookrightarrow X,$$
where the images of the connected components of $\hat Q$   serve as the  $(n-1)$-faces $Q_i$ of $Q$;

$\bullet_\partial$  \hspace {1mm} the immersion $h$  is one-to one on the  boundary of $\hat Q$;
 $$h: \hat Q\setminus \partial \hat Q \hookrightarrow X; $$
 ;
 
 $\bullet_\subset$  \hspace {1mm}the image of  the boundary  of $ \hat Q$ is contained in the image of its interior, 
 $$h(\partial \hat Q)\subset h(\hat Q\setminus \partial \hat Q),$$
where  the images of the connected components of $\partial \hat Q$   serve as the  $(n-2)$-faces or corners of $Q$;

$\bullet_{min}$ \hspace {1mm} {\it \color {blue!60!black}$Q$ is minimal}: it  has zero mean curvature and it is normal to itself along the corners.\vspace{1mm}

{\it Sub-Example 8.}  The most transparent instance of this is where  $Q$ is the  union  of two faces  that are  smooth submanifolds in $X$ with boundaries, 
 $$Q=Q_1\cup Q_2\subset X,$$
  such that the boundary of one is contained in the  interior of another,
$$\partial Q_1\subset int(Q_2)\mbox  {  and } \partial Q_2\subset int(Q_1).$$  
Thus, the 
 the corner of $Q$ is equal the intersection of the two faces of $Q$,
$$Q_{12}=Q_1\cap Q_2,$$
and where one may think of $Q_1$ as the solution of the Plateau problem with free boundary in $Q_2$ and, similarly,  $Q_2$
is minimal with boundary in $Q_1$.

 \vspace {1mm}

{\it Proposition/Example 9: {\it Finite dimensionality of  Deformations and Codeformations}.}  It {\color{red!40!black} seems obvious}, (I didn't check  this carefully)  that, by the standard elliptic estimates, the space of the above compact minimal $Q$ in a
 small  $C^\infty$-neighbourhood of a given minimal  $Q_0$ is finite dimensional. 
 
 It is slightly less obvious that, given an above  minimal   $Q\subset X=(X,g)$, there exists a {\it finite dimensional}{\footnote {This  dimension  can be  bounded by the index of the second variation operator for $Q$.} linear space $\Delta$ of $C^\infty$-smooth  quadratic forms $\delta$ on $X$, such that, for all Riemannin metrics $g'$ sufficiently $C^\infty$-close to  $g$  \vspace {1mm}
 
 {\sf there exit a small $\delta \in \Delta$ and   $C^\infty$-small perturbation $Q'$
  of $Q$ such that   \vspace {1mm}
  
\hspace {10mm}  $Q'$ is 
   {\it  minimal with respect to  the Riemannin metric $g'+\delta$}.  }

 \vspace {1mm}

 Let us explain this in he simplest case where 
 $Q=Q_1\cup Q_2\subset X=(X,g) $ as in the above sub-example, where we assume that both   $ Q_1$ and $Q_2 $ are {\it strictly locally minimizing} with the free boundary conditions 
 $\partial Q_1\subset Q_2$ and $\partial Q_2\subset Q_1$.

 Slightly $C^\infty$-perturb the Riemannin metric in $X$, say $g\leadsto g'$,  
 and  show that $Q$ can be accordingly deformed to $Q'\subset X$, which  is   strictly $(n-1)$-volume minimizing with respect to $g'$ with
 similar  free boundary conditions . \footnote{In the classical  case,  where $Q\subset  X$ is a smooth closed  strictly locally minimizing submanifold (no boundaries),  it is not hard to show that it is  stable under {\it $C^0$-small} perturbations of $g$; {\color{red!30!black} probbaly} the same applies to   $Q$  with smooth edge(s) $Q_{12}$ and, {\color{red!60!black}possibly} to  general   {\it semi-regular} $Q$ presented later in this section.}

 The simplest way to do it is by consecutively   minimizing $g'$-volumes of $Q_1$ with free boundary 
 $\partial Q_1\subset Q_2$, then of the volume of $Q_2$ with boundary in the new $g'$-minimal $Q_1$, etc
  
  Then, for sufficiently small $g-g'$,  the  strict minimality of $P$ implies the convergence of this process 
to   $Q'$ which lies $C^\infty$-close to $Q$ (for the obvious $C^\infty$-topology in the space of our $Q\subset X$) and  $ g'$-volume  minimizing with free boundary  positioned on the non-singular part  of $Q$. \footnote {This $Q'$ can be {\it defined} as a  fixed point of a self mapping in the space of  $P$ behind this iteration process, where the strictness of minimality makes this self-mapping 
 (which,  by the way, is compact) {\it contracting.} }\vspace {1mm}

{\it Example 10: Higher Order Corners.} Let us generalize the above  sub-example by allowing piecewise smooth 
$$Q=\cup Q_i\subset X$$ 
  where all $Q_i\subset X$, $i=1,2,...,k$, are submanifolds with corners, such that the boundary of each of them  is contained in the union of others,
  $$\partial Q_i\subset  \bigcup_{j\neq i}Q_j.$$
 
 More general $Q$ of this kind is where such a  decomposition exits only locally at all points in $Q$:
  
  {\sf given a point $q\in Q\subset X$, there exists a neighbourhood of this point in $X$, say $U(x)\in X$,
  such that the intersection $Q\cap U(q)$ admits the above kind of decomposition
  $$ Q\cap U(q)=\bigcup _i Q_i(q),
  \mbox  {  where  }
   \partial Q_i(q) \subset  \bigcup_{j\neq i}Q_j(q).$$}
  
  It still make sense here to speak of {\it minimal} $Q$, i.e. with all $mean.curv(Q_i(q))=0$ and where,  minimality with free boundary in $\bigcup Q_j(q)$ is also well defined for all $Q_i(q)$.
  \vspace {1mm}
 
 {\it Question.} {\sf What is the most general assumption(s) on local topology of such $Q$ that would imply the above kind
  Deformations and Codeformations   finite dimensionality  properties?} 
 
   \vspace {1mm}
 
  {\it Example 11: Semi-regular $P$ and $Q$.}  Recall that a {\it simple cone} in  $\mathbb R^{n-1}$ is the intersection of at most $n-1$ half spaces, with  mutually transversal  boundary hyperplanes.
  
 Now, call  a piecewise linear linear cone $P\subset \mathbb R^n$  {\it semi-regular} if it is equal to the union of $k\leq n$ 
 mutually transversal simple cones $P_i$ in some hyperplanes in $\mathbb R^n$, 
 $$P=\bigcup_i P_i,$$
  such that  the boundary of each $P_i$ is contained in the union of the remaining ones,
  $$\partial P_i\subset \bigcup_{j\neq i  } P_j,$$
and,  moreover,  such that   the interior of each $(n-2)$-face in $P_i$, for all $i$, is contained in the interior of some $Q_j$.  
  
 \vspace{1mm}
 
 {\it Example 12:  Cones of rank k=1, 2  and 3}. The cones of rank 1 are just hyperplanes in $\mathbb R^n$.
 
 A cones of rank 2 is a union of a hyperplane  $P_1\subset \mathbb R^n$ and a half-hyperplane  $P_2\subset \mathbb R^n$ with
 its boundary (an $(n-2)$-subspace) in $P_1$.

 If $k=3$, then   there are two possibilities for the position of the third face  $P_{3}\subset P$:
 This can  be either   a  half-hyperplane  positioned in the halfspace  bounded by $P_1$ on  the other side from $P_2$ or    be an 
 $(n-1)$-cone with two $(n-2)$-faces which is positioned  in one of the two convex cones 
 bounded by $P_1$ and $P_2$.
   \vspace{1mm}

  {\it Semi-Regular $Q\subset X$.}   A piece wise smooth $Q\subset  X$ is called {\it semi-regular} if, locally, at each point it is diffeomorphic to a   semi-regular cone.
  
    \vspace{1mm}

 { \color {red!50!black} Conjecture.} {\it Minimal Semi-regular} $Q\subset X=(X,g)$ {\sf enjoy the   deformations and codeformations   finite dimensionality  properties.}\vspace{1mm}

   {\it Remark} (a) This conjecture, is  probbaly, not hard to prove  but 
     the semi-regularity condition is too restrictive and much of it seem unneeded, such as 
   the  transversality condition between  the half-hyperplanes $P_2$ and  $P_3$ attached along their boundaries on two different side to a hyperplane $P_1\subset \mathbb R^n$ in the above rank 3 example.
   
  More seriously, semi-regularity excludes singular minimal $Q$ in dimensions $N\geq 8$.
  
  (b) It remains unclear if our minimal $Q$ have any global  geometric significance.

    (c) The full transversality condition, albeit, probbaly, redundant,    implies the following convenient  (irrelevant?) simple property.    \vspace{1mm}

   {\sf Let $Q\subset X$ be  compact semi-regular and let $Q'\subset X$ be $C^\infty$-close to $Q$, which means all, locally defined, 
   $(n-1)$-faces of $Q'$  are close to the corresponding faces of $Q$. (This is formulated more carefully below.)  Then there is a diffeomorphism $\Phi':X\to X$
   that moves $Q$ to $Q'$; moreover, there is a $C^\infty$-continuous  map $\Psi': Q'\to \Phi'$ from the space of $Q'\subset X$ to $Diff$ such that $\Psi'(Q')(Q)=Q'$  and $\Psi'(Q)=Id$.}

   \vspace {1mm}

 {\it Bubble-Spaces $Q\subset X$ with Variable Mean Curvatures.} The basic properties, including 
 deformations and codeformations   finite dimensionality  properties  for semi-regular 
 minimal Plateau spaces  $Q$ extend  verbatim to to  bubble spaces with  constant mean curvatures $M$,
 including stability of strictly minimizing ones under variations of $M$ keeping $M$ constant.

But one  needs be more careful with variable  mean curvature  of $Q$, since it is not and is  {\it not supposed to be continous} as function    on  $Q$ with the topology 
 induced by the  
 the {\it embedding}  $Q\subset X$. 
 
 Another problem is comparing the mean curvatures of two different spaces, $Q$ and  $Q'$ in $X$, 
 let them even be very close one to another.
 
 To handle this, we recall that $Q$ is the image of a smooth  manifold with boundary under a smooth  immersion, 
 $$h: \hat Q\to X.$$
 Accordingly, the mean curvature is required to be continuous, smooth if you wish,  as a function on this 
 $\hat Q$.
 
Also  the $C^r$-distance for all $r<\infty$  between $Q$ and $Q'$ is defined  as the infimum of  the numbers $d_r$, 
such that there exists a diffeomorphism $\phi:  \hat Q\to \hat Q'$ for which the $C^r$ distance 
between the immersions $h:\hat Q\to \mathbb R^n $ and $\phi\circ h': \hat Q\to \mathbb R^n $
for $h': h:\hat Q'\to \mathbb R^n$ is $\leq d_r$.

 Finally, the measures $\mu$ behind  (the existence theorems for) $\mu$-bubbles are not defined on $X\supset Q$ or on a neighbourhood 
 $U=U_Q \supset Q$ in $X$ but  on an $n$-manifold 
 $\hat U= \hat U_Q$,  which comes with   an {\it immersion $\alpha: \hat U\to X$} and  an {\it embedding 
 $\beta: \hat Q\to \hat U$}
 such that 
 $$ \alpha \circ \beta=h:\hat Q\to Q\subset X.$$
  
  Granted that one sees  the same picture of {\it small}  deformation and codeformation of bubble-spaces 
   as for constant 
  mean curvature,  where  codeformations refer  here  to finite dimensional families of functions (or measures) on 
  $\hat U$. 
  But understanding global properties of such $\mu$-bubbles remains even more  limited   than 
that for constant $M$.
  
 {\sf  {\color {red!20!black} But can one prove (or just conjecture)  nevertheless something non-trivial  about these $Q$  
 in relation  to the  scalar curvature problems? 
  
  \vspace {0.6mm}

Wouldn't it be, perhaps,  more sensible  to switch to a reasonably regular class of {\it minimal   varifolds},  e. g.  $V^{n-1}\in X$, where the singular locus of such a $V^{n-1} $ is a smooth $F^{n-2}\subset V^{n-1} $,  where  there  are  three local branches of $V^{n-1} $ meeting along this  $F^{n-2}$ and where the dihedral angles between these branches are  $\frac {2\pi}{3}$? \vspace {0.6mm}

 Can, one, alternatively,   "rigidify" bubblehedra by minimizing a single functional, a weighted  combination of volumes of faces of different dimensions and /or involving also dihedral angles and (directions of  vectors of) means curvatures of low dimensional faces?}}

%%%%%%%%%%%%%%%

\subsection {\color {blue}Stability of Geometric Inequalities, Metrics and Topologies in Spaces of  Manifolds,   Limits  and Singular Spaces with   Scalar Curvatures bounded from Below}. \label {stability3}

%%%%%%%%%%%%%%%%%

Inequalities relating geometric quantities $\mathcal A$   and   $\mathcal B$ of geometric objects $Ob$ progress along the following lines.

1. {\it Rough Inequalities.}   This says that $\mathcal A(Ob)$ is bounded by {\it some} function of $\mathcal B(Ob)$.

               \hspace {10mm }  For instance,   {\sf volumes  of Euclidean domains $V\subset \mathbb R^n$ are bounded by

   \begin{turn}{-90} 
$\leadsto$
\end{turn}     \vspace{-2mm}     \hspace {10mm } the (n-1)-volumes of their boundaries.}
    
    \hspace {10mm } (There is  about a dozen  of   "direct elementary" proofs, of this which 
          
          \hspace {10mm } generalize to a variety of situations, e.g. to Riemannian  manifolds
          
          \hspace {10mm } with certain restrictions to their curvatures.)\vspace {1mm}

2.  {\it Sharp  Inequalities.}  These specify   the  maximal  values  of $ {\mathcal A}(Ob)$   among  all $Ob$ with  a given bound on $ {\mathcal  B}(Ob)$, say, in the form  $ {\mathcal A}(Ob) \leq E_{sharp} ({\mathcal  B}(Ob))$.

    \hspace {10mm } For instance, the Euclidean domains satisfy  the sharp {\it isoperimetric

     \hspace {10mm }  inequality} 
        $vol(V)\leq \gamma_n \cdot(vol_{n-1}\partial V)^{\frac {n-1}{n}}$, where $\gamma_n$ is equal

  \begin{turn}{-90} 
$\leadsto$
\end{turn}  \vspace{-2mm}\hspace {10mm }  to the volume of the $n$-ball with  unit $(n-1)$-volume of the  boundary.

    \hspace {10mm } (There is no direct elementary proof of this, except for $n=2$ and $4$, 
    
     \hspace {10mm }  and  the present  day "non-elementary"  proofs don't    generalize to the 
     
   \hspace {10mm }  expected cases,    such as  complete manifolds  with non-positive 
   
      \hspace {10mm }  sectional  curvatures.)
 \vspace {1mm}

3. {\it Rigidity.}  This is a description of all   {\it extremal} $Ob$ that maximize $ {\mathcal A}(Ob)$ with  a given bound on  $ {\mathcal  B}(Ob)$,
that is where $ {\mathcal A}(Ob) = E_{sharp} ({\mathcal  B}(Ob))$.

  \begin{turn}{-90} 
$\leadsto$
\end{turn} \vspace{-5mm}

  \hspace {10mm }  For instance,  the balls in $\mathbb R^n$ are "isoperimetrically extremal": 
  
    \hspace {10mm } they are the {\it only} Euclidean  domains, where  the isoperimetric 
  
   \hspace {10mm } inequality  becomes equality, 
      $vol(V)=\gamma_n (vol_{n-1}\partial V)^\frac {n-1}{n}$.\vspace {1mm}

4. {\it Stability.}   An extremal object  $Ob_{extr}$  is  {\it stable} if convergence  $ {\mathcal  A}(Ob_\varepsilon ) \to   {\mathcal A}(O_{extr})$ and 
${\mathcal  B}(Ob_\varepsilon)\to {\mathcal B}(Ob_{extr}))$  implies that $Ob_\varepsilon\to Ob_{extr}$ in a "suitable sense", where determination of this 
"sense" is the main problem here.

 \hspace {10mm }  For instance, the the balls $B$ in $\mathbb R^n$ are  {\it isoperimetrically stable (modulo translations) with respect to the flat topology, }:
 if $vol(V_\varepsilon )\to vol(B_0)$ and $vol_{n-1}(\partial V_\varepsilon)\to vol_{n-1} (\partial B_0)$,
 then translates $V'_\varepsilon$ of $V_\varepsilon$ converge to $B_0$ in the {\it flat topology}.
 This means  in the present case that  that 
 
 $\bullet$   $vol(B_0\cap V'_\varepsilon)\to vol(B_0)$

 $\bullet$   $vol(B_0\setminus V'_\varepsilon)\to 0$

 $\bullet$ the $\delta$-neighbourhoods of $B_0$  for  $\delta\to_{\varepsilon\to 0} 0$   contain almost all of the boundary of $V'_\varepsilon$, that is 
  $Vol_{n-1}(U_\delta(B_0)) \cap \partial V'_\varepsilon))\to Vol_{n-1}( \partial V'_\varepsilon).$

\vspace {1mm}

Turning to scalar curvature, observe,  that poofs of sharp inequalities, be they   Dirac theoretic or relying on the 
$\mu$-bubble,  are  easily adaptable in most known   cases,  at least   for compact manifolds, for identification of rigid objects, such as

(i) Riemannian {\it flat metrics} $g_{extr}=g_{fl}$  for the inequality $inf Sc(g)\geq 0$ on the torus, 

(ii)  metrics $g_{extr}=g_{sph}$  with constant curvature one  on the $n$-sphere
 $S^n$   for the inequality 
 $\inf Sc(g)\geq n(n-1) $ for metrics $g\geq g_{sph}$ on $S^n$.

However, the following two questions remain unsettled.

{\it Problem  1.} Fully  describe in the case (1)   metric $g_\varepsilon$, $\varepsilon>0$, on the $n$-torus with 
$Sc(g_\varepsilon) \geq -\varepsilon\to 0$, and, in the case (ii), metrics $g\geq g_{sph} $    on $S^n$
with $Sc(g_\varepsilon)\geq n(n-1)-\varepsilon.$

{\it Problem  2.} Find a minimal set of reasonable  additional conditions on $\varepsilon$, such that the metrics $g_\varepsilon$ would converge to $g_{extr}$

The following example indicates what can be expected in regard to  problem  1.

 {\it  Bubble-Convergence.}  
   Let $X$ be a Riemannian $n$-manifold, $n\geq 3$,  and let $X_i=X_{N_i, \varepsilon_i}$,  $\varepsilon_i>0$, be the connected sum of $X$ with  closed Riemannian manifolds $X_{i,j}$, $j=1,2,,,,.N_i$,    where the connected sum is  realized by $\varepsilon_i$-thin surgery  localised at $N_i$  disjoint $\varepsilon_i$-balls 
  $B_{i,j}=B_{ij} (\varepsilon_i)\subset X$, $j=1,...N_i$.

  If  $\varepsilon_i\to 0$, then, this is geometrically clear,  that {\vspace {1mm}
  
   \hspace{6mm} $X$ {\it "emerges" from the sequence  $X_i$ in the limit for} $i\to \infty$,  {\vspace {1mm}
  
  \hspace{-6mm}where     "emerges" 
 becomes   "{\it Hausdorff converges}" if $diam(X_i)\to 0$  and  "{\it converges to $X$ in the intrinsic  flat topology}",\footnote{The definition of this   metric, introduced by  Christina Sormani and  Stefan Wenger in  [Sormani-Wenger(intrinsic flat) 2011], is given in later on in this section.} if
 $$\sum^{N_i}_{i=1}vol(X_i)\to 0. \mbox { }$$

  We explained in section \ref {thin1}  that if $Sc(X)\geq \sigma$ and   $Sc(X_{i,j})\geq \sigma$, then
  the manifolds $X_i$  "naturally" carry metrics with  $Sc(X_i)\geq \sigma -\epsilon_i$,  where   $ \epsilon_i\to 0$.

More interestingly,  the argument indicated in section \ref{C0-limits3} can be used to show\footnote {If $dim(X)\geq 9$ or if some among manifolds $X_{i.j}$ are non-spin, then one needs  new  not formally published results by Lohkamp and/or by Schoen and Yau on "desingularization" of minimal hypersurfaces.}  that 
\vspace {1mm}

{\sl if the set of the centers $x_{i,j}\in X$  of  all these   balls  is   dense in $X$, i.e. all open sets 
  $U\subset X$  contain some  balls $B_{i,j}$,

  if the distances between the balls are much larger than their radii,  
  
   $$dist(B_{i, j_1},B_{i, j_2})/\varepsilon_ i\to \infty   \mbox { for }  i\to \infty$$
   and if  the scalar curvatures of  the   manifolds $ X_{\\i,j}$ are bounded from below, 
  $Sc(X_{i,j})\geq \sigma$,

\hspace {39mm}{\it \color {blue!49!black}then      $Sc(X)\geq \sigma$.} }

  \vspace {1mm}

(One doesn't   need here any bounds  on the diameters and/or volumes of $X_i$, and, {\color {red!49! black}\large \sf probably}, the lower bound on the distances between   $B_{i,j}$ is  redundant.) 
    \vspace {1mm}

{\it Intrinsic Flat Distance.} Given two compact  oriented  $n$-dimensional pseudomanifolds with peace-wise  Riemannian metrics $X_1$ and $X_2$ define $dist_{if}(X_1,X_2)$ as the infimum of the numbers $D\geq 0$ such that there exists an oriented $(n+1)$-dimensional piecewise Riemannian   pseudomanifold $W$ with a boundary, such that

$\bullet$  the oriented boundary of $W$ is $\partial W=X_1\sqcup -X_2$, where the imbeddings 
$X_1, X_2\hookrightarrow W$  are isometric with respect to the {\it distance functions} associated to the Riemannian structures in these spaces;
 
 $\bullet$ $vol_{n+1}(W)\leq d$.    \vspace {1mm}
 
 {\it Remark.} If $X_1$ and $X_2$ are Riemannian {\it manifolds}, then one  can also  take a     Riemannian manifold for $W$, but now with a larger boundary 

 $\partial W=X_1\sqcup -X_2\sqcup X_3$
    and with the condition
     $vol_{n+1}(W)+vol_n(X_3)\leq d$.
 \vspace {1mm}

The following conjecture, in agreement with the Penrose inequality, gives an idea of how wild metrics
with $Sc\geq -\varepsilon$ can/can't be.

{\it Sormani {\color {red!50!black} Conjecture}.}  {\sf Let $X_i$ be a sequence of   Riemannian manifolds homeomorphic to the torus $\mathbb T^3$, such that $$Sc(X_i)\geq -\varepsilon_i\underset {i\to\infty}\to 0.$$ 

{\it If the {\sf volumes} and the  {\sf diameters} of all $X_i$ are bounded by a constant and the {\sf areas of all closed minimal surfaces} in $X_i$ are bounded from below by a {\sf positive} constant},

 {\sf then a subsequence of $X_i$ converges to a flat torus with respect to the intrinsic flat distance in the space of 
Riemannian 3-manifolds}.\footnote {See [Sormani(scalar curvature-convergence)  2016],  [AH-VPPW  (almost non-negative)  2019], [Sormani(conjectures on convergence) 2021],  [Allen(conformal to tori) 2020],  [Pa-Ke-Pe(graphical tori) 2020].} \vspace {1mm}}

{\it Exercise.} Show that the above  condition $\sum^{N_i}_{i=1}vol(X_i)\to 0$ does  imply the intrinsic flat 
convergence $X_i\to X$ as it is claimed in the above example. \vspace {1mm}

{\it Hint.} Use the {\sf \large  {\color {blue!50!black}filling volume inequality:}}\footnote {See [G(filling) 1983],  [Wenger(filling) 2007], [Katz(systolic geometry) 2017].} \vspace {1mm}

{\sl Given a compact Riemannian $n$-manifold $X=(X,g)$, there exists a Riemannian metric $g_\circ$  on the cylinder $W_\circ=X\times (0,1]$, such that:

(i) the metric $g_\circ$ is conical near $0$, 
$$g_\circ(x,t)=t^2dx^2+dt^2, \mbox { for } t\leq \varepsilon=\varepsilon_X>0;$$
 
 (ii) the distance function $dist_{g_\circ}$ on $X=X\times \{1\}\subset W_\circ$ is equal  to $dist_g$;
 
 (iii) the volume of $W_\circ$ is universally bounded by that of $X$ 
$$vol_{n+1}(W_\circ)\leq const_n\cdot vol_n(X)^{\frac {n+1}{n}}.$$}

.\vspace {1mm}

{\it Other kinds of  convergence.} 
Besides the  intrinsic flat there are other distances in the spaces of Riemannian manifolds 
(more or less} adapted  to scalar curvature such as  the {\it directed Lipschitz metric} in section 10 in [G(Hilbert)  2012] and the $d_{p,g}$-distance introduced in  in [Lee-Naber-Neumayer](convergence) 2019] which well goes along with  $Sc\geq -\sigma$ under  a lower bound on Perelman's 
$\nu$-functional.

\vspace{1mm}

Once you have a metric in the space  $\cal X$ of Riemannian manifolds, you are inclined to complete 
this space and study the resulting singular spaces  $X$ from this completion.

Then you isolate the essential properties of these  $X$ and define more general " {\it singular spaces $X$
with $Sc(X)\geq\sigma$"}

Then you dream  of  an abstract category of "objects" with  $Sc\geq \sigma $  that carry the essence of what we know (and don't know) about the scalar curvature.

%%%%%%%%%%%%%%%%%%%%%%%%%%%%%%

\subsection {\color {blue}Who are you, Scalar Curvature?}\label {who3}

%%%%%%%%%%%%%%%%%%%%%%%%%%%%%%
There are two  issues here.

{\color{magenta}1.} What are most general geometric objects that display features similar to these of manifolds with 
positive and  more generally, bounded from below, scalar curvatures?

{\color{magenta}2.}  Is there a  direct link between   Dirac operators and minimal varieties or their joint appearance in the ambience  of scalar curvature is purely accidental? 

\vspace{1mm}

Notice in this regards that  there are  two  divergent branches  of the growing  tree  of scalar curvature. 

{\color {teal} A.} The first one is concerned with  the effects of $Sc>0$ on the {\it differential structure} of  spin (or spin$^\mathbb C$) manifolds $X$, such as their  $\hat\alpha$  and Seiberg-Witten invariants.

{\color {teal} B.}  The second aspect is about coarse geometry and topology of $X$  with $Sc(X)\geq \sigma$, the (known) properties of which are derived by means of minimal varieties  and twisted Dirac  operators;   here the spin condition, even when it is present,  must be  redundant.\vspace {1mm}
  
To better  visualise   separation  between {\color {teal} A} and to {\color {teal} B}, think of possible 
 {\it singular  spaces $X$ with $Sc(X)\geq 0$} corresponding to  {\color {teal} A} and to {\color {teal} B} -- these  must be grossly different.

For instance,  if $X$ is  an  {\it Alexandrov space} with (generalised) sectional curvature  $\geq \kappa>-\infty$ then the inequality  $Sc\geq 0$ makes perfect sense  and, {\color{red!50!black}probably} most (all?) of {\color {teal} B} can be transplanted to these spaces. \footnote{ It seems, much of the geometric measure theory extends to Alexandrov spaces but it is unclear what would correspond to twisted Dirac operators on these spaces.}

But nothing    of  spin related results  makes  sense  for singular Alexandrov spaces. \vspace {1mm}

And  if you start from the position of {\color{magenta}2} you better  go away  from conventional spaces and   start 
 dreaming of   geometric magic  glass ball    with  ghosts  of   harmonic   spinors and of minimal varieties 
dancing   within.  (See  section \ref {beyond manifolds} for continuation of this discussion.)
\vspace {1mm}
 
In  concrete terms one formulates two problems.\vspace {1mm}

{\large \sf \color{blue} A.} {\sf  What is the  largest class of  spaces (singular, infinite dimensional ...) 
which display the basic features of manifolds with $Sc\geq 0$ and/or with  $Sc\geq \sigma>-\infty$
and, more generally, of spaces  $X$, where  the  properly understood   $-\Delta+\frac {1}{2}Sc(X)$ is positive or, at least  not too negative?

For instance, which (isolated) conical singularities and which singular   volume minimising hypersurfaces  belong to this class?
\vspace {1mm}}

\vspace {1mm}

{\large \sf \color{blue} B.} {\sf  Is there a partial  differential equation, or something more  general,  the solutions of which would  mediate between twisted harmonic spinors and minimal hypersurfaces (flags of hypersurfaces?) and which would be  non-trivially linked to scalar curvature? 
 
 Could one  non-trivially couple the twisted Dirac ${\cal D}_{\otimes  L}$ with some  equation $\cal E_L$ on the connections in the bundle $L$ the Dirac    operator  in the spirit of the Seiberg-Witten equation?}\footnote{Natural candidates for $\cal E_L$ are  equations for    critical points  of energy-like  functional on spaces of connections,  where, observe, $L$-twisted harmonic spinors $s:X\to \mathbb S\otimes L$ themselves   minimize $s\mapsto   \int_X\left\langle\mathcal D_{\otimes L}( s(x)) \mathcal D_{\otimes L}(s(x)\right\rangle dx$.}
 %%%%%%%%%%%%%%%%%%%%%%%%%

%%%%%%%%%%%%%%%%%%%%%%%%%
 
\section {Dirac  Operator Bounds on the Size and Shape of Manifolds $X$ with $Sc(X)\geq \sigma$} \label {Dirac4}

%%%%%%%%%%%%%%%%%%%%%%%%%

%%%%%%%%%%%%%%%%%%%%%%
\subsection {\color{blue}  Spinors, Twisted Dirac Operators, and Area Decreasing maps}\label {twisted4}
%%%%%%%%%%%%%%%%%%%%%%%%%%%%

The Dirac  $\cal D$ on a Riemannian manifold $X$   tells you by itself preciously little  about the geometry of   $X$, but the same   $\cal D$ twisted with  vector bundles $L$ over $X$ carries the following message:

\vspace {1mm}

\hspace {18mm} {\sf \large \color {blue}  manifolds with  scalar curvature $Sc\geq \sigma>0$  

\hspace {28mm} can't be too 
large 
area-wise. }\vspace {1mm}

Albeit  the best possible  result of this kind  (due to Marques and Neves, see {\color {blue} B} in section 
\ref {slicing3D.3}, which is known for  $X$ homeomorphic to $S^3$ and  which says that \vspace {1mm}

 \hspace {-3mm}{\it if $Sc(X)\geq 6=Sc(S^3)$, then $X$ can be "swept over" by $2$-spheres of areas $\leq 4\pi$,}

 \vspace {1mm}

 \hspace {-6mm} was proven by means of minimal surfaces,  all  known bounds on "areas" of Riemannian manifolds of  dimensions $\geq 4$  depend on Dirac operators  $\cal D$ twisted (or "non-linearly coupled" for n=4)   with complex vector bundles $L$  over $X$ with unitary connections in $L$, where,  don't forget it,   the very definition of $\cal D$ needs $X$ to be spin.\footnote {Recently,  	Jintian Zhu [Zhu(rigidity) 2019]  and Thomas  Richard [Richard(2-systoles) 2020]  established new kind of  bounds on areas of surfaces applicable to   higher dimensional  non-spin manifolds by using  geometric calculus of variations, but these bounds depend on   lower distance bounds (that may be hidden in the topological assumptions, such as in the Zhu paper)  and  are not sufficient, for instance, to show that the unit sphere  $S^n$ for $n\geq 4 $  admits no metric $g$ with $Sc(g)\geq n(n-1)$ and such that the $g$-areas of all surfaces $\Sigma$ in $S^n$  satisfy $area_g(\Sigma)\geq C\cdot area_{S^n}(\Sigma)$  for arbitrary large $C$. }

Recall %(compare  with section \ref {3.11.1})  
 that the {\it twisted} Dirac  operator, denoted  
$${\cal D}_{\otimes L }:C^\infty(\mathbb S\otimes L )\to C^\infty(\mathbb S\otimes L),$$
acts on the tensor product of the spinor bundle  $\mathbb S \to X$ \footnote{All  you have to know about $\mathbb S(X)$  is that  it is a vector bundle associated with the tangent bundle $T(X)$,  which can be  defined for spin manifolds $X$, where "spin" is needed, since the structure group of  $\mathbb S(X)$ is the double cover of the orthogonal group $O(n)$ rather  than $O(n)$ itself. }  with $L\to X$, where it is related to the (a priori positive Bochner Laplace operator)  $\nabla^2_{\otimes L}=\nabla^2_{\otimes L}=\nabla^\ast_{\otimes L}\nabla_{\otimes L} $  in the bundle $\mathbb S \otimes L$,   by the {\it twisted} Schroedinger-Lichnerowicz-Weitzenboeck formula
$$ \mbox { ${\cal D}_{\otimes L}^2=\nabla^2_{\otimes L} + \frac{1}{4}Sc(X) + {\cal R}_{\otimes L}$},$$
where  {$\nabla_{\otimes L}$ denotes the covariant derivative  in $\mathbb S\otimes L$ and 
${\cal R}_{\otimes L}$ is a certain (zero order)    which acts in the fibers of the  twisted spin bundle 
$\mathbb S\otimes L$  and which is derived  from the curvature of the connection in $L$. 

If we are not  concerned with the  sharpness of   constants, all we have to know is that ${\cal R}_{\otimes L}$  is controlled by 
  $$||  {\cal R}_{\otimes L}|| \leq const\cdot ||curv(L)||$$
for $const=const(n, rank(L))$, where 
a little thought (no computation is needed) shows that, in fact, this constant depends only on $n=dim(X)$. (The actual formula for $\mathcal R_{\otimes L}$ is written down in   the next section, also  see [L-M(spin geometry) 1989]  and   [MarMin(global riemannian) 2012] for further details and references.) 

 We regard a closed orientable  {\it even dimensional} Riemannian manifold  $X$  {\it area wise \color {blue} large}, if  it carries a {\it \color {blue} homologically substantial}  or {\it  essential}  {\it \color {blue} bundle} $L$ over it with {\it {\color {blue} small} curvature}, where
"homologically substantial" signifies that some Chern number  of $L$ doesn't vanish.
 It is easy in this case\footnote { See {\color {blue}$(L^\wedge)$} in section \ref {fredholm4} and references therein.} that there exists  an associated bundle $L^\wedge$, such that  
$$|curv|(L^\wedge)\leq const_n |curv|(L)$$   
and such that  the Chern character  in the index  formula 
guaranties   non-vanishing 
  of the cup product
$ \hat A(X)\smile Ch(L^\wedge)$
evaluated at $[X]$,
$$ (\hat A(X)\smile Ch(L^\wedge))[X]\neq 0$$
and, thus, by Atiyah-Singer theorem, the presence of {\it non-zero harmonic twisted spinors}: 
{\sf sections $s$ of the bundle $\mathbb S\otimes L^\wedge$  for which ${\cal D}_{\otimes L^\wedge}(s)=0.$ }

\vspace {1mm}

 If the dimension $n$ of $X$ is odd, the above applies to $X\times  S^1$  for a sufficiently long circle $S^1$.
 
For instance,  $n$-manifolds, which 
 admit  area decreasing non-contractible   maps to  spheres $S^n(R)$ of  large radii $R$  are area-wise large, where 
the relevant bundles $L$  are  induced from non trivial bundles over the spheres. (One may take $L^\wedge= L$ for these $L$.)

 \vspace{1mm}

But   if the scalar curvature of $X$ is $\geq \sigma$ for a  {\color {blue}\it  large} $\sigma>0$, where this   {\color {blue}\it "large"} properly matches the above  {\color {blue}\it "small"}, then by  the Schroedinger-Lichnerowicz-Weitzenboeck formula the   ${\cal D^2}_{\otimes L^\wedge}$  is positive and no such harmonic twisted spinor  exists; therefore, a suitably defined "area"$(X)$ must be bounded
by $\frac {const } {\sigma}$. (See the  sections \ref {K-area3}, \ref{cowaist4}}  for a definition of this "area" 
called  K-{\it area} and K-{\it  cowaist}.) 
\vspace{1mm}

Next, recall that the  $\hat  A$-genus, 
$$\hat A(X)=1-\frac {1}{24}p_1+\frac{1}{5760}(-4p_2+7p_1^2)+... \in H^\ast(X)$$ 
is a certain  polynomial in Pontryagin classes $p_i\in H^{4i}(X)$  of $X$ and 
$$Ch(L)=rank_{\mathbb C}(L)+c_1(L)+\frac{1}{2}(c_1(L)^2-2c_2(L))+...\in H^\ast(X)$$
 is a polynomial in  Chern classes $c_i(L)\in H^{2i}(X)$ of $L$, while $[X]\in H_n(X)$ denotes the fundamental class of $X$.

If $n=dim(X)$ is even, the spin bundle $\mathbb S$ naturally splits,   
$\mathbb S= \mathbb S^+\oplus \mathbb S^-$, the   ${\cal D}_{\otimes L}$ also splits:
${\cal D}_{\otimes L}={\cal D}^+_{\otimes L}\oplus {\cal D}_{\otimes L}^-$, for 
$${\cal D}^\pm_{\otimes L}: C^\infty(\mathbb S^\pm\otimes L)\to C^\infty(\mathbb S^\mp\otimes L)$$
and the index formula reads:
$$ind({\cal D}^\pm_{\otimes L})=\pm (\hat A(X)\smile Ch(L))[X].)$$
\vspace {1mm}

{\it \color {blue} Relative Index Theorem on Complete Manifolds.} Let  $X$ be a  {\it complete} Riemannian manifold the scalar curvature of  which is  {\it uniformly positive at infinity}.\footnote {It is shown in [Zhang(Area Decreasing)  2020] that   "uniformly"can be dropped --  "positive at infinity" suffices.}   Then the Schroedinger-Lichnerowicz-Weitzenboeck formula implies that

\hspace {5mm}  {\it the Dirac   operator is positive at infinity, i.e.  outside some  compact subset $V\subset X$}: 
$$\int_X \langle \mathcal D^2s(x), s(x)\rangle dx  \geq \varepsilon \int_X \mathcal ||s(x)||^2 dx$$
for some $\varepsilon= \varepsilon (X)>0$ and all $L_2$-spinors $s$ supported outside $V$. 
This  (easily)  implies, in turn, that the operators  $\mathcal D^\pm$ are Fredholm but  the indices of these operators  depend on delicate   information on geometry of $X$ at infinity  and no simple formula for 
$ind(\mathcal D^\pm)$ is available. 

However if there are two operators  $\mathcal D_1$ and  $\mathcal D_2$, which are {\it equal  at infinity}, e.g. $\mathcal D_1=\mathcal D^+_{\otimes L}$,   and $\mathcal D_2=\mathcal D^+_{\otimes |L|}$, where $L\to X$ is a bundle with a unitary connection, where  $|L|$  is    the trivial bundle of rank $k=rank_{\mathbb C}L$ over $X$ and where $L$ comes with an isometric connection preserving isomorphism with $|L|$ at infinity, 
as in section \ref {relative3}, then  the difference of their indices  -- both are Fredholm for the same reason as $\cal D^\pm$-- satisfy the Atiyah-Singer formula: 
$$ind({\cal D}^+_{\otimes L})-ind({\cal D}^+_{\otimes |L|})= (\hat A(X)\smile (Ch(L)-Ch|L|))[X].$$
where, 
$$Ch(L)-Ch|L|=c_1(L)+\frac{1}{2}(c_1(L)^2-2c_2(L))+...$$
is understood as a cohomology class with compact supports and $[X]$ is the fundamental homology class with infinite supports.

More generally, if  $\mathcal D_i  =\mathcal D_{\otimes L_i}$, $i=1,2$, where $L_1$ is equated with $L_2$ at infinity, then 
$$ind({\cal D}^+_1)-ind({\cal D}^+_2)= (\hat A(X)\smile (Ch(L_1)-Ch (L_2))[X],$$
where one needs the operators $\mathcal D_i$ be positive at infinity. 
\vspace {1mm}

{\it The proof } of   this can be obtained by adapting any version of the {\it local proof} of the compact 
Atiyah-Singer theorem (see (see [GL(complete) 1983], [Bunke(relative index) 1992],  [Roe(coarse geometry) 1996]).

Namely,   the index is represented by the difference of the traces of  families of auxiliary operators 
  $ K^+_{1,t} -K^+_{2,t}$  and $K^-_{1,t} -K^-_{2,t}$, $t>0$, where \vspace {1mm}
  
  \hspace {-6mm}(i)\hspace {3mm}  these $K_{..., t}$-s are given by   continuous kernels $K_{...,t} (x, y)$ which are
  
  \hspace {2mm}supported in the $t$-{\it neighbourhood of the diagonal}  in $X\times X$, i.e.  where $dist (x,y)\leq t$;\vspace {1mm}
  
   \hspace {-6mm}(ii)\hspace {3mm} $K^\pm_{1,t} (x,y) =K^\pm_{2,t}(x,y)$  for $x$ and  $y$ in the {\it complement of a compact subset} 
   
    \hspace {2mm}$V_t\subset X$, where $V_{t_1}\subset V_{t_2} $ for $t_2>t_1$ and where $\bigcup_tV_t=X$;\vspace {-4mm}

 $$ trace (K^+_{1,t} -K^+_{2,t})-  trace (K^1_{1,t} -K^1_{2,t}) =(\hat A(X)\smile (Ch(L_1)-Ch (L_2))[X];\leqno {\rm (iii)}$$
for all $t>0$;\vspace {1mm}
 
   \hspace {-6mm}(iv)\hspace {3mm} the operators $K^\pm_{i, t}$, $i=1,2$, {\it weakly converge}\footnote{ The corresponding functions $K_{...,t}(x,y)$ uniformly converge on compact subsets in $X\times X$.} for $t\to \infty $ to  the projection 
   
    \hspace {-5mm}operators on the kernels of 
$\mathcal D\pm_i$.

The  quickest  way to get such $K_{..., t}$ is by taking suitable functions $\psi_t$ of the corresponding Dirac  operators, where the Fourier transforms of $\psi_t$ have {\it compact supports}, and where (as in all arguments of this kind)  the essential issue  is the proof   
of  {\it uniform bounds} on the traces of  the operators 
$K^\pm_{1,t} -K^\pm_{2,t}$ for $t\to \infty$.

Specific  bounds for particular  $K_{...,t}$  are crucial for an (approximate)  extension of the  index theory   to non-complete manifolds, %(e.g. needed for the problem discussed in  section  \ref {4.6}) 
 but these  bounds are often buried  in  the $K$-theoretic  formalism  of the recent papers. Also, I must admit, this point was not explained (overlooked?) in the exposition of Roe's argument  in my  paper [G(positive) 1996].

%%%%%%%%%%%%%%%%%
\subsubsection {\color {blue} Negative Sectional Curvature against Positive Scalar Curvature}
 \label {negative sectional4}
%%%%%%%%%%%%%%%%%%%%%%

A characteristic  topological corollary of the above  is as follows.\vspace{1mm}

{\large \color {blue} $[\kappa\leq 0] \leadsto [Sc \ngtr 0]$:} {\sl If a  closed  orientable spin $n$-manifold $X$ admits a map to a complete Riemannian  manifold
 $\underline X$ with {\color {blue} $sect.curv(\underline X)\leq 0$}, 
 $$f: X\to \underline X,$$
such that the  homology  image $f_\ast[X]\in H_n(\underline X;\mathbb Q)$
doesn't  vanish, then $X$ admits no metric with $Sc(X)>0$.} \vspace {1mm}

{\it Two Words about the Proof.} All we need of $sect.curv\leq 0$ is the existence  of   distance decreasing maps from the universal covering of ${\underline X}$ to  (large)  spheres,
 $$F_ {\underline x}:  \underline X \to S^{\underline n}(R), \hspace{1mm} \underline n=dim (\underline X),\hspace{1mm}  {\underline x}\in {\underline X},$$
which can be (trivially) obtained with a use of the inverse exponential maps 
$$\exp^{-1}_{\underline x}:  \tilde {\underline X} \to T_{\underline x}(\underline X),\hspace{1mm}  {\underline x}\in {\underline X}.$$

To make the idea clear,  let $\underline X$ be compact,    the fundamental group of $\underline X$ be residually finite, (e.g. $\underline X$ having constant sectional curvature or, more generally  being a  locally symmetric space)   and  $X$ be embedded to $\underline X$. 

Let $X^\perp \subset \underline X$ be a closed oriented submanifold of dimension $m=\underline n-n$ for $\underline n=dim(\underline X)$, which has {\it non-zero intersection index} with $X\subset  \underline X$. 

Also assume that the restriction of the tangent bundle of  $\underline X$ to  $ \underline X^\perp \subset \underline X$
is trivial.

Then --  this is rather obvious -- there exist finite covers $\tilde {\underline X_i}\to \underline X$,  such that the products of the lifts (i.e. pull-backs) of  $X$ and of $X^\perp $ to $\tilde {\underline X}_i$, denoted  $\tilde X_i\times  \tilde X^\perp_i $, admit smooth maps to the spheres of radii $R_i$, 
$$F_i:\tilde X_i\times  X_i^\perp\to S^{\underline n}(R_i),$$ 
where

$\bullet_1 $ $R_i\to \infty $,

$\bullet_2 $   $deg(F_i)\neq 0$,

$ \bullet$  the maps $F_i$ are {\it distance  decreasing on the fibers $\tilde X_i\times x^\perp$}  for all  $x^\perp\in   X_i^\perp$
for the Riemannian  metric  in these fibers induced by the embedding 
$\tilde X_i\times x^\perp=\tilde X_i\subset \tilde {\underline X_i}$.

It follows  that for {\it arbitrary} Riemannian metrics   $g$ and $g^\perp$ on $X$ and on  $X^\perp$ there exists   (large)  constants $\lambda $ and $C$
independent of $i$, such that 
 
\vspace{1mm}

{\sf  the maps $F_i$ are {\it $C$-Lipschitz} with respect to the sum of the lift of the  metric $g$ to $\tilde X_i$  and the lift of  $\lambda\cdot g^\perp$ to  $\tilde X^\perp_i $ that is the metric 
$$\tilde g_i\oplus \lambda \cdot\tilde g_i^\perp\mbox { on } \tilde X_i\times  \tilde X^\perp_i.$$}

If $Sc(g)\geq \sigma>0$, then also  $Sc(\tilde g_i\oplus \lambda \cdot\tilde g_i^\perp)\geq \sigma'>0$ for all sufficiently large $\lambda$, which, for large $R_i$,  rules out non-zero harmonic spinors on 
$\tilde X_i\times  \tilde X^\perp_i$ twisted with the bundle $L^\ast=F_i^\ast(L)$ induced from any given bundle $L$ on  $S^{\underline n}$.

But if $\underline n=2k$ and  the  Chern class $c_k(L)$ is non-zero, then  non-vanishing of $deg(F_i)$  implies  non-vanishing of
of $ind({\cal D}_{\otimes L})$  via the index formula and the resulting contradiction delivers the proof for even $\underline n$ and the odd case follows with $\underline X\times S^1$.

\vspace {1mm}

{\it Remarks.}   This argument, which is rooted in Mishchenko's proof of Novikov conjecture for the fundamental group of the above $\underline X$,  which was  adapted to scalar curvature in [GL(complete) 1983]  and further   generalized/formalised in  [Rosenberg($C^\ast$-algebras - positive  scalar) 1984],  and  [CGM(Lipschitz control) 1993], doesn't really need compactness of $\underline X$,   residual finiteness of $\pi_1(\underline X)$ and  triviality of  $T(\underline X)|X^\perp$.  Beside,  the spin condition for $X$ can be   relaxed to that for the  universal cover of $X$.
  
  Moreover,   since   the bound on the size of $\tilde X_i\times \mathbb  T^{\underline n-n}$ by $\frac {const}{\sqrt\sigma}$
can be obtained with the use of minimal hypersurfaces (see \S 12 in  [GL(complete) 1983]), [G(inequalities) 2018] and 
section \ref{separating5}) the spin condition can be dropped altogether.

{\it\large \color {blue} Question.} {\sf Are there other topological non-spin obstructions to $Sc>0$?} 

For instance, is the following  true?\vspace{1mm}

 {\it\large \color {blue} {\color {red!50!black} Conjecture}.}  {\sf  Let $X$ be a  closed orientable Riemannian  $n$-manifold, such that no closed orientable $n$-manifold $X'$ which admits a map $X'\to X$ with non-zero degree carries a metric with $Sc>0$. Then there exists an integer $m$ and 
 a sequence of maps 
$$F_i :\tilde X_i\times \mathbb R^m\to S^{n+m}(R_i),$$
 where $\tilde X_i$ are   (possibly infinite) coverings of $X$, such that

$\bullet$  the maps $F_i$  are constant at infinity and they  have non-zero degrees, 

$\bullet$ $R_i\to \infty$,

$\bullet$ the  maps $F_i$ are distance decreasing on the fibers $\tilde X+i\times x^\perp$ for all $x^\perp \in \mathbb R^m$.}\footnote{See   [Dranishnikov(asymptotic) 2000],    [Dranishnikov(macroscopic) 2010],  [DFW)flexible) 2003],    [Dranishnikov( hypereuclidean) 2006]  [BD(totally non-spin) 2015]  and references therein  for relations between various   largeness conditions (e.g. of  universal covering of compact  manifolds) and their roles in the proofs of  the   Novikov conjecture  and of non-existence of metrics with  $Sc>0$.}

\vspace{1mm}

Apparently,  there is  no  instance of a {\it specific} homotopy class $\cal X$  of closed manifolds $X$ of dimension $n\geq 5$,  
where a Dirac theoretic  proof of non existence of metrics with $Sc>0$ on all $X\in \cal X$ couldn't be  replaced by  
a  proof via minimal hypersurfaces. 

\vspace {1mm}

(This  seems  to disagree  with what was said concerning  the "quasisymplectic theorem" {\color {blue} $\otimes_{\wedge \omega}$}   in section \ref{SY+symplectic2}.

In fact the {\it general} condition for $Sc\ngtr 0 $  in  {\color {blue} $\otimes_{\wedge \omega}$},   can't be treated, not  as it stands,  with minimal hypersurfaces, but this may be possible in all 
{\it \it specific examples}, where this condition  was {\it proven to be}  fulfilled.)

And  it is conceivable when it comes to the Novikov conjecture,  that its validity in all  proven specific examples, can be derived by an elementary argument from  the  invariance of rational 
Pontryagin classes under  $\varepsilon$-homeomorphisms.\footnote{The original proof  of topological invariance of Pontryagin classes by Novikov, as well as  simplified versions and modifications   of his proof in [G(positive) 1996)  automatically apply to $\varepsilon$-homeomorphisms and, sometimes,  to homotopy equivalences})

But even though the relevance of twisted   Dirac theoretic methods is questionable as far as {\it topological} non-existence  theorems are concerned, these methods seem irreplaceable  when  it comes to {\it geometry} of  $Sc\geq \sigma$.

%%%%%%%%%%%%%%%

\subsubsection {\color {blue} Global  Negativity of the  Sectional  Curvature,  Singular Spaces with $\kappa\leq 0$,  and Bruhat-Tits Buildings}\label {Tits4}

%%%%%%%%%%%%%%%%%%%%

The essential feature of complete spaces $\kappa\leq 0$ (these   often come under heading of CAT(0)-spaces) needed for   $[Sc \ngtr 0]$  is as follows.

  \vspace {1mm}

{\color {blue} \large [$\circlearrowleft_\varepsilon$]  \it Self-contraction Property.}  {\sl  $\underline X$ admits a family of   proper  $\varepsilon$-Lipschitz selfmaps  $\phi_{\varepsilon}:  \underline X\to  \underline X$,  for all $\varepsilon>0$, where these 
maps are properly homotopic to the identity map $id$. }\footnote {See [G(large)  1986]  for more about such manifolds.}\vspace {1mm}

If  $ \underline  X$ is a topological $n$-manifold, than  this property implies the existence of proper   Lipschitz maps $ \underline  X\to \mathbb R^n$ of degree one, but unlike the latter it 
makes sense for singular spaces that are not topological manifolds or pseudomanifolds.

On the other hand, if a possibly singular, say  finite dimensional polyhedral space $X$ satisfies {\color {blue} \large $\circlearrowleft_\varepsilon$}, then there exists a manifold $ \underline  X^+\supset  \underline  X$, which also satisfies {\color {blue} \large $\circlearrowleft_\varepsilon$},  where the most transparent case is that of spaces 
$  \underline X$ which come with   free isometric actions by  discrete groups $\Gamma$ with compact quotients   $\underline  X$.

To derive   $\underline  X^+$  from $\underline X$ in this case, embed
$\underline X/\Gamma \hookrightarrow \mathbb R^N$, take a small 
   regular neighbourhood $U\subset \mathbb R^N$ of $\underline  X^+/\Gamma\subset \mathbb R^N$ and let $\tilde U \to U$.  be   the universal covering of $U$.  
   
   Then this $\tilde U$ with  a {\it suitably blown-up metric}  serves for  $\underline  X^+$, where the simplest such blow up is achieved by multiplying the (locally Euclidean) metric in $\tilde U$  by the  function 
   $\frac {1}{dist(\tilde u, \partial  U)}.$

In fact, what is truly needed for the non-existence argument, and what is satisfied by complete simply connected spaces $\underline  X  $  with $\kappa<0$ is the following parametric version of  {\color {blue} \large [$\circlearrowleft_\varepsilon$]}.

\vspace {1mm}

{\color {blue} \large [$\circlearrowleft_\varepsilon${\color {blue} \large $\circlearrowleft_\varepsilon$]}}. {\sf   There exist a continuous map    $\Phi_{\varepsilon}:  \underline X\times  \underline X \to  \underline X$  with the following properties.}

{\sl $\bullet_\varepsilon $ the  maps $\phi_{\varepsilon, \underline x_0} = \Phi_{\varepsilon}:  \underline X=\underline x_0\times  \underline X\to  \underline X$
are proper $\varepsilon$-Lipschitz  for all $\underline x_0\in \underline X$  and all $\varepsilon>0$;\vspace {1mm}

$\bullet_n $ the restrictions of these maps $\phi_{\varepsilon, \underline x_0}:\underline X\to \underline X$ to the $n$-skeleton   $\underline X^{(n)}\subset \underline X$ are proper homotopic to the inclusions $\underline X^{(n)}\subset \underline X$;\footnote{Here we assume  that $ \underline X$ is  triangulated and $n$ denotes the dimension of a manifold $X$ we are going to map to $ \underline  X$;}
\vspace {1mm}

$\bullet_\Gamma $  the family  $\phi _{\varepsilon,\underline x_0} $ is  equivariant under the isometry group of $\underline X$: \vspace {1mm}

{\it  if $\gamma:\underline X\to \underline X$ is an isometry, then
  $$\phi_{ \varepsilon,  \gamma(\underline x_0)}   =\gamma\circ \phi _{\varepsilon,\underline x_0}. $$}}

The above argument combined with that in the previous section yields the  following  generalization of the non-existence theorem {\large \color {blue} $[\kappa\leq 0] \leadsto [Sc \ngtr 0]$.}
\vspace {1mm}

{\large \color {blue} $[\kappa\leq 0]_{global} \leadsto [Sc \ngtr 0]$:} { {\sl If  a complete Riemannian spin manifold $\tilde X$  of dimension $n$  with a discrete (not necessarily free) co-compact  isometric  action of a group $\Gamma$ admits a proper $\Gamma$-equivariant  map to an $\underline X$ which satisfies {\color {blue} \large $\circlearrowleft_\varepsilon${\color {blue} \large $\circlearrowleft_\varepsilon$}},
then $\inf_x(Sc(X,x)\leq 0$.}
\vspace {1mm}

{\it \color {blue} Corollary.} {{\sl Let $\Gamma$ be a finitely generated  subgroup in the linear group 
 $GL_N(\mathbb C)$,\footnote {One may place here any field   instead of $\mathbb C$.} let $X$ be a compact oriented  Riemannian spin $n$-manifold 
with $Sc(X)>0$ and let $f: X\to {\sf B}(\Gamma)$ be a continuous map,  where ${\sf B}(\Gamma)$ denotes the classifying (Eilenberg MacLane) space of $\Gamma$. 

Then the  image 
$$f_\ast [X]_{\mathbb Q} \in H_n({\sf B}(\Gamma); {\mathbb Q})$$  
of the rational fundamental class 
$$\mbox{ $[X]_{\mathbb Q} \in   H_n(X; \mathbb Q)$    for $ f_\ast :  H_\ast(X; \mathbb Q)  \to H_\ast({\sf B}(\Gamma); \mathbb Q)$} $$ 
is zero.\footnote{A more sophisticated  theoretic  version of this in the context  of the Novikov conjecture appears in [Kasp-Scan (Novikov) 1991]. }}

\vspace {1mm}

{\it Proof.} A finite index subgroup 
in $\Gamma$  freely,\footnote{Finite index was needed  for his "freely"} discretely and isometrically acts on the product 
$\underline X$   of Riemannian symmetric spaces and {\it Bruhat-Tits buildings}, where such  products, according to Bruhat-Tits are   \vspace {1mm}

\hspace {6mm}{\it complete  simply connected polyhedral space  with $\kappa(X)\leq 0$.}\vspace {1mm}

\hspace {-6mm} Since  \large $\circlearrowleft_\varepsilon${\color {blue} \large $\circlearrowleft_\varepsilon$}}  apply to such spaces,   the proof of the corollary  follows.

\vspace {1mm}

 {\it \color {blue}  Historical Remark.}  Around 1950,  A.D. Alexandrov,  H. Pedersen and Busemann who suggested (two somewhat different) definitions of $\kappa\leq 0$  applicable to {\it  singular} metric  spaces, and  their followers   focused on essentially local geometric properties of these spaces $X$, and tried to {\it alleviate effects of singularities} 
by adding extra assumptions on $X$. \footnote  {A  brief overview of this circle of  ideas is given in section 2.3 of [G(hyperbolic)2016]   and contributions by the Alexandrov's school are presented  in [AKP(Alexandrov spaces) 2017].}

The theory of  $\kappa\leq 0$  has acquired a global mathematical status  in early seventies  with the discoveries of    Bruhat-Tits buildings (1972)\footnote {Bruhat and Tits independently developed the local and global  theory of their  spaces  being  unaware of definitions of $\kappa\leq 0$   suggested  by differential geometers.} 
and the   link of   $\kappa\leq 0$  with the index theory and the Novikov conjecture by   Mishchenko(1974).

This  has eventually  led to the modern perspective on  CAT(0)-spaces, i.e. those with $\kappa\leq 0$,  the main  interest in  which is due to a  multitude  of significant    examples of   singular CAT(0)-spaces with interesting  fundamental groups  inspired by the ideas behind the construction(s) and applications  of  the
 Bruhat-Tits buildings.

\vspace {1mm}

{\it \color {blue} Hyperbolic Remark}.  "$\varepsilon$-Lipschitz" in the theorem {\large \color {blue} $[\kappa\leq 0]_{global} \leadsto [Sc \ngtr 0]$}  is only needed on the {\it large scale}, that is expressed by 
the inequality $$dist(f_{\varepsilon, x_0}( x_1),   f_{\varepsilon, x_0}(x_2))\leq \varepsilon dist(x_1,x_2)+const.$$}

Thus, for instance, \vspace {1mm}

{\sf the non-existence conclusion for  metrics with $Sc>0$  on $X$ applies, where $\underline X$ is the {\it Vietoris-Rips complex }of a {\it hyperbolic group}}. 

It follows, that the conclusion of the above corollary holds for hyperbolic groups $\Gamma$:   \vspace {1mm}

{\sl Let $X$ be a closed orientable Riemannian spin manifold with  $Sc(X)>0$  and let  $\Gamma$ be  a hyperbolic group.  Then the class   $f_\ast [X]_{\mathbb Q}\in H_n({\sf B}(\Gamma);\mathbb Q)$ vanishes for all continuous maps $f:X\to  {\sf B}(\Gamma)$.}

 \vspace {1mm}

{\it \color {blue}   "$\varepsilon$-Area"  Remark.}  Instead of "$\varepsilon$-Lipschitz" one may  require   "{\it $\varepsilon$-area contracting}" or some  large scale counterpart to this condition.
 
 This may be  significant, because the    $\varepsilon$-area version of  {\large \color {blue} $[\kappa\leq 0]_{global} \leadsto [Sc \ngtr 0]$}  is  not-approachable  with the (known) techniques of minimal hypersurfaces and/or of stable  $\mu$-bubbles, 
 while the above  "$\varepsilon$-Lipschitz" {\large \color {blue} $[\kappa\leq 0]_{global} \leadsto [Sc \ngtr 0]$} can be proved  in many, probably in all, cases with these techniques  having   an advantage of not requiring manifolds $X$ to be spin.
 
 On the other hand, for all I know,  there is no example of an $\underline X$, say with a cocompact  action of an isometry group $\Gamma$, 
 which satisfies  a version of {\color {blue} \large $\circlearrowleft_\varepsilon${\color {blue} \large $\circlearrowleft_\varepsilon$}}  with  the $\varepsilon$-contracting area  property but not with the $\varepsilon$-Lipschitz one.\footnote {Neither, it seems, there are examples of $\underline X$ with compact quotients  $\underline X/\Gamma$, which   satisfy {\color {blue} \large $\circlearrowleft_\varepsilon$} but not {\color {blue} \large $\circlearrowleft_\varepsilon${\color {blue} \large $\circlearrowleft_\varepsilon$}}.}
  
 \vspace {1mm}

%%%%%%%%%%%%%%%%%%%%%%
\subsubsection {\color {blue} Curvatures of Unitary  Bundles, Virtual Bundles and Fredholm Bundles}\label {fredholm4}
%%%%%%%%%%%%%%%%%%%

Let us try to formalise the concept of    \vspace{1mm}

{\sf "area",  of a Riemannian manifold $X$, where this "area"  is associated with curvatures of vector bundles over $X$ and which  has the property of being bounded   \vspace{1mm}
 by $const\cdot \frac {1}{\sigma}$, for  $\sigma=\inf_xSc(X,x)>0$.} \vspace {1mm}

{\large \sf \color {blue} $||curv(L)||.$ } Given a   vector bundle $(L,\nabla)$ with an orthogonal (unitary in the complex case)  connection, over a Riemannian manifold
$X$,   let  
$$ ||curv(L)||(x)=||curv(\nabla)||(x)=||curv(L, \nabla)||(x)$$ denote\vspace {1mm}

{\it the infimum of positive functions  $C(x)$, such that the maximal {\it rotation angles} $\alpha\in [-\pi, \pi]$ of the   parallel transports  along the boundaries of  smooth discs $D$ in $X$ satisfy
 $$ |\alpha|=|\alpha_D|\leq\int_D C(d).\footnote{This definition is adapted to vector bundles over rather  general metric spaces, e.g. polyhedra with piecewise smooth metrics.}$$}

(The holonomy  splits into the direct sum of  rotations $z\mapsto \alpha_iz$, $z\in \mathbb C$, $\alpha_i\in \mathbb T\subset \mathbb C$,  $i=1,2,...,rank(L)$, and  our $\alpha =\max_i\alpha_i$.)

For instance, if $D$ is  a geodesic   digon in $S^2$ with the angles $\beta \pi$, $\beta \leq 1$,
then the holonomy of the tangent vectors around the boundary of $D$ satisfies:
$$ |\alpha_D|=2 \beta\pi =area (D),$$
 which agrees with the equality $|curv|(T(S^2))=1$.

It follows that the curvature of the tangent bundle (complexified if you wish) of the product of spheres, 
 satisfies
$$\left |\left |curv\left (T\left (\bigtimes_iS^{n_j}(R_j)\right )\right)\right |\right |=\frac{1}{\min_jR_j^2}.$$
 
What is more amusing is that the even dimensional spheres  $S^n$, $n=2m$, support unitary bundles $L$
with with twice smaller curvatures and  {\it non-zero} top  Chern classes, 
$$\mbox {  $|curv|(L) =\frac  {1}{2}$  and $c_m(L)\neq 0$.}$$  

For instance, if $n=2$, then the  Hopf bundle, that is the square root of the tangent bundle,  has these  properties and  in general, the positive $\mathbb C$-spin bundle $\mathbb S^+$ can be taken for such an $L$. 

This is the smallest curvature a non-trivial bundle over $S^n$ may have:\vspace{1mm}

\hspace {5mm}  {\sl Unitary vector bundles over $S^n$ with $|curv|<\frac {1}{2}$ are trivial. }

 \vspace{1mm}
{\it Proof.}  Follow the parallel transport of  tangent  vectors  from the north to the south pole.

More generally    \vspace{1mm}

{\sf there are bundles $L$ on 
the products of  even dimensional spheres  $\bigtimes_iS^{n_j}(R_j)$,  which  are induced by $\lambda$-Lipschitz maps 
to $S^n$, $n=\sum n_j$, $\lambda = \frac{1}{\min_jR_j^2}$, such that $|curv|\leq \frac{1}{2\min_jR_j^2}$
 and such that  {\it some Chern numbers  of these $L$ are non-zero}, and this is the best one can do. }

In fact, \vspace{1mm}

{\sl If a unitary   vector bundle $L=(L, \nabla) $ over a product manifold  $S^n\times Y$ has  $|curv|(L)<\frac {1}{2}$, then all Chern
numbers  of $L$ vanish.} (see \S 13 in [G(101) 2017]).  \vspace{1mm}

The role of the  Chern numbers here is motivated by the following observation (see [GL (spin) 1980,  [G(positive) 1996]).

Let $X$ be a closed orientable  spin manifold of dimension $n=2m$  and $L=(L, \nabla)$  a unitary vector  bundle, such that  {\it some Chern  number of $L$ doesn't vanish}. Then \vspace{1mm}

{\large \color {blue}$(L^\wedge)$ } {\sl there exists an associated bundle $L^\wedge$, which  is a polynomial 
in the  exteriors powers of $L$,  such that 
$$ind(\mathcal D_{\otimes L^\wedge})\neq 0$$.}
\vspace{1mm}

Since (it is easy to see) the degree and the coefficients of such a polynomial must be bounded by a constant depending only on $n$, the curvature of  $L^\wedge$ satisfies 
$$ |curv|( L^\wedge) \leq const_n \|curv|(L);$$
 Therefore, \vspace{1mm}

{\Large \color {blue} $\bullet$} {\sl if the  scalar curvature of a closed orientable $2m$-dimensional spin manifold satisfies $Sc(X)\geq \sigma >0$, then}  -- this is explained in the previous section  --  {\sl non-vanishing  $c_m(L)\neq 0$,  implies the following lower bound on
the curvature of the bundle $L$: 
$$|curv|(L)\geq {\epsilon}\cdot{\sigma},\mbox { } \epsilon=\epsilon(n)>0.$$}

{\it \large \color {magenta!50!black} Open  Problem}. {\sf Prove {\Large \color {blue} $\bullet$} without the  spin condition.} 

\vspace{1mm}

 The above suggest the definition of "area"$(X) $ of a Riemannian manifold $X$  as the supremum of $\frac {1}{|curv|(L)}$  over all  unitary vector bundles $(L=L, \nabla)$ with non-zero Chern numbers.
 
 However, the  "area" terminology we  introduced in [G(positive)  1996], despite several natural/functorial properties of this "area"  (see  [G(positive)  1996] and  [G(101 2017]),   seems inappropriate, since this "area" is {\it by no means additive}. A more  adequate word , which we prefer to use from now on   is   {\it K-cowaist.}\vspace {1mm}
 
 {\it Virtual  Hilbert and Fredholm.} To define this, we  represent the {\it (Grothendieck) classes}  $\mathbf h$ of vector bundles   over $X$, which are also called {\it virtual (Fredholm)  bundles}, 
 by {\it Fredholm homomorphisms} between Hilbert bundles with unitary connections $\mathcal L_i= (\mathcal L_i,\nabla_i)$, $i=1,2$, $$h: \mathcal L_1\to \mathcal L_2,$$
  where these $h$ must {\it almost  commute}, i.e.  {\it commute modulo compact s},  with the parallel transports in 
 in  $ \mathcal L_1$ and $\mathcal L_2$ along smooth paths in $X$. 
 
 (This  idea for flat bundles goes back to { [Atiyah(global) 1969],  [Kasparov(index) 1973], [Kasparov(elliptic)   1975],  
 [Mishchenko(infinite-dimensional) 1974]  and  where  non-flat generalizations and applications  are discussed in \S$9\frac {1}{6}$  of [G(positive) 1996].)
  
  (Such an $\mathbf h$  
 represents the finite dimensional virtual (not quite)  bundle  $ker(h)-coker (h)$.)
  
  Define
  
   $$|curv| (h) =\max (|curv|\|(\mathcal L_1),|curv|(\mathcal L_2))$$ 
and let 
 $$|curv|(\mathbf h)=\inf  |curv| (h)$$
 where the infimum is taken over all $h$ in the class $\mathbf h$.\vspace{1mm}
 
 {\it Why Hilbert?} If one  limits the  choice of representatives of $\mathbf h$  to virtual  {\it finite dimensional}  bundles $L\to X$, then the resulting curvature function on $K^0(X)$ may only increase:  
 $$|curv|(\mathbf h)_{fin.dim}\geq |curv|(\mathbf h).$$
Apparently, this must be standard,  the  Hilbert spaces in  the definition of   Fredholm bundles can be approximated by finite dimensional Euclidean ones,
\footnote {This is an exercise that the author delegates to the reader.} that implies that
  $$|curv|(\mathbf h)_{fin.dim}= |curv|(\mathbf h),$$ 
 but even so "Hilbert"   allows  greater flexibility of certain constructions, example of which   we shall see below. 
 
 \vspace{1mm}

{\color {blue} \it Naive (Strong Novikov)
 {\color {red!50!black} Conjecture}.}  Let $Y$ be a compact {\it  aspherical}\footnote {The universal covering of $X$ is contractible.}  Riemannian manifold, possibly with a boundary.
Then {

{\sl all} (classes of complex  vector bundles)   {\sl $\mathbf h\in K^0(Y)$ satisfy:
$$\inf_N|curv|(N\cdot\mathbf h)=0,\mbox { } N=1,2,3,...,\mbox { }  .$$}

{ \it Exercises.} (a) Show that the equalities  $|curv|(\mathbf h)=0$ and $\inf_N|curv|(N\cdot\mathbf h)=0$ are {\it homotopy invariants} of $Y$.

(b) Show that if $Y$  satisfies this naive conjecture and $X$ is a closed Riemannian orientable spin $n$-manifold with $Sc(X)>0$, then all continuous maps $f:X\to Y$ send the fundamental rational homology class $[X]_\mathbb Q \in H_n(X, \mathbb Q) $ to zero in $H_n(Y, \mathbb Q)$.

 %%%%%%%%%%%%%%%%%%%%%%
 \subsubsection {\color {blue} Area, Curvature  and $K$-Cowaist}\label {cowaist4}

 %%%%%%%%%%%%%%%%%%%%%
 \vspace{1mm}
 
 {\it $K$-cowaist$_2$.} Given a Riemannian manifold $Y$ (or a more general space, e.g. a polyhedral one  with a piecewise smooth metric), define the {\it $K$-cowaist }  on the  homology classes $h_\ast  \in H_\ast(Y)$, denoted  
K-$cowaist_2(h)$ \footnote{Subindex 2  is to remind that curvature of  bundle $L$ over $Y$  is seen on restrictions of $L$ to surfaces in $Y$.} as the infimum of $|curv|(\mathbf h)$ over all $\mathbf h\in K^0(Y)$, such that $\mathbf h(h_\ast)\neq 0$, where this equality serves as an  abbreviation for the value of the Chern character of $\mathbf h$ on $h_\ast$,
$$\mathbf h(h_\ast)=_{def} Ch(\mathbf h)(h_\ast).$$

In these terms the above {\Large \color {blue} $\bullet$} can be reformulated as follows.\vspace{1mm}

 {   \color {blue} \it $K$-cowaist Inequality for Closed Manifolds.} {\sl The K-cowaists of  (the fundamental classes of) closed orientable $2m$-dimensional spin manifolds $X$ with $Sc(X)\geq \sigma>0$ satisfy:
$$\mbox {K-}cowaist_2[X] \leq \frac {const_m}{\sigma}.\leqno {\mbox {\Large \color {blue} $\bullet_{wst}$}}$$}
 
Notice, that {\it \large conjecturally,} a similar  inequality also holds for the {\sl ordinary  2-waist}, (see  [Guth(waist) 2014]  for an exposition of  this "waist")  where it is confirmed  for 3-manifold by   the Marques-Neves theorem (see  section \ref {slicing3D.3})  
 
\vspace{1mm} 

{\it Exercises.} Show that the K-cowaist is bounded by the hyperspherical radius defined in section \ref {filling+hyperspherical+asphericity3} as follows, 
$$ \mbox{K-}cowaist_2[X]\leq 4\pi Rad_{S^{2m}}^2(X)$$

(b)  Show that K-$cowaist_2(S^n)=4\pi$.\vspace {1mm}

{\it Almost Flat Bundles Over Open Manifolds.}  If $X$ is a non-compact manifold, then we deal with the  K-theory 
with compact support that is represented by Fredholm homomorphisms 
$$h:\mathcal L_1\to  \mathcal L_2$$
which are {\it  isometric and connection preserving isomorphisms at infinity}, i.e. away from  compact subsets in $X$ where the corresponding $K$-group is denoted $K^0(X/\infty)$.
 (If $X$ is compact then $K^0(X/\infty)=K^0(X)$

 Here the Hilbertian nature of  "Fredholm" allows a painless  (and obvious by deciphering terminology)   definition  of the {\it pushforward homomorphism} for possibly {\it infinitely} sheeted  covering maps $F:X_1\to X_2$,  
 $$F_\star: K^0(X_1/\infty) \to K^0(X_2/\infty),   $$
where, clearly,
$$ |curv|(F_\star(\mathbf h)) \leq  |curv|(\mathbf h)$$
for all $\mathbf h\in K^0(X_1/\infty).$\vspace{1mm}

It follows that 
 {\sl $$K\mbox {-}cowaist_2[X_1]\leq K\mbox {-}cowaist_2[X_2]$$
 for  coverings $X_1\to X_2$ between orientable Riemannian manifolds.}

\vspace {1mm}

{\it On K-cowast Contravariance.} The {\it compact support} property of (virtual) bundles $L\to X_2$  is preserved under pullbacks by  {\it proper} maps $F:X_1\to X_2$, e.g. by finite coverings, but it fails, for instance, for {\it infinitely sheeted} coverings  
$F:X_1\to X_2$.

This makes the  inequality 
 $$K\mbox {-}cowaist_2[X_1]\geq  K\mbox {-}cowaist_2[X_2]$$
  (that is  obvious  for  { \it finitely sheeted} coverings)  
   {\it problematic}  for infinite covering maps  $F:X_1\to X_2$.

This should be compared with the  {\it covariance problem} for  {\it max-scalar curvature} which is defined in section \ref {max-scalar5}  and  which obviously lifts under covering maps,
 $$Sc^{\sf max}_{prop}[X_1]\geq Sc^{\sf max}_{prop}[X_2],$$
while the   opposite inequality  causes a problem (see section \ref {max-scalar5}).

{\it \large Question.} Can one match the covariance of $Sc^{\sf max}$ by a somehow generalized  K-$cowaist_2$ that would be invariant under (finite and infinite) covering maps $F:X_1\to X_2$?

Specifically, one looks  for {\it almost flat} (virtual)   {\it infinite dimensional} Hilbert bundles in  a suitable $K$-theory, which would be  
compatible with the index theory and with the Schroedinger-Lichnerowicz-Weitzenboeck formula 
in the spirit of Roe's $C^\ast$-algebras. 
c\vspace {1mm}
 
 {\it Amenable  Cutoff Subquestion.} Let $X_2$ be a closed orientable Riemannian manifold of dimension $n=2k$  and let $L\to X_2$ be a vector bundle induced by an $\varepsilon$-Lipschitz map $f: X_2\to S^n$
 from the positive spinor bundle $L=\mathbf S^+=\mathbf S^+(S^n)\to S^n$. 
c
 Suppose that the fundamental group $\pi_1(X_2)$ is {\it amenable}, let 
 $X_1=\tilde X_2\to X_2$ be the universal covering map  and let 
 $$\tilde  L =F^\ast (L)\to  X_1$$ 
be the pullback of $L$. 
 
 When do there exist  unitary  bundles $\tilde L_i\to X_1$, $i=1,2,...$, with unitary connections, such that
 
 $\bullet_\infty$ the bundles  $\tilde L_i$ are flat  trivial at infinity;  
 
  $\bullet_{|\tilde L}$ there is an exhaustion of $X_1$ by compact {\it F\o lner  subsets} 
  $$V_1\subset ... \subset V_i\subset ... \subset X_1,$$
   such that the restrictions of  $\tilde L_i$ to  $V_i$ are equal to the restrictions of $\tilde L$,
 $$( \tilde L_i)_{|V_i}=\tilde L_{|V_i};$$
 
 $\bullet_{\int}$ the integrals of the $k$-th powers of the  curvatures of $L_i$ are dominated by such integrals for $\tilde L$ over $V_i$, 
 $$\frac {\int_{X_1} |curv|^k(\tilde  L_I) dx_1} {\int_{V_i} |curv|^k  (\tilde L) dx_1} \underset {i\to \infty}\to 0;$$
 
 $\bullet_\epsilon$ the curvatures of all $\tilde L_i$ are bounded by 
  $$|curv|(\tilde L_i) \leq \epsilon, $$ 
where  $\epsilon=\epsilon_n(\varepsilon)\to 0 $ for $\varepsilon\to 0$.
 
 (The  Federer Fleming isoperimetric/filling inequality in the rendition of [MW(mapping classes) 2018]
may be useful here.)
 
 \vspace {1mm}
 
{\it  Non-Amenable Cutoff  Example.} Let $X(=X_2)$ be a closed orientable  Riemann surface of genus $\geq 0$ and 
  $L\to X$  a complex line bundle with a  unitary connection, e.g. $L$ is   the tangent bundle  $T(X)$, the Chern number of  which   $c_1(T(X))[X]=\chi(X)$ doesn't vanish  for $genus(X)>0$.  
  
 Let $\tilde L\to \tilde X$ be the lift (pullback)  of $L$ to the universal covering $\tilde X(=X_1)$ of $X$ and observe that there exit disks $\tilde D^2(R)\subset \tilde X$,  such that the parallel translates over the boundary circles $\tilde S^1(R)=\partial \tilde D^2(R)$ are  a multiples of $2\pi$ and where the radii $R$ of such disks can be arbitrary large.  \vspace {1mm}
 
{\sf Then the restriction of $\tilde L\to \tilde X$ to such a disk 
$\tilde D^2(R)\subset \tilde X$ {\it extends to a bundle},  call it $\tilde L_R\to \tilde X$, which  is  {\it trivial outside $\tilde D^2(R) $}  and such that 
$$c_1(\tilde L_R/\tilde S^1(R))\sim area(\tilde D^2(R))\underset{R\to \infty}  \to \infty,$$
provided the curvature of $L$ (that is a closed  2-form on $X$)  {\it doesn't vanish.}}

   \vspace {1mm}
  
  {\it Problem for  $n>2$.} The main difficulty in   similarly trivializing at infinity bundles over $n$-dimensional Riemannian manifolds $X$ for  $n=dim(X)\geq 3$
   seems to be associated with the following {\large \it \color {blue}questions.} 
  
  {\sf  Let $\mathcal U_b(k)=\mathcal U_b(k, X)$, $b\geq 0$, be the space    of the unitary connections $\nabla$ on a trivial bundle $L\to X$ of rank $k$, such that $|curv|(\nabla)\leq b$. }

{\color {blue} (a)}   {\sl For which values $b_1$ and $b_2>b_1$  are  the  connections from $ \mathcal U_{b_1}(k)$  homotopic in 
  $\mathcal U_{b_2}(k)\supset \mathcal U_{b_1}(k)$?

  {\color {blue} (b)}   When do  
  the  homomorphisms of the   homotopy groups 
  $$\pi_i(\mathcal  U_{b_1}(k))\to \pi_i(\mathcal  U_{b_2}(k)), \mbox { } i\geq 1, $$
   induced by   the inclusions
  $\mathcal U_{b_1}(k) \hookrightarrow \mathcal U_{b_2}(k)$ vanish?
  
 {\color {blue} (c)}  How do the Whitney sum homomorphisms
 $$ \mathcal U_b(k_1)\times  \mathcal U_b(k_2)\to  \mathcal U_b(k_1+k_2)$$
 behave in this respect?
  
   In particular, what happens to the  homomorphisms
  $\pi_i(\mathcal  U_{b_1}(k))\to \pi_i(\mathcal  U_{b_2}(k))$
  under stabilization  
  $$\underset {N}{\underbrace {\mathcal U_b(k)\times...\times \mathcal U_b(k)}} \leadsto   (\mathcal U_b(Nk))$$ 
   for $N\to \infty$?}

  \vspace {1mm}

 {\it Exercise.} Let $X$ be a complete orientable  even dimensional Riemannian  manifold  with nonpositive  sectional curvature.  
Show that there exists a $K$-class $\mathbf h\in K^0(X/\infty)$, such that  
 $$|curv|(\mathbf h)=0 \mbox {  and } \mathbf h[X]\neq 0,$$
 where $[X]$ denotes the fundamental homology class of $X$ with {\it infinite supports}.
 
%%%%%%%%%%%%%%%%%%%%
\subsubsection {\color {blue} Sharp Algebraic  Inequalities for the  $L$-Curvature in the  Twisted SLW(B) Formula }\label {sharp algebraic4}
%%%%%%%%%%%%%%%%

{\it \color {blue}  Normalization of Curvature.} In so far as the scalar curvature is concerned we are interested not in the curvature $|curv|(L) $ per se but rather in the norm of the endomorphism()
$$\mathcal R_{\otimes L}:\mathbb S\otimes L\to \mathbb S\otimes L $$
 in the  Schroedinger-Lichnerowicz-Weitzenboeck formula for the twisted Dirac  operator, 
 $$ \mbox { ${\cal D}_{\otimes L}^2=\nabla^2_{\otimes L} + \frac{1}{4}Sc(X) + {\cal R}_{\otimes L}$},$$
 (see the  previous section)
where this $\mathcal R_{\otimes L}$ is as following linear/tensorial  combination of the values of the curvature of $L$ on the  tangent bivectors in the manifold $X$, 
(see [GL(spin) 1980],[Lawson\&Michelsohn(spin geometry) 1989] and section \ref {Clifford3})
 $$\mathcal R_{\otimes L}(s\otimes l)=\frac {1}{2}\sum_{i,j}(e_i\circ e_j\circ s) \otimes R^L_{e_i\wedge e_j}(l),$$
where

{\sf $e_i\in T_x(X)$, $i=1,...n=dim(X)$ is an orthonormal frame of tangent vectors at a point $x\in X$,
 
 $s\in \mathbb S$, are spinors, 

$l\in L$ vectors in the bundle $L$, 

$R^L(e_i\wedge e_j):L\to L$ is  the  curvature of $L$ (written down  as  the  valued $2$-form on $X$)

and 

"$\circ$" denotes the Clifford multiplication.}

This suggest the definition  of
$$\lambda_{\min}[curv]_{\otimes \mathbb S}(L)$$ as the   smallest (usually negative)  eigenvalue of the  $||\mathcal R_{\otimes L}||. $\vspace {1mm} 

{\Huge   \color {blue}  $\ast$ }{\large  \color {blue} \sf Example: Llarull's algebraic  inequality}. [Llarull(sharp estimates)  1998] 
{\sl Let $f:X\to S^n$ be a smooth 1-Lipschitz, or more generally, an area non-increasing map and let $L\to X$  
be the pullback  the  spinor bundle $\mathbb S(S^n)$. Then 
 this minimal eigenvalue  of the   $\mathcal R_{\otimes L}$ satisfies:
 $$\lambda_{\min}[curv]_{\otimes \mathbb S}(L)=-\frac {1}{4}(n(n-1)= -\frac {1}{4}Sc(S^n).$$}
 (We return to this   in 
 %%%%%%%%%%%%%%%
 {\color {red}corrected}%% ???
 the  next section,)  \vspace {1mm}

 Using this  $\lambda_{min}[curv]$ instead of the  $|curv|$ one defines 
 $$\lambda_{min} [curv]_{\otimes \mathbb S}(\mathbf h), \mbox {  $\mathbf h \in K^0(X)$,} $$ 
 as the {\it supremum} of $\lambda_{min} [curv]_{\otimes \mathbb S(L)}$ for all (virtual) bundles $L$
in  the class of $\mathbf h$,

 Accordingly  one  modifies the above   $K$-$cowaist_2(\mathbf h)$ and  define the corresponding {\it K-cowaist coupled with spinors}, denoted 
  K-$cowaist_{\otimes \mathbb S, 2} (h_\ast)$, $h_\ast\in H_\ast(X)$, as the supremum of $\lambda_{min} [curv]_{\otimes \mathbb S}(\mathbf h)$ over 
  over all $\mathbf h\in K^0(Y)$, such that $\mathbf h(h_\ast)\neq 0$.

Then, for instance, the above  {\Large \color {blue} $\bullet_{wst}$} for spin manifolds $X$ takes more elegant form:

$$\mbox {K-}waist_{\otimes \mathbb S, 2} [X] \leq \frac {4}{\sigma}  \mbox { for }  \sigma =\inf_x Sc(X,x)>0.$$

Notice that  this inequality, combined with the above {\Huge   \color {blue} $\ast$},  implies  Llarull's {\it geometric inequality} $Rad_{S^n}(X)\leq \sqrt\frac {n(n-1)}{\sigma}$, 
which we discuss at length in 
the next section.

  Also this may give better formulae for {K-}cowaists of product of manifolds.

(See section \ref {max-scalar5} and also [G(positive)  1996] and [G(101)  2017]
for   other known  and conjectural   properties  of  $|curv|(\mathbf h)$ formulated in these papers  in the language of the K-area.)

%%%%%%%%%%%%%%%%

\subsection {\color {blue} Llarull's and Goette-Semmelmann's {\sf Sc}-Normalised Estimates for Maps to  Convex Hypersurfaces in Symmetric Spaces. } \label {Llarull4}

%%%%%%%%%%%%%%%%%%%%%%%

Let us now  look  closer at  the above
$${\cal R}_{\otimes L}(s \otimes l)= \frac{1}{2}\sum_{i,j}(e_i\circ e_j\circ s)\otimes  R_{ij}(l), $$
that is the 
 endomorphism of ( on) the bundle $\mathbb S\otimes L\to X$,

which appears in  the zero order term in the twisted Dirac   
 $$  {\cal D}_{\otimes L}^2=\nabla^2_{\otimes L} + \frac{1}{4}Sc(X) + {\cal R}_{\otimes L},$$
for    
$${\cal D}_{\otimes L }:C^\infty(\mathbb S\otimes L )\to C^\infty(\mathbb S\otimes L).$$

{\it Example  of $ L=\mathbb S$ on $S^n$.} Since the norm of the curvature  of (the Levi-Civita connection on) 
the tangent  bundle is one,  the norm of the curvature  operators $R_{ij}:\mathbb S \to \mathbb S $ are  at most ( in fact, are to)  $\frac{1}{2}$,
$$||R_{ij}(s)||\leq \frac{1}{2},$$
since the spin bundle $\mathbb S(X)$ serves as the "square root" of the tangent bundle $T(X)$, where this  is literally true for  $n=dim(X)=2$, that formally implies the inequality $||R_{ij}(s)||\leq \frac{1}{2}$ for all $n\geq 2$.

And since the Clifford multiplication operators  $s\mapsto e_i\cdot e_j\cdot s$ are unitary, 
$$||{\cal R}_{\otimes L}(s \otimes l)||\leq \frac {1}{4}n(n-1) =  \frac {1}{4}Sc(S^n)$$

This doesn't, a priori,  imply this inequality for all (non-pure)  vectors $v$ on the tensor product
 $\mathbb S\otimes L$ for $L=\mathbb S$, but,  by diagonalising  the Clifford multiplication operators in a suitable basis
and by employing the {\it essential constancy}\footnote{Some eigenvalues of this  are $\pm1$
 and some zero.}  of the curvature $R_{ij}$ of $S^n$, [Llarull(sharp estimates)  1998]  shows that
$$||\langle {\cal R}_{\otimes L}(\underline \theta ), \underline \theta\rangle ||\geq - \frac {1}{4}n(n-1)$$
for all unit vectors $\underline \theta\in \mathbb S(S^n)\otimes  \mathbb S(S^n)$. 

 This inequality for twisted spinors on $S^n$ trivially yields the corresponding inequality on all 
 manifolds $X$ mapped to $S^n$, where the bundle  $L\to X$ is the induced from the spin 
 bundle $\mathbb S(S^n)$.  
 
 Namely, let $X=(X,g)$ be an $n$-dimensional Riemannian manifold, 
 $f: X\to S^n$ be a smooth map,  $L=f^\ast(\mathbb S(S^n))$, let $df:T(X)\to T(S^n)$ be the differential of $f$ and 
 $$\wedge^2df: \wedge^2 T(X)\to \wedge^2 T(S^n)$$  be the exterior square of $df$.\footnote {Recall that the norm $||\wedge^2df||$ measures by how $f$ contracts/expands surfaces in $X$. For instance the inequality  $||\wedge^2df||\<1$ signifies that $f$ decreases the areas of the surfaces in $X$.}

Then the  $${\cal R}_{\otimes L}: \mathbb S(X)\otimes L\to  \mathbb S(X)\otimes L$$
satisfies
$$||\langle {\cal R}_{\otimes L}( \theta ),  \theta\rangle ||\geq -||\wedge^2df||\frac {n(n-1)}{4}, \mbox  { } L=f^\ast(\mathbb S(S^n)),$$
 for all unit vectors $\theta \in \mathbb S(X)\otimes f^\ast(\mathbb S(S^n)).$

Moreover,  -- this is formula  (4.6) in [Llarull(sharp estimates)  1998] -- 

$$||\langle {\cal R}_{\otimes L}( \theta ),  \theta\rangle ||\geq -\frac {1}{4}|trace \wedge^2df|,$$
where $trace \wedge^2df$ at a point $x \in  X$
stands for
$$\sum_{i\neq j} \lambda_i\lambda_j,$$
for  the differential $df :T_x(X) \to  T_{f(x)}(S^n)$   diagonalised to the orthogonal sum of
multiplications by $\lambda_i.$

This inequality, restricted to $L^+= f^\ast(\mathbb S^+(S^n))$ together with the index formula, which says for this 
$L_+$ that 
$$ind( {\cal D}_{\otimes L^+}) = \frac{|deg(f)|}{2} \chi(S^n),$$
provided $X$ is a closed  oriented spin manifold.

Thus we arrive  at the proof  of  Llarull's theorem in the Sc-normalized trace form    suggested by Mario Listing in  [Listing(symmetric  spaces) 2010]. 
\vspace{1mm}  

{\huge \color {blue}$\star$}  {\large \it \color {blue}  $trace \wedge^2df$-Extremality of $S^n$.}\footnote {This is {\it Spherical Trace Area Extremality Theorem} from  section \ref{area extremality3}.} {\sf Let $X$ be a closed orientable  Riemannian spin $n$-manifold  and $f:X\to S^n$ a smooth map of nonzero degree.}

{\it If 
$$ Sc(X,x) \geq \frac {1}{4}|trace \wedge^2df(x)|  $$
at all points $x\in X$n then, in fact,  $ Sc(X) = \frac {1}{4}|trace \wedge^2df|  $ everywhere on $X$.}

\vspace{1mm}

{\it In fact,} if $n$ is even and $\chi(S^n)=2\neq 0$, this follows from the above. And if $n$ is odd, there are  (at lest) three different 
reductions  to the even dimensional case (see [Llarull(sharp estimates)  1998], [Listing(symmetric  spaces) 2010], [G(inequalities) 2018]),  but these are artificial and conceptually unsatisfactory.
  
Also see see [Llarull(sharp estimates)  1998]  and  [Listing(symmetric  spaces) 2010] for   characterisation  of maps  $f$, where $ Sc(X) = \frac {1}{4}|trace \wedge^2df|  $. 
 
\vspace {1mm} 

Llarull's estimate for the bottom of the spectrum of the  curvature operator in spin  bundle $\mathbb S(S^n)$, was generalized by  Goette and Semmelmann    [Goette-Semmelmann(symmetric) 2002]
to the other  Riemannian manifolds $\underline X$ with {\it non-negative curvature operators}, 
and (in  the  Sc-normalized form suggested  Listing)  resulted in the following.

\vspace {1mm}

{\huge \color {blue}$\star \star$}  {\large \it \color {blue}  $ \wedge^2df$-Extremality Theorem.}\footnote {This generalizes {\it Spin-Area Convex Extremality Theorem} from section \ref{area extremality3}.} {\sf Let $\underline X=(\underline X,\underline g)$  and $X=(X.g)$  be a closed orientable Riemannian {\it spin} $n$-manifolds.  where $\underline X$ has  {\it non-negative curvature operator}  and let $f:X\to \underline X$ be a smooth map of {\it non-zero} degree.

  {\it If $\underline X$ has {\it non-zero} Euler characteristics, then
 this map can't be strictly area decreasing with respect to the $Sc$-normalised metrics
  $g^\circ =Sc(g)\cdot g$ and 
  $\underline g ^\circ=Sc(\underline g)\cdot \underline g$. }
  
  This means that 

{\sf  if  of the exterior square of the differential of $f$ with respect to the original
 metrics  $g$ and $\underline g$ is related to the scalar curvatures of the two manifolds by 
 the  {\it inequality} 
$$Sc(g,x)\geq ||\wedge ^2df(x)|| Sc(\underline g, f(x))$$
 at all $x\in X$, where $||\wedge ^2df(x)||$ stands fo the sup-norm with respect to the metrics $g$ and $\underline g$,}

\hspace {30mm} {\it then the   equality holds}:
$$Sc(X,x)= ||\wedge ^2df|| Sc(\underline g, f(x)).$$}

 {\it Examples, Remarks, Conjectures.} (a) All {\it compact symmetric} spaces  have 
 non-negative curvature operators.
 
 Also 
 
(a$_1$)  {\it the  induced metrics in  convex hypersurfaces in these spaces 
 also have the  curvature  operator non-negative},\footnote {This was explained to me by Anton Petrunin,  who  introduced  a    class of   metrics inherited by convex hypersurfaces, see [Petrunin(convex) 2003].}
 
 and

 (a$_2$) {\it Riemannian products of manifolds with $curv.oper\geq 0$ have  $curv.oper\geq 0$.}

  \vspace {1mm}

(a$_3$)  By a theorem of Alan Weinstein  [Weinstein(Positively curved) 1970],

{\it submanifolds $\underline X^n\subset  \mathbb R^{n+2}$
with non-negative  sectional curvatures have non-negative curvature  operators.}

(b) Llarull, Goette-Semmelmann, and Listing also analyzed the equality cases in their  papers  and 
  proved the corresponding  {\it rigidity theorems.}

(c) Goette  and Semmelmann also state  in their paper  an extremality/rigidity result for  odd dimensional $\underline X$, that was  scrutinized  and generalized  in [Goette(alternating torsion)2007].

(d) Besides  symmetric spaces, Goette  and Semmelmann proved  $ \wedge^2df$-extremality 
for  was proven  for K\"ahler manifolds 
with positive Ricci curvature.\footnote {See   [Goette-Semmelmann(Hermitian) 1999] and the earlier "symmetric" paper  [Min-Oo(Hermitian) 1998].}
\vspace {1mm}

(e) {\color {red!40!black}Conjecturally,} neither  {\it spin} nor $\chi(\underline X)\neq 0$-condition are 
 necessary for the $ \wedge^2df$-extremality.

In fact, Goette  and Semmelmann (as well as Min-Oo) prove their theorems  not only for spin manifolds but also for certain {\it spin$^c$-manifolds} and also for  for  { \it spin maps}
$f:X\to \underline X$  between non-spin manifolds, i.e. where $f$ pulls back the 
 Stiefel-Whitney class $w_2(\underline X)$ to $w_2(X)$.

(f) The above extremality theorems were generalized in the original papers to maps  
$f:X\to \underline X$, where $dim(X)=n=\underline n+4k$,  $\underline n= dim( \underline X)$,  and where $f$ has non-zero
$\hat A$-degree, i.e.  where the pullback $f^{-1}(\underline x) \subset X$ of a generic point  
$\underline x \in \underline X$ has {\it non-zero $\hat A$-genus.}

 \vspace {1mm}
 
{\it $\mathbb T^\rtimes$-Stabilization  of Extremality Theorems and   Generalizations}.
The above (f) suggests the following.c\vspace{1mm}
 
{\color {red!50!black} \textbf {Conjecture}}. {\sf  If  a  Riemannin manifold  $\underline X$ is $ \wedge^2df$-extremal, then, for all $X$ and all smooth maps $f:X\to \underline X$, such that 
$$Sc(g,x)> ||\wedge ^2df(x)|| Sc(\underline g, f(x)),$$
the  generic  pullback $f^{-1}(\underline x) \subset X$  is homologous (even bordant) in $X$ to a submanifold $Y$, which {\it supports a metric with positive scalar curvature}}.\vspace{1mm}
 
As it stands, this seems not very realistic.

 However, if the extremality  of $\underline X$ follows by the above kind of argument relying on a {\it sharp}   SLW(B)-{\it inequality} for the Dirac operator on $X$
twisted with the pullback $L^\ast=f^\ast(\underline L) $ of some  bundle  $\underline L\to\underline X$, with a unitary connection, then, as we shall explain below, \vspace {1mm}

{\huge \color {blue}$\star \star \star$}   {\it  the inequality  $Sc(g,x)> ||\wedge ^2df(x)|| Sc(\underline g, f(x))$ implies vanishing not only of 
$\hat A(f^{-1}(\underline x))$ but of more general (all?)  Dirac theoretic obstructions for $Sc>0$ on 
$(n-\underline n)$-dimensional manifolds.}\footnote {
The simplest instance of this, where $\underline X= S^{\underline n}$ and where $X$ is a warped extension $X_0\rtimes, \mathbb T^1$, was observed   in \S$5\frac {4}{9}$ in [G(positive) 1996], 
and  used for the proof of  a special case of  $C^0$-{\it closure theorem} from section 
 \ref {C0-limits3}.} \vspace{1mm}
 
 The basic (and fairly general) instance of this is where $X$ supports $\varepsilon$-flat bundles $L_\varepsilon \to X$ for all $\varepsilon>0$   (i.e.  $L_\varepsilon$ are endowed with  unitary connections the  curvatures of which are bounded in norm by
$\varepsilon$), such that that the indices of the Dirac operator $\cal D$ on $X$ twisted with $L^\ast\otimes L_\varepsilon$, as expressed by the index formula,  don't vanish for $\varepsilon \to 0$.

Since the norm of the connection curvature term in SLW(B)-formula for the operator 
$\mathcal D_{\otimes (L^\ast\otimes L_\varepsilon)}$ converges,  for $\varepsilon\to 0$, to that 
for $\mathcal D_{\otimes L^\ast}$, the inequality
 $$Sc(g,x)- ||\wedge ^2df(x)|| Sc(\underline g, f(x))\geq \delta>0$$
implies vanishing of of the index of $\mathcal D_{\otimes (L^\ast\otimes L_\varepsilon)}$
 for  $\varepsilon << \delta$ and the proof follows. \vspace{1mm}

{\it Remarks and Examples.} (a) The bundles $L_\varepsilon$ can be understood in a
 fairly general way, e.g. is  virtual Fredholm bundles, as families of such bundles or, more 
 generally as moduli over the (reduced)  $C^\ast$-algebra of a quotient group of the fundamental group of $X$. 
 
 (b) If $\underline X=   \underline X_0 \times \underline X_1 $, where $X_1$ is a  compact orientable Riemannin spin manifold with $curv.oper(X_0)\geq 0$ as  in  {\huge \color {blue}$\star \star$}, if $X$ is orientable spin, and if $f:X\to \underline X$  
 is a map of non-zero degree, such that $Sc(g,x)> ||\wedge ^2df(x)|| Sc(\underline g, f(x)),$
 then,  {\color {red!50!black} probably,}  {\it the rational Rosenberg index}(see  [Zeidler(width) 2020])
  $$\alpha (\underline X_1) \in (KO_{\underline n_1}(C^\ast \pi_1(\underline X_1)))\otimes \mathbb Q$$
 vanishes.}
 \vspace{1mm}
 
I feel shaky in these matters (this must be obvious to the readers  well versed in the  K-theory of $C^\ast$-algebras) but the proof of this  is transparent in many cases.  \vspace{1mm}

For instance, this is so   

{\it if the universal covering  $\tilde {\underline X_1}$ of  $\underline X_1$ is $\wedge^2$-hyper-Euclidean} i.e. there exists a smooth  proper area non-decreasing map $\tilde {\underline X}_1\to\mathbb R^{n_1}$, $n_1=dim(X_1)   $,

In fact,  the above considerations and the relative index theorem yield the following more general 
proposition. 
\vspace {1mm}

{\color {green!50!black}$\bigstar$}  \textbf {Non-compact Extremality Theorem.} {\sf  Let $X$ and
 $\underline X_0$  be connected orientable Riemannin {\it spin} manifolds of dimensions $n$ and $\underline n_0$, where  $X$ is {\it complete}  and $\underline X_0$ is {\it compact}, and let {\it the curvature operator  of  $\underline X_0$ be non-negative}. 

  Let
$$f =(f_0, f_1): X\to \underline X_0\times \mathbb R^{ m}, \mbox { }  m=n- n_0,$$ 
be a  smooth proper map  with non-zero degree. }

{\it Let    $\underline X_0$  be simply connected, let     the Euler characteristics of $\underline X_0$ don't   vanish, $\chi(\underline X_0)\neq 0$, and let the map 
$f_1: X\to  \mathbb R^{ m}$ be  area non-increasing.}
 Then  
$$\inf_{x\in X}\left( Sc(g,x)- ||\wedge ^2df(x)|| \cdot Sc(\underline g, f(x)\right)\leq 0.$$}}\vspace {1mm}

{\it  Proof.}
 Here, the  relevant bundle 
 $L^\ast\to X$ is  the $f_0$-pull back of the positive spin bundle $\mathbb S^+(  \underline 
 X_0)$  (as in {\huge \color {blue}$ \star$} and in {\huge \color 
 {blue}$\star \star$})  from $ \underline X_0$ to $X$, while the bundles  $L_\varepsilon\to X$ are the 
$f_1$ pullbacks of the complex   {\it $\varepsilon$- flat bundles of ranks $l=\frac {m}{2}$} (if $m$ is odd, multiply    $X$ and $ \underline X_0$  by $\mathbb R^1$)
 on  $\mathbb R^m$ with compact supports    (e.g. flat split at   infinity) and such that 
 {\it they have their relative 
Chern numbers  $c_l\to \infty  $  for $\varepsilon \to 0$.}

This and non-vanishing of $\chi( \underline X_0$ imply (half line computation)  the non-vanishing of the relative index of 
 $\mathcal D_{\otimes (L^\ast\otimes L_\varepsilon)}$ by the relative index theorem.
%The computation of the relative  indices of $\mathcal D_{\otimes (L^\ast\otimes L_\varepsilon})$ 
%and their non-vanishing is straightforward in this case
  and the proof concluded is with  the ($\varepsilon$-perturbed) SLW(B)-formula  for  $X_0$ as in  the Goette-Semmelmann theorem {\huge \color {blue}$ \star \star$}.

\vspace{1mm}

{\it Remarks/Corollaries} (a) {\color {red!40!black} Probbaly},  it is  is not hard to prove rigidity of 
$\underline  X _0$ in this case.

(b)  Instead  of "{\it simply connected and } $\chi(\underline X_0)\neq 0$}"   one could requirer that {\it the universal covering $\tilde {\underline X_0}$ has non-zero Euler characteristics}.

 Indeed,  by the Gromoll-Meyer  theorem,}  $\tilde {\underline X_0}$
isometrically splits, 
$$\underline X_0=\underline X_0' \times \mathbb R^k,$$
where $\underline X_0'$ is compact simply connected 
and the theorem apples to  $\underline X_0'\times \mathbb R^{k+m}$.

(c) The above proof, similarly to these of {\huge \color {blue}$\ \star$} and {\huge \color {blue}$\star \star$},  easily  generalizes to maps  $f$ with  {\it non-vanishing  $\hat A$-degrees.}

\vspace {1mm}

{\it Question.} {\sf Can one approach  the above conjecture from the  opposite angle by actually {\it  constructing 
$(n-m)$-submanifolds in $X$ with positive scalar curvatures} in the homology class of 
$f^{-1}(\underline x)$?}

(Application of $\mu$-bubbles, as we know,  allows  such constructions, but these {\it fail to deliver   sharp inequalities} of this kind).
 
%%%%%%%%%%%%%%%%%%%%%%%%%%%%%%

 \subsection {\color {blue}Bounds on Mean Convex Hypersurfaces}\label {mean4}

%%%%%%%%%%%%%%%%%%%%%%%%%%%

Recall that the {\it spherical radius} $Rad_{S^{n-1}}(Y)$ of a  connected  orientable Riemannian manifold of dimension $(n-1)$ is the supremum of the radii $R$ of the spheres $S^{n-1}(R)$, such that $X$ admits a   distance decreasing  map $f:Y\to S^{n-1}(R)$ of non-zero degree, where this $f$ for non-compact $Y$ this map is supposed to be constant at infinity.\footnote{Alternatively, one might require $f$ to be {\it locally} constant at infinity, or more generally, to have the limit set  of codimension$\geq 2$ in $S^{n-1}(R)$. }

We already indicated in  section \ref {mean convex3} also see [G(boundary) 2019] that  Goette-Semmlenann's theorem (above {\huge \color {blue}$\star \star$}), applied to smoothed doubles  \DD$X$ and  \DD$\underline X$ yields the following  corollary.

\vspace{1mm}

{\Large \color {blue}$\Circle^{n-1}$}{\it Let  $X$ be  a compact  orientable Riemannian manifold with boundary  $Y=\partial X$.}

{\it If $Sc(X)\geq 0$ and the mean curvature of $Y$ is bounded from below by $mean.curv(Y)\geq \mu>0$,  then the hyperspherical radius of $Y$  for the induced Riemannian metric is bounded  by
$$ Rad_{S^{n-1}}(Y)\leq \frac {n-1}{\mu}.$$}

In fact, the proof of this indicated in section \ref {mean convex3} (also see [G(boundary) 2019]) together with the above {\color {green!50!black}$\bigstar$} yields the following 
more general theorem. \vspace{1mm}

{\color {green!50!black}$\bigstar_{mean}$} \textbf {Non-Compact Mean Curvature Inequality}. {\sf  Let $X$ and
 $\underline X_0$  be connected orientable Riemannin {\it spin} manifolds of dimensions $n$ and $\underline n_0$ with boundaries, where  $X$ is {\it complete}  and $\underline X_0$ is {\it compact}, and let {\it the curvature operator  of  $\underline X_0$ be non-negative}. 

  Let
$$f =(f_0, f_1): X\to \underline X_0\times \mathbb R^{ m}, \mbox { }  m=n-\underline n_0,$$ 
be a  smooth proper map, which sends $\partial X\to \partial  \underline X_0\times \mathbb R^{ m}$ and which has  {\it non-zero} degree. }

{\it Let    $\underline X_0$  be simply connected, let     the Euler characteristics of $\underline X_0$ don't   vanish, $\chi(\underline X_0)\neq 0$,  let the map 
$f_1: X\to  \mathbb R^{m}$ be  area non-increasing and let the restriction of this  map to the boundary of $X$, be distance non-increasing, i.e. 
 $$\mbox {$||df_1(x)||\leq1$ for $x\in \partial X$.}$$  
If
$$Sc(g,x)- ||\wedge ^2df(x)||\cdot Sc(\underline g, f(x))\geq  0\leqno {\color {blue} 
\sf scal_{\geq }}$$
then

$$\inf_{x\in \partial X} \left (mean.curv (\partial X,x)- ||df(x)|| \cdot mean.curv(\partial\underline X_0\times \mathbb R^m , f(x))\right )\leq 0.\leqno { \color {blue} \sf mean_{\leq}}$$}

\vspace{1mm}

{\it Remarks.}  (a)  If $Sc(\underline X_0)=0$, (hence,  $\underline X_0$ is Riemannian flat) then the condition {\color {blue} \sf scal$_\geq$} reduces to $Sc(X)\geq 0$.

(b) The inequality  {\color {blue} \sf mean$_\leq$} also  yields some information for manifolds $X$ with  {\it negative scalar curvatures bounded from below}. 

For instance,  if $X$ is compact and  $Sc(X)\geq -2)$, then {\color {blue} \sf mean$_\leq$}, this 
is achieved by  applying  {\color {green!50!black}$\bigstar_{mean}$} to
maps from  $X\times S^2$ to the unit balls $B^{n+2}\subset \mathbb R^{n+2}$ (see [G(boundary) 2019]).

However, the sharp inequalities for $Sc(X)<0$, such, for instance, as {\it optimality} of the hyperspherical  radii of the boundary spheres of balls $B^n(R)$ in the hyperbolic spaces $\mathbb H^n_{-1}$, remain {\it \large conjectural}.\footnote{ This "optimality" means that if  $Sc(X)\geq -n(n-1)$ and $mean.curv (\partial X)\geq mean.curv( \partial B^n(R))$ than $Rad_{S^{n-1}}(\partial X) \leq  Rad_{S^{n-1}}(\partial B^n(R))$.}

(c) It is unknown if the spin condition on $X$ is necessary, but it can be relaxed by  requiring the universal cover of $X$, rather  than $X$ itself is spin.(This done with the $L_2$-version of  
the Goette-Semmelman  theorem}Goette-Semmelmann

And if one is content with a non-sharp bound
$$ Rad_{S^{n-1}}(Y)\leq \frac {const_n}{\inf mean.curv(Y)},$$
then one  and  can prove this without the spin assumption by the by  a capillary version of the (iterated)  warped product argument for manifolds  with boundaries \ref {warped boundary5}. \ref{capillary warped5}.

(d) Unavoidable  {\it approximation error terms}  in   the smoothing   of the corners in the doubles \DD$X$ and  \DD$\underline X_0$ make our proof of  {\color {green!50!black}$\bigstar_{mean}$}
  poorly adjusted for {\it establishing rigidity}  of   $\underline X_0\times \mathbb R^m$.
  
  For this purpose,  it  would be  better to use Lott's index theorem for manifolds with boundary.

  In fact, Lott himself  proves in [Lott(boundary)2020] a {\it non-normalized}   rigidity theorem for {\it compact}\ manifolds  $\underline X_0$ of {\it even} dimension $n$.
  
   Apparently, Lott's argument extends to the Sc- and mean.curv- normalized  case and non-compactness of $\underline X_0$  also causes no serious problem. But it is unclear how handle
   the case of odd $n$ without an approximation argument.
   
   The simplest case, where this difficulty arises is for maps from compact manifolds  $X$ to  {\it odd dimensional}  balls $B^n\subset \mathbb R^n$  and to products of such balls by tori,
$\underline X_0=B^{2k+1}\times \mathbb T^{n-2k-1}$, where {\color {green!50!black}$\bigstar_{mean}$}
applies to the universal coverings of these manifolds.

{\color {red!50!black} Possibly} on can resolve the problem with a generalized  {\it Bourguignon-Kazdan -Warner perturbation theorem}  or with (also generalized)  {\it Burkhart-Guim's regularized Ricci flow argument.}

%%%%%%%%%%%%%%%%%%%%%%%%%%%%%

\subsection {\color {blue}Lower Bounds on the Dihedral Angles of Curved Polyhedral Domains}\label {dihedral4}

%%%%%%%%%%%%%%%%%%%%%%%%%%%%%%%%

We want to generalise the above {\color {green!50!black}$\bigstar_{mean}$}  to manifolds 
$X$ with non-smooth boundaries with suitably defined mean curvatures bounded from below, where we limit ourself to manifolds with rather simple singularities at their boundaries. 

Namely,  let  $X$  and $\underline X$  be  Riemannian $n$-manifolds  {\it with corners}, which means that their boundaries $Y=\partial X$ and $\underline Y=\partial \underline X$ are  decomposed into $(n-1)$-faces $F_i$ and $\underline F_i$ correspondingly, where, locally, at all points $y\in Y$,  and $\underline y\in \underline Y$  these   decompositions are  is diffeomorphic to such decomposition of the boundary of  a convex 
$n$-dimensional polyhedron (polytope) in  $\mathbb R^n$.

Let  $f:X\to   \underline X$ be a smooth map, which is  compatible with the corner structures in $X$ and   $\underline X$: 

{\sf  $f$ sends the  $(n-1)$-faces  $F_i$ of $X$ to faces    $\underline F_i$ of $\underline X.$}
 
Assume as earlier that

$$Sc(X,x) \geq ||\wedge^2df||\cdot Sc(\underline X, f(x))\mbox {  for all $x\in  X$}\leqno {  {\color {blue} \sf scal_\geq}}$$
and replace  {\color {blue} \sf mean$_\leq$} by the opposite inequalities   applied to    {\it for all faces} $F_i\subset Y$ individually, 
$$mean.curv(F_i,y) \geq ||df||\cdot mean.curv (\underline F_i, f(y))\mbox {  for all $y\in F_i$ }.\leqno {     {\color {blue} \sf mean_\geq}}$$

 Let $\angle_{i,j} (y)$  be the dihedral angle  between the faces  $F_i$ and $F_j$ at   $y\in  F_i\cap F_i$ and let us 
 impose our  main inequality between these    $\angle_{i,j} (y)$ for all $F_i$ and $F)j$  and the dihedral   angles between the corresponding faces  faces  $\underline F_i$ and $\underline F_j$ at  the points $f(y)\in   \underline F_i\cap \underline F_j$: 
$$\angle_{i,j} (y)\leq \angle_{i,j} (f(y))\mbox  { for all } F_i,  F_j \mbox { and } y\in  F_i\cap F_j.\leqno {\mbox {\large \color {blue} $[\leq]^{\angle _{ij}}$}} $$

Besides the above,  we need to  add the following condition the relevance of which remains unclear.

Call a point $ y\in    Y=\partial X $ {\it suspicious} if  one of the following two conditions is satisfied

(i)   the corner structure of  $X$ at $y$ is {\it non-simple} (not cosimplicial), where simple means that a neighbourhood of $y$   is diffeomorphic to a neighbourhood of a point in the $n$-cube, which is equivalent to transversality of the intersection of the $(n-1)$-faces which meet 
at $\underline y$; 

(ii) there are two $(n-1)$-faces  in $X$ which contain $y$, say  $F_i\ni y$ and $F_j\ni y$ , such that
the dihedral angle $\angle_{ij}=\angle(F_i,F_j$ is $>\frac {\pi}{2}$;

\vspace {1mm}

Then out final condition says that 
$$\angle_{i,j} (y)= \angle_{i,j} (f(y))\mbox.\leqno {\mbox {\large \color {blue} $[\equiv]^{\angle _{ij}}$} }$$
for all  suspicious points $y$. 
\vspace{1mm}

  {\color {green!50!black}\DiamondSolid$\angle_{ij}$}  {  \textbf {Compact Dihedral Extremality Theorem}.} {\sf.  Let $X$ and
 $\underline X$  be compact connected orientable Riemannin {\it spin} manifolds of dimension $n$  with corners, where  {\it the curvature operator  of  $\underline X$ be non-negative}, all faces $\underline F_i\subset \partial \underline X$  are mean convex.} 
  Let  $f: X\to \underline X$ be a smooth
{\it proper corner} map, (it   respects the corner structure)  of {\it non-zero degree} and let $f$  
satisfy the four conditions 
 {\color {blue}   {\color {blue} \sf scal$_\geq$} and  {\color {blue} \sf mean$_\geq$},  {\large \color {blue} $[\leq]^{\angle _{ij}}$}} and  {\large \color {blue} $[\equiv]^{\angle _{ij}}$}

{\it If  the universal covering of   $\underline X$  has non-zero   the Euler characteristics $\chi(\tilde {\underline X})\neq 0$, then $f$ satisfies the equalities corresponding to the inequalities  
{\color {blue}   {\color {blue} \sf scal$_\geq$}, {\color {blue} \sf scal$_\geq$},   {\large \color {blue} $[\leq]^{\angle _{ij}}$}}:
  $$Sc(X,x) = ||\wedge^2df||\cdot Sc(\underline X, f(x))\mbox {  for all $x\in  X$},$$
$$mean.curv(F_i,y) = ||df||\cdot mean.curv (\underline F_i, f(y))\mbox {  for all  $y\in F_i$} ,$$
  $$\angle_{i,j} (y)= \angle_{i,j} (f(y))\mbox { for all } F_i,  F_j \mbox { and } y\in  F_i\cap F_j. $$}
\vspace {1mm}

\vspace {1mm}

{\it About the Proof.} This is shown by smoothing the boundaries of $X$  and of   $\underline X$ and applying    {\color {green!50!black}$\bigstar_{mean}$} from the previous section, to the universal covering of  $\underline X$  and the corresponding (induced)  covering of $X$ \footnote{ If, instead of "$X$ is spin"  we only assume  "the  universal covering of $X$ is spin", then we pass to this universal covering of $X$ and use there the $L_2$-index theorem. }  where an essential feature of 
non-suspicious points  follows from the following\vspace {1mm}

{\it \large Elementary Lemma.} {\sf Let $\Delta\subset S^n$ be a spherical simplex with all edges of length $\geq l \geq \frac{\pi}{2}$. Then there exists a continuous family of simplices 
 $\Delta_t\subset S^n$, $t\in [0,1]$  with the following properties.} 
 
{\sl  $\bullet$  $\Delta_0=\Delta$ and  $\Delta_1$ is a regular simplex with the edge length $l$;

$\bullet$ all $\Delta_t$  have the  edges of length $\geq l$;

$\bullet$ $\Delta_{t_2}\subset \Delta_{t_1}$  for $t_2\geq t_1$;

$\bullet$ for each $t<1$ there exists an $\varepsilon >0$, such that $n$ (out of $n+1$) vertices
 of $\Delta_{t+\varepsilon}$ coincide with those of $\Delta_{t}$.}\footnote {This Lemma explains the role of the condition   {\large \color {blue} $[=]^{\angle _{ij}}$} in our proof. The conclusion of the Lemma fails, in general to be true for obtuse angles (it seems OK if there is a single obtuse angle at each vertex , e.g, as it is  for products of convex polygons) 
 but it remains unclear if this condition is needed  for the validity of the theorem itself.
 }\vspace {1mm}

{\it The proof} of the lemma is a high school exercise while {\it construction of  adequate smoothing  of $X$ }
with the help of this lemma, which  is straightforward and boring, will be given elsewhere.

\vspace {1mm}

Notice that  the {\Large \color{blue} $\times\hspace {-0.6mm} \blacktriangle^i$}-{\large \it \color{blue} Inequality} from section \ref {corners3}, which says that\vspace {1mm}

 {\sf convex  polyhedra $\underline X\subset \mathbb R^n$ with the dihedral angles $\leq \frac {\pi}{2}$ admit no deformations which would {\it decrease their dihedral angles} and simultaneously {\it increase the mean curvatures of their faces}, }\vspace {1mm}

 \hspace {-6mm}is an immediate corollary of  {\color {green!50!black} \DiamondSolid$\angle_{ij}$}. \vspace {1mm}

{\it Two Problems.} 1.  There is {\color {red!50!black}little doubt} that the above  extremal
 manifolds with corners $  \underline X$  are {\it rigid}, but our argument, as  we explained this in the previous section is, technically,  not good  enough for proving it,   and no index theorem theorem for general manifolds 
with corners is available, at  least not at the present moment.

2. It remains unclear what is the {\it full  class} of extremal polyhedra and manifolds with corners in general, but the following  generalization of  {\color {green!50!black} \DiamondSolid$\angle_{ij}$} is 
easily available.

\vspace {1mm}

{\it  Fundamental Domains of    Reflection Groups.} 
What  underlies
the double \DD-construction, $X\leadsto $ \DD$X$ in  the proof of the  {\color {green!50!black} \DiamondSolid$\angle_{ij}$  theorem}  is   the  doubling  $S^n$= \DD$S^n_+$, which is 
associated with the  reflection of $S^n$ with respect to the equatorial  subsphere.

With this in mind, one can generalise everything from this section to general reflection groups, 
including spherical, Euclidean and  hyperbolic  ones, (such as 
 we met     in section  \ref {reflection3}) and also to  products of these.\vspace{1mm}

{\it \textbf {Example of Corollary.}} {\sf  Let $X$ be a compact  manifold with corners, where the  (combinatorial) corner structure is isomorphic to that of the  product of an $(n-m)$-simplex $\blacktriangle$ with the rectangular  fundamental domain $\blacksquare$ 
(orbifold) of a cocompact reflection group in an aspherical $m$-manifold}.\footnote {These exist for all $m\geq 4$ by Michael Davis 1983 theorem, see his lectures  [Dav(orbifolds) 2008] and references therein.}  \vspace{1mm}

{\sl If $X$ is spin, than it admits {\sf no Riemannian metric $g$},  such that 
$Sc(g)\geq 0$,  where all faces have $mean.curv_g\geq 0$ and where the  dihedral angles are  
smaller than the corresponding angles 
 in the  product of the  regular Euclidean simplex $\blacktriangle$ by  $\blacksquare$  with $\frac {\pi}{2}$ dihedral angles.}

\vspace{1mm}

{\it Problem with Rigidity.} If not for  

%%%%%%%%%%%%%%%%%%%%%%%%

\subsection {\color {blue} Stability of Geometric   Inequalities with $Sc\geq \sigma$  and  Spectra of Twisted Dirac Operators.} \label {stability4}

%%%%%%%%%%%%%%%%%%%%%%%%%%%%%%

Sharp geometric inequalities,  as we explained in section \ref {stability3},  beg for  a company of  their nearest neighbours. 

For instance, the {\color {blue} Euclidean isoperimetric inequality} for bounded  domains $X\subset \mathbb R^n$,  which says that 
$$vol_n(X) \leq \gamma_n vol_{n-1}(\partial X)^{\frac {n}{n-1}}\mbox {  for }\gamma_n= \frac {vol(B_n)}{vol_{n-1}(S^{n-1})^{\frac {n}{n-1}}},$$
goes along with the following.\vspace {1mm} 

{\large {\sf A}. {\color {blue} \it Rigidity.} {\sl  If $vol_n(X) = \gamma_n vol_{n-1}(\partial X)^{\frac {n}{n-1}}$, then $X$ is  a ball.}}

\vspace {1mm} 

{\large  {\sf B}. \color {blue} \it  Isoperimetric Stability.}  {\sf Let $X\subset \mathbb R^n$ be  a bounded  domain 
with  $vol_n(X) = vol_n (B^n)$ and $vol(\partial X)\leq vol_{n-1}(S^n)+\varepsilon$.}

{\sl Then there exists a  ball $B=B_x^n(1+\delta)\subset \mathbb R^n$ of radius $\delta$  with center $x\in X$, where  $\delta\underset {\varepsilon\to  0} \to 0$, such that the volume of the difference satisfies
$$vol_n(X \setminus B) \leq  \delta_1, $$
and, moreover, 
$$ vol_{n-1}(\partial B\cap X)\leq \delta_2, \mbox { and } vol_{n-2}(\partial B\cap \partial X)\leq \delta_3,$$
where $$ \delta_1, \delta_2, \delta_3\underset {\varepsilon \to 0} \to 0.$$}

(Unless $n=2$ and $X$ is connected,  there is no bound on the diameter of $X$, but 
the constants $\delta,\delta_1,\delta_2,\delta_3$ can be explicitly evaluated even  for moderately large $\varepsilon$.)

\vspace{1mm}

In the case of sharp  scalar curvature inequalities, their poofs by Dirac theoretic methods\footnote {See  [Llarull(sharp estimates)  1998], [Min-Oo(Hermitian)  1998], 
  [Goette-Semmelmann(symmetric) 2002],  [Listing(symmetric spaces) 2010],  [Zeidler(width) 2020],  [Zhang(area decreasing) 2020], [Lott(boundary) 2020], [Guo-Xie-Yu(quantitative K-theory) 2020].}(more or less)  automatically deliver   rigidity. 
For instance,\vspace{1mm}

{\color {blue} \huge $\star$}  {\sl if a manifold  $\underline X $ homeomorphic to $S^n$,  besides having $curv.oper(\underline X )\geq 0$ has $Ricci(\underline X )>0$ and if  $X$ is a closed orientable spin Riemannian manifold with  $Sc(X)\geq n(n-1)$ then, 
all  {\sf smooth 1-Lipschitz }  maps $X\to \underline X$ of non-zero degrees are {\sf  isometries}.} \footnote{Even if Ricci vanishes somewhere, one still may  have a satisfactory description of the extremal cases. For instance, if  $\underline X =(S^{n-m} \times \mathbb R^m)/\mathbb Z^m$, e.g.  $\underline X  =S^{n-m}\times \mathbb T^m$, then all (orientable spin) $X$ with $Sc(X)\geq Sc( \underline X )=(n-m)(n-m-1)$, which admit maps $f: X\to \underline X $ with  $deg(f)\neq 0$, are {\it locally isometric} to $\underline X $ (albeit the map $f$ itself doesn't have to be a  local isometry. }

\vspace{1mm}

What we want to understand next is what happens if the inequality  $Sc(X)\geq n(n-1)$
is relaxed to  $Sc(X)\geq n(n-1)-\varepsilon$ for a small  $\varepsilon>0$, where an application of thin surgery\ref{thin1} delivers the  following.\vspace {1mm}

{\color {blue}\it Example.}\footnote {Compare with  [GL(classification) 1980], 
[BaDoSo(sewing Riemannian manifolds) 2018]  and section 2 in [G(101) 2017].}
  Let $\Sigma\subset S^n$ be a compact smooth submanifold of dimension $\leq n-3$.
Then there exists an arbitrary small $\varepsilon $-neighbourhood  $U_\varepsilon=U_\varepsilon (\Sigma)\subset S^n$  with a smooth boundary $\partial_\varepsilon=\partial U_\varepsilon$ and a family of smooth metrics $g_{\varepsilon, \epsilon}$ on the double \vspace {1mm}

\hspace{19mm} \DD$(S^n\setminus U_\varepsilon) =(S^n\setminus U_\varepsilon) \cup_{\partial_\varepsilon} (S^n\setminus U_\varepsilon)$, 

\vspace {1mm}

\hspace {-6mm}
where $Sc(g_{\varepsilon, \epsilon})\geq n(n-1)-\varepsilon -\epsilon$   and which, for  $\epsilon \to 0$,  {\it uniformly converge} to the natural continuous Riemannian metric on    \DD$(S^n\setminus U_\varepsilon (\Sigma)$.

Moreover, if $\Sigma\subset S^n$ is contained in a hemisphere, then -- this follows from the {\it spherical Kirszbraun theorem} -- the (double) manifolds \DD$(S^n\setminus U_\varepsilon, g_{\varepsilon, \epsilon})$ admit 1-Lipschitz maps to the sphere $S^n$ with degrees one, for all sufficiently small $\varepsilon >0$ and , $ \epsilon= \epsilon(\varepsilon)\underset {\varepsilon \to 0}\to 0$.

For instance, if $n\geq 3$ and  $\Sigma$ consists of a single point, then 
  \DD$(S^n\setminus U_\varepsilon)$, that is the connected sum $S^n\#S^n=S^n\#_{S^{n-1}(\varepsilon)}S^n$ of the sphere $S^n$ with itself (where the $\varepsilon$-sphere  $S^{n-1}(\varepsilon)$ serves as $\partial_\varepsilon$ and $S^n\#S^n$ is  homeomorphic to $S^n$), admits, for small $\varepsilon$,  a 1-Lipschitz map to $S^n$ with degree 2.
  
Furthermore,  iteration of the connected sum construction,  delivers manifolds (topologically spheres)
$$(S^n)^{k\#_\varepsilon}=\underset {k}{\underbrace {S^n\#_{S^{n-1}(\varepsilon)}S^n\#...\#_{S^{n-1}(\varepsilon)}S^n}}m$$
which carry metrics with $Sc(S^n)^{k\#_\varepsilon}\geq n(n-1)-\varepsilon-\epsilon$ and, at the same time,  admit maps to $S^n$ of degree $k$, where   these maps  are  1-Lipschitz everywhere and which are  locally isometric  away from 
$\sqrt\varepsilon$-neighbourhoods of  $k-1$ $\varepsilon$-spherical  "necks" in $(S^n)^{k\#_\varepsilon}$.
  
  (For general $\Sigma$ and even $k$ one has such maps  $f$ with  $deg(f)=k/2.$

\vspace{1mm}

{\it \large Conjecturally,} this example faithfully represents possible geometries of closed  Riemannian $n$-manifolds 
$X$ with  $Sc(X)\geq n(n-1)-\varepsilon $, which admit 1-Lipschitz maps to the unit sphere $S^n$,  but only the following two, rather superficial, results of this kind are available.\vspace{1mm}

{\large  \sf \color {blue}\textbf 1.} {\sf Let $X =(X,g)$ be a closed oriented Riemannian  spin $n$-manifold with $Sc(X)\geq n(n-1)-\varepsilon $ and let $f: X\to\underline X= S^n$ be a smooth 1-Lipshitz  map of degree $d\neq 0$ and 
let $J_f(x)=\wedge^ndf$  denote the Jacobian of $f$}.

Let $X_{\leq\lambda}\subset X $ denotes the subset, where $|J_f(x)|\leq \lambda$, for some $\lambda <1$. }

{\sl Then the signed 
$f$-volume of   $X_{\leq\lambda} $  satisfies
$$ |vol_f (X_{\leq\lambda})|=_{def} \left |\int_{X_{\leq\lambda}}J_f(x)dx \right| \leq c_{\lambda,n, \tilde V}(\varepsilon)\underset {\varepsilon \to 0}\to 0. \footnote{This was incorrectly stated in an earlier  version of this text  for {\it non-signed} volume of $f$, that is 
$\int |J_f(x)|dx$; the error was pointed out to me by Bernhard Hanke.}
 \leqno {\color {blue} [|X_{\leq\lambda}|\leq] }$$}
 
(Observe that  since $f$ is 
1-Lipschitz, $|J_f|\leq 1$ and $1-|J_f(x)|^{-1}$ measures the the distance from    the differential  $df(x):T_x(X)\to T_{f(x)}(\underline X)$ to being    isometry.)

 {\it Sketch of the Proof.} Since  the twisted Dirac  $D_\otimes$ in Llarull's rigidity argument from  [Llarull(sharp estimates)  1998]  has non-zero kernel, its  square  $D^2_\otimes$ is non-positive (we assume here that  $n=dim(X)=dim (\underline X)$ is even), and, by the  Bochner-Schr\"odinger-Lichnerowicz-Weitzenb\"ock formula  (that is above  $[D_\otimes^2]_f$), this  implies non-positivity of 
$$\nabla^2 +\frac {1}{4} Sc(X) + {\cal R}_\otimes.$$

Consequently, $-\Delta_g -\frac {1}{4} (\varepsilon + (1-\underline l(x)))$, where $\Delta_g$  is an ordinary Laplace  on $X =(X,g)$,  also non-positive, since the  coarse (Bochner)
Laplacian $\nabla^2$ is "more positive" than the  (positive) Laplace(-Beltrami)  $-\Delta$ as it follows from  the {\it Kac-Feynman} formula and/or from  the {\it Kato inequality.}

(In general,  this applies in the  context of the above   rigidity theorem {\color {blue} \huge $\star$} and yields non-positivity of $-\Delta_g -\frac {1}{4} (\varepsilon + \underline C(1-\underline l_f(x))$
with $\underline C$ depending on the smallest eigenvalue of $Ricci(\underline X)$.)

\vspace{1mm}
 
 In order to   extract required  geometric information concerning the metric $\tilde g$  from this property of the metric $g$, we observe that the essential  part of $X$, that is the one, where we  need
 to bound from below the $L_2$-norms of  the $g$-gradients of functions $\phi(x)$ (to which the above  $\Delta_g$ applies) is 
 where 
 $$\lambda \geq \underline l_f(x) \geq \lambda_{\tilde V}>0$$ 
 for some $\lambda_{\tilde V}>0$, and  where the geometries of  $g$ and of $\tilde g$  are mutually  $(\lambda_{\tilde V})^{-1}$-close. 
 
Thus,  the relevant  lower $g$-gradient estimate for $\phi(x)$ comes from the isoperimetric inequality for $\tilde g$ which, in turn, follow from such an inequality  in $\underline X$, that is the sphere in the present case. (Filling in the details is left to the reader.)
 
 \vspace {1mm}
{\it Remark.} (a) The above example shows that the  $g$-volume of   $X_{\leq\lambda}\subset X $ can be large and that the bound on $\tilde V$ concerns not only the subset  $X_{\leq\lambda}$ but  its complement $X\setminus  X_{\leq\lambda}$ as well.
\vspace{1mm}

{\it Corollary } + {\it Question.} (a)  {\sf Let $X$ be a closed orientable Riemannian spin $n$-manifold   with $Sc(X)\geq n(n-1)$
and let  $f:X\to S^n$ a (possibly non-smooth!)  1-Lipshitz map of degree $\neq 0$.}

{\it If the map $Y$ is a homeomorphism, then it is an isometry.}\vspace{1mm}

(b) {\sf Is this remain true for {\it all} 1-Lipshitz maps?} \vspace {1mm} \vspace {1mm}

The  inequality {\color {blue}$ [|X_{\leq\lambda}|\leq]$ }  doesn't take advantage of $deg(f)$ when this is large, but the following proposition does just that.

 \vspace{1mm}
 
{\large  \sf \color {blue} \textbf  2.} {\sf Let $X$ be a compact oriented  Riemannian spin $n$-manifold   with a boundary $Y=\partial X$, such that $Sc(X)\geq n(n-1)+\varepsilon$, $\varepsilon >0$.

Let $f:X\to S^n$ be a smooth map, which is constant on $Y$,   which  is {\it area contracting away from the a neighbourhood $U\subset X$  of $Y=\partial X\subset X$, 
$$ || \wedge ^2df(x)||\leq 1\mbox { for  all } x\in X\setminus U,$$}
 and where
$$ || \wedge ^2df(x)||\leq C_o\mbox { for  all  $x\in X\setminus U$ and some constant $C_o>0.$} $$}
 {\sl Then the degree of $f$ is bounded by a constant $d$ depending only on $U$ and on $C_o$,
 $$|deg(f)|\leq d= const_{U,C_o}.$$}
 
{\it Sketch of the Proof.} (Compare with \S\S$5\frac{1}{2}$ and  6 in [G(positive) 1996].)
Let $s(x)$ be the  (Borel) function on $X$ which equals to $\varepsilon$ away from $U$ and is equal to $E=-C_n\times C_o$ on $U$ for some  universal $C_n\approx n^n$.

Then arguing (essentially)  as in the first part of the above proof, we conclude that the spectrum of the  
$-\Delta +s(x)$ on the (smoothed) double  \DD$(X)$  contains at least $d=deg(f)$ negative  eigenvalues.

This an easy argument  would deliver $d$ eigenvalues $\lambda_i$  of the  $-\Delta$ on  \DD$(U)$, where the corresponding eigenfunctions vanish on the two copies  of the boundary 
of $U$ in $X$ (but not, necessarily on $Y$), and such that 
$\lambda_i \lessapprox E$. 

This would yield the required bound on $d$. (Here again, the details are left to the reader.)

\vspace{1mm}

{\it Remark} + {\it Example} + {\it Two Problems.} (a) If the boundary of $Y=\partial X$ admits an orientation reversing involution, then the constancy of $f$ on $Y$ can be relaxed  to  $ dim(f(Y))\leq n-2$, where the constant  $d$ will have to depend on the geometry of this involution and of the map $Y\to S^n$.
 
 (It is unclear if the existence of such an involution is truly necessary.)

(b) This (a)    apply, for instance,  to coverings  $X=\Sigma^2_{d, \delta}$ of the 2-sphere minus two $\delta$-discs  as well as to the products of these $\Sigma^2_{d, \delta}$ 
with the  Euclidean ball $B^{n-1}(R)$  of radius $R>\pi$.
 
(c) {\sf What are the sharp and/or comprehensive versions of these {\large  \sf \color {blue} \textbf  1} and  {\large  \sf \color {blue} \textbf  2}}?

(d) {\sf Let $Y$ be a homotopy sphere of dimension $4k-1$, which bounds a Riemannian manifold  $X$ with $Sc\geq \varepsilon >0$. Give an {\it effective} bound on the $\hat A$-genus of $X$ in terms of the geometry of $Y$ and its second fundamental form $h={\rm II}(Y\subset X)$ and study the resulting invariant 
$$Inv_\varepsilon(Y,h)=\sup_X|\hat A(X)|,\mbox {  where } \partial X=Y, \mbox { }  Sc(X)\geq \varepsilon, 
  \mbox { }  {\rm II } (Y\subset X)=h. $$ }

%%%%%%%%%%%%%%%%%%%%%%

\subsection { \color {blue} Dirac Operators  on Manifolds with Boundaries} \label {boundaries4}

  %%%%%%%%%%%%%%%%%%%%%%%%%%%%%%%%%%%
 
  When I was delivering these lectures in the Spring 2019,    all known  relevant for us  index  theorems for (twisted) Dirac operators  $\cal D$ directly  applied only to {\it complete} Riemannian manifolds. \footnote {This is not quite true: Roe partitioned index theorem and its generalization do allow boundaries, see [Roe(partial vanishing) 2012],  [Higson(
cobordism invariance) 1991], [Schick-Zadeh(multi-partitioned) 2015], [Karami-Zadeh-Sadegh(relative-partitioned) 2018]  and section \ref {Roe3}.}
  But 
  then Cecchini, Guo-Xie-Yu  and 
Zeidler 
 \footnote {[Cecchini(long neck) 2020], [Guo-Xie-Yu(quantitative K-theory) 2020],
[Zeidler(bands)  2019],  
[Zeidler(width)  2020], [Cecchini-Zeidler(generalized Callias) 2021], [Cecchini-Zeidler(scalar\&mean) 2021].} 
have developed \vspace{1mm}

  {\sl \color {blue!59!black} an index theory for manifold with boundary including the solution of the
  
   long neck problem for spin manifolds  by Cecchini} (see  section \ref {Roe}).\vspace{1mm} 
 
 Even though, much(all?) what is presented  in this and the following  sections 4.6.1-4.6.5  may  follow from   the recent results    of these authors, 
 we  keep it as it was originally written, since  this  suggests an additional perspective on the role of the Dirac  operator in the geometry   of scalar curvature.\vspace{1mm}

As far as   the scalar curvature  is concerned, all the   index theorems  are needed for  is delivering   
 {\it \color {blue!40!black} non-zero  harmonic} or {\it approximately harmonic} (often twisted)  spinors  on Riemannian manifolds $X$ under certain  certain geometric/topological conditions on $X$, which, a priori,  have nothing to do with the  scalar curvature  but which are eventually  used to obtain upper bounds on $Sc(X)$ via the (usually twisted) Bochner-Schr\"odinger-Lichnerowicz-Weitzenb\"ock formula.

The index theorems for Dirac operators on  closed manifolds 
can  yield  a non-trivial information on existence of approximately harmonic spinors on non-complete manifolds as well as on manifolds with boundaries, where the main issue, say for manifolds with boundaries,  can be formulated as follows.

\vspace {1mm} 

{\it \color {blue}\large   Spectral  $\mathcal D^2$-Problem.}   {\sf Let $X$ be a compact  Riemannian spin manifold with a boundary
and $L\to X$ be a (possibly infinite dimensional Hilbert)  vector bundle with a unitary connection.

Under which  geometric/topological 
conditions does  the first eigenvalue of the twisted Dirac  $\mathcal D_{\otimes L}$  on $X$ with the zero boundary condition is  $\leq\lambda>0$?

In other words, when does $X$  support a smooth  {\it non-zero} twisted spinor $s:X\to \mathbb S(X)\otimes L$, which {\it vanishes on  the boundary} of $X$ and  such that 
$$    \int _X\langle  \mathcal D_{\otimes L}^2(s(x)), s(x)\rangle dx  \leq \lambda^2\int_X ||s(x)||^2dx\leqno {\mbox {\color {blue}\large   \ScissorHollowRight\hspace{-0.2mm}$_{\lambda}$}}$$ 
for a given constant $\lambda\geq 0$?}\footnote {Recall that the first eigenvalue of the Dirichlet problem is
the infimum of $ \int _X||\mathcal D_{\otimes L}(s(x))||^2dx$ taken over all $L$-twisted  spinors $s(x)$, such  that  $s|\partial X=0$ and 
 $\int _X||s||^2dx=1$.}

\vspace {1mm} 

{\it  Motivating Example.}  If  $X$ is  obtained from  a complete  manifold $X_+\supset X$ by cutting away $X_+\setminus X$, and 
if  $X_+$  carries a non-vanishing (twisted) $L_2$-spinor $s_+$ delivered by applying the relative index theorem, then  the cut-off  spinor $s=\phi\cdot s_+$,  for a "{\it slowly} decaying" positive function $\phi$  with supports in $X$ satisfies  {\color {blue}\large   \ScissorHollowRight\hspace{-0.2mm}$_{\lambda}$} with "{\it rather small }" $\lambda$.

\vspace {1mm}

\vspace {1mm} 

{\it \color {blue} Potential  Corollary.} Since   $$\mathcal D_{\otimes L}^2(s)\geq \nabla^2_{\otimes L}(s) + \frac {1}{4} Sc(X)(s) - 
const'_n |curv|(L)$$
by the Bochner-Schr\"odinger-Lichnerowicz-Weitzenb\"ock formula   
and since 
 $$\int\langle  \nabla^2_{\otimes L}(s), s\rangle= \int_X\langle \nabla_{\otimes L}(s),  \nabla_{\otimes L}(s)\rangle\geq 0$$
for $s_{|\partial X}=0$, the inequality  {\color {blue}\large \ScissorHollowRight \hspace{-0.2mm}$_\lambda$} implies 
\vspace {1mm}\vspace {1mm}

$$ \inf_xSc(X,x)\leq \frac {4const_n}{\rho^2} + 4const'_n|curv|(\nabla).\leqno {\mbox {\color {blue}\large   \ScissorRight\hspace{-0.2mm}$_{Sc}$}} $$
for some universal  positive  constants $const_n$ and  $const'_n$.

\vspace {1mm}

From a geometric perspective, the role of above is to advance the solution of the following.

{\it \color {blue} \large  Long Neck Problem.} Let $X$ be an orientable (spin?)  Riemannian  $n$-manifold with a boundary  and $f:X\to S^n$ be a  smooth  area decreasing map. \vspace {1mm}

{\sf What kind of a {\sl lower bound on} $Sc(X,x)$   and a {\sl lower bound on}  the "length of the neck" of $(X,f)$, that is 

\hspace {-4mm}{\it the distance between the 
support of the differential of $f$ and the boundary of $X$,}

\hspace {-6mm}would make $deg(f)=0$?}

An instance of a desired result\footnote {This is settled for spin manifolds in [Cecchini(long neck) 2020].} would be 
{\color {blue} $$ [Sc(X)\geq n(n-1)] \&  [dist(supp(df), \partial X)\geq const_n]\Rightarrow deg(f)=0, 
$$} 
but it is more realistic to expect a weaker implication
{\color {blue} $$ [Sc(X)\geq n(n-1)] \&  [dist(supp(df), \partial X)\geq 
const_n{\color {green!50!blue}\cdot \sup_{x\in X}| ||df(x)||}]\Rightarrow deg(f)=0. $$}

In fact, Roe's proof  of the partitioned index theorem as well as 
the proof of the relative index theorem, e.g. via the {\it finite propagation speed} argument, combined with Vaffa-Witten kind spectral  estimates (see $6\frac {1}{2}$ in [G(positive) 1996]) suggest that 

{\sf if a compact    orientable Riemannian  spin manifold of even dimension $n$  with boundary admits a 
 a   smooth  map $f:X\to S^n$, which is locally constant on the boundary of $X$  and which has   {\it non-zero degree}, then there exists a {\it non-zero spinor} $s$,  twisted  
 with the pullback bundle $L= f^\ast(\mathbb S(S^n))$ such that   $s$ {\it vanishes on the boundary} $\partial X$ and which satisfies  {\color {blue}\large   \ScissorHollowRight\hspace{-0.2mm}$_{\lambda}$},
$$    \int _X\langle  \mathcal D_{\otimes L}^2(s), s\rangle  \leq \lambda^2\int_X ||s||^2dx,$$
where 
$$\lambda\leq const_n \frac  {\sup_{x\in X}||df(x)||}{dist(supp (f), \partial X)}.$$}

\vspace {1mm}

This still remains problematic, but we prove in the sections below some inequalities  in this regard 
for manifolds  $X$ with certain restrictions on their local geometries.\footnote {An influence   of the metric geometry of a Riemannian manifold  $X$ on the spectra of  twisted Dirac operators on $X$ is briefly duscussed in 
\S6 of [G(positive) 1996].}

 %%%%%%%%%%%%%%%%%%%%%%%%%%

 \subsubsection{ \color {blue} Bounds on  Geometry  and Riemannian Limits}\label {limits4}

Some properties of manifolds $X$  with boundaries trivially follow  by a limit argument from the corresponding properties  of complete  manifolds  as follows.

A sequence of manifolds $X_i$ marked   with distinguished points $\underline x_i\in X_i$ is said to 
{\it Lipschitz  converge} to a marked Riemannian manifold $(X_\infty, \underline  x_\infty)$, if 

 there  exist  {\it $(1+\varepsilon_i)$-bi-Lipschitz}  maps \footnote{Here and below  
 "$\lambda$-bi-Lipschitz"  is understood as the $\lambda $-bound on the {\it norms of the differentials}  of our maps and their inverse.}  from the balls 
$B_{\underline x_i}(R_i)\subset X_j$  to the balls $B_{\underline x_\infty}(R_i)\subset X_\infty$, say 
$$\alpha_i: B_{\underline x_i}(R_i)  \to B_{\underline x_\infty}(R_i+1),$$
which send $\underline x_i\to \underline x_\infty$ and 
where 
$$\mbox { $\varepsilon_i\to 0$ for $i\to\infty$}.$$

Observe that if 
$$dist(\underline x_i, \partial X_i)\to \infty\mbox { for }  i\to \infty, $$
then {\it the limit manifold $X_\infty$ is complete.}

\vspace {1mm}

{\color {blue} \small $\bigstar$} {\it Cheeger Convergence Theorem.}   {\sf If the (local)  $C^k$-geometries  of Riemannian  manifolds $X_i$ at the points  
$x_i \in B_{\underline x_i}(R_i)$ for $R_i\to \infty$   are bounded (as defined below)  by $c(dist(x_i, \underline x_i))$ for some continuous function $ b(d)$, $d \geq 0$ independent of $i$, 
then  some subsequence of $X_i$ converges to a $C^{k-1}$-smooth Riemannian manifold $X_\infty$.}
 
 See   [Boileau(lectures)  2005]  for  the proof and further  references.

\vspace {1mm}

{\it Definition of Bounded Geometry.}  The $C^k$-geometry of  a smooth Riemannian $n$-manifold $X$ is  bounded    by a constant $\b geq 0$  at a point $x\in X$, if the $\rho$-ball  $B_x(\rho)\subset  X$ for $\rho=\frac {1}{b}$  admits a  smooth $(1+b)^2$-bi-Lipschitz  map   $\beta: B_x(\rho)\to \mathbb R^n$,  such that the  norms of the $k$th  covariant derivatives of $\beta$ in  $B_x(\rho)$ are bounded by $b$.

\vspace {1mm}

Notice that  the {\it traditionally  defined     bound} on geometry  in terms of the curvature and the injectivity radius of $X$, implies the above one:\vspace {1mm}

  {\sl if the norms of the curvature tensor of $X$ and its $k$th-covariant derivatives are   bounded by $\beta^2$ and there is {\it no geodesic loop} in  $X$ based at $x$ 
of length $\leq \frac {1}{\beta}$, then  (the proof is very easy)  the $C^{k+1}$-geometry of $X$ at $x$ is bounded by $b(\beta)$ for some universal continuous function  $b(\beta)=b_{n,k}(\beta).$}

\vspace {1mm}

{\small \color {blue} \FiveStarOpenCircled}  {\it Application of {\color {blue} \small $\bigstar$} to Scalar Curvature.} {\sf Let $ b=b(d)\geq 0$, $d>0$, be a continuous function  and let $(X, \underline x\in X)$ be a  marked compact Riemannian $n$-manifold  with a  boundary, such that 
the local geometry of $X$ at $x\in X$ is bounded by $b(dist(x, \underline x))$  and let
$$R=dist (\underline x, \partial X).$$}

{\sf Let $d_0$ be a positive number  and  let  $f: X \to S^n$ be a  smooth    {\it area decreasing}   map which is constant  within distance $\geq d_0$ from  $\underline x\in X$ and which has {\it non-zero} degree. }

\vspace {1mm}

{\color {blue} A.}  {\it If $X$ is spin and $n=dim(X)$ is even,  then  there exists a spinor $s$ on $X$ twisted with the induced spinor bundle $L=f^\ast(\mathbb S(S^n))\to X$, such that $s$ vanishes on the boundary $\partial X$ of $X$  and such that  
$$    \int _X\langle  \mathcal D_{\otimes L}^2(s), s\rangle  \leq \lambda(R)^2\int_X ||s||^2dx$$

 where $\lambda= \lambda_{n, b, d_0}(R)$ is a certain universal function in $R$, which asymptotically  vanishes  at infinity,
$$\lambda(R)\underset {R\to \infty} \to 0.$$}

 {\color {blue} B.} {\it The scalar curvature of $X$ is bounded by
 $$\inf_{x\in X}\leq n(n-1) + \lambda'_{n, b, d_0} (R),$$
where, similarly to the above $\lambda$, this  $\lambda'( R)\to 0$ for $R\to \infty$.}(One can actually  arrange  $ \lambda'=\lambda$.)

\vspace{1mm} 

{\it Proof.} According to Cheeger's theorem, if $R=dist (\underline x, \partial X$) is sufficiently large , then $X$ can be well   approximated by a complete manifold $X_\infty$,   
where such an $X_\infty$   supports a non-zero $L$-twisted harmonic spinor $s_\infty$ by the relative index theorem. 

Then  this $s$ can be truncated to $s_i$ by multiplying it with  a slowly decaying function on $X$ with compact  support and then  transporting it   to the required  spinor on $X$.

This takes care of {\color {blue} A}  and  {\color {blue} B}  follows by Llarull's inequality.

\vspace {1mm}

{\it Remarks.} (a)   The major  drawback of {\small \color {blue} \FiveStarOpenCircled} is an  excessive presence and   non-effectiveness of the   bounded geometry condition. 

We {\sl \color {red!30!black} don't know what the true dependence of $\lambda$ 
on the geometry of $X$ is}, but we shall prove several inequalities in the following sections that suggest what one may expect in this regard.

(b) If the "area decreasing"  property  of the above  map $f:X\to S^n $ is strengthened to  "1-Lipschitz", then a version of  {\color {blue} B} follows from the {\sf double puncture theorem} (see sections  \ref {punctured3} and  \ref {log-concave5}), which needs neither    spin nor  the  bounded geometry 
conditions.
%%%%%%%%%%%%%%%%%%
%%%%%%%%%%%%%%%%%%%

%%%%%%%%%%%%%%%%%%%

 \subsubsection{\color {blue} Construction of Mean Convex Hypersurfaces and Applications to $Sc>0$} \label {4.6.2}%%%%%%%%%%%%
%%%%%%%%%%%%%%%%%%%%%%
%%%%%%%%%%%%%%%%%%%
%%%%%%%%%%%%%%%%%%%

Since doubling  of  manifolds with mean convex boundaries preserves positivity of the scalar curvature (see section \ref {mean1}), 
 some problems concerning $Sc>0$ for manifolds $X$ with  boundaries can be reduced to the corresponding ones for closed manifolds by doubling {\it mean convex}  domains $X_{\hspace {-0.7mm}\partialvartoiint} \subset X$   
across their boundaries $\partial X_{\hspace {-0.7mm}\partialvartoiint}.$

To make use of this,  
we shall present below some  a simple criterion for the existence  of such $X_{\hspace {-0.7mm}\partialvartoiint}$ and apply this for  establishing 
 effective versions of the above {\color {blue} B}.
 \vspace {1mm}

 Let $X$ be a compact $n$-dimensional {\it Riemannian band} (capacitor),  that is  the boundary of 
    $X$ is divided into two disjoint subsets, that are certain unions of boundary components of $X$,
    $$\partial  X=\partial_- \cup \partial _+$$
  and let us give a condition for the existence of a    domain $X_{\hspace {-0.7mm}\partialvartoiint} \subset X$
 which {\it contains} $\partial_-$  and the boundary of which is smooth and has {\it positive} mean curvature. 
  
 {\it \color {blue} \textbf {Lemma}.} {\sf Let the boundaries of all domains $U\subset  X$,  which {\it contain the $d_0$-neighbourhood of $\partial X_-$} for a given $d_0< dist(\partial_-, \partial _+)$, satisfy
 $$ vol_{n-1}(\partial U)>vol_{n-1}(\partial _-)\leqno {\color {blue}   [\ast_1]}$$ 
  and let all {\it minimal\footnote {Here "minimal" means "volume minimizing" with a given boundary.} hypersurfaces} $Y\subset  X$, the  boundaries of which are  contained in $\partial _+$ and which themselves  contain points $y\in Y$  far away from $\partial_+$, namely, such that   
  $$dist (y,\partial_+) \geq dist(\partial_-, \partial _+)-d_0,$$ satisfy
  $$ vol_{n-1}(Y)> vol_{n-1}( \partial_-).\leqno {\color {blue}   [\ast_2]}$$}
  
  {\sl  Then there exists  a    domain $X_{\hspace {-0.7mm}\partialvartoiint} \subset X$
 which {\it contains} $\partial_-$  and such that  the boundary of which is smooth with {\it positive} mean curvature.}
  \vspace{1mm}
  
 {\it Proof.} Let $X_0\subset X $  minimises $vol_{n-1}(\partial X_0)$ among all domains in $X$ which contain $\partial_-$ and observe that,  because of {\color {blue}  $ [\ast_1]$}, the boundary of $X_0$ contains a point $y\in \partial X_0$
 with $dist (y,\partial_+)\geq dist(\partial_-, \partial _+)-d_0 $  and, because of  {\color {blue}  $ [\ast_2]$},
this $X_0$ doesn't intersect $\partial_+$. 

Then, by an elementary  argument (see [G(Plateau-Stein) 2014])  the hypersurface $\partial X_0$ can be smoothed and its mean curvature made everywhere positive. \vspace {1mm}
 
{\color {blue} $[\star\star]$}   {\it Two Words  about  {\color {blue} $ [\ast_2]$.}} There  are several well known  cases of manifolds where the lower bound on the volumes of minimal hypersurfaces  $Y\subset X$, where $\partial Y\subset partial X$
  and where $dist(y, \partial) X\geq R$ for some $y\in Y$, 
   are available.

   For instance  if $X$ is $\lambda$-bi-Lipschitz to the $R$-ball in the simply connected space $X^n_\kappa$
 with   constant curvature $\kappa$, then the volume of $Y$ is bounded from below in terms of  the volume of the $R$-ball $B^{n-1}_0(R)\subset X^{n-1}_\kappa$ as follows.

  Let $g=dr^2 + \phi^2(r)ds^2$, $r\in [0,R]$, be the metric in  the ball $B(R)=B^{n-1}_0(R)\subset X^{n-1}_\kappa$ in the polar coordinates where $ds^2$ is the metric on the unit sphere $S^{n-1}$ and let  $g_\lambda=dr^2 + \phi_\lambda^2(r)ds^2$ be the metric (which is typically singular at $R=0$),    such that  the volumes of the concentric balls and of their boundaries satisfy 
 $$ \frac  {vol_{g_\lambda, n-1}  B(r)}  {vol_{g_\lambda, n-2} (\partial  B(r))}= \Psi_\lambda(r)=\lambda^{2n-3}  \frac  {vol_{g, n-1}  B(r)}  {vol_{g, n-2} (\partial  B(r))}. \leqno {\color {blue} [\star]}$$
 
Then the   standard  relation between  $vol (Y)$ and the filling volume bound in $X$ says that,\vspace {1mm}

\hspace {15mm}{\it the volume of the above $Y$  is bounded by $vol_{g_\lambda, n-1}  (B(R))$.}\footnote{The quickest  way to show this is with a use of {\it Almgren's sharp isoperimetric inequality}. But  since this   still remains unproved for $\kappa<0$,  one needs a slightly  indirect argument in this case,   which, possibly  -- I didn't check it carefully  -- gives a  slightly weaker inequality, namely  $Vol(Y)\geq c_n\cdot vol_{g_\lambda, n-1}  (B(R))$  for some  $c_n>0$.}

\vspace {1mm}

 Notice that {\color {blue} $[\star]$} uniquely and rather explicitly defines the function  $\phi_\lambda$. 

In fact, since 
 $$vol_{g_\lambda, n-2} (\partial  B(r))= \phi_\lambda ^{n-2} \sigma_{n-2}$$
  for  
 $\sigma_{n-2}=vol(S^{n-2})$, and since 
   $$\frac {d vol_{g_\lambda, n-1}  (B(r)}{dr}= vol_{g_\lambda, n-2} (\partial  B(r))$$
 this {\color {blue} $[\star]$}   can be written as the following  differential equation on $\phi_\lambda$
$$ \phi_\lambda^{n-2} = \frac {d(\phi_\lambda^{n-2}\Psi_\lambda)}{dr},$$
 where our $\phi_\lambda$ satisfies  $\phi_\lambda(0)=0$.
   \vspace {1mm}.

\vspace {1mm}

\hspace {35mm} {\large  \color {blue} \sf Examples of Corollaries.} \vspace {1mm}

{ \color {blue} \textbf A.} {\sf Let $X$ be a complete Riemannian $n$-manifold with {\it infinite $(n-1)$-volume at infinity}, which means that the boundaries of compact domains which exhaust $X$,
$$U_1\subset U_2 \subset...\subset U_i\subset... \subset  X,$$
have $vol_{n-1}(U_i)\to \infty$. }

{\sl If $X$ contains no complete non-compact  minimal hypersurface with finite $(n-1)$-volume,
then $X$ can be exhausted by compact smooth domains the boundaries of which have positive mean curvatures.}

Notice that according to {\color {blue} $[\star\star]$}, 

{\sl no such minimal hypersurface exists in manifolds with uniformly bounded, or even, slowly growing, local geometries.}

 Also notice that
 
 {\sf  infinite non-virtually cyclic coverings $\tilde X$ of  {\it compact} Riemannian manifolds $X$, besides having {\it uniformly bounded} local geometries, also have {\it infinite $(n-1)$-volumes at infinity}; hence they
 can be exhausted by compact smooth mean convex domains.}

And even  the {\it virtually cyclic coverings} $\tilde X$ admit such exhaustions unless they are isometric cylinders $Y\times \mathbb R$.

Also notice that  if $\tilde X$ is a  Galois (e.g. universal) covering with {\it non-amenable}   deck transformation (Galois)  group,
then it can be exhausted by $U_i$ with 

\hspace {-6mm}$mean.curv(\partial U_i)\geq \varepsilon >0$.
(See 1.5(C) in [G(Plateu-Stein) 2014].)
\vspace {1mm}

{\it Exercises.} (a) Show that 
 if a complete connected non-compact   Riemannian $n$-manifold  $X$ has  uniformly bounded local geometry,    then $ X\times \mathbb R$ has infinite $n$-volume at infinity.
 
 (b) Show that if   $X$ has $Ricci(X)>-(n-1)$, then 
 $X\times \mathbf H^2_{-1}$ has infinite $(n+1)$-volume at infinity and that it can be exhausted by compact smooth mean convex domains.

\vspace {1mm}

{ \color {blue} \textbf B.}  
 {\sf Let $A$ be  $\lambda$-bi-Lipschitz  to the annulus $\underline A=\underline  A(r, r+R)$   between two concentric spheres of radii $r$ and $r+R$  in the Euclidean space $\mathbb R^n$.
\footnote{This means the existence of a $\lambda$-Lipschitz homeomorphism from $\underline A$ onto $A$,  the inverse of which $A\to \underline A$
is also $\lambda$-Lipschitz.}}

{\sl If $R\geq 100 \lambda r$, then $A$ contains a hypersurface $Y$ which separates the two 
boundary components of $A$ and  such that
$$mean.curv (Y)\geq \frac {100}{r}.$$}

{ \color {blue} \textbf C.} {\sf Let $\underline X$ be  a complete simply connected $n$-dimensional manifold with  non-positive   sectional curvature and such that $Ricci(X) \leq -(n-1)$, e.g.  an irreducible   symmetric space  with $Sc(X)= -n(n-1)$.

 Let $A$   be a compact Riemannian manifold which is  $\lambda$-bi-Lipschitz   to the annulus between
two concentric balls  $B(r)$ and  $B(r+R) $ in $ \underline X$.}

  {\it There exists a (large) constant $const_n>0$,  such that  if $R\geq const_n \cdot \log \lambda$, then  there exists a smooth closed  hypersurface  $Y\subset A$, which separates the two boundary components in $A$ and  such that
$$mean.curv(Y)\geq \frac {n-1}{\lambda +const_n(\lambda-1)}.\footnote {The sign convention for the mean curvature is such that the {\it mean convex} part of $V$ bounded by $Y$  is the one which contains the boundary component {\it corresponding to the sphere  $\partial B(r)$} in $\underline X$.}$$ }

{\it About the Proof.} If $\kappa(X) \leq -1$ this follows from {\color {blue}$[\ast \ast]$}, while  the general case 
needs a minor generalization of this.

 \vspace {1mm}

{\it \color {blue!40!black}\large  First Application to Scalar Curvature.}  Since $$Rad_{S^{n-1} }(Y)\geq \lambda^{-1}  Rad_{S^{n-1}}(\partial B(r)) \gtrapprox \exp r,$$ 
the above inequality 
 together with Remark (b) after  {\Large \color {blue}$\Circle^{n-1}$} from section \ref {mean4}. yields the following. \vspace{1mm}

 {\sl If  a Riemannian manifold $X$ is $\lambda$-bi-Lipschitz to the ball  $B(R)\subset \underline X$, where   $R\geq const_n \log \lambda$, then 
 the  scalar curvature of  $X$ is bounded by:
$$\inf _{x\in X}Sc(X,x)\leq - \frac {1}{const_n\cdot \lambda^2}. $$ }

\vspace {1mm}

{\it \color {blue!40!black}\large  Second Application to Scalar Curvature.} It may happen that a manifold $X$ with $Sc(X)>0$ itself contains no mean convex domain, but it may acquire such domains after  a   modification of its metric that doesn't change  the sign of the scalar curvature. Below is an instance 
of this.

 Let $X =(X, g)$ be a compact $n$-dimensional Riemannian band, as in the above 
{\it \color {blue} \textbf {Lemma}}, where   the boundary of a compact Riemannian manifold $X=(X,g)$ with $Sc(X)\geq 0$  is decomposed as earlier,
$\partial  X=\partial_ -\cup \partial_+$.

Let $Sc(X)>0 $ and let us  indicate possible modifications of the Riemannian  metric $g$, that would enforce the conditions  {\color{blue} $[\ast_1]$}  and      {\color{blue} $[\ast_2]$}  in the {\it \color {blue} \textbf {Lemma}}, while keeping the scalar curvature positive.

We will  show below that  this  can be  achieved in some  cases by  multiplying  $g$ by a positive function $e=e(x)$ , which is equal one near $\partial_-\subset X$ and  which is  as  large far  from  $\partial_-$  as is needed for {\color {blue} $[\ast_1]$}   and where 
 we  also  need  the Laplacian of  $e(x)$ to be bounded from above by $\varepsilon_n Sc(X,x)$  in order to keep $Sc>0$  in agreement with the Kazdan-Warner conformal change  formula from section \ref {conformal2}. 

The simplest case, where there is no need for any particular formula, is where
the sectional curvatures of $X$ are pinched between $\mp b^2$, no geodesic  loop in $X$ of 
length<$\frac{1}{b}$  exists, while   the scalar curvature of $X$ is bounded from below by $\sigma>0$. 

In this case, let $$e_0(x)= c\frac {\sqrt \sigma}{b+1} dist_g(x, \partial _-0)$$
and observe that if $c=c_n>0$ is sufficiently small, then  $e_0(x)$ has a {\it small} (generalized)  gradient $\nabla(e_0)$ and,  because the the geometry of $X$ is suitably bounded, the function $e_0$ can be   approximated by a smooth function 
$e(x)$ with second derivatives significantly smaller than $\sigma$,

\hspace {20mm} thus, {\it ensuring the inequality $Sc(eg)>0$.}

On the other hand, if
$$dist (\partial_-, \partial_+)\geq C (b+1) ||\nabla (e)||^{-1} vol(\partial _-)^{\frac {1}{n-1}},$$
for a large $C=C_n$, 

 then \vspace {0.6mm}
 
\hspace {6mm} the {\it condition {\color {blue} $[\ast_1]$} is satisfied}, say  with $d_0=\frac {1}{2}dist(\partial_-, \partial_+),$  \vspace {0.6mm}

\hspace {-6mm}and, due to the bound on the geometry of $X$,   \vspace {0.6mm}

\hspace {6mm} the {\it condition  {\color {blue} $[\ast_2]$} is satisfied} as well.\vspace {1mm}

Now let us  look  closer at what kind  $e(x) $  we need and observe the following

{\color {blue} [1]} The  bound on the geometry of $X$ is needed only, where the gradient of $e$ doesn't vanish.

Thus, it suffices to have the geometry of $X$  

{\it bounded only in the $\frac{1}{b}$-neighbourhoods of the boundaries of domains $U_i$},
$$\partial_-\subset U_1\subset ... \subset U_i\subset... \subset U_k\subset X,$$
where $dist (U_i, \partial U_{i+1})\geq \frac {1}{b}$ and where  $\frac {k}{b}$ is sufficiently large.

{\color {blue} [2]}  Since, the by the standard comparison theorem(s),

{\sf Laplacians  of the distance-like functions are bounded from above  in terms of the Ricci 

curvature,} \vspace {1mm}

{\it the $b$-bound on the full local geometry  can be replaced by $Ricci(X,x) \geq -b^2g$.}\vspace {1mm}

Summing up, this yields the following refinement of  {\color {blue} B} in 
{\small \color {blue} \FiveStarOpenCircled}  from the previous section. \vspace {1mm}

{\sf Let $X =(X,g) $ be a, {\it possibly non-complete}   Riemannian $n$-manifold, 
such that 
$$Sc(X)\geq 0,$$ 
and let 
$$f :X\to S^n$$
 be an {\it area non-increasing map}, such that the support of the differential of $f$  is compact and the scalar curvature of $X$ in this support is {\sf bounded from below by that of $S^n$}, 
$$\inf_{x\in supp(df)}Sc(X,x)\geq n(n-1).$$

Let $A_i$  
 be disjoint   "bands" in  $X$, that are {\it $a_i$-neighbourhoods of the boundaries} of compact  domains $U_i$, such that  
$$supp(df)\subset U_1 \subset ...\subset U_i ... \subset U_k\subset X.$$} 

Let us give an effective   {\it criterion for vanishing of the degree of the map} $f$ in terms of the geometries of  $A_i$.

{\it \color {blue} \textbf {Proposition}}. {\sf Let    {\sf  the scalar and the Ricci curvatures of  $X$ in $A_i$ for 
$i= 2,... k-1$  be bounded from below} by 

$$  Sc(A_i) \geq \sigma_i \mbox { and } Ricci(A_i) \geq -b^2 g, \mbox {  } 2\leq i\leq k-1,   $$
and set 
$$\beta_i=\frac {\sqrt \sigma_i}{b_i}.$$

Let the {\it sectional curvatures} of $U_k$ outside $U_{k-1}$  be {\sf bounded from above} by 
$$\kappa(U_k\setminus U_{k-1} )\leq c^2,\mbox { } c>0,$$
 and let the complement   $U_k\setminus U_{k-1} $ contains {\it no geodesic loop  of length} $\leq \frac {1}{c}.$}

{\it If the following weighted sum  of $a_i$} {\sf (that are half-widths of the bands $A_i$)} {\it  is sufficiently large,
$$ \sum_{1<i<k} \beta_ia_i   \geq const_n \frac {(vol_{n-1} (\partial U_1))^{\frac {1}{n-1}}}{\frac {a_k}{c} },$$

and if $X$ is orientable spin, then 
$$ deg(f)=0.$$}

{\it Proof.} Arguing as above, one finds a smooth  function $e(x)$, the differential of which is supported in the union of $A_i$, $1<i<k$,  such that $Sc(e\cdot g)$ remains nonnegative (and even can be  easily made everywhere  positive)  and such that  $U_k$ satisfy the assumptions  {\color {blue} $[\ast_1]$}  and  {\color {blue} $[\ast_2]$} of the above {\it \color {blue} \textbf {Lemma}}, that yields a  subdomain  
$$X_{\hspace {-0.7mm}\partialvartoiint} \subset  U_k,$$ which is mean convex with respect to the metric $eg$ and  to a smoothed double of which compact
Llarull's theorem applies.

{\it Remarks.}   (a)  Even in the case of {\it complete} manifolds $X$,  this  doesn't  (seem to)  directly follow from Llarull's theorem, since the latter, unlike the former, needs {\it uniformly positive} scalar curvature at infinity. 

(b) The above  proposition, as well  construction of mean-convex hypersurfaces in general, doesn't  advance,  at least not directly, the solution of the   {\it \color {blue}\large   spectral  $\mathcal D^2$-problem} formulated in section \ref {boundaries4}.
 
%%%%%%%%%%%%%%%%%%%%%%%%%%%%

%%%%%%%%%%%%%%%%%%%

{\sf Let  $X =(X, g)$ be a complete Riemannian $n$-manifold, let  $f:X\to S^n$ be a smooth {\it area contracting} map  the differential $df$  of which has {\it compact} support.

Let 
$$  |d|=\sup_{x\in X}  ||df(x)||$$ and $$ r=r(x)=dist(x, supp(df)).$$

Let the Ricci curvature of $X$ outside $supp(df)$ be bounded from below by 
$$Ricci(x)\geq -b(r(x))^2g(x)$$ 
for some continuous  function $ b(r)$, $r\geq 0$.}\vspace{1mm}

{\sl If the function $b(r)$ grows  sufficiently slowly for $r\to \infty, $ e.g. $\sigma(r) \leq \sqrt [3]r$ for large $r$,  
then there is an effective lower  bound 
$$Sc(X,x)\geq \sigma(r(x)),$$ 
which implies that }

\hspace {30mm}  {\it the map $f$ has zero degree, \vspace {1mm}

  where  $\sigma(r)$,  $r\geq 0$, is  a certain "universal"  function,  which is "small negative" 
at infinity.}

\vspace {1mm}

More precisely, there exists a universal effectively computable family of    
functions in $r$, {
$$\mbox {$\sigma(r)=\sigma_{b, |d|,  N, }(r)$, $r\geq 0$, $N=1,2,....$,}$$ 
with the following five properties 
  
 \hspace {-7mm} (i) \hspace {6mm} the functions  $\sigma(r)$  are {\it monotone decreasing} in $r\geq 0$,

 \hspace {-7mm} (ii) \hspace {5mm}  $\sigma_{b, |d|,  N, }(r)$ is  {\it monotone decreasing} in $N$,}

 \hspace {-7mm} (iii)   \hspace {4mm}  {\it $\sigma_{b, |d|,  N, }(r)$ is {\it monotone increasing} in $b$ and in $|d|$},

$$\mbox {$\sigma(0)=N(N-1)$, while     $\sigma(r)\underset {r\to \infty}  \to -\infty$}\leqno {(\rm iv)}$$

$$\mbox { $\sigma_N(r)=\sigma_{b, |d|,  N, }(r)\underset {N\to \infty}  \to -\infty$ for fixed $b$, |d| and $r>0$,}\leqno {(\rm v)}$$
  {\large \color {blue} such that }

 \vspace {1mm}

{ \color {blue} [$ \bigtimes \Circle^{N-n}] $}    \hspace {1mm}  {\it if   $Sc(X, x)\geq \sigma_{b, |d|,  N, }(r(x))$  for all $x\in X$  and some $N\geq n+2$, 
then, assuming $X$ is orientable and spin, the degree of $f$ is  zero.\footnote{Compare with "inflating  balloon" used in 7.36 of [GL(complete) 1983].}}

\vspace {1mm}

{\it Proof.}  The     bound on $\Delta  \varphi(x)$ for $Ricci \geq -b^2$ (compare with {\color {blue} [2]} from the previous section)  shows that there exists $\sigma_{b, |d|,  N, }(r)$ with the above properties 
(i)-(v) and a  positive function $\varphi(x) $ on $X$, such that

(a) {\sf  $\varphi$ is  equal to $|d|$ on the support $supp(df)\subset X$ 

and such that 
$$\sigma(r(x))+\frac {m(m-1)}{ \varphi(x)^2}    -\frac  {m(m-1)}{ \varphi^2(x )}
 || \nabla  \varphi(x)||^2-\frac {2m}{\varphi(x)}\Delta  \varphi(x)\geq \varepsilon>0\mbox {  for }  r(x)>0.\leqno {(b)}$$

Therefore, by the formula
{(\Large ${\color {blue}\star\star}$})}
  from section \ref {warped+2}  for the scalar curvature of the warped product metrics $g_\varphi =g+\varphi^2ds^2$
  on $X\times S^m$, $m=N-n$, 
$$Sc(g_\varphi )(x,s) =Sc(g)(x)+\frac {m(m-1)}{ \varphi(x)^2}    -\frac  {m(m-1)}{ \varphi^2(x )}
 || \nabla  \varphi(x)||^2-\frac {2m}{\varphi(x)}\Delta  \varphi(x),$$
  the metric $g_\varphi$ has uniformly positive scalar curvature and because of $(a)$ the map $f:X\to S^n$ suspends to an area decreasing map $(X\times S^m, g_\varphi)\to S^{n+m}$  of the same degree as $f$.
Then Llarull's theorem applies and the proof follows.\vspace {1mm}

{\it On Manifolds with Boundaries. } If $X$ is a compact manifold with a  boundary, the above can be applies to the smoothed double $X\cup_{\partial X}X$ , where the scalar curvature of such a double near the smoothed boundary can be bounded from below by the geometry of $X$ near the boundary  and the (mean) curvature of the boundary  $\partial X\subset X$. 

Thus, the above yields a  condition for $deg(f)=0$ in terms of the lower bound on $Sc(X,x)$ and on 
$dist (x, supp(df))$, which is  similar to, yet is different from such a condition from the previous section.

\vspace {1mm}

{\it \color {blue}   Dirac operators with Potentials.}  The recent 

\hspace {10mm}{\it  relative index theorem for the Dirac operators  with  potentials}

\hspace {10mm} by   Weiping Zhang\footnote {See [Zhang(area decreasing) 2020],
[Zhang(deformed Dirac) 2021].}, 

\hspace {-6mm} which applies to complete manifolds $X$   with non-negative scalar curvatures at infinity and which is more efficient in many (all?) cases than multiplication of $X$ by spheres, makes most (all?) of the above redundant.
 
\vspace {10mm}

%%%%%%%%%%%%%%%%%%%%%%%

\subsubsection {\color {blue} Amenable Boundaries} \label {amenable4}
%%%%%%%%%%%%%%%%%
 
  {\sf  If the {\it volume of the  boundary} of   a  compact manifold $X$ is {\it significantly smaller than the volume of $X$} and if it is  additionally supposed that  the manifold is {\it not very much curved near the boundary},  then  we shall see in this section that   } \vspace{ 1mm}

{\color {blue!49!black} \sl the  index theorem  applied to  the double of such an  $X$  with a  smoothed metric,  yield   geometric bounds on the area-wise size of    $X$ in terms of the lower bound on the scalar curvature of $X$.}

\vspace{ 1mm}

{\it Elliptic Preliminaries.} Let $V$ be a  (possibly non-compact)  Riemannian manifold with a boundary,  and let $l$ be a 
section of a bundle $L\to V$
with a unitary connection $\nabla$, such that  $l$ satisfy the following (elliptic) {\it G\aa rding $(\delta_\circ, C_\circ)$-inequality}:
{\sf the $C^1$-norm of  $l$   
  at $v\in V$   is  bounded at  by the 
$L_2$-norm  of $l$ in the $\delta_\circ$-ball $B=B_v(\delta_\circ)\subset V$ as follows
$$ ||l(v)|| + ||\nabla l(v)||\leq C_\circ \sqrt{\int_B||(l)||^2dv}$$ 
for all points $v\in V$, where $$dist (v, \partial V)\geq \delta_\circ.$$

Let $$\rho(v)=dist (v, \partial V)\mbox { and }  \beta=\sup_{v\in V} vol(B_v(\delta_\circ))$$}

{\it Lemma.}  {\sl  If $l$ vanishes on an $\varepsilon$-net $Z\subset V$, then 
$$ || l(v)|| + ||\nabla l(v)||\leq  \left( {10C_\circ \varepsilon\beta}\right)^{\rho(x)-2\delta_\circ}\sqrt{\int _V l^2(v)dv }$$
Moreover, if $V$ can be covered by $2\delta_\circ$-balls with the multiplicity of the covering at most $m$, then
the $L_2$-norms of $l$ and $\nabla l$  on the subset $V_{-\rho}\subset V$ of the points  $\rho$-far from the boundary, that is 
$$V_{-\rho}=  V\setminus U_{\rho}(\partial V)=\{v\in V\}_{dist(v,\partial  V)\geq \rho},$$ 
satisfies
$$\sqrt{\int _{V_{-\rho}}|| l||^2(v)dv}
\leq \epsilon\sqrt{\int _{V}|| l||^2(v)dv}\leqno {\mbox {\color{blue} \rm[\NibLeft]}}$$ 
for $\epsilon= m  \left( {10C_\circ \varepsilon\beta}\right)^{\rho(x)-2\delta_\circ}$.

{\it Proof.} Combine  {\it G\aa rding's inequality} with the following obvious one:
$$ ||l||\leq \varepsilon ||\nabla||l$$
and iterate the resulting inequality 
$i$ times insofar as  $\rho-i\delta_\circ$ remains positive.}
\vspace{1mm}

{\it Remark.} A  single  round of iterations suffices   for  our immediate applications.    

\vspace{1mm} 

{\it\color {blue} Corollary.} {\sf Let $X$ be a complete orientable  Riemannian manifold  of dimension $n$ with compact boundary (e.g. $X$ is compact or homeomorphic to $X_0\times \mathbb R_+$, where $X_0$ is a closed manifold), 
and let, y
{\it for some $\rho>0$ and $0<\delta_\circ<\frac {1}{4}\rho$, \vspace {1mm}

the $\rho$-neighbourhood  of the boundary of $ X$, denoted 
$U=U_\rho(\partial X)\subset X$,   has }}

{\it (local) 
geometry bounded by}  $\frac {1}{\delta_\circ}$,   

{\rm \hspace {-6mm}where we succumb to  tradition and  define this bound 
on geometry 
  as follows}: \vspace {1mm}

{\sf \hspace {-2mm}  the sectional  curvatures $\kappa$  of $U$ are  pinched between $-\frac {1}{\delta^2_\circ}$ and 
 $\frac {1}{\delta^2_\circ}$ and the injectivity radii are bounded from below by  ${\delta_\circ}$ at all points $x\in U$, for which  $dist(x,\partial X)\geq \delta_\circ $, that is, in formulas, 
$$\mbox { $|\kappa(X,x)|\leq \frac {1}{\delta_\circ^2}$ for $dist(x, \partial X)\leq \rho $  and $injrad (X,x)\geq \delta_\circ   $   for $ \delta_\circ  \leq dist (x, \partial X)\leq \rho$}.$$

Let the scalar curvature of $X$   be non-negative $\frac {1}{2}\rho$-away from the boundary,
$$Sc(X, x)\geq 0\mbox {  for  } dist(x, \partial X)\geq \frac {1}{2}\rho.$$

Let $f:X\to S^n(R)$, where $S^n(R)$ is the sphere of radius $R$,   be  a  smooth  {\it area decreasing}  map ,  which is constant on $U_\rho$, and, if $X$ is non-compact,  also locally  constant at infinity.

 Let the degree of this map  be {\it bounded from below by the volume of $U_\rho=U_\rho(\partial X)$} as follows.
$$d> C vol(U_\rho)  \mbox { for some } C\geq 0. $$

 If $\delta_\circ$,  $\rho$ and $C$ are {\it sufficiently large},  then, provided  $X$ is {\it  \color {blue!59!black}spin, }  the  scalar curvature of 
the complement 
$$X_{-\rho} =X\setminus U_\rho=\{x\in X\}_{dist(x,\partial X>\rho}$$ 
{\it can't be everywhere  much greater than}  $Sc(S^n(R))= \frac {n(n-1)}{R^2}$. 
Namely
$$\inf_{x \in X_{-\rho}}Sc(X,x)\leq \sigma_+\frac {n(n-1)}{R^2}+ 
\sigma,\leqno{ \mbox {\color {blue}[\NibSolidLeft]}}$$
where $\sigma= \sigma_n(\delta_\circ, \rho, C)$ is  a  positive  function, which may be {\sf infinite} for {\sf small}   $\delta_\circ$ and/or  $\rho$ and/or $C$ and which has the following properties.\vspace{1mm}

 $\bullet $ {\it the function $\sigma$  is {\it monotone decreasing  in}  $\delta_\circ$,  $\rho$ and $C$;
 
 $\bullet $     $\sigma_n(\delta_\circ, \rho, C)\to 0$  for $C\to \infty$ and arbitrarily fixed $\delta_\circ>0$ and 
 $\rho>\delta_\circ$.}}
 
 \vspace{1mm}

{\it Proof.} Let $2X=$\DD$X$ be a smoothed double of $X$ and $L\to 2X$ 
the vector bundle induced from $\mathbb S^+(S^n)$ by $f$ applied to a copy 
(both copies, if you wish)  of $X\subset 2X$.

Assume $n=dim(X)$ is even, apply the  index theorem and conclude that  the dimension of the space of $L$-twisted harmonic spinors on $2X$  is  $\geq d$.

Therefore, there exists such a non-zero spinor $l$ that vanishes at given $d-1$ points in $2X$.  

Let  such points 
make  a  $ \varepsilon$-net on the subset  $2U_{\rho_\circ}=$ \DD$U_{\rho_\circ}\subset  2X$ with a minimal possible $\varepsilon$.

If $d$ is much  larger then  $vol  (2U_\rho)\approx 2 vol(U_\rho)$, then this $\varepsilon$ becomes small and, consequently,
$\epsilon$ in the above inequality   {\color{blue} \rm[\NibLeft]} also becomes small.  Then, the inequality  {\color{blue} \rm[\NibLeft]} applied to  the domain $2U_\rho\subset 2X$, shows that the integral
$$\int_{2U_\rho} ||l||^2(x)dx$$  is much smaller then the integral of $||l||^2$ over the   complement 
$2X_0=2X\setminus 2U_\rho$. 

Therefore, if $\sigma_+$ is large then the sign of the full  integral 
$$ \int_{2X} Sc(X,x) ||l||^2(x)dx= \int_{2X_{\rho}} Sc(X,x)|| l||^2(x)dx +   \int_{U_\rho} Sc(X,x) ||l||^2(x)dx$$
is equal to the sign  of  $\int_{2X_{\rho}} Sc(X,x)|| l||^2(x)dx$, which contradicts the Schroedinger-Lichnerowicz-Weitzenboeck formula for harmonic $l$. 

Thus, modulo simple verifications and evaluations of constants left to the reader, the proof is completed. 
\vspace{1mm}

{ \color {blue} {\it Example} 1}.  {\sf Let    a   complete non-compact  orientable spin Riemannian  $n$-manifold $X$ with {\it compact boundary} admits smooth {\it area decreasing} maps 
$f_i: X\to S^n$  of {\it non-zero  degrees},\footnote{Here as everywhere in this paper, when you you speak of $deg(f)$
the map $f$ is supposed to be locally constant at infinity as well 
as on the boundary of $X$.}  such that the  "supports" of $f_i$,   i.e.  the  subsets  where these maps are {\it  non-constant}, may lie  arbitrarily far   from the boundary of $X$,
$$\mbox {  dist ("supp"$f_i, \partial X) \to \infty$ for $i\to \infty$}.$$}
 
 {\it \color {blue!50!black} Then the scalar curvature of $X$ {\it can't be uniformly positive at infinity}:
$$ \liminf_{x\to \infty}Sc(X,x)\leq 0.$$}

Moreover, the same conclusion holds, if 

{\sf there exist $i$-sheeted  coverings $\tilde X_i\to X$, which admit  smooth area decreasing maps  
$f_i: \tilde X_i\to S^n$, such that 
$$ \frac  {deg(f_i)}{i}\to \infty \mbox  { for } i\to \infty.$$}

\vspace {1mm}

{\it Example} 2. Let $ Y_k$ be a $k$-sheeted covering of the unit 2-sphere $S^2=S^2(1)$  minus two opposite balls of radii $\frac {1}{k^m}$, for some $m\geq 1$.

Then  the product manifold $X_0=Y_k\times S^{n-2}(k)$ admits an area decreasing map $f:X_0\to S^n(R)$ constant on the boundary  and  such that  
$$deg(f)\geq \frac {k}{10d}$$ 
and it follows from  the above corollary that the  Riemannian metric on $X_0$ can't be extended    to a larger manifold $X\supset X_0$, with bounded geometry  and $Sc\geq 0$ without adding much  volume to $X_0$, say in the case 
$m=n-1$, although  $vol_{n-1}(\partial X_0)$ remains bounded for $R\to \infty$.

\vspace {1mm}

{\it Melancholic Remarks.} Rather than indicating the richness of the field, the diversity of the results  in the above sections  {4.6.1}- {4.6.4}  is due to our inability to formulate and to prove  the true general theorem(s).

%%%%%%%%%%%%%%%\%%%%%%%
%%%%%%%%%%%%%
%%%%%%%%%%%%%%%%%%%%%%%%%
\subsubsection {\color {blue} Almost Harmonic Spinors on Locally Homogeneous and  and Quasi-homogeneous  Manifolds with Boundaries}\label {4.6.5}%%%%%%%%%%%%%

%%%%%%
%%%%%%%%%%%%%%%%%%%%

Let $X$ be a complete Riemannian manifold with a transitive isometric action of a group $G$, let $L\to X$ be a vector bundle with a unitary connection $\nabla$  and let the action of $G$ equivariantly lifts to an action on $(L, \nabla)$.

Let the $L_2$-index  of the twisted Dirac operator $\mathcal D_{\otimes L}$ (see [Atiyah($L_2$) and  
[Connes-Moscovici($L_2-index$ for homogeneous)  1982],    be non zero. For instance, if $X$ admits a free discrete isometry group $\Gamma \subset G $ with compact quotient, then this is equivalent to this index to be non-zero on $X/\Gamma$.\vspace{1mm}

{\sf The main class   of examples of such $X$ are {\it symmetric spaces with non-vanishing "local Euler characteristics}  (compare with  [AtiyahSch(discrete
series) 1977]) i.e. where the corresponding (G-Invariant) $n$-forms, $n=dim(X)$ don't vanish.  

The simplest instances of these are   hyperbolic spaces $\mathbf H ^{2m}_{-1}$,  where the indices  of the Dirac operators twisted with  the positive spinor  bundles don't vanish.
 In  fact, such an  index for  a  compact quotient   manifold $\mathbf H ^{2m}_{-1}/\Gamma$  is equal to  $\pm$one half of the Euler characteristics of this manifold by the Atiyah-Singer formula  (compare [Min(K-Area) 2002]).}

\vspace{1mm}

Let $(X,L)$ be an above  homogeneous pair with  $ind(\mathcal D_{\otimes L})\neq 0$  and let 
$X_R\subset X$  be a ball of radius $R$. Then the restrictions of  $L_2$-spinors on $X$ (delivered by the $L_2$-index theorem) to $X_R$ can be perturbed (by taking products with  slowly decaying cut-off functions) to {\it $\varepsilon$-harmonic} spinors that {\it vanish on the boundary} of $X_R$, where 
$\varepsilon\to 0$   for $R\to\infty$ and where
"$\varepsilon$-harmonic" means that %%%
$$    \int _{X_R}\langle  \mathcal D_{\otimes L}^2(s), s\rangle  \leq \varepsilon^2\int_{X _R}||s||^2dx$$ 
as in  {\color {blue}\large c \ScissorHollowRight\hspace{-0.2mm}$_{\lambda}$} in section \ref {boundaries4}.

In fact, it follows from  the  local proof of the $L_2$-index theorem in [Atiyah($L_2$)  1976] or, even better,  from   its later version(s) relying on the finite propagation speed, that these $\varepsilon$-harmonic spinors can be constructed internally in  $X_R$ with no reference to
the ambient $X\supset X_R$.

Moreover, a trivial  perturbation (continuity) argument shows that 

{\it similar spinors exist on manifolds $X'_R$  with these metrics close to these on $X_R$.} 

but it is unclear "how close" they should be. Here is a specific problem of this kind.\vspace {1mm}

{\sf Let $X_R$ be a compact Riemannian   spin manifold with a boundary, such that 
$$\sup_{x\in X}dist(x, \partial X_R)\geq R$$ 
and let  
 the sectional curvatures of $X$ are everywhere  pinched between 
$-1$ and $-1-\delta$. }
\vspace {1mm}

{\color {blue} (A)}   {\sl \color {blue!34!black}Under what conditions on $R,\delta$ and $\varepsilon$ does $X_R$ support a non-vanishing $\varepsilon$-harmonic spinor twisted with
the  spin bundle $\mathbb S(X_R)$?} 

\vspace {1mm}

Besides, one wishes to have 

{\color {blue} (B) } {\sf  similar spinors on manifolds $\overline X$ mapped to $X_R$ with {\it non-zero} degrees and with 

 {\it controlled  metric  distorsions}}  

 \hspace {-0mm}in order to get {\sl bounds on the scalar curvatures} of such $\overline X$
 
  (See section \ref {spin harmonic6})for  continuation  of this discussion to {\it fibrations} with  quasi-homogeneous fibers.)

%%%%%%%%%%%%%%%%%%%%%%%%

\subsection {\color {blue}Topological Obstructions to Complete Metrics with  
 Positive Scalar Curvatures Issuing from the Index Theorems  for Dirac Operators}\label {obstructions4}
%%%%%%%%%%%%%%%%%%%%%%%%%%%%

Obstruction on {\it homotopy types} of compact manifolds $X$ implied by the  existence of metrics of positive scalar curvature on  $X$, obtained by Dirac theoretic methods 
usually 
(always?)  generalize to {\it non-compact complete} manifolds, where   "homotopy"  means "proper homotopy",
i.e. the maps being  "homotopies"  as well as the maps establishing homotopies   must be proper: {\sf infinity-to-infinity.}

Moreover, such obstructions    not only rule out metrics with positive scalar curvatures on   $n$-manifolds $X'$ which are homotopy equivalent to $X$, but also on $n$-manifolds $\hat X$  that {\it dominate} (the fundamental
homology class of) $X$, i.e. admits maps $f: \hat X \to X$  with  $deg(f)=\pm 1$ to $X$ in the orientable  cases, and often, even with any  $deg(f)\neq 0$.\vspace {1mm}

{\it Dimension$+m$-Domination.}  The above  also applies  to smooth  proper  maps of $(n+m)$-dimensional manifols to $n$-dimensional $X$,
say $f: \hat X^{+m}\to X$, such that  the pullbacks of generic points  under  $f$ and by  all smooth maps  $X^{+m}\to X$  {\it homotopic to} $f$ -- these pullbacks (but not necessarily all $m$-manifolds homotopy equivalent to these pullbacks)  
admit no metrics with $Sc>0$.\vspace {1mm}

 \textbf {Example 1}: {\sl Maps of  non-zero  $\hat A$-degree to Enlargeable\footnote {A compact Riemannian  $n$-manifold  $X$ is enlargeable if it admits (finite or infinite)  coverings $\tilde X$  with {\it arbitrarily large} hyperspherical radii, i.e. for all $R>0$, there exists a covering  $\tilde X$,  which admits a locally constant at infinity distance decreasing map $\tilde X\to  S^n$ with non-zero degree.
 
 Notice that this condition doesn't depend on the Riemannian metric in $X$, moreover it is a homotopy (even domination) invariant. }
 Manifolds   and Similar Maps.}  If
 a compact {\it spin}  $(n+m)$-manifolds  $\hat X^{+m}$  admits a  smooth map $f$  to  compact  {\it enlargeable} $n$-manifolds $X$, (see section \ref{filling+hyperspherical+asphericity3}   e.g. to the   torus
  $\mathbb T^{n}$, or, more generally, to a 
 Riemannian manifold with non-positive sectional curvature,  such that the pullback $f^{-1}(x)\subset \hat X^{+m}$ of a generic point $x\in X$ has {\it non-zero} $\hat\alpha$-invariants, e.g.  $\hat A(f^{-1}(x))\neq 0$ in the case 
 $m=4k$, then $\hat X^{+m}$  can't carry a metric with $Sc>0$.\vspace {1mm}

{\it About Relevance of  Spin.}  {\color{red!50!black}Probably},  the same non-existence conclusion holds if only 
the pullback  "$\hat X^{+m}$ is spin, for instance, where $\hat X^{+m}$ is diffeomorphic to $X\times X^m$, where $X^m$ (but not necessarily $X$) is spin. 

In fact, if $m+n\leq 8$ this follows from the the  $\mu$-bubble separation theorem in section \ref{separating3},  and if  $m_n\geq 9$, this might follow from Lohkamp's desingularization  results. 
(Schoen-Yau's 2017  theorem  is non-sufficient for this purpose.)

 On the other hand, the  Dirac theoretic method has an advantage of being applicable to {\it  $\wedge^2$-enlargeable} 
 manifolds $X$ defined in example 4  below.

 Also Dirac operators   serve well if  the underlying  $X$ is a  quasisymplectic {\color {blue}$\otimes_{\wedge^k\tilde\omega}$}-manifolds as in section \ref{SY+symplectic2}, e.g. a  closed aspherical $4$-manifolds $X$ with $H^2(X;\mathbb Q)\neq 0$.\footnote {One should note, however, that   no example is known of a {\it
  compact} 
 {\it non-enlargeable} manifold that is 
 $\wedge^2$-enlargeable or  quasisymplectic {\color {blue}$\otimes_{\wedge^k\tilde\omega}$}  or a manifold  with infinite $K$-area.}

 \vspace {1mm}

 {\it "Positive"  versus "Uniformly Positive".} If $X$ is non-compact, one has to distinguish   "just  (strict or not)  {\it positivity}" of the scalar curvature,   $Sc(X)>0$  along with $Sc(X)\geq 0$  -- the existence of the former implies the existence of the latter  except for a few exceptional "rigid"  examples, such  as Riemannian flat manifolds and Ricci flat  K\"ahler (Calabi-Yau)   manifolds, 
  from  "{\it uniform uniform positivity"},  where  
$Sc(X)\geq \sigma>0$.\vspace {1mm} 

 \textbf {Example 2}:  {\it Metrics  with Positive Curvatures in  the Plane  and their  High Dimensional   Warped Descendants.}     The products of tori $\mathbb T^{n-2}$   by the plane $\mathbb R^2$ (obviously) admit metrics with $Sc>0$, but {\it no metric with $Sc\geq \sigma>0,$}  where 
the latter follows from Roe's partitioned index theorem.

Also one can do it with   Zeidler's-Cecchini's  Dirac theoretic $\frac {2\pi}{n}$-inequality for Riemannian   spin bands, while  our   non-Dirac theoretic proof needs Lohkamp-Schoen-Yau desingularization theorem(s)  for $n\geq 9$.

More generally, by the same  token  the product manifolds  $X= X_0 \times  \mathbb R^2$,  support  complete  metrics with $Sc>0$,  
but if $X_0$ admits {\it no domination} by a manifold with a complete  metric with {\it positive} scalar curvature, then $X$ admits {\it no domination}  by a manifold with a complete  metric with  with {\it uniformly positive} scalar curvature.\footnote{This follows from the $\mu$-bubble separation theorem (section \ref{separating3}) that relies on Lohkamp-Schoen-Yau desingularization   for $n\geq 9$ , but  
I am not certain how much of this can be proven for by Dirac theoretic methods in the case of  spin manifolds.}\vspace {1mm}

{\it Exercise.}\footnote {I haven't  solved this exercise.} Show that  products $X_1\times X_2$ of  {\it non-compact} manifolds $X_1$ and $X_2$  admit complete metrics with $Sc>0$, while such triple  products, $X_1\times X_2\times X_3$  admit complete metrics with $Sc\geq \sigma >0$.\vspace {1mm}

\textbf {Example 3:} {\it Simply Connected Manifold Dominated  by $Sc>0$.} There only instance of a Dirac theoretic  obstruction for $Sc>0$ on topology of compact {\it simply connected} manifolds, which  is (this is an accident) a homotopy theoretic one, is  
Lichnerowicz' $\hat A[X]\neq 0$  for $n=dim(X)=4$.  (If $n\geq 5$ there is no  constraints on rational Pontryagin 
classes of $X$ except for the signature and none of higher  $\hat\alpha$-invariants used in Hitchin's theorem is homotopy invariant either.)

But even this obstruction is not "domination invariant":   connected sums $X_{\#\pm}=X\#-X$ have $\hat A[X\#-X]=0$ 
for all $X$  and, by Milnor' homotopy classification theorem, these  $X_{\#\pm}$ are homotopy equivalent to manifolds which admits metrics with $Sc>0$, namely to  connected sums  of  $CP^2$  and $S^2\times S^2$ by Milnor's 1958 theorem  and 
by adding more copies of $S^2\times S^2$ these become {\it diffeomorphic} to 
 connected sums  of  $CP^2$  and $S^2\times S^2$ by Wall's 1964 theorem.\vspace {1mm}

\vspace {1mm}

All (known) Dirac  theoretic   non-domination results  of compact $n$-manifolds $X$ by compact   $\hat X$ with $Sc(\hat X)>0$ apply only  to spin manifolds  $\hat X$\footnote {In all known examples it suffices  that  the universal  covering of $\hat X$  is spin.}  and rely on existence of flat or almost flat (generalized,  e.g. virtual Fredholm) unitary  vector  bundles 
over $X$ (or over $X\times \mathbb T^1$)  with non-zero Chern numbers. 

In fact, the limit of applicability of such results would be  (essentially) reached if one could resolve the following. 
\vspace{1mm}

{\it \textbf {Problem A}.}  {\sf   Let $B$ be a Riemannian manifold,  let $X\subset  B$ be a compact {\it relatively aspherical}  submanifold, i.e. the inclusion homomorphisms of the higher homotopy groups, $\pi_i(X)\to\pi_i(B)$, 
 vanish  for all $i\geq 2$.}

Prove (or disprove)  that
 {\sl for 
all   complex vector bundles $L\to B$ and all 
 $\varepsilon>0$ there exist  vector  bundles $L_\varepsilon\to X$ with unitary connections, such that

(i) the bundles $L_\varepsilon$ are isomorphic to multiples  of $L$ restricted to $X$
$$L_\varepsilon=k\cdot L_{|X}\mbox { 
for }
k\cdot L=\underset {k}{\underbrace {L\otimes L\otimes ...\otimes L}};$$

(ii) The curvature  operators $R_{L_\varepsilon}$  of $L_\varepsilon$ satisfy
$$||R_{L_\varepsilon}||\leq \varepsilon.$$}

In fact, as we know,  that if $X$ is spin, then the   index theorem applied to the twisted Dirac operators 
$\mathcal D_{\otimes L_\varepsilon}$ (that act on spinors with values in the bundles $L_\varepsilon$)
shows that the (untwisted)  Dirac operators $\cal D$  on certain  covering manifolds $\tilde X_\varepsilon$  contain zero in their spectra; thus    $Sc(X)\ngtr 0$ by the Schroedinger-Lichnerowicz-Weitzenboeck-(Bochner) formula.

In this in mind, one asks   another question.

\textbf B.  {\sf Suppose,  an even dimensional  compact spin  submanifold $X$ in an aspherical space $B$ represents a {\it  non-torsion} homology class in $B$. \footnote {One  has little idea of  what to  expect for non-zero  torsion classes.}}

{\it Does  then the spectrum  of the Dirac operator on  some  covering of $X$ contain zero in the  spectrum?}

\vspace {1mm}

Now, let us look more systematically at   what of the above generalizes to complete  manifolds with $Sc>0$ and with $Sc\geq \sigma>0$.

Originally, the results for  $Sc(X)>0$ were derived from these, where  $Sc\geq \sigma>0$, namely applied to  
$X\times S^2(R)$ for suitably large $R$. 

Nowadays, one has at one's disposal index theorems  for {\it Dirac operators with   potentials} proved in  
 [Cecchini(Callias) 2018], 
[Cecchini(long neck) 2020]  and in [Zhang(area decreasing)  2020].\vspace {1mm}
\vspace {1mm}

\textbf {Example 4}: {\it  $\wedge^2$-Enlargeability against $Sc>0$.}  A  Riemannian metric $g$ on  a manifold $X$  is  $\wedge^2$-enlargeable
if, for all $R>0$,  there exists    coverings $\tilde X$, which  admit locally constant at infinity  $g$-area non-increasing maps with non-zero degrees to the $R$-spheres $S^n(R)$, 
 where, a  priori, such a covering may depend on $R$. 

A smooth manifold $?$ is  $\wedge^2$-enlargeable  if     all  Riemannian metrics on it are enlargeable.

For instance, 

{\it metrics with  infinite areas on connected surface are  $\wedge^2$-enlargeable},

while 

{\it  connected surfaces are enlargeable if they have infinite fundamental groups.}
\vspace {1mm}

{\it Exercises.}  (4a) Show that the products $X=X_1\times X_2$, where both $X_1$ and $   X_2$  are  connected  non-compact 
are   not $\wedge^2$-enlargeable.

(4b) Show that the products of enlargeable manifolds by $\wedge^2$-enlargeable are $\wedge^2$-enlargeable.

(4c) Show that  the product $X=X_0\times\mathbb R$, where $X_0$ is enlargeable, is  
$\wedge^2$-enlargeable.

{\it \color {red!40!black}Probably} the converse is also true: {\sf if  $X_0\times\mathbb R$ is $\wedge^2$-enlargeable, then $X_0$ is enlargeable.} (This is  close  in spirit to   {\color {red!40!black}  stabilization conjecture} in section 7.3)

Also it is   {\it \color {red!70!black} \it not  impossible} that (b) also admits a converse: {\sf if the product $X_1\times X_2$ is $\wedge^2$-enlargeable, then one of the two manifolds is   $\wedge^2$-enlargeable and another one is  enlargeable.} 
 
(4d) Show that if $X$ dominates (a multiple of the fundamental class of) a $\wedge^2$-enlargeable manifold $\underline X$, i.e. if there is a quasi-proper map $f:X\to \underline X$ of non-zero degree,\footnote{A map $f$ is quasi-proper if it extends to a continuous map between the compactified spaces, 
from $ X^{+ends}\supset X$ to $ \underline X^{+ends} \supset \underline X$.} then $X$  $\wedge^2$-enlargeable.

For instance, 
  complements to Cantor (closed zero-dimensional) subsets   in enlargeable manifolds $X$ and  connected sums of  $X$ with arbitrary manifolds are $\wedge^2$-enlargeable.
\vspace {1mm}

\textbf {Theorem 4e.} {\it  $\wedge^2$-Enlargeable manifolds $X$, the universal coverings $\tilde X$  of which are spin, admit no metrics with $Sc>0$.}

This is proven in \S6 in [GL(complete)1983] for spin manifolds $X$  with a use of the relative index theorem
applied to $X\rtimes S^2(R)$, where in  the case of $\tilde X$ spin, one does it   with  relativized 
 Atiyah's  $L_2$-index theorem. 

\vspace {1mm}

\textbf {Example 5}: {\it Obstruction on $Sc>0$ of Complete Metrics  for Manifolds with Infinite relative   $K$-areas and  for Quasisymplectic
 {\color {blue}$\otimes_{\wedge^k\tilde\omega}$}-Manifolds. } 
Let us formulate   two special cases of general non-existence theorems  for  complete metrics with  $Sc>0$ from  [Cecchini-Zeidler(generalized Callias) 2021]and from  [Zhang(deformed Dirac) 2021]  proved   with a use of Dirac operators  with potentials.\footnote {I want to thank  Simone Cecchini  and  Weiping Zhang  for explaining their results to me.} 
\vspace {1mm}

\textbf {Theorem 5a.} {\sf Let $X$ be an orientable     manifold of even dimension $n$ and let 
$X_0\subset X$ be a compact subset, such that $X$ has {\it infinite $K$-area relative to the complement $X\setminus X_0$}.}

This means that  for some,  hence for every, Riemannian metric $g_0$ on $X$ the following holds.

{\sf For all $\varepsilon>0$, there exist complex  vector bundles  $L_1, L_2\to X$ with unitary connections, such that:

$\bullet$
the norms of the curvature operators of these connections  with respect to $g_0$ are  everywhere $\leq \varepsilon$;

$\bullet$  these norms vanish outside  $X_0$, i.e. the connections are flat over $X\setminus X_0$;

$\bullet$ There exists  a parallel,  i.e,   connections preserving, isomorphism between the bundles $L_1$ and $L_2$ over 
 $X\setminus X_0$.

$\bullet$ some (relative)  Chern number $c_I[X]$,  $c_I\in  H^n(X, X\setminus X_0)$, of the virtual bundle $L_1-L_2$ 
doesn't vanish.}

{\it If the universal covering of $X$ is spin, then $X$ admits no complete Riemannian metric $g$  with $Sc(g)>0$.}
\vspace {1mm}

\textbf {Theorem 5b.} {\sf Let $X$ be an orientable     manifold of  dimension $n=2k$ and let 
$X_0\subset X$ be a compact subset.

Let $h\in H^2(X, X\subset X_0)$ be a relative cohomology class,   such that $h^k\neq 0$, while
the   lift of $h$ to the   universal covering of $X$, say $\tilde h\in H^2(\tilde X; \widetilde{X\setminus X_0})$, 
 {\it vanishes.}}

{\it If  the universal  covering $\tilde X$ of $X$ is spin, then $X$ admits no complete  metric with $Sc>0$.}

\vspace {1mm}\vspace {1mm}

\textbf {Example 6.} {\it  Topology at Infinity of  Complete Manifolds with Uniformly Positive Scalar Curvatures. } If instead of $Sc>0$  we want to rule out complete metrics  $Sc\geq \sigma >0$,  we need the above topological conditions on $X$ 
satisfied only at  infinity.
Below is a specific formulation of this.

\vspace {1mm}

\textbf {Theorem/{\color {red!50!black} Conjecture} 6a.}  {\sf  Let $X$ be an orientable   manifold of even dimension $n$,  
 let $X^\circ\subset X$ be an open subset with a {\it compact} complement in $X$ and  let $X_i^\circ\subset  X^\circ$, $i=1,2,...$ ,  be a sequence of compact subsets that tend to infinity in $X$, i.e.
every  compact subset in $X$ intersects only finitely many  $X_i^\circ$.

Let one of the following two conditions be satisfied.

$\bullet_{area}$ The relative  $K$-areas of $X^\circ$ with respect to $X^\circ\setminus X^\circ_i$ are  infinite  for all $i$.

$\bullet_{sympl}$ There exists cohomology  classes  $h_i\in H^2(X^\circ , X^\circ \setminus X^\circ_i)$,  such that $h_i^k\neq 0$, while  
the   lifts  $\tilde h_i\in H^2(\tilde X^\circ; \widetilde{X^\circ\setminus X^\circ_0})$, where $ \tilde X^\circ$ denotes the universal covering of $X^\circ$, 
 vanish.}

{\it If the universal  covering $\tilde X$ of $X$ is spin, then $X$ admits no complete Riemannian metric with
$Sc\geq \sigma>0$.}\vspace {1mm}

{\it The proof of this} must follow from a suitable version of Cecchini's long neck principle and   from [Guo-Xie-Yu(quantitative K-theory) 2020]  but I haven't carefully checked this. Nor am I certain that that  the same conclusion holds  under more general condition(s), where  
the subset $X^\circ$ is not fixed but dependent on (decreasing with)  $i=1,2,...$ .

But  we do  know for sure that  Roe's partitioned  index theorem  shows (in agreement with what follows from the
  $\mu$-bubble separation theorem) that {\sf if a spin manifold $X$ is {\it enlargeable at infinity}}, i.e. 
  \vspace {1mm}
  
  {\sl  if  there exists  an exhaustion of $X$ by 
  compact domains  $X_i\subset X$  with smooth boundaries $Y_i\subset \partial X_i$, such that the {\it complements} of all $X_i$ admit sequences of coverings, say $\tilde X _{ij}^\perp\to X\setminus X_i $, where the hyperspherical radii of the corresponding coverings  $\tilde Y_{ij} =\partial \tilde X _{ij}^\perp$ of $Y_i$  tend to infinity  for $j\to\infty$, }
  
 \hspace {10mm} {\sf  then $X$ admits no  complete metric with $Sc\geq \sigma>0$.}

    \vspace {2 mm}

          \hspace {22 mm}    {\sc Remarks, Problems, {\color {red!50!black} Conjecture}s.}    \vspace {1 mm}

 \textbf {Question 7.}  {\sf Let $X$ be an open aspherical $n$-manifold. 
   Does  {\it non-contractibility of $X$  to the $(n-1)$-dimensional
skeleton $X^{[n-1]}\subset X$} imply that $X$ is $\wedge^2$-enlargeable?}

(It is not even clear, in the case where $X$   admits a complete metric with non-positive sectional curvature,
whether $X$ admits a metric with positive scalar curvature.)
\vspace {1mm}

  \textbf {Question 8.} Are there {\it "topological conditions at infinity"}, which prevent complete  metrics with $Sc>0$?
 
 Or, conversely,  given an open $n$-manifold $X$, there exists a $n$-manifold $X'$,  such that

(i)  {\sf  $X'$  "contains $X$ at infinity",  i.e.  a complement in $X$ to a {\it relatively compact} open subset,  admits a {\it  proper} imbedding $X\setminus U  \hookrightarrow X'$;

 (ii) $X'$ admits a complete metric with $Sc>0$, or, at least, can be dominated by a complete manifold with $Sc>0$.}
 
 \vspace {1mm}
 
 {\color {red!30!black}Probably}, the minimal  surface argument from  [Wang(Contractible)  2019] shows that 
 
 {\sl 3-manifolds $X'$,
 which  "contains ends"  of  contractible non-simply connected at  (their single ends at) infinity  manifolds $X$,  
  can't be dominated by manifolds with positive scalar curvatures. }
 
 But no such result is in sight for manifolds of dimensions $n \geq 4$. \vspace {1mm}

  \textbf {(Over)Optimistic Existence {\color {red!50!black} Conjecture}  9{\Large\textsubscript\smiley}}.  {\sl All open simply connected manifolds of dimensions $n\geq 4$ admit complete  metrics  with $Sc>0$.}\vspace {1mm}

{\it \color {red!40!black}Questionable Case.}  If $X^{n-1}$ is a simply connected 
manifold, which admits no metric with $Sc>0$, e.g, where $n=4k+1$ and   $\hat A[X^{4k}]\neq 0 $
or where $X^{n-1}$ is Hitchin's  sphere, the results by Cecchini, Zeidler and Zhang may imply 
that  $X=X^n= X^{n-1}\times \mathbb R^1$,
admits no complete  metric  with $Sc>0$.
(Unquestionably, these $X$ admit {\it  no metrics} with $Sc>\sigma>0$, )

In view of this,  it is safer to reformulate   9{\Large\textsubscript\smiley} as follows.
\vspace {1mm}

   \textbf {More Realistic {\color {red!50!black} Conjecture} 9{\Large\textsubscript\frownie}}. {\sl  All open simply connected manifolds $X$ of dimensions $n\geq 4$ with
 $H_{n-1}(X)=0$ (which is equivalent to  "connected at infinity"  for $\pi_1(X)=0$) admit complete metrics with $Sc>0$.}

{\it Example.} The products $X=Y\times \mathbb R^2$, as we know  do admit complete metrics with
$Sc>0$  for all $Y$  and these can be  made simply connected by thin surgery for $dim(X)\geq 4$.

\vspace {1mm}

{\it Non-Example  \textbf {10{\Large\textsubscript\frownie}}. .}  There is no instance of a compact  {\it contractible} manifold $\bar X$ with {\it aspherical}  boundary, where we know whether the interior $X$ of $\bar X$ admits a complete metric with $Sc> 0$.\vspace {1mm}
 
\vspace {1mm}

\textbf { Codimension one Optimistic Reduction {\color {red!50!black} Conjecture}  10}{\Large\textsubscript\smiley}. {\sl Let  $X$  be a complete orientable $n$-manifold with $Sc(X)>0$. If $X$ is orientable, then all   $(n-1)$-dimensional homology classes in   $X$ are realizable by  by smooth closed oriented  hypersurfaces $Y\subset X$,  which  support metrics with 
$Sc>0$.}

But this    { \color {red!40!black} \it  contradicts} to  {\sf  4B{\Large\textsubscript\smiley}}  in the  { \color {red!40!black} questionable case}. Maybe,  it  would be   better to  stick to a  weaker conjecture, e.g. as follows.

 \vspace {1mm}

\textbf { Codimension one more Realistic  {\color {red!50!black} Conjecture}  10}{\Large\textsubscript\frownie}. {\sl  All   $(n-1)$-dimensional homology classes in   $X$ are realizable by the images of the fundamental homology classes of smooth closed ${n-1}$-manifolds $Y$  under continuous maps  $Y\to X$,  where these $Y$ support metrics with 
$Sc>0$.}

\vspace {1mm}

   \vspace{1mm} 
  
  If $X$ is compact,   
  one knows  that  deg $\pm1$ dominants of SYS-manifolds $X_{SYS}$ and manifolds $X_{\kappa\leq 0}$  with non-positive sectional curvatures, as well as their products  $X_{SYS}\times X_{\kappa\leq 0}$ have this $Sc \ngtr 0$ property: they have no   dominants with $Sc>0$;\footnote {I am inclined to think  that  products of    SYS-manifolds may, in general,  carry metrics with $Sc>0$, but I am not certain about it.} we shall prove  in  section {obstructions5}  a similar property for open manifolds, thus confirming the following conjecture in special cases.

  \vspace{1mm}

\textbf { Non-compact Domination   {\color {red!50!black} Conjecture}  11}{\Large\textsubscript\frownie} {\sl If      a compact  orientable $n$-manifold (or pseudomanifold)  $X_0$   can't be dominated (with maps of degree 1) by {\it compact} manifolds  
 with $Sc>0$, then it can't be dominated  by  {\it complete} manifolds with $Sc>0.$}

 \vspace{1mm}
 
Despite the validity of this is known in a variety of {\it specific cases}, including  complete manifolds $X_0$, where non-domination by complete $X$ with $Sc(X)>0$  via {\it proper} maps   implies this property with {\it quasi-proper} ones(see section \ref{domination1}),
 one can't even rule out in general  domination by  complete  manifolds with $Sc\geq \sigma>0.$

%%%%%%%%%%%%%%%%%%%%%%%%%

\section {Variation, Stabilization and Application of   $\mu$-Bubbles} \label {bubbles5}

%%%%%%%%%%%%%%%%%%%%%%%%%

Given a  a Borel measure $\mu$   on an $n$-dimensional  Riemannian manifold $X$, 
  {\it $\mu$-bubbles}  are   critical points of the  following functional on a topologically  defined  class   
  of domains $U\subset X$ with boundaries called  $Y=\partial U$: 
 $$(U, Y)\mapsto vol_{n-1}(Y)-\mu(U).$$

  Observe that  in our examples,   $\mu(U)=\int _U\mu(x)dx$ for (not necessarily positive)  continuous  functions $\mu$ on $X$ and
 that $\mu(U)$  can be regarded as a {\it closed 1-form} on the space of cooriented hypersurfaces  $Y\subset X$. Then  $vol_{n-1}(Y)-\mu(U)$ also comes as such an 1-form
 which we denote $vol^{[-\mu]}_{n-1}(Y)(+const)$.\vspace{1mm}

%%%%%%%%%%%%%%%%%%

\subsection{\color {blue}Second Variation Formula and Pointwise    Scalar Curvature Estimates for 
 $\mathbb T^\rtimes$-Stabilized Bubbles} \label{variation5}

%%%%%%%%%%%%%%%%%%%

The first and the second variations of  $vol^{[-\mu]}_{n-1}(Y)(+const)$ are the sums of these
for $Vol_{-1}(Y) $  and of $vol(U)$ where the former were already computed in section \ref {2nd variation2}. 

And turning to the latter, it is obvious that  the first derivative/variation of   $\mu(U)$  under $\psi\nu$, where $\nu$ is the outward looking  unit normal normal field to $Y$  and $\psi(y)$ is a function on $Y$, is 
$$\partial_{\psi\nu} \int_U \mu(x)dx = \int_Y\mu(y)\psi(y)dy$$
and the second derivative/variation is 
$$ \partial^2_{\psi\nu} \int_U \mu(x)dx =\partial_{\psi\nu} \int_Y\mu(y)\psi(y)dy=\int_Y(\partial_\nu\mu(y)+M(y)\mu(y))\psi^2 (y)dy,$$
where the field $\nu$ is extended along normal geodesics to $Y$, (compare section \ref {2nd variation2})  and where 
$M(y)$ denotes the mean curvature of $Y$ in the direction of $\nu$.

It follows that $\mu$-bubbles $Y$, (critical points of $vol_{n-1}^{[-\mu]}(Y)= vol_{n-1}(Y)-\mu(U)$) have 
$$mean.curv(Y)=\mu(y)$$
and that  

 {\sf second variation   of {\it locally minimal bubbles $Y\subset X$,
$$\partial_{\psi\nu}(vol_{n-1}^{[-\mu]}(Y))= \partial_{\psi\nu}\left (vol_{n-1}(Y)- \int_U \mu(x)dx )\right),$$
is non-positive. }}

Then we recall,  the formula {\Large \color {blue}$[ \circ\circ]$} from section \ref {2nd variation2}
$$\partial^2_{\psi\nu} vol_{n-1}(Y)=\int_Y ||d\psi(y)||^2dy +R_-(y) \psi^2(y)dy$$
for 
$$R_-(y)=-\frac {1}{2}\left ( Sc(Y,y)-Sc(X,y)  +M^2(y) -\sum_{i=1}^{n-1}\alpha_i(y)^2\right ),$$where $\alpha_i(y)$ are  the principal curvatures of $Y$ at $y$, 
and where $\sum \alpha_i^2$ is related to the mean  curvature $M=\alpha_1+...+\alpha_{n-1}$, by the inequality
$$ \sum \alpha_i^2\geq \frac {M^2}{n-1}.$$
Thus, summing up all of the above, observing that 
$$\partial_\nu\mu(x)\geq -||d\mu(x)||$$
and letting
 $$R_+(x)=\frac {n\mu(x)^2}{n-1} -2||d\mu(x)|| +Sc (X,x),\leqno {\mbox {\large\color {blue} $[R_+=]$}}$$
we conclude that

{\it  if $Y$  locally minimises $vol_{n-1}^{[-\mu]}(Y)(=vol_{n-1}(Y)-\mu(U))$, then
$$\int ||d\psi||^2dy +\left (\frac {1}{2} Sc(Y)-\frac {1}{2} R_+(y)\right)\psi^2(Y)dy \geq \partial_{\psi\nu}vol_{n-1}^{[-\mu]}(Y)
\geq 0$$
 for all functions $\psi$ on $Y$.}
 
 Hence, \vspace{1mm}
 
 {\color {blue}\small\PlusCenterOpen$_{\geq 0}$} {\it the   $-\Delta +\frac {1}{2}Sc(Y,y) -\frac {1}{2}R_+(y)$, for  $\Delta=\sum_i\partial_{ii}^2$  is  positive on $Y$.}
\vspace{1mm}

 {\it Examples.} (a)  Let $X=\mathbb R^n$ and $\mu(x) =\frac{n-1}{r}$, that is the mean curvature of the sphere of radius $r$.
Then 
$$R_+(x)=\frac {n(n-1)}{R}-2\frac {n-1}{r^2}+0=\frac {(n-1)(n-2)}{r^2}= Sc(S^{n-1}(r)).$$

(b) Let $X =\mathbb R^{n-1}\times \mathbb R$ be the hyperbolic space with the metric 
$g_{hyp}=e^{2r}g_{Eucl} + dr^2$ and let $\mu(x)=n-1$. Then 
$$R_+(x)=n(n-1) -0 +(-n(n-1))=0= Sc(\mathbb R^n). $$

(c) Let  $X=Y\times \left ( - \frac {\pi}{n}, \frac {\pi}{n}\right)$ 
 with  the metric $\varphi^2h+dt^2$, where the metric  $h$ is a metric on  $Y$ and where 
  $$\varphi(t) =\exp \int_{-\pi/n}^t -\tan \frac {nt}{2} dt, \mbox { }  -\frac {\pi}{n}<t < \frac {\pi}{n}.$$
Then a simple computation shows that  
$$R_+(x)=\frac {n(n-1)}{R}-2\frac {n-1}{r^2}+0=\frac {(n-1)(n-2)}{r^2}= Sc(S^{n-1}(r)).$$
 $$\frac {n\mu(x)^2}{n-1} -2||d\mu(x)||+n(n-1)=0.$$
 Furthermore, if  $Sc(h)=0$, than $Sc(X(=n(n-1)$ and $R_+=0$.
\vspace{1mm}

Two relevant corollaries to  {\color {blue}\small\PlusCenterOpen$_{\geq 0}$} are as follows.\vspace{1mm}

 {\sf Let $X$ be a Riemannian manifold of dimension $n$, let $\mu(x)$ be a continuous function 
 and  $Y$ be a smooth minimal $\mu$-bubble in $X$}.
 \vspace{1mm}

 {\color {blue}\small\PlusCenterOpen$_{conf}$} {\sl If 
 $$R_+(x)=\frac {n\mu(x)^2}{n-1} -2||d\mu(x)|| +Sc (X,x)> 0,$$ 
 then by Kazdan-Warner  conformal change theorem} (see section \ref {conformal2})  {\sl $Y$ admits a metric with $Sc>0$.}

\vspace{1mm}

{\color {blue}\small\PlusCenterOpen$_{warp}$} {\sl 
   There exists a metric $\hat g$ on the  product $Y\times \mathbb R$ of the form $g_Y+\phi^2dr^2$
for the metric $g_Y$ on $Y$ induced from $X$, such that
$$Sc_{\hat g}(y,r)\geq R_+(y).$$}
\vspace{1mm}
which implies that $Sc_{\hat g}(y,r)\geq R_+(y)$, since $\lambda \geq 0$. QED.\vspace{2mm}

%%%%%%%%%%%%%%%%%%%%%%%%
\subsection {\color {blue} On Existence and Regularity of Minimal Bubbles}\label {existence5}
   
   %%%%%%%%%%%%%%%%%%%%%%%%%%%%%%%%%

    {\sf Let $X$  be  a  compact  connected Riemannian manifold  of dimension $n$  with boundary $\partial X$ and let $\partial_-\subset \partial X$  and $\partial_+\subset \partial X$  be disjoint compact domains  in $\partial X$.}
    \vspace{1mm}
  
  {\it Example.}  Cylinders $Y\times [-1,1]$ naturally come with such a  $\partial_\mp$-pair for
  $ \partial_-=Y\times \{-1\}$ and $ \partial_+=Y\times \{1\}$, where, observe,  $ \partial_-\cup \partial_+
  =\partial (Y\times [-1,1])$ if and only if $Y$ is a manifold without boundary.    
      \vspace{1mm}
  
  Let us agree that the mean curvature of $\partial _-$ is evaluated with the incoming normal field 
  and $mean.curv(\partial_+)$ is evaluated with the outbound field.

   For instance, if the boundary of $X$ is {\it concave}, as for instance for $X$ equal to  the sphere minus two small disjoint balls, t
    then $mean.curv(\partial_-) \geq 0$ and   $mean.curv(\partial_+)\leq 0$.       \vspace{1mm}

   {\it Barrier  {\color {blue}$[\gtrless\mp mean]$}-Condition.} {\sf A continuous  function $\mu(x)$ on $X$  is said to satisfy { \color {blue}$[\gtrless\mp mean]$}-condition if 
  $$\mu(x) \geq  mean.curv(\partial_-, v)\mbox { and } \mu(x) \leq  mean.curv(\partial_+, x) \leqno { \color {blue}[\gtrless\mp mean] }$$
 for all  $x\in \partial_-\cup \partial_+. $}

  \vspace{1mm}
 
  It follows by  the maximum principle in the   geometric measure theory that  \vspace{1mm}

  {\huge \color {blue} $\star$}  {\sl the  { \color {blue}$[\gtrless\mp mean]$}-condition ensures the existence of a minimal 
  $\mu$-bubble $Y_{min}\subset X$.
  which separates $\partial_-$ from $\partial_-+$.}  \vspace{1mm}

   If this condition is {\it strict}, i.e. if $\mu(x)> mean.curv(\partial_-)$  and   $\mu(x)< mean.curv(\partial_+)$  and if 
 $X$ has no boundary apart from  $\partial_\mp$, then $Y_{min}\subset X$ doesn't intersect $\partial\mp$; in general, the intersections $Y_{min}\cap \partial_\mp$ are contained in the 
 {\it side boundary} of $X$ that is the closure of the complement 
 $\partial X\setminus (\partial_-\cup \partial_-)$.
 (This, slightly reformulated, remains   true for non-strict { \color {blue}$[\gtrless\mp mean]$}.)

 \vspace{1mm}

If $dim(X) = n \leq  7$,  then, (this well known and easy to see) Federer's regularity theorem(see  section \ref {SY+symplectic2}) applies to minimal bubbles as well  as to minimal subvarieties and the same can be said about Nathan Smale's theorem on non- stability of singularities for $n = 8$. Thus, in what follows we may assume our minimal bubbles smooth for $n \leq 8$.
  
  Then, by the stability of $Y_{min}$ (see section \ref {variation5} above), \vspace {1mm}

 {\color {blue} \Large  $\bullet_{\varphi_\circ}$}:  {\sl there exits  a function $\phi_\circ =\phi_\circ(y) >0$  defined in the interior $^\circ Y$  of $Y$, i.e. on
$Y\setminus \partial X$, 
such that the 
metric 
$$g_{\varphi_\circ}= \varphi_\circ^2g_Y +dt^2\mbox {  on  the cylinder } ^{\circ}Y\times \mathbb R,$$ 
where $g_Y$ is the Riemannian metric on $Y$  induced from $X$, satisfies
$$Sc_{g_{\varphi_\circ}}(y, t)\geq Sc(X,y)+\frac{ n\mu(y)^2}{n-1}-2||d\mu(y)||\leqno {\mbox {\color {blue}\EllipseShadow}}$$ 
for all $y\in \hspace {0.3mm} ^{\circ}\hspace {-0.3mm}Y$}.\footnote {Since the metric $g_{\varphi_\circ}$ is $\mathbb R$-{\it invariant} its scalar curvature is {\it constant} in $t\in \mathbb R$.}

\vspace {1mm}

      \hspace {32mm}{\sc \color {blue} \textbf {What if $n\geq 9$?}.}
\vspace{1mm}

The overall  logic of the proof indicated  in  [Lohkamp(smoothing) 2018]   leads one to believe that, assuming strict   { \color {blue}$[\gtrless\mp mean]$},  there always  exists a smooth $Y_o \subset X$, which   separates $\partial_\mp$ and and which admits a function  $\phi_\circ$ with the property 
{\color {blue}\EllipseShadow}.

  The proof of this, probably, is    automatic,  granted a full   understanding 
 Lohkamp's arguments. But since I have not seriously studied these arguments, everything which follows in sections  5.3-5.8   should be regarded as {\it \color {magenta} \large conjectural} for $n\geq 9$.\footnote {In some cases, a generalization of Schoen- Yau's  theorem 4.6 from [SY(singularities) 2017] can be used instead of Lohkamp's theory; namely, this is possible in those applications,
which  don't depend on   the Dirac operators on these bubbles,  but  can  be obtained  by  relying only  on the geometric measure theory.}\vspace {1mm}

  {\it Barrier  {\color {blue}$[\gtrless mean=\mp\infty]$}-Condition.} Let $X$ be  a non-compact, possibly non-complete,  Riemannian manifold $X$ and let the set of the ends   of  
  $X$ is subdivided to $(\partial_\infty)_- = (\partial_\infty)_-(X)$  and  $(\partial_\infty)_+=(\partial_\infty)_+(X)$, where this can be  accomplished, for instance,  with a proper map from $X$ to an open (finite or infinite) interval $(a_-,a_+)$ where  "convergence"
  $x_i\to (\partial_\infty)_\mp$, $x_i\in X$, is defined as $e(x_i)\to a_\mp$. 

For example,  if $X$ is the open cylinder, $X=  Y\times  (a,b) $, where $Y$ is a compact manifold, possibly with a boundary, this is done  with the projection  
$Y\times  (a_-,a_+) \to   (a_-,a_+)$.\vspace{1mm}

{\it\large  Obvious Useful Observation.}  {\sf If  a   function $\mu(x)$ satisfies
$$\mu(x_i) \to \pm \infty 
\mbox {  for } x_i\to (\partial_ \infty)_\mp$$ 
then $X$ can be exhausted by compact manifolds $X_i$  with distinguished
 domains $(\partial_\mp)_i\subset \partial X_i$, such that
 
$\bullet $ these   $(\partial_\mp)_i$ separate  $(\partial_ \infty)_-$ from  $(\partial_ \infty)_-$  for all $i$  and 
$$(\partial_\mp)_i\to  (\partial_ \infty)_\mp;$$
 
 $\bullet $ restrictions of $\mu$ to    $(X_i, (\partial_\mp)_i)$
 satisfy the  barrier  {\color {blue}$[\gtrless\mp mean]$}-condition.}

 {\it This  ensures the existence of locally  minimising $\mu$-bubbles in $X$  which  separate  $(\partial_ \infty)_-$ from  $(\partial_ \infty)_+$.}

%%%%%%%%%%%%%%%%%%%%%%%%%%%%%

\subsection {\color {blue} Bounds on Widths of Riemannian Bands and on Topology of Complete Manifolds with $Sc>0$}\label {bands5}

%%%%%%%%%%%%%%%%%%%%%%%%%%%%%%

Let us  prove  the following version of the {\color {blue}$\frac {2\pi}{n}$-inequality}  from section \ref {bands3}.

{\it \large \color {blue}$\frac {2\pi}{n}$-Inequality$^\ast$.} {\sf Let $X$ be an open, possibly non-complete Riemannian manifold of dimension $n$ and let 
$$f: X\to (-l,l) $$
be a  proper (i.e. infinity $\to$ infinity) smooth  distance non-increasing  map, such that the pullback
$f^{-1}(t_o)\subset X$ of  a generic point $t_o$ the interval $(-l,l)$ is non-homologous to zero in $X$.}

{\it If $Sc(X)\geq n(n-1)=Sc(S^n)$  and if the following condition {\color {blue}  \textSFliii$_{Sc \ngtr 0}$ }  is satisfied, 
then 
$$l\leq \frac {\pi}{n}.$$}

{\color {blue}  \textSFliii$_{Sc \ngtr 0}$ } {\sf No smooth closed cooriented  hypersurface in $X$ homologous to  $f^{-1}(t_o)$  admits a metric with $Sc>0$. }

\vspace {1mm}

{\it Proof.}  Assume   $l>  \frac {\pi}{n}$. and  let $\underline \mu(t)$ denote the mean curvature of the hypersurface 
$\underline Y\times \{t\} $  in the warped  product metric  $\varphi^2h+dt^2$. on $\underline Y\times \left ( - \frac {\pi}{n}, \frac {\pi}{n}\right)$ for 
$$\varphi(t) =\exp \int_{-\pi/n}^t -\tan \frac {nt}{2} dt, \mbox
 { }  -\frac {\pi}{n}<t < \frac {\pi}{n}$$
as in example (c) from the previous section.

Since $\underline \mu(t)\to \pm\infty$ for $t\to \mp \frac {\pi}{n}$, the   barrier  {\color {blue}$[\gtrless mean=\mp\infty]$}-condition  from the  section \ref {existence5} guaranties the existence of 
a locally minimizing $\mu$-bubble in $X$ for $\mu$ being a slightly modified $f$-pullback of 
$\underline \mu$ to $X$. 

Let us spell it out in detail.

Assume without loss of generality that the pullbacks $Y_{\mp}= f^{-1}\left (\mp\frac{\pi}{n}\right)\subset X$ are smooth, and let $\mu(x) $ be a smooth  function on $X$ with the following properties.

$\bullet_1$ $\mu(x) $ is   constant on $X$ on   the complement of $f^{-1} \left ( - \frac {\pi}{n}, \frac {\pi}{n}\right)$ for 
$ \left ( - \frac {\pi}{n}, \frac {\pi}{n}\right)\subset (-i,i)$;

 $\bullet_2$ $\mu(x) $ is   equal to  $\underline \mu\circ f$ in the interval  $\left ( - \frac {\pi}{n}+\varepsilon , \frac {\pi}{n}-\varepsilon \right)$ for a given (small) $\varepsilon >0$;

 $\bullet_3$ the absolute values of the mean curvatures of the hypersurfaces  $Y_{\mp}$ are everywhere smaller than the absolute values of $\mu$;

 $\bullet_4$ $ \frac{n \mu(x)^2}{n-1}-  2||d\mu(x) ||+ n(n-1)\geq 0$ 
  at all points $x\in X$.\vspace {1mm}

In fact, achieving  $\bullet_3$ is possible, since $\underline \mu(t)$ is infinite at $\mp\frac {\pi}{n}$, while the mean curvatures of the hypersurfaces  $Y_{\mp}$ and  what is needed for 
  $\bullet_4$ are  the inequality   $||df||\leq 1$ and 
the equality 
$$\frac{n\underline \mu(t)^2}{n-1}-\vert \frac{d\underline \mu(t)}{dt}\vert +n(n-1)=0$$ indicated in   example (c) from section {5.1}). 

Because of  $\bullet_3$, the submanifolds $Y_{\mp}$ serve as barriers for $\mu$-bubbles (see  the previous section) between them;  this implies the  existence of  a minimal $\mu$-bubble $Y_{min}$ in the subset 
$f^{-1} \left ( - \frac {\pi}{n}, \frac {\pi}{n}\right)\subset X$  homologous to $Y_o$.
 by  {\huge \color {blue} $\star$} in section \ref {existence5}. 

Due to  $\bullet_4$, the  $\Delta + \frac{1}{2}Sc(Y)$
 is positive by  {\color {blue}\small\PlusCenterOpen$_{\geq 0}$} from the  section \ref {variation5}.

Hence, by {\color {blue}\small\PlusCenterOpen$_{conf}$} the manifold  $Y_{min}$ admits a metric   with $Sc>0$ and the inequality $l\leq \frac {\pi}{n}$ follows.\vspace {1mm}

{\color {blue} \large \it On Rigidity.}  A a close look at
minimal $\mu$-bubbles (see section \ref {rigidity5}) shows that\vspace {1mm}

{\sl if $l= \frac {\pi}{n}$, then $X$ is isometric to a warped product , $X=Y\times \left ( - \frac {\pi}{n}, \frac {\pi}{n}\right)$ 
 with  the metric $\varphi^2h+dt^2$, where the metric  $h$ on $Y$ has $Sc(h)=0$ and where 
  $$\varphi(t) =\exp \int_{-\pi/n}^t -\tan \frac {nt}{2} dt, \mbox { }  -\frac {\pi}{n}<t < \frac {\pi}{n}.$$}

 \vspace {1mm}
 
 {\it Exercises.} (a) Let $X$ be an open manifolds with two ends, Show that if no  closed hypersurface in $X$ that separates  the ends admits a metric with positive scalar curvature then $X$ admits no metric with $Sc>0$ either.\footnote 
 {This, for  a class of spin manifolds $X$,   was  shown in [GL(complete) 1983]   by applying  a relative index theorem for suitably  twisted Dirac operators on $X\times S^2(R)$.}
   
  \vspace{1mm}
  
(b)  Let $X$ be a complete Riemannian manifold,  and let
 $$S(R)=\min_{B(R)} Sc(X)$$ denote the minimum of the scalar curvature (function)  of $X$ on the ball  $B(R)=B_{x_0}(R)\subset X$  for some centre point $x_0\in X$.  Show that \vspace {1mm} 
 
 {\sf if $X$ is  homeomorphic to $\mathbb T^{n-2}\times \mathbb R^2$, then there exists a constant $R_0=R_0(X, x_0)$, such that} 
 $$  S(R)\leq \frac{4\pi^2}  {(R-R_0)^2}  \mbox  {  for all } R\geq R_0. \footnote {
 A rough version of this for a class of spin manifolds $X$  can be proved by   Dirac operator methods.}\leqno {\color {blue}[\asymp \frac {4\pi^2}  {R^2}]}$$

{\it Hint.}  Since  the bands  between the concentric  spheres of   radii $r$ and $r+R$,  call them 
$X(r, r+R)=B(r+R)\setminus B(r)$,
 are,  for large $r$, quite similar to   the cylinders  $\mathbf T^{N-1}\times [0,R]$,    the   {\it \large \color {blue}$\frac {2\pi}{n}$-Inequality$^\ast$} applies to them and  says that 
 their scalar curvatures  satisfy

$$S(R)=\inf Sc _x (X(r, r+R),x)  \leq  \frac {4(n-1)\pi^2}{nR^2}.$$

 \vspace {1mm}

 {\it Question.}  {\sf What  are topological obstructions, if any,  for the existence  of a complete Riemannian metric $g$ on an open manifold     $X$ (possibly with a boundary),  such that  $Sc(g)>0$ and/or   $Sc(g)\geq 0$   at infinity, i.e. outside a given  compact subset in $X$. }
  \vspace {1mm}
mple
 
We describe obstructions  on  topology implied by positivity of the scalar curvature   in sections  \ref{obstructions4}  and \ref{obstructions5}; {\footnote {Also see  
[GL(complete) 1983], 
[Cecchini(Callias) 2018], [Wang(Contractible)  2019]) for the existence of  complete metrics with {\it everywhere} positive scalar curvatures, where the techniques from [Cecchini(long neck) 2020]  and/or from  [Zhang(Area Decreasing) 2020] may be(?) applicable to $Sc>0$ {\it at infinity.}}
here we make a couple of preparatory remarks  concerning this issue.
 \vspace {1mm}

(a)  If  the  sectional curvature of  a complete Riemannian manifold $X$ is non-negative at infinity, then, by the standard argument, $X$ admits a proper continuous  function  $f:X\to [0,\infty)$, where the  levels $f^{-1}(r)\subset X$ are {\it convex hypersurfaces  for all $r\in [1,\infty)$.} Thus, 

\hspace {20mm}{\it $X$ is topologically cylindrical at infinity. }\vspace {1mm}

(b) Similarly to (a),  if  $Ricci (X)\geq 0$ at infinity, then  $X$ admits a proper continuous  function  $f:X\to [0,\infty)$, where the  levels $f^{-1}(r)\subset X$ have {\it non-negative (generalized) mean curvatures   for all $r\in [1,\infty)$.} Thus,

\hspace {30mm} {\it $X$ has finitely many ends.}\vspace {1mm}

(c)  Let $X =Y_0\times \mathbb T^{n-2}$, where $Y_0$ is connected surface of infinite topological type, e.g. with infinitely many ends.

Then "the first "-Torical symmetrization  from section \ref {quadratic3}, at least if $n\leq 8$,  brings us to a complete  surface $Y\subset X$ of infinite topological type, such that
a (generalized)   warped product $g^\circ$ metric on $Y\times  \mathbb T^{n-2}$ has $Sc>0$ at the points $x\in Y$, where $Sc(X,x)>0$.

Now, to prove that $g^\circ$ can't have $Sc(g^\circ)>0$ at infinity, one needs to find a  complete  $\mathbb T^{n-2}$-invariant volume minimizing hypersurface that doesn't intersect a given compact subset   $K\subset Y\times  \mathbb T^{n-2}$, 
where such a hypersurface can be seen as a minimizing geodesic in  the surface $Y$ with the metric,  which is conformally equivalent to the metric $dy^2$  induced  from $X$, and where the conformal factor is equal to 
$(vol(\mathbb T_y^{n-2})^2$.
 
In general, the volumes of the tori  $\mathbb T_y^{n-2}=\mathbb T^{n-2}\times \{y\}$  may grow very fast  for $y\to \infty$,  such that all minimal  hypersurfaces intersect the  subset $K$, where $Sc\leq 0$,   but there is a limit to such a growth due to the differential inequality satisfied by   the conformal factor (see section \ref{warped+2}. Besides, a significant  growth of this factor, may allow  stable $\mu$-bubbles  away from $K$.

But it is unclear if this can be made rigorous and

{\sf  non-existence of complete metrics with $Sc>0$  on the above $Y_0\times \mathbb T^{n-2}$ remains {\color {blue}conjectural.}}

On the other hand, this kind of  reasoning  rules out complete metrics with $Sc>0$  everywhere on many manifolds  with sufficiently complicated topologies, which may include manifolds considered in [Cecchini(Callias) 2018]  and/or in [Wang(Contractible)  2019].\footnote{Cecchini's proof, which  applies to spin manifolds of all dimensions,  depends on the index theory for  Dirac-type 
operators, while     Wang's argument, which    relies  on  specifically 2-dimensional   properties of minimal surfaces,  shows that certain contractible 3-manifolds admit no metrics with $Sc>0$.}

%%%%%%%%%%%%%%%%%%%%%
 \subsection {\color {blue} Equivariant Separation and  Bounds on Distances Between Opposite  Faces of Cubical Manifolds with $Sc>0$}\label {separating5}
%%%%%%%%%%%%%%%%%%%%%%%%%%%

Recall  the following general purpose proposition  from section \ref {separating3]}. \vspace {1mm}

 {\color {blue} {\large \sf I{\large \sf \color {red!80!black}\sf I}\sf I}$_\circlearrowleft$} \textbf{Equivariant  Separation Theorem}. {\sf Let $X$ be an $n$-dimensional,  Riemannian band, possibly non-compact and non-complete.

   Let  $$Sc(X,x)\geq \sigma(x)+\sigma_1, \mbox {  },$$
   for a continuous function $\sigma=\sigma(x)\geq 0$ on $X$ and a constant $\sigma_1> 0,$
   where $\sigma_1$  is related to  $d=width(X)=dist_X(\partial_-,\partial_+)$ by the inequality
 $$\sigma_1 d^2>   \frac {4(n-1)\pi^2}{n}.$$    
(If scaled  to $\sigma_1=n(n-1) $, this   becomes $d> \frac {2\pi}{n}$.)

{\it Then   there exists a smooth  hypersurface $Y  \subset X$,   which separates
  $\partial_-$ from  $\partial_+ $, and a smooth positive  function $\phi$ on   $Y$, 
such that the scalar curvature of the  metric $g_\phi = g_{Y_{-1}}+\phi^2dt^2$ on $Y\times \mathbb R$ is bounded from below by 
$$Sc(g_\phi,x)\geq \sigma(x).$$}

Furthermore} 

{\it if $X$ is isometrically acted upon by a compact  connected  group $G$,
then the separating hypersurface $Y\subset X$ and the function $\phi$ on $Y$ can be chosen invariant under this action.}

{\it Proof.} The    general case of this reduces to that of $\sigma=n(n-1)$ by on obvious scaling/rescaling argument and when $\sigma=n(n-1)$ we use the same $\mu$ as above associated with 
$\varphi(t)= \exp \int_{-\pi/n}^t -\tan \frac {nt}{2} dt, \mbox { }  -\frac {\pi}{n}<t < \frac {\pi}{n}.$
Then, as earlier, since
$$Sc_{g_{\varphi_\circ}}(y, t)\geq Sc(X,y)+\frac{ n\mu(y)^2}{n-1}-2||d\mu(y)||$$
by 
{\color {blue}\EllipseShadow}   from the previous section,  the above  equality  $\frac{n\underline\mu(t)^2}{n-1}-\vert \frac{d\underline\mu(t)}{dt}\vert +n(n-1)=0$ implies the  requited bound  $Sc(g_o)\geq \sigma_1$. QED.

\vspace {1mm}

{\it Example of Corollary} {\sf Let  $X$ be an {\it orientable spin} manifold, let   $\partial_-\cup \partial_+=\partial X$ and let  $f:X\to S^{n-1} \times [-l,l]$ be a smooth map, such that $\partial_\mp\to  S^{n-1}\times \{\mp l\}$.} 

Let the following conditions be satisfied.

{\sf $\bullet $  $deg(f)\neq 0$, 

$\bullet $ the map $X\to S^{n-1}$, that is the  composition of $f$ with the projection
$S^{n-1} \times [-l,l]\to S^{n-1}$, is {\it area decreasing};

 $\bullet $ $Sc(X)\geq (n-1)(n-2)+\sigma_1$ for some $\sigma_1\geq 0$.}
 
  {\sl Then  in conjunction with the (stabilised) Llarull theorem shows that 
$$dist (\partial_-,\partial_+)\leq \frac {2\pi}{n} \frac {n(n-1)}{\sqrt \sigma_1}=\frac {2\pi(n-1)}{\sqrt \sigma_1}.$$}

{\it Remark.}  This inequality if it looks   sharp, then only for  $\sigma_1\to 0$, while sharp(er) inequality of this kind need different   functions $\mu$.

 \vspace {2mm}

 \vspace {1mm}
 
 {\it {\large $\mathbf \square^{n-m}$-\textbf {Theorem. }} {\sf Let $X$ be a compact connected   orientable Riemannian manifold of dimension $n$ with a boundary   and let $ \underline X_\bullet$ is a closed orientable manifold of dimension $n-m$, e.g. a single point $\bullet$  if $n=m$.
 
 Let
 $$  f:X\to [-1,1]^m \times \underline X_\bullet$$
   be a continuous map,  which sends the boundary of $X$ to the boundary of $[-1,1]^m \times \underline X_\bullet$  and which has  {\it non-zero degree}.

Let $\partial_{i \pm}\subset X$, $i=1,...,m$, be the pullbacks of the pairs of the opposite faces of the cube $[-1,1]^m$ under the composition of $f$ with the projection  
$[-1,1]^m \times \underline X_\bullet \to [-1,1]^m$.

Let $X$ satisfy the following condition:\vspace {1mm}

{\color {blue}  \textSFliii$_{Sc \ngtr 0}^m$ }  \hspace{1mm} No transversal intersection $Y_{-m\pitchfork}\subset X$ of $m$-hypersurfaces $Y_i\in X$ which separates $\partial_{i -}$ from $\partial_{i +}$, admits a metric with $Sc>0$; moreover, the products $Y_{-m\pitchfork}\times \mathbb T^m$ admit no  metrics with $Sc>0$ either.\footnote {This "moreover" is  unnecessary, since the  relevant for us case  of  stability of the  $Sc \ngtr 0$ condition under multiplication by tori   is more or less automatic. (The general case needs some effort.)} 

{\it If   $Sc(X)\geq n(n-1)$, then the distances  $d_i=dist(\partial_{i-},\partial_{i +})$ satisfy the following inequality} (which generalise  {  \color{black}  $\square^n$-\textbf {inequality}}  from  section \ref{multi-width3}).
{\it $$\sum_{i=1}^m \frac{1 }{d_i^2}\geq  \frac  {n^2}{4\pi^2} \leqno {\color {blue}\square_{\sum}}$$
Consequently 
$$min_i dist (\partial_{i -},\partial_{i +})\leq \sqrt m \frac { 2\pi}{n}.\leqno {\color {blue}\square_{\min}}$$}}}

{\it Proof.} Let  $$\sigma'_i=\left (\frac {2\pi}{n}\right)^2\frac {n(n-1)}{d^2}=\frac {4\pi^2(n-1)}{nd^2}$$

and rewrite {\color {blue}$\square_{\sum}$} as
$$ \sum_i\sigma'_i\geq n(n-1).$$

Assume $ \sum_i\sigma'_i< n(n-1)$ and let  $\sigma_i>\sigma'_i$ be such that $\sum_i\sigma_i<n(n-1)$.

 Then, by  induction on $i=1,2,...,m$ and using  $\mathbb R^{i-1}$-invariant  $\square$-Lemma  on the $i$th step, construct   manifolds $X_{-i}=Y_{-i}\times \mathbb R^i$ with $\mathbb R^i$-invariant metrics $g_{-i}$, such that
$$Sc(X_{-i})>  n(n-1)-\sigma_1-..-\sigma_i.$$
The proof is concluded by observing that this  for $i=m$ would contradict  to {\color {blue}  \textSFliii$_{Sc \ngtr 0}^m$.} 

\vspace {1mm}

{\it Remarks.} (a) As we mentioned earlier, this inequality is non-sharp starting from $m=2$, where    the sharp inequality 
$$min_{i=1,2} dist (\partial_{i -},\partial_{i +})\leq\pi.\leqno {\square^2_{\min}}$$
for squares with Riemannian metrics on them with $Sc\geq 2$ follows by an elementary argument.

(b) One can show  for all  $n$  that
$$min_i dist (\partial_{i -},\partial_{i +})\leq \sqrt m \frac { 2\pi}{n}-\varepsilon _{m,n},$$
where  $\varepsilon _{m,n}>0$  for $m\geq 2$.
\vspace{1mm}

(c) A possible way for  sharpening  ${\square_{\sum}}$,  say for the case $m=n$, is by using $n-2$   inductive steps instead of $n$ and then generalizing the elementary proof of $\square^2_{\min}$  to  $\mathbb T^{n-2}$-invariant metrics on $[-1,1]^2\times \mathbb T^{n-2}$.

In fact, all theorems for surfaces $X$ with  positive (in general, bounded from below) sectional curvatures beg for their  generalisations to  $\mathbb T^{m-2}$-invariant metrics on $X\times \mathbb T^{m-2}$  with  positive (and/or  bounded from below) scalar curvatures. 

%%%%%%%%%%%%%%%%%
\subsubsection  {\color {blue} Max-Scalar Curvature  with  and without Spin} \label {max-scalar5}

%%%%%%%%%%

It remains a {\large \sf \color {red!30!black} big open problem} of making sense of  the inequality $Sc(X)\geq \sigma$,  e.g. for $\sigma=0$, for {\it non-Riemannian} metric spaces, e.g. for piecewise smooth polyhedral spaces $P$.

But   lower bounds on Lipschitz constants of homologically  substantial maps $X\to P$ entailed by the inequality $Sc(X)\geq \sigma>0$,   that, for a fixed $P$, tell you  something about the    geometry of $X$, can be used the other way around for the definition of  scalar curvature-like  invariants of general metric spaces $P$ as follows.

  \vspace {1mm}

Given a metric space $P$  \footnote {To be specific  we  assume that $P$ is locally compact and locally contractible, e.g. it is locally triangulable space} and a homology class $h\in H_n(P)$ define $Sc^{\sf max}({ h})$ as the supremum  of the numbers $\sigma\geq 0$, such that $H$ can be {\it dominated with $Sc\geq \sigma$}.  Here (slightly unlike how it is in  \ref {domination1})
this means that

{\sf  there exists a closed 
 orientable   Riemannian $n$-manifold $X$ and a 1-Lipschitz 
map $f: X\to  P$,  such that  the fundamental homology class $[X]$ goes to $h$,
$$f_\ast[X]=h.$$}

Similarly,  one defines   $Sc^{\sf max}_{sp}({ h})$ by allowing only {\it spin} manifolds $X$, where, for instance, the discussion in section 4.1.1 shows that
$$Sc^{\sf max}_{sp}({ h})\leq const_n\cdot  \mbox {K-}waist_2(h).$$

Below are a few observations concerning these definitions.  

 $\bullet_1 $ $Sc^{\sf max} [X]\geq \inf_xSc(X,x)$   for all closed {\it Riemannian} manifolds $X$, where the equality 
 $Sc^{\sf max} [X]=Sc(X,x)$, $x\in X,$ holds  for what we  call  {\it extremal} manifolds $X$.
 
$\bullet_2 $  More generally, the product  homology  class $ h\otimes [X]\in H^{n+m}(P\times X)$, $m=dim(X)$, where $P\times X$ is endowed with the Pythagorean product metric,  satisfies 
  $$ Sc^{\sf max} (h\otimes [X]) \geq  Sc^{\sf max} (h)+ \inf_xSc(X,x).$$

{\color {blue} $\bullet_3 $
 \it Possibly,}  $$ Sc^{\sf max} (h\otimes [S^m])= Sc^{\sf max} (h)+m(m-1),$$ 
 but even   the rough inequality $$ Sc^{\sf max} (h\otimes [S^m]) \leq  Sc^{\sf max} (h)+ const_m.$$
remains beyond   splitting  techniques from  section 5.3. \footnote{These techniques deliver such an inequality for     the   stabilized max-scalar curvature: 
 $Sc^{\sf max stab}(h)= \lim_{m\to\infty} (h\otimes [\mathbb T^m])$, where  one may additionally require the  manifolds $X$ mapped to $P\times \mathbb T^m$ to be isometrically acted upon by the $m$-tori }

\vspace{1mm}
 
$\bullet_4 $  If $F:X_1\to X_2$ is a finitely sheeted  covering  between closed orientable Riemannian manifolds, 
then  \vspace {1mm} 

\hspace {10mm}$Sc^{\sf max}_{sp}[X_1] \geq Sc^{\sf max}_{sp}[X_2]$ as well as 
$Sc^{\sf max}_{sp}[X_1] \geq Sc^{\sf max}_{sp}[X_2]$, \vspace {1mm} 

\hspace {-6mm} but }the equality may fail to be true, e.g. for SYS-manifolds $X_2$   defined in section \ref{SY+symplectic2}

(It is less clear when/why this  happens to {\it infinitely sheeted} coverings, where the problem  can be  related to    possible failure  of  contravariance of $K$-$waist_2$, see 
 section \ref{cowaist4})
 \vspace {1mm} 

{\it Non-Compact Spaces and $Sc^{\sf max}_{prop}$.} The above definitions naturally extends to homology with infinite supports in non-compact spaces , e.g. to the fundamental classes $[P]$ of open manifolds and  pseudomanifolds  $P$,
  where the Riemannian manifolds $X$ mapped  to these spaces  are now non-compact and not even complete. 
  
  Also we use the  notation $Sc^{\sf max}_{prop}$ for fundamental classes of (psedo)manifolds $P$ {\it with boundaries}, where {\it proper maps} $X\to P$   are those sending $\partial X\to \partial P$.

 \vspace {1mm} 

{\it \color {blue} Stabilized $max$-Scalar Curvatures.} These for a space $P$ are defined as 
$$stabSc^{max}_{...}(P) =  Sc^{max}_{...}(P\times \mathbb T^N)$$
where  $\mathbb T^N$ is flat torus that may be assumed arbitrarily large (this proves  immaterial at the end of day), where $N$
 is also large and where
the implied metric  in the product is the Pythagorean one: 
$$dist((p_1,t_1),(p_2,t_2))=\sqrt {dist(p_1,p_2)^2+dist(t_1,t_2)^2 }. $$\vspace {1mm}

{\it Examples. }  (a)   Llarull's and Goette-Semmelmann's inequalities from section 4.2  can be regarded as  sharp  bounds on  
$Sc^{\sf max}_{sp}$  for (the fundamental homology  classes of)  spheres and  convex hypersurfaces.

  \vspace {1mm} 

(b)  The $\square$-inequalities from the previous   section provide similar  bounds  on stabilised $Sc^{\sf max}_{prop}(P)$  for the fundamental homology classes  of the  rectangular solids $P=\bigtimes _{i=1}^n[0,a_i]$. 

  \vspace {1mm}

(It seems, there are interesting examples in the spirit of  SYS-spaces from section \ref{SY+symplectic2} where   one needs to  allow $ f_\ast [X]_{ \mathbb  Z/l\mathbb  Z}\neq 0$,  at least for  for odd $l$.

Also one may ask in this regard if  
$Sc^{\sf max}_{prop}$ of the universal covering of a  closed orientable manifold $X$ with a residually finite fundamental group   is equal to the limit of $Sc^{\sf max}_{prop}$ of the finite coverings of $X$.)\vspace {1mm}

(c) {\it Spaces with S-Conical Singularities and $Sc\geq \sigma$.}    Let us define  classes $\mathscr S^n_{\geq\sigma}$, $n=2,3,...   $  of piecewise Riemannian spaces with 
$Sc\geq \sigma>0$  by induction on dimension $n\geq 2$  as follows.

Let $Y=Y^{n-1} $ 
from $\mathscr S^{n-1}_{\geq\sigma}$  be isometrically realized 
 by a piecewise smooth $(n-1)$-dimensional   subvariety in a $(N-1)$-dimensional sphere, $N>>n$, that serves as the boundary of   the $N$-dimensional hemisphere, 
 $$Y \subset S^{N-1}(R)=\partial S^{N}_+(R),$$  where the radius of the sphere satisfies,
  $$R\geq\sqrt {\frac{ (n-1)(n-2)}{\sigma}}$$
and  where "isometrically" means preservation of the lengths of piecewise smooth curves in $Y$. 
  
 Then   the {\it spherical cone } of $Y$, that is the union of the geodesic segments which the center of 
 the spherical $n$-ball    $S^{N}_+ \subset S^N$  to all $y\in Y$ is, by definition, belongs to $\mathscr S^{n}_{\geq\sigma'}$ for 
 $$\sigma'=\sigma\frac {n}{n-2}$$
 and, more generally, a piecewise smooth $Y$ is in  $\mathscr S^{n}_{\geq\sigma'}$ if its scalar curvature at all non-singular points is $\geq \sigma'$ and near singularities $Y$ is isometric to a spherical cone over a space from $\mathscr S^{n-1}_{\geq\sigma}$.

To conclude the definition, we  agree to start the induction with $n-1=1$, where our  admissible spaces are 
circles of length $\leq 2\pi$ and, if we allow boundaries, segments of any length.

$Y \subset S^{N-1}$ be a closed  submanifold of dimension ${n-1} \geq 2$,   and let $S(Y)\subset S^{N} \supset S^{N-1}$ be the {\it spherical  suspension} of $Y$, that is the union of the geodesic segments which go from the north and the south poles of  $ S^N$ to $Y$.

Notice that this $S(Y)$ with the induced Riemannian metric is smooth away from the poles, where it is singular unless   the induced Riemannian metric in $Y$  has constant sectional curvature +1 and $Y$ is simply connected (hence, isometric to $S^{n-1}$).

\vspace{1mm}
 
{\sl Let $Y$ be a space from  $\mathscr S^n_{\geq\sigma}$  with $k$  isolated singular points   $y_i\in Y$ where $X$ is locally isometric to  S-cones 
over $(n-1)$-manifolds, call them $V_i$, $i=1,...k$ such that every such $V_i$ 
 bounds a Riemannian manifold $W_i$, where 
$Sc(W_i)>0$ and the mean curvature of $V_i=\partial W_i$ is positive. Then 
$$Sc^{\sf max}_{prop}(Y)\geq \sigma.$$}

{\it Sketch of the  Proof.} Arguing as in [GL(classification) 1980], one can, for all  $\varepsilon>0$, deform the metric in $X$ near  singularities  
keeping $Sc\geq \sigma -\varepsilon $,   such that the resulting metric on $Y$ minus the singular points $y_i$ becomes complete, where its   $k$ ends are isometric  to the cylinders $\varepsilon V_i\times [\infty)$,  where $\varepsilon V$ stands for an $V$  with  its Riemannian metric multiplied by $\varepsilon^2$.

This complete manifold, call it $Y_\varepsilon$, admits a locally constant at infinity  1-Lipschitz map 
  $Y_\varepsilon\to Y$ of degree 1, and then the closed  manifold $\bar Y_\varepsilon$,   obtained  from  
  $Y_\varepsilon$ by attaching $\varepsilon W_i$ to $\varepsilon V_i\times \{t_i\}$,   for large $t_i\in[0,\infty]$ admits a required 1-Lipschitz map to $Y$ as well. QED

{\it Remark.} Instead of filling  $V_i$ by $W_i$ {\it individually} it is sufficient to fill  in their (correctly oriented!) disjoint union $V= \bigsqcup_i V_i$ by $W$.
For instance,  if there are only two singular points, where $V_1$ and $V_2$ are isometric and admit orientation reversing isometries then $V_1\sqcup -V_2$ bounds the cylinder $W$ between them.

This kinds of "desingularization by surgery" also applies to $Y$, where the  singular loci $\Sigma\subset Y$  have  dimensions $dim(\Sigma)\geq 1$, similarly to  how it is done to manifolds with corners (see section 1.1 in [G(billiard0 2014]) but the filling condition becomes less manageable.

In fact even if $dim(\Sigma)=0$,  it is unclear how essential our  filling truly is,   especially for evaluation $Sc^{\sf max} $ of a {\it multiple} of the fundamental class of an $Y$; 
yet, the  spaces  $Y\in \mathscr S^n_{\geq\sigma}$ with isolated singularities  seem to enjoy the same  metric properties as smooth manifolds with $Sc\geq \sigma$ filling or no filling.

For instance, if the non-singular locus of such an $Y$ is spin then 
the hyperspherical radius  $Y$ is  bounded in the same way as it is for smooth manifolds: 
$$Rad_{S^n}(Y)\leq \sqrt \frac {n(n-1)}{\sigma},$$
 as it follows from Llarull's theorem for complete manifolds.
 
 In fact,  the construction from    [GL(classification) 1980] for  connected sums of manifolds with $Sc>0$, when   applied to $Y\setminus \Sigma$,   achieves  a blow-up of the metric $g$ of $Y$ on  $Y\setminus \Sigma$  to a complete one, say  $g_+$, such that $g_+\geq g$  and   
 $\inf_xSc(g_+,x)\geq \inf_xSc(g,x) -\varepsilon$ for an arbitrarily small $\varepsilon>0$. 
 
Also mean convex  cubical domains $U$  in $Y$ with none of  the  singular $y_i\in Y$ lying on the   boundary $\partial U$ satisfy the constraints  on the dihedral angles similar to those for smooth Riemannian manifolds with $Sc\geq \sigma$

But the picture becomes   less transparent for $dim(\Sigma)>0$, as it is exemplified by the following.
\vspace {1mm}

{\it Question}. {\sf Does the inequality $Rad^2_{S^n}(Y)\leq const_n\frac {n(n-1)}{\sigma}$ hold true for all 
 $Y\in \mathscr S^n_{\geq\sigma}$?}
\vspace {1mm}

{\it Perspective.} In view of  [Cheeger(singular) 1983],  [GSh(Riemann-Roch) 1993]  and[AlbGell(Dirac operator on pseudomanifolds) 2017 ], it is tempting to use the Dirac operator on the non-singular locus $Y\setminus \Sigma$ with a controlled behavior for $y\to \Sigma$, but it remains unclear if one can  actually  make this work for $dim(\Sigma)>0$. 

The only realistic approach at the present moment  is offered by the method of minimal hypersurfaces
 (and/or of stable $\mu$-bubbles), which may be additionally  aided by   surgery  desingularization,  such as   multi-doubling
 similar to that described in   [G(billiards) 2014] for manifolds with corners.

\vspace {1mm}

{\it Max-Scalar  Curvature Defined via {\sf Sc}-Normalized Manifolds .} Given a Riemannian  manifold $X=(X,g)$  with positive scalar 
curvature, let   $ g_\sim=Sc(g)\cdot g$, consider   Lipschitz  maps $f$   of  closed oriented Riemannian manifolds $X=(X,g)$ with $Sc(X)>0$  to $P$, such that $f_\ast[X]=h$, for a given $h\in H_n(P)$,  let 
$\lambda_\sim^{min}$ be the infimum of the Lipschitz constants of these  maps with respect to the metrics $g_\sim$
and let  $$Sc_\sim^{\sf max }(h)= \frac {1}{(\lambda_\sim^{min})^2}. $$  

And if $P$ is a {\it a piecewise  smooth polyhedral space}   (e.g.  a Riemannian  manifold), define   $Sc_{\wedge^2_\sim} ^{\sf max }(h)$ by taking the infimum 
 $\inf_f\sup_{x\in X}||\wedge ^2df(x) ||$ instead of the $\lambda_\sim^{min}$ (as in {\color {blue} $\wedge^2$-inequality} from  
section 4.2\footnote{The definition of $||\wedge ^2df(x) ||$ makes sense   for Lipschitz maps (at almost all $x$) but the arguments with  Dirac operators need smoothness of the maps. But it may be interesting to go beyond smooth manifolds and maps to  general continues maps with {\it  bounded area dilations}, where, probably, the most adequate definition of "area" in non-smooth metric spaces $P$  is the Hilbertian one in the sense of [G(Hilbert)  2012].}): 
$$Sc_{\wedge_\sim^2} ^{\sf max }(h)=\frac {1}{\inf_f\sup_{x\in X}||\wedge ^2df(x) ||}.$$
Clearly, 
$$Sc^{\sf max} \leq Sc_\sim^{\sf max}\leq Sc_{\wedge_\sim^2}^{\sf max}.$$

(Similar inequalities are satisfied by the spin and  by proper versions of $Sc^{\sf max}$), where 
{\it most bounds on    $Sc^{\sf max}$} we prove and/or conjecture below \vspace {1mm}
{\it can be} {\sf more or less automatically} {\it sharpened to their   $Sc_\sim^{\sf max}\l$  and   $  Sc_{\wedge_\sim^2}^{\sf max}$} ({\sf as well as to their  spin and   proper}) {\it  counterparts}.)

\vspace {1mm}\vspace {1mm}

{\it \large  Problem.} {\sf Evaluate $Sc^{\sf max}_{prop}$ of (the fundamental classes of)  "simple" metric space,  such as products of $m_i$-dimensional balls of radii $a_i$ where $\sum _im_i=n$  and the product  distance is $l_p$, i.e.
$dist_{l_p}( (x_i), (y_i))=\sqrt [p]{\sum_idist(x_i, y_i)^p}$, e.g. for $p=2$.}

\vspace {1mm}

This is related to  the problem of a  general nature of evaluating $Sc^{\sf max}(h_1\otimes h_2)$ of 
$h_1\otimes h_2\in H_{n_1+n_2}(P_1\times P_2)$ in terms of $Sc^{\sf max}(h_1)\in H_{n_1}(P_1)$ and $Sc^{\sf max}(h_2)\in H_{n_2}(P_2)$.

 It follows from the additivity of the scalar curvature (see section 1) that 
 $$Sc^{\sf max}(h_1\otimes h_2)\geq  Sc^{\sf max}(h_1)+ Sc^{\sf max}(h_2),$$
 but it is  unrealistic (?) to expect that, in general 
 $$Sc^{\sf max}(h_1\otimes h_2)\leq  const_{n_1+n_2}
\cdot(Sc^{\sf max}(h_1)+ Sc^{\sf max}(h_2)),$$
albeit
the geometric method  from the   section 5.4 does deliver non-trivial  bounds on  $Sc_{prop}^{\sf max}$ of {\it product spaces} whenever  lower bounds on  the {\it hyperspherical radii of the  factors} are available.\footnote {One may define   $Rad_{S^n}(h)$, $h\in H^n(P)$,  as the suprema of the radii  $R$ of the $n$-spheres, for which 
$P$ admits a 1-Lipschitz  map $f:P\to S^n(R)$, such that $f_\ast(h)\neq 0$.}

%%%%%%%%%%%%%%%%%%%%%

\subsection {\color {blue} Extremality and Rigidity of $\log$-Concave Warped products}\label {log-concave5}
   %%%%%%%%%%%%%%%%

The inequalities proven in  section 5.3 say, in effect, that 
the metric  
$$\mbox {$g_\phi=\phi^2g_{flat} + dt^2 $ 
on $\mathbb T^{n-1}\times\mathbb R$ for  $\phi(t)= \exp \int_{-\pi/n}^t -\tan \frac {nt}{2} dt$}$$ 
is  {\it extremal}:\vspace{1mm}
{\it \color{blue} one can't increase    $g_\phi$ without decreasing  its scalar curvature,}\footnote{To be precise, one should say that 

{\sf one can't modify the metric, such that the {\it scalar curvature increases} but the metric itself {\it doesn't decrease}.} 

The relevance of this formulation is seen in the example of $X=S^n\times S^1$, where one can stretch the obvious product metric $g$  in the $S^1$-direction without changing the scalar curvature, but one {\it can't increase} the scalar curvature by deformations that increase $g$.}\vspace{0mm}

 \hspace{-5.5mm}where the essential feature   of   $\phi$
  
   (implicitly) 
 used for this purpose was 
 {\it $\log$-concavity of  $\phi$}:
$$\frac {d^2\log\phi(t)}{dt^2}<0.$$

We show in this section that the same kind of extremality (accompanied by rigidity) holds for other $\log$-concave functions, notably for $\varphi(t)=t^2$,  $\varphi(t)=\sin t$ and $\varphi(t)=\sinh t$ which results in 

\hspace {7mm}{\it \color {blue}  rigidity of  punctured  Euclidean, spherical and hyperbolic spaces.}
\vspace{1mm}

More generally,  let $X=Y\times \mathbb R$ comes with the warped product metric
$g_\phi  =\phi^2 dg_y+dt^2$. Then the  mean curvatures   of the hypersurfaces $Y_t=Y\times \{t\}$, $t\in \mathbb R$,  satisfy (see 2.4)
$$mean.curv (Y_t)=\mu(t)=(n-1)\frac {d\log \phi(t)}{dt}=\frac {\phi'(t)}{\phi(t)},$$
and, obviously, are these $Y_t\subset X$ are  locally (non-strictly) minimizing $\mu$-bubbles. \footnote{If $Y$ is non-compact, the minimization is understood here for variations with compact supports.}

Now, clearly, {\it $\phi$ is log-concave, if and only if 
$$ \frac {d\mu}{dt}= -\left | \frac {d\mu}{dt}\right |.$$}
Thus,  $R_+$ defined (see section 5) as 
$$R_+(x)=\frac {n\mu(x)^2}{n-1} -2||d\mu(x)|| +Sc (X,x)$$
is equal in the present case to 
 $$\frac {n\mu(t)^2}{n-1} +2\mu'(t) +Sc (g_\phi(t))= \frac {2(n-1) \phi''(t)}{\phi(t)} +(n-1)(n-2)\left(\frac {\phi'}{\phi}\right)^2+Sc(g_\phi(t)) $$

which implies (see section  5)
that $$(R_+)_{Y_t}=\frac {1}{\phi^2}Sc(g_{Y_t})=Sc(g_{Y_t})\mbox { for } g_{Y_t}=\phi^2g_Y.$$

Here  our  $-\Delta_{Y_t} +\frac {1}{2}Sc(g_{Y_t}) -(R_+)_{Y_t}$  from section \ref{variation5}) is
equal to  $-\Delta_{Y_t} $,  the  lowest eigenvalue  which is
  {\it zero} with  {\it constant}  corresponding eigenfunctions and the corresponding ($\mathbb T^1$-invariant warped product) metrics on $Y_t\times \mathbb T^{1}$ are  (non-warped)   $g_{Y_t}+dt^2$ for $Y_t =Y\times \{t\} \subset X=Y\times \mathbb R$ and all $t\in \mathbb R$.

\vspace{1mm}

This computation together with {\color {blue}\small\PlusCenterOpen$_{warp}$} in section \ref {variation5}
yield the following.
\vspace{1mm}

{\it \large \color {blue}  Comparison Lemma.} {\sf Let $\underline X=\underline Y \times [\underline a, \underline  b]$ be an $\underline n$-dimensional   warped product manifold with the metric  

$$g_{\underline X}= g_{\underline \phi}=  \underline \phi^2g_{\underline Y}+dt^2,\mbox { }\underline t\in  [\underline a, \underline b],$$
 where  $\underline \phi(\underline t)$ is a smooth positive  log-{\it concave} function on the segment $[\underline a,\underline b]$.

Let $X$ be an $n$-dimensional Riemannian  manifold, 
with a smooth function $\mu(x)$ on it  and let $Y =Y_{\mu}\subset X$ be a stable, e.g. locally minimising 
$\mu$-bubble in $X$.

Let $g^{\rtimes } =g_{\phi_\rtimes} =\phi_\rtimes ^2g_{Y}+ dt^2$
be the metric on $Y\times \mathbb T^1$ where $g_Y$ is the metric on $Y$ induced from $X$, and where $\phi_\rtimes $ is the first eigenfunction of the  
$$\mbox {$  -\Delta +\frac {1}{2}Sc(g_Y,y) -R_+(y)$ for $R_+(x)=\frac {n\mu(x)^2}{n-1} -2||d\mu(x)|| +Sc (X,x)$}$$
 ($\phi_\rtimes $ is not assumed positive at this point).

Let $f:X\to \underline X$ be a smooth map let $f_{\underline Y}:X\to \underline Y$ denote 
the  $\underline Y$-component of $f$, that is the composition of $f$ with the projection 
$\underline X=\underline Y \times [\underline a,\underline b]\to \underline Y$.

Let 
$$f_{[\underline a,\underline b]}: X\to [\underline a, \underline b]$$ 
be the $[\underline a,\underline b]$-component of $f$, let
$$\underline \mu^\ast(x) =\underline \mu\circ f_{[\underline a,\underline b]}(x)\mbox  { for  } \underline     
 \mu(\underline t)=(\underline  n-1)\frac {d\log \underline \phi(\underline t)}{d\underline t}=mean.curv(\underline      Y_{\underline t}), \mbox { } \underline t=f_{[\underline a, \underline b]}(x)$$
and let 
$$\underline \mu'^\ast=\underline \mu'\circ  f_{[\underline a,\underline  b]}(x)\mbox { where } \mu'=\mu'(\underline t)=\frac {d\underline\mu(\underline t)}{d\underline t}. $$

Let 
$$\underline R_+^\ast(x)= \frac {n\underline \mu^\ast (x)^2}{\underline n-1} -2||d\underline \mu^\ast (x)|| +Sc (\underline X,f(x))$$

{\it If  $$  R_+(x) \geq \underline R_+^\ast(x),$$ 
then the function $\phi_\rtimes$ is positive and the scalar curvature of the metric $ g^{\rtimes}=g_{\phi_\rtimes}$ on $Y \times \mathbb  T ^1$ satisfies 
$$Sc_{g^{\rtimes}}(y,t)   \geq  \frac {1}{ ||df_{[\underline a,\underline b]}(y)||^2}Sc(\underline Y, f_{\underline Y}(y))=  Sc(\underline Y_{\underline t}, f(y))\mbox { for $\underline Y_{\underline t}\ni f(y)$}.$$ }

\vspace{1mm}

{\large \it The  main case} of this lemma, which  we use below, is where\vspace{1mm}

({\it $\bullet_{df_{[\underline a,\underline b]}}$})  {\it \hspace{1 mm}the function $f_{[\underline a,\underline b]}:X\to [\underline a,\underline b]$ is 1-Lipschitz, i.e. $||df_{[\underline a,\underline b]}||\leq 1$,

\hspace{-6 mm}}and 

({\it $\bullet_{\mu}$})  {\it  \hspace{1 mm}$\mu(x)=\underline \mu\circ f_{[\underline a,\underline b]} $, that is $\mu(x)=mean.curv (\underline Y_{\underline t}, f(x))$ 
for $\underline Y_{\underline t}\ni f(x)$}}

and where the conclusion reads:
$$Sc_{g_\phi}(y,t )   \geq  \frac {1}{ (f_{[\underline a, \underline b]}(y))^2}Sc(\underline Y, f_{\underline Y}(y))+Sc(X,y)-Sc(\underline X,f(y)).\leqno {\color{blue} [Sc\geq]}.$$ \vspace{1mm}

{\large \it Corollary.} {\sf Let $X^{\rtimes }$ denote the above Riemannian (warped product) manifold  
 $(Y\times \mathbb T^1, g^\rtimes=g_{\phi_\rtimes})$ 
and let $ f_\rtimes : X^{\rtimes }\to \underline Y$  be defined by $(y,t)\mapsto f_{\underline  Y}(y).$}

{\it If besides {\it $\bullet_{df_{[a,b]}}$ and  $\bullet_\mu$,
$$ || \wedge^2df||\leq 1,\mbox { e.g. } ||df||\leq 1$$
and if 
$$Sc(X,y)\geq Sc(\underline X,f(y)),$$
then the map $ f_\rtimes $ satisfies 
$$ Sc(X^{\rtimes },  x_\rtimes )\geq ||df_\rtimes ||^2 Sc(\underline  Y,  f_\rtimes  (x_\rtimes))\geq || \wedge^2d f_\circ||Sc(\underline  Y, f_\rtimes (x_\rtimes )).$$}}

Now, the existence of minimal bubbles under the barrier  {\color {blue}$[\gtrless mean=\mp\infty]$}-condition (see section 5.2) and 
 a  combination of the above  with the Llarull 
 {\color {blue} $trace \wedge^2df$-inequality} from section 4.2 yields the following.\vspace {1mm}

{\Large \color {blue} $\odot_{S^n}$.} {\it \large Extremality of Doubly  Punctured Spheres.} {\sf Let $X$ be an oriented  Riemannian spin $n$-manifold, let $\underline X$ be the  $n$-sphere with two opposite points removed
and let $f:X\to  \underline X$ be  a {\it smooth 1-Lipschitz} map of {\it non-zero} degree.}

{\it If $Sc(X)\geq n(n-1)=Sc(\underline X)=Sc(S^n),$ then

(A) the scalar curvature of $X$ is constant $=n(n-1)$;

(B)  the map $f$ is an isometry.} 

\vspace{1mm}

{\it Proof.}  The spherical metric on $ \underline X=S^n\setminus \{s,-s\}$ is the  warped 
product $S^{n-1}\times \left (-\frac {\pi}{2},\frac {\pi}{2}\right )$ where the warping factor
$\underline \phi(t)=\cos t$ which  is logarithmically concave, where 
$\underline \mu(t)=\frac {d\log \underline  \phi(t)}{dt}  \to \pm\infty $  for $t\to \mp\frac {\pi}{2}$. \footnote{If a log-concave function $\phi$ on the segment $[-l,l]$ is positive for $-l <t<l $  and it  vanishes at $-l$,  then the  logarithmic derivative of $\phi$ goes to $\infty$  for $t\to -l$;  similarly, $$\frac {\phi'}{\phi}\underset {t\to l}\to -\infty,$$ if $\phi$ vanishes at $t=l$.}

This implies  (A) while  (B) needs a little extra (rigidity) argument indicated in section 5.7.

\vspace{1mm}

1-{\it Lipschitz Remark.} As it is clear from the proof, the 1-Lipshitz condition can be relaxed to the following one.

{\sf The radial component $f_{\left [-\frac {\pi}{2},\frac {\pi}{2}\right]}:X\to \left [-\frac {\pi}{2},\frac {\pi}{2}\right]$ of $f$, which corresponds to the signed  distance function from the equator  in
$S^n\setminus \{s,-s\}$ is 1-Lipschitz and (the exterior square of) the  differential of the $S^{n-1}$ component
$f_{S^{n-1}}: X\to S^{n-1}$ satisfies}

$$ \wedge^2df_{S^{n-1}}\leq\frac {1}{ \left (\cos f_{\left [-\frac {\pi}{2},\frac {\pi}{2}\right]}(x)\right)^2}.$$

\vspace{1mm}
{\it Non-Spin Remark. }  If $n=4$, one can drop the spin condition,
 since  $\mu$-bubbles $Y\in X$, being  3-manifolds, are spin. 
\vspace{1mm}

Similarly to {\Large \color {blue} $\odot_{S^n}$} one shows the following.\vspace{1mm}

{\Large \color {blue} $\odot_{\mathbb R^n}.$} {\sf Let  Let $X$ be as above, let $\underline X$ be $\mathbb R^n$ with a point  removed
and let $f:X\to  \underline X$ be  a {\it smooth 1-Lipschitz} map of {\it non-zero} degree.

{\it If $Sc(X)\geq n(n-1)\geq 0$ and if $X$ is an isometry at infinity, then

(A) $Sc(X)=0$;

(B)  the map $f$ is an isometry.} \vspace{1mm}

{\Large \color {blue} $\odot_{\mathbf H^n}.$} Let  Let $X$ be as above, let $\underline X$ be the hyperbolic space with a point  removed
and let $f:X\to  \underline X$ be  a {\it smooth 1-Lipschitz} map of {\it non-zero} degree.

{\it If $Sc(X)\geq  -n(n-1)$ and if $X$ is an isometry at infinity, then

(A) $Sc(X)=-n(n-1)$;

(B)  the map $f$ is an isometry.} }

 \vspace{1mm}
 
 {\it Question.} 
 Let $d_0(\underline x) =dist (\underline x, \underline x_0)$    be the distance function in $\underline X$ (used in  {\Large \color {blue} $\odot_{\mathbb R^n}$}  and/or in   {\Large \color {blue} $\odot_{\mathbf H^n}$}) to the point $\underline x_0$, which was removed from  $\mathbb R^n$ or from $\mathbf H^n$, and let $d_f(x)= d_0(f(x))$.

 {\sf Can one relax the 1-Lipschitz condition in the propositions   {\Large \color {blue} $\odot_{\mathbb R^n}$}  and in {\Large \color {blue} $\odot_{\mathbf H^n}$}  by requiring that not $f$ but only 
 the function $d_f(x)$ is $1$-Lipschitz?}

\vspace {1mm}

We conclude this section with the following proposition, that  is proven (in a different form) in 
[Richard(2-systoles) 2020] (compare [Zhu(rigidity) 2019])  and   which  provides  a useful geometric  information on manifolds with scalar curvature $\geq \sigma>0$  on the scale $\sim\frac {1}{\sqrt \sigma}$.

{\it  \color {blue} Richard{'}s  Lemma.} {\sf Let $X$ be an   oriented $m$-dimensional Riemannian   manifold (possibly non-compact and non-complete)  with compact  boundary 
and $X_0\subset X$  be an open subset  with smooth boundary  such that the complement $X\setminus X_0$ is compact.
Let $h\in H_{m-2}( \partial X)$ and $h_0\in H_{m-2}(X_0)$ be homology classes, which have {\it equal  images}  under  the  homomorphisms      induced by the inclusions\ $\partial X\hookrightarrow X \hookleftarrow X_0$,  
that are
$$h\in H_{m-2}(\partial X)\rightarrow   H_{m-2}(\partial X) \leftarrow H_{m-2}(X_0)\ni h_0.$$}

{\it Let   $$Sc(X)\geq \sigma>0, $$ 
and
$$ dist^2 (X_0, \partial X)  \geq \frac  {m(m-1)\pi^2}{\sigma}.$$}

 Then (we can  vouch 100\%  here, as everywhere in this text, only for $n\leq 8$)

{\it the image of the homology class $h$ in $H_{m-2}(X)$ can be realized by a  closed smooth $(m-2)$-dimensional submanifold $Y\subset X$, on which  there exists a  smooth positive  function $\phi(y)$,  such that  the metric $g_\ast=dy^2+\phi(y)^2dt^2$ on the product manifold  $Y\times \mathbb R^2$
satisfies
$$Sc(g_\ast)\geq \frac  {m-2}{m}\sigma, $$}  
where $dy^2$  denotes the Riemannian metric on $Y$ induced  from $X\supset Y$ and $dt^2$ is the Euclidean metric in the plane $\mathbb R^2$.
\vspace{1mm}

{\it Proof.}   Use the     codimension 2 argument   as in the proof of the {\it quadratic decay theorem} in section 1 in [G(inequalities) 2019] (see  also  section 7 in  {[GL(complete) 1983]   and \S $9\frac {3}{11}$ in  [G(positive) 1996]) together with 
 the above comparison lemma  combined   with the equivariant  separation theorem from section \ref{separating5}.

\vspace{1mm}

(A version of Richard{'}s lemma is also established in  [Chodosh-Li(bubbles) 2020] in the course of their proof of non-existence of metric with $Sc>0$ on aspherical 4- and 5-manifolds; also, this lemma is used for a similar purpose in   [G(aspherical) 2020].)

%%%%%%%%%%%%%%%%%%%%%%%%%%%%%%%

 \subsection {\color {blue} On
  Extremality of   Warped Products of Manifolds with Boundaries and with Corners} \label {warped boundary5}
 
%%%%%%%%%%%%%%%%%%%%%%%%%%

We explained in section 4.4  how reflection+ smoothing  allows an extension of the Llarull and Goette-Semmelmann theorems from section  4.2 to manifolds  with smooth boundaries and to a class of manifolds with corners. This, combined with the above,  enlarges the class of manifolds with corners to which the conclusion of the extremality  {\color {blue} \DiamondSolid$\angle_{ij}$}  {\large \it theorem}    applies.
example.

Let $\triangle^{n-1}\subset S^{n-1}$  be the regular spherical simplex with flat faces and  the   
dihedral angles $\frac {\pi}{2}$ and let ${\sf S^\ast_\ast}\triangle^{n-1}\subset S^n\subset S^{n-1}$ be the spherical  suspension of  $\triangle^{n-1}$ and let  $\underline X={\sf S}_{a}^{b}(\triangle^{n-1})\subset  {\sf S^\ast_\ast}\triangle^{n-1}$, $a,b \in (-\frac{\pi}{2}, \frac{\pi}{2})$,
be the region of  $\sf S^\ast_\ast \triangle^{n-1}$ between a pair of $(n-1)$-spheres concentric to our equatorial 
$S^{n-1}\subset S^{n}$.\vspace {1mm}

{\sf Let $X$ be an $n$-dimensional orientable Riemannian spin manifold with corners and let 
$f:X\to \underline X$ be a smooth  1-Lipschitz map which respects to the corner structure and which has non-zero degree.}  \vspace {1mm} 

{ \color {blue} \it \large Spherical ${\sf S}_{a}^{b}(\triangle)$-Inequality.}  {\it If $Sc(X)\geq Sc(\underline  X)=n(n-1)$, if all $(n-1)$-faces $F_i\subset \partial X$  have their mean curvatures bounded from below by those of the corresponding faces in $\underline X$, \footnote{All these but two have zero mean curvatures.}
$$mean.curv(F_i)\geq mean.curv (\underline F_i),$$
and if all dihedral angle of $X$ are bounded by the corresponding ones of $\underline X$,
$$\angle_{ij}\leq \underline \angle_{ij}=\frac {\pi}{2},$$
then 
$$Sc(X) = n(n-1),$$
$$mean.curv(F_i)= mean.curv (\underline F_i)$$
and 
$$\angle_{ij} =\frac {\pi}{2}. $$}\vspace{1mm}

{\it Exercise.}  (a)  Recall  {{ {\color {blue}  $\blacksquare$}}-hyperbolic comparison theorem for cubical manifolds  diffeomorphic to 
$$\underline V=[0,1]\rtimes [0,1]^{n-1}\subset \mathbb H^n=(\mathbb R^1\times\mathbb R^{n-1},  dt^2+e^{2t}dx^2)$$ 
from section \ref{Sc-criteria3} and generalize it to all compact cubical manifolds $V$ (to be sure, of dimension $n\leq 8$).

(b)  Formulate and prove (for $n\leq 8$)  the Euclidean and hyperbolic versions of the  
{\color {blue} ${\sf S}_{a}^{b}(\triangle)$}-inequality for spin manifolds $V$ with corners.\footnote {See [Li(parabolic) 2020) for further  results in this direction.}

\vspace {2mm}

{\it \color {blue} Question.} Do  the counterparts to the  ${\sf S}_{a}^{b}(\triangle)$-inequality hold for other simplices 
and polyhedra?
%%%%%%%%%%%%%%
 
 \subsection {\color {blue}On Rigidity of  Extremal Warped Products}  \label {rigidity5}

%%%%%%%%%%%%%%%%%%%
Let us explain,  as a matter of example, that \vspace {1mm}

\hspace {10mm} {\sl doubly punctured sphere 
$\underline X=S^n\setminus\{\pm s\}$  is spin-rigid.} \vspace {1mm}

This means   that\vspace {1mm}

{\it  if an oriented Riemannian  spin $n$-manifold     $X$ with $Sc(X)\geq n(n-1)=Sc(\underline X=Sc(S^n)$  admits a smooth proper 1-Lipschitz map  $f:X\to \underline X$
such that $deg(f)\neq 0$, then, in fact, such an  $f$ is an isometry.}
\vspace {1mm}

{\it Proof.} We know (see the  the proof of  {\Large \color {blue} $\odot_{S^n}$} in \ref{log-concave5}) that $X$ contains a minimal $\mu$-bubble $Y$, which separates the two (union of) ends of $X$, where
$\mu(x)$ is the $f$-pullback of the mean curvature function of the concentric $(n-1)$-spheres  in
$\underline X=S^n\setminus\{\pm s\}$ between the two punctures and that this $m$-bubble must be umbilic, where we assume at this point that $Y$ is non-singular, e.g. $n\leq 7$.

What we want  to prove now is that these bubbles {\it foliate all of} $X$, namely they come in a  continuous family of mutually disjoint minimal $\mu$-bubbles  $Y_t$, $ t\in\left( -\frac {\pi}{2}, \frac{\pi}{2}\right)$, which together  cover $X$. 

Indeed, if the {\it maximal} such family $Y_t$ wouldn't cover $X$, then the  would exists a small perturbation $\mu'(x)$ of $\mu(x)$ in the  gap between two $Y_t$ in the maximal family, such that 
$|\mu'|>|\mu|$ in this gap, while $||d\mu'||= ||d\mu||$ in there and such that there would exist a minimal $\mu'$-bubble $Y'$ in this gap. 

But  then, by calculation in section  \ref{log-concave5}, the resulting  warped product metric on $Y'\times S^1$ would have  $Sc>n(n-1)$, thus proving "no gap property" by contradiction.

Therefore, $X$ itself is the warped product, $X=Y\times \left( -\frac {\pi}{2}, \frac {\pi}{2}\right)$  with the  metric $dt^2=(\sin t)^2 g_Y$, where $Sc(g_Y)=n(n-1)$ and which by Llarull's rigidity theorem, has constant sectional curvature. QED. 

\vspace {1mm}

{\it Remarks} (a) On the positive side,  this argument is quite robust, which makes it  compatible with approximation of bubble and metrics. For instance it nicely works  for  $n=8$ in conjunction with Smale's generic regularity theorem and, probably, for all $n$ with Lohkamp's smoothing theorem.

But it is not quite  clear how to make this work for  non-smooth limits of smooth 
metrics.

For instance, 

{\sf let $g_i$ be a sequence  of Riemannian metrics on the torus $\mathbb T^n$ , such that 
$$Sc(g_i)\geq -\varepsilon_i\underset {i\to \infty} \to 0$$  
and such that $g_i$ uniformly   converge to a continuous metric $g$. 

{\sl Is this $g$, say for $n\leq 7$,  Riemannian flat?}\footnote {Yes, according to 
[Burkhart-Guim(regularizing Ricci flow) 2019].} \vspace {1mm}}

(The above argument shows that,  given an  indivisible $(n-1)$-homology  class in  $\mathbb T^n$, there exists  a  foliation of $\mathbb T^n$  by $g$-minimal submanifolds from this class.
 But it is not immediately clear how to show  that these
submanifolds are totally geodesic.)

{\sf let $g_i$ be a sequence  of Riemannian metrics on the torus $\mathbb T^n$ , such that 
$$Sc(g_i)\geq -\varepsilon_i\underset {i\to \infty} \to 0$$  
and such that $g_i$ uniformly   converge to a continuous metric $g$. 

The above argument shows that,  given an  indivisible $(n-1)$-homology  class in  $\mathbb T^n$, there exists  a  foliation of $\mathbb T^n$  by $g$-minimal submanifolds from this class, 
 but it is not  clear how to show  that these
submanifolds are totally geodesic, that is needed   for the proof of  flatness of $g$,

Yet, 

{\it {\color{blue!50!black}  the Ricci flow argument from 
[Burkhart-Guim(regularizing Ricci flow) 2019] 

does show    the metric  $g$ }  is  flat. }

%%%%%%%%%%%%%%%%%%%%%%%%%%%%%%c
\subsection {\color {blue}
Capillary Surfaces:  $\mu$-Bubbles with Measures $\mu_\partial$ Supported on  Boundaries}\label {capillary5}
%%%%%%%%%%%%%%%%%%%%%%%%%
In order to extend extremality and other  results to more general manifolds with boundaries,  such, e.g. as  conical domains in $\mathbb R^n$, one shouldn't  limit oneself to 
the definition of a  $\mu$ bubble from  section 5,  where the admissible 
measures $\mu$   on $X$  are of the form $\mu(x)dx$ for {\it continuous functions} $\mu(x)$.

In fact the definition of $\mu$-bubbles  makes sense for more general  measures, where
 a  geometrically  interesting case  is that of a manifold  $X$ with boundary, here denoted $S=\partial X$,
and   our measure is of the form $\mu_\bullet(x)dx+\mu_\partial(c)ds$,  where $\mu_\bullet$   and 
$\mu_\partial $ are  continuous (or measurable)  functions on $X$ and on $S=\partial X,$
and where we let
$$|\mu_\partial|<1$$ 
for a  reason  that becomes clear later on.

Let $\cal Y$ be the set of    cooriented hypersurfaces $Y\subset X$ with boundaries contained in $S$,
 $$Z=\partial Y\subset S=\partial X,$$
 where the unit    field   normal to $Y$, which defines the coorientation is called 
the {\it upward} field  and denoted $\nu=\nu_{Y\uparrow}$  and let 
$$\mu=\mu_\bullet(x)dx+\mu_\partial(s)ds.$$

Then  a hypersurface  $Y\in \mathcal Y $ is called a {\it $\mu$-bubble} (compare 5.1), if it is
extremal or, at least,  stationary  for 
$$Y \mapsto  vol^{[-\mu]}_{n-1}(Y)=_{def} vol_{n-1}(Y )-  \mu(X_{<}),$$
where  $X_{<}\subset X$  is the region in $X$  "below" $ Y\subset X,$
where
 $$\mu(X_{<})=\int_{X_{<}}\mu_\bullet(X)dx+ \int_{S_{<}}\mu_\partial(s) ds$$
for $$S_<= S\cap X_<\subset S=\partial X.$$

This kind  of (2-dimensional) $Y$ for constant functions $\mu_\bullet$ and $\mu_\partial$ are called {\it capillary surfaces}.

  An essential for our  geometric  purposes feature of
    such surfaces and of $(\mu_\bullet+\mu_\partial$)-bubbles in general is a particular algebraic  property of the second  variation formula for  stationary $Y$,  similar to that for $\mu$-bubles with continuous $\mu(x)$ on  manifolds without boundary that is proved and used at  the beginning of section \ref  {bubbles5}; 
 the  derivation of this formula for capillary surfaces  was given in   [Ros-Souam(capillary) 1997]  and    used in 
  [Li(comparison) 2017] for  the proof of extremality of certain polyhedra.
  \footnote{Necessary existence and regularity of capillary hypersurfaces follow from  
[Simon-Spruck(capillary) 1976], [Gerhard(capillarity) 1976], [Liang(capillarity) 2005], [Philippis-Maggi(capillary) 2015]  as it is indicated in   [Li(comparison) 2017] and [Li(rigidity) 2019].}

In what follows we present a geometrically  transparent  derivation of this formula with an eye on  further  applications.
\footnote{The first version of this manuscript contained a  computational error  that lead to a most   disappointing  conclusion. I am thankful to  Mike Anderson who   encouraged me  to double check my computation.}    \vspace {1mm}

{\it\color {blue}    Example.} 
 Let $$ X=B^n\subset \mathbb R^n=\mathbb R^{n-1}\times \mathbb R, \mbox { }
  S=\partial B^n=S^{n-1},$$ be the unit ball, 
with the boundary sphere $S=S^{n-1}$ 
and let $Y_t= Y_t^{n-1}\subset B^n$, 
be  the horizontal discs, that are the 
intersections 
 $$Y_t =B^n\cap  (\mathbb R^{n-1}\times \{t\}),
 \mbox { } -1<t<1.$$ 

Let  
$$\angle_t= \angle_{Z_t}  (Y_t, S^{n-1})$$
 be the dihedral  angles between the hypersurfaces  $ Y_t$ and $S^{n-1}$ along their intersection 
$$  Z_t =\partial Y_t= Y_t\cap S^{n-1}$$
where, 
  we agree  that  this  angle is measured   "below" $Y_t$;  thus
 $\angle_{-1}=0$ and $\angle_1=\pi$, i.e.
 it is related to  the height $t=t(s)$ of $Z_t\subset S^{n-1}$ by 
  $$t=\cos (\angle_t-\pi/2).$$

Next,  let $\mu_\bullet=0$ and  let 
$$\underline \mu_\partial=\underline\mu_\partial(s)=cos\angle_{t(s)},$$
where $t$ is regarded here as the height function $t: S^{n-1}\to [-1,1]$ for 

\hspace {-6mm}$S^{n-1} \subset  \mathbb R^{n-1}\times [-1,1]\subset \mathbb R^{n-1}\times  \mathbb R$.
  
 Then the {\it normal derivative} $\partial _\nu=\frac {d}{dt}$     of the volume of the discs $Y_t\subset B^n$ is expressed in terms of 
 $$\mbox{$|Z_t|=vol_{n-2}(Z_t)$,  and the angle $\angle_t\in(0,\pi)$ }$$ 
 as follows
$$\partial_\nu vol_{n-1}( Y_t)=|Z_t|\cot \angle_t,$$

 while the derivative of the $ \underline\mu_\partial$-measure of the region $S_{\leq t}\subset S^{n-1}$  below $Z_t=\partial Y_t\subset S^{n-1}$  for   $\mu=\underline\mu_\partial(s)=cos \angle _{t(s)}$ is 
 $$\partial_\nu \underline\mu_\partial(S_{\leq t})) = \frac {|Z_t|\underline\mu_\partial(s)}{\sin \angle_t}= |Z_t|\cot \angle_t .$$
  Thus,
 
 {\sf the   derivatives of $vol_{n-1}(Y_t)$ and of  $ \underline\mu_\partial(S_{\leq t})$  by the field $\psi\nu$ for all   $C^1$-smooth functions 
$\psi=\psi(y)$, $y\in Y_t$ satisfy
$$\partial_{\psi\nu}vol_{n-1}(Y_t)=\partial_{\psi\nu}\underline\mu_\partial(S_{\leq t})$$
and 
$$\partial_{\psi\nu}vol^{[-\mu]}_{n-1}(Y)=0,$$}
which says that   

\hspace {20mm} {\sl $Y_t$ are $\mu$-bubbles for this  $\mu=\underline\mu_\partial(s)$,} 

\hspace {-6mm} since they {\it are}   stationary for the the functional $Y\mapsto vol^{[-\mu]}_{n-1}(Y)$.

\vspace{1mm}

{\it Exercise.}  Let  
$$Y_\rho\subset B^n(1)=\{x\in \mathbb R^n\}_{||x\leq1||}, \mbox { } 0<\rho<2,$$
be   the intersections of the   concentric $\rho$-spheres $S^{n-1}_{x_0}(\rho )  \subset \mathbb R^n$ around  $x_0=(0, 0,..,-1)\in \mathbb R^n$ with the unit ball $B^n\subset \mathbb R^n$.

Determine  the measure $\mu=\mu_\bullet(x)dx+ \mu_\partial(s)dy$, for which these $Y_\rho$  
serve as $\mu$-bubbles.

\vspace {1mm}

Let us return to the general Riemannian manifold  $X$  with boundary $S=\partial X$,  a hypersurface  $Y \subset X$,   where  $Z=\partial Y \subset S=\partial X$ and let $\mu$ 
be a  measure of the form $\mu=\mu_\bullet(x)dx+\mu_\partial(s)ds$
for  continuous  functions $\mu_\bullet (x)$ on $X$  and  $\mu_\partial(s)$ on $S$. 

\vspace {1mm} 

{\color {blue} \sf\large  First Variation Formula for  $vol^{[-\mu]}_{n-1}(Y)$.}  In order to define $\partial_\nu$ 
one needs to choose {\it a vector field} extending the "upward" normal field to $Y$, denoted,  as earlier,
 $\nu$,  {\it from $Y$ to a neighbourhood of} $Y$ in $X$, that we do as follows.

Smoothly extend $X$ beyond its boundary by a slightly greater Riemannian manifold $X_+$ of the same dimension and extend 
$Y$ by a hypersurface $Y_+\subset X_+$. 

Move $Y_+$ in both normal directions by distance $|t|$ to $Y_{+, t}\subset X_+$, $-\varepsilon\leq t\leq \varepsilon$ for a small $\varepsilon>0$  and let $Y_t\subset X$ be the intersection of the so moved $Y_+$ with $X\subset X_+$
$$Y_t=  Y_{+, t}\cap X\subset X,$$
where, observe, $Y_0=Y$, and where $Y_{\pm t}$ are $t$-equidistant hypersurfaces  to $Y_0$ in $X$, except, maybe, for the $|t|$-neighbourhood of $S=\partial X$, where 
we  agree that $Y_t$ with $t<0$ lies below $Y$ i.e. in the domain  $X_< \subset X$ and $Y_{t>0}$
are positioned over $Y$ in $X$.

Now $\partial_\nu$ is understood as $\frac{d}{dt}|_{t=0}$  and  $\partial_{\psi\nu}$ and 
$\partial^2_{\psi \nu}$ are understood accordingly.

Let $\angle_z\in (0,\pi)$  denote the angle between  (the tangent spaces of)   $Y$ and $S$ at $z\in Z=\partial Y=Y\cap S$ measured below $Y$, i.e. in $X_<$.

Then, clearly,  for all smooth functions $\psi=\psi(y)$, 
$$ \partial_{\psi\nu}vol_{n-1}(Y)= \int _Y  \psi(y)\cdot  mean.curv(Y,y) dy +\int_Z \psi(z) \frac{\cos \angle_z}{\sin \angle_z}  dz, $$
$$\partial_{\psi\nu} \mu_\bullet (Y)=\int_{Y}\psi(y)\cdot \mu_\bullet(y)dy$$
and 
$$\partial_{\psi\nu} \mu_\partial(S_-)= \int_Z\frac {\psi(z)} {\sin \angle_z} \mu_\partial(z)
dz,$$
where $\mu(S_<)$ stands for $\int_{S_<} \mu_\partial(s)ds$, where the "$Y$"-integrals are the  ones we met earlier in section 5  for $X$ without boundary and where the shape  of the  $Z$-integrals  can be seen  by looking at the above example.

\vspace {1mm}

Thus, the first variation {\color {blue}$\partial_{\psi \nu} vol^{[-\mu]}_{n-1}(Y)$} equals the sum of two integrals, one over $Y$ and the other one over $Z=\partial Y$, {\color {blue}  
$$ \int _Y  \psi(y)  (mean.curv(Y,y) -\mu_\bullet(y))dy
+\int_Z \psi(z)\left  (\frac{\cos \angle_z}{\sin \angle_z}   -  \frac {1} {\sin \angle_z} \mu_\partial(z)\right)
dz,$$}

Therefore, $\partial_{\psi \nu} vol^{[-\mu]}_{n-1}(Y)$ vanishes for all smooth functions $\psi(y)$ and 
{\it $Y$ is a (stationary)  $\mu$-bubble,  if and only if
{\color{blue}$$\mbox {$mean.curv(Y)=\mu_\bullet $ and $ \cos  \angle_z=\mu_\partial(z) $},\mbox { } z\in Z=\partial Y. $$}}

\vspace {1mm}

{\color {blue}$\partial Y$-Contribution to the  \sf\large Second Variation Formula for  $vol^{[-\mu]}_{n-1}(Y)$.} 
Let us compute the second derivative (variation)  $\partial_{\psi\nu}\partial_{\psi\nu}vol^{[-\mu]}_{n-1}(Y)$ on a {\it stationary}  $Y=Y_0$, where the first variation vanishes and, thus, 
$$\cos\angle_z=\mu_\partial(z).$$

To make it clear, we do it for the normal deformation $Y_t$ of $Y=Y_0$ with $\psi=1$,
we ignore the contribution from the $Y$-integral and observe, that because of the identity  $\cos\angle_z=\mu_\partial(z)$ on $Y_0$, the 
only non-zero term in the (Leibniz formula for the ) derivative $\partial_\nu -\frac {d}{dt}$ 
of the above $Z$-integral is 
$$\int_Z \frac{1}{\sin \angle_z} (\partial_\nu \cos \angle_z  -  \partial_\nu \mu_\partial(z))
dz=-\int_Z \partial_\nu \angle_zdz+  \partial_\nu \mu_\partial(z)dz,  $$
where the derivative of the angle $\angle_z=\angle_z(t)$ is determined  as follows.

 Intersect $S=\partial X$  and $Y=Y_0$ with (a  germ at $z$ of )  a surface  $E_z\subset X_+\supset X$, which is  normal to $Z=S\cap Y=\partial Y$ at $z\in Z$ and is geodesic at 
 $z$, e.g. being the image of the local exponential  map from the  normal plane 
 $T_z^\perp(Z)\subset T_z(X)=T_z(X_+)$ to $X_+$, where $X_+$ is the above extension of $X$.
 
 Let  $$\mbox {$\underline Y=\underline Y(z)=Y\cap E_z$ and $\underline S=\underline S(z)=S\cap E_z$}$$  
 be the   intersection curves 
  in this   surface  $E_z$,  where   we  identify $E_z$ with a small ball in the Euclidean plane $\mathbb R^2=T^\perp_z(Z)$
  and let $\underline Y_t\subset \mathbb R^2$ be $t$-equidistant curves to 
$\underline Y=\underline Y_0$.

Let  $z_t= \underline Y_t \cap \underline S$, (thus $z_0= z_{t=0}=z$)
 let $y_t\in \underline Y_0$ be  the {\it normal projection of $z_t\in \underline Y_t$}  to $\underline Y_0$, which means that  the straight  segment $[y_t, z_t]\subset \mathbb R^2$   is normal  to the curve  $ \underline Y_0$, (also normal  to $ \underline Y_t$ and having length $|t|$, since $Y_t$ is  equidistant to $Y_0$).
 
 There are two summands that contribute to the difference $\angle_{z_t}-\angle_{z_0}$  between the angles between our curves at their intersection points. 
 
 (1)  The first summand  is due to the turn of the tangent lines to $\underline S$  along the segment   $\underline S_{z_0,z_t}\subset \underline S$  between the points $z_0$ and $z_t$,  which is  
  equal  to the integral of the curvature $\kappa_S$  of $S$ over this (curved)  segment, where 
 $$\int_ {\underline S_{z_0,z_t}} \kappa_ {\underline S}( \underline s)d \underline s= \kappa_ {\underline S}( \underline z_0) |z_0-z_t| +o(t)= \kappa_ {\underline S}\cdot(\underline z_0) \frac {1}{\sin \angle_{z_0}}+o(t),  \mbox {  } t\to 0. $$

(2) The second contribution to $\angle_{z_t}-\angle_{z_0}$ comes  from the curvature of the curve $\underline Y=\underline Y_0$ integrated over the segment $\underline Y_{y_t,z_0}\subset \underline Y$, 
$$\int_{\underline Y_{y_t,z_0}} \kappa_{\underline Y}(\underline y)d\underline y= 
 \kappa_{\underline Y}(z_0) \cot \angle_{z_0}+o(t).$$

Summing up,  the normal derivative of the the $Z$-integral in the first variation formula is expressed in terms of the curvatures of  $Y$ and $S$  and the  angle between them as follows.
{\color {blue}$$\partial _\nu\int_Z \left  (\frac{\cos \angle_z}{\sin \angle_z}   -  \frac {1} {\sin \angle_z} \mu_\partial(z)\right)dz=
-\int_Z\frac {\kappa_{\underline S(z)}(z)}{ \sin \angle_{z}}dz +\kappa_{\underline Y(z)}(z) \cdot\cot \angle_{z} dz+\partial_\nu \mu_\partial(z)dz,$$}
where, recall, 

(i)  {\sf $S$ is the boundary the ambient $n$-manifold $X$, 

(ii) $Y\subset X$ is a hypersurface with boundary $Z=\partial Y=Y\cap S$,

(iii)  $\underline S(z)\subset S$ and $\underline Y(z)\subset Y$ are intersections of $S$ and $Y$ with a germ of a surface $E_z\subset X$   normal to $Z$ at $z$ and geodesic at $z$\footnote{It is better,  as  earlier, to think of $E_z$ in $X_+\supset X$  and  take the image of the exponential map from a small disc in the normal plane $T^\perp_z(Z)\subset T_z(X_+)$ to $X_+$  for  this $E_z$.}, where $\kappa_{\underline S}$ and $\kappa_{\underline Y}$ denote the curvatures of these curves in $E_z$  with the following sign convention:
 
({\Large $\mathbb \pm$})  {\it if the boundary $S= \partial X$ is convex, then $\kappa_{\underline S}\geq 0$; 
 if the  the   "lower region" $X_<\subset X$  bounded by $Y$ is convex, then $\kappa_{\underline Y}\geq 0.$}
 
 (iv) $\angle_z$ is the angle between the tangent spaces $T_z(S)$ and $T_z(Y)$ in $T_z(X)$, which is measured in  $X_<$  { \it under} $Y$,\footnote{
 This  "under", together with ($\pm$), determines the signs of the integrants  in the above formula that is crucial for our (potential) applications.}
 
 (v) the above formula is supposed to hold if $Y$ is  {\it stationary} for the functional
$$Y\mapsto vol^{[-\mu]}_{n-1}(Y)=vol_{n-1}(Y)-\int_{X_<}\mu_\bullet(x) dx- \int_{S_<}\mu_\partial(s)ds$$
where  $\mu_\bullet(x)$  and $\mu_\partial(s)$  are continuous functions on $X$ and on $S$ and where  $X_<\subset X$ is the just mentioned   "lower region" in $X$ bounded by $Y$ and $S_<=S\cap X_<$.}

The above formula  for {\color {blue}$\partial _\nu\int_Z...$} can be neatly rewritten in terms of the mean curvatures $M_S$ of $S$, $M_Y$ of  $Y$  and $M_Z$  of $Z$ in $Y$  by invoking   the following. 

{\it Algebraic Identity.}\footnote{Compare with  3.8 in [Li(comparison) 2017]. }
{\color {blue}$$M_Z(z) =\frac {M_S(z)}{\sin\angle_z}+(\cot\angle_z)M_Y(z)- \frac {\kappa_{\underline S(z)}(z)}{ \sin\angle_z }-(\cdot \cot \angle_z)\cdot \kappa_{\underline Y(z)}(z). $$}
(What is {\color {magenta!60!black} significant} is that  the coefficients  on the right hand side here are the same as in  the above expression for the $Z$-term in the second variation formula.)

To prove this identity, let us express everything in terms of the traces second fundamental forms $\mathrm {II}_S$, $\mathrm {II}_Y$ and $\mathrm {II}_Z$, where

$\mathrm {II}_S$ is taken with outward normal field denoted $\nu_S=\overleftarrow\nu^\perp_{ S}$

$\mathrm {II}_Y$  is taken  with the upward  field $\nu=\nu_Y=\nu^\perp_{Y\uparrow }$

and where

 $\mathrm {II}_Z$ will be evaluated with two unit normal fields to $Z$, one of them $\nu_Y$
 restricted to $Z$ and the other one is  tangent to $Y$  and facing outward, call it 
 $\nu_Z=\overleftarrow\nu^\perp_{ Z}$. (If $Y$ is {\it normal} to $S$ at $z$, i.e. $\angle_z=\pi/2$, then $\nu_Z(z)=\nu_S(z)$.)

Observe that 

   $ \nu_Y$  
is {\it normal} to $\nu_{Z}$   and that 

 the angle between   $\nu_{S}$  and $\nu_Y$ is {\it complementary} to 
 the angle $\angle_z$  between $S$ and $ Y$ at all  $z\in Z=S\cap Y$,
 $$\angle_z (\nu_{S}, \nu_Y)=\pi-\angle_z.$$
Therefore,  
 $$\nu_{ S}(z)= (\sin \angle_z)\cdot  \nu_Z(z)- (\cos \angle_z) \cdot \nu_Y(z) $$
and, by the linearity of  the form $ \mathrm {II}_Z$ and its trace in the normal vectors, 
 
$$ trace_Z (\mathrm {II}_S)
=(\sin \angle_z)\cdot trace_{\nu_Z}( \mathrm {II}_Z)-(\cos\angle_z) \cdot trace_Z( \mathrm {II}_Y),  $$
or

$$M_Z=trace_{\nu_Z}( \mathrm {II}_Z)=\frac {1}{ \sin\angle_z} trace_Z (\mathrm {II}_S)
+(\cot\angle_z) \cdot trace_Z( \mathrm {II}_Y),\leqno{\color {blue}[1/\sin]}$$ 
where  $trace_Z (\mathrm {II}_S)$ denotes the trace of the restriction of  the form $\mathrm {II}_S$ to (the tangent bundle of)  $Z\subset S$, that is the same as   the trace of  the form  $\mathrm {II}_Z$ with respect to the vector field $\nu_S$ restricted to $Z$, where   $trace_Z( \mathrm {II}_Y)$ is understood similarly  and where 
$trace_{\nu_Z}( \mathrm {II}_Z)$
denotes the trace with respect to  the field $\nu_Z$, where, indeed, it is equal to the mean curvature $ M_Z$  of $Z$ in $Y$, 
$$trace_{\nu_Z}( \mathrm {II}_Z)=M_Z.$$

Finally, we observe that the curvatures of the curves $\underline S(z)$ and $\underline Y(z)$ in $E_z$
are equal to the traces of the form  $\mathrm {II}_S$ and $\mathrm {II}_Y$ restricted to  $\underline S(z)\subset S$ and  $\underline Y(z)\subset Y,$
$$\kappa_{\underline S(z)}(z)=trace_{\underline S(z)}\mathrm {II}_S(z) \mbox { and } 
\kappa_{\underline Y(z)}(z)=trace_{\underline Y(z)}\mathrm {II}_{Y}(z),$$
while 
$$trace_{\underline S(z)}\mathrm {II}_S(z)+ trace_{Z}\mathrm {II}_S(z)=trace_S\mathrm {II}_S(z)=M_S(z)$$ 
and 
$$trace_{\underline Y(z)}\mathrm {II}_Y(z)+ trace_{Z}\mathrm {II}_Y(z)=trace_S\mathrm {II}_Y(z)=M_Y(z)$$
These, combined with {\color {blue}$[1/\sin]$} yield the required algebraic  identity which we write now as 
$$- \frac {\kappa_{\underline S(z)}(z)}{ \sin\angle_z }-(\cdot \cot \angle_z)\cdot \kappa_{\underline Y(z)}(z)=M_Z(z)-
\frac {M_S(z)}{\sin\angle_z}+(\cot\angle_z)M_Y(z)$$

Substitute this  into the above  formula for the  {\color {blue}$\partial_\nu$}-derivative of the integral
{\color {blue}$\int_Z...dz$} in the first variation formula for {\color {blue} $vol^{[-\mu]}_{n-1}(Y)$} and express this derivative   by the following

\hspace {20mm} {\color {blue}\large  \it Mean Curvature Stability Relation.}

{\color {blue!40!black}}$$\partial _\nu\int_Z...dz=\int_Z \left(M_Z(z)-\frac {M_S(z)}{\sin\angle_z}+(\cot\angle_z)n\cdot M_Y\right)(z)-\partial_\nu \mu_\partial(z)dz.$$}

Thus, for instance, if  $\mu_\partial(z)$ is {\it constant} and $Y$ is a {\it local minimizer} for $vol^{[-\mu]}_{n-1}(Y)$, then, this formula, which necessarily holds for the integrals over all subdomains $U\subset Y$,  shows that  
{\color {blue} $$M_Z(z)-\frac {M_S(z)}{\sin\angle_z}+
(\cot\angle_z)n\cdot M_Y\geq 0,\leqno {\color {blue}  [\geq]}$$}
which us most informative (and quite useful) if $\mu_\bullet=M_Y=0$.

\vspace {1mm}
{\it \color  {magenta} About the Signs.} Consistency in the choices of  signs in the definitions of the curvatures, and/or of normal vectors  for the second fundamental forms  is crucial for applications.  

There is hardly a   problem here with  the sign of the $M_S/\sin$-term, since it is  clearly visible by looking at the case where $Y$ is {\it normal} to $S$, i.e.   $\angle_z=\pi/2$;  it is also instructive to go through the full calculation in the following.

{\it Example/Exercise 1.} Let $f(t)>0$, $0<t<\infty$ be a smooth function and $X\subset \mathbb R^{n-1}\times \mathbb R_+$ be the rotation body of the subgraph of the function $f$ (i.e. the region below the graph)  around the $t$-axes. 

Let 
$\mu_\bullet=0$ and $\mu_\partial $ is a constant, say  $\mu_\partial(t)=c$.

Let $Y_t=\subset X$ be the $(n-1)$-balls of radii $f(t)$ normal to the $t$-axes  and let  us compute the the second variation,  of  $vol^{[-\mu]}_{n-1}(Y_t)$, at $t$ where this ball is stationary, that is the second derivative
$$(bR^{n-1}- c\cdot vol_{n-1}(S_{<t})'',$$
where $b=b_{n-1}$ denotes the volume of the unit  ball $B^{n-1}\subset \mathbb R^{n-1}$ and where
$S_{<t}\subset S=\partial X$  is the part of this boundary below (or to the left from) $t$, where the stationary $Y_t$ is where the first derivative vanishes, i.e. 
$$(bR^{n-1}- c\cdot vol_{n-1}(S_{<t}))'=(n-1)b R'R^{n-2}-bc (n-1)R^{n-2} \sqrt {1+R'^2}=0,$$
that is
$$\sqrt {1+R'^2}=R'/c.$$
Then an elementary calculation show that  
$$(bR^{n-1}- c\cdot vol_{n-1}(S_{<t})''=-vol(S^{n-2}(R))\frac {k_f(t)}{\sin \angle_t} =-vol_{n-2}(\partial Y_t)\frac {k_f(t)}{\sin \angle_t},$$
for $k_f(t)$ being the curvature of the graph of $f(t)$ which, according to our convention, is {\it negative} for $f(x)=x^2,$ since the subgraph of $x^2$ is concave. 
\vspace {1mm}

The computation becomes messier if  $Y$ is non-flat  but it still 
durable in simple cases.\vspace {1mm}

{\it Example/Exercise 2.} Let $f$ and $X$ be as above and let  $Y_t \subset X$ be the intersections of $X$ with the spheres of radii $t$  centered at  $\mathbf 0$ that is the  zero on the $z$-axes, where we assume that  the intersections $Z_t$ of these $t$-spheres with $S=S_f=\partial X$ are  {\it  non-empty connected  transversal} for $t\geq 1$.

Let $\mu_\bullet(t) $  be equal to the  mean curvature  of $Y_t$, i.e. $\mu_\bullet(t)=\frac  {n-2}{t}$ and
let  $\mu_\partial$ be constant, denoted $\mu_\partial=c$.

We invite  the reader to
  evaluate the second variation of $ vol^{[-\mu]}_{n-1}(Y_t)$ at stationary $Y_t$. 
 %%%%%%%%%%%%%%%%%%%%%%%%%%%% 

\subsubsection {\color {blue}  Capillary Warped Products Inequalities}\label{capillary warped5}

%%%%%%%%%%%%%%%%%%%%%%%%
Most (all?) extremality/rigidity  properties of warped  products proved in the earlier
 sections,  as well as Gauss-Bonnet kind inequalities from the following section, have their  counterparts 
 for manifolds with boundaries,  which are proven with "capillary" $\mu$-bubbles with a use of the  above inequality 
 {\color {blue}  $[\geq]$}.

 We formulate  below a  a few  examples  and postpone a more thorough analysis and   applications, e.g. to manifolds with corners {\footnote {See [Li(comparison) 2019] in this regard.} 
 until another occasion.\vspace{1mm}
 
 \textbf  {Spherical Suspension Inequality.}  Let $\underline X_0\subset S^{n-1}\subset S^n$ be a smooth convex domain 
   in the equatorial sphere  $S^{n-1}\subset S^n$  and let $\underline X_{\pm 1}= \underline X_{\pm 1}(r)\subset S^n$ be the the union of the  geodesic segments between the north and the south           
  poles  of $S^n$ through  this domain. Remove the poles  $\pm 1\in \underline X_{\pm 1}\subset S^n$ from $ \underline X_{\pm 1}$ and denote 
 $$ \underline X= \underline X_{\pm 1}\setminus \{-1,+1\}\subset S^n\setminus  \{-1,+1\}$$

{\sf Let $X$ be a (non-compact) Riemannian manifold with a  boundary and let 
$f:X\to \underline X$
be a smooth  proper} (boundary-to-boundary, infinity-to-infinity)  {\sf  1-Lipschitz map 
$f:X\to \underline X$ of non-zero degree.}

{\it Let  the scalar curvature of $X$ is bounded from below by
$$Sc(X)\geq n(n-1)=Sc(S^n).$$
If $X$ is spin, if $n\leq 7$, 
and if either $n$ is  odd or $\underline X_0$ is a ball,\footnote {Probably these "if"  are unnecessary.}
then there exists a point $x\in\partial X$, where the mean curvature of $X$ at $x$ is bounded by that of   $\underline X$ at $f(x)$,
$$mean.curv(X,x)\leq mean.curv(\underline X,f(x)).$$}
Moreover,

{\it  if  $mean.curv(X,x)\geq mean.curv(\underline X,f(x))$  at all $x\in X$,
then $mean.curv(X,x)= mean.curv(\underline X,f(x))$ and the map $f$ is an isometry.}

\vspace {1mm}

{\it About the Proof.} The  "right sign" {\color {blue}  $[\geq]$}  at the boundary, allows   carrying through 
the $\mu$-bubble 
 argument  from \ref{log-concave5} in the proof of the extremality/rigidity of double punctured spheres.  This    
reduces the problem to  the  { \color {blue}\EllipseSolid}-comparison theorem 
 for $Sc(V) \geq 0$ of log-concave warped products from section \ref{Sc-criteria3}  and the proof follows.

\vspace {1mm}

{\it About $n\geq 8$.} Probably, the case $n=8$ follows  by a Natan Smale's kind of perturbation argument 
and if $n\geq 9$ generalizations of  Lohkamp's and/or of Schoen-Yau's arguments would work, but  
  singularities at  capillary boundaries need additional care. 
\vspace {1mm}

{\it About  Corners.} The above theorem remains valid for non-smooth domains $\underline X_0\subset S^{n-1}$  with properly understood 
generalized mean curvature, e.g. for convex $k$-gons in $S^2$, where the proof can be obtain by properly smoothing the corners.(See the next section.) 

In general,  
the density function $\mu\partial(s)$ on the boundary $S=\partial X$ from  the previous section may be (and typically is) 
discontinuous along the corners.  In this case the smoothing argument introduces an unpleasant error and the behaviour of $\mu_\partial$-bubbles at the corners needs an additional care.\footnote {In the  case 
of $\mu_\partial(s) $ constant on the faces of 3-dimensional domains,  the proof of  $C^{1,\alpha}$-regularity of  capillary surfaces at the corners is indicated in [Li(rigidity) 2019], see next section.}  

\vspace {1mm}

 \textbf  {Spin-Extremality of Doubly Punctured Balls.}  {\sf Let $X$ be a compact manifold with {\it non-negative scalar  curvature} and  a {\it mean convex boundary} $S=\partial X$ ,   let $P_-,P_+\subset \partial X$ be two closed subsets and  let $f: S\to \underline S= S^{n-1}$ be a 1-{\it Lipschitz} map  of {\it non-zero degree} , such that the subsets {\it $P_\pm$ go to the North and the South poles}  of the  unit  sphere $\underline S=S^{n-1}=\partial B^n\subset\mathbb R^n$.}

{\it If $X$ is spin and if $n=dim(X)\leq 8$,\footnote {Both  conditions are, probbaly,  redundant, where dropping the  the latter could be possible with the recent Lohkamp's techniques, while removing  the former remains beyond the range of the present day knowledge.} then the mean curvature of 
 $S=\partial X$ outside $P_-$ and $P_+$ can't be greater than that of $S^{n-1}$,
$$\inf_{s\in S\setminus (P_-\cup P_+)}mean.curv( S, s)\leq n-1.$$}

\vspace {1mm}

{\it Sketch of  the Proof.} Let $\underline\mu_\partial(\underline s)=\cos \angle_{t(\underline s)}$, $\underline s\in \underline S=S^{n-1}$,  be the function on $S^{n-1}$ as in 
 {\color {blue}example} from the previous  section} (where  $\angle_{t(\underline s)}$ are the  angles between $S^{n-1}$ and the  (parallel) hyperplanes $\mathbb R^{n-1}\times\{t\}$)  and let 
 $\mu_\partial=\mu_\partial(s)$ be the composed function 
 $$S\overset {f}\to \underline S\overset{\underline \mu_\partial}\to \mathbb R\mbox {  on $S=\partial X$, that is 
 $\mu_\partial(s)= \underline \mu_\partial (f(s) $, $  s\in S$}.$$.
 
Then,  arguing as in the proof of the 
double puncture  theorem  for spheres $S^n$ in sections  \ref{punctured3} \ref {log-concave5},
we conclude to the existence  of a stable $\mu$-bubble $Y\subset X$  for $\mu=\mu_\partial ds$,
which is smooth up to the boundary for $n\leq 8$. (if   $n=8$ one  needs a version of Nathan Smale's generic regularity  theorem.)

Next,  the 
mean curvature stability relation from the previous section show that the mean curvature   
$M=M_Z$ of the boundary $Z=\partial Y\subset Y$ in $Y$
  and the 
norm of the differential of the natural map $\phi: Z\to S^{n-2}$ satisfy the inequality 
$$\frac {M(Z, z)}{||d\phi(z)||}\geq n-2.$$ 
(Hopefully, there is no silly error   in the computation)

Finally,  the mean curvature spin-extremality theorem  \footnote {One needs here 
  this theorem for maps   to the convex hypersurfaces not only in Euclidean spaces but also  in other Riemannian flat manifolds, specifically in  $\mathbb R^m\times \mathbb T^N$ in the present case.} 
from section \ref{mean convex3} applies to $\mathbb T^\rtimes$-stabilized $Y$, that is to  
$Y\rtimes \mathbb T_1$, and the proof follows. \vspace {1mm}

{\it Remarks.} (a) It is not hard to prove, as in  the cases we encountered earlier, that the balls are rigid in this regard: 

{\it if $$\inf_{s\in S\setminus (P_-\cup P_+)}mean.curv( S, s)= n-1,$$ 
then $X$ is isometric to $B^n$.}

(b) The above argument generalizes to {\it complete} non-compact manifolds $X$, but, probbaly, 
 completeness can be replaced by a weaker condition.

\vspace {1mm}

{\textbf {Corollary to the Proof:] Multi-Width  Mean Curvature Inequality for non-Spin Manifolds.}

{\sf Let $X$ be a compact Riemannin $n$-manifold  with a boundary,  
and let $f$ be a continuous map  from $\partial X$ to the boundary of the $n$-cube with {\it non-zero degree,}
$$f:\partial X\to\partial [-1,1]^n, $$
such that the {\it distances between the pullbacks of the opposite faces of the cube are all 
$\geq \pi$,}
$$dist (\partial_{-,i}, \partial_{+,i})\geq \pi,\mbox { } i=1,...,n.$$}

{\it If $X$ has non-negative  scalar curvature, $Sc(X)\geq 0$, and if $n=dim(X)\leq 8$\footnote{One can drop this if one  extend  Schoen-Yau's "desingularization" theorem for capillary hypersurfaces}   then 
$$\inf_{x\in \partial X} mean.curv(\partial X,x)\leq n-1.$$}

{\it Proof. }   Apply the above argument to the $f$-pullbacks  of a pair of opposite faces of the cube,
 say to 
 $$P_{\mp}=f^{-1}(\partial_ {\mp,1})\subset \partial X$$ 
let $Y\subset X$ be the  corresponding  stable $\mu$-bubble which separates $P_-$ from $P_+$.
Then apply the same argument to  $Y\rtimes \mathbb T^1$  and continue  inductively as in the  proof of 
 the multi-width {  \color{black}  $\square^n$-{Inequality} $\sum_{i=1}^n \frac{1 }{d_i^2}\geq  \frac  {n^2}{4\pi^2}$ in section\ref 
 {multi-width3} thus reducing the problem to the  case of $X^2\rtimes \mathbb T^{n-2}$, where 
 $X^2$ is a surface with $curv(\partial X^2)\geq 1$, where $Sc(X^2\rtimes \mathbb T^{n-2})\geq 0$,
 and where the proof follows by the proof of  the mean curvature spin-extremality theorem.
 \vspace {1mm}

{\it Remarks}  (a)  The above inequality improve non fill-is results 4(A)   and (5) in section  \ref {fill-ins}

(b) If $X$ is spin this inequality follows from the mean curvature spin-extremality theorem.

(c) One can improve this inequality  with a better (iterated) warped product model manifold $\underline S$,
and, probbaly, with the best such $\underline S$ the improved inequality, will not follow from the mean curvature extremality of spheres, even for spin manifolds $S$.
(It is not impossible, that the the sphere $S^{n-1}$ radially mapped to $\partial [-1,1]^n$ is extremal 
 this inequality.)

\vspace {1mm}

\textbf {Capillary   Mean Curvature Separation Theorem.}  {\sf Let $X$ be a compact manifold  with boundary $S\partial X$  and let  ,   let $P_-,P_+\subset \partial X$ be two closed subsets
 sich that,
 
 $\bullet_\sigma$ the scalar curvature of $X$ is bounded from below by a given non-positive number,
 $$Sc(X)\geq  \sigma,\mbox  { } \sigma\leq 0;$$ 
 
 $\bullet_M$ the values of mean curvature of $S$ in  the subset   $P_-\subset S$ and in its complement  are bounded from below as follows
$$\mbox{ $mean.curv (S, s) \geq M_-$,  $s\in P_{-}$
   and  $mean.curv (S, s) \geq M_+$, 
   $s\in S\setminus P_-$,}$$
 for some   $M_-$, $ M_+$, where  $M_+$ is  {\it positive} while  $M_-$ can be {\it negative}. 
 
 {\it Let $M_+$  be bounded from below in terms of  $-\sigma $  and $-M_-$ according to the following inequality
$$M^2_+\geq \max \left (      \frac {n-1}{-n\sigma},  -M_-\right),$$
and let  
  the  distance $D$ between $P_-$ and $P_+$ measured  in $S$,   with respect to  the induced Riemannian metric  in $S\subset X$  be bounded in terms of $M_+$ as follows,}  
$$D\geq const_n \frac{1}{M_+}\mbox { for } const_n\geq 100\pi. $$}

Let  $M'>0$ be a given number  and let 
 the numbers $=M_+$ and $D$  be  {\it sufficiently large} depending on $n$,  $\sigma$,  $M_-$ and $M'$  -- specific inequalities are indicated below.}
 
 {\it Then, assuming $n=dim X\leq 8$,  there exists  a smooth compact hypersurface $Y\subset X$ with boundary $\partial Y\subset S=\partial X$, such that
 the mean curvature of the boundary of $ Y$ is bounded from below by $\frac {1}{2} M_+$,
$$ mean.curv(\partial  Y)\geq \frac {1}{2} M_+,
$$and the scalar curvature of some  warped $\mathbb T^\rtimes $-extension of $Y$ is non-negative, 
$$Sc(Y\rtimes \mathbb T^1)\geq 0.$$} \vspace {1mm}

{\it Sketch of the Proof.} Let 
$$\mu=\mu_\bullet(x) dx + \mu_\partial(s)  dx $$
where  $\mu_\bullet(x) = M_+$   and where  $\mu_\partial(s)$ is "induced" as earlier   from the   function $\underline \mu_\partial(\underline s)=\cos_{t(\underline  s)}$ on $S^n$ by a  $\frac {\pi}{D}$-Lipschitz map $S\to S^{n-1}$, 
 which sends $P_+$  to the North pole and $P_+$ to the South pole  of $S^{n-1}$. (The existence of such a map is obvious.) 

Then our conditions on the mean curvatures guarantee the existence of a stable $\mu$-bubble $Y\subset X$, which  separates $P_-$ and $P_+$ and  has (free) boundary in $S$ and   where the second variation formula along
with the mean curvature stability relation from the previous section imply the desired properties of this $Y$.

\vspace{1mm}

{\it Remarks. } (i) Our bound   on $D$ is very rough. We suggests the reader would find a better estimate.

(ii) If one takes  into account, besides $D=dist_S((P_-,P_+)$, the distance $d=dist_X(P_-,P_+)$    then, one can prove, with some $\mu_\partial ds+\mu_\bullet(x)dx$  for  a suitable function $\mu_\bullet(x)$,     a comprehensive   separation theorem incorporating the above with 
{\color {blue} {\large \sf I{\large \sf \color {red!80!black}\sf I}\sf I}} from section \ref{separating3} for closed manifolds.

{\it Problem.} Find adequate version of the  "log-convexity"  condition  on  $\mu_\partial ds+\mu_\bullet(x)dx$  and find all "interesting"  sharp  capillary   extremality/rigidity  inequalities including all such inequalities presented in the previous sections.
\vspace{2mm}

%%%%%%%%%%%%%%%%%%%%%%

\subsection {\color {blue}3D Gauss Bonnet Inequalities}\label {Gauss-Bonnet5}

%%%%%%%%%%%%%%%%%%%%%%%
The simplest inequality of this kind, which  applies  to closed connected cooriented stable minimal surfaces $Y$ in orientable   Riemannian 3-manifolds $X=(X,g)$, is a bound  on  the integral of the scalar curvature of $X$ over $Y$, that reads:
$$\int_y Sc(X,y)dy\leq 8\pi,\leqno {(A)}$$
where the equality holds for Riemannian products 
$Y_0\times S^1$, for sufaces  $Y_0=(Y_0, h_0)$ homeomorphic to $S^2$.
(Compare with {\it area exercises} in section \ref{SY+symplectic2}.)
{\it Proof}. Combine the inequality  {\Large \color {blue}$[ \star\star]$} involved in the second variation formula (section \ref {2nd variation2}) with the Gauss-Bonnet theorem.  

\vspace{1mm}

{\it Corollary}. {\sf If  $X$ is compact with $Sc(X)>0$, then the 
2-systole of $X$ with the metric $g^\ast(x)=Sc(g,x)g(x)$  satisfy  is bounded by $8\pi$.}
Moreover

\vspace {1mm}

{\it The 2-dimensional homology of $X$ admits   a basis represented   by closed surfaces $Y\subset X$
with $area_{g^\ast}(Y) \leq 8\pi$}.\vspace {1mm}

{\it Question.} Can one directly bound the areas of $g^\ast$-minimal surfaces in $X$?
\vspace {1mm}

{\it Using the  Dirac Operator.} Let us give a Dirac theoretic proof of this corollary, where, observe,
this is the only known case where the Dirac operator  goes in parallel with  minimal surfaces.

To simplify, let $X$ be homeomorphic to $S^2\times S^1$  and to  keep track of constants let us   compare the metric $g_\ast$ on this $X$ with the Riemannian  product $\underline X = (\underline X,\underline g)$ of the circle  $S^1$  the unit sphere $S^2$ with it{'}s  usual metric with $Sc=2$.

Let $X^4=X\times S^1$  and $\underline X^4=\underline X\times S^1$ be the corresponding 
$4$-manifolds, where the Dirac operator will be employed,  let $\underline L^\ast \to \underline X^4$ be the line bundle induced from the Hopf bundle by the natural map $ \underline X^4\to  S^2$ and let $ \underline L^\circ_\varepsilon \to X^4$ be an $\varepsilon$-flat bundle induced by the natural map $X^4\to S^1\times S^1$ 
from an $\varepsilon$-flat bundle  $ L_\varepsilon \to S^1\times S^1$, such that the first Chern class of 
$ L_\varepsilon$ doesn't vanish

Then the twisted Dirac  $\mathcal D_{\otimes L^\ast\otimes L^\circ}$ on $X^4$ has 
{\it non-zero index},   and this nonvanishing of $ind(\mathcal D_{\otimes L^\ast\otimes L^\circ})$ persists for all
Riemannian metrics $c{'}$ on $X^4$.  

On the other hand, if a line bundle $L \to X^4=(X^4, g{'}_\ast=g{'}Sc(g{'}))$   has the   the norms of its   curvature  $\omega$  bounded by curvature $\underline\omega$ of $\underline L^\ast $ according to the inequality
$$  ||\omega||_{g{'}_\ast}=||\omega||_g{'}(x) /Sc(g,x)< ||\underline\omega||_{\underline g}(x) /Sc(\underline g,x)=\frac {1}{8\pi},$$

then $ind(\mathcal D_{\otimes L^\ast\otimes L^\circ})=0$ as it follows   from the twisted  Lichnerowicz-Weitzenboeck-formula (and a little computation). 

Now let us  assume that  the $g_\ast$-areas of all non-homologous to zero $2$-cycles in $X$ are bounded from below by $8\pi+\epsilon$.

 Then,  by the {\it Morse lemma for mass in codimension 1}, the mass of the generator of 
 $H_2(X,\mathbb R$ is also bounded from below by $8\pi+\epsilon$, which, by duality,   bounds  the comass of the corresponding generator of $H^2(X;\mathbb R$ by 
 $(8\pi+\epsilon)^{-1}.$
 Hence, there exists a 2-form $\omega_0$ on $X$  with $g_\ast$-norm $\leq (8\pi+\epsilon)_{-1}$ in the cohomology class of the curvature form of the line bundle induced from the Hopf bundle. 
 
 $\omega_0$  by the curvature of a line bundle over $X$, lift this bundle $X^4=X\times S^1$   and,  confront its   properties with the above discussion.
 
   Then, 
  by contradiction, we conclude  that  the $g_\ast$-areas of  {\it certain  non-homologous to zero } $2$-cycles in $X$ must be arbitrarily close to $8\pi$. (One could go to the limit and get  such cycles 
with areas $\leq 8\pi$, but doing this, which  needs an additional, let it be a   well known, argument,  is unnecessary for our purpose.)
\vspace {2mm}

Let us return to minimal surfaces and formulate  a version of the above (A)
 for    (compact  orientable Riemannian) $3$-manifolds $X$  {\it with boundaries}, denoted $S=\partial X$, which   involves the integral of 
the mean curvature $M(S)$ over boundary curves of surfaces $Y\subset X$ with
$Z=\partial Y=Y\cap S$. Namely,\vspace {1mm}

{\it connected cooriented   cooriented    surfaces $Y\subset X$ with non-empty boundaries $Z=\partial Y=Y\cap Z$  which are stable  minimal for the free boundary condition, satisfy:
$$\frac {1}{2}\int_Y(Sc(X,y) dy+\int_Z M(S,z)dz\leq 2\pi.$$}

{\it About the Proof}. This can be obtained by applying the above (A) to the double of $X$, or,  alternatively,  with a use of the second variation formula for manifolds with boundaries  from section\ref{}, (where only the simplest case of $\mu=0$ is needed here).
\vspace {1mm}
 
 {\it Corollary. } {\sf Let $S\subset \mathbb R^3$ be a smooth {embedded 
  {\it non-simply   connected} closed surface. Then there exists  a closed {\it non-contractible} curve
 $Z\subset S$. such that 
 $\int_ZM(S,z)dz\leq 2\pi$.}

\vspace {1mm}

{\it Question.} Can one find such a curve in $S$ without using minimal surfaces in the domain bounded by $S$?}

\vspace {1mm}

{\it Gauss-Bonnet Extremality of Truncated Cones.}  Let $\underline X\subset  \mathbb R^3$
be a round truncated cone, the essential  invariant of which is the angle $\beta$ between the side surface $\underline S$ of this cone and the bottom, where  as in section {\ref {capillary5} we prefer do deal with the complementary angle $\alpha=\pi-\beta$ and where the inequality
{\color {blue}$[\geq ]$} from section \ref {capillary5} 
{\color {blue} $$M_Z(z)-\frac {M_S(z)}{\sin\angle_z}+
(\cot\angle_z)n\cdot M_Y\geq 0,\leqno {\color {blue}  [\geq]}$$}
 become an equality, where
the  curvature $M(\underline  Z)$ of the  horizontal circles $\underline Z  \subset\underline S$  is related to the mean curvature of  of $S$ along these circles by
$$M(\underline  Z) = \frac {M(\underline S)}{\sin\alpha}.$$

Now let $X$ be a compact Riemannian 3-manifold  with boundary which is divided in 3 parts
bottom $B$, top $T$ and the side surface $S$, which separates $B$ from $T$ and such that  

{\it the angle between  $B$ and $S$ is everywhere $\leq \beta$  and the angle between $S$ and 

$T$ is everywhere $\leq \alpha$.}

 Let, moreover,  $B$ and $S$ be  {\it mean convex}, i.e. their mean curvatures with respect to the outward normals are {\it positive}. 

Under this condition the functional $Y\mapsto area(Y)- \cos(\alpha) area (S_< )$ defined on surfaces $Y\subset X$ with boundaries $Z=\partial Y\subset S$ and separating $B$ from $T$ 
necessarily assumes minimum at a  surface $Y\subset X$ with $\partial Y\subset S$, which satisfies
according to  {\color {blue}$[\geq ]$}:
 $$\int_Y Sc(X,y)dy+ \frac {1}{\sin\alpha} \int_Z M(S,z)dz\leq \pi.$$

As earlier,  the  most interesting case is  for $Sc(X)\geq 0$ and $M(S)\geq 0$, already  for domains 
$X\subset \mathbb R^3$, where the existence of a   curve $Z\subset S$  separating 
the top from the bottom and having  $\frac {1}{\sin\alpha}\int_Z M(S,z)dz\leq \pi$ seems  non-obvious. (Am I missing a direct obvious proof?)

Also   note  that similar inequalities hold for manifolds $X$ with  more complicated corners (see  section 5.4 in [G(billiards) 2014] and  [Li(comparison) 2017]) but many such inequalities  still reman conjectural.

Besides  manifolds with $Sc\geq 0$, the above type Gauss-Bonnet inequalities yield geometric information for manifolds with  scalar curvatures bounded from below by  {\it negative constants $\sigma$,} where this information is somewhat opposite to that for manifolds $X$ with $Sc(X)\geq \sigma>0$.

Namely, in the later case one conclude that $X$ must have representatives of non-zero homology classes by surfaces of area bounded by $const
\cdot\sigma$. On the contrary,  the bound  $Sc(X)\geq \sigma$  for $\sigma<0$, implies, under additional {\it topological conditions}, that $X$ can{'}t have such  surfaces with small area.

{\it Example.} {\sf Let  $X$ be  homeomorphic to $S_\chi\times S^1$, where    $S_\chi$ is a closed connected orientable  surface  with the Euler  characteristic $\chi<0$}. 

{\it If $Sc(X)\geq -2$,  then all  surfaces $Y\subset X$ in the homology class of  $S_\chi\times \{s_0\}\subset X$ have 
$$area(Y)\geq 2\pi |\chi(Y)|.  $$}

This is, of course, obvious. What is slightly more interesting is a similar inequality for "area minimizing"  families of 2d-foliations in $X$, but these inequalities are inherently non-sharp  in the key example
of hyperbolic manifolds $X$ (see  [G(foliated) 1991]) for more about it). 

What looks more promising are foliations by $\mu$-bubbles using  horospherical 
foliations for models,  but the corresponding inequalities here is yet to be properly formulated and proved.

%%%%%%%%%%%%%%%%%%%%%%%

\subsection {\color {blue}  Topological Obstructions to $Sc>0$  Issued from  Minimal Hypersurfaces and $\mu$-Bubbles}\label{obstructions5}
%%%%%%%%%%%%%%%%%%%%%%%

Start with recalling the proof of    Schoen-Yau's Non-Existence\&Rigidity Theorem  by $\mathbb T^\rtimes$-stabilization argument (see section \ref{bubbles1}) applied to complete non-compact manifolds  as follows.

\vspace{1mm}

\textbf 0.   {\sf Let a smooth open  orientable   manifold $X$ contain  a decreasing chain  (flag) of oriented   {\it properly embedded}  (infinity-to-infinity)  submanifolds
 $$X\supset X_{-1}\supset ...\supset X_{-i} \supset... \supset X_{-(n-2)}, \mbox { }  dim( X_{-i})=n-i,$$ 
such that 
the homology classes   $[X_{-i}\in H_{i, inf}(X)$ of $X_{-i}$   with infinite supports  are  {\it non-zero}  
for all $i$ and the class  $[X_{-(n-2)}] \in H_{2, inf}(X)$  is not representable by a simply 
connected surface} (i.e. by $S^2$ or $\mathbb R^2$).

{\it If $X$ supports a  complete metric with $Sc\geq 0$, then $X$ is isometric to the product 
$X=X_0\times\mathbb R^1$, where $X_0$ is a flat manifold.}\footnote{As usual, if $n\geq 8$, 
one has to  appeal   to "desingularization" results  from  [Lohkamp(smoothing) 2018]  or from [SY(singularities) 2017]. 
(If $X$ is spin and $H_1(X)$ has no torsion, then  the results from section 6 in [GL(complete) 1983] apply.)}

\vspace {1mm}

{\it Proof.} By  Jerry Kazdan's  perturbation theorem and Cheeger-Gromoll splitting theorem, the case $Sc\geq 0$
reduces to that of $Sc>0$, where the $\mathbb T^\rtimes$-symmetrization shows that $X$ contains a properly embedded surface  $Y\subset X$ in the( infinite) homology class of $X_{-(n-2)}$, such that some warped  
product $Y\rtimes \mathbb T^{n-2} = (Y\times \mathbb T^{n-2}, g^\rtimes = dy^2+\phi^2(y)dt^2)$ has 
positive scalar curvature, 
 $$Sc(Y\rtimes \mathbb T^{n-2})=Sc(g^\rtimes)>0.$$

Hence, $Y$ must be simply connected. Otherwise  a covering $\tilde Y$ of $Y$ with
 infinite cyclic fundamental group $\pi_1(\tilde Y)=\mathbb Z$     would allow an  extra  $\mathbb T^\rtimes$-symmetrization, and turn into  a complete manifold $\mathbb R\times \mathbb T^{n-1}$ with a 
(warped product) metric $\tilde g^\rtimes =dx^2+\varphi(x)^2d\tilde t^2$, for $x\mathbb R^1$  
and $t \in \mathbb T^{n-1}$,   on $\mathbb R\times \mathbb T^{n-1}$ {\it invariant under the action of the torus}  and such that  $Sc(\tilde g^\rtimes)=0$.

 Thus, impossibility of this follows by the  formula
 $$Sc(g)(x,\tilde t) =  -\frac  {(n-1)(n-2)}{ \varphi^2(y )}
\Bigg\lVert  \frac{ d \varphi(x)}{dx}\cdot \frac{1}{\varphi}\Bigg\rVert^2-\frac {2(n-1)}{\varphi(y)}\frac{ d^2  \varphi(x)}{dx^2},$$
 since no   function  $\varphi>0$ can have  negative second derivative.\vspace {1mm}

\vspace {1mm}

{\it Manifolds with Spines.}  Let us now turn to  more general open 
manifolds  $X$,  including {\sf infinite coverings of}  {\it compact  enlargeable} 
(e.g.  admitting metrics with non-positive sectional curvatures) manifolds {\it with punctures}  and 
on products of   SYS-manifolds by enlargeable ones, where  geometry depends on 
the distance to a distinguished closed subset $S\subset X$ called the {\it spine of $X$.}

\vspace {1mm}

{\it Example.} If  $X$  comes with a covering map   to a   compact manifold minus a point, $X\to X_0\setminus x_0$, 
then relevant spines $S\subset X$  are the pullbacks of compact subsets 
$S_0\subset  X_0\setminus x_0$.  

\vspace {1mm}

\vspace {1mm}

{\it S-Quasiproper Maps.}  Given a spine $S$ in $X$, a continous map from $X$ to a metric space, say
 $f:X\to\underline X$, is called {\it uniformly  $S$-quasi-proper} if it 
is {\it constant on the connected components} of the
 complement $X\setminus S$ and if the restriction of $f$ to $S$,
 $$f_{|S}: S\to \underline X$$ 
is {\it uniformly proper},  i.e.  the diameters of the $f_{|S}$-pullbacks of subsets from $\underline X$ are bounded in terms of the diameters of  these subsets,
$$diam(f^{-1}(\underline U)\cap S) \leq \xi(diam(\underline U)),$$
  for   some continuous  function  $\xi(d)$, $d\geq 0$  and all $  \underline U
\subset \underline X$.
\vspace {1mm}

{\it Bounded Geometry along  Spine.} A  Riemannian manifold $X$ with a spine $S$ is said to have    {\it bounded $C^\infty$-geometry along $S$} if there are continuous functions $\xi_i(d) $  and $\xi_o$ such the $i$-th
covariant derivatives of the curvature tensor  of $X$ satisfy
$$ ||\partial_icurv(X,x)||\leq \xi_i(dist(x,S))\mbox {   and } \frac {1} {inj.rad(X_i,x)}\leq \xi_o(dist(x,S)).
 \leqno   {\square_{bnd}}$$

\vspace {1mm}

\textbf {Lemma: $\mathbb R^\rtimes$-Symmetrization of Manifolds with Spines}. {\sf  Let $X$ be a complete connected orientable Riemannian $n$-manifold with a spine
$S\subset X$ 
and let $f:
X\to \underline X =\underline Y\times \mathbb R^1$ be a uniformly   S-quasi-proper  1-Lipschitz  map. 

Let the  scalar curvature of $X$ be bounded from below in terms of the  distance function 
$d(x)= dist(x,S),$ 
$$Sc(X,x)\geq \sigma(d(x))$$
for some continuous monotone decreasing function $\sigma(d)$ $d>0$.}

{\it  If $n=dim(X)\leq 7$ and if $X$ has bounded $C^\infty$-geometry along $S$,\footnote{This $C^\infty$  is a minor technicality:  the geometry which is actually used in the proof  below  is  that of the curvature itself and of  the  injectivity radius,  where even these maybe redundant.} then there exists  a smooth connected complete Riemannian  warped product  $n$-manifold $X_1=(Y_1\times \mathbb R^1, dy^2+\phi(y)^2dt^2)$ with a  $\mathbb R^1$-invariant  spine  $S_1\subset X_1$ and with a  uniformly   $S_1$-quasi-proper $\mathbb R^1$-equivariant  1-Lipschitz  map 
$$f_1:X_1=Y_1\times\mathbb R\to\underline X=\underline  Y\times\mathbb R^1$$ 
 for the obvious action of the group $\mathbb R^1$ on both spaces, 
such that   } {\sf 

$\bullet_{bnd}$ $X_1$ has bounded $C^\infty$  geometry along $S_1$;

 $\bullet_{Sc}$ the scalar  curvature of $X_1$ is bounded from below by the same function $\sigma(d) $ as the the scalar curvature of $X$, 
$$Sc(X_1,x_1)\geq \sigma(dist(x_1,S_1)).$$

$\bullet_{f_1}$ the topology of the map $f_1$ is "essentially the same" as that of $f$, where, in our   case,
we shall  need two specific  instances of this:

$\bullet_{deg}$  if  $dim(X)=dim(\underline X)$, then the map $f_1$ has the {\it same degree} as $f$;

$\bullet _{SYS}$ if  $dim(X)=dim(\underline X)+2$ and if  the homology class of the $f$-pullbacks of  generic points, $f^{-1}(\underline x)\subset X$,  $\underline x\in \underline X$,  is {\it spine detectably non-spherical}, i.e. all surfaces $\Sigma\subset X$ in this class contain closed curves in the {\sf  intersection 
$ \Sigma\cap  S$},  which are non-contractible in $X$ then the    the homology class of $f$-pullbacks of  generic points, $f_1^{-1}(\underline x)\subset X_1$, is also    {\it spine detectably non-spherical.}}

\vspace{1mm}

 {\it Proof.} Apply  $\mu$-bubble separation theorem  from section \ref{separating3} to the bands 
 $X_{[-d,d]}\in X$ that are the pullbacks of the bands $\underline Y\times [d,d]\subset \underline Y\times \mathbb R^1 $,  
 $$X_{[-d,d]}=f^{-1}(\underline Y\times [d,d]$$ 
  for the segments 
 $[d,d]\subset \mathbb R^1$, $d>0$,  and  thus obtain hypersurfaces $Y=Y(d)\subset X_{[-d,d]}\subset X$
 and warping functions $\phi_d(y)$, such  that the manifolds $X^\rtimes(d)=(Y(d)\times \mathbb R^1,  dy^2+\phi_d(y)^2dt^2)$ (obviously) satisfy all requirements of the  lemma, except for $\bullet_{Sc}$ which is replaced by an $\varepsilon_d$-weaker inequality,
 $$Sc(X_1,x_1)\geq \sigma(dist(x_1,S_1))-\varepsilon_d, $$
where $\varepsilon_d\to 0$  for $d\to \infty.$

 Now, the  $C^\infty$ -geometry of $X$ is bounded along the spine $S\subset X$, the standard elliptic estimate implies that the  $C^\infty$-geometries  of all  $X^\rtimes(d)$  are uniformly, (i.e. independently  of $d$) bounded  along the spines of these manifolds; hence, some sequence $X^\rtimes(d_i)$ Hausdorff converges to the required $X^1$.
 QED.
 \vspace {1mm}
 
Then we recall  the  "symmetry appendix" to the separation theorem and conclude that the $\mathbb R^\rtimes$-symmetrization is also compatible with extra symmetries and with the warper product structures as follows.
  \vspace {1mm}
 
\textbf {$\mathbb R^\rtimes$-Symmetrization in a Presence of a  Group Action.} {\it If the manifolds $X$ and 
$\underline Y$ are isometrically acted upon  by a group $G$,  and if the map $f:X\to \underline X=\underline Y\times \mathbb R^1$ is $G$-equivariant, then $X_1$ comes with an isometric  action 
of $G\times \mathbb R^1$ and the map $f_1: X_1\to X=\underline X=\underline Y\times \mathbb R^1$ is $G\times \mathbb R^1$-equivariant.}\vspace {1mm}

Furthermore, if 

{\it $\bullet $ $G=\mathbb R^m$;  
 
$\bullet $   $\underline Y= \underline Z\times\mathbb R^m$;

 $\bullet$  $X=(Z\times \mathbb  R^m, dz^2+ \psi(z)^2dt^2$, 
 
 then 
 $X_1=Z\times \mathbb  R^{m+1} dz^2 + \varphi(z)^2dt^2$.}

{\sf($dt^2$ stands for the  Riemannian metric in both Euclidean spaces $\mathbb R^m$  and $\mathbb R^{m+1}$.)}
\vspace {1mm}

\textbf {Corollary A.} {\sf Let $X$ be a complete Riemannian $n$-manifold with a spine $S\subset X$ and  $f:X\to \mathbb R^n$ be a uniformly $S$-quasi-proper 1-Lipschitz map. Then the  scalar curvature of 
$X$ {\it can't be  uniformly positive along $S$}}, i.e. 

{\it there is no positive  function $\sigma(d)>0$, such that 
 $Sc(X,x)\geq \sigma(dist (x,S))$, $x\in X$.} 
\vspace {1mm}

\textbf {Non-Existence/Rigidity Sub-Corollary A$'$.}  {\it If a complete orientable 
-$n$-manifold $\hat X$, $n\leq 7$,   dominates with non-zero degree a  compact orientable enlargeable manifold $X_0$,{sf  e.g.  $\hat X$ is homeomorphic to $X_0$ minus a point,} then $\hat X$  is a compact flat manifold.}
\vspace {1mm}

\textbf {Corollary B.}  {\sf Let $X$ be a complete Riemannian $n$-manifold with a spine $S\subset X$ and  $f:X\to \mathbb R^{n-2}$ be a smooth  uniformly $S$-quasi-proper 1-Lipschitz map, such that  the homology class of the $f$-pullbacks of  generic points, $f^{-1}(\underline x)\subset X$,  $\underline x\in \underline X$,  is {\it spine detectably non-spherical}.

{\it Then the scalar curvature of 
$X$ {\it can't be  uniformly positive along $S$}}.

\vspace {1mm}

\textbf {Non-Existence/Rigidity Sub-Corollary B$'$.}  {\it If a complete orientable 
$n$-manifold $\hat X$, $n\leq 7$,   dominates with  degree one  the product of an orientable enlargeable manifold by a SYS-manifold, then $\hat X$  is a compact flat manifold.}}

\vspace {2mm}

{\it \textbf   {Remarks  on Rigidity, $n>7$, and on  SYS-Enlargeable manifolds.}} 

(a) Corollaries  \textbf A and  \textbf B also extend to the case of $Sc\geq \sigma(dist(x, S)$, where  the function $\sigma$ is not strictly positive, $\sigma(d)\geq 0$,  where the conclusion is that $X$ is Riemannian flat: it is isometric to $\mathbb R^n$ , in the case   \textbf A, and to $\mathbb R^{m-2}\times T^2$, for  a flat (possibly non-split)  torus $ T^2$  for  \textbf B.

This can be proven either by adapting Jerry Kazdan's perturbation  argument or arguing as in the proof of the rigidity  of warped products in section \ref{rigidity5}

(b) If $n=8$, then the conclusion of   \textbf A and  \textbf B remain intact, since the perturbations  from Nathan Smale's 
argument  are controlled by the bound on the $C^\infty$ geometry.

Also, the rigidity  sharpening of  \textbf A and  \textbf B   remans valid,  since  the  warped product  rigidity  proof 
compensates for Smale's  perturbations. 

(c) Probably -- IF I understand the logic of Schoen-Yau's  "desingularization"  proof correctly  --  it, similarly to Smale's proof, extends to the present case and implies as much of the lemma  as is  needed for $A$ and $B$, but proving rigidity  for $n\geq 9$ seems technically more involved. 

(d) 
It is {\color {red!40!black} unclear}   if the non-domination corollary \textbf B$'$  for   {\it SYS-enlargeable} manifolds (defined below) that  are significantly more general then those in \textbf B'  
  follows from the $\mathbb   R^\rtimes $-symmetrization lemma, because of  to the   "spine detectability" condition in this lemma that can (can it?)  fail to be satisfied  in the general case.
\vspace{1mm}

{\it Definition.}  A Riemannian manifold $X$  is {\it SYS-enlargeable},   if, for all $d>0$, there exists a proper compact $n$-dimensional Riemannian  band    $X_d$   with width $width(X_d)>d$, which admits a locally isometric immersion $X_d\to X$  and such that all compact hypersurfaces $Y\subset X$, which separate $\partial _-(X_d) \subset \partial X_d$  from $\partial_ +(X_d) \subset \partial X_d$, are  
SYS, i.e.  Schoen-Yau-Schick manifolds.  

(A more   general  class of such manifolds is defined in [g(inequalities] 2018], but I  
admit, finding the  "true definition" remains {\color {red!40!black} problematic}.) 
definition.)
%%%%%%%%%%%%%%%%%%%%%%%%%

{\section {Generalisations, Speculations} \label {speculations6}

%%%%%%%%%%%%%%%%%%%%%%%%%

The most tantalizing aspect  of  scalar curvature is  that  it serves as a meeting point between  two different 
branches of analysis: {\sf the index theory and the geometric measure theory}, 

Each of the this theories, has its own domain of applicability to the scalar curvature problems 
 (summarized 
 below)  with a  significant overlaps and distinctions between the two domains.  
 
This  suggests, on the one hand,  

\hspace {20mm} {\sf a possible  unification of these two theories}

 \hspace {-0mm}and, on the other hand, 
 
  \hspace {10mm}{\sf a radical generalization, or  several such generalizations,   
  
  of the concept of a space 
 with the scalar curvature 
  bounded from below.} \vspace{1mm}

This is a dream. In what follows, we  indicate what seems  realistic, something lying  within the reach of the currently used  techniques and ideas. 
 
%%%%%%%%%%%%%%%%%%%%%%%%%%

\subsection {\color {blue} Dirac Operators versus Minimal Hypersurfaces}\label {versus6}
Let us briefly outline the relative  borders of the domains of applicability  of the two methods.\vspace{1mm}

%%%%%%%%%%%%%%%%%%%%%%%%%%%%%

1. {\color {blue} \sf \Large Spin/non-Spin.} There is no single   instance of    {\it topological obstruction}
 for a metric with $Sc>0$ on a closed   manifold $X$, the {\it universal coverings $\tilde X$} of which is {\it non-spin}\footnote {Relaxing   the condition "$X$ is spin"  to "$\tilde X$ is spin" is achieved  with (a version of) the  Atiyah $L_2$-index theorem from [Atiyah($L_2$) 1976], as  it is explained in \S\S$9\frac {1}{9 }, 9\frac {1}{8 }$ of [G(positive)  1996].}  that is obtainable  by the (known)  Dirac operator methods.\footnote{Never mind  Seiberg-Witten equation for $n=4$}  

But the minimal hypersurface method delivers  such obstructions for a class   manifolds $X$,  which admits  continuous maps $f$  to {\it aspherical spaces} $\underline X$,  such that such an  $f$   doesn't annihilate the fundamental 
class $[X]\in H_n(\underline X)$, $n=dim(X)$, i.e.  where the image  $f_\ast [X]\in H_n(\underline X)$ doesn't vanish.  

{\it Example.}  
The connected sum $X=\mathbb T^n\#\Sigma$, where $\Sigma$ is a  simply connected non-spin manifold are instance of such $X$  with the universal coverings  $\tilde X$ being non-spin.)  
\vspace{1mm}

2.  {\color {blue} \sf \Large Homotopy/Smooth Invariants.}   The  minimal hypersurface method alone can only deliver {\it homotopy theoretic} obstructions for the existence of metrics with $Sc>0$ on $X$.

But  $\hat \alpha(X)$,  non-vanishing of which  obstructs $Sc>0$ according to the results by  Lichnerowicz  and Hitchin    proven with 
 {\it untwisted} Dirac operators  is not  homotopy invariant.  (Non-vanishing of $ \hat \alpha$ is {\it the only} obstruction for $ Sc>0$  for simply connected manifolds of dimension $ \geq 5$, see section \ref{spin index3}.)

Here, observe, the spin condition is essential, but when
 it comes to  twisted Dirac operators,  those obstructions for the existence of metrics with $Sc>0$, which are {\it essentially due to twisting}  are also {\it homotopy invariant}, and, for all we know, the spin condition is redundant there.  

Furthermore,    minimal hypersurfaces can be applied together  with that    Dirac  operators.

For example the product manifold  $X=X_1\times X_2$, where 
$\hat \alpha (X_1) \neq 0$ and $X_2=\mathbb T^n\#\Sigma$, doesn't carry metrics with $Sc>-0$ , which for $dim (X)\leq 8$ follows from   Schoen-Yau's  [SY(structure) 1979] (with a use Nathan Smale's generic non-singularity  theorem for $n=8$), while the general case needs Lohkamp's  [Lohkamp(smoothing) 2018]. 

 Notice that the  twisted Dirac operator  method also applies to these, $X=X_1\times X_2$, provided that $\Sigma $ is spin, or at least, the universal covering 
$\tilde \Sigma$ is spin. \vspace{1mm}

3. {\color {blue} \sf \Large  SYS-Manifolds.} The most challenging for the Dirac operator methods is Schoen-Yau's  proof of    non-existence 
of metrics with $Sc>0$  on   Schoen-Yau-Schick manifolds (see section \ref{SY+symplectic2}}),  where the known Dirac  operator methods, even in the spin case, don't apply.  \vspace{1mm}

And as far as the topological non-existence theorems go, the  minimal hypersurface method remains silent on  the issue of metrics with $Sc>0$ on quasisymplectic manifolds $X$ as in section \ref{SY+symplectic2}}, (e.g.   closed aspherical $4$-manifolds $X$ with $H^2(X;\mathbb Q)\neq 0$.) 
And we can't
rule out  metrics with $Sc>0$  on the connected sums  $X\#\Sigma$  with any one  of the present day methods, if 
the universal coverings $\tilde \Sigma$ are non-spin.\vspace{1mm}

4. {\color {blue} \sf \Large Area Inequalities.} The main   advantage of the twisted Dirac  operator over minimal hypersurfaces is that geometric  application of the  latter to $Sc>0$   depend on lower bounds on the sizes of Riemannian  manifolds $X$, where these sizes are expressed
 in terms of the {\it distance functions} on $X$, while the twisted Dirac relies on the 
 {\it area-wise  lower bounds} on $X$.   

The simplest (very rough) result in this regard says that every (possibly non-spin) smooth manifold  $X$ admits a Riemannian metric $g_0$, such that every {\it complete}\footnote{ "Complete" is essential as it is seen already for $dim(X)=2$. But if  $area_g(S)\geq area_{g_0}(S)$ is strengthened to $g\geq g_0$ 
one can drop "complete", where the available proof goes via minimal hypersurfaces and where there is a realistic possibility of a Dirac  operator proof as well.} 
 metric $g$ on $X$, for which
$$ area_g(S)\geq area_{g_0}(S)$$
for all smooth surfaces $S\subset X$, satisfies:
$$\inf_{x\in X} Sc(g,x)\leq 0$$
(see  section 11 in [G(101) 2017]). 
More interestingly, there are  better, some of them  sharp,  bounds on  the area-wise size of manifolds with $Sc \geq \sigma>0$,  such as sharp  area inequalities in section \ref{trace-norms3} and  Cecchini's long neck theorem for maps  of manifolds with boundaries to spheres  \ref {Roe3}.
   
These     can't be obtained, in general,  with the (present day) techniques of  minimal hypersurfaces and stable $\mu$-bubbles, but   
the following    area bounds do follow by these techniques, yet they are  unapproachable  with Dirac operators.

(a)  {\it Marcus-Neves{'} $S^3$ by $S^2$-Sweeping Theorem} [Marques-Neves(min-max spheres in  3d)  2011]]  (section \ref{slicing3D.3}).

(b)  {\it Zhu{'}s $S^2\times T^n$-Systole  Theorem}  [Zhu(rigidity)  2019], (see footnote  in section \ref{twisted4})

(c) Richard{'}s $S^2\times S^2$-Systole  Theorem   [Richard(2-systoles) 2020],
(same footnote in section \ref{log-concave5}).\vspace {1mm}

5.  {\color {blue} \sf \Large Inequalities for   Metrics  Normalized by $Sc$.} Dirac  operator arguments  that yield   geometric bounds on Riemannian manifolds $X=(X,g)$  with $Sc(X)\geq \sigma>0$, e.g. on their spherical radii, 
 in terms of $\sigma$, automatically deliver in most (all?) cases similar bounds on 
 $Sc(X)\cdot X=(X, Sc(X,x)\cdot g(x))$.
 
 For instance, Llarull's  algebraic inequality  (see section \ref {Llarull4}) not just implies that 
 $$Rad_{S^n}(X/\sigma) =  Rad_{S^n}(X)/ \sqrt \sigma \leq   1/ \sqrt {n(n-1)} $$ for 
 $ \sigma =\inf_ {x\in X} Sc(X,x),$ but in fact, that 
 $$Rad_{S^n}(Sc(X)\cdot X)\leq   \sqrt {n(n-1)} =Rad_{S^n}(Sc(S^n)\cdot S^n)) $$
for {\it all}  compact spin manifolds  $X$ with positive scalar curvatures.

But it is unclear if  such  inequalities,  let them be non-sharp ones,   can be obtained  with techniques of minimal hypersurfaces and stable bubbles  and, the  bound 
$ Rad_{S^n}(Sc(X)\cdot X)\leq  const_n $ remains {\large \it \color {red!40!black}problematic} for {\it non-spin} manifolds $X$, while 
 the  inequality  $Rad_{S^n}(X/\sigma)\leq const_n$ follows with   minimal hypersurfaces (see section 12 in
 [GL(complete) 1983]  and section  \ref {log-concave5},   augmented    by the regularity results from  [Lohkamp(smoothing) 2018]  and/or   [SY(singularities) 2017]  for $n\geq 9$.

\vspace {1mm}

6. {\color {blue} \sf \Large   Families of Manifolds,  Foliations  and  Homotopies of Metrics with $Sc>0$.}  Individual index  formulas typically (always?)   extends to families of operators and deliver harmonic spinors on  members of 
appropriate families.   But there is no (apparent?)   counterpart  of this for minimal hypersurfaces and/or for stable $\mu$-bubbles that  is partly due to discontinuity of minimal subvarieties  under deformation of metrics in the ambient manifolds. 

Consequently, non-triviality of    homotopy groups (except for $\pi_0$)   of spaces of metrics with $Sc>0$ is undetectable by   minimal hypersurfaces.   
Also   the Sc-normalized (in the sense of \ref{warped stabilization and Sc-normalization2}) distance inequalities, as well 
as topological and geometric obstruction for $Sc>\sigma$ on foliations,   escape the embrace of minimal hypersurfaces.\footnote{Possibly, this can be remedied by  an  extension of  the Schoen-Yau inductive  descent  method to a class of discontinuous families.} 

\vspace {1mm}

7.  {\color {blue} \sf \Large  Non-Completeness and Boundaries.} Until recently, 
the major drawback of the Dirac operator methods  was   reliance on completeness of manifolds $X$ it applied to,\footnote {Our attempts to alleviate this limitation  in section \ref{boundaries4},  remains    unsatisfactory.} but recent results by Zeidler, Cecchini,  Lott and 
 Guo-Xie-Yu on index theorems for manifolds with boundaries\footnote{See 
  [Cecchini-Zeidler(Scalar\&mean) 2021],  [Guo-Xie-Yu(quantitative K-theory) 2020.}
 have effectively  extended the Dirac  operator  index theory  to such manifolds.

Also  minimal hypersurfaces and especially  stable $\mu$-bubbles in conjunction with twisted Dirac operators,  fare better  in non-complete  manifolds, especially  in manifolds with controlled mean curvature of their boundaries,  as it is demonstrated in section \ref{5} of this paper, but the recent  articles by John Lott [Lott(boundary 2020]  and   Christian B\"ar  with Bernhard Hanke [B\"ar]-Hanke(boundary) 2021]  open here new  possibilities for Dirac operators.

\vspace {1mm}

8. {\color {blue} \sf \Large $Sc\geq \sigma$ for $ \sigma<0$.} Both methods have more limited applications here than for $\sigma\geq 0$, where the most impressive performance  of the  Dirac operator  is in the proof of the 
Ono-Davaux spectral inequality  (stated in section \ref{negative3}}),  which also may be  seen from a more geometric  perspective of  stable $\mu$-bubbles, as it is suggested  by the {\it Maz'ya-Cheeger inequality}.
\vspace {1mm}

 9. {\color {blue} \sf \Large Singular Spaces.} Unlike Dirac operators,  minimal varieties and  $\mu$-bubbles cn be defined  for many relevant  singular spaces, such as 

(i) {\it pseudomanifolds} with piecewise linear or piecewise smooth metrics, 

(ii) {\it Alexandrov spaces} with sectional curvatures bounded from below,

(iii) {\it singular minimal hypersurfaces} and related spaces, e.g. doubles of smooth manifolds over such hypersurfaces. 

However, despite the recent progress in the  papers  [SY(singularities) 2017]  and   [Lohkamp(smoothing) 2018],  there is neither a  concept  of  $Sc\geq \sigma$ for such spaces $X$ nor comprehensive theory of minimal hypersurfaces in $X$.

And it is not clear at all if there is  room for Dirac operators on this kind of singular  spaces $X$.

%%%%%%%%%%%%%%%%%%%

\subsubsection {\color {blue} 13  Proofs of non-Existence of Metrics with  $Sc>0$  on Tori}
%%%%%%%%%%%%%%%%%%%

The present-day   proofs can be divided according the techniques they are achieved   with; these are 

A. Dirac operators. 

B. Minimal  hypersurface and stable $\mu$-bubbles. 

C. Combination of A and B.

D. Harmonic maps in dimension 3.

E. Ricci flow in dimension 3.

(I am not certain  if one can  do something   with  the Seiberg-Witten equations.)

\vspace {1mm}

In what follows,  $X$ is a Riemannin manifold diffeomorphic to $\mathbb T^n$. 
We agree that 
 two A-proofs of $Sc \ngtr 0$ on $X$ are 
 {\it different}  if they rely on different variants of the  index theorems and which   deliver  
different harmonic spinors for generic metrics  in $X$. Similarly, B-proofs are regarded  different if the relevant minimal surfaces or $\mu$-bubles are, generically, 
different.
\footnote  {Difference between  spinors and  minimal hypersurfaces  often  disappears for flat metrics on tori and also two seemingly  different spaces of spinors may, in fact, be canonically isomorphic, such as  the space of
spinors on the universal covering $\tilde X$ of an $X$ and the space of spinors on $X$ twisted with the flat bundle over $X$ with the fiber $L_2(\pi_1(X)) $ associated with the covering  $\tilde X\to X$ via the  regular representation of $\pi_1(X)$. 

I must admit I haven't systematically traced such isomorphisms in all cases and some proofs in our list can be not different after all.}

Here one notice that  all proofs based on  index theorems on  compact manifolds $X$ and  relative index theorems on complete manifolds 
have their $L_2$-{\it counterparts on Galois coverings}  $ X_\ast \to X$ that result in {\it different }harmonic  spinors\footnote {To compare spinors on different coverings of $X$ we lift them all to the universal covering $\tilde X$ of $ X$. (For general $X$, this $L_2$ has an advantage of allowing one to   relax the spin condition on $X$ 
to that on $\tilde X$.)} if the fundamental  groups $\pi_1(X_1)$ and $\pi_1(X_2)$ are {\it non-commensurable.}

Shall we regard such proofs different?\vspace {2mm}

\hspace {27mm} {\sc {\color {blue}Six}  A-Proofs with Variations }\vspace {2mm}

\textbf 1.   {\sf    Lusztig's Kind of  Proof}.   This, for $n$ even,  goes with the  family of Dirac operators 
$ \mathcal D$ on  $X$   twisted with unitary  line bundles  $l_\tau\mathbb T^n$ parametrized by the dual torus
 $hom (H_1(X)\to \mathbb T^1)\ni \tau$.

This proof can be rendered in the language of $C^\ast$-algebras (here this is the algebra of continous function on the dual torus) but, probbaly, the harmonic spinors will be the same.

If $n$ is odd, besides the reduction to the even case,  either for  $X\times \mathbb T^1$ or 
or  $X\times X$ (are these two proof different?) 
one,  probably   can proceed with the odd dimensional spectral flow argument. (I am not certain if, in a general 
$C^\ast$-algebraic K-theoretic setting, there is a distinction between what happens to  even and to  odd $n$.)

\textbf  2. {\sf    $\wedge ^2$-Hypersphericity of $\tilde X$.}  Here,  never mind odd $n$,  one uses the relative index theorem for the Dirac 
operator   on the universal covering $\tilde X$ twisted with almost flat bundles $L_\varepsilon\to\tilde X$ on $\tilde X$ with compact supports .

\textbf  3. {\sf    Infinite $K$-Area/Cowaist$_2$.} Since K-$cowaist_2(X)=\infty$, can use the ordinary index theorem on $X$ for $\mathcal D$  
 on  $\tilde X$ itself twisted with almost flat bundles over $X$.

This is close to but different from  $\mathcal D$ twisted with Mishchenko's infinite dimensinal Fredholm bundles which also yields $Sc \ngtr 0$ on tori.

\textbf  4 {\sf    Quasi-symplectic Proof.} This depends ($n$ is even) on the $L_2$-index theorem applied $\mathcal D$ on $\tilde X$   twisted  with fractional powers of a lift of a line bundle from $X$ to $\tilde X$/
(I am not certain how to arrange a spectral  flow argument for odd $n$  in this case.)

\textbf  5. {\sf    Roe's  Index Theorems.}  Since $\tilde X$ is hypereuclidean, the Roe's algebra index theorem applies to $\tilde X$. Also one may use Roe's partitioned index theorem applied to the half-cyclic cover of $X$ (homeomorphic to $ \mathbb T^{n-1}\times \mathbb R^1_+$ in our case)

\textbf 6.  {\sf   Bounds on Widths of Bands.} All  infinite  coverings $X_\ast$ of  $X$   contains  arbitrarily wide torical bands to which  the index theorem by  Zeidler-Cecchini  and by Guo-Xie-Yu apply and yield
$Sc \ngtr 0$

\vspace {2mm}

\hspace {40mm} {\sc { \color {blue} Fife } B-Proofs }\vspace {2mm}

All these proofs rely on  
{\color{magenta} S}choen-{\color{magenta} Y}au's {\color{magenta} I}nductive  {\color{magenta} D}escent with {\color{magenta}   M}inimal  {\color{magenta}H}ypersurfaces or with 
{\color{magenta} $\mu$-B}ubbles.  This is a list of these. \vspace {1mm}

 \textbf  1. {\sf    {\color{magenta} S}-{\color{magenta}Y} {\color{magenta}I}{\color{magenta}D} with  {\color{magenta}   M}inimal  {\color{magenta}H}ypersurfaces in $X$  and with Conformal 

Modification  of Metrics
via Kazdan-Warner's Theorem.}

\textbf  2.  {  \sf   {\color{magenta} MH} Inductive Descent  in $X$  with $\mathbb T^\rtimes$-Symmetrization.}

(Both 1 and 2 apply to all SYS manifolds.)

{\textbf  3.}  {  \sf {\color{magenta}I}{\color{magenta}D} in $\tilde X$ with $\mathbb T^\rtimes$-Symmetrization   of {\color{magenta}   M}inimal  {\color{magenta}H}ypersurfaces  

in large balls in $\tilde X$ with prescribed  boundaries.}

(This proof applies  to  all, possibly non-complete) manifolds with large hyperspherical radii.)
 
 \textbf  4. {  \sf  Proofs  via Bounds on Width of Bands in infinite  Coverings of $X$
 
  either with  {\color{magenta} MH}
 or with
 {\color{magenta} $\mu$-B}.} (Compare with A8.)
 
\textbf  5. {  \sf Exhaust $\tilde X$  by domains $U_i\subset \tilde X$  bounded by {\color{magenta} $\mu$-B}ubbles $Y_i=\partial U_i$ and Apply either 3 or 4  to $Y_i\rtimes\mathbb T^1$.}
(One can also  use  here  {\color{magenta}I}{\color{magenta}D} with 5 itself applied in all dimension.)

\vspace {2mm}

\hspace {40mm} {\sc { \color {blue} Two } C-Proofs }\vspace {2mm}

\textbf 1. {\sf   The above  $Y_i$ have their hyper-spherical radii 
$Rad_{S^{n-1}}(Y_i)\to\infty$, 
 that is incompatible with  $Sc(Y_i\rtimes \mathbb T^1)\geq \sigma>0$ by the index theorem for the Dirac operators on    $Y_i\rtimes \mathbb T^1$)  twisted with bundles induced from a complex vector bundle  $\underline L\to\mathbb S^{n-1}$ with  a non-zero top Chern class.}
 (You know what to do if $n$ is odd.) 

\textbf 2  {\sf   Exhaust $\tilde V$ by domains $Y'_i$ with $mean.curv (\partial U'_i)>0$  and apply Lott's
index theorem for maps from these  $U'_i$ to the hemisphere $S^{n}_+$, or use the Goette-Semmelmann's theorem for smoothed doubles of $U'_i$ mapped to $S^n$.}

(There are also variations of these proofs with exhaustions of $\tilde V$ by cubical domains 
but these, albeit especially useful in dimension $9$, where 1 and 2 don't apply,  are unbearably  artificial.)\vspace {1mm}

All these proofs, have different  possibilities for generalizations to non-torical $X$ and different ranges of applications. It would be pleasant  to find a unifying  framework  for  them.

%%%%%%%%%%%%%%%%%%%

{\subsubsection   {\color {blue}On Positivity  of $-\Delta +  const\cdot Sc$,
Kato's Inequality  and  Feynman-Kac Formula} \label{Kato6}
 %%%%%%%%%%%%%%%%%%%%%

 \textbf 1. {\it\large  \color {red!30!black} Question.} {\sf What are  effects on the topology and/or  metric geometry of 
a  Riemannian manifold $X$ played by {\it positivity} of the   
$$L_\gamma: f(x) \mapsto -\Delta f(x)  +  \gamma\cdot Sc(X,x)f(x)$$   for a given constant $\gamma>0$?}

Observe that the greater the constant $\gamma$ is, the stronger this effect should be.  Indeed, since
$-\Delta$ is a positive , 
$$-\Delta  + \gamma_1\cdot Sc(X)\geq 0 \Rightarrow  -\Delta   + \gamma_2\cdot Sc(X), \mbox { for }   \gamma_1\geq \gamma_2.$$

\vspace{1mm}

 If $\gamma= \frac {1}{2}$, then    the product $X\times \mathbb T^1$ admits a $\mathbb T^1$-invariant metric   with $Sc\geq 0,$ namely the warped product metric $g_\rtimes  (x,t)=\phi^2dx^2+dt^2$, where $\phi$  is the lowest eigenfunction of 
the  $-\Delta f(x)  + \frac {1}{2}Sc(X)$ (see section \ref{FCS1}).

Thus, all we know about geometry and topology of $\mathbb T^\rtimes$-stabilized manifolds with $Sc\geq 0$ applies to manifolds with positive $-\Delta f(x)  + \frac {1}{2}Sc(X)$.

Yet,  there can be (maybe not?)   a  difference between metric  geometries of manifolds $X$ with positive $-\Delta f(x)  + \gamma\cdot Sc(X)$ for different $\gamma\geq \frac {1}{2}$.

\vspace{1mm}

 Now, turning to small $\gamma$,  observe the following. \vspace{1mm}

  \textbf 2.   {\it All  compact smooth manifolds $X$ of dimension $n\geq 3$ admit Riemannian metrics $g$ for which the  $-\Delta_g+\varepsilon Sc(g)$ is positive for some 
$\varepsilon=\varepsilon(X)>0$.}

\vspace {1mm}

 {\it Idea  of the Proof.}  Make a "thin connected sum" of $(X,g_0)$  with a huge (volume-wise huge) topologically spherical manifold $X_\circ$, where   $Sc(X_\circ) \geq 1$ and apply the following.
\vspace{1mm}

{\it Lemma/Exercise.} {\sf Let  $s(x)$ be a continuous  function on a compact connected  manifold $X$, such that  $\int_X s(x)dx>0$, then the  $-\Delta +\varepsilon s$
is positive for all sufficiently small $\varepsilon>0$.}
\vspace {1mm}

\textbf 3. {\it {\color {red!50!black} Conjecture}.}     {\sl There is a universal $\bar\varepsilon=\bar\varepsilon _n>0$, such that all compact $n$-manifolds admit Riemannian metrics $g_\varepsilon$, for all $0\leq \varepsilon<\bar\varepsilon$,  such that  
 $$-\Delta_{g_\varepsilon}+\varepsilon  Sc(g_\varepsilon)\geq 0.$$}

One knows in this respect that if such $\bar\varepsilon _n$ does exist, then  it can't be greater  than the 
conformal  Kazdan-Warner constant, 
$$\bar\varepsilon _n \leq \gamma_n=\frac{n-1}{4(n-2)}$$
and, for all we know,  $\varepsilon _n$, may be equal to this $\gamma_n$.

But it would be more interesting to have  $\bar\varepsilon=\bar\varepsilon(X)$ as a topological invariant which 
takes infinitely many different  values on $n$-dimensional  manifolds $X$.

 If the   $-\Delta_g+\gamma_n Sc(X)$, where $\gamma_n=\frac{n-1}{4(n-2)}$
then,  by Kazdan-Warner theorem 
$X$ admits a (conformal) metric with $Sc\geq 0$. hen 

 Moreover, there may exist a universal   $\varepsilon>0$ that serves all manifolds $X$ or at least all $X$ of a given dimension  $n$, but all one can say at this point is that  
this$ \varepsilon$ must be $<\frac{n-2}{4(n-1)}.$\vspace {1mm}

{\it Remark.} If  $-\Delta+\gamma_n Sc(X)> 0$, 
then,  by Kazdan-Warner theorem,  
$X=X,(g_0)$ admits a metric $g$  (conformal to $g_0$ with $Sc(g)>0$. 
In particular, if $X$ is spin, it admits no  $g$-harmonic spinors by Lichnerowicz-Hitchin vanishing  theorem;
 thus,  $\hat\alpha(X)=0$ by the Atiyah-Singer index theorem. 

In fact,  regardless of the sign of the scalar curvature,   the existence of harmonic spinors is a conformal  invariant by Hitchin's theorem, and   
 this, applied to Dirac operators twisted with infinite dimensional  unitary (almost) flat bundles, allows an extension of most (all) Dirac operator topological obstructions to $Sc>0$  to manifolds with 
positive operators $-\Delta+\gamma_n Sc(X).$

But it feels  a bit  strange (have I confused the values of the  constants?) that  a natural alternative argument with the {\it refined Kato's inequality} (see below) and the Schroedinger-Lichnerowicz-Weitzenboeck-(Bochner) formula delivers such a conclusion {\it only}  with
 $\gamma'_n=  \frac{n-1}{4n} >\gamma_n= \frac{n-2}{4(n-1)}.$

\vspace {1mm}

\textbf  4. {\it Kato's Inequality}  [Hess-Schrader-Uhlenbrock(Kato) 1980]. {\sf Let $V\to X$ be a vector bundle with a unitary connection $\nabla$  over a Riemannian manifold and 
 let $f:X\to V$ be a smooth section.}

 Then, an elementary calculation shows that the {\it gradient of the norm of $f$ is bounded by the norm of the  covariant derivative of $f$}, 
$$\left  \arrowvert d|f|\right \arrowvert \leq |\nabla f|,$$
where this inequality at the zero points of $f$ is understood in the distribution sense.
\vspace {1mm}

 \textbf 5. {\it Corollary}. {\sf  Let  $S: V\to V$  be  a selfadjoint endomorphism,  i.e. a family of 
selfadjoint operators in the fibers, $S(x):V_x\to V_x$,  and let $s(x)$ be the lowest eigenvalue of $S(x)$.}

{\it Then the lowest eigenvalue   of the  $f\mapsto  \nabla^2 f(x)+S(x)f(x)$ is bounded from below 
 by the lowest eigenvalue of the scalar   $-\Delta +s(x)$ on $X$.}

{\sf Equivalently,} 

{\it if the  $-\Delta +s_\lambda(x) $ on $X$ for $s_\lambda=s-\lambda$ is positive for 
some real number $\lambda$,    then the 
$\nabla^2+S_\lambda(x)$  for $S_\lambda = S-\lambda$ is also positive.}

In fact, {\it non-positivity} of a  selfadjoint     means that there exists a {\it test  vector} $\phi$, such that 
$\langle A\phi,\phi \rangle<0$.

Thus , if $\nabla^2+S_\lambda$ is
  non-positive, then there exists a section $f:X\to V$,  such that
  $$0>\int_X(\langle \nabla^2f(x), f(x)\rangle+\langle S_\lambda f(x),f(x)\rangle)dx=\int_X(|\nabla f(x)|^2  +\langle S_\lambda f(x),f(x)\rangle)dx\geq 
  $$$$ \geq \int_X (|\nabla f(x)|^2  +s_\lambda |f(x)|^2) dx, $$
 and, by Kato's inequality, the norm $|f(x)|$ serves as the test function for non-positivity of 
 $-\Delta +s_\lambda$, for
$$\int_X(-\Delta |f(x)|+s_\lambda(x)|f(x)|)dx =\int_X \left  \arrowvert d|f(x)|\right \arrowvert^2+s_\lambda(x)|f(x)|^2)dx
 \leq$$$$\leq \int_X|\nabla f(x)|^2  +s_\lambda(x)|f(x)|^2)dx <0. $$

Bundles  $V$ relevant for applications to scalar curvature are spinor bundles $\mathbb S(X)\to X$ twisted with unitary bundles $L$ with  "small" curvatures, such as the following.

\textbf 6. {\it \color {blue} Example}. {\sf Let $X$ be a compact orientable Riemannian $n$-manifold, which admits a distance (or area)  decreasing map to the unit sphere $S^n$ with {\it  non-zero degree. }}

{\it If $X$ is spin, then  the lowest eigenvalue $\lambda_1(X, \frac {Sc}{4})$
of the  $f(x) \mapsto -\Delta f(x) +\frac {1}{4}Sc(X,x)f(x) $  is bounded by that for $S^n$,
$$\lambda_1\left (X, \frac {Sc}{4}\right)\leq \frac{1}{4} n(n-1).$$} \vspace {1mm}

{\it  Exercise.} Formulate and prove similar generalizations of other bounds on the size of Riemannian manifols $X$ by $\inf Sc(X)$, such as area (non)-contraction  inequalities from sections 3. 3 and 3.4. \vspace {1mm}

\textbf 7.  {\it {\color {red!50!black} Conjecture}.}  Probably, the  Dirac  operator proofs of  geometric inequalities for non-complete manifolds  with $Sc\geq \sigma>0$  
 by Cecchini,  Zeidler and Guo-Xie-Yu also extend to manifolds $X$ with lower bounds on (properly defined) eigenvalues $\lambda_1\left (X, \frac{Sc}{4}\right)$. \vspace {1mm}

{\it Remark.} No single results of this kind is available by the methods of the geometric measure theory, where one 
 faces  the following \vspace {1mm}
 
\textbf 8. {\it Open Problem.}  Find counterexamples to the following claim. 

Let $X=(X, g_0)$  be a  {\sf compact Riemannian  $n$-manifolds $X=(X,g_0)$, 
$n\neq 4$.

If  {\it the universal covering of $X$ is  non-spin}, then for all $\gamma <\frac {1}{2}$ and $\lambda>0$, there exists a  Riemannian metrics $g=g_{\gamma, \lambda}$ on $X$, such that $g\geq  g_0$,  and such that    the  $-\Delta_g+\gamma\cdot Sc(g)-\lambda$ is {\it positive}.} \vspace {1mm}

\textbf 9.{\it Exercises.} Denote by  $\lambda_1(X, \gamma Sc)=\lambda_1(X, \gamma g)$ the bottom of the spectrum of the  
  $-\Delta_g+\frac{1}{2}Sc(g)$ on a Riemannian manifold $X=(X,g)$.\vspace{1mm}

(a) Show that  $\lambda_1(X, \gamma Sc)$ is invariant under finite covering $\tilde X\to X$ of 
compact Riemannian manifolds, $$\lambda_1(\tilde X, \gamma Sc)=\lambda_1(X, \gamma Sc).$$

(b) Show that $\lambda_1(X, \gamma Sc)$ is additive under Riemannian products of manifolds, 
$$\lambda_1(X_1\times X_2, \gamma Sc)=  \lambda_1(X_1, \gamma Sc)+\lambda_1(X_2, \gamma Sc).$$

(c) {\sf Let $X=(X,g_0)$ be a  compact (possibly non-spin)  Riemannian manifold.  Show  that  there is a constant $\lambda=\lambda(X)$,  such that 
all Riemannian metrics $g$   on $X$, which are  {\it area wise greater than $g$},   i.e. such that $area_g(S)\geq area_{g_0}(S)$ for all smooth surfaces $S\subset X$ satisfy
$\lambda_1(g, \frac {1}{2}Sc) \leq \lambda.$}
  
  \vspace {1mm}

 (d) Show that  $\lambda_1(g, \gamma Sc)$ is semicontinuous under $C^0$-limits of Riemannian metrics:   
 
 {\sl if Riemannian  metrics $g_i$ uniformly converge to $g$ then 
  $$\lambda_1(g, \gamma Sc)\geq \limsup  \lambda_1(g_i, \gamma Sc)\mbox { for all }\gamma.$$}
 Also prove this for other eigenvalues of the operators $-\Delta_g+\gamma Sc$.\vspace {1mm}

\textbf {10}. {\it Questions.}  (a) Is there a geometric definition of  $\lambda_1(g, \gamma Sc)$ (and of higher eigenvalues of  $-\Delta_g+\gamma Sc$) applicable to {\it continuous} Riemannian metrics similarly to 
{\color {blue} $ \blacksquare$}  and 
 { \color {blue}\EllipseSolid}  from section \ref{Sc-criteria3}.

(b) Is there any kind of semicontinuity of the spectra of Dirac operators $\cal D:   \mathbb S(X)\to   \mathbb S(X) $  under "weak" limits of Riemannian metrics $g$ on $X=(X, g)$  and/or under  "weak" limits of connections of "twisting" vector bundles $L$ in $\mathcal D_{\otimes L}: \mathbb S(X)\otimes L\to  \mathbb S(X)\otimes L$?

 \vspace {1mm}
 
 \textbf {11}. {\it Refined Kato's Inequality} [Herzlich(Kato) 2000)], [Davaux(spectrum) 2003].  This  improves   \textbf 4 in the case where $V$ is a twisted 
 spin bundle and $f$ is in the kernel of the twisted Dirac operator as follows.   
 $$\left  \arrowvert d|f|\right \arrowvert \leq \sqrt \frac {n-1}{n} |\nabla f|.$$
 
Accordingly, the above \textbf 6 and \textbf 7 {\it \color {red!30!black} should hold} for  $\lambda_1(X, \frac {(n-1)Sc}{4n}$, but I didn't check it carefully.

 \textbf {12}. {\it  Feynman- Kac Formula.}  The Kato inequality, implies further bounds on the spectrum of  the  $\nabla^2$ that acts  on  sections of the bundle $V\to X$, by the spectrum  of  $-\Delta$, namely the inequality
 $$trace (\exp-t \nabla^2)\leq rank (V)\cdot trace (\exp  t \nabla^2),\mbox  { }t>0, $$
which follows from the point-wise inequality between the corresponding heat kernels.

 Remarkably, the latter (trivially!) follows from an identity -- {\it \color {blue!50!black}   Feynman-Kac formula}, that says that
 
 {\it the heat value }  $H_{\nabla^2} (x_0,x_1): V_{x_0}\to V_{x_1}$  {\sf is equal to the {\it average of  
the parallel transport s} from the fiber $V_{x_1}$ to   $V_{x_1}$ along "all" paths between these points, where this average is taken with respect to the {\it Wiener measure} on the space of paths between $x_0$ and $x_1$ in $X$. }\vspace{1mm}

\textbf {13.} {\it \large \color {red!50!black} Question.} {\sf Can  geometric inequalities on  scalar curvatures of Riemannian manifolds $X$, at least those proven with Dirac operators, be derived from {\it integral identities} for natural measures in spaces of maps from graphs to $X$?}

%%%%%%%%%%%%%%%%%%%

\subsection {\color {blue}  Logic  of Propositions about the  Scalar Curvature}\label {logic6}
%%%%%%%%%%%%%%%%%%%%%

Propositions/properties   $\mathcal  P|Sc$  concerning  the scalar  curvatures of Riemannian  manifolds or related invariants,  makes a kind of an "algebra", vaguely similar to how it is in algebraic topology,  where properties of  invariants    $\mathcal P|Sc$ can be  modified, generalized, stabilized in a systematic manner, e.g.  those concerning  $X$ and $Y$, can be coupled to 
 corresponding propositions, let them be only conjectural,  concerning  
   the {\sl \color {blue!30!black} Riemannian   products $X\times Y$.}

Then these hybridised propositions can be developed/generalized to statements on

 \hspace {24mm}   {\sl \color {blue!30!black} fibrations over $ Y$ with  $X$-like fibers} 
 
  \hspace {-6mm} and then further to

 \hspace {0mm}{\sl \color {blue!30!black}foliations with $X$-leaves, where a properly 
understood  (non-commutative?) 

space of leaves is  taken for $Y$.}

\vspace{1mm}

{\it \large {\color {red!30!black}Conjectural} Example}:  {\color {blue!50!black} \sf \large Lichnerowicz $\times $  Llarull  $ \times$  Min-Oo.} Let $\underline X$ be the product of the  the hyperbolic space by the   unit sphere, 
$$\underline X= \mathbf H^n\times S^n.$$
Let $X$ be a complete  orientable spin Riemannian  manifold, such that  $Sc(X)\geq 0$.
Let $f:X\to \underline X$ be a smooth proper map  with the following two  properties.\vspace{1mm}

{\sf $\bullet_{S^n}$   The $S^n$-component  $f_{S^n}: X\to S^n$  of $f$, that is the composition of $f$ with the  projection
 $\underline X=\mathbf H^n \times S^n\to S^n$,  is an  {\it area contracting}, e.g. 1-{\it Lipschitz} map.

$\bullet_{\mathbf H^n}$    { \sf The  $\mathbf H^n$-component    of $f$ is a {\it Riemannian submersion at infinity}:\vspace {0.8mm}

 the map $f_{\mathbf H^n}:X\to \mathbf H^n $ is a {\it submersion} outside a compact subset in $X$,  where the differential $ df_{\mathbf H^n}: T(X)\to T( \mathbf H^n)$  is {\it isometric} on the orthogonal complement to the kernel of  $ df_{\mathbf H^n} $.}\vspace{1mm}}

{\it Then either $Sc(X)=0$, or   the $\hat A$-genera of the pullbacks  $f^{-1}(\underline x)  \subset X$ of generic points $\underline x\in \underline X$ vanish.

 In particular, if $dim(X)=2n$ and  $Sc(X,x_0)>0$ at some $x_0\in X$,  then $deg(f)=0$.}\vspace {1mm}

 {\it Theorem  Generalisation}.  There   
  other avenues  for generalizations of  results on scalar curvature. Below we indicate directions of  some of these 
  "avenues"  mentioned  in the previous sections.

$\bullet$  from manifolds to  distance and area controlled maps between manifolds 
  $\bullet$   from closed  manifolds to manifolds with boundaries, where   the mean  curvature is bounded from below; 
  
$\bullet$   from   manifolds to manifolds with boundaries to manifolds with corners;

$\bullet$    from $(X,g)$, where $Sc(g)>0$ to  to $(X, Sc(g)\cdot g$; 
 
$\bullet$    from $X$ to $X\times \mathbb R^N$;
  
 $\bullet$   from complete     to  non-complete manifolds   with long necks;
 
 $\bullet$ from  properties of compact manifolds $Y$ with  $Sc (X)\geq \sigma$ to similar  properties of  generic point-pullbacks $Y=f^{-1} (\underline x)$ of smooth proper  distance decreasing maps $f: X\to \underline X$,  $Sc(X)\geq \sigma$  and   $\underline X$ is a "large" manifold, e.g.   $\underline X=\mathbb R^m$.
   
\vspace {1mm}

{\it Suggestions to the Reader}. Hybridize/generalize various  theorems/inequalities from  the previous as well as of the  following   sections.
 More specifically, formulate and prove whenever possible counterparts of results for $n$-dimensional manifolds  with $Sc\geq  \sigma$ to $n-N$ dimensional ones with $Sc\geq  \sigma$ and which admit an isometric (possibly non-free) action of
the torus $\mathbb T^N$.

%%%%%%%%%%%%%%%%%%%%%%%%

\subsection {\color {blue}  Almost flat  Fibrations, K-Cowaist  and  {\it max}-Scalar Curvature }\label {almost flat6}
%%%%%%%%%%%%%%%%%%%%%%%%%%%

{\sf \color {red!30!black} Much of what follows in this section and in 6.4 and 6.5 represents an attempt to  find geometric counterparts to the foliated 
 $Sc\geq 0$  non-existence theorems   based on the {\it Connes'    fibration} idea.}\footnote{See [Connes(cyclic cohomology-foliation) 1986],   [Bern-Heit(enlargeability-foliations) 2018],  [Zhang(foliations) 2016]  and  also [Zhang(foliations:enlargeability) 2018],   [Su(foliations) 2018] 
 and [Su-Wang-Zhang(area decreasing foliations) 2021] for a definite  results in his direction.}

\vspace {1mm}

Let  let $P$  and $Q$ be  Riemannian manifolds, let $ F: P\to Q$ be a smooth  fibration.
and let $\underline \nabla$ be the connection defined by 
   the   {\it horizontal tangent (sub) bundle} on $P$  that is the orthogonal complement  
to the {\it vertical}  subbundle of $T(P)$, where "vertical" means  "tangent to the fibers" called  $S_q=F^{-1}(q)\subset P$, $q\in Q$.\vspace {1mm}

{\it \large \color {blue} Problem.} {\sf Find relations between the K-cowaists$_2$ and   between  {\it max}-scalar curvatures of $P$, $Q$ and the fibers $F^{-1}(q)$ for fibrations with "small" curvatures  $|curv|(\underline \nabla)$.}\footnote{Recall that the K-cowaists$_2$ defined in section {4.1.4}  measure  area-wise sizes of spaces, e.g.  K-$cowaist_2(S)=area(S)$ for simply connected surfaces and 
K-$cowaist_2(S^n) = 4\pi$, while {\it max}-scalar curvature of a metric space $P$ defined in section \ref{max-scalar5}  is the supremum of scalar curvatures of Riemannian manifolds $X$ that are in a certain sense are   greater than $P$.} 
\vspace {1mm}

We already know in this regard the following\vspace {1mm}

(A)  If $P\to Q$  is a {\sl unitary vector bundle with a non-trivial Chern number}, then,
by its very definition, {\sl K-cowaist$_2(Q)$ is bounded from below} by $\frac {const_n}{|curv|(\underline \nabla)}$.

\vspace {1mm}

(B) There is a fair bound on  $Sc^{\sf max}$ of product spaces $P=Q\times S$, such as the  rectangular solids, for instance,  as is shown by methods  of minimal hypersurfaces and of  stable $\mu$-bubbles  in section \ref {separating5}.

\vspace {1mm}

In what follows, we say a few words  about (A) for non-unitary bundles   in the next section and then  turn to 
several  extensions of (B) to non-trivial fibrations.

%%%%%%%%%%%%%%%%%%%%%%%

\subsubsection {\color {blue} Unitarization of Flat and  Almost Flat Bundles.}\label {unitarization6}
%%%%%%%%%%%%%%%%%%%%%%%%%

Let $Q$ be a closed oriented manifold and start with the case where   $L\to Q$ is  a {\it flat} vector bundle with a structure group $G$,
e.g.  the orthogonal group $O(N_1,N_2)$.

Let some characteristic number of $L$ be non-zero, 
which means that the classifying map $f:Q\to {\sf B}(G)$   sends the fundamental class $[P]_\mathbb Q$ to a non-zero element   in  $H_n({\sf B}(G);\mathbb Q)$.\footnote {If $G$ is compact, or if $G=GL_N(C)$, then $H_n({\sf B}(G);\mathbb Q)$, then the homology  homomorphism $f_\ast:  H_i(Q,\mathbb Q)\to H_i ({\sf B}(G);\mathbb Q)$, $i>0$, for {\it  flat } bundles $L$,
 but  it is not so, for  instance, if    $G=O(N_1,N_2)$ with $N_1,N_2>0$. }

{\it Then $X$ admits no metric with $Sc >0$.} 

\vspace {1mm}

{\it First Proof.} Let $\Gamma\subset G$ be the monodromy group of $L$ and recall (see section  \ref{Tits4})
that $\Gamma$ properly and discretely acts on a product $\underline X$ of Bruhat-Tits building. Since 
this $\underline X$ is  $CAT(0)$ and $Sc(P)>0$, the homology homomorphism  $H_n(P;\mathbb Q) \to H_n({\sf B}(\Gamma);\mathbb Q)$ induced by the classifying map 
 $f_\Gamma:  P\to {\sf B}\Gamma$ is zero.

Since the classifying  map $f:Q\to {\sf B}(G)$ factors through  $f_\Gamma:  P\to {\sf B}\Gamma$ via the embedding $\Gamma\hookrightarrow G$,
 the homomorphism $H_n(P;\mathbb Q) \to H_n({\sf B}(G);\mathbb Q)$
is  zero as well and the proof follows. 
\vspace {1mm}

{\it Second  Proof?}    Let $K\subset G$  be the  maximal compact subgroup and let  $S$ be the quotient space, $S=G/K$ endowed with a $G$-invariant Riemannian metric.  

Let $\mathcal S_\ast$ be the space of $L_2$-spinors on $S$ twisted with some bundle $L_
 \ast\to S$ associated with the tangent bundle of $S$ and let $\mathscr S_\ast\to Q$ be the corresponding  Hilbert bundle over $Q$ with the fiber $\mathcal S_\ast$.
 
 Apparently,  an     argument by  Kasparov (see below) implies that, 
at least under favorable conditions on $G$, a certain generalized {\it index of the Dirac operator} on  $Q$ twisted with $\mathscr S_\ast\to Q$ is {\it non-zero}; hence,  $Q$ carries a {\it non-zero harmonic}  (possibly almost harmonic)  {\it spinor} and the proof follows by revoking the Schroedinger-Lichnerowicz-Weitzenboeck formula.\vspace {1mm}

 {\it Kasparov $KK$-Construction.} Let  
 $G$ be semisimple, 
and observe that the quotient space $S=G/K$ carries a $G$-invariant metric with non-positive sectional curvature.

Take a point $s_0\in S$ and let $\tau_{0} (s)= \tau_{s_0} (s) $ be the gradient of the distance function 
$s\mapsto dist(s,s_0)$ on $S$ regularized at $r_0$ by smoothly interpolating between  $r\mapsto dist(s,s_0)^2$ in a small ball around $
s_0$ with 
$dist(s,s_0)$ outside  such a ball.

Let  $\tau^{\tiny \bullet}_{0}  :  \mathcal S_\ast\to  \mathcal S_\ast$ be the Clifford multiplication by 
$ \tau_{0}(r), $ that is
$\tau^{\tiny \bullet}_{0}   : s\mapsto \tau_0(r) {\tiny \bullet} s$, $s\in  \mathcal S_\ast$.

\vspace {1mm}

{\it \color {blue!40!black} Discreetness  Assumption. }  Let the monodromy subgroup  $\Gamma\subset G$ be  discrete  and let us restrict the space  $\mathcal S_\ast$ and the  
$\tau^{\tiny \bullet}_{0}$ to a $\Gamma$ orbit $\Gamma(s)\subset S$  for a point $r\in R$  different from $r_0$

Then, according to an {\it observation by Mishchenko}  [Mishchenko(infinite-dimensional) 1974] the resulting  
 on the space of spinors restricted to $\Gamma(s)$,
 $$\tau^{\tiny \bullet}_{s_0,\Gamma}=\tau^{\tiny \bullet}_{s_0|\Gamma(s)}:\mathcal S_{\ast |\Gamma(s)}\to   \mathcal S_{\ast |\Gamma(s)},$$
  has the following  properties:\vspace {1mm}

{\Large  \color {blue} ($\star$)}  \hspace {2mm} {\sl $\tau^{\tiny \bullet}_{s_0, \Gamma}$ is  Fredholm;}

 \hspace {-3mm} {\Large  \color {blue} ($\star\star$)}  \hspace {1mm} {\sl $\tau^{\tiny \bullet}_{s_0,\Gamma}$ commutes with the action of $\Gamma$ modulo compact operators in the following  sense:
{\sl the  operators $$\tau^{\tiny \bullet}_{\gamma(s_0), \Gamma}-\tau^{\tiny \bullet}_{s_0,\Gamma} : \mathcal S_{\ast |\Gamma(s)}\to   \mathcal S_{\ast |\Gamma(s)}$$
are compact for all $r\notin \Gamma(s_0)$ and all  $\gamma\in \Gamma.$}}

\vspace {1mm}

These properties and the contractibility of $S$, show,  by an  elementary extension by skeleta argument  [Mishchenko(infinite-dimensional) 1974], that \vspace {1mm}

  {\Large  \color {blue} ($\star\star\star$)}    {\sl  the (graded) Hilbert bundle  $ \mathcal S_{\ast|\Gamma} \to Q$ admits a Fredholm endomorphism homotopically compatible with   {\sl $\tau^{\tiny \bullet}_{s_0, \Gamma}.$ \vspace {1mm}}}

Finally,  a  {\it  K-theoretic  index  computation}  in    [Kasparov(index) 1973], [Kasparov(elliptic)   1975]  and/or in [Mishch 1974]
yields  
\vspace {1mm}

 {\Large  \color {blue} ($\star\star\star\star$)}  {\sl non-vanishing of the index of the  Dirac  operator on $Q$ twisted with $\mathscr S_{\ast|\Gamma}$}   in relevant cases (which delivers non-zero harmonic spinors on $Q$ and the issuing $Sc(Q)\ngtr )$ conclusion in our case).
  \footnote{ The properties {\Large  \color {blue} ($\star$)} and {\Large  \color {blue} ($\star\star$)}, however simple,  establish the key link  between geometry and the index theory. These  were discovered and used by  Mishchenko in the ambience of the Novikov higher signatures conjecture and the  Hodge, rather than   the Dirac, operator on  manifolds with non-positive sectional  curvatures. 
  
  It seems,  no essentially new  geometry-analysis connection has be been discovered since, while    {\Large  \color {blue} ($\star\star\star\star$)}  grew into a fast field of the  KK-theory of   $C^\ast$-algebras  in the realm of the non-commutative geometry. }

\vspace {1mm}

Now, let us {\it drop the discreetness assumption} and  make the above  ($\Gamma$-equivariant) construction(s)    fully  $G$-equivariant.

The (unrestricted to an orbit  $\Gamma(s)\subset S$)   $\tau^{\tiny \bullet}_{0}:\mathcal S_\ast \to \mathcal S_\ast$  seems at the first sight no good for tis purpose:  

the properties    {\Large  \color {blue} ($\star$)} and   {\Large  \color {blue} ($\star\star$)} fails to be true for it, since the space   $\mathcal S_\ast$ of  $L_2$-spinors on $S$ is too large and "flabby". 

On the positive side, the  space  $\mathcal S_\ast$ may contain  a $G$-invariant subspace, roughly as large as $\mathcal S_{\ast|\Gamma}$, namely the subspace of {\it harmonic} spinors in it.
 But the   $\tau^{\tiny \bullet}_{0}$  doesn't, not even approximately, keeps this  space invariant.
However  -- this is an idea of Kasparov, I presume,  -- one can go around this problem by invoking the full Dirac  $\mathcal D:\mathcal S_\ast \to \mathcal S_\ast$, rather than   its kernel alone.

Namely, we add  the following extra structure to  $\mathcal S_\ast$:\vspace {1mm}

(A) {\sf the action of the Dirac operator  $\cal D$  or rather of   the   technically  more convenient  first order operator
 $$\mathcal E=\mathcal D(1-\mathcal D^2)^{\frac {1}{2}}: \mathcal S_\ast\to \mathcal S_\ast$$:}

(B)  {\sf the action  of continuous functions $\phi$ with compact supports in $S$.}

 These functions $\phi(s)$ act on spinors by multiplication, where this action, besides {\sl commuting with the action by $G$},
 
 \hspace {16mm}{\it  commute with $\mathcal E$ modulo compact s. }\vspace {1mm}

\vspace {1mm}

Now, because of (A) and(B), a suitably generalized index theorem applies, I guess,   and, under suitable topological conditions, yields non-zero  (almost) harmonic spinors on $Q$.\footnote{I  couldn't find any explicit 
 statement of this kind in the literature,  but it must  be buried somewhere  under several layers of   KK-theoretic formalism, which fills pages of  the  books and articles I looked into. 
 
 (In my article [G(positive) 1996]), \S$8\frac{1}{2}$,   I  mistakenly use a simplified argument
  of composing $\tau^{\tiny \bullet}_{0}$ with
a  projection on $ker(\mathcal D)$) }

\vspace {1mm}

{\it \large Problem.} {\sf Does the above (assuming it is correct) generalises to non-flat bundles $L\to Q$?}

Namely, 

{\sf is there a natural  Hilbert bundle $\mathscr S\to Q$   associated with $L$ and having its curvature bounded in terms of  that of $L$ and such that  $\mathscr S$ carries an additional structure, such as a (graded) Fredholm endomorphism, that would yield, under some topological conditions, {\it non-zero harmonic} (or almost harmonic)   
$\mathscr S$-twisted spinors on $Q$ via a suitable index theorem?}\footnote {"Almost flat"  generalizations of the "flat" Lusztig  signature  theorem are  given in   \S\S$8\frac {3}{4}, 8\frac {8}{9}$ of [G(positive) 1996]. }

\vspace {1mm} 

{\it \large Generalized  Problem.} {\sf Does the above generalizes further to fibrations with variable fibers
 with nonpositive sectional  curvatures?}\vspace {1mm} 
 
   Namely, let $F:P\to Q$ be a smooth fibrations between complete Riemannian manifolds, where
the fibers $S_q=f^{-1}(q)\subset P$  are simply connected and the induced metrics in which  have  non-positive sectional curvatures. 

Let  a connection  in this fibration be given by  a horizontal  subbundle 
$T^{hor} \subset T(P)$,   that  is the orthogonal complement to the vertical bundle -- the kernel of the differential  $dF:T(P)\to T(Q) $.

Let  $[q, q']\subset Q$  be a (short) geodesic segment between $q, q' \in Q$ and let 
$[p, p']^\sim \subset P$ be a horizontal lift of $[q, q']$.

We don't assume that the holonomy transformations $S_{q}\to R_{q'}$ are isometric 
and let

(1)  $maxdil_p(\varepsilon)$  be the supremum of the norm of the differentials of  the transformations 
$S_{q}\to S_{q'}$ at $p\in S_q$ for all horizontal path $[p, p']^{\sim}  \subset P$ of length $ \leq \varepsilon$ 
issuing from $p\in P$;

\hspace {-6mm} and 
 
(2)  $maxhol_p(\varepsilon, \delta)$ be the supremum of $dist(p,p')$   for all horizontal paths  $[p, p']^\sim$ of  length$\leq \varepsilon$, where $p'$ lies in the fiber of $p$, i.e. $F(p')=F(p)=q$ and where there is a smooth surface   $S\subset P$ the boundary of which is contained in the union of the  path $[p, p']^\sim$ and the fiber $F_q$ which contains $p$ and $p'$ and 
such that $area((S)\leq \delta^2$.

\vspace {1mm} 

{\sf Can one bound $\inf_q Sc( Q,q)$, or, more generally, {\it max-$Sc(Q)$}  in terms of bounds on the functions 
$\log maxdil_p(\varepsilon)$     and  $maxhol_p(\varepsilon, \delta)$, for all (small) $\varepsilon,\delta>0$ and all $p\in P$?}

\subsubsection{\color {blue} Comparison between  Hyperspherical Radii and $K$-cowaists of Fibered Spaces}\label {fibered spaces6}

{\color {red!40!black}\textbf A. } The  methods  of minimal hypersurfaces and of  stable $\mu$-bubbles 
 from section \ref{separating5}  that deliver fair bounds on $Sc^{\sf max}$ of product spaces $P$, such as the  rectangular solids, for instance, 
 dramatically fail (unless I  miss something obvious)  for
fibrations with  {\it non-flat connections} because of the following. 

\vspace {1mm}

{\it \color {red!29!black}   Distortion Phenomenon.} What may happen, even  for (the total spaces of) {\it unit  $m$-sphere bundles} $P$  with orthogonal connections $\underline \nabla$ over closed Riemannian manifolds $Q$, where {\it the hyperspherical radius is large, and the curvature is small}, say
$$\mbox { $Rad_{S^{n}}(Q)=1$, $n=dim(Q)$,  and $|curv|(\underline\nabla)\leq \varepsilon$,  }$$
is that,  at the same time, 
$$Rad_{S^{m+n}}(P)\leq \delta, \mbox { }  m+n=dim(P),$$
where  $ \varepsilon>0 $   and $ \delta>0$ can be {\it arbitrarily  small.\footnote {This doesn't happen if the action of the structure group on the fiber  of our fibration has bounded displacement, see  (2) in section \ref{unitarization6}.  }\vspace {1mm}}

This  possibility is due to   the fact that,  in general,  $P$ admits {\it no Lipschitz controlled retractions} to the spherical    fibers of our fibration,   even if the fibration is topologically  trivial and continuous    retractions 
(with uncontrollably large Lipschitz constants) do exits, where

{\sf  non-triviality of monodromy,  say at $q\in Q$ can  make the distance function $dist_P$ on  the  fiber  $S^{m}_q\subset P$  {\it significantly smaller} than the (intrinsic) spherical metric.}
\vspace {1mm}

{\it \color {red!29!black} Example.}  Let  $Q$ be obtained from the unit sphere $S^2$ by adding $\varepsilon$-small  handles at   finitely many points
which  are together  $\varepsilon $-dense  in $S^2$ and such that  $Q$ goes to $S^2$ by a 1-Lipshitz map  of degree one.\footnote{E.g.  let  the  handles lie outside (the ball bounded by) the sphere  $S^2\subset \mathbb R^3$  and let our  map  be the normal projection $Q\to S^2$.}

Let $P\to Q$ be  a topologically trivial  flat  unit  circle bundle, such that 
 the monodromy   rotations  $\alpha \in \mathbb T^1$ of he fiber $S_q=S^1$ around the 
 loops at $ q\in Q$ of length $\leq \delta$  are   $\delta$-dense in the group $\mathbb T^1$  for all $q\in Q$. 
 
 Then, clearly, $Rad_{S^3}(P) \leq 10\delta$, where $\delta $  can be made arbitrarily small for $\varepsilon\to 0$, whilst
 the trivial fibration  has large hyperspherical radius, namely, $Rad_{S^3}(Q\times S^1) =1$.

\vspace {1mm}

 {\color {blue} \textbf B.}  Metric distortion of the fibers of the fibration $P\to Q$    has, however,  little effect on the  K-cowaist of $P$, that can be used, instead of the hyperspherical radius,  as a measure of the size of $P$ and 
that  allows non-trivial bounds on $Sc^{\sf max}(P)$ for {\it spin} manifolds $P$ with a use of twisted Dirac operators.

In practice, to make this work, one needs vector bundles with unitary connections over  the base $Q$ and over the manifold   $S$ isometric to  the fibers  $S_q\subset P$, call these bundles 
$L_Q\to Q$ and  $L_S\to S=S_q$, where the following properties of these bundles are essential.  \vspace {1mm}

{\large \color {blue}$\bullet_{\rm I}$} {\it Monodromy Invariance of $L_S$}. The bundle $L_S\to S$, where $S$ is isometric to the fibers $S_q$ of the fibration $P\to Q$,  must be  {\it equivariant} under the action of the monodromy group $G$ of the connection $\underline \nabla$ on the fibers $S_q$  of the fibration $P\to Q$.

(Recall  that an equivariance structure   on a bundle $L$ over a $G$ space $S$ is an
 {\it equivariant  lift}  of the action of $G$ on $S$  to an action of $G$ on $L$.)

If a bundle $L_S\to L$ is $G$-equivariant,  it extends fiberwise to a bundle over $P$, call it $L_\updownarrow\to P$.

 (An archetypical example of this  is the tangent bundle $T(S)$ which extends to what is called  call the vertical  tangent  bundle  for all fibration with $S$-fibers. But, in general, actions of     groups $G$ on $S$ do not lift to vector  bundles $L\to S$. 
However,  such  lifts may become possible  for  suitably  modified   spaces $S$ and/or bundles over them.)

 \vspace {1mm}

{\large \color {blue}$\bullet_{\rm II}$} {\it Homologically Substantiality of the two Vector Bundles}. Some Chern numbers.
of the bundles $L_S$ and $L_Q$ must  be {\it non-zero}. 

{\large \color {blue}$\bullet_{\rm III}$}  {\it Non-vanishing of $F^\ast[Q]_\mathbb Q^\circ\in H^n(P;\mathbb Q)$}.  The  image     of the fundamental 
cohomology class  $[Q]^\circ\in H^n(Q)$, $n=dim(Q)$,  under  the  rational cohomology
homomorphism induced by $F:P\to Q$ doesn't vanish,
 $$F^\ast[Q]^\circ\neq 0.$$
(This is satisfied, for instance, if the fibration $P\to Q$ admits a section $Q\to P$.)   \vspace {1mm}

Granted {\large \color {blue}$\bullet_{\rm I}$}-{\large \color {blue}$\bullet_{\rm II}$}-{\large \color {blue}$\bullet_{\rm III}$}, there exists a vector bundle $L^\rtimes \to P$,  which  is equal to  a tensor product of  exterior powers of the "vertical bundle" $L_\updownarrow\to P$ and $F^\ast(L_Q)\to P$  (that is  $F$-pull back of $L_Q$) and  such that a {\it suitable} Chern number of $L^\rtimes$ doesn't vanish.

 Here  "suitable" is what ensures  {\it non-vanishing of   the index} of the twisted Dirac operators $\mathcal D_{\otimes f^\ast(L^\rtimes)}$ on manifolds $X$  mapped to $P$ by  maps  $f:X\to P$ with {\it non-zero degrees}.  (Compare with $5\frac {1}{4}$ in 
[G(positive) 2016].)

Then bounds on curvatures of the bundles  $L_S$ and $L_Q$ together with such a bound for $\underline \nabla$ 
and also a bound on  {\it parallel displacement of the $G$ action on} $S$ (see below)
yield a bound on $|curv| (L^\rtimes)$, which implies a bound on $Sc^{\sf max}(P)$ according to 
the twisted  Schroedinger-Lichnerowicz-Weitzenboeck formula applied to the operators  $\mathcal D_{\otimes f^\ast(L^\rtimes)}$ on manifolds $X$ mapped to $P$ by (smoothed) {\it  1-Lipschitz maps} 
$f:X\to P$,
used in the definition of $Sc^{\sf \max} (P)$.
\vspace {1mm}

 {\color{blue} \it  Parallel Displacement.} The geometry of  a $G$-equivariant unitary bundle $L=(L, \nabla)$ over a Riemannian  $G$-space  $S$  is characterized, besides the (norm of the) curvature of $\nabla$,  by the difference between the parallel transform and  transformations by small $g\in G$.
 
 To define this, fix a norm in the Lie algebra of $G$ and let $|g|$, $g\in G$ denote  the distance from $g$ to the identity in the corresponding "left" invariant Riemannian metric in $G$.

Then, given a  transformation $g: S\to S$ and    a  lift  $\hat g:L\to L$ of it to $L$,  compose it 
with  the parallel translate of it  back to $L$ along shortest curves (geodesics for complete $S$) between all pairs  $s, g(s)\in S$. 
Denote by $\hat g\div\nabla:L \to L$ the resulting endomorphism 
and let 
$$ |\hat G\div\nabla|= \limsup_{|g|\to 0}\frac { || (\hat g\div\nabla)-\mathbf1||}{|g|},$$ 
 where  $\mathbf 1: L\to L$ is the identity  endomorphism () and  $||...||$ denotes the  norm.

Notice at this point that  the curvature of the connection $\underline \nabla$ takes values in the Lie algebra of $G$ and the  norm $|curv|(\underline \nabla)$, similarly to the above "parallel displacement",   depends on a choice of the norm in this Lie algebra.

If $S$ is compact,  we agree to use the norm equal to the sup-norms of the corresponding vector fields on $S$, but one must be careful in the case of   non-compact $S$. (Compare with  (2) in section \ref{unitarization6}.)

%%%%%%%%%%%%%%%%%%%
  \subsubsection{ \color {blue} $Sc^{\sf max}$ and  $Sc_{sp}^{\sf max}$ for Fibrations with  Flat Connections}\label {Scmax6}
  %%%%%%%%%%%%%%%%%%%%

Let $P$ and $Q$ be closed orientable Riemannian manifolds  and let us observe that  what happens to the non-spin and spin {\it max}-scalar curvatures and of  the K-cowaists\footnote{
 K-cowaist$_2(P)$ is    the reciprocal of the {\it infimum of the norms of the  curvatures} of unitary bundles over $P$ with {\it non-zero} Chern numbers.

$Sc^{max}(P)$  is  the {\it supremum } of $\sigma$, such that  $P$ admits  an equidimensional  {\it 1-Lipschitz}  map with {\it non-zero} degree  from a closed Riemannian manifold $X$ with $Sc\geq \sigma$, and where $X$ in the definition of $Sc^{max}_{sp}$  must be {\it spin}).

The hyperspherical radius $Rad_{S^N}(P)$, $N= dim P$, the {\it supremum} $R_{max}$  of 
{\it radii of the spheres} $S^N(R)$, which receive  {\it 1-Lipshitz} maps  from $P$ of {\it non-zero} degree.  

It is (almost) 100\% obvious that $Rad_{S^N}(S^N)=1$, it is not hard to show  that
K-cowaist$_2(P)$ is $4\pi$, that
the equality $Sc_{sp}^{max}(S^N)=Sc(S^N)=N(N-1)$ follows from Llarull's' inequality for twisted Dirac operators and it remains unknown if 
$Sc^{max}(S^N)=Sc_{sp}^{max}(S^N)=Sc(S^N)=N(N-1)$ for $N\geq 5$ (see   section \ref{log-concave5} for $N=4$).}
  of     fibrations $P\to Q$ with flat connections,  follows from what we know for   trivial 
  fibrations over  covering  spaces  $\tilde Q\to Q$.
   \footnote{ {\it A flat  structure (connection)} in a  fibration $F:P\to Q$ with $S$-fibers is defined for  arbitrary topological spaces $Q, S$ and $P$, as a $\Gamma$-equivariant splitting $\tilde F:\tilde P=\tilde Q\times S\to \tilde Q$ for some  $\Gamma$-covering 
 $\tilde Q \to Q$ and the induced covering $\tilde P\to P$.
 
 In the present case we assume that our $Q$ and $ S$, hence   $P$, are  compact {\it orientable pseudomanifolds}  with piecewise smooth Riemannian metrics, where $\tilde P=\tilde Q\times S$ carries the (piecewise) Riemannian product metric and the action of $\Gamma$ on  $\tilde P$ is isometric.}

  {\color {blue}  (\textbf A)} {\sl If the monodromy group of 
 a flat fibration of)   $F:P\to Q$ is {\sf \color {blue!60!black}finite} and the map $F$ is {\sf 1-Lipschitz}, then 
 $$Sc_{sp}^{\sf max} (P)\leq const_{m+n}\cdot\max \left (\frac {1}{K\mbox {-}cowaist_2(Q)}, \frac {1}{K\mbox {-}wast_2(S) }\right),\leqno{ {\mbox {\Large $\star$}_{waist_2}}}$$
$$Sc_{}^{\sf max} (P)\leq const'_{m+n}\cdot\max \left (\frac {1}{Rad^2_{S^n}(Q)}, \frac {1}{Rad^2_{S^m}(S)}\right)\leqno{ {\mbox {\Large $\star$}}_{Rad^2}}$$
and 
$$Sc_{sp}^{\sf max} (P)\leq (m+n)(m+n-1)\cdot\max \left (\frac {1}{Rad^2_{S^n}(Q)}, \frac {1}{Rad^2_{S^m}(S)}\right)\leqno{ {\mbox {\Large $\star$}}_{sp, Rad^2}}$$
for   $n=dim(Q)$  and $ m=dim (S), $} where $S$ is the fiber of our fibration $P\to Q$.

\vspace{1mm}

In fact, these reduce to the corresponding inequalities for the product $\tilde P=\tilde Q\times S$ for the finite(!)  covering  $\tilde P$ of $ P$,  induced from  the 
{monodromy covering $\tilde Q\to Q$}, where 
 
 $\bullet$ in the case  {\Large $\star$}$_{waist_2}$, one uses the tensor product of the  relevant vector  bundles  over $\tilde Q$ and $S$ and  where the $\otimes$-product bundle  can be pushed forward from $\tilde P$ back to $P$, if one wishes so;
 
  $\bullet$ in the case   {\Large $\star$}$_{Rad^2}$,   the (obvious) inequalities  
  $$Rad_{S^{n+m}}(\tilde P)\geq Rad_{S^{n+m}}( P)$$ -- the finiteness of  monodromy   is crucial in this one  --   and 
   $$Rad_{S^{n+m}}(\tilde Q\times S)\geq min (Rad_{S^{n}}(\tilde Q),Rad_{S^{m}}( S))$$ 
allows a use of the  "cubical bounds" from the previous section,  which need no spin condition, while 
  the corresponding sharp inequality {\Large $\star$}$_{sp, Rad^2 }$ for   spin manifolds $P$
 follows from Llarull's theorem.

\vspace{1mm}

 {\color {blue}(\textbf B)} {\sl If the  monodromy group $\Gamma$ of the fibration $P\to Q$ is {\sf \color {blue!60!black}infinite}, then the above argument yields  the following modifications   of the inequalities {\Large $\star$}$_{sp, Rad^2 }$,   {\Large $\star$}$_{sp, Rad^2 }$ and {\Large $\star$$_{waist_2}$}.} \vspace{1mm}

{\Large $\star $}$^\infty_{Rad^2 }$   The two $Rad^2$  inequalities  {\Large $\star$}$_{Rad^2}$ and  {\Large $\star$}$_{sp, Rad^2 }$ for spin manifolds $P$  remain  valid for infinite  monodromy,  if {\sl $Rad_{S^n}(Q)$      is replaced  in these inequalities by $Rad_{S^n}(\tilde Q)$ for a (now infinite)  $\Gamma$-covering $\tilde Q$  of $Q$}. 

(The universal covering of $ Q$ serves this purpose  but the monodromy covering  gives an a priori   sharper result.)\vspace{1mm}

 \vspace{1mm}

 {\Large $\star$}$^\infty_{ waist_2}$ {\it One keeps  {\Large $\star$$_{waist_2}$}   valid for infinite  $\underline \nabla$-monodromy by  replacing  K-$cowaist_2(Q)$ by   K-$waist_2(\tilde Q)$.}\footnote {It is known [Brun-Han(large and small) 2009] that the hyperspherical radius can drastically  decrease under infinite coverings
 but the situation with    K-$cowaist_2$ remains unclear.}

\vspace{0mm}

{\it Remarks.}  (a) {\it Sharpening the   Constants.}  Our  argument allows improvements of  the  above   inequalities  as we shall see, at least for {\Large $\star_{sp, Rad^2}$} , in the following sections. 

(b) {\it On Displacement and Distortion.} None of  the above  inequalities contains corrections terms for  {\color {blue}parallel displacement} defined earlier in section \ref{fibered spaces6}, albeit  
it may result in a decrease of the hyperspherical radii of $P$ due to  distortion of the fibers $S\subset P$ as  the {\color {red!29!black} example}  in section \ref{fibered spaces6} shows.

Notice at this point that the presence of large distortion is inevitable for fibrations  with non-compact fibers, where the monodromy along short loops has unbounded displacement. 

{\it Example.}   Let  $Q$ be a surface and $P\to Q$ an $\mathbb R^2$-bundle with an  orthogonal  connection, 
the curvature form of which  doesn't vanish, and let $g$ be  a Riemannian metric on $P$ which agrees with the Euclidean metrics in the $\mathbb R^2$-fibers and such that the map $P\to Q$ is a Riemannian fibration, i.e. it is isometric on the horizontal subbundle in $T(P)$ corresponding to the connection.

The the Euclidean distance  between points  in the fibers,  
$$\mbox {$p_1, p_2\in \mathbb R^2_q\subset P$,  $q\in Q$}$$
 is related to the $g$-distance in $P$ as follows
  $$ dist_{\mathbb R^2}(p_1,p_2)\sim (dist_P(p_1,p_2))^2 \mbox { for } dist_{\mathbb R^2}(p_1,p_2)\to \infty.$$

(This is the same phenomenon as the distortion of central subgroups in  two-step nilpotent groups.)

%%%%%%%%%%%%%%%%%%%%%%%
\subsubsection {\color {blue} Even and Odd Dimensional  Sphere Bundles} \label {even odd6}
 
%%%%%%%%%%%%%%%%%%%%%%%%%

{\it $Sc_{sp}^{\sf max}$-Bound for Sphere Bundles.} {\sf  Let $P$ and $Q$ be  closed orientable  {\sf spin} manifolds,  where $P$ serves as the total space  
of a unit  $m$-sphere bundle $F:P\to Q$ with an orthogonal connection $\underline\nabla$.}

{\it If the map $F:P\to Q$ is 1-Lipschitz\footnote
{The role of this  "1-Lipschitz" is   seen by looking at  the trivial fibrations $P=Q\times S\to Q$ and also  at {\it Riemannian} fibrations $F:P\to Q$ (the differentials of) which are {\it isometric} on the horizontal (sub)bundle. In general, when the metrics in the horizontal tangent spaces may vary,   estimates on $Sc^{\sf max}(P)$ should incorporate along with , besides $curv(\underline \nabla)$, (a  certain function of) these  metrics. (Observe, that the scalar curvature of $P$ itself is influenced  by the  first and second "logarithmic  derivatives" of these metrics.)} 
and if  the cohomology class 
 $$\mbox { $F^\ast[Q]^\circ_\mathbb Q \in H^n(P;\mathbb Q )$, $n=dim(Q)$,}$$
   doesn't vanish  (as in  {\large \color {blue}$\bullet_{\rm III}$} in section \ref{fibered spaces6}),  {\it  then the spin max-scalar curvature of $P$} (defines with spin manifolds $X$ mapped to $ P$) {\it  is bounded in terms of the hyperspherical radius  
$R=Rad_{S^n}(Q)$
and  of the norm of the curvature of $\underline\nabla$} as follows:
  $$Sc_{sp}^{\sf max}[P]\leq const\cdot(1+\underline \epsilon) \cdot\left(Sc(S^n(R))+ Sc(S^m)\right), \leqno {\color {blue} [\rtimes S^m]}$$}
where, recall, $Sc(S^n(R))=\frac {n(n-1)}{R^2}$,   $Sc(S^m)=m(m-1)$, where $const=const_{m+n}$ is a universal constant (specified later) and where
{\it $\underline \epsilon$ is a certain positive  function $\underline \epsilon=\underline \epsilon_{m+n}(\underline c)$, for 
$\underline c= |curv|(\underline\nabla)$, such that \vspace {1mm}

\hspace {27mm} $\underline \epsilon_{m+n}(\underline c)\to 0$ for $\underline c\to 0.$}

\vspace {2mm}

{\it Proof.} Start by observing that if either $m=0$ or $n=0$, then {\color {blue} $[\rtimes S^m]$}
with $const=1$ reduces to Llarull's  inequality,  which says in these terms, e.g. for $Q$,  that 
 $$Sc^{\sf max}(Q) \leq \frac {n(n-1)}{Rad^2_{n}(Q)}=Sc(S^n(R)).$$

What we need in the general  case if we want $const=1$  is a complex vector bundle $L\to P$ with non-zero top Chern number and such that the normalised curvature (defined in section \ref{warped stabilization and Sc-normalization2}.) satisfies 
  $$|curv|_{\otimes \mathbb S}(L) \leq \frac {Sc(S^n(R))+Sc(S^m(1)}{4} +const'\cdot\underline\epsilon.$$ 

Now, let $m=dim(S=S^m)$  and $n=dim(Q)$ be even and   observe  that the non-vanishing condition $F^\ast[Q]^\circ_\mathbb Q \neq 0$ always holds for {\it even dimensional} sphere bundles.

Also observe that  $S^m$  and $Q$ support   bundles needed for our purpose, call them $L_S$  and $L_Q$, where $L_S$ is the positive spinor bundle  $\mathbb S^+(S^m)\to S=S^m$  and $L_Q\to Q$
 is  induced from the spinor bundle 
$\mathbb S^+(S^n(R))$ by a 1-Lipschitz map $ Q\to S^n(R)$ with non-zero degree. 

One knows that the top  Chern numbers  of these bundle don't vanish  and,  according to Llarull's calculation, 
 $$|curv|_{\otimes \mathbb S}(L_S) =\frac {1}{4}Sc(S^m)=\frac {1}{4}m(m-1)$$
 and 
$$|curv|_{\otimes \mathbb S}(L_Q) \leq \frac {1}{4}(Sc(S^n(R))=\frac {n(n-1)}{ 4R^2}.$$

 Since the  (unitary) bundle   $L_S\to S^m$   
 is   {\it  invariant} under the action
of the spin group, that is the double covering of $SO(m)$, \footnote {This bundle is {\it not} $SO(m)$ -invariant, but I am not certain if this is truly relevant.} it  defines  a bundle $L_{\updownarrow }\to P$,   the curvature of which satisfies
$$|curv|( L_\updownarrow) =    |curv| (L_S) +O(\underline \epsilon).$$

Then  all one needs to show   is  that  the tensor product of  
$$L=L^\rtimes =L_\updownarrow\otimes F^\ast(L_Q),$$
satisfies 
 $$|curv|_{\otimes \mathbb S}(L) \leq \frac {Sc(S^n(R))+Sc(S^m(1)}{4} +const'\cdot\underline\epsilon.$$

This  follows by a multilinear-algebraic computation    similar to what goes on in  the paper by Llarull, where,   I admit,  I didn't  
carefully check this computation.

But if one  doesn't care for sharpness of  $const$,  then a  direct  appeal to the  {\color {blue}  $\bigotimes_\varepsilon$-Twisting Principle} formulated in section \ref{twisted3}
%%%%%%%%%%%%%%%%%%
{\color {red}corrected}%???
  suffices.\vspace {1mm}

{\it Remark.} Even the non-sharp version of {\color {blue}\large [{\Large$\rtimes$}$S^m$]}, unlike how it is with a non-sharp bound   $Rad_{S^n}(X) \leq const_n (\inf_x Sc(X,x))^{-\frac{1}{2}}$, $n=dim(X)$,      {\it can't be proved}   at the present moment  
{\it without  Dirac operators},  which necessitate  spin  as well as compactness (sometimes completeness)  of our manifolds. \vspace {1mm}

{\it\color {blue}  Odd Dimensions.} If $n=dim (Q)$ is odd, multiply $P$ and $Q$ by a long circle, and then either of the three arguments, used in the odd case of  Llarull's theorem which are mentioned in section \ref{area extremality3} and referred to [Llarull(sharp estimates)  1998], [Listing(symmetric  spaces) 2010] and [G(inequalities) 2018], applies here.

Now let  $n$ be even and   the dimension $m$ of the  fiber be odd. Here  we multiply  the fiber $S$, and thus $P$ by  $\mathbb R$,  
and endow the new fiber, call it $S'=S^m\times \mathbb R$ with the bundle  $L_{S'}$  over  it, which is induced by an     $O(m+1)$-{\it equivariant}   $1$-Lipschitz  map $S^m\times \mathbb  R\to S^{m+1}$, which is  {\it locally constant at infinity}.
Since the curvature  of the new fibration $P'=P\times \mathbb R \to Q$ is   equal to that of 
the original one  of $\underline \nabla$  in $P\to Q$,  the proof follows via the relative index theorem.\vspace{1mm}

{\it Remarks/Questions.} (a) Is there an alternative argument, where, instead of $\mathbb R$, one multiplies the fiber   $S$  with the circle $ \mathbb T$,  and uses,  in the spirit of Lusztig's argument,   the  obvious $ \mathbb T$-family of flat connection in it.

(b) Is there a version of   the inequality  {\color {blue}$ [\rtimes S^m]$}, which is  sharp for $|curv|(\underline \nabla)$ far from zero?  

(c)  What are $Sc_{sp}^{\sf max}$ of the Stiefel manifolds of orthonormal 2-frames in the Euclidean  $\mathbb R^n$,   Hermitian  $ \mathbb C^n$ and quaternion $\mathbb H^n$?\footnote{Notice that  $St_2(\mathbb C^2)=S^3$ and   $St_2(\mathbb H^2)=S^7$, but not all invariant metrics on   Stiefel manifolds are symmetric.
 
 Also  notice that the corresponding (Hopf) fibrations $F: P=S^3\to Q=S^2$ and  $ F:P=S^7\to Q=S^4$
 have $F^\ast[Q]^\circ =0$ in disagreement with the above condition {\large \color {blue}$\bullet_{\rm III}$}; this  makes one wonder whether  this condition is essential.}

%%%%%%%%%%%%%%%%%%
 
  \subsubsection{ \color {blue}K-Cowaist and $Sc^{\sf max}$ of Iterated Sphere Bundles, of   Compact Lie Groups and of
   Fibrations with Compact Fibers} \label {iterated sphere bundle6}

  %%%%%%%%%%%%%%%%%%%%%

  Classical compact Lie groups are  equivariantly homeomorphic to  iterated sphere bundles.
  
  For instance, $U(k)$ is equal to the complex Stiefel manifold  of Hermitian  orthonormal  $k$-frames $St_k(\mathbb C^k)$, where 
   $St_i(\mathbb C^n)$ fibers over    $St_{i-1} (\mathbb C^n)$ with fibres $S^{2(k-i)-1}$  for all $i=1,...,k$.
  
Since the rational  cohomology of  $U(k)$ is the same as of the product $S^1\times S^3\times ...\times S^{2k-1}$, these fibrations 
  satisfy the above  non-vanishing condition   {\large \color {blue}$\bullet_{\rm III}$}, which  implies by the above  
  {\color {blue}$ [\rtimes S^m]$}  that  \vspace{1mm}
  
  {\it the product
  $U(k)\times \mathbb R^k$ carries a $U(k)$-invariant bundle,  which is  trivialized at infinity,  such that the  top Chern number of it   is non-zero.}\vspace{1mm}

  This, by the argument from the previous section, delivers   \vspace{1mm}

{\sf complex  vector  bundles  with  {\it curvature controlled}  unitary connections 
  and {\it non-vanishing} top Chern classes  over  total spaces $P$ of principal $U(k)$-fibrations  $F:P\to Q$, {\it provided
    $F^\ast[Q]_\mathbb Q^\circ\neq 0$} (that is the above  {\large \color {blue}$\bullet_{\rm III}$}). }\vspace{1mm}
   
   This yields 
  
 \hspace {20mm} {\it  a  lower bound on the $K$-cowaist 
  of $P\times  \mathbb T^k$,} 
  
   which, in turn, implies,
  the following.\vspace {1mm}
  
{\it Corollary 1.} {\sf Let $F:P\to Q$ be a principal $U(k)$-fibration with a unitary connection $\underline \nabla$, where the map $F$ is 1-Lipschitz   and   $F^\ast[Q]_\mathbb Q^\circ\neq 0$.\footnote {For a  principal fibration, this is a very strong condition, saying, in effect, that 
the fibration is "rationally trivial".}

Then  
$$Sc_{sp}^{\sf max}[P]\leq const_{m+k}\cdot (1+\underline \epsilon)\cdot \left(\frac {n(n-1)}{Rad_{S^n}(Q)^2}+ const_k\right), \leqno {\color {blue} [\rtimes U(k)]},$$}

{\it where  $\underline \epsilon$ is a certain positive  function $\underline \epsilon=\underline \epsilon_{k+n}(\underline c)$, for 
$\underline c= |curv|(\underline\nabla)$, such that \vspace {1mm}

\hspace {27mm} $\underline \epsilon_{k+n}(\underline c)\to 0$ for $\underline c\to 0.$}
\vspace {1mm}

Now let us state and prove  a similar inequality for topologically trivial fibrations with {\it arbitrary compact holonomy groups $G$}. 
 
\vspace {1mm}

{\it Corollary 2.} {\sf Let $S$ and $Q$   be compact  connected orientable Riemannian manifolds  of dimensions $m=dim(S)$ and $n=dim(Q)$  and  let $G$ be a compact   isometry group of $S$ endowed with a biinvariant Riemannian metric.\footnote {If $G$ is disconnected "Riemannian" refers to the connected components of $G$.}

Let  $F_{pr}: P_{pr}\to Q$ be a  principal $G$-fibration with a 
$G$-connection $\underline \nabla$  and with a Riemannian metric on  $P_{pr}$, which agrees with our metric on the $G$-fibers,   for which the action of $G$ is isometric  and for which   the differential of the map $F_{pr}$ is  isometric on the  
$\underline \nabla$-horisontal tangent bundle $T_{hor}(P_{pr})\subset T(P_{pr})$.

Let $F:P\to Q$ be an associated  $S$-fibration  that is 
$$ P= (P_{pr}\times S)/G$$
where the quotient is taken for the diagonal action of $G$. 

Endow  $P$ with with the   Riemannian  quotient metric.}

\vspace{1mm}

{\color {blue}$ [\rtimes S_G]$}    {\it Let $F_{pr}: P_{pr}\to Q$ be a  topologically} (but not, in general geometrically) {\it trivial  fibration } (i.e. $P_{pr}=Q\times G$ with the obvious action by $G$).

 {\it    There exists a positive constant $\underline c_0$ and a function $\underline \varepsilon=\underline \varepsilon_{m+n}(\underline c)$,
 $0\leq \underline c\leq \underline c_0$, } {\it where $\underline \varepsilon\to 0$ for $\underline c \to 0$, and  such that if 
$|curv|(\underline \nabla)=\underline c\leq  \underline c_0$, then the spin max-scalar curvature of $P$ is bounded by

$$Sc_{sp}^{\sf max}[P]\leq const_\ast\cdot (1+\underline \epsilon)\cdot \left(\frac {n(n-1)}{Rad_{S^n}(Q)^2}+ \frac {m(m-1)}{Rad_{S^m}(S)^2} +const_G\right).\footnote{I apologise for the length of this statement that  is due to  so many, probably redundant,  conditions  needed for the proof.}$$}

{\it Proof.}  Embed $G$ to a unitary group $U(k)$ and let $F_U: P_U\to Q$  be the  fibration with the fiber $U=U(k)$ associated to $F_{pr}  :P_{pr}\to Q$. 

Let $P^U\to Q$ be the  fibration with the  fibers $S_q\times U_q$, $q\in Q$ and observe that 
this $P^U$ fibers over $P$ with $U$-fibers and  over $P_U$ with $S$-fibers, where the latter is a {\it trivial fibration.}

To show this it is enough to consider  the case, where $P$ is  the principal fibration $P_{pr}$  for which 
 $P^U=P_{pr}\times U$  and $P_U$ is the quotient space,  $P_U=(P_{pr}\times U)/G$
for the diagonal action of $G$. 

Then the triviality of the principal $G$-fibration  $P^U\to P_U$ is seen with the  map
 $P^U \to U= U(k)$ for  $\{G_q\times U_q\} \mapsto U_q=U $  which sends the diagonal $G$-orbits 
from all $G_q\times U_q$  to $G\subset U(k)=U$.

Thus, assuming $m=dim (S)$ is even (the odd case is handled by multiplying by the circle  as earlier)   we obtain  an {\it \color {blue!70!black} upper   bound }
on spin max-scalar curvature of $ P^U=P_U\times S$ in terms of the  $K$-cowaist of $P_U$ and $Rad_{S^m}(S)$.

On the other hand, if the fibration $P\to Q$ has   curvature  bounded by $ \underline c$, the  same applies to the induced fibration
 $P^U\to P$ with $U$-fibers, and since the (biinvariant metric in  the)  unitary group $U=U(k)$ has positive scalar curvature, the max-scalar curvature of  $P^U$ is {\it {\it \color {blue!70!black} bounded  from below}  by 
one half of that for $P$} for all  sufficiently small  $\underline c$ and when $\underline c\to 0$ these estimate converge to what happens to Riemannian product  $P=Q\times S$.

Confronting  these upper and lower bounds yields  a qualitative version of {\color {blue}$ [\rtimes S_G]$}, while completing the 
proof of   the full  
quantitative 
statement  is left to the reader.

{\it About the  Constants.} A Llarull's kind of computation  seems  to show that the above inequalities hold with   $const_{m+n}=const_\ast=1$. 

%%%%%%%%%%%%%%%%%%%%%%%%%%%%%%%%%

 \subsection {\color {blue} K-Cowaist and Max-Scalar Curvature  for  Fibration with Non-compact Fibers}\label {max-scalar  non-compact6}
%%%%%%%%%%%%%%%%%%%%%%%%%%%%%%%%%

 Let  $P\to Q$ be a Riemannian fibration where the fiber $S$ is a complete contractible  manifold with non-positive sectional curvature and such that the monodromy of the natural connection $\underline \nabla$ in this fibration (defined by the horizontal tangent subbundle $T^{hor}\subset  T(P)$) {\it isometrically} acts on $S$. \vspace {1mm}

{\it \color {blue}  Problem.} (Compare with "Generalized Problem" in section  \ref{unitarization6}.) {\sf  Is there a lower  bound on the
K-$cowaist_2(P)$ in terms of such a bound on K-$cowaist_2(Q)$ and on an upper bound on the norm of the curvature of $\underline \nabla$ that can be represented by the function   $maxhol_p(\varepsilon, \delta)$} as in (2)  of section \ref{unitarization6}?

%%%%%%%%%%%%%%%%%%%%%%%%%%%

 \subsubsection {\color {blue} Stable Harmonic Spinors and Index Theorems.} \label {stable harmonic6}

%%%%%%%%%%%%%%%%%%%%%%%%%

Our primarily interest in  such a lower  bound  is that it would yield an {\it upper bound} on the {\it proper spin max-scalar} curvature of $P$. %\footnote

  {This "proper spin max-scalar"  is defined via proper 1-Lipschitz  maps of open  spin manifolds $X$ to $P$,  section \ref{max-scalar5}  where  following recipes 
{\large \color {blue}$\bullet_{\rm I}, \bullet_{\rm II},\bullet_{\rm III}$},  from  {\color {blue} \textbf B}  
in section \ref{fibered spaces6} one has to construct 
 a (finite or infinite dimensional   graded)   with a unitary connection  vector  bundle  $ {\cal L}\to S$, which is 

{\large \color {blue}$\star_{\rm I}$} {\it invariant (modulo compact operators?) under isometries of $S$} (compare with {\large \color {blue}$\bullet_{\rm I}$} in section \ref{fibered spaces6}).

and 

 {\large \color {blue}$\star_{\rm II}$} {\it homologically substantial}, where this substantiality must generalize that of   {\large \color {blue}$\bullet_{\rm II}$} by   properly  {\it incorporating the action of the isometry group $G$ of} $S$.
  (An inviting possibility is 
  the above $L^{\otimes_N}$.) \vspace {1mm}
 
 What one eventually  needs is not  such a bundle $ {\cal L}\to S$ per se, but rather  some Hilbert space  of sections for a class of related bundles over $P$, where

(i) {\sf  a suitable {\it index theorem}, e.g. in the spirit of our the second "proof" in section \ref{unitarization6} (with   a {\it Hilbert $C^\ast$-module} $\mathscr H$  over 
the {\it reduced $C^\ast$-algebra} of  the group $G$ being utilized)},

 and  where

  (ii) {\it the Schroedinger-Lichnerowicz-Weitzenboeck formula} {\sf applies to twisted   harmonic  $L_2$-spinors delivered by such a theorem  
and provides a bound on the scalar curvature of $P$.}}

\vspace {1mm}

{\it Who is Stable?} Harmonic  spinors delivered by index theorems (and also spinors with a given asymptotic behaviour as  in 
Witten's and Min-Oo's  arguments) 
are stable under  certain deformations  (and some discontinuous modifications, such as surgeries) of   the metrics and bundles in questions, albeit the exact range of these perturbation on non-compact manifolds is not fully understood.

But the Schroedinger-Lichnerowicz-Weitzenboeck formula doesn't use, at least not in a visible way, this stability,  which is unlike 
how it is
with   stable  minimal  hypersurfaces  and stable $\mu$-bubbles.

One wonders, however, {\sf \color {blue!60!black} 

 \hspace {3mm}whether there is a common ground for these two stabilities in our context}.
%%%%%%%%%%%%%%%%%%%%%%

 \subsubsection {\color {blue} Euclidean Fibrations}\label {Euclidean fibrations6}

%%%%%%%%%%%%%%%%%%%%%%%%

Let us  indicate an elementary approach to the above   {\it \color {blue}  problem}  in the case where the fiberes $S$  of the fibration  $F:P\to Q$    are  isometric to the {\it Euclidean space.}\vspace {1mm}

(1)   Start with the case where  the  (isometric!) action of the (structure) group $G$ on the fiber $S$  of the fibration  $P\to Q$  has a fixed point, then
assume $m=dim (S) $  is even  and observe that   radial maps $S\to S^m$, which are constant  at  infinity and have 
 degrees one,  induce  homologically substantial $G$-invariant bundles $L=L_S$ bundles on $S$.

Since  $S=\mathbb R^m$,   such  maps can be chosen with arbitrarily small Lipschitz constants, thus making  the curvatures of these bundles arbitrarily small, namely, (this is obvious)  with the  supports in the  $R$-balls  $B_{s_0}(R)\subset S$,   around the fixed point $s_0\in S$ for the $G$-action and with curvatures of our (induced from $\mathbb S(S^m)$) bundles $L_S=L_{S, s_0,R}\to S$ bounded
  by $\frac {1}{R^2}$.\footnote {It suffices to  have the universal covering  $\tilde S$ of $S$    isometric to $\mathbb R^m$, where radial bundles on $\tilde S$ can be
   pushed  forward   to Fredholm bundles on $S$.}

 Then we see as earlier that  in the limit for $R\to \infty$,   the curvature of the  bundle $L_\updownarrow\to P$, which is on the fibers $S=S_q\subset P$ is equal to
  $L_S\to S$, (see  {\ {\large \color {blue}$\bullet_{\rm I}$} in { \color {blue} \textbf B} of section \ref{fibered spaces6}) will be bounded by the curvature of the connection $\underline \nabla$  on $P\to Q$, provided the map $P\to Q$ is 1-Lipschitz.\footnote{The parallel displacement  contribution to the curvature of  $L_\updownarrow$ (see  { \color {blue} \textbf B} of section \ref{fibered spaces6}))  cancels   away by an easy argument.}
\vspace {1mm}

Consequently,

{\it the K-$cowaist_2$ of $P$ is bounded from below by the minimum of the  K-$cowaist_2$ of $Q$ and the reciprocal of  the curvature $|curv|(\underline  \nabla)$}
\vspace{1mm}

(2) Next, let us deal with  the  opposite case, where the structure group $G=\mathbb R^m$, i.e. the Euclidean space  $\mathbb R^m$  acts on itself by parallel translations. 

 Then, topologically speaking, the fibration $F:P\to Q$ is trivial, but the above doesn't,  apply since  this $P\to Q$ typically admits {\it no
parallel section}.

But since the $\underline \nabla$-monodromy transformations,  that are parallel translations on the fiber
$S=\mathbb R^m$, have bounded displacements, there  exists a continuous {\it trivialization map}  
$$G: P\to Q\times \mathbb R^n,$$
 which, assuming $Q$ is compact, (obviously)  has  the following properties.

(i) The fibers $\mathbb R^m_q\subset P$ are {\it isometrically} sent by $G$ to $\mathbb R^m= \{q\}\times \mathbb R^m\subset Q\times \mathbb R^m$ for all $q\in Q$.

(ii) The composition of $G$ with the projection $Q\times \mathbb R^m \to \mathbb R^m$,
call it $$G_{\mathbb R^m}:P\to  \mathbb R^M$$
is {\it 1-Lipshitz on the large scale},
$$ dist (G_\mathbb R^m(q_1,q_2))\leq   dist (q_1,q_2)=cost.$$

It follows  by a  standard {\it Lipschitz extension}   argument, that, for an arbitrary $\varepsilon>0$,  there exists 
 a smooth  map
  $$G_\varepsilon': P\to Q\times \mathbb R^m, \mbox { } \varepsilon>0,$$
which is properly homotopic to $G$ and such that the corresponding map
 $$G'_{\varepsilon,\mathbb R^m}:P\to  \mathbb R^{m}$$
is $\lambda$-Lipschitz for $\lambda\leq m+n +\varepsilon$, where, recall, $m+n=dim(P)$

Now, the concern expressed in  {\color {red!40!black}\textbf A } of section \ref{fibered spaces6}  notwithstanding,
the $\mu$-bubble splitting argument from section 5.3 applies and shows that\vspace {1mm}

{\color {blue}(a)}  {\it the stabilized max-scalar curvature of $P$} defined  via products of  $P$ with flat tori {\it is bounded, up to a  multiplicative 
constant,  by that of $Q$.} \vspace {1mm}

Besides,  the existence of fiberwise   contracting  scalings of $P$, which  fix a given section $Q\to P$, show that 

{\color {blue}(b)} {\it if $Q$ is compact  and if $m$ is even, then the K-$cowaist_2$ of $P$  is bounded from below,  by that of $Q$.}\vspace {1mm}

Notice here, that \vspace {1mm}

{\sf  unlike most previous occasions, neither a bound on the curvature of the fibration $P\to Q$ is required, nor   the manifold $X$  in the definition of the max-scalar curvature mapped to $P$ 
need to be spin.}\vspace {1mm}

And besides dispensing of the spin condition,   one may allow here  \vspace {1mm}

{\sf non-complete manifolds  $Q$ and $X$ and/or manifolds in {\color {blue}(a)} and compact manifolds $Q$  with boundaries in {\color {blue} (b).} }\vspace {1mm}

(3) Finally, let us turn to the general case where the structure group of a fibration $P\to Q$ with the fiber 
$S=\mathbb R^m$ is the full isometry group $G$ of the Euclidean space $\mathbb R^m$. 

Recall that $G$ is a the semidirect  product, $G=O(m)\rtimes \mathbb R^m$, let $P_G\to Q$  be the principal bundle  with fiber $G$ associated with $ P\to Q$ and let $P_O\to P$ be the associated 
$O(m)$ bundle. Let 
$$ P_O\leftarrow P_G\rightarrow P$$
be the obvious fibrations.

Now,   granted a bound on the Lipschitz constant   of $F:P\to Q$  and the curvature of this fibration, we obtain

(i) {\sl  a bound on the max-scalar curvature of the space  $P_G$ in terms of such a bound on $P$}

In fact,  the curvature of the fibration  $ P_G\to P$ as well as its Lipschitz 
constant are bounded by those of  $F:P\to Q$ and our  bound (i) follows from  non-negativity of the  scalar curvature of the fiber $O(m)$  of this fibration by the (obvious) argument used  in section \ref{iterated sphere bundle6}.

Then we look at the fibrations $P_G\to P_O\to Q$ and observe that

(ii)  the fibration  $P_O\to Q$  has $O(m)$-fibers 
and, thus   the K-$cowaist_2(P_O)$ is bounded from below by  that of  $Q$  as it was shown in section \ref{iterated sphere bundle6};

(iii) the fibration $P_G\to P_O$ has $\mathbb R^m$-fibers  and the structure group $\mathbb R^m$ and, by the above (2), the K-$cowaist_2$ of $P_G$ is bounded from below by  that of  $Q$; hence 

{\it  K-$cowaist_2(P_G)$ of $P_G$ is bounded by  K-$cowaist_2(Q)$.} \vspace {1mm}

We recall   at this point  the basic bound on $Sc_{sp}^max(P_G)$ by the reciprocal of the K-$cowaist_2(P_G)$, confront (i) with (iii)  and conclude (similarly to how it was done in  section \ref{iterated sphere bundle6}) to the 
 final result of this section. \vspace {1mm}

\PencilRightDown\PencilLeftDown \hspace {1mm} {\sf Let $F:P\to Q$ be a smooth  fibration between Riemannian manifolds with  fibers $S_q=\mathbb R^m$ and a connection $\underline  \nabla$, the monodromy of which isometrically acts  on the fibers.
If  the map $F$ is 1-Lipschitz, then}

{\it  the proper  spin max-scalar curvature of $P$ is bounded in terms of the  curvature $|curv|(\underline  \nabla)$ and the reciprocal to 
K-$cowaist_2(Q)$.}

\vspace {1mm} 

{\it \large Corollary}.  {\sf Let  $Q$ admit  a  constant at infinity area decreasing map to 
$S^n$, $n=dim(Q)$, of non-zero degree.

Let
the norm of the  curvature of (the connection  $\underline\nabla$ on)  a bundle $P\to Q$ with $
\mathbb R^m$-fibers  is bounded by 
$\underline c$.}

{\sl Let   a  complete orientable  Riemannian  spin manifold $X$ of dimension $m+n$  admit a proper area decreasing map to $P$.} 

{\it Then  
$$\inf_{x\in X } Sc(X,x)\leq  \Psi (\underline c),$$
where,  $\Psi= \Psi_{m+n}$ is an effectively describable   positive function;  in fact,
 the above  proof of \PencilRightDown\PencilLeftDown \hspace {1mm} shows  that 
one may take 
$$\Psi (\underline c)= (m+n)(m+n-1)+const_{m} \underline c$$
and where, probably, $(m+n)(m+n-1)$ can be replaced by $n(n-1)$.}

%%%%%%%%%%%%%%%%%%%

\subsubsection {\color {blue} Spin  Harmonic Area of Fibrations With  Riemannian Symmetric Fibers}  \label 
{spin harmonic6}

%%%%%%%%%%%%%%%%%%

Let $S$ be a complete Riemannian manifold with a transitive isometric action of a group $G$  which equivariantly lifts  to
a  vector bundle $\mathbf L_S\to S$ with a unitary connection, such that the  integrant  in the local formula for the index of the twisted  Dirac  $\mathcal D_{\otimes L}$  doesn't vanish.
Then a certain generalized analytic $L_2$-index of $\mathcal D_{\otimes L}$ doesn't vanish as well,\footnote {As we have  already mentioned in section \ref{4.6.5},  if $S$ admits a free discrete cocompact isometric action of 
a group $\Gamma$, this is equivalent to the  non-vanishing of the index of the  corresponding  on $S/\Gamma$
[Atiyah (L2) 1976];  in general, this index is defined by Connes and Moscovici in[Connes-Moscovici($L_2-index$ for homogeneous)  1982].}  which implies the existence of non-zero harmonic $L$-twisted square summable spinors on $S$.

\vspace {1mm}

\hspace {10mm} {\it Example:  Hyperbolic and Hermitian Symmetric spaces.} \vspace{1mm}

(a) The hyperbolic space  $S=\mathbf H_{-1}^{2m}$  admits a non-zero harmonic $L_2$-spinors twisted with the spin bundle 
  $L_S =\mathbb S^+(\mathbf  H_{-1}^{2m})$  (compare with section \ref{4.6.5}).\vspace{1mm}

(b) Hermitian symmetric spaces $S$, e.g.  products of hyperbolic planes or the quotient space of 
the symplectic group $Sp(2k, \mathbb R)$ by $U(k)\subset Sp(2k, \mathbb R)$,  admit   non-zero harmonic 
$L_2$-spinors twisted with     tensorial powers  of the canonical line bundles. 
\vspace {1mm}

{\it \color {blue} Questions.} (a)  {\sf  What are (most) general {\it local} conditions on pairs $(X, L)$, where $X$\footnote { We return to the notation $X$ instead of $S$, since, in general, this $X$ doesn't have to be  anybody's  fiber.} is  complete Riemannian manifold  and  $L\to X$ is a vector bundle with a unitary connection, such that $X$ would support non-zero  $L$-twisted harmonic  $L_2$-spinors, or, at least,  the  $\mathcal D^2_{\otimes L}$ would  contain zero in      its spectrum?\footnote {Possibly, the answer is in  [NaSchSt(localization)   2001], but I haven't read this paper and the book with the same title.  }}
\vspace {1mm}

(b)  {\sf What happens, for example,  to non-vanishing (twisted) harmonic $L_2$-spinors on homogenous  spaces $X$  under (small and/or big)
non-homogeneous deformations of the  metrics on $X$? }

\vspace {1mm}

(c) {\sf Do  non-vanishing  harmonic $L_2$-spinors twisted with the spinor bundle $\mathbb S(X)$   exist on Riemannian manifolds $X$ which are bi-Lipschitz homeomorphic to even dimensional hyperbolic  spaces with constant sectional curvatures?\vspace {1mm}

(d) Do complete simply connected  Riemannian manifolds $X$  of even dimension $n$ with their sectional curvatures pinched between $-1$ and $ -1-\varepsilon_n$ for a small $\varepsilon_n>0$ carry such spinors?\vspace {1mm}}

 \vspace {1mm}

\vspace{1mm}

{\it Example of an Application.} {\sf Let a  complete   oriented Riemannian $n$-dimensional spin  manifold $X$  admit a smooth area decreasing  map to the unit sphere,
$f:X\to S^n$, such that the pullback of the oriented volume form $\omega_{S^n}$ is {\it non-negative} on $X$,
$$\frac{f^\ast(\omega_{S^n})}{\omega_X}\geq 0,$$
and the pullback of the Riemannian metric from $S^n$ to $X$ is {\it complete}, that is the $f$-images of {\it unbounded connected} curves from $X$ have {\it infinite lengths} in $S^n$.}

{\it \it {\color {red!50!black} Conjecture}. }  {\sf The scalar curvature of $X$ is bounded by that of $S^n$,  
$$ \inf _{x\in X}Sc(X,x)\leq n(n-1).$$}

{\it Remark.}  If  {\it area-decreasing} is  strengthened to {\it $\varepsilon$-Lipschitz}  for a small $\varepsilon=\varepsilon_n>0$,
then   this conjecture (without the {\it spin} assumption)  might follow in many (all?) cases by the geometric techniques of section \ref{bubbles5}.

\vspace{1mm}
\vspace{1mm}

{\it Back to Fibrations}.  Let  $F: P\to Q$  be a fibration with the fiber $S$ and the structure group $G$, let $P$ be endowed with a complete Riemannian metric  
and let $L_\updownarrow\to P$ be the natural extension of the ($G$-equivariant!)  bundle $L_S$ to  $P$ (compare with section \ref{fibered spaces6}).

Let  $L_Q\to Q$ be a vector bundle with a unitary connection.
and let 
$$L^\rtimes=F^\ast(L_Q) \otimes  L_\updownarrow\to P$$.

Conceivably there must exist (already exists)  an index theorem for the Dirac  operator on $P$ twisted with the bundle 
$L^\rtimes$ that would ensure the existence of non-zero twisted harmonic $L_2$-spinors on $P$ under favorable topological and geometric conditions.

For instance, if $Q$ is a complete Riemannian of {\it even}    dimension $n$, if the bundle $L_Q$ is induced from the spin bundle $\mathbb S^+ (S^n)$ by a smooth  constant at infinity  map $Q\to S^n$ of positive degree, if $P$ is spin and if the map $F:P\to Q$ is isometric on the horizontal subbundle in $T(P)$, then, {\it \color {blue} \large  conjecturally},\vspace {1mm}

{\sf the manifold $P$ supports a non-zero $L^\rtimes$-twisted harmonic $L_2$-spinor.}

\vspace {1mm}

In fact this easy if the fibration is flat, e.g. if  the fibration  $P=Q\times S$ and, if the curvature  of  this fibration is (very) small, then 
 a trivial perturbation argument as in section \ref{4.6.5}
yields almost harmonic spinors on large domains  $P_R\subset P$.

But what we truly wish is the solutions of the following counterparts to 
{\color {blue}(A)} and {\color {blue}(B)}  from section \ref{4.6.5}.

Let  $F: P_R\to Q_R$  be a submersion between compact Riemannian manifolds with boundaries,
where 
$$R=\sup_{p\in P}dist (p, \partial P)$$
and where the local geometries of the fibers  are $\delta$-close (in a reasonable sense) to the geometry of an above homogeneous $S$ and let $ L_{\rtimes,R} \to P_R$  be a vector bundle, also $\delta$-close (in a reasonable sense) to an above $ L_{\rtimes}$.\vspace {1mm}

{\sf {\color {blue} (A$_F$)}  When does $P_R$ support a $\varepsilon$-harmonic  $ L_{\rtimes,R}$-twisted spinor which vanishes on the boundary of $P$?

{\color {blue} (B$_F$)} When does a similar spinor exist on a manifold $\overline P_R$, which admits a  map to $P_R$ with  non-zero degree and with   a
controlled metric distorsion?}
 (See  section \ref{4.6.5} for a specific conjecture in this direction.)

%%%%%%%%%%%%%%%%%%%

\subsection {\color {blue} Scalar Curvatures of  Foliations}\label {foliations6}

%%%%%%%%%%%%%%%%%%%%%%%
Let $X$  be a smooth $n$-dimensional  manifold and $\mathscr L$ a smooth foliation of $X$  that is a smooth partition of $X$  into $(n-k)$-dimensional leaves, denoted $\mathcal L$.

Let $T(\mathscr L)\subset T(X)$  denote the tangent bundle of $\mathscr L$ and Recall that the {\it transversal (quotient)  bundle} $T(X)/T(\mathscr L)$ carries a  natural {\it leaf-wise flat  affine  connection} denoted 
$\nabla_\mathcal L^\perp$, where the parallel transport  is   called {\it monodromy.} 

This  $\nabla_\mathcal L^\perp$  can be (obviously but non-uniquely)  extended  to an actual (non-flat) connection on the bundle  $ T(X)/T(\mathscr L)\to X$, which is  called {\it Bott connection}. \vspace {1mm}

{\it Two Examples} (1)  Let  $ \mathscr L$ admit a transversal $k$ dimensional foliation, say   $ \mathscr  K$  and observe that the bundle  $ T(X)/T(\mathscr L)\to X$ is canonically (and obviously)  isomorphic to the tangent bundle $ T(\mathscr  K)$.

Thus, every  $ \mathscr  K$-leaf-wise connection in the tangent bundle $ T(\mathscr  K)$, e.g. the  Levi-Civita
 connection for a leaf-wise Riemannian metric in $ \mathscr  K$, defines a  $ \mathscr  K$-leaf-wise connection, say  $\nabla_{\mathscr  K }$ of   $ T(X)/T(\mathscr L)$ 

Then there is a unique connection on the bundle $ T(X)/T(\mathscr L)\to X$, which agrees with $\nabla_\mathcal L^\perp$ on the 
$\mathscr L$-leaves and with $\nabla_{\mathscr  K}$  on the $\mathscr  K$-leaves, that is the Bott connection.

(2)   Let  the bundle $ T(X)/T(\mathscr L)\to X$ be  topologically trivial and let  $\partial  _i :X\to  T(X)$, $i=1,...,k$, be linearly independent vector fields transversal to $\mathscr L$. Then there exists a unique Bott connection,  for which the projection
 of  $\partial_i$ to    $ T(X)/T(\mathscr L)$ is   parallel  for  the translations  along the orbits of the field $\partial  _i$ for all $i=1,...,k$.

\vspace {1mm}

In what follows, we choose  a Bott connection on the bundle   $ T(X)/T(\mathscr L)\to X$   and  denote it  $\nabla^\perp_X$.

Also we  choose  a subbundle  $T^\perp\subset T(X)$  complementary  to $T(\mathscr L)$, which, observe, is  canonically  isomorphic to $T(X)/T(\mathscr L)$, where this isomorphism is implemented  by the quotient homomorphism  
$T^\perp\subset T(X) \to T(X)/T(\mathscr L).$

With this isomorphism, we transport the connections  $\nabla_\mathcal L^\perp$ and  $\nabla^\perp_X$ from $T(X)/T(\mathscr L$ to $\nabla^\perp_X$
 keeping  the  notations unchanged. (Hopefully, this will bring no confusion.)

%%%%%%%%%%%%%%%%%%%%%%

\subsubsection {\color {blue} Blow-up of Transversal Metrics on Foliations}\label {blow-up transversal6}

%%%%%%%%%%%%%%%%%%%%

Let $g=g_\mathscr L$ be a leaf-wise Riemannian metric on the foliation $\mathscr L$, that is a positive quadratic form on the bundle $T(\mathscr L)$, let $g^\perp$ be such a form on $T^\perp$ and observe that the sum of the two
$g^\oplus=g\oplus g^\perp$ makes a Riemannian metric on the manifold  $X$.

This metric  itself  doesn't tell you much about  our foliation $\mathscr L$, but the family
$$g^\oplus_e=g\oplus e^{2}g^\perp, \mbox { }  e>0,$$    
is more   informative in this respect, especially for 
$e \to \infty$.
For instance, \vspace {1mm}

 {\color {blue} \large  [\textbf a]} {\sf  if the metric $g=g _\mathscr L$ has {\it strictly positive scalar curvature}, i.e. $Sc_g(\mathcal L)>0$
for all leaves $ \mathcal L$ of $\mathscr L$, and}, {\sf this is essential,} {\sf if  the metric $g^\perp$ 
is {\it invariant under the monodromy along the leaves  $ \mathcal L$} -- foliations which comes with such a  $g^\perp$  are called {\it transversally Riemannian}, -- then, assuming $X$ is {\it compact}, 
$$ Sc(g^\oplus_e)>0$$
{\sf  for all sufficiently large $e>0$.}}\vspace {1mm}

{\it Proof of} {\color {blue} \large  [\textbf a]}. Let $x_0\in X$,  let $\mathcal L_0=\mathcal L_{x_0}\subset X$ be the leaf which contains $x_0$ and observe   that the pairs   pointed Riemannian manifolds $(X_e,\mathcal L_0\ni x_0)$ for 
$X_e=(X,g^\oplus_e)$  converge to the (total space of the) Euclidean  vector bundle $T^\perp$ restricted to $\mathcal L_0$ 
 with the metric 
$$g_{\lim}=g_{\mathcal L_0}\oplus g_{Eu}^\perp,\leqno {\color {blue}[\oplus]}$$ 
where $g_{\mathcal L} =g_{\mathscr L}|\mathcal L_0$, where  $g_{Eu}^\perp=g_{Eu}^\perp(l)$, $l\in \mathcal L_0$, is  the  a family of the  Euclidean metrics in the fibers of  the bundle  $T^\perp|\mathcal L_0$ corresponding to $g^\perp$ on $\mathcal L_0$, and where
"$\oplus$" refers to the local splitting  of this bundle via the (flat!) connection $\nabla^\perp_{\mathscr L}|\mathcal L_0$.\footnote{The limit space $(T^\perp, g_{lim})$ can be regarded as the {\it tangent cone of $X$ at}  $\mathcal L_0\subset X$,
where the characteristic feature of this cone is its  scale invariance under  multiplication  of the metric $g_{lim}$ normally to $\mathcal L_0$ by constants.}
The  scalar   curvature of the metric $g_{\mathcal L_0}\oplus g_{Eu}^\perp$ is determined by\vspace {1mm}

{\sf the scalar curvature of the leaf $\mathcal L_0$ and the first and second (covariant)

 logarithmic derivatives of $g_{Eu}^\perp(l)$,}\vspace {1mm}

\hspace {-6mm} where  $g_{Eu}^\perp(l)$ is regarded as a function on $\mathcal L_0$ with values in the space of (positive)  quadratic forms on $\mathbb R^k$, which in the case   $g_{Eu}^\perp(l) =\varphi(l)^2g_0$  reduces to  the "higher   warped product formula" from section \ref{warped+2}:
 $$Sc(\varphi(l)^2g_0)(l,r) =Sc(\mathcal L_0)(l)  -\frac  {k(k-1)}{ \varphi^2(l )}
 || \nabla  \varphi(l)||^2-\frac {2k}{\varphi(l)}\Delta  \varphi(l),\leqno {\mbox {\Large ${\color {blue}(\star\star_{\mathcal L})}$}}$$
  where $(l, r) \in\mathcal L_0\times \mathbb R^k$ and $\Delta= \sum \nabla_{i,i}$ is the Laplace  on 
  $\mathcal L_0$.

Since, in general, these  "logarithmic derivatives" denoted $g_{Eu}^\perp(l)'/g_{Eu}^\perp(l)$ and $g_{Eu}^\perp(l)''/g_{Eu}^\perp(l)$
 are the same as of the original (prelimit) metric $g^\perp(l)$, 
it follows, that 
 $$Sc(g_{\mathcal L_0}\oplus g_{Eu}^\perp )\geq Sc( g_{\mathcal L_0})- const_n\left ( ||(g^\perp(l)'/g^\perp(l))^2||  +||g^\perp(l)''/g^\perp(l)|| \right).\leqno {\mbox {\Large ${\color {blue}(\star\star_{Sc})}$}} $$

 In particular, if $g^\perp$ is constant with respect to  $\nabla^\perp_{\mathscr L}|\mathcal L_0$, then 
 the limit metric  $g_{\lim}$ locally is the Riemannian product $(\mathscr L, g_\mathscr L)\times \mathbb R^k$ with the scalar curvature equal to that of $\mathscr L$. QED.

However obvious, this immediately implies\vspace {1mm}

 {\color {blue} \large  [\textbf a1]}  {\sl vanishing 
of the $\hat A$-genus as well as  of its products with the Pontryagin classes of $T^\perp$} for {\it transversally Riemannian} foliations on closed spin  manifolds $X$, where
 the "product part"   of this claim 
follows from the {\it twisted}  Schroedinger-Lichnerowicz-Weitzenboeck formula for the  Dirac operator
$\mathcal D_{\otimes T^\perp}$, since  the {\it curvature of the (Bott connection in the) bundle $T^\perp  \to  X$ converges to zero for 
$e \to \infty$.}

(This is not {\it formally} covered by Connes' theorem stated in section \ref{foliations3}, where the spin condition must be satisfied by 
  $\mathscr L$ rather than  $X$ itself  as it is required here;
   but it can be easily derived from  Connes' theorem.)
\vspace {1mm}

Another   equally obvious corollary of {\color {blue}$[\oplus]$} is as follows.\vspace {1mm}

 {\color {blue} \large  [\textbf a2]} {\sl  If $Sc(\mathscr L)> n(n-1)$  and if $X$ is closed orientable spin, then 
$X$ admits   no  map $f:X\to S^n$,  such that $deg(f)\neq 0$ and such that the restrictions of $f$ to the leaves of $\mathscr L$ are 1-Lipschitz.} \vspace {1mm}

But this  is not fully satisfactory,  since it it {\it \color {blue}\large  remains unclear } \vspace {1mm}

{\sf\large  if one  {\color {red!39!black}  truly  needs}  the inequality $Sc(\mathscr L)> n(n-1)$ or \vspace {1mm}

\hspace {-2mm} $Sc(\mathscr L)> (n-k)(n-k-1)$ for $n-k=dim (\mathscr L)$  will suffice}?\vspace {1mm}

{\it Exercise.} Show that $Sc(\mathscr L)>2$ does   suffice for 2-dimensional foliations.\vspace {1mm}

{\it \large Flags of Foliations.}  
Let $$\mathscr L= \mathscr L_0\prec   \mathscr L_1 \prec...\prec \mathscr L_j,$$ 
where the relation  $ \mathscr L_{i-1} \prec  \mathscr L_{i}$  signifies that  $ \mathscr L_{i} $  {\it refines}   $ \mathscr L_{i-1}$,  which means  
 the  inclusions between  their leaves,
$$\mathcal L_{i}\subset \mathcal L_{i-1},$$
and where $\mathscr L_0$ is the bottom foliation with  a single leaf  equal $X$.: 

Let $T^\perp_{i}=T^\perp_{i}\subset T(\mathscr L_{i-1})$, $i=1,2,...,j$ be transversal subbundles isomorphic to $T(\mathscr L_{i-1})/T(\mathscr L_{i})$,  let $g_j=g_{\mathscr L_j}$  be a  ${\mathscr L_j}$-leaf-wise  Riemannian metric, 
let $g^\perp_i$, $i=1...j,$,   be Riemannian metrics on $T_{i}^\perp$  and let 
$$g^\oplus_{e_1,...,e_j}=g_0\oplus e_1g^\perp_1\oplus...\oplus  e_jg^\perp_j,  \mbox { }  e_i>0, .$$\vspace {1mm}

{\color {blue} \large  [\textbf b]}   {\sl If the metrics in the quotient  bundles  $T(\mathscr L_{i-1})/T(\mathscr L_{i})$, $i=1,...,j$, 
which corresponds to 
$g^\perp_i$, are invariant under holonomies along the leaves of  $\mathscr L_j$, if  $e_i\to \infty$,  
then 
$$Sc(g^\oplus_{e_1,...,e_j})\to Sc(g_j),$$
where this convergence is uniform on compact subsets in $X$.}

\vspace{1mm}

{\it Proof}. Since 

{\it the logarithmic derivatives of  maps from Riemannian manifolds to the 

Euclidean spaces tend to zero as the metrics in these manifolds are 

scaled by 
constants $\to \infty$,}

 the above  {\Large ${\color {blue}(\star\star_{Sc})}$} implies the following. \vspace{1mm} 

{\color {blue} \large  [\textbf b$_{lim}$]}   {\sf The  pair of   pointed Riemannian manifolds $ (X_{e_1,...,e_j},  \mathcal L_j \ni  x_j)$, for all leaves $\mathcal L_j$ of $\mathscr L_j$ and all $x_j\in L_j $,
converges to the (total space of the) flat Euclidean vector bundle $T^\perp_{1}\oplus ... \oplus T^\perp_{j}\to\mathcal L_j$, where 

 the limit metric on (the total space of)  $T^\perp_1\oplus ... \oplus T^\perp_j$ locally splits as
$$g_{\lim}=g_{\mathcal L_j}\oplus g_{Eu} \otimes g_{Eu,1(l)} , \leqno {\color{blue}[\oplus_i\hspace {-1mm}\perp]}, $$  where  $g_{Eu}$ is the Euclidean metric  on $\mathbb R^{k_2+...k_i+...+k_j}$ for $k_i=rank (T^\perp_i)$  and 
$g_{Eu,1(l)}$, $l\in \mathcal L_j$
is a family of Euclidean metrics in the fibers of the bundle   $T^\perp_1\to X$ restricted to  $\mathcal L_j$,
where the logarithmic derivatives of these metrics are equal these for the original (prelimit) metrics in the bundle  
 $T^\perp_1$ over $\mathcal L_j$.}

\vspace {1mm}

Now, we see, as earlier, that {\color {blue} \large  [\textbf b$_{lim}$]} $\Rightarrow$ {\color {blue} \large  [\textbf b$_{}$]} 
and the proof follows.

\vspace {1mm}

Thus, the  above  {\color {blue} \large  [\textbf a1]}  and  {\color {blue} \large  [\textbf a2]}  generalize to  transversally  Riemannian flags of foliations 

%%%%%%%%%%%%%%%%%%%%

\subsubsection {\color {blue} Connes' Fibration}\label {Connes6}

%%%%%%%%%%%%%%%%%%%%%%%
Let the "normal"  bundle $T^\perp\to X$ to a foliation $\mathscr L$ on $X$ admits a smooth $G$-structure for a subgroup 
$G$ of the linear group $ GL(k)$, $k=codim(\mathscr L)$, which (essentially)  means that  the monodromy transformation   for the above canonical  flat leaf-wise connection $\nabla^\perp_{\mathscr L}$ are contained in $G$.

For instance, being  Riemannian for a  foliation is the same as to admit $G=O(k)$  and $G=GL(k)$  serves all foliation. 

Let $G$ {\it isometrically} act on a Riemannian manifold $S$ and let
$P\to X$ be a fibration associated to $T^\perp\to X. $

Then the monodromy of $\nabla^\perp_{\mathscr L}$ is isometric on the fibers $S_x\subset P$.\vspace {1mm}
 
{\it Principal Example.}[Con(cyclic cohomology) 1986] Let $$\mbox {$G= GL(k)$ and $S=GL(k)/O(k)$}$$
and let us identify 
the fiber $S_x$, for all $x\in X$,  with  the space of Euclidean structures, i.e. of positive definite quadratic forms, in the linear space $T^\perp_x$.

Clearly, this $S$
 canonically  splits as 
 $$\mbox {$S=R\times \mathbb R$  for   $R=SL(k)/SO(k)$},$$ 
 where, observe,  $R$  carries a unique up to scaling $SO(k)$-invariant  Riemannian (symmetric)  metrics with non-positive sectional curvature  and where the   $\mathbb R$-factor  is the logarithm of the central multiplicative subgroup
 $\mathbb R ^\times _+\subset GL(k)$.

Thus, $S=R\times \mathbb R$ carries an invariant Riemannian  product metric, call it $g_S$, which is unique up-to scaling of the factors.  

Next, observe that the tangent bundle $T(P)$ splits as usual
$$T(P)=T^{vert}\oplus T^{hor}$$

where $T^{vert}$ consists of the vectors tangent to the fibers $S_x\subset P$,  $x\in X$, 
and where  $T^{hor}$

is the horizontal subbundle corresponding to the Bott connection, and where the splitting
 $T(X)=T(\mathscr L)\oplus T^\perp$ lifts to a splitting of  $T^{hor}$, denoted 
 $$T^{hor}=\tilde T(\mathscr L)\oplus \tilde T^\perp.$$
 Thus, the tangent bundle $T(P)$ splits into sum of three bundles,
 {\color {blue}$$T(P)=T^{vert}\oplus \tilde T(\mathscr L)\oplus \tilde T^\perp,$$}
where, to keep track of things, recall that
  $$  rank (\tilde T(\mathscr L))=dim (\mathscr L)=n-k, \mbox { } rank  (\tilde T^\perp)=codim(\mathscr L)=k$$
 and $$rank (T^{vert})=dim (GL(k)/O(k))=\frac{ k(k+1)}{2}.$$
 
Let us record the  essential features  of these three  bundles  and their roles in  the geometry of the space $P$ (see [Connes(cyclic cohomology-foliation) 1986] and compare with   \S$1\frac {7}{8} $ in [G(positive) 1996]). \vspace {1mm}

  {\color {blue}  (1)   \it Metric $\tilde g^\perp$ in $\tilde T^\perp$.}
The (sub)bundle $\tilde T^\perp \subset T(P)$ carries a {\it tautological  metric} call it $\tilde g^\perp$, 
which, in the fiber $\tilde T^\perp_p\subset \tilde T^\perp$ for $p\in P$ over $x\in X$,   is equal to this very $p\in  P_x$ 
 regarded as  a metric in $T^\perp_x\subset T^\perp\to X$.\vspace {1mm}

  {\color {blue} (2) \it Foliation ${\mathscr L}^+$ of $P$.} The  leaves $\mathcal L^+\subset P$  of this foliations are   the pullbacks of the leaves
$\mathcal L$ of $ {\mathscr L}$ under the map $P\to X$.
These  $\mathcal L^+$ have dimensions $n-k+\frac{ k(k+1)}{2}$ and the tangent bundle $T({\mathscr L}^+)$ is canonically isomorphic to  $\tilde T(\mathscr L)\oplus T^{vert}.$ \vspace {1mm}

  {\color {blue} (3) Foliation $\tilde {\mathscr L}$  of $P$.} This is the   natural lift of the original foliation $\mathscr L$ of $X$: 
  
\vspace {1mm}

{\sf  
the leaf $\tilde {\mathcal L}_p$  of   $\tilde {\mathscr L}$  through a  given point $p\in P$ over an $x\in X$ is equal 
to the  set of the Euclidean metrics in the fibers $T^\perp_l\subset T^\perp\to X$ for all  $l\in \mathcal L_x \subset X$, which are obtained from $p$, regarded as such a metric in $T^\perp_x\subset T^\perp\to X$,  by the monodromy along the leaf $\mathcal L_x$ of the foliation $\mathscr L$ of $X$.}\vspace {1mm}

This foliation    can be equivalently  defined via its tangent (sub)bundle,  that is   

$$T(\tilde {\mathscr L})= \tilde T(\mathscr L)\subset T(P).$$
Also observe that this $\tilde {\mathscr L}$ refines  ${\mathscr L}$, written as $\tilde {\mathscr L}\succ {\mathscr L}^+$, 
where, in fact,  the leaves of $ {\mathscr L}^+$ are products of the  monodromy covers of the leaves of  ${\mathscr L}$ by $S$.

\vspace {1mm}

 {\color {blue} (4)   $\tilde {\mathscr L}$-Monodromy Invariance of the Metric  $\tilde g^\perp$.}
 The bundle $\tilde T^\perp \subset T(P)$, where the metric $\tilde g^\perp$  resides, is naturally isomorphic to the "normal" bundle $T(P)/T(\mathscr L^+)$, but this metric is {\it not invariant} under the monodromy of the foliation $ \mathscr L^+$. 
 
   However,  $\tilde g^\perp$ is {\it invariant  under the monodromy of the sub-foliation} $\tilde {\mathscr L}\succ \mathscr L^+$  with the leaves  $\tilde {\mathcal L}\subset  \mathcal L^+$  as it follows from the above description of the leaves $\tilde{\mathcal L}_p$ of $\tilde{\mathscr L}$.
\vspace{1mm}

  {\color {blue} (5)   $\tilde {\mathscr L}$-Monodromy Invariance of   $\tilde g_S$ in the Bundle  $T^{vert} $}.  Since the fibration $P\to X$  with the fiber $S=GL(k)/O(k)$ is associated with $T^\perp\to X$, every 
$ GL(k)$ metric $g_S$ on $S$ gives rise to a {\it monodromy invariant metric in the fiberes of this   fibration}, which is   denoted $\tilde g_S$ and regarded as the metric in the subbundle  $T^{vert}\subset T(P)$, which made of the vectors tangent to the $S$-fibers and which is canonically isomorphic to $T(\mathscr L^+)/T(\tilde { \mathcal L})$.

Clearly, 

{\it this metric  $\tilde g_S$ is invariant under the monodromy along the leaves of the foliations  $\tilde {\mathscr L}$ on $P$}.

\vspace {1mm}

  {\color {blue} (6) \it  Scalar Curvature   under Blow-up of  Metrics in $T(P)$.}  Let $g=g_{\mathscr L}$   be a Riemannian metric in the tangent bundle $T(\mathscr L)\subset T(X)$ of  a foliation $\mathscr L$ of $X$ as earlier  and let $\tilde g$ be its lift to the bundle 
$ \tilde T(\mathscr L)=T(\tilde {\mathscr L})\subset T(P)$.

Let $\tilde g_{e_S, e_\perp}$, $e_S, e_\perp>0$, be the Riemannian metric on the manifold $P$ that is the metric in the bundle 
$$T(P)= \tilde T(\mathscr L)\oplus T^{vert} \oplus \tilde T^\perp.$$ 
where this  $\tilde g_{e_S, e^\perp}$ is split into the sum  
of the  metrics  from th above {\color {blue}(5)} and  {\color {blue} (4)}. which are  taken here 
with (large)  positive $e$-weights  as follows. 
{\color {blue} $$\tilde g_{e_S,e^\perp}=\tilde g+ e_S\tilde g_S+e_\perp \tilde g^\perp.$$ }

  Then it follows from the above {\color {blue} \large  [\textbf b]}, that
  if $$  e_S,  e_\perp  \to \infty,$$
  
 then

{\color {blue} [$\Uparrow_{Sc}$]}   {\it the scalar curvature of the metric  $\tilde g_{e_S, e_\perp}$ at $p\in P$ over $x\in X$ converges to that of $g$ on the leaf $\mathcal L_x\ni x$ at $x$, where

this convergence is uniform on the compact subsets in $ P$.}\footnote {This convergence property, which is  implicit in [Connes(cyclic cohomology-foliation) 1986], is used   in \S$1\frac{7}{8}$ of [G(positive) 1996] and    in "adiabatic" terms in Proposition 1.4 of [Zhang(foliations) 2016], where  
 it is required  that  $ e_\perp/ e_S$ is large,  since the  shape of the  
compact domain in $P$ where    {\sf the scalar curvature of the metric  $\tilde g_{e_S, e_\perp}$  becomes $\varepsilon$-close to that of $g$,} depends on the ratio $ e_\perp/ e_S$, (see section   \ref{Hermitian Connes6}.)}

\vspace {1mm}

{\it Generalizations.} Much of the above    {\color {blue} (1)} -  {\color {blue} (6)} applies to   foliations with monodromy groups $G$ not necessarily equal to $GL(k)$ and with fibrations with the fibers   that may be  different from $G/K$, which we
will approach in the following sections on the case-by-case basis.

%%%%%%%%%%%%%%%%%%%%%%
\subsubsection {\color {blue} Foliations with Abelian Monodromies}\label {Abelian Monodromies6}

 %%%%%%%%%%%%%%%%%%%%
 {\sf Let a foliation $\mathscr L$ of  an orientable  $n$-dimensional  Riemannian  manifold  $X$ admit a smooth $G$-structure   invariant under the monodromy, where the group $G$ is Abelian  and let the scalar curvatures of the leaves with the indices Riemannian metrics are bounded from below by $\sigma>n(n-1)$.}

{\color {blue} \Large  $\square$}{\small $\Circle $}. {\it The hyperspherical radius of $X$ is bounded by one,
$$ Rad_{S^n}(X)\leq 1.$$}
That is,  if $R>1$, then  

{\it $X$ admits no 1-Lipschitz map to the sphere $S^n(R)$, which  is constant 

at infinity and which  has non-zero degree.}

\vspace {1mm} 

Prior to turning to the proof, that is an easy corollary of what we discussed  about  $\mathbb R^k$-fibration in section \ref{Euclidean fibrations6}, we'll clarify a couple of  points.  \vspace {1mm}

1. We don't assume here that the manifold $X$ is compact or  complete, nor do we require it is being spin.\vspace {1mm}

2. We don't know if our Abelian assumption on $G$ is  essential.
It is {\color {red!50!black} conceivable} that \vspace {1mm}

{\sf {\color {blue} \Large  $\square$}{\small $\Circle $} holds for all   foliations,  i.e. for  $G=GL(k)$,  $k=codim(\mathscr L)$, and, moreover,  with the  
 bound  $Sc(\mathscr L)\geq (n-k)(n-k-1)$.} \vspace {1mm}

2. Examples of foliations with Abelian $G$, include: \vspace {1mm} 

 {\sf  foliations with   transversal conformal structure, e.g. (orientable) foliations of codimension one, where $G$ is the multiplicative group $\mathbb R^\times$;

 flags of codimension one  foliations  (where $G= (\mathbb R^\times)^k$)  and/or of  foliations with   transversal conformal structures.}

\vspace {1mm} 

{\it Proof of {\color {blue} \Large  $\square$}{\small $\Circle $}}. Let $P\to X$ be the principal fibration associated with  
the bundle $T(X)/T(\mathscr L)$ and by blowing up  the metric of $P$ transversally to the lift $\tilde{\mathscr L}$ to $P$ as in the previous section, make the scalar curvature of $P$ on a given compact domain  $P_\varepsilon\subset P$  greater than $n(n-1)-\varepsilon$ for a given $\varepsilon>0$.

Also with  this blow-up,    make the Lipschitz  constant of the map $P\to X$ as small as you want.

(A possibility of this formally follows from the above  {\color {blue} (1)} -  {\color {blue} (6)} for foliations of codimension one, while  the proof in  general case amounts to    replaying   {\color {blue} (1)} -  {\color {blue} (6)} word-for-word in the present case.)

Next, let   $G=\mathbb R^m$, observe as in 
 (2) in section \ref{Euclidean fibrations6}  that   $P_\varepsilon$ admits a $(1+\varepsilon)$-Lipschitz map  of degree one from   $P_0$ to $X\times [0,L]^k$ for an  arbitrary large $L$  and apply the maximality/extremality theorem for punctured spheres from 
sections \ref{punctured3} and \ref{log-concave5}. 

This concludes the proof for $G=\mathbb R^m$  and the case of the  general Abelian $G$ follows by passing to the quotient of $G$ by the maximal compact subgroup.

\vspace {1mm}

To get an idea  why one can control the geometry of the  blow-up  only on compact subsets in P,
 look at the following.

\vspace {1mm} 

{\it \color {blue} Geometric Example.} Let $(Y,g)$ be  a Riemannian manifold and let $P_Y\to Y$ be the fibration, with the fibers $S_y$, $y\in Y$,  equal to the spaces of quadratic forms in the tangent spaces $T_y(Y)$ of the form $c\cdot g_y$, $c>0$.
Thus, $P_Y=Y\times \mathbb R$, for $\mathbb R  =\log \mathbb  R^\times_+$ with the metric  $e^{2r}dy^2+dr^2$. 

When  $r\to +\infty$  and the curvature of  $e^{2r}g$ tends to zero, then the   metric $e^{2r}dy^2+dr^2$ converges to the hyperbolic one with constant curvature $-1$, but  when  $r\to -\infty$, then the curvatures of  $e^{2r}g$ and of $e^{2r}dy^2+dr^2$
blow up at all  points $y\in Y$, where the curvature of $g$ doesn't vanish. 

And 
 if apply this to the fibration $P=P_Y\times \mathcal L \to X=Y\times \mathcal L$ with the same $\mathbb R$-fibers,  then  we see that the convergence of the scalar curvatures of the  blown-up $P$ to those of   $\mathcal L$ is by no means uniform.

%%%%%%%%%%%%%%%%%%%%

\subsubsection {\color {blue} Hermitian Connes' Fibration}\label {Hermitian Connes6}
%%%%%%%%%%%%%%%%%%%%%%%
  Let  
  $\mathscr L$ be a foliation on $X$  of codimension $k$ as earlier with 
  a  transversal (sub)bundle $T^\perp\subset T(X)$ and a Bott connection in it. 
  Let  $T^{\bowtie}$  be the sum of  $T^\perp$ with its dual bundle  and endow   $T^\bowtie$ with the natural, hence monodromy invariant,  symplectic structure.
  
   Let $S_x$ denote the space of Hermitian structures in the space $T_x^\bowtie$, for all $x\in X$, and let 
   $P\to X$ be the corresponding fibration, that is the  fibration associated with  $T_x^\bowtie$ with the fiber
   $S=Sp(2k,\mathbb R)/U(k)$.
   
    Equivalently, this fibration  $P\to X$  is associated to $T^\perp\to X$ via the action of $GL(k) $ on $S=Sp(2k,\mathbb R)/U(k)$ for the natural embedding of the linear group $GL(k)$ to the symplectic $Sp(2k,\mathbb R)$
   
  Besides  sharing  the  properties {\color {blue}  (1)}-{\color {blue}  (6)} of the original Connes' bundle formulated  in section \ref{}, this new $P\to X$ has,  
   a lovely    additional  feature:  >??
   
   {$S$ is  a {\it Hermitian} (irreducible) symmetric space,  which implies (see section  \ref{spin harmonic6})
 non-vanishing of the index  
   of some  twisted Dirac  on $S$ that is invariant under the isometry group (that is $Sp(2k,\mathbb R)$ of $S$. 
   
      This, as it was  {\sf \color {blue} }  stated in section  \ref{spin harmonic6} must imply the existence of twisted harmonic $L_2$-spinors on fibrations with $S$-fibers, which we formulate below in the form relevant to foliations positive   scalar 
   curvatures and which, besides being   interesting in its own right, would simplify the proof by Connes in [Connes(cyclic cohomology-foliation) 1986] as well as the arguments from [Zhang(foliations)  2016].

    {\sf L et $Y=(Y, \omega)$ be a closed symplectic manifold of dimension $2k$  and let $F: P_Y\to Y$  be the fibration associated with the tangent bundle $T(Y)$  with the fiber $S=Sp(2k, \mathbb R)/U(k)$.

Observe that the quotient bundle $T(P)/T^{vert}$ carries a tautological Hermitian metric $g_\bowtie$, and   a granted  $Sp(2k, \mathbb R)$-connection
   in the tangent bundle  $T(Y)$, that is a horizontal subbundle $T^{hor}\subset T(P)$, one obtains a Riemannian metric $g_P$
   in the tangent bundle $T(P)=T^{vert}\oplus T^{hor}$ that is $$g_P=g_S+g_\bowtie$$
where $g_S$ is a $Sp(2k, \mathbb R)$-invariant Hermitian  metric in $S$, which is unique  up to scaling.}
   
     \vspace{1mm}
 {\sf Let the symplectic  form $\omega$ be integer and thus  serves as the curvature of a unitary line bundle $L\to Y$.  }
   
     \vspace{1mm}
 
   {\sf \color {blue} \large {\color {red!50!black} Conjecture} 1}  {\sl The  bundle of spinors  on $P$ twisted with   some tensorial power of  the bundle  $F^\ast(L)\to P$
  admits a non-zero harmonic $L_2$-section on $P$. }
  
   \vspace{1mm}

  {\it Remarks and Examples.} (a) The geometry of this  $P$, unlike of what we  met in section \ref{spin harmonic6}, is as far from being a product as  in $P_Y$ from the { \color {blue} geometric example} in section  \ref{Abelian Monodromies6}.

(b) The simplest instance of $Y$ is that of an even dimensional  torus $\mathbb T^{2k}$ with an invariant symplectic form $\omega$  and trivial flat symplectic connection.

In this case, the universal covering $\tilde P_Y$ of the  manifold $P_Y$ is Riemannian homogeneous;  moreover, the (local) index  integrant is homogeneous as well. It is probable, that a version of  the Connes-Moscovici theorem applies in this case and yields twisted  harmonic $L_2$-spinors on $\tilde P_Y$ and, eventually,  on $ P_Y$.

(c) It would be most amusing to find a link between  the symplectic geometry of $(Y, \omega).$ and   twisted Dirac operators on $P_Y$ or their non-linear modifications.

\vspace{1mm}
  
 Let us modify the above conjecture  to make it applicable to foliations.

 Let $S=G/K$  be a symmetric space, where the index of the Dirac  twisted with some bundle $L_S\to S$ associated with the $K$-bundle $G\to S$ doesn't vanish, e.g.  $S=Sp (2k)/U(k)$.

Let $F:P\to X$ be a smooth $S$-fibration with a smooth  $G$-connection $\underline \nabla$ and let $T^{hor}\subset T(P)$ be the corresponding horizontal subbundle.

Let $L_{\updownarrow}\to P$ be the  bundle   the restriction of which to the fibers $S=S_x \subset P$, $x\in X$ are equal to $L_S$.

Let $g^{hor}$ be a  smooth Riemannian metrics (positive quadratic forms) in $T_{hor}$ and $g_P=s^{hor}+ g_S$ be the sum of this metric with a $G$-invariant metric in the fiber.

    Let $L_X\to X$ be a vector bundle with a unitary connection $\nabla_X$ trivialized at infinity (which is relevant for  non-compact manifolds $X$) and let $L^\ast =(L^\ast, \nabla^\ast)\to P$ be the  bundle pulled back by $F$ from $L_X$ along with the connection $\nabla_X$, that is $L^\ast=F^\ast( L_X, \nabla_X)$. \vspace {1mm}

   {\sf \color {blue} \large {\color {red!50!black} Conjecture} 2.} {\sf If $X$ is complete and if some Chern number of $L_X$ doesn't vanish,
   then $P$ supports a non-zero harmonic $L_2$-spinor $s =s(p)$ twisted with $L_{\updownarrow}$ and   with (i.e. (tensored with) some  bundle associated with $L^\ast$.
   
   Moreover, there exists such a non-zero spinor $s(p)$,  the rate of  the  decay   of which  at infinity is  independent of the metric $g^{hor}$}:
   
  { \sl given an exhaustion of $P$ by compact domains,  $P_1\subset ... \subset P_i\subset ... P$, then
   $$\frac {\int_P||s(p)||^2dp} {\int_{P_i}||s(p)||^2dp} \geq 1-\varepsilon(i)\underset {i\to \infty}\to 1,
    \leqno {\mbox {\color {blue} {\large (}$L_2$\hspace {-0.8mm}\large{$\searrow)$}}} $$
   where the function  $\varepsilon(i)>0$ may  depends on $P_i$, $X$, $F$, $\underline \nabla$, $ \nabla_X$ and   $g_S$, but       which is independent of the metric  $g^{hor}$. }   
   \vspace {1mm}
  
\subsubsection{ \color {blue} Hermitian Connes' Fibrations over Foliations  with Positive Scalar Curvature}\label {over  positive scalar6}

   \vspace {1mm}
  
  Let    $F:P\to X$ be the Hermitian Connes' fibration with the $S$-fibers, $S=Sp(2k)/U(k)$, 
  over a Riemannian manifold $X=(X,g)$ with a foliation
   $\mathscr L$  of codimension $k$ on it , as in the previous section, let  
    the subbundle  $T^{hor}\subset T(P)$ corresponds to a Bott connection on the   bundle $T^\perp\to X$ normal to the tangent subbundle $T(\mathscr L)\subset T(X)$ and let us lift the splitting 
   $T(X)= T(\mathscr L)\oplus T^\perp$ lifts to the corresponding splitting 
  $T^{hor}= \tilde T(\mathscr L)\oplus \tilde T^\perp$.

 Recall that the points   $p\in P$ correspond to Hermitian structures in the symplectic   spaces $T_{F(p)}^\bowtie \supset T^\perp_{F(p)}$, the real parts of which give  Riemannian/Euclidean structures to $T^\perp_{F(p)}$  and which  then pass to the spaces 
  $T _p^{hor}$ via the  the differentials $dF_p :T_p(P)\to T_{F(p)}(X)$ which  are  isomorphic on the fibers of $T_p^{hor}\subset  T_p(P)$.  
 
 Thus, the bundle $\tilde T^{\perp}\to P$ carries a canonical Riemannian metric, which we  call $\tilde g^\perp$.
 
Next, let $\tilde g_{\mathscr L}$ be the metric on  the bundle 
  $\tilde T(\mathscr L)$ 
 that is induced from the Riemannian metric $g_{\mathscr L}$ that is the metric $g$ 
 on $X$ restricted to the bundle $T(\mathscr L)$
 and  let us endow the bundle $T(P)=\tilde T(\mathscr L)\oplus \tilde T^\perp \oplus T^{vert}$, where   $T^{vert}\subset T(P)$ is the bundle tangent to the $S$-fibers of the fibration $F:P\to Q$, with the metrics  
 $$\tilde g_{e_S,e^\perp}=\tilde g_{\mathscr L}+e_\perp\cdot  \tilde g^\perp+ e_S\cdot g_S, \mbox { }  e_S,e_\perp>0, $$ }
as this was done in {\color {blue} (6)} of section \ref{Connes6}, except that now the fibers $S$ are isometric to  $Sp(2k)/U(k)$ with a (unique up to scaling) $Sp(2k)$-invariant metric $g_S$, rather than  to $GL(k)/O(k)$ as in section  \ref{Connes6}.

\vspace {1mm}
 
 Now as in {\color {blue} [$\Uparrow_{Sc}$]} of section  \ref{Connes6}, we observe  the following. \vspace {1mm}

{\it \color {blue} Effect of  $g_S $-Blow-up }. If 
the constant $e_S$ is  much greater than $\frac {1}{\sigma} $, then 
 the scalar curvatures of the leaves  
${\mathcal L}^+ \subset P$, which are     the pullbacks of the leaves
$\mathcal L\subset X$  under the map $P\to X$ (see   {\color {blue} (2)} in section \ref {Connes6}),
become close to those of the underlying leaves  $\mathcal L$, hence {\color {blue}$\geq \sigma-\varepsilon>0$}, while the norms of the  logarithmic covariant  derivatives $\nabla_{\cal L^+}(\log \tilde g^\perp)$ 
 of the transversal metric  $\tilde g^\perp$ along the leaves     $\mathcal L^+$ with respect to  the metrics 
 $ \tilde g_{\mathscr L}+ e_S\cdot  g_S,$   becomes {\color {blue} $\leq \varepsilon$} on the leaves
   $\mathcal L^+\subset P$, where observe, that, given an $\varepsilon>0$,  \vspace {1mm}
   
   {\it  if  $e_S = e_S(\varepsilon) $ is sufficiently large, then   these two $\varepsilon$-bounds hold  on {\color {red!60!black} all}  of $P$.}\vspace {1mm}

 {\it \color {blue}  $ \tilde g_S^\perp$-Blow-up}. This has two   effects  on the geometry of  $P$.\vspace {1mm}

{\color {blue} $\mathbf 1_{Sc} $}  {\sf   If the scalar curvatures of the leaves $\mathcal L^+$  are  bounded from below by $\sigma-\varepsilon$ and if 
 the norm of  $\nabla_{\cal L^+}(\log \tilde g^\perp)$ is bounded by $\varepsilon$, then 
 $$Sc(P)\underset {e_\perp \to \infty}  \to \sigma -\epsilon$$
where   
 $$ \epsilon\underset {\varepsilon \to 0 }\to 0$$
 and where the convergence $Sc(P)  \to \sigma -\epsilon$ is uniform on compact subsets in $P$} (but not on all of $P$). 
 \vspace {1mm}

{\color {blue} $\mathbf 2_{Lip} $} {\sf If $e_\perp\to \infty$, then (regardless  of $e_s$)  the Lipschitz constant of the map $F:P\to X$ tends to zero  uniformly on  compact subsets in $P$.}

  \vspace {1mm}

{\it \large \color {blue} Conclusion.}    Granted   {\sf \color {blue} \large {\color {red!50!black} Conjecture} 2} from the previous section,  we see that, as far as the the Dirac operators  are concerned,  the   positivity of the scalar curvature of  $g|\mathscr L$ has the same effect   as of the metric $g$ on $X$ itself. 

For example,      {\sf \color {blue} \large {\color {red!50!black} Conjecture} 2}  implies the following.\footnote{Here we use the fact that
if $X$ is spin then $P$ is also spin, where, if you are in doubt,  this implication  can be achieved by taking $S\times S$ instead of $S$.}

 \vspace {1mm}

 {\color {blue}\small $\bigstar$  }  {\sl  Let $X=(X,g) $ be a complete orientable spin Riemannian $n$-manifold and let   $\mathscr L$, be a  smooth foliation 
on $X$of codimension $k$, such that the scalar curvature of 
 $g$    restricted to the leaves of  $\mathscr L$ satisfies 
$$Sc(g_{| \mathscr L})> n(n-1).$$ 
Then every $1$-Lipschitz map  $X\to S^n$ locally constant at infinity has zero degree.} \vspace {1mm}

Let us spell it  all out again.
  
   Let $\mathscr D$ be a Dirac operator on $P$,  which is

(a) twisted with an $S$-fiber bundle that makes the local  index of the corresponding  Dirac operator on $S$ non-zero;

(b) on the top of that the  $\mathscr D$  is twisted with a bundle induced from a bundle $L_X$ on $X$,  with a non-zero Chern number, where one may assume that the integrant in the local formula for the index of $\mathscr D$ doesn't vanish on $P$.

\vspace{1mm}

The above shows, that   if $Sc(\mathscr L)\geq \sigma>0$  and the curvature of $L_X$ is small, then,  given a compact  domain $P_0\subset P$, 
 the spectrum of the  $\mathscr D_{e_S,e_\perp}^2$  on $P_0$ (that is   $\mathscr D^2$ for the metric 
  $\tilde g_{e_S,e_\perp}$ on $P_0$) with the zero boundary condition can be made uniformly separated away from zero,
 that is  $\lambda_1\geq \delta=\delta(n, \sigma)>0$
  by taking   sufficiently large $e_S$ and $e_\perp$.\vspace{1mm}

But this contradicts  {\sf \color {blue} \large {\color {red!50!black} Conjecture} 2}, which  implies that \vspace{1mm}

 {\sl given $e_S\geq 1$ and $\delta>0$, there exists a compact domain $P_0\subset P$, such that the first eigenvalue of $\mathscr D_{e_S,e_\perp}^2$ on $P_0$ satisfies
$$\lambda_1=\lambda_1(P_0, \mathscr D_{e_S,e_\perp}^2)\leq \delta.$$}

To prove such a bound on $\lambda_1$,   one  needs to construct an almost , $\mathscr D$-harmonic spinor with support in $P_0$, where  a natural pathway   to  this end  goes along  the lines of the local proof of the index theorem,  roughly, as follows.

Let $\cal K$ be an  that is a function of the (unbounded self-adjoint)  $\mathscr D$ (we suppress the subindices $e_S$ and $e_\perp$),  which can be represented by a smooth kernel  $\mathcal K(p,q)$, $p,q\in P$ supported in a $d$-neighbourhood of the diagonal, where
this the values of  $K(p,q)$ at all pints $(p,q)\in P\times P$,  depends on the metric in $P$  only  in the ball of  $B_{p,q}(4d) \subset P\times P$ and such that
the super-trace of $\mathcal K$ and of all its powers $\mathcal K^i$ is equal to the index of $\mathscr D$, whenever this construction is applied to compact manifolds $P$. 

Then, as $i\to \infty$, this  converges to the projection to the kernel of $\mathscr D$, that is the space of harmonic 
$L_2$-spinors,  and  the only issue to settle in the present case is certain uniformity of this convergence, for $e_\perp\to\infty$.

What may facilitate the estimates needed for the  proof of this is the uniform bound on  geometries  of  the metrics $\tilde g_{e_S,e_\perp}$ for $e_S,e_\perp\to infty$, probably, where, possibly,  one can get a fair  representation/approximation  of  $\mathcal K^i_{e_s,e_\perp}$ by a (singular) perturbation   argument at $e_\perp=\infty$. 

\vspace{1mm}

{\it Remark.} Even if the above argument is carried thought it, as we mentioned in section  \ref{foliations3},  it will be not deliver 
what, probably,  follows by Connes' argument:

  "{\sl no $1$-Lipschitz map  $f: X\to S^n$ with $deg(f)\neq 0$}" for $Sc(g_{| \mathscr L})> (n-k)(n-k-1). $ \footnote {According to what  was explained to me by Jean-Michel Bismut, the same  may apply to Zhang's argument.} 

But it seems beyond the present day methods to drop  the spin assumption in    {\color {blue}\small $\bigstar$}
for $k\geq 2$.
\vspace{1mm}

{\it On Non-Integrable Generalization.} Let  $X=(X,g)$ be  a Riemannian manifold, let   $\Theta\subset T(X)$ be a smooth subbundle of codimension $k$ and let 
$$\Lambda=\Lambda(\Theta): \wedge^2\Theta\to \Theta^\perp $$
 be its    curvature, that is the 2-form on $\Theta $ with values in the normal subbundle  (identified with the quotient bundle $T(X)/\Theta)$, which is defined by the normal components  of  commutators of pairs of  tangent fields $X\to\Theta$.
x $Sc(g{|\Theta}, x)$, $x\in X$,  be the sum of the sectional curvatures over the pairs of   vectors in an orthonormal  basis in $\Theta_x$.

It seems probable, that most (all?) we know and/or conjecture about (tangent subbundles of) foliations with $Sc\geq \sigma>0$  \vspace {1mm}

{\sf extends to  $\Theta$ with $Sc(g{|\Theta})\geq \sigma >0$, if $||\Lambda(\Theta)||$ is much smaller than $\sigma$.} \vspace {1mm}

(Homogeneous $\Theta$, e.g. on  spheres of dimensions $2m+1$ and $4m+1$  may serve as extremal cases of the corresponding inequalities.) 

\subsubsection {\color {blue} Geometry and  Dynamics of Foliations with Positive  Scalar Curvatures} \label {dynamics6}

 Let us formulate a few versions of  the width/waist conjecture  
from section \ref{slicing3D.3} for foliations with  $Sc >0$.\vspace {1mm}

{\sf  Let $X$ be a complete Riemannian $n$-manifold with a foliation $\mathscr L$ of codimension $k$,  where $n-k\geq 2$ and where the scalar curvatures of  the induced Riemannian metrics in the leaves  satisfy $Sc\geq (n-k)(n-k-1)$.}  \vspace {1mm}

 {\huge $\star \star$} {\large \color{blue} (Unrealistically?)  \it Strong Foliated  Width/Waist {\color {red!50!black} Conjecture}.}  {\sf There exists a continuous map from $X$ to an $(n-2)$-dimensional polyhedral space,  say $F:X\to P^{n-2}$, such that }

{\sl the pullback $F^{-1}(p)\subset X$ is contained in a single leaf of $\mathscr L$ for all $p\in P$  and 
$$ diam (F^{-1}(p))\leq const_n  \mbox    { and } vol_{n-2}(F^{-1}(p))\leq const'_n \mbox { for all }p\in P^{n-2},$$
where, conceivably,  these constants  don't even depend on $n$, e.g.  with a  possibility   $const'=4\pi$.}
 \vspace {1mm}

{\it An Impossible Proof of {\huge $\star \star$.} }Ideally, one would like to have a continuous family of metrics in the leaves that would eventually 
simultaneously collapse all leaves this factorizing $X$ to something $n-2$-dimensional.

But the obvious candidate for this -- Hamilton's Ricci flow, even if it is defined for all time, doesn't collapse $X$ fast enough to bound $width_{n-2}$ or $waist_{n-2}$. 

In fact, what happens is better    seen for the mean curvature flow, where the collapsing map for ellipsoids with the principal axes of the lengths $1, 1, d$, moves this ellipsoid by distance $ \sim d$, which may be arbitrary large.

\vspace {1mm}

Here is another way to look at this problem.
\vspace {1mm}

(Provisional)   {\it  Bounded Distance Deformation Question. } Let $(X,g_1)$ be  a complete Riemannian manifold with $Sc\geq 1\sigma>0$.

Does there exist a Riemannian  metric $g_2$ on $X$, such that $Sc(g_2)\geq 2$  and  
$$|dist _{g_1}(x,y)-dist _{g_2}(x,y)| \leq const_n$$
  for all $x,y\in X$ and some $const_n<\infty$?

\vspace {1mm}

{\it \color {blue} 3D-Foliations.} Let $X=(X,g)$ be a complete  Riemannian manifold and   $\mathscr L$ be a smooth 3-dimensional  foliation such that the restrictions of $g$ to all leaves $\mathcal L$  of   
{\it \color {blue} 3D-Foliations.} Let $X=(X,g)$ be a complete  Riemannian manifold and   $\mathscr L$ have $Sc_g(\mathcal L)\geq \sigma >0$.
 
 Then there exist continuous maps from $\mathcal L$  to  locally finite 1-dimensional  polyhedra,  say  $F:\mathcal L\to P^1$,
 such that $diam (F^{-1}(p))\leq  c<\infty$ for all leaves $\mathcal L$ of $\mathscr L$.
 
(Such maps exist with $c=2\pi\sqrt{\frac  {36}{\sigma}}$, see section \ref{slicing3D.3}, but what is relevant at the moment is the universal  bound 
$c=c(\sigma)<\infty$.)
 
Moreover, since such maps also exist for the {\it universal  coverings of the leaves}, the maps $F$ can be coherently chosen on all leaves simultaneously by adapting 
{\sf \color {blue!50!black} the geometric proof of  \it the  Stallings Ends of the Groups Theorem} to foliations (compare with  \S$ 2\frac {2}{3}$ of  [G(positive)   1996]).

Let us recall this proof in the  simplest case where   $X$  is a smooth  manifold with  a discrete cocompact action of a group $\Gamma$ \footnote{Contrary to the statements found  sometimes  in the literature, all versions of Stallings' theorem, as well as of its refinements and generalizations, effortlessly   follow with a  {\it proper} use  of  minimal hypersurfaces.} and show that  \vspace {1mm}

{\it Proof. } (Compare [G(infinite)  1984])  Let $Y\subset X$ be a connected {\it volume minimizing}  compact hypersurface 
which  {\it separates two ends} of $X$. Then, because of minimality,  no  $\gamma$-translate   $\gamma(Y)\subset X$,  intersects $Y$ unless  $\gamma(Y)=Y$.

Then the proof easily follows, since, due to  the    the end separation property,  the complement to the set  of the translates of $Y$, that is  
$$X\setminus \bigcup_{\gamma\in \Gamma}\gamma(Y),$$
is the union of {\sl  mutually non-intersecting subsets of diameters $\leq const$.}\vspace {1mm}

Similar argument applies to foliations where one similarly achieves  invariance of the  relevant maps under the monodromy 
groupoid instead of $\Gamma$ (see  \S$ 2\frac {2}{3}$ of  [G(positive)   1996])
\footnote{ The argument from [G(positive)   1996]  becomes more transparent,  if one  makes the metric  $g$ on $X$ {\it generic} by a small perturbation, for which  all compact  locally minimizing hypersurfaces $Y$ in the 
leaves are isolated;  hence stable under small transversal deformations of leaves.

  This allows a sufficient quantity of   compact  volume minimizing  hypersurfaces $\tilde Y$ with $diam(\tilde Y)\leq const$   in the leaves  
$\tilde {\mathcal L} $  of the lift   $\tilde {\mathscr L}$  of  ${\mathscr L}$ to the universal covering $\tilde X$, such that
the  
intersections of the leaves   $\tilde {\mathcal L}$ with  the complement of the union of these  $\tilde Y$ have the diameters of all their connected components   also uniformly bounded, say by  $\leq const'$.}

which implies, for instance the following.\vspace {1mm}

{\Large \color {blue} ($\star_3$)}   {\it  If a compact orientable Riemannian $n$-manifold $X=(X,g)$ carries a 3-dimensional foliation, where the leaves have positive scalar curvatures,
 $$Sc_g(\mathcal L)>0,$$
 then $X$ admits no maps of non-zero degrees to aspherical $n$-manifolds.}\vspace {1mm}

Conclude with the following purely foliation theoretic question, the positive answer to which, that, I think, is unlikely,  motivated the above conjectures.

{\sf Is it true that  {\it no smooth  foliation} on  $\mathbb R^n$ of positive dimension invariant under the action of $\mathbb Z^n$\footnote{Ideally, one  would like to drop this invariance condition.} can have   the diameters of all leaves bounded by a common constant $C<\infty$?}

(An approach to  counterexamples may be found in   [EM(wrinkling) 1998].)

%%%%%%%%%%%%%%%%%%%%
\subsection {\color {blue}Moduli Spaces Everywhere}\label {moduli spaces6}
%%%%%%%%%%%%%%%%%%%%%

{\sf \large  \color {red!50!blue} {\it \large All }  topological and geometric constraints on metrics with $Sc\geq \sigma$  are accompanied by non-trivial homotopy theoretic properties of spaces of  such metrics.} 
\vspace{1mm}

A manifestation of this  principle    is  seen in how    topological obstructions for the existence of metrics with $Sc>0$ on closed manifolds $X$ of dimension $n\geq 5$ give rise to 

{\sf pairs $(h_0, h_1)$ of    metrics with $Sc\geq \sigma >0$ on closed hypersurfaces $Y\subset X$ which 

{\sl can't be joined by homotopies $h_t$ with $Sc(h_t)>0$}.}

\vspace{1mm}

The elementary argument used for the proof of this (see section \ref{spaces of metics3}) also shows that    (known)  constraints on {\it geometry}, not only on topology,  of manifolds with  $Sc\geq  \sigma$ play  a  similar role.  

For instance,   assuming for notational simplicity, $\sigma=n(n-1)$, and recalling the {\color {blue} $ \frac {2\pi}{n}$-inequality} from section \ref{bands3}, we see that \vspace{1mm}

(a)   {\sl if  $l\geq \frac {2\pi}{n}$, then   the pairs of metrics $h_0 \oplus dt^2$  and  $ h_1 \oplus dt^2$ on the cylinder  $Y\times [-l,l]$,  for the above $Y$ and $l\geq \frac {2\pi}{n}$, {\it can't be joined} by homotopies of metrics $h_t$  with $Sc(h_t)\geq n(n-1)$ and with
$dist_{h_t}(Y\times \{-l\}, Y\times \{l\}) \geq \frac {2\pi}{n}$.}\vspace{1mm}

This  phenomenon is also observed   for manifolds with {\sl controlled  mean curvatures of their boundaries}, e.g.  for Riemannian bands $X$ 
with $mean.curv (\partial_\mp X)\geq \mu_\mp$ and with $Sc(X)\geq \sigma$, whenever these inequalities imply that $dist(\partial_- X, \partial_+X )\leq d=d(n, \sigma ,\mu_\mp)$. 
 (One may  have $\sigma<0$ here  in some cases.)

Namely,\vspace{1mm}

(b)  {\sl  certain  sub-bands  $Y\subset X$ of codimension one with $\partial_\mp(Y)\subset \partial_\mp(X)$   admit pairs of metrics $(h_0,h_1)$, such that  $mean.curv_{h_0,h_1} (\partial_\mp Y)\geq \mu_\mp$ and  $Sc_{h_0,h_1}(Y)\geq \sigma$ while  
  $dist_{h_0,h_1}(\partial_-,\partial_+)\geq  D$  for a given $D\geq d$.
  But these metrics   can't be joined  by homotopies $h_t$ , which would keep these inequalities  on the scalar and on the mean curvatures and  have $dist_{h_t}(\partial_-,\partial_+)\geq  d$  for all $t\in [0,1]$.}\vspace{1mm}

(c) This  seems to persist (I haven't carefully  checked it)  for manifolds with corners, e.g. for cube-shaped manifolds $X$: these, apparently contain hypersurfaces $Y\subset X $, the boundaries of which $\partial Y\subset \partial X$ inherit the corner structure from that in $X$, and which admit pairs of  "large" metrics $h_0, h_1$, which also have "large"  scalar curvatures, "large" mean curvatures of the codimension one faces $F_i$ in $Y$   and  "large"  complementary ($\pi-\angle_{ij}$)  dihedral angles along the codimension two faces $F_{ij}$, but where these $h_0, h_1$ can't be joint by homotopies of metrics $h_t$ with comparable  "largeness" properties.\vspace {1mm}

It is  unclear, in general,  how to  extend the $\pi_0$-non-triviality (disconnectedness) of our spaces of metrics to the higher homotopy groups, since the techniques currently used for this purpose   rely entirely on the Dirac theoretic techniques (see 
 [Ebert-Williams(infinite loop spaces) 2017] and references therein), which are poorly adapted to manifolds with boundaries.  But some of this is possible for closed manifolds. 

For instance, let $Y$ be a smooth closed spin manifold, and 
$h_p$,  $p\in P$,  be a homotopically non-trivial family of metrics with $Sc(h_p)\geq \sigma> 0$, where, for instance, $P$ can be a $k$-dimensional sphere  and  non-triviality means non-contractibility.

Let ${\cal S}^m_\sigma(S^m\times Y)$ denote the space  of pairs $(g, f)$, where $g$ is a Riemannian  metric on  $ S^m\times Y $ with $Sc(g)\geq \sigma$ and $f:(S^m\times Y,g)\to S^m$ is a distance  decreasing map   homotopic to the projection $f_o:S^m\times Y \to S^m$. 

 \vspace {1mm}
 
{\sf If    non-contractibility of the 
 family $h_p$  follows from non-vanishing of the index of some Dirac operator, then (the proof of)  Llarull's theorem suggests that the corresponding family $(h_p+ds^2, f_o)\in {\cal S}^m_{\sigma_+}(S^m\times Y)$ for $\sigma_+=\sigma+m(m-1)$ is non-contractible in the space 
$$ {\cal S}^m_{m(m-1)}(S^m\times Y) \supset {\cal S}^m_{\sigma_+}(S^m\times Y). $$}

 This is quite transparent in many cases,  e.g. if $h_p = \{h_0, h_1\}$ is an above kind of   pair of  metrics  with $Sc>0$, say an embedded codimension one  sphere 
 in a Hitchin's  homotopy sphere. 
\vspace {1mm}

{ \it Remarks.}  (i)  If  {\it "distance  decreasing"} of $f$  is strengthened to {\it "$\varepsilon_n$-Lipschitz"}   for a sufficiently small $\varepsilon_n >0$, then the above  disconnectedness of the space of pairs  $(g,f)$
follows for all $X$ with a use  of minimal hypersurfaces instead of Dirac operators.

\vspace {1mm}

(ii) The above definition of the space ${\cal S}_\sigma^m$ makes sense for all manifolds $X$ instead of $S^m\times Y$, where one may allow   $dim(X)< m$ as well as $>m$.

However,  the following   remains   problematic in most cases.\vspace{1mm}

{\sl For which closed  manifolds $X$ and numbers  $m$, $\sigma_1$ and 
$\sigma_2>\sigma_1>0$ is the inclusion ${\cal S}^m_{\sigma_2}(X) \leq {\cal S}^m_{\sigma_1}(X)$ homotopy equivalence?}

\vspace{2mm}

{\it \large Suggestion to the Reader.} Browse through all theorems/inequalities  in the previous as well as in the following  sections, formulate their possible homotopy parametric versions and try to prove some of them.

%%%%%%%%%%%%%%%%%%%
\subsection {\color {blue}Corners, Categories and Classifying Spaces.}\label {corners categories6}
%%%%%%%%%%%%%%%%%%%%%%%%%%

It seems (I may be mistaken)  that all known results concerning  the homotopies  of spaces with metrics $Sc>0$  are about  {\it iterated (co)bordisms}  of manifolds  with $Sc>0$ and/or about  {\it cobordism categories with $Sc>0$} in the spirit of  [EbR-W(cobordism categories) 2019], rather than about spaces of metrics per se.\footnote {See  [Kaz(4-manifolds) 2019] for a computation of such cobordisms in dimension $3$.}

To explain this, start with  thinking of  morphisms $a \to b$  in  a  category     as members of  class of {\it labeled} (directed) {\it edges/arrows } [0,1]  with the   0-ends labeled by  $a$ and the 1-end labeled  by $b$.  

Then define  a {\it cubical category}  $\cal C$ (I guess there  is a standard  term but I don't know it)
as a class of {\it labeled combinatorial cubes} of all dimensions, $[0,1]^i$, $i=1,2,...$, where  all faces are labeled  by members  of some class  and which satisfied  the obvious generalisations of the axioms of the ordinary  categories: associativity and the presence of the identity morphisms.

{\it Example.}   Let    ${\cal C}=A^\square$ consist of 
continuous maps from  cubes to a topological space $A$, e.g. to the space $A=G_+=G_+(X)$ of metrics with positive scalar curvature on a given manifold $X$, where these maps are regarded as labels on the cubes they apply to.
 
 If we glue all such cubes along faces with equal labels, we obtain a cubical complex, call it $|{\cal C}|$,  which is  (weekly) 
 homotopy equivalent to $A$, where possible degeneration of cubes.e.g. gluing two faces of the same cube, is offset by possibility of unlimited subdivision of   cubes by means of  cubical identity morphisms.

Next,
  given a smooth closed manifold $X$,  consider "all" Riemannian manifolds  of the form  $(X\times [0,1]^i, g)$, $i=0,1,2,...$, such that $Sc(g)>0$, and such that the  metrics $g$    in  small neighbourhoods of all  "$X$-faces" $X\times F_j$, where $F_j$ is are  ($(i-1)$-cubical)  codimension one faces in the cube  $[0,1]^i)$, split as  Riemannian products: $g=g_{X\times F_j}\otimes dt^2$.
  Denote the resulting cubical category by $XG_+^\square$ and observe that there is a natural cubical map 
   $$\Xi:| G_+(X)^\square|\hookrightarrow   |XG_+^\square|.$$

Now we can express the above  "iterated cobordism" statement  by saying that the only part of the homotopy invariants of $G_+(X)$  (which is homotopy equivalent to $G_+(X)$), e.g of its homotopy groups, which is  detectable by the present  methods is what remains non-zero  in   $  |XG_+^\square|$ under $\Xi$.

Similarly one can  enlarge  other spaces of Riemannian  metrics on non-closed manifolds  from the previous section with lower bounds on their curvatures and their sizes, where the latter can be  expressed with  maps $f: (X,g)\to \underline X$, with controlled Lipschitz constants with respect to $g$, or with respect to the  {\it {\sf Sc}-normalised  metric   
$Sc(X)\cdot g$}.  

\vspace {1mm}

There is yet another way of enlarging the cubical category $XG_+^\square$, namely by 
$B^{\ast} G^\square ( {\sf D})$, where  {\sf D} is topological, e.g. metric space  and where  

$\bullet _0$ closed oriented Riemannian  manifolds $X$ of all dimensions $n$ along with continuous maps $X\to   {\sf D}$ stand for 0-cubes -  "vertices",

$\bullet _1$   "edges" ; i.e 1-cubes are cobordisms $W^{n+1}$  between $X_0, X_1$,  with Riemannian metrics split near their boundaries $\partial W^{n+1}= X_0\sqcup - X_1$,  and continuous maps to  {\sf D} extending those from $X_0 and X_1$,

$\bullet _2$   "squares",  are  (rectangularly cornered $(n+2)$-dimensional)  cobordisms between $W$-cobordisms  with maps to {\sf D}, etc.\vspace {1mm}

The actual  cubical subcategory of $B^{\ast} G^\square ( {\sf D})$, which is relevant for the study of the space $|XG_+^\square|$ (that is, essentially,  the space of metrics with $Sc>0$ on $X$) is where 
 all manifolds in the picture are spin, the scalar curvatures of their metrics are positive, 
 $D$ is the classifying space of a group $\Pi$  and where one may assume the fundamental groups of all $X$ to be coherently (with inclusion homomorphisms) to be isomorphic $\Pi$ \footnote {This "assume"  relies on the codimension two surgery of manifolds with $Sc>0$, which is possible for making  the fundamental groups of $n$-manifolds isomorphic to $\Pi$   if   $n\geq 4$ and  where 
 more serious  topological conclusions need  $n\geq 5$.}   (compare   [Ebert-Williams(infinite loop spaces) 2017],  [BoEW(infinite loop spaces) 2014],  [HaSchSt(space of metrics)2014] , 
\vspace{1mm}

{\large \it \color {blue} Question.}  {\sf What are possible generalizations of the above to manifolds with corners,  which are far from being either cubical or rectangular?}

For instance, prior to speaking of spaces of metrics and of  categories of cobordisms,  let $X$ be an individual manifold with corners, say a (smoothly) topological $n$-simplex or a dodecahedron, let   $(\infty<\sigma<\infty)$,  let $(\infty<\mu_i<\infty)$  be numbers assigned to  the codimension one  faces $F_i$ of $X$ and 
$0<\beta_ij< \pi $ be assigned to the  codimension two faces of the kind $F_i\cap F_j$. \vspace {1mm}

\hspace {15mm}{\color {blue}\sf When does $X$ admit  a Riemannian metric $g$ such that
$$\mbox { $Sc_g\geq \sigma$, $mean.curv_g(F_i) \geq \mu_i$ and  
$\angle_g(F_1,F_j)\leq \pi-\beta_{ij}?$}$$}

Let moreover, $D\subset \mathbb R_+^N$, where the $N$ Euclidean coordinates  are associated with the faces  $F_i$ of $X$, be a closed convex subset, introduce  
 the following additional condition on $g$: \vspace {1mm}

{\color {black!59!blue} \sl the} $N$- {\color {black!59!blue} \sl  vector of distances} 
$\{d_i(x)=dist_g(x, F_i)\}$   {\color {black!59!blue} \sl  is in } $D$   {\color {black!59!blue} \sl   for all} $x\in X$.\vspace {1mm}

\hspace {-6mm} We   ask  when does there  exist  a $g$   with this additional condition and also\vspace {1mm}

{\color {blue} \sf  what is the homotopy type of the space of metrics} $g$ {\color {blue} \sf on} $X$,
{\color {blue} \sf  such that} 

$$\mbox { $Sc_g\geq \sigma$, $mean.curv_g(F_i) \geq \mu_i$,  
$\angle_g(F_1,F_j)\leq \pi-\beta_{ij}$  and  $\{d_i(x)\}\in D$?}$$

(For instance,  if $X$ is a  topological $n$-simplex, then an "interesting" $D$ is 
defined by $\sum _id_i(x)\geq const$.)\vspace {1mm}

One may also  try to  generalize the concept of cubical category  by allowing all kinds of combinatorial types of manifolds $X$  with corners and of attachments of $X$ to $X'$ along isometric  codimension one faces $ X\supset F \leftrightarrow   F'\subset X'$,  where the isometries $F \leftrightarrow F'$,
must match  the mean curvatures of the faces: 

 $mean,curv(F')=-mean.curv (F)$ which is equivalent to the natural metric on 
 $$X\underset{F \leftrightarrow F'}  \cup X$$
 being $C^1$-smooth.\vspace {1mm}

\hspace {4mm}{\sf \color {blue!65!red}  Is there a coherent category-style  theory along  these lines of thought?}

%%%%%%%%%%%%%%%%%%%%
%\subsection {\color {blue}  Limits  Singular Spaces, and Beyond}\label {limits  singular6}
 
%%%%%%%%%%%%%%%%%%%%

\subsection {\color {blue} Scalar Curvature under  Weak Limits of  Manifolds}\label {weak limits6}

 {\sf We show  in this section by means of  examples how  the scalar curvature may behave under limit of sequences of Riemannian manifolds.}\footnote {Our examples a similar to these  from  [Sormani(scalar curvature-convergence) 2016] and    [Lee-Naber-Neumayer](convergence) 2019].}
  \vspace {1mm}

 We saw in section \ref{stability3}   how a Riemannian manifold $X$  "emerges" as  "bubble-limit" from a "foam" (sequence)  $X_i$ obtained by taking thin connected sums of $X$ with compact Riemannin  manifolds $X_{i,\circ}$,  where this "emergence" becomes  {\it Hausdorff} or {\it intrinsic flat Sormani-Wenger}  convergence under suitable conditions imposed on $X_{i, \circ}$  and where
the scalar curvatures of $X_i$ subconverge to that of $X$ in these cases.

\vspace {1mm}

{\it Counter examples.} The inequality $Sc(X_i)\geq \sigma$ is {\it not  always  preserved} 
by   the Hausdorff and by  the  intrinsic  flat limits.

 In fact,   \vspace {1mm}
 
 {\sf all Riemannian manifolds $X$ of dimensions $n\geq 3$ can be approximated by $n$-dimensional $X_i$ with $Sc(X_i)\geq 1$ 
 
 (a) in the Hausdorff metric,
 
 (b)  in the intrinsic flat  metric.}
 (Here on speaks of closed  oriented  manifolds.) 
 \vspace {1mm}
  
 {\it Proof of} (a). All Riemannian manifolds $X$ can be  Hausdorff approximated by graphs $\Gamma$ and 
 boundaries of suitable small neighbourhoods of these  graphs embedded to $\mathbb R^{n+1}$ for $n\geq 3$,   have arbitrarily large scalar curvatures  (see section \ref {thin1}).
 
  {\it Proof of} (b) Assume $X$ bounds an orientable $(n+1)$-manifold $V$ (otherwise take the connected sum of $X$ with a small copy of $X$ with reverse orientation) and endow $V$ with a Riemannian metric and let  for which the embedding $X\to V$ is distance  preserving.

Let    $U\subset V$ be a union of small  balls, or just of small "sufficiently convex" subsets $U_j$,    the scalar curvatures of the boundaries  of which  satisfy $Sc(\partial U_j)\geq 2$.
   
 If  $vol(U)\geq vol(V)-\varepsilon$, that is easily achievable,  then the {\it flat distance} between $X$ and the boundary $\partial U=\cup_j \partial U_j$ is also $\leq \varepsilon$.
 
 What remains in order  to satisfy the definition of the  intrinsic flat distance from 
 [Sormani-Wenger(intrinsic flat) 2011]  is to modify $\partial U$ and the metric in $V$ in order to have the embedding from $\partial U$ to the complement of the interior of $U$, denoted $W=  V\setminus int(U)$,
  {\it isometric}.
  
  To do this,  let $\delta =dist_V(X, U)>0$ and  this,  take a finite  $\delta'$-dense subset   $K\subset \partial U$ for 
  $\delta'$ much smaller than $\delta$ and let $\{[k, k']\}\subset V$,  be the set of   those 
  geodesic segments in $V$ between the points $k\in K$ and$K\neq k' \in K$  which  don't 
  intersect   the interior of $U$.
  
 Assume   that 
 the segments $ [k_i, k_j]$ are mutually disjoint  and that their length are much smaller that $\delta$, say of order $\delta'$; other wise , add extra small ball to $U$.

Now,   perform the (very)  thin surgery along $ [k_i, k_j]$, that is attach thin 1-handles to $U$,  keeping the scalar curvature of the boundary of the resulting $U'$ essentially  as positive as that of $\partial U$, let $W' =V\setminus U'$ and observe  that the oriented boundary of $W'$ is
$$\partial W'= X-\partial U'$$
and   that  $vol(W')\leq \varepsilon $.

Since $\delta'<< dist (U', X)\approx \delta$,   and since the additive  difference between the "intrinsic"  metrics in $\partial U'$ and the "extrinsic" one, both  defined by shortest paths, the former in  in $\partial U'$ and in $W$ respectively, is of order $\delta'$, 
 one  can enlarge  the metric of $W$ that would make it equal to the intrinsic metric in $\partial U'$ 
 without changing the metric on $X\subset W$, and  also  only slightly changing the volume of $W$.
 
  This makes the intrinsic  flat distance  between $X$ and $\partial U' $ smaller than  $2\varepsilon$ and the proof of (b) is concluded.
   
\vspace{1mm}

The examples  (a) and (b) suggest   the following.    \vspace{1mm}

{\sf \large Definitions.}  {\it  A  Riemannian $(\alpha,\beta)$-cobordism} between closed oriented Riemannian $n$-manifolds $X_1$ and $X_2$ is an  oriented Riemannian $(n+1)$-manifold 
$$W=W_{\alpha,\beta}=\overleftrightarrow W_{\alpha,\beta}$$
with oriented boundary
$\partial W=X_1-X_2$,such that the Hausdorff  distance between  $X_1$ and  $X_2$ in $W$ satisfies 
   $$dist_{Hau}(X_1,X_2)\leq\alpha$$
and the volume of $W$ is 
$$vol(W)\leq\beta.$$
Such a cobordism can be regarded  as a morphism $W: X_1\to X_2$ with an obvious
composition for 
$X_1\overset {W_{\alpha_1,\beta_1}}\to X_2\overset {W_{\alpha_2,\beta_2}}\to X_3$:
$$W_{\alpha_1,\beta_1} \circ W_{\alpha_2,\beta_2}=W_{\alpha_1+\alpha_2,\beta_1+\beta_2}:X_1\to X_3.$$

 {\it A  Riemannian  $(\alpha,\beta, \lambda)$-cobordism} between $X_1$ and $X_2$, denoted 
$$W=W_{\alpha_1,\beta_1,\lambda}= \overrightarrow W_{\alpha_1,\beta_1,\lambda},$$
is an
$(\alpha,\beta)$-cobordism with a $\lambda$-Lipshitz retraction $W\to X\subset W$.

Here, the arrows are not invertible and the composition for $X_1\to X_2\to X_3$
 is multiplicative in $\lambda$,
$$W_{\alpha_1,\beta_1, \lambda_1} \circ W_{\alpha_2,\beta_2,\lambda_2}=W_{\alpha_1+\alpha_2,\beta_1+\beta_2, \lambda_1\cdot\lambda_2}:X_1\to X_3.$$

{\sf  Two Observations.}  

(i) Given a Riemannian manifold  $X$ (which corresponds to $X_2$  from our definition),
there exists an $\varepsilon=\varepsilon(X)>0$ such that 
the $\varepsilon$-neighbourhood $U_\varepsilon(X)\subset W$ of $X$ in $W$ admits a continuous retraction $F: U_\varepsilon \to X$,  which is  
$(1+4\varepsilon)$-Lipschitz  on the scale $>>\varepsilon$.  Moreover, 
$$ dist (F(u_1), F(u_2))\leq  dist(u_1,u_2)+ 5\varepsilon  \mbox { for all  } u_1,u_2\in U_\varepsilon(X). $$  

Indeed, there is such an $F$ which sends each $u\in U_\varepsilon$ to an almost nearest point  $x=x(u)\in X$, namely, such that    $dist (u, x(u))\leq 2\varepsilon.$\vspace {1mm}

 ({\it Probably irrelevant}){\it  Remark.} It is not hard to show (an exercise to the reader) that there exist  such retractions  $F: U_\varepsilon(X)\to X$, that are 
$(\sqrt {N} + C_X\cdot\varepsilon)$-Lipschitz on all scales,  
where $C_X$ is a constant which depends only on $X$.  (If $\varepsilon$  were allowed  to depend on $W\supset X$,the map $F$ could  be made $1+C_W\varepsilon$-Lipschitz.)

\vspace {1mm}

(ii) {\it  From $(\alpha,  \beta)$ to $(\alpha\leq \varepsilon, \beta)$}. A regularized 
$\varepsilon$-neighbourhood  $W_\varepsilon\subset W$ of $X\subset W$ is not quite a
$(\varepsilon, \beta)$-cobordism, since the embedding of the new boundary component to $W_\varepsilon$, say 
$X_\varepsilon\subset W_\varepsilon$ is not isometric. 

But if $\beta$ is much smaller than $\varepsilon$, this error can localized,  by making $\varepsilon$ smaller if necessary,  on a  small part $X'_\varepsilon\subset X_\varepsilon$, namely on the difference 
$X'_\varepsilon =X_\varepsilon\setminus \partial W=\partial W_\varepsilon \setminus \partial W $, since, by the coarea formula 
   $$\int _0^\varepsilon vol (X'_\varepsilon) d\varepsilon \leq \frac {\beta}{\varepsilon}.$$

 Moreover, if most of the volume of $X_1=\partial W\setminus X(=X_2)$ is concentrated near $X$,
namely,   
$$ vol(X_1\setminus W_\delta) <<\varepsilon \mbox {  for } \delta<< \varepsilon,$$
 e.g. 
if 
$$vol(X_1)-vol(X_2) << \delta,$$
 then, 
  by the coarea inequality,  
the boundary of $X'_\varepsilon$ can be also  made small. Then, by    filling in  $\partial X'_\varepsilon$ by a  $X''_\varepsilon$ of small volume   according to {\it the filling inequality} and then  by applying the filling inequality to  $X'_\varepsilon \cup X''_\varepsilon$, 
one   modifies   the metric in $W_\varepsilon$ such that  the embedding $X_\varepsilon\to W_\varepsilon$ becomes   distance preserving.\footnote {I didn't check the  details.}
\vspace {1mm}

{\it $(\alpha, \beta, \lambda, \sigma)$-Problem.} {\sf Given a closed oriented  Riemannian $n$-manifold 
$X$ and numbers $(\alpha> 0, \beta>0, \lambda\geq 1, \sigma>-\infty )$. 
Does there exist a cobordism  $W_{\alpha, \beta}:X_1\to X$ or   $W_{\alpha, \beta,\lambda }:X_1\to X$, where $Sc(X_1)\geq \sigma$?} 

\vspace {1mm}

{\it Open   Manifolds} The definitions of $(\alpha...)$ cobordisms $W: X_1\to X_2$ generalize to open manifolds and manifolds with boundaries, where in the latter case  $W$ comes with  a corner structure, organized as that of   cylinders $X\times [1,2]$ regarded as cobordisms between
$X\times 1$ and $X\times 2$, where  
the  flat distance  between  $X_1$ and $X_2$  defined by such a $W$     incorporates, besides 
  $vol_{n+1}(W)$, the $n$-volume of the "side boundary" of $W$,  that is $\partial_{side}W=\partial W\setminus (X_1\cup X_2)$.

\vspace {1mm}

{\it Two} {\sf  {\it {\color {red!50!black} Conjecture}s.} (1) Let the sequence $W_{\alpha_i, \beta_i}: X_i\to X$
defines a $\alpha, \beta$-convergence of $X_i$ to $X$, for
$$ \alpha_i, \beta_i\to _{i \to  \infty}0.$$  
If the scalar curvatures of all $X_i$ satisfy $Sc(X_i)\geq \sigma$, then also $Sc(X)\geq \sigma$.

(2) Let $X_i$ converge to $X$ via  $(\alpha, \beta, \lambda)$-cobordisms, that is a sequence 
$W_{\alpha_i, \beta_i, \lambda}: X_i\to X$, 
$$  \beta_i\to 0 \mbox { and } \lambda_i\to 1.$$
(The  Hausdorff distance $dist_{Hau} (X, X_i)$ and its bound $\alpha$ play no  role here.)

If  $Sc(X_i)\geq \sigma$ for all $i=1,2,...,$  then $Sc(X)\geq \sigma$ as well.\vspace {1mm}

{\it How to prove and how to improve, how to modify, and how to generalize.}  A natural approach to the proof of (1) and (2) could be as follows.

Let $\sigma\geq 0$,  assume $Sc(X,x_0)<0$,  take a $ \blacksquare$}-neighbourhood  $ \blacksquare^n\subset X$  of $x_0$ that violates the  {\color {blue} $ \blacksquare$} criterion for $Sc\geq 0$,  and then approximate $\blacksquare^n$ by    neighbourhoods $ \blacksquare_i^n \subset X_i$, which  violate the {\color {blue} $ \blacksquare$} criterion as well.

 To appreciate  the issue, let  $Y\subset X$ be a closed  volume  minimizing hypersurface
and try  to find  minimizing  hypersurfaces in  $ X_i$ that converge to $Y$ for
 $i\to \infty$ .  
 
To do this,  start with $Y_i\subset X_i$ that approximates $X_i$ for $i\to \infty$ and which can't be {\it fully}  moved  away from  their  small neighbourhood in $X_i$,    but, in the course of volume minimization, 
these $Y_i$ may, a priori,  develop"thin   fingers" protruding far away from the original $Y_i$ and carrying tiny,  yet definite positive, amounts of volume. 

The latter  problem  can be ruled out by imposing additional geometric condition(s) on $X_i$ (which is  automatic in the case of $C^0$-convergence, as in section 10 of  [G(Hilbert) 2012] and  section 4
of [G(billiards) 2014]), but in general, one has to accept these fingers  that  would  allow only   {\it weak approximation} of $Y$ by $Y_{i, min}$. (This doesn't seem to create   a serious problem for closed manifolds $X$, but may need a modification of the {\color {blue} $ \blacksquare$} criterion
for open ones.)
 
Possibly,  the validity of these   conjectures  needs additional conditions on $X_i$.
e.g. the convergence of volumes $vol(X_i)\to vol(X)$ as in section 10 of   [G(Hilbert) 2012].
 
 On the other hand,  the  bubble  example  suggests that even  a more general 
 convergence may  preserve positivity of the scalar curvature.   \vspace {1mm}

 {\it About Singular  $X$ and $W$.} The above  {\sf "convergence assisted by cobordisms"}     makes  sense for {\it pseudomanifolds} $X$, $X_i$ and $W_i$   with piecewise smooth metrics on them.\footnote {When it come to proofs, one  needs to  deal with  {\it integral current spaces}, (see [Allen-Sormani(convergence) 2020], [Sormani(conjectures on convergence) 2021] and references therein) but as far as our geometric statements are concerned, pseudomanifolds will do.}  
 
This  suggests  a provisional  definition of  $Sc^{?}(X)$, where 
$Sc^{?}(X,x)> \sigma$, $x\in X$, if and only if

{\sf there exists  a closed  neighbourhood $U\in X$ of $x$ with piecewise smooth boundary,  where $U$ admits an $(\alpha,\beta)$-approximation by Riemannian manifolds  $U_i$, $i\to \infty$, with $Sc(U_i)\geq \sigma'>\sigma$.}

 Namely,  
 
 {\sf there exist
 cobordisms  $W_i=W_{\alpha_i,\beta_i}=W_{U,\alpha_i,\beta_i}$, which are pseudomanifolds with "cornered"  boundaries 
  $$\mbox {$\partial W_i =X\cup X_i \cup \partial_ {side}W$ with 
$\partial \partial_ {side}W=\partial X\cup \partial X_i$}, $$
where 
  
$\bullet$ $Sc(X_i)\geq \sigma_i\to\sigma$,

$\bullet$ $dist_{Hau}(X_i, X)\leq  \alpha_ i\to 0$,
 
$\bullet$ $dist_{flat}(X_i, X)= vol_{n+1}(W)+vol_n(\partial _{side}W)\leq \beta_i\to 0$.}\vspace {1mm}

 Observe  that  psedomanifold $X$,  obtained,   by   $\varepsilon$-thin surgery  with  $\varepsilon\to 0$  (see      
  [BaDoSo(sewing Riemannian manifolds) 2018] and [BaSo(sequences) 2019]) may have  nasty  singularities of  codimensions $k\geq 3$, such e.g. as in   joins $X_1\vee X_2$, which  are 
 limits of thin connected sums of manifolds of dimensions $n\geq 3$.
\vspace {1mm}
 
Nevertheless, singular     $X$ with $Sc^?\geq \sigma$  in these example satisfy the {\color {blue} $ \blacksquare$}-criterion and for all we know,  
  enjoy  all  essential   geometric 
 properties  known for  smooth manifolds with $Sc\geq \sigma$.
 
But none of this is known at the present moment for general limit spaces $X$ with the  following questions remaining unresolved.\vspace {1mm}

1. {\sf Does the inequality $Sc^?(X)\geq \sigma$, that is   {\it local approximability} of  $X$ at all points $x\in X$  by (small open) manifolds $X_i=X_i(x)$ with $Sc(X_i)\geq \sigma$,   imply the  existence of  {\it global approximation} of   $X$ by  manifolds  with $Sc\geq \sigma$? 

2. Is  $Sc^?$ satisfy the additivity relation $Sc^?(X\times Y)= Sc^?(X)+Sc^?( Y)$?

3. Do $X$ which  admit   global approximations}     by  manifolds  with $Sc\geq \sigma$ 
satisfy the {\color {blue} $ \blacksquare$}-criterion? \vspace {1mm}

Notice that spaces   $X$, which do  admit global  approximation by manifolds $X_i$ with $Sc\geq \sigma>0$, satisfy (essentially) the same 
geometric bounds  as   $X_i$, because the retractions $W_i\to X$ indicated in he above (i)  defines maps $X_i\to X$ of degrees 1, which are $\lambda_i$-Lipschitz on the scale $\geq \varepsilon_i$, where $\lambda_i\to1$  and $\varepsilon_i\to 0$ for $i\to \infty$.  
(If not for "scale>0",    these $X$  would have $Sc^{\max}\geq \sigma$, see section  \ref {max-scalar5}.)

  \vspace {1mm}
 
% (?????We shall be  returning several times to these   and related    questions, e.g. in the context of the       stability problems of {\it geometric inequalities} with    
  %  $Sc\geq \sigma$, see sections \ref {3.2},  \ref{4.5} \ref{6.8}.)

   \vspace {1mm}  \vspace {1mm}  \vspace {1mm}  \vspace {1mm}  \vspace {1mm}

  It is  unproven at the present moment   that  the limits $g$ of   {\it in measure convergence} sequences $g_i\to g$  inherit positivity of scalar curvature  from $g_i$, but, probably, the  {\color {blue} $ \blacksquare$}-criterion can be used  to do this.
 (This must be easy,  if   $|\log g/g_i|\leq const<\infty$, that is  if the  Lipschitz constants $Lip_g(g_i)$ and $Lip_{g_i}(g)$  of  the identity maps $(X, g_i)\to X,g)$  and $(X, g)\to X,g_i)$ are uniformly bounded.)

%%%%%%%%%%%%%%%%%%%%%%%%%%%%%%%
\subsection{\color {blue} Scalar  Curvature beyond Manifolds Limits}\label{beyond manifolds}

%%%%%%%%%%%%%%%%%%%%%%%%%%%%%

There (at least) three  different  avenues of thought on  generalization   the concepts $Sc\geq \sigma$. 

{ \color {magenta} \textbf I.}  Finding  workable classes of  (singular) metric spaces that share their properties with smooth manifolds with $Sc\geq \sigma$, e.g. for $\sigma=0$.

{ \color {magenta} \textbf {I.A}.}   An attractive  class of such spaces $X$, that have been already  mentioned in section \ref {stability3},  is  where the  generalized sectional curvatures in the sense of Alexandrov  satisfy $sect.curv(X)\geq \kappa >-\infty$,  and where $Sc\geq \sigma$  at all
$C^2$-smooth points of these $X$.\footnote{Alexandrov spaces with $sect.curv(X)\geq \kappa$  {\color{red!50!black} seem} to provide  a perfect playground for the geometric measure theory in all dimensions and codimensions as examples with conical singularities show. Thus,   for positive $\kappa$,   they {\color{red!50!black} probably} enjoy {\it Almgren's sharp isoperimetric inequality} in all codimensions and {\it Almgren's waist estimate}.

 And as far as the scalar curvature and  minimal hypersurfaces are concerned one may try  more general  singular spaces with the Ricci curvatures bounded from below.}

{\it \color {blue}Conjecturally}, the basic properties of minimal subvarieties of all codimensions extend to these spaces, where such subvarieties of codimension one, as well as stationary $\mu$-bubbles,  serve for proving  geometric inequalities similar to the ones we have for smooth manifolds.

In fact, this is  not hard to prove  for   spaces  $X$ with   {\it isolated conical singularities}, where, as far as minimal hypersurfaces are concerned, the positivity of the sectional curvatures of the the links (bases)  of the singularities can be relaxed to positivity of the Ricci curvatures.
\vspace {1mm}

{ \color {magenta} \textbf {I.B}.} Another  class, that immediately jumps to one's mimd  is that of piecewise smooth, e.g.  spaces with iterated conical singularities, such as piecewise flat spaces, where the key issue is working out a condition for
$Sc\geq \sigma $ at conical singularities, where it may prudent to to require these spaces to be {\it rational homology manifolds.}  

\vspace {1mm}

{ \color {magenta} \textbf {I.C}.} It seems plausible that (stationary?)  minimal hypersurfaces in smooth manifolds have some generalized scalar curvatures $\geq -\infty$. 

Also doubles \DD of domains bounded by (possibly singular) minimal hypersurfaces in smooth manifolds $X$  must have (generalized) scalar   curvatures bounded from below by 
 
\hspace {40mm}  Sc(\DD) $\geq Sc X$.

\vspace {1mm}

{ \color {magenta} \textbf {I.D}.}  {\color {blue}Topologies in  Spaces of Riemannian Manifolds  Associated with Scalar Curvature }
It remains unclear  what is the weakest topology in the space of isometry classes of Riemannian manifolds   for which  the condition $Sc\geq \sigma$ is closed under limits. \footnote {See   [Sormani-Wenger(intrinsic flat) 2011],  [Sormani(scalar curvature-convergence)  2016], [Allen-Sormani(convergence) 2020], [Sormani(conjectures on convergence) 2021] and section 10.1 in [G(Hilbert)  2012] for something about it.}

Besides, properly defined weak limits  $X_\infty$ of spaces $X_i$  with $Sc\geq \sigma$, even for singular  $X_\infty$ must have a suitably defined scalar curvature 
$\geq \sigma $ as well, where  the following  may be instructive.\vspace {1mm}

{\it Example 1.} Infinite geometric connected sums $$X_\infty = \lim_{i\to \infty}  X_i \#Y_{i+1},$$
where $Y_i$ are closed Riemannian $n$-manifolds with $Sc(Y_i)>\sigma$, 
 such that
  $$\sum_{i=1}^\infty diam (Y_i)<\infty,$$ 
must have (possibly under extra conditions on the geometries of $Y_i$)
$$Sc(X_\infty)\geq \sigma.$$

{\it Recollection.}  {\sf  A  {\it thin}  (geometric)    connected sum (see section \ref {thin1})
$X_1\#X_2 $ is an abbreviation for a family   of  Riemannian manifolds $X_\varepsilon=X_1\#_\varepsilon X_2 $, for small positive  $\varepsilon\to 0$, which Hausdorff converge  to the join $X= (X,x)=(X_1, x_1) \vee (X_2, x_2)$, where the tube $T=T_\varepsilon \subset X$ (homeomorphic to $S^{n-1}\times [1,2]$)   joining the two manifolds is based on small, say  
of radii $\frac {\varepsilon}{10}$,  spheres in $X_1$ and $X_2$ around  $x_1$ and $x_2$, and such that

 $\bullet$ the complement to the $\frac {\varepsilon}{2}$-neighbourhood of $T$  in $X$,
$$X\setminus U_{\frac {\varepsilon}{2}}(T)\subset X$$ 
is isometric to the disjoint union of the complements to the $\varepsilon$-balls
 $B_{x_1}(\varepsilon) \subset X_1$ and  $B_{x_2}(\varepsilon) \subset X_2$,

  and where -- this can be arranged -- 

 $\bullet$ the scalar curvature of $Sc (X_1\#_\varepsilon X_2)\geq \sigma -\varepsilon$, in this neighbourhood     
is almost  bounded from below by the scalar curvatures of $X_1$ and $X_2$ at the  points  $x_1$ and $x_2$,
$$Sc(U_{\frac {\varepsilon}{2}}((T)\geq \min (Sc(X_1,x_1), Sc(X_2,x_2))-\varepsilon.$$

Accordingly  $X_i \#Y_{i+1}$ stands for $X_i \#_{\varepsilon_i}Y_{i+1}$ say with $\varepsilon_i=\frac {1}{2^i}$.}

\vspace {1mm}

{\it Question 2.} Let $X$ be a compact  smooth  Riemannian n-manifold and let $X_i$ be a sequence of Riemannian n-manifolds with $Sc(X_i)\geq \sigma$ and let $U_i\subset X_i$ be domains with smooth boundaries $\partial U_i$, such that \vspace {1mm}

{\sf $U_i$ admit $(1+\varepsilon_i$)-bi-Lipschitz embeddings to $X$, where $\varepsilon_i\to 0$ for $i\to\infty$.}\vspace {1mm}

{\sl What bound on the sizes of the boundaries $\partial U_i$ for $i\to 
\infty$  would imply that $Sc(X)\geq \sigma$?}\vspace {1mm}

{\it Partial Answer.}  If $\sigma=0$, then the following bound on the diameters of the connected components $comp_{ij}\subset\partial U_i$, which says that the {\it limit  Haussdorf dimension} of of $\partial U_i$ is <1,  is sufficient:
$$\sum_j diam(comp_{ij})\to 0\mbox {  for } i\to \infty.$$

In fact, this follows from the  $ \blacksquare$-{\it criterion}  in section \ref {Sc-criteria3}.

 If  $n=3$, this may be close to the necessary condition,
but if $n\geq 3$ the sufficiency of  the similar  strict  bound on the  {\it limit  Haussdorf dimension} by $n-2$, 
which says that all $\partial U_i$ can be covered by subsets $B_{ij}$, such that
$$\sum_j diam(B_{ij})^{n-2}\to 0\mbox {  for } i\to \infty,$$
remains problematic.

{\it Remark} (a)  We  make no assumption on   geometries of the complements $X_i\setminus U_i$. Thus, the relationship between $X$ and$X_i$  are, unlike any kind of distance,  non-symmetric.
(If $n=3$, one imagines  $X_i\subset U_i$  as kind of   white holes  universes emerging from $X$, where they are   seen as black holes.)

{\it  Remark} (b) The Penrose inequality suggests that  if $n=3$, then requiring that 
{\sf the areas of $\partial U_i$ tend to zero,   for $ i\to \infty$,h  would have  little effect (if ar all)  on the geometry of the limit space $X$. } But it is unclear what should be the corresponding condition for $n\geq 3$.
(Could the  areas be replaced by something like   $2$-waists of $\partial  U_i$?)

\vspace{1mm}

\vspace {1mm}

{ \color {magenta} \textbf {I.E}.}  Since the scalar curvature is additive under finite    Riemannian products it is tempting 
to extend the idea to infinite products and iterated fibrations, and to find  geometric  meaning of the inequality $Sc\geq \sigma$ for infinite dimensional Hilbertian manifolds, such as spaces of maps between Riemannian manifolds.  But no plausible conjecture is known in this direction.  

\vspace {2mm}

{ \color {magenta} \textbf {II}.} Instead of the spaces one may focus on  analytic techniques used for the study of
$Sc\geq \sigma$, 
in particular  the index theory for the Dirac  operator and the geometric measure theory and  search for  generalisations (unification?)  of these that  would be applicable to singular spaces.\vspace {2mm}

{ \color {magenta} \textbf {III}.} One may think of manifolds with $Sc\geq \sigma$ and the methods used for 
their study as  as geometric/analytic  embodiment of  certain algebraic formulae behind these, such as the {\color {blue} GaussTheorema Egregium coupled with the second variation formula} and the {\color {blue}  Schroedinger-Lichnerowicz-Weitzenboeck-(Bochner}) formula coupled with the  formula(s) involved in the local proof of the index theorem.

{ \it Conceivably},  there may exist alternative  implementations of these formulas  in 
 categories  which are quite different from those of manifolds and/or of metric spaces, where, e.g.  the objects  are  represented by functors from a category of "decorated graphs" 
to that of measure spaces as in the last section of [G(billiards) 2014] or

 {\color {blue}\large\sf  something else, something  far removed from the present day idea 
 
 of  what  a geometric space is}.

%%%%%%%%%%%%%%%%%%%%

\section { Metric Invariants Accompanying Scalar Curvature}\label{metric7}

%%%%%%%%%%%%%%%%%%%

Many  invariants of metric spaces $X$ can be expressed in (quasi)-category theoretical language,   e.g.  in terms of $\lambda$-Lipschitz maps between   
$X$  and a  "measuring rod"  space (or spaces) $\underline  X$. \footnote  {There  are  also some  questionable,   albeit  sometimes  coming in  handy,   ad hoc  invariants,   such as the "injectivity radius",   but these are useless as far as the scalar curvature is concerned.}

In fact, the distance function in  $X$ is fully encoded   by the sets of  1-Lipschitz functions,  i.e. distance non-increasing maps $X\to \mathbb R$:\vspace {1mm}

(dist$_{cnrt}$)  {\sf  the distance  $dist(x_0,x_1)$ is  (obviously) equal to the {\it supremum} of  numbers $d\geq 0$, such that $X$ admits a 1-Lipschitz map $f:X\to 
\mathbb R$, such that $f(x_0)=0$  and $f(x_21=d$.} \vspace {1mm}

Alternatively,  $dist(x_0,x_1)$ can be defined {\it covariantly} via maps  of two point subsets from  $ \mathbb R$  to
$X$, as   follow:\vspace {1mm}

(dist$_{cov}$){\sf the distance  $dist(x_0,x_1)$ is equal to the  {\it infimum} of $d\geq 0$, such that $\{0,d\}$ admits a 1-Lipschitz map to $X$  with the  {\it  image }   $\{x_0,x_1\}\subset X$.}

Similarly,  one can define the volume of a   connected Riemannian $n$-manifold  $X$ as  \vspace {1mm}

(vol$_{cov}$) the {\it infimum}  of numbers $v=d^n$, such that that
{\sf $X$ receives  a   smooth {\it locally volume non increasing} map $f$  (i.e. $||\wedge^ndf||\leq 1$)  from the cube $[0,d]^n$ {\it onto} $X$.}\vspace {1mm}

And  -- this is closer to invariants  used   in the study of  scalar curvature  -- one can   define $vol(X)$ of a closed connected manifold $X$ {\it contravariantly} as \vspace {1mm}

(vol$_{cntr }$)  {\sf the supremum of volumes $v=(2n)d^n$ of the boundaries of  the cubes $[0,d]^n $,   which receive  {\it non-contractible} 
 locally volume non increasing piecewise  smooth maps from $X$.}

\vspace {1mm}

{\it Exercise.} Prove (vol$_{cov}$)  and (vol$_{cntr }$).
%%%%%%%%%%%%%%%%%%%%%%%%%%%%

\subsection{\color {blue}  Multi-Spreads of Riemannianan Manifolds:   $\square^\perp$  and $\tilde\square^\perp$} \label{multi-spreads} 

%%%%%%%%%%%%%%%%%%%%%%%

Let $\tilde U$ be a  compact $n$-dimensional  manifold possibly with a boundary,   let $g$ be  Riemannian metric on $\tilde U$ and let $\tilde h\in H_{n-k}(\tilde U)$ be a homology class of codimension $k$.

The  $\square^n$-inequality   from  section \ref{multi-width3}  for widths  of cubes  with  metrics with $Sc\geq \sigma$  motivates the following.

{\it Definition.} {\it The $\square^\perp$-spread} of  a homology class $\tilde h\in H_{n-k}(\tilde U)$, denoted  
$$\square^\perp (\tilde h)=\square_g^\perp (\tilde h),$$  
is  the supremum   of the  numbers $d\geq 0$, for which  there exists

{\sf a continuous proper (boundary-to-boundary) map $\psi=(\psi_1,...\psi_i,...\psi_k): U\to  [-1,1]^k$,   
$\psi_i: \tilde U\to [-1,1]$,
such that 

(a)   the homology class of the  $\psi$-pullback of a point\footnote{If $\psi$  is smooth       this an actual pullback of a generic point.}
 is equal to $\tilde h$, symbolically
$$\psi^\ast[t]=\tilde h, $$
where  $ [t] \in H_0([-1,1]^k)$,  $t\in   [-1,1]^k$, is the homology class of a point in  $ [-1,1]^k$;
  
(b)  the distances  between the pullbacks of the opposite faces in the cube $[-1,1]^k$, 
$$d_i=dist_g(\psi_i^{-1}(-1),\psi_i^{-1}(1)), i=1,...,k,$$
are bounded from below by the following inequality 
$$ \left (\frac {1}{k}\sum_{i=1}^k \frac {1}{d_i^2}\right)^{-\frac {1}{2}}\geq d,$$
that is 
$$\frac {1}{k}\sum_{i=1}^k \frac {1}{d_i^2}\leq \frac {1}{d^2}.$$}

(Equivalently, one could require the maps $\psi_i$ to be $d_i^{-1}$-Lipschitz.)\vspace{1mm}
 
Next define $\tilde\square^\perp(h)\geq \square^\perp(h) $ of a homology class  of codimension $k$ in a Riemannian $n$-manifold $X$, possibly non-compact and with a boundary, denoted 
 $h\in H _{n-k}(X)$, as the supremum of the  numbers $d\geq 0$, such there exist 
 
 (i)  a Riemannian manifold  $\tilde U$,  
 
 (ii)  a homology class $\tilde h\in H_{n-k}(\tilde U)$ with $\square^\perp( \tilde h)=d$, 
 
(iii)  a  locally isometric map $\phi:\tilde U\to X$,  for which the induced homology homomorphism 
$\phi_\ast :H_{n-k}(\tilde U) \to H_{n-k}(\tilde X)$
sends $\tilde h$ to $h$, in writing:
$\phi_\ast(\tilde h)=h.$

(This definition make sense for an arbitrary metric  $dist$  on $\tilde U$.) \vspace{1mm}

{\it Topological Remark.}  The  $ \tilde\square^\perp$-spread of $h$ vanishes if and only if 
none  of $\tilde h$ is homologous to any point-pullback, that is the case, for instance,
if $h$ has non-zero self intersection $h\frown h\neq 0$.

On the other hand, by a theorem of Serre on cohomotopy groups,
If $k$ is odd, or if $h\frown h=0$, then some non-zero multiple of $h$, say $Nh$, has  $\square^\perp(Nh)>0.$ \vspace{1mm}

 Say that $h$ is  {\it $\tilde \square^\perp$-spread infinite} or that $h$ has infinite $\tilde\square^\perp$-spread if  $\tilde\square^\perp  ( h)=\infty.$\vspace{1mm}
 
 {\it Define  $\square^\perp$-spread and  $\tilde\square^\perp$-spread} of a compact  connected  orientable $n$-dimensional Riemannian manifold  $X$,   possibly with a boundary,  denoted $\square^\perp(X) $ and   $\tilde\square^\perp(X)$, as the $\square^\perp$- and  $\tilde\square^\perp$-spreads of the {\it zero dimensional homology class} $[x]$ of a single point   $x\in X$.

 {\it If $X$ is non-compact}   define  $\square^\perp(X)$ as 
 $\limsup  \square^\perp(X_i)$, $i\in I$, for all  compact $n$-submanifolds $X_i$ exhausting $X$  and 
   let 
   
   $\square^\perp(j\cdot X)$ and $\tilde \square^\perp(j\cdot X)$  denote these spreads of the {\it $j$-multiple $j\cdot [x]\in H_0(X)$}
 of the homology class of $x\in X$.\vspace{1mm}

 Observe that if $X$ is  {\it compact without boundary} then  $\square^\perp(X)=0$
  and  that
the {\sl $\tilde\square^\perp$-spread}, unlike the $\square^\perp$-spread, {\sl of the universal covering} $\tilde X$ is {\it equal} to that
  of $X$.
 
 Thus, for instance, the $n$-torus $\mathbb T^n$ is   $\tilde\square^\perp$-infinite,
     $$\tilde\square^\perp(\mathbb T^n)=\square^\perp(\mathbb R^n)=\infty, \mbox{ 
     while } \square^\perp(\mathbb T^n)=0.$$
  
And,  in general, if   a homology class $h$ is representable by a simply connected cycle in $X$, that is 
if $h$ is equal to the image of a class $\tilde h\in H_{n-k}(\tilde X)$ under  (the homology homomorphism induced by)   the  universal covering map  $\tilde X\to X$, then 
$$\tilde\square^\perp  \tilde (h)= \tilde\square^\perp  ( h).$$

   \vspace{1mm}

Say that  $X$ is {\it $\tilde\square^\perp$-spread infinite} or that $X$ has infinite $\tilde\square^\perp$-spread if  $\tilde\square^\perp  ( X)=\infty,$ and observe that this property is equivalent to {\it iso-enlargeability} from [G(inequalities)  2018];

 \vspace {1mm}

As far as the scalar curvature is concerned,  we are interested in  {\it lower bounds} on $\tilde\square^\perp (h)$,  which are usually easily available, e.g.  in the examples  (1)-(4)  below.  \vspace{1mm}

(1) The ball $B^n(R)\subset \mathbb R^n$ has 
$$\square(B^n(R)) \geq  \frac {2R}{\sqrt n}.$$

(2) Closed connected surfaces $X$ with {\it infinite} fundamental  groups $\pi_1(X)$  are (obviously) {\it $\tilde\square$-spread infinite},  i.e.  $\tilde \square  ( X)=\infty.$
\vspace{1mm}

(3)  The spread of an  $n$-manifold $X$ with non-empty boundary is (obviously) related to the inradius   
 $inrad(X)=\sup_x(dist, (x,\partial X)$ by the following inequality.
 $$ \tilde \square (X)\leq 2\sqrt {n} \cdot inrad (X),$$
 where the equality holds for $X= [0, 2r]\times \mathbb R^{n-1}.$
 
 Furthermore, if $n=2$,  then 
$$ \tilde \square (X) (X)\geq \sqrt 2\cdot inrad (X).$$
This is seen with the universal covering $\tilde U$ of $X$ minus the furthest point from the boundary,  where 
$inrad (\tilde U)=\frac {1}{2}\cdot inrad (X)$.

In particular,  {\it complete non-compact} surfaces $X$ are  {\it $\square$-infinite}.

\vspace {1mm}

(4) Surfaces $X$ {\it homeomorphic to the 2-sphere}  have   $\tilde \square\perp  ( X)\geq  \sqrt 2\cdot diam(X)$,  that is seen by evaluating $\square$   of the universal cover $\tilde U$ of $X$ minus two furthest points in it.

And  since the cut  loci to all  points  $x$ in this $X$ contain  {\it conjugate points} to $x$,  the inradii of surfaces  $\tilde U$, which locally isometrically immerse to $X$ are bounded by $diam (X)$; hence, 
$\tilde \square\perp  ( X)\leq  2 \cdot diam(X)$.

\vspace {1mm}

(All compact simply connected manifolds $X$ have $\tilde \square\perp  (X)\leq C<\infty$, but  no bound on $C$ by the diameter is possible for $n>2$.

 In fact, an arbitrary closed  $n$-manifold $X$,   $n\geq 3$,   e.g.  $X=S^3$, admits,  by geometric  surgery argument,    
   Riemannian metrics $g_C$ for all $C>0$,  with $diam_{g_C}(X)= 1$,   with
$sect.curv(g)\leq \frac {1}{100n^2C^2}$ and, thus,   with  $\tilde \square\perp  (X)>C$, where 
 $\tilde U$ is the $R$-ball in the tangent space $T_{x_0}(X)$ for $R=nC$, sent to $X$ by the exponential map 
and endowed with the Riemannian  metric induced from that on $X$.)

(5)  \textbf {Product  Inequality}.  {\sf Let $\underline X_i$,  $ i=1,2,...,m$, be Riemannian manifolds of dimensions $n_i$ , possibly with boundaries,  
non-compact and  non-complete  and let  $f_i:X\to  \underline X_i$  be proper
 (infinity-to-infinity boundary-to-boundary) maps and let 
$$f=(f_1,...,f_m) :X\to  \underline X=\underline X_1\times ...\times \underline X_m.$$}

{\it Then the $\tilde\square^\perp$-spread of the homology class $h= f^\ast[\underline x])  \in H_{n-k}(X)$, $k=dim(X)-dim(\underline X)$ of the point-pullback $f ^{-1}(\underline x)$ of $f$ is bounded from below by the 
$\tilde\square^\perp$-spreads of $\underline X_i$ as follows.
$$\tilde\square^\perp(f^\ast[\underline x])\geq \left ( \frac {1}{n}\sum_{i=1}^m \frac {n_i}{\tilde\square(\underline X_i)^2}\right)^{-\frac {1}{2}}$$
that is
$$\frac {1}{\tilde\square^\perp(f^\ast[\underline x])^2}\leq \frac {1}{n}\sum_{i=1}^m \frac {n_i}{\tilde\square(\underline X_i)^2}.$$}

{\it In fact},  the $\square^\perp$-spread of intersection of cycles $h_1$ and $h_2$ of codimensions $k_1$ and $k_2$, denoted $d=\square^\perp(h_1\frown h_2)$, satisfies 
$$\frac{1}{d^2}\leq\frac {1}{k_1+k_2} \left( \frac{k_1}{(\square^\perp(h_1))^2}+\frac{k_2}{(\square^\perp(h_2))^2}\right).$$
Then the  the proof   for  $\square^\perp $   follows by induction on $m$  and the corresponding inequality for 
$\tilde\square^\perp$  follows.

For instance,   the  $\square$-spread of the rectangular solid,   $d=\square \left (\bigtimes_{i=1}^n [0,d_i]\right)$ satisfies 
  $$\frac{1}  {d^2}\leq  
   \frac {1} {n}\sum_i \frac  {1} {d_i^2}.$$

\vspace {1mm}

(6)  Connected sums of compact  connected {\it $\tilde\square^\perp$-spread infinite} manifolds $X$ with  complete manifolds are
(obviously) {\it $\tilde\square^\perp$-infinite.} 

In particular, complete metrics on  a   $\tilde \square^\perp$-spread infinite  compact  manifold $X$ minus a point  are    $\tilde\square^\perp$-infinite.

(7)  Let $X$ be a complete connected non-compact manifold and $Y\subset X$ be a compact connected  submanifold of {\it codimension 1.}

If  the  inclusion homomorphism $\pi_1(Y)\to \pi_1(X)$ is {\it injective}, then
$$\tilde \square^\perp(X) \geq \tilde \square^\perp(Y).$$
In particular, if 
  $Y$ is  {\it $\tilde\square^\perp$-spread infinite} then also  $Y$ is  {\it $\tilde\square^\perp$-spread infinite}.

(8) Let $X$ be a  connected complete non-compact  manifold and $Y\subset X$ be a compact 
connected {\it  $\square$--spread infinite} submanifold of {\it codimension 2},  such that  the  inclusion homomorphism $\pi_1(Y)\to \pi_1(X)$ is {\it injective.}

If the real homology class  of the $\varepsilon$-circle $S^1_y(\varepsilon) \subset X\setminus Y$   in the  normal plane to $Y$ doesn't vanish, then 
$$\tilde \square^\perp(X) \geq \tilde \square^\perp(Y).$$

\vspace {1mm}

Unlike  lower bunds,  upper bounds  on $\tilde \square$ find no, at least, no immediate, applications to  scalar curvature. What make them   amusing is an  unexpected complexity of sharp evaluation of $\tilde \square$,   
and even of  $\square\leq \tilde \square$,
  in  simple  examples  indicated
below, where there   are more questions than answers. \vspace {1mm}

 (A) The $d$-cube $[0,d]^n$ satisfies
 $$ \square [0,d]^n=d.$$
 
 {\it Proof.}  The inequality $ \square [0,d]^n\geq d$ is obvious. (It is the simplest case of  the product inequality.)
 
The lower bound follows  from 
    {\it Besicovich-Derrick}  \& geometric/arithmetic means inequalities, which shows that 
    $(\square (X))^n\leq vol(X)$, $n=dim(X)$, for all Riemannin manifolds $X$.
 
{\it {\color{red!50!black}Probably},}   
$ \tilde \square [0,d]^n$ is equal to $d$ as well,  and 
 the following  more general property of $\tilde \square$  also looks plausible.

(B)  {\it {\color {red!50!black} Conjecture}}.   {\sf All {\it convex} domains $X\subset \mathbb R^n$ satisfy   $\tilde \square (X)=\square (X)$.}
\vspace{1mm}

(C) The  universal covering $\tilde U$ of the  2-{\it ball $B(r)\subset  \mathbb R^2$ minus the center}  (obviously) satisfies:
            $$ \tilde \square (\tilde U)=  \square (\tilde U)=\sqrt 2 r,$$
which is equal  to the $\square$-spread of the (inscribed) square $\left[-\frac {r}{\sqrt 2},\frac {r}{\sqrt 2}\right]\subset B(r). $

(D)  {\it {\color {red!50!black} Conjecture}}. {The  rectangular solid $\bigtimes_{i=1}^n[0,d_i]\subset \mathbb R^n$ satisfies 
$$\tilde \square\left  (\bigtimes_{i=1}^n([0,d_i]\right)\leq\left( \sum_i\frac {1}{d_i^2}\right)^{-\frac {1}{2}}.$$
(It is obvious that  $\tilde \square\left  (\bigtimes_{i=1}^n([0,d_i]\right)\geq\left( \sum_i\frac {1}{d_i^2}\right)^{-\frac {1}{2}}.$)

(E) {\it {\color {red!50!black} Conjecture}.}  The $\tilde \square$-spread of the ball $B(r)\subset \mathbb R^n$ is equal to the (conjectural) 
 $\tilde \square$-spread  of the 
 inscribed  cube  (where,  obviously, the latter is bounded by the former): 
 $$\tilde\square(B(r))=\square(B(r))=  \frac {2r}{\sqrt n}.$$
Moreover,  if a $\tilde U$  admits a locally isometric immersion to $B(r)$,  if
$$ \tilde\square(B(r))=\square(B(r))=  \frac {2r}{\sqrt n}= \frac {2r}{\sqrt n}$$
and if $n\neq 2$, then $\tilde U=B(r)$.

The following is a  step toward   a  weaker, also conjectural,   inequality   

$ \square\left (\bigtimes_{i=1}^n([0,d_i]\right)\leq \left( \sum_i\frac {1}{d_i^2}\right)^{-\frac {1}{2}}.$

(F)  {\it \textbf {Proposition}.} {\sf Let $X=\bigtimes_{i=1}^n[0,d_i]\subset \mathbb R^n$   and $X'=\bigtimes_{i=1}^n[0,d'_i]\subset \mathbb R^n$ be rectangular solids, such that     either}

 (i) {\sl there exists a proper (boundary to 
 boundary)  1-Lipschitz map of {\it odd} degree\footnote{{\it Non-zero}  degree should be OK, but I only vaguely see how to  prove this.}
 $X'\to X$,}

or 

(ii) {\sl there exists a smooth   locally expanding  (non-decreasing the lengths of
 smooth curves)   embedding $X \to  X'$.}

Then 
$$\left (\sum_i \frac {1}{d_i^2}\right)^{-\frac {1}{2}}\leq \left (\sum_i \frac {1}{(d'_i)^2}\right)^{-\frac {1}{2}}.$$}

  {\it This   is an immediate   corollary} of the following  (a) and (b)   pointed out to me by 
Roman Karasev,\footnote  {Besides (a) and (b) Roman has  made a few other illuminating remarks, including counter examples to  some of my naive suggestions on this subject matter.}  where (a)   depend on 
 the concept of a $k$-dimensional {\it $\mathbb Z_2$-waist},  denoted  $waist_k(X)$,  of a Riemannian manifold $X$  (possibly with a boundary). This   is a numerical invariant, which is (almost by definition,
see next section) 
is

{\it non-increasing 
under    proper 1-Lipschitz map of {\it odd} degree $X'\to X$  and  non-decreasing under  locally expanding     embedding $X \to  X'$,i.  e. $waist_k(X')\geq waist .(X)$ in these cases}

(The $\square$-spread may decrease under  locally isometric      embedding $X \to  X'$. For instance, concentric R-balls in the unit sphere, satisfy:   $\square(B(R))\geq  \square(B(R'))$ for $\frac {\pi}{2}\leq R \leq R'\leq\pi$.)

(a){\sl The {\it $\mathbb Z_2$-waists}   of the solids $\bigtimes_i [0,d_i]$,  $d_1\leq...\leq d_i\leq ...\leq d_n$,
satisfy
$$waist_k(\bigtimes_i [0,d_i])=d_1\times  ...\times  d_k\mbox {  for all  } k=1,...,n.$$}

{\it This is stated}, in a slightly different form  in    corollary 5.3 in  [Klartag(waists) 2017] (also  see  [Akopyan-Karasev(non-radial Gaussian) 2019] ).

(b) {\sl If positive  numbers  $d_1\leq...\leq d_i\leq ...\leq d_n$, and $d'_1\leq...\leq d'_i\leq ...\leq d'_n$ satisfy
$$d_1\times ...\times d_k\leq d'_1\times ...\times d'_k$$   for all $k=1,...,n$,
 then 
 $$\sum_i \frac {1}{d_i^2}\geq \sum_i \frac {1}{(d'_i)^2}.$$}

{\it Indeed,}  since  the numbers  $l_i=-2\log d_i$ dominate  $  l'_i=-2 \log  d'_i$,
i.e.
$$\sum_i^k l_i\geq  \sum_i^k l'_i, \mbox { } k=1,...,n.$$
 the {\it Karamata inequality}  applied to to (convex !) function $\exp l$ yields the required inequality: 
$$\sum_i \frac {1}{d_i^2}=\sum_i^n \exp l_i \geq \sum_i^n \exp l'_i  =\sum_i \frac {1}{(d'_i)^2}.$$

(G) {\it Generalizations.} The above argument yields similar monotonicity of $\Sigma_\alpha=\left (\sum_i  d_i^\alpha\right)^\frac {1}{\alpha}$ for all {\it negative}  
$\alpha$, but  it  is  unclear which (if any) of these $\Sigma_\alpha$ is increasing under (globally) {\it non}-one-to-one locally expanding maps between solids. (This monotonicity may fail for small $|\alpha |$.)

Also,  waist evaluation in 
in [Akopyan-Karasev(tight estimate) 2016] (corollary 4)
yields similar monotonicity for maps between  ellipsoids  with principal axes $d_i$ and $d'_i$,  and between solids and ellipsoids.

(H) {\it Questions.} Is there an "effective" set of inequalities between the numbers $d_i$ and $d'_i$  necessary and sufficient for the  existence of an (affine) isometric  embedding  from solid to solid,  $\bigtimes_i [0,d_i]\to   \bigtimes_i [0,d'_i]$,  and from  ellipsoid to solid  or to  ellipsoid?  

(From a  geometric perspective, one would rather have  this kind of inequalities for (non-affine)  locally injective and/or non-injective locally expanding maps, while from the convexity point of view it is   natural to study the convex set of all 
affine embeddings between  convex sets,  with a  special consideration of affine self-embeddings of such sets.)

(I)  {\it {\color {red!50!black} Conjecture}.}   {\sf The manifolds  $\underline X_i$ from the above  product inequality  satisfy the equality:
$$\tilde\square^\perp(f^\ast[\underline x])= \left ( \frac {1}{n}\sum_{i=1}^m \frac {n_i}{\tilde\square(\underline X_i)^2}\right)^{-\frac {1}{2}}.$$}

(This generalizes the above  conjectural formula $\tilde \square\left  (\bigtimes_{i=1}^n([0,d_i]\right)\leq\left( \sum_i\frac {1}{d_i^2}\right)^{-\frac {1}{2}}$,  and,  in fact,  may follow from such a formula.)

%%%%%%%%%%%%%%%%%%%%%%%%%%%%%%

 \subsection  {\color {blue}Manifolds with Distinguished Side    
 Boundaries and  Gauss-Bonnet/Area Inequalities} \label{Gauss-Bonnet7}
 
 %%%%%%%%%%%%%%%%%%%%%%%%%%%%%%%

 Let $\partial_{side}\subset \partial X$ be an open subset in the 
 boundary of our Riemannian $n$-manifold $X$ and let us generalize the definitions of   $\square^\perp$  and  $\tilde\square^\perp$ for  a {\it relative homology} class $h\in H_{n-k}(X, \partial_{side})$,  as earlier but with 
 {\it side-proper}  rather than just proper  maps $\psi$. 
 
 Namely:  
 
$\bullet$  {\sf the auxiliary $n$-manifold $\tilde U$ also comes with a distinguished {\it side boundary},  denoted 
 $\tilde \partial_{side}\subset \partial \tilde U$;
 
 $\bullet$   continuous  maps 
 $$\psi=(\psi_1,...\psi_i,...\psi_k): U\to  [-1,1]^k,   \mbox { } 
\psi_i: \tilde U\to [-1,1],$$ must be {\it side-proper},   which means that they   send the complement $\partial X\setminus \partial_{side}$ to the boundary of the cube,
 $$\psi (\tilde U\setminus \tilde \partial_{side})\subset \partial [-1,1]^k;$$
 
   $\bullet$ locally isometric maps $\phi:\tilde U\to X$ must send $\tilde \partial_{side}\to  \partial_{side}\subset \partial  X$.}
 
 \vspace{1mm}
 
 \textbf {Theorem: G-B-Inequality.} {\sf  Let $X$ be   a Riemannian manifold of dimension $n$ and  with a distinguished open subset 
 $\partial_{side}\subset \partial X$ and let $h \in H_2(X, \partial_{side} )$ be a relative  homology class.}
 
  {\it Then $h$  can be represented by an immersed  smooth surface $\Sigma\subset X$,  the   boundary of which is  contained in
$\partial_{ side}$  and   such that  the integrals of the scalar curvature  of $X$ over  all connected components $S$ of $\Sigma$ and of the mean curvature\footnote{Our sign convention is such that the boundaries of {\it convex} domains have {\it positive} mean curvatures.}   of $\partial_{side}$ over $\Theta=\partial S $ are related to the
 multi-spread $\tilde \square$ of the pair  $(X, \partial_{side})$  by the following inequality.

satisfy:
$$
 \int_S  Sc(X,s)ds+ 2\int_{ \Theta} mean.curv(\partial_{side}, \theta)d \theta \leq 4\pi\chi( S)+ C_{\tilde\square}\cdot area(S),$$
where $\chi( S) $ is the Euler characteristics of $S$ and
$$C_{\tilde\square}=\frac {4(n-1)(n-2)\pi^2}{n}\left (\tilde\square(X,\partial_{side})\right)^{-2}. \leqno{\mbox{\eye}}$$}
The proof of this given in  [G-Z(area) 2021])   combines (a version of) the  the  $\square^n$-inequality    for widths  of cubes (see   sections \ref{multi-width3}  and \ref{separating5})  with the argument similar to that in  [Zhu(rigidity) for the proof of  the  sharp  equivariant area inequality (see  section \ref {warped stabilization and Sc-normalization2}), where we consider only the  case of $n\geq 7$.

The case $n=8$, which needs a version of Natan Smale's  generic regularity result we postpone till another paper, 
while $n\geq 8$ needs a generalization of Lohkamp's or of   Schoen-Yau's regularization theorems.

\vspace {1mm}

{\it Remarks.}(a)  This  theorem, as stated, is non-vacuous  {\it only if}   $Sc(X)\geq 0$ and $mean.curv(\partial X)\geq 0$; otherwise,    all relative homology classes can be represented by surfaces with {\it arbitrarily small}  integrals,
$\int_S  Sc(X,s)ds$   or $\int_{ \Theta} mean.curv(\partial_{side}, \theta)d \theta.$   

(b) If $Sc(X)\ngeqslant 0$ or 
$mean.curv(\partial X) \ngeqslant 0$,  a rough, (but  meaningful)  \eye kind inequality is possible if, for instance,  $Sc(X)\geq -1$, $mean,curv (\partial X)\geq -1$ and 

the sectional curvature of $X$ is bounded by  $+1$ in the 1-neighbourhood of the region in $X$, where $Sc(X,x)<0$ and/or $mean.curv (\partial X  x )<0$;

 the principal  curvatures of $\partial X$  are bounded by 1  in the 1-neighbourhood of the region in $\partial X$, where $mean.curv (\partial X  x )<0$.

(c) If  the boundary of $X$ is mean convex,   $mean.curv (\partial X  x )\geq 0$  and $Sc(\geq -1$, then   the proof of  \eye, which delivers  {\it area minimizing} surface in the class $h\in H_2(X, \partial_{side}$, provide a non-trivial {\it lower bound} on the area-norm of this class.

But , it is unclear what should be a {\it correct} version of \eye  for $n\geq 3 $, where     $Sc(X)\ngeqslant 0$ and/or 
$mean.curv(\partial X) \ngeqslant 0$.\vspace {1mm}

(d) If  we allow  $C_\square =\frac {4(n-1)(n-2)\pi^2}{n}\left (\square(X,\partial_{side})\right)^{-2}$instead of $C_{\tilde\square}$ in \eye, then the  required  surface $\Sigma\hookrightarrow X$ may be assumed {\it embedded}.\vspace {1mm}

(e) A  kind of  \eye (systolic)   inequality  for metrics with $Sc> 0$ on 
$S^2\times S^2$ was established  in  [Richard(2-systoles) 2020; also a version of  Zhu's    sharp  equivariant area inequality  for manifolds with $Sc\geq 0$  and with  mean convex boundaries  is proven in [Barboza-Conrado](disks) 2019].\vspace {1mm}

 {\sf \large Examples of Corollaries.} \textbf  A.  {\sf Let $ X$ be  a Riemannian manifold  diffeomorphic   to the product $\pentagon \times \mathbb R^{n-2}$, where
$\pentagon$} is a planer $j$-gon, or, more generally,  let  $X$  be  a {\it manifold with $j$  corners}, which  admits a proper (boundary-to-boundary, infinity-to-infinity) map of  positive degree $ f:X\to \pentagon\times \mathbb R^{n-2}$, such  that 
the images of these corners in $\partial \pentagon \times \mathbb R^{n-2}$  have  non-zero intersection indices with the circles $\partial \pentagon \times t\subset \partial \pentagon\times  \mathbb R^{n-2} $, $ t\in \mathbb R^{n-2}.$

{\sl  If $Sc(X)\geq 0$, if the  mean curvature of $\partial X$ away from the corners is $\geq 0$ and if the dihedral angles $\angle_i$, $i=1,...j$, of $X$ at the corners satisfy
 $\angle_i\leq \alpha_i\leq \pi$, where 
 $$\sum_{i=1}^j \pi-\alpha_i  > 2\pi,$$}
then 
\hspace  {20mm} {\it the map $f$ can't be (globally)  Lipschitz.}

  Moreover, 
  
 \hspace  {15mm} {\it there exist  sequences  of   points $x_i, y_i \in X$, such that  
  $$dist (x_i,y_i)\leq const<\infty,\mbox {  and 
$dist (f(x_i),f(y_i))\to \infty$ for $i\to\infty.$}$$}

\vspace {1mm}

{\it In fact} the inequality \eye  holds for manifolds with {\it non-smooth} boundaries (here $\partial X= \partial_{side}X$), if  the mean curvature understood in a suitable distribution way. But to  fully make sense of this  one  needs additional data on regularity of the boundary $S=\partial \Sigma$. 

However, for just  keeping track of the inequality $\sum_{i=1}^j \pi-\alpha_i  > 2\pi$  in the integral $\int_{ \Theta} mean.curv(\partial_{side}, \theta)d \theta$,     one can simply  smooth the boundary 
$\partial X$ in an obvious manner and  thus approximate $X$ by  domains $X_\varepsilon\subset X$ with   smooth boundaries.  Then  \eye,  applied to these 
$X_\varepsilon$, 
yields  the corollary  for $\varepsilon\to 0$. (We suggests the reader would fill in the details of this argument.)

\textbf B. {\sf  Let $\underline S$ be  compact connected surface with a boundary, $\hexagon$  be  a
 planar $k$-gon,  $X$ be a Riemannian $n$-manifold and let 
$$f: X\to \underline X=\underline S\times \hexagon \times \mathbb R^{n-4} $$}
 be a diffeomorphism. (A continuous proper  map of degree one   will do).
 
 {\sf Define  the side boundary of $X$ as the one corresponding  to the boundary  $\partial \underline S$, 
 $$\partial_{side}(X)=f^{-1}( \partial \underline S\times  \hexagon \times \mathbb R^{n-4}),$$ }
(this $\partial_{side}$ is smooth)  and let $\partial_\angle X\subset \partial X$  be the   "cornered part" of the boundary of $X$, that is  $$\partial_\angle X =f^{-1} (\underline S\times \partial \hexagon\times \mathbb R^{n-4}),$$ 
 where the faces and the "corners" of  $\partial_\angle X$ correspond to the edges  and the vertices of $\hexagon$.

 Let  the following four conditions be     
  satisfied.

  $(\bullet)$ {\sl The map $f$ is  roughly  asymptotically  Lipschitz-like: 
$$dist(f(x) f(y))\leq \mathcal L(dist(x,y))$$ 
for some continuous function $\mathcal L(d)=\mathcal L_f(d)$, $d\geq 0$, and all $x,y\in X$ with
$dist(x,y)\geq 1,$  e.g. $||df||\leq const <\infty$.}

 $(\bullet\bullet)$ {\sl The faces of $\partial_\angle X$ are mean convex, i.e have positive mean curvatures.}

  $(\bullet\bullet\bullet )$   {\sl The  
 dihedral angles  $\angle_i$, $i=1,2, ...k$, between the faces of $\partial_\angle X$ at all points in  the 
 "corners"   are all bounded as follows.
$$\angle_i\leq \frac {2\pi}{l}$$  
where $l$ is a positive integer, such that

if $k=3$, then  $l\geq 6 $ i.e.  $\angle_i\leq \frac {\pi}{3},$

 if $k=4,5$, then $l\geq 4$  i.e.  $\angle_i\leq \frac {\pi}{2},$

 if $k\geq 6$,  then $l\geq 3$,  i.e.   $\angle_i\leq \frac {\pi}{2}$.}

$(\bullet\bullet\bullet\bullet)$  {\sl  Either $l$ is  even or let,  for every  pair of adjacent $(n-1)$-faces in $\partial_\angle$,
say $\partial_i$ and $\partial_{i+1}$   there exist  
an isometric, i.e. preserving the induced Riemannian metric,  involution of $\partial_\angle$, which   interchanges these faces,  $\partial_i\leftrightarrow\partial_{i+1}$, and fixes the corner  $\partial_i\cup \partial_{i+1}$ between them.\footnote   {Although,    the existence of this involution   is probbaly unneeded in the present  case,   it suggests a generalisation   of the  
  {\color{red}$ \pentagon$}-problem from section \ref{reflection3}  by adding to the structure of  the manifold $V$ an action of a compact group $G_\partial$  on its boundary, where  the metric $g$ in this problem must be required to be  $G_\partial$-invariant.

In fact, the persistence  of $\mathbb T^\rtimes$-stabilization also  suggests an  addition  of an action of  a  compact group $G$ on $V$ with requirement of $g$ being $G$-invariant as well.}} \vspace {1mm}

{\it Then $X$ contains a surface $\Sigma$ as in the above theorem. In fact, there exists a smooth compact  connected oriented  surface  $S\subset X$  with $\partial_{side}(X)$  which represent a {\sf non-zero} homology class in $H_2(X, \partial_{side}(X)$ and such that 
 $$\int_S  Sc(X,s)ds+ 2\int_{ \Theta} mean.curv(\partial_{side}, \theta)d \theta \leq 4\pi\chi( \underline S).$$}

\vspace {1mm}

{\it About the Proof.}  Develop $X$ by reflections in the faces, divide the resulting manifold $\tilde X$ (diffeomorphic to $S\times \mathbb R^2\times \mathbb R^{n-4}$) by a non-torsion subgroup  $\Gamma_0$ of finite index in the reflection group $\Gamma$ (that isometrically acts on 
$\tilde X$ with $\tilde X/\Gamma=X$)  and smooth the (natural continuous) Riemannian 
 metric on $X/\Gamma_0$ with almost no decrease  of  its scalar curvature. 

This  reduce the problem to the case, where $\hexagon$ is replaced by a closed surface of positive genus and where  G-B-inequality applies.   

{\it Exercises} {\sf  (a) {\sf Fill in the details in this argument.}}
   
%{\it Hint} See section ??? and references therein.

(b) {\sf Extend the proof to the case of higher dimensional "reflection polyhedra"instead of $\hexagon$,e.g.
for  $m$-cubes $[0,1]^m$.}

(c) Apply (b) to $\underline X=\underline S\times  [0,1]^{n-2}$  and work out yet another criterion for 
$Sc\geq 0$  additionally to these in    section \ref{Sc-criteria3}.

(d) Formulate and proof the hyperbolic version (i.e. for $Sc\geq \sigma<0$) of this  criterion in the spirit $(2_{\leq 0})$  in  section \ref{reflection3}.

(e)  Formulate and proof the  version of this for $Sc\geq \sigma<0$  by taking into account the area of $S\subset X$ .}

%%%%%%%%%%%%%%%%%%%

\subsection {\color {blue} Width, Waist and other Slicing Invariants}\label{width7}

%%%%%%%%%%%%%%%%%%%%%

Given  numerical  invariant ${\color {blue}INV}$ of $k$-dimensional spaces  $Y$,  one  defines the "slicing version " of ${\color {blue}INV}$  for $n$-dimensional $X$, $n\geq k$, as the infimum of the numbers ${\color {blue}I}$, such that 
$X$ can be  "sliced" into $k$-dimensional subspaces $Y=Y_{\underline x} \subset X$, parametrized by an $(n-k)$-dimensional space $\underline X\ni \underline x$, such that  ${\color {blue}INV\leq I}.$
\vspace{1mm} 

{\it \textbf {Example}} 1:  {\it Uryson's Width.}  If  ${\color {blue}INV}$ stands for "diameter" then the corresponding slicing invariant of a, say locally compact metric space $X$, called   $k$-width is defined, via slicings of  $X$   by, where   pullbacks of points under continuous
maps $f: X\to\underline X$  for polyhedral (triangulated) spaces  $\underline X$ of dimensions  
$dim(\underline X)=m=  dim(X)-k$..

If the dimension of $X$ is unspecified of if $X$ is infinite dimensional, we speak of {\it codimension $m$}  width.)

{\it Exercises}.\footnote {I haven't done these  exercises.}
 (a) Evaluate the widths of balls, ellipsoids simplices and  rectangular solids in Euclidean spaces.

(b) Decide whether the $k$-width is (essentially) non-increasing under proper 1-Lipschitz maps of non-zero degrees between Riemannian $n$-manifolds  for all  $k\leq n$: the existence of such a map $X_1\to X_2$  should(?)  imply that 
$width_k( X_2)\leq  width_k( X_1)$, or at least,  that $width_k( X_2)\leq  const_n\cdot width_k( X_1)$.

\hspace{1mm}

Below, as a matter of  instance,   we formulate a quantified  version of  the  classical  bound on Lebesgue covering dimension by the Hausdorff dimension
  conjectured in  [Guth(volumes of balls-width) 2011],    proved  in  [Lio-Li-Na-Ro(filling) 2019] and  refined in  [Papasoglu(width) 2019],  where also  a  (relatively)   direct proof was found.

\vspace {1mm}

\textbf {Theorem A.}  {\sf There exists a universal constant  $\epsilon = \epsilon_n>0$, such that all {\it proper}}
 (closed bounded subsets are compact)   {\sf  metric spaces  $X$ admit the following bound on the codimension $ n-1$   Uryson 
 width.}

{\sl Let, for some  $R=R_X>0$,  all   pairs of concentric balls  $R$-balls, 
$$B_x(R)\subset B_x(10R)\subset X,\mbox   { } x\in X,$$   
admit  closed subsets $S$ pinched between the boundaries of these balls,
$$S\subset B_x(10R)\setminus B_x(R),$$ 
  such that 

$\bullet_1$  $S$ separates  the ball $B_x (R)$ from the complement $X\setminus B_x(10R)$, i.e. no connected component   of this complement intersects both $B_x(R)$   and $X\setminus B_x(10R)$;

$\bullet_2$   $S$ can be covered  by  countably many balls 
$$S\subset \bigcup_i  B_{x_i}(r_i),$$
such that 
$$\sum_i r_i^{n-1}\leq \epsilon R.$$}

{\it Then $X$ admits a continuous map into an $(n-1)$-dimensional  polyhedral space,
$f:X\to \underline X,$
such that 
$$diam (f^{-1}(\underline x))\leq R, \mbox  { for all } \underline x \in \underline X.$$}

(A significant instance  is\\of this, proven by Guth, is that  of Riemannian $n$-manifolds  $X$, where the inequality  
$vol_x(B(1))\leq\epsilon_n$ for sufficiently  small $\epsilon=\epsilon_n>0$, implies that $width_1(X)\leq const_n\varepsilon.$ 

For example, 

{\it all  Riemannian $n$-manifolds $X$ satisfy: $width_1(X) \leq const_n vol(X)^\frac {1}{n}$.})

\vspace {1mm}

Another  basic property of Uryson width, now  in relation  to curvature, is the following.\vspace {1mm}

\textbf {Theorem B.}  [Perelman(width) 1995] {\it The the volumes  of all Riemannian $n$-manifolds (and singular Alexandrov spaces)  
$X$ with non-negative sectional curvatures are bounded by their Uryson width (essentially) the same way   as it  is for  rectangular solids 
 $$\frac {1}{const_n}\prod_{k=1}^n width_k(X)\leq vol(X)\leq const_n\prod_{k=1}^n width_k(X).$$}

({\it \color {red!50!black} Probably}, there are   similar bounds for the waists of these manifolds:
\hspace {-0mm}$$\frac {1}{const_n}\prod_{k=1}^l width_k(X)\leq waist_l(X)\leq const_n\prod_{k=1}^lwidth_k(X),
 l=1,...,n.)$$

\vspace {1mm}

{\it \textbf {Example}} 2:  {\it From Volumes to  Waists.}   If  ${\color {blue}INV}$ represents  the $k$-volume of $k$-dimensional Riemannian manifolds, 
 then the corresponding slicing invariant of  Riemannian $n$-manifolds is called  the  {\it $k$-wast}, denoted $waist_k(X)$, which, in  the simplest case, can be defined  with slicings of $X$ by pullbacks of points under 
 continuous maps $f:X\to\underline X=\mathbb R^{n-k}$ with $vol_k(f^{-1}(\underline x))$  understood as {\it $k$-dimensional Hausdorff measure.}
 
 It is known  (see section 1.3 in [G(singularities)  2009]) that all Riemannian $n$-manifolds have strictly positive $k$-waists for $k\leq n$:
 
 {\it Every continuous map $f: X\to \mathbb R^{n-k}$ admits a point $\underline x\in \mathbb R^{n-k}$, such that
 $$Hau_k (f^{-1}(\underline x))\geq \delta=\delta_X>0.$$}
 However,   (non-trivial) 
  sharp bounds on the  waist, such as $waist_k(S^n)=vol_k(S^k)$  for unit spheres,
  have been proved only under annoying, {\color {red!40!black} probably unnecessary}, assumptions  on $f$, such, e.g. as  being smooth {\it generic} or {\it piece-wise  real analytic}.\footnote{ See [G(filling) 1983], [G(waist) 2003], [Guth (waist) 2014], [Akopyan-Karasev( tight estimate) 2016], [Akopyan-Karasev(non-radial Gaussian) 2019],  [Klartag(waists) 2017].}\vspace {1mm}

{\it $\mathbb Z_2$-Waist and the  Even Degree Problem.} The known  proof of the  lower bounds on waists 
of manifolds  $X$, such as rectangular solids, for example, which depend on a Borsuk-Ulam topological lemma,    apply to the {\it $\mathbb Z_2$-waists}
defined in terms of the {\it Morse spectrum of the $k$-volume function} on the space of $\mathbb Z_2$-cycles of  dimension $k$ (see [Guth(Steenrod) 2007],  [G(Morse Spectra) 2017]), where this waist is monotone decreasing under  smooth  maps $f:X_1\to X_2$ of {\it odd degree}:  

{\sl if the map $f$ is k-volume non-increasing, $||\wedge^k df||\leq 1$, and $deg(f)$ is {\it odd}, then  
$\mathbb Z_2$-$waist_k(X_2)\leq  
\mathbb Z_2$-$waist_k(X_1)$.}

But it is {\it \color {red!50!black!}unclear}  if such monotonicity holds for all k-volume non-increasing maps with { \it \color {red!50!black!} non-zero} {\it degree.}

\vspace {1mm}

{\it Almgen's Min-Max Theorem.} There is an  alternative proof of the sharp lower bound on $waist_k(S^n)$, that  
relies on Almgen's min-max theorem, which  delivers minimal subvarifolds of volume $\leq v$ in Riemannian manifolds sliced 
into cycles of volumes $\leq v$ (see  [Guth (waist) 2014]). 

this Although  proof, doesn't (seem to) apply to rectangular solids,  it  does  yield  

{\sf sharp  lower bounds  for the $k$-waists of   compact manifolds with {\it sectional curvatures $\geq \kappa>0$} }(see {G(singularities)  2009]). 

However, the following remains unsettled.

\textbf {Problems.}  A.  {\sf Extend Almgen's method to {\it singular} Alexandrov spaces with $sect.curv\geq\kappa$.}

B.  {\sf Develop  a unified  method that would yield, for instance,  sharp inequalities for products of spaces $X_i$
with $sect.curv(X_i)\geq \kappa_i>0$.}\vspace {1mm}
 
 {\it Spherical Waists with the the Dirac operator.} The   sharp parametric area contraction theorem  from section \ref{parametric3}  implies the sharp lower bound
 on the {\it spherical waists} of $N$-spheres:
 
{\sf the space of smooth  {\it strictly area decreasing}  maps $f:S^2\to S^{N}$ is contractible  in the space of all continuous maps $S^2\to S^{N}$  for all $N\geq 2$.}

Moreover,

{\sf Let $\underline X=(\underline X\underline g)$ be a compact Riemannian $N$-manifold with {\it positive curvature operator}, e.g. a convex  
hypersurface in $\mathbb R^{N+1}$  and  let   $\underline g_\circ= 
\underline g_\circ( \underline x)= \frac {1}{N(N-1)}Sc( \underline X,  \underline x)\underline g_\circ( \underline x).$ }
Then  the argument used in the proof of the sharp parametric area contraction theorem yields the bound on the spherical waist of $\underline X_\circ=(\underline X, \underline g_\circ)$  from below:

{\sf  $Sc(\underline X)>0$, then   the space $\mathcal F_\circ$ of smooth  {\it strictly area decreasing}  maps $f:S^2\to  \underline X_\circ$ is {\it contractible  in the space of all continuous maps} $S^2\to \underline X_\circ$  for all $N\geq 2$.}
\vspace {1mm}

{\it Questions.} (a) Does the space  $\mathcal F_\circ$ is contractible?

(b) Is the ordinary 2-waist of  $\underline  X_\circ$ is similarly  bounded  from below as 
$$waist_2(\underline  X_\circ)\geq 4\pi?$$
 In particular, is the space of  maps  $f: \Sigma\to \underline X_\circ $, where $\Sigma$ is surface of genus>0
and where $area (f(\Sigma))<4\pi$, also contractible in the space of all continous maps      $\Sigma\to \underline X_\circ$?

(c)  Is there a counterpart of the above for $n>2$, e.g.  for maps  $S^n\to (X,g_{?})$, $n>2$, in  the spirit of Almgren's style proof of the lower waist bound for  manifolds with $sect.curv>0$?

%%%%%%%%%%%%%%

\subsection {\color {blue}  Hyperspherical Radii, their Parametric and $k$-Volume  Multi-contracting Versions} \label{hyperspherical7} 

%%%%%%%%%%%%%%%%%%%%%%%%

From a category/homotopy  theoretic point of view  the main role of Riemannian metrics on  manifolds $X$ and $Y$ is a definition of a "norm" on smooth maps $f:  X\to Y$, where we distinguish the following.

$\bullet^k_{sup} $ The sup-norm on the $k$th exterior power of the differential of $f$, denoted 
$$||\wedge^kdf||=\sup_{x\in X} ||\wedge^kdf(x)||^\frac {1}{k}.$$
For instance,  the inequality $||\wedge^kdf||< 1$ means that $f$ strictly decreases the $k$-volumes of smooth
$k$-submanifolds in $X$.

$\bullet^k_{trace} $ The normalized trace norm on $\wedge^kdf(x)$, 
$$||\wedge^kdf||_{trace}=\sup_{x\in X} \frac {1}{\binom{n}{k}}(trace\wedge^kdf(x))^\frac {1}{k},$$
(In terms  of an orthonormal frame $e_1, ...   e_n\in T_x(X)$, for which the vectors $df(e_i)\in T_y(Y)$, $y=f(x)$ are orthogonal
$$trace \wedge^kdf(x)=\sum_{{1\leq i_1}<{i_2}<...<{i_k}\leq n }\lambda_ {i_1}\cdot\lambda_{i_2}\cdot...\cdot\lambda_{i_k}.$$
for $\lambda_i=||df(e_i)||$.)

Such a "norm"  defines a "norm" on homotopy and/or other classes $[f]$ of maps $f$, by
$${\sf "norm"}[f]=\inf_{f\in [f]} {\sf "norm"}(f), $$
where a relevant example is where 
$[f]= [f]_{h,\underline h}$   consists of the maps that send a given homology class  $h\in H_\ast(X)$ to a given set  $\{\underline h\}$ or a set of classes  
$\underline h\in H_\ast(Y)$.

For instance   if $Y=S^n$  and  the set  $\{\underline h\}$  consists of non-zero multiples of the fundamental class $[S^n]\in H_n(S^n)$,  we define various  {\it  hyperspherical radii} of $h$ as the reciprocals of such norms,
 $$Rad _{S^n}^{{\sf norm}}(h)= \frac {1}{{\sf "norm"}[f]},$$ 
where "{\sf norm}"  may stand for $||\wedge^kdf||$  and $||\wedge^kdf||_{trace}$.  

And if  the class $[f]$ consists of the maps $f:X\to S^n$ with {\it non-zero homology homomorphism
$:H_n(X)\to H_n(S^n)=\mathbb Z$,  } 
we write 
$$Rad _{S^n}^{{\sf norm}}(X)= \frac {1}{{\sf "norm"}[f]}=\sup_{0 \neq h\in H_n(X)} Rad _{S^n}^{{\sf norm}}(h).$$

In particular, if $X$ is  a connected   orientable $n$-manifold  and $[f]$  is the class of {\it locally constant at infinity maps $f:X\to S^n$ of non-zero degrees}, i.e. which {\it dominate} non-zero-multiples of the fundamental class $[S^n]\in H_n(S^n)$, 
we speak of  {\it hyperspherical radii of $X$},
$$Rad _{S^n}^{{\sf norm}}(X)=Rad _{S^n}^{{\sf norm}}[X]= \frac {1}{{\sf "norm"}[f]},$$ 
with an emphasis on the  norms"  $Lip=||df||$, $||\wedge^kdf||$  and $||\wedge^kdf||_{trace}$, $k=1,2$,
where non-trivial bounds on these radii for  manifolds $X$ with $Sc(X)\geq \sigma>0$,
 are  given (in different terms) in  section  \ref{area extremality3}.

{\it Exercise} (a)  Show that the hyperspherical radius  of   the $R$-sphere $S^n(R)$ defined with any of  the norms   $\bullet^k_{sup} $ and  $\bullet^k_{trace} $ is equal to $R$.

(b) Evaluate these radii for the (open) Euclidean  ball $B^n(R)$  and the cube $[0,R]^n$.\footnote{I haven't done this exercise for the cube.}

(c) Show that if $X$ is a (Riemannian product)  cylinder, $X=X_0\times\mathbb R^1$,  and $[f]_{\underline h}$ is the class of maps $f:X\to Y$, which  send the fundamental class of $X$ to $\underline h\in H_n(Y)$, $n=dim(X),$
then the  "norms"  of the multiples of   $\underline h$ are bounded by the corresponding "norms" of  $\underline h$,
$${\sf "norm"}(j\cdot \underline h)\leq  {\sf "norm"}( \underline h),  \mbox { }  j=0, \pm 1, \pm 2,...\mbox { } . $$

\vspace {1mm}

{\it Spaces  of Maps and Parametric Radii.}  A norm  on maps $f:X\to Y$ can be regarded as  a function on the space $\cal F$ of maps $X\to Y$ (not only on the set of homotopy classes of maps).  

Call  such a function $\Psi: \mathcal F\to [0,\infty)$  and define a (Morse-kind) filtration on the homology $H_\ast(\mathcal F)$, by the {\it  images of the homology    homomorphisms} induced by the sublevels of $\Psi$ to $\cal F$, 
 $$\mbox{ $H_\ast (\Psi^{-1}[0,\lambda]) \to H_\ast(\mathcal F)$, 
 $0\leq\lambda <\infty$,}$$ 
 where  these images are denoted 
 $$H_\ast(\mathcal F|{\tightunderset {\Psi}\leq \lambda})\subset H_\ast(\mathcal F).$$ 

Equivalently, $\Psi$ defines a function on $H_\ast(\mathcal F)$, call  it 
$$\Psi_\ast:  H_\ast(\mathcal F)\to [0,\infty),$$  
where $\Psi_\ast(h)$ is the infimum of $\lambda|geq 0$ for which $h\in H_\ast(\mathcal F{\tightunderset {\Psi}\leq \lambda})$.

 In other words, 
 
 {\sf the inequality $\Psi_\ast(h)\leq \lambda$  for  $h\in H_i\mathcal F)$ signifies that }

 {\sl $h$ is representable by a  family $P$ of maps 
 $f_p: X\to Y$, $p\in P$,  where $P$ is an oriented $i$-pseudomanifold  and 
 $\Psi(f_p)\leq\lambda$ for all $p\in P$. }\vspace{1mm}

{\it Example: Stabilized Radii.} Let  $X$ be an orientable $n$-manifold and $Y$ be the unit sphere $S^{n+m}$. Then
the homology of the space $\mathcal F$ of continuous maps $f:X\to S^{n+m}$ vanishes for $0<k<m$ and 
$H_m(\mathcal F)=\mathbb Z.$
Define the (stabilized)  spherical radii of $X$, by
$$Rad^{\sf {norm}}_{S^{n+m}}(X)=Rad^{\sf {norm}}_{S^{n+m}}(X\times \mathbb R^m),$$
 and observe that such a  radius is equal to the infimum of the "norms"  of the non-zero classes in $H_m(\mathcal F)$.

\vspace{1mm}

{\it Remark.}  If  the  above  "norm" is associated with $||\wedge^ndf||$, $n=dim(X)$, (where the  maps with
 $"norm"(f)=||\wedge^ndf||\leq 1$ are volume non-increasing),  then, according the the sharp  waist inequality from the previous section, the stabilized radii are equal to the basic one:
 $$Rad^{\wedge^n}_{S^{n+m}}(X)=Rad^{\wedge^n}_{S^{n}} \mbox { for all } m=1,2,...  \mbox { }.  $$

\vspace{1mm}

{\it \color {red!49!black} Stabilization {\color {red!50!black} Conjecture} .} If $k<n$ then  the stabilized radii satisfy:
   {\sf $$Rad^{\wedge^k}_{S^{n}}(X)\geq Rad^{\wedge^k}_{S^{n}} 
   \geq c_{n,m,k}Rad^{\wedge^k}_{S^{k}}(X),$$
  for all $k=1,...,n-1$ and   universal constants $c_{n,m,k},$
such that 
  $$1>c_{n,1,k}>c_{n,2,k}>...>c_{n,m,k}>...\geq c_n>0.$$}

{\it Admission.} I haven't proved that  either $c_{2, 3,1}>0$ or that  $c_{2, 3,1}=1$,  that  is a  possible decrease (if any)  of minimal  Lipschitz constants for maps $X\times \mathbb R^1\to S^3$ with non-zero degrees  compared to such  maps
$X\to S^2$  of  oriented surfaces $X$. 

\vspace{1mm}

{\it Diagrams and Multiple Norms.} All of the above definitions can be generalized by replacing single maps between Riemannian manifolds by diagrams $\mathcal D=f_I$  of maps $f_i$  with  homotopy commutativity relations imposed on some sub-diagrams in  $\cal D$.

We have met simple instances of such diagrams for  distance  and area multi-contracting  
 maps to products,
$$f=(f_1,f_2,...f_k):  \underline  X\to \underline X_1\times  \underline X_2\times ...\times  \underline X_k$$
(see section \ref{multi-contracting3}),
where a  "total norm" of such an $f$ related  to scalar curvature 
is 
$$ \left (\frac {1}{k}\sum_{i=1}^k \frac {1}{("norm"(f_i))^2}\right)^{-\frac {1}{2}}.$$

{\it Problem} {\sf Find constraints on  norms  $Lip(f_i)=||df_i||$ and on $||\wedge^2df_i||$  for  more complicated diagrams  $f_I=\{f_i\}$  of maps between manifolds with $Sc$-normalized and/or 
$\mathbb T^ \rtimes$-stabilized (see section \ref{warped2}) manifolds with positive scalar curvatures.}

%%%%%%%%%%%%%%%%%%

\subsection {\color {blue} $m$-Radii of    Uniformly Contractible Spaces}\label {m-radii of    uniformly contractible7}

%%%%%%%%%%%%%%%%%%%%%%%%%%

Define the $m$-radius  with an above   "norm"  of a Riemannian manifold $X$ as 
the supremum of the hyperspherical radii of all $m$-cycles   in $ X$, or, more  formally as 
$$Rad_m^{norm}(X)=\sup_{V\subset X} Rad^{norm}_{S^m}(V).$$
where the supremum is taken over  all relatively compact open (not to  worry about pathologies)   subsets $V$ in $X$.

\vspace {1mm}

{\it Exercise.} Show that if $X$ is an $n$-dimensional Riemannian manifold with $H_{n-1}(X)=0$, then 
$$Rad^{Lip}_{n-1}(X)\leq const_n Rad^{Lip}_{S^n}$$
for $const_n< 10 n$.\vspace{1mm}

If $X$ is uniformly contractible   (see  section \ref{5D.3})  then -- it is (almost) obvious -- that  $Rad^{Lip}_1(X)=\infty$.  
But it is unclear, in general, if this true for  $Rad^{Lip}_m$, for $m\geq 2$.

 Below is a partial result in this direction slightly generalizing that in  \S$9\frac{3}{11}$ [G(positive) 1996].

  \vspace {1mm}

{\it Lipschitz Suspension Lemma.} 
{\sf Let a Riemannian manifold  $X$ contain a double  sequence of triples of disjoint $(m-1)$-cycles,   i.e. of oriented  
$(m-1)$-dimensional  sub-pseudomanifolds, $A_{ij}, B_{ij}, C_{ij}\subset X,$ and let  
$$\mbox {$[0, A_{ij}],  [A_{ij}, B_{ij}],  [B_{ij}, C_{ij}]\subset X$ 
}$$ 
be  $m$-chains represented by  oriented  $m$-sub-pseudomanifolds with boundaries $-A_{ij}$,   $A_i\cup- B_{ij}$ and $B_{ij}\cup- C_{ij}$.}

Let  

{\sl $\bullet_1$ $R_{ij}=Rad_{S^{m-1}}(A_{ij})\geq R_i\to \infty$  for $i\to \infty$;  

$\bullet_2$ the diameters $diam[0, A_{ij}]\leq d_j$  for some constants $d_j$  and all $i$;

$\bullet_3$ $[A_{ij}, B_{ij}]$ is  contained in the $\delta\cdot  R_i$-neighbourhood of  $B_{i,j}$, i.e. 
$dist (x, A_{ij})$, $ x\in[A_{ij}, B_{ij}]$, is bounded by n  $  \delta\cdot R_i$, for   $\delta \leq \frac {1}{10n}$, $n=dim(X)$;

$\bullet_4$ $dist ( [B_{ij}, C_{ij}], [0, A_{ij}])\geq r_i\to\infty $  for $i\to \infty$.

$\bullet_5$  $dist (C_{i,j},[0,  A_{i,j}]\geq d^+_j\to \infty$ for $j\to \infty$.}

Then,  if  $X$ is  uniformly contractible (or  uniformly rationally acyclic) and $n=dim(X)>m$,
then
$$Rad^{Lip}_m(X)=\infty.$$ }

{\it Proof.}  To keep track of these    $\bullet_1$-$\bullet_5$, visualize  $A_{i,j}$, $B_{ij}$  and $C_{ij}$ as  concentric  circles of radii $i$,  $\frac {41}{40}i$ and $2i+j$ in 
$\mathbb R^2\subset X=\mathbb R^3,$  where the chains $[0, A_{ij}],  [A_{ij}, B_{ij}],  [B_{ij}, C_{ij}]\subset X$  are the annuli between these circles.

Then proceed with the proof by observing that 
 uniform contractibility of $X$ implies that  the  cycle  $C_{ij}$ for $j>>i$  (much greater) bounds  chain,  call it
$[C_{ij}, \varnothing_{ij}]$ with the the support far  from  $[0, A_{i,j}]$,  say 
$$\mbox {$dist([C_{ij}, \varnothing_{ij}], [0, A_{i,j}])\geq 10R_i$}.$$

Then the union $D^m_{ij}$ of these  four chains
$$D^m_{ij} = [0, A_{ij}]\cup  [A_{ij}, B_{ij}]\cup   [B_{ij}, C_{ij}]\cup [C_{ij}, \varnothing_{ij}]$$
makes a  $m$-cycle, such that    
$$Rad^{Lip}_{S^m}(D^m_{ij})\geq \varepsilon R_i,$$ 
say,  for   $ \varepsilon= \frac {1}{1000n^3}.$

This is shown by constructing  a  $\lambda$-Lipschitz 
map $D^m_{ij}\to S^{m}(R_i)$ with non-zero degree and with    $\lambda< 100n^2$,  such that 
$ [0,A_{ij}]\subset D^m_{ij}$ goes to the south pole of $S^{m}(R_i)$ and $  [B_{ij}, C_{ij}]\cup [C_{ij}, \varnothing_{ij}]$
to the north pole and 
 where the two main ingredients of  this    construction   are the following:

(i) a $\lambda_1$-Lipschitz  extension of the 1-Lipshitz map $A_{ij}\to S^{m-1}(R_i)$ to 
the $\delta R_i$-neighbourhood of $A_{ij}\subset X$, for $\delta=\frac{1}{10n}$;

(ii) distance function  $x\to dist(x, A_{ij})$, $x\in X$.

(Fitting all this together  is left to the reader.)

\vspace {1mm}

{\it Large Scale Lipschitz Uniform Embeddings.}  A map between metric spaces, $\phi:Z\to X$ is LSL if 
there exist positive constants $\lambda$ and $e$, such that 
$$dist (\phi(z_1z_), f(z_2))\leq \lambda\cdot dist(z_1,z_2)+e.$$

A map $\phi:Z\to X$  is LSUE if 
there exists a function $\Delta(D)= \Delta_\phi(D)$, such that $\Delta(D)\to \infty$ for $D\to \infty$ and 
$$dist (\phi(z_1), f(z_2))\geq  \Delta(dist(z_1,z_2)).$$ 

A map $f:Y\to X$ is LSLU  embedding if it is LSL  as well as  LSU.

{\it LSLUE-Lemma.} \vspace {1mm} {\sf Let  $X$ and $Z$ be Riemannian manifolds, $A_{ij}, B_{ij}, C_{ij}\subset Z$ and 
$[0, A_{ij}],  [A_{ij}, B_{ij}],  [B_{ij}, C_{ij}]\subset Z$ be  $(m-1)$-cycles and $m$-chains  satisfying the above conditions $\bullet_1$-$\bullet_5$
and let $Z\to X$ be an LSLU  embedding.}

{\it If $X$ is uniformly  contractible  (uniformly rationally acyclic will do) then $X$ also contains  $(m-1)$-cycles and $m$-chains cycles, which satisfy  $\bullet_1$-$\bullet_5$.}

Thus, for instance, 

{[\Large\color {blue} $\star$}]  {\it if $dim(Z)=m$, if
 $Rad^{Lip}_{S^m}(Z)=\infty$,    if   $X$ is  uniformly  contractible and if   $dim(X)>m$, then
$$Rad_m^{Lip}(X)=\infty.$$}

{[\Large\color {blue} $\star\star$}] {\it \textbf {Example of Corollary}.} {\sf Let  $X$  be a  compact aspherical manifold of dimension six.}

{\it If the fundamental group 
$\pi_1(X)$ contains a  surface group $\Gamma$ (e.g. $\Gamma=\mathbb Z^2$)  as a subgroup, 
then $X$ admits no metric with $Sc>0$.}

\vspace {1mm}

{\it Proof.} The inclusion $\Gamma\subset \pi_1(X)$  implies that the universal covering $Z$  of the surface with the fundamental group $\Gamma$ admits  an LSLU  embedding to the universal covering $\tilde X$ of $X$.
Hence, $Rad_2^{Lip}(\tilde X)=\infty$.

On the other hand, an easy argument  (see \S$9\frac{3}{11}$ [G(positive) 1996]  and  [G(aspherical) 2020] shows that if a uniformly contractible $n$-manifold $\tilde X$ ) satisfies $Rad^{Lip}_m=\infty$, then it contains  
compact submanifolds $Y$ of dimension $n-m-1$, which have arbitrarily large filling radii,  while, 
if $Sc(X)\geq \sigma$, then 
 $\mathbb T^{\rtimes}$-stabilizations  $Y_{\rtimes}$ of $Y$  have their scalar curvatures bounded from below 
by $\sigma/2>0$.   

This,  in the present 6d-case, contradicts to the bound $fillrad (Y)\leq const\cdot \sigma$ for $dim(Y)=3$. QED.\vspace {1mm}

{\it Exercise.} Extend all $5d$-results from   section \ref{5D.3}  to   $(5+m-1)$-manifolds  $X$,  which admit maps of non-zero degree to uniformly contractible (and uniformly rationally acyclic) manifolds $\underline X$ (and pseudomanifolds with at most 2-dimensional singularities),  where the fundamental groups $\pi_1(\underline X)$ contains subgroups $\Gamma$,  which  serve as fundamental groups of compact $m$-pseudomanifolds the universal coverings $Z$ of which have $Rad^{Lip}_{S^m}(Z)=\infty$.

%%%%%%%%%%%%%%%%\%%%%%%%%%%
%%%%%%%%%%%%%%%%%%%

\section { References }  

%%%%%%%%%%%%%%%%%
%%%%%%%%%%%%%%%%%%%

\hspace {5mm}[AbVaWa{Clifford  Salingaros Vee)2018]  R. Abramowicz, M. Varahagiri, A. Walley,  {\sl A Classification of Clifford Algebras as Images of Group Algebras of Salingaros Vee Groups.} Adv. Appl. Clifford Algebras 28, 38 (2018). \url{https://doi.org/10.1007/s00006-018-0854-y}

  \vspace{1mm}   \vspace{1mm}

[AH-VPPW  (almost non-negative)  2019] B. Allen, L. Hernandez--Vazquez, 
D. Parise,  A. Payne, S. Wang, {\sl Warped tori with almost non-negative scalar curvature,}
Geom Dedicata (2019) 200:153?171.
\vspace {1mm}\vspace {1mm}

[Akopyan-Karasev( tight estimate) 2016]    A. Akopyan, R. Karasev,    {\sl  A tight estimate for the waist of the ball,}
arXiv:1608.06279 [math.MG]
 	
  \vspace{1mm}   \vspace{1mm}

[Akopyan-Karasev(non-radial Gaussian) 2019]       A. Akopyan, R. Karasev, {\sl   Gromov's waist of non-radial Gaussian measures and radial non-Gaussian measures},  arXiv:1808.07350v4.

  \vspace{1mm}   \vspace{1mm}

[AKP(Alexandrov  spaces) 2017]  S. Alexander, V. Kapovitch, A. Petrunin, {\sl 
Invitation to Alexandrov geometry:
CAT[0] spaces} arXiv:1701.03483v3 [math.DG] 
   \vspace{1mm}   \vspace{1mm}

[AlbGell(Dirac operator on pseudomanifolds) 2017 ]  P. Albin, J. Gell-Redman, {\sf  The index formula for families of Dirac type operators on pseudomanifolds}, 
arXiv:1712.08513v1. 
   
   \vspace{1mm}   \vspace{1mm}

[Allen(Sobolev Stability) 2019] B. Allen,{\sl  Sobolev stability of the PMT and RPI using IMCF}
arXiv:1808.07841 [math.DG]

\vspace{1mm}   \vspace{1mm}

[Allen(conformal to tori) 2020]  B. Allen, {\sl Almost non-negative scalar curvature on Riemannian manifolds conformal to tori,}  	arXiv:2010.06008 
 \vspace {2mm}

[Allen-Sormani(convergence) 2020]  B. Allen, C. Sormani, {\sl Relating Notions of Convergence in Geometric Analysis}, arXiv:1911.04522v3.
\vspace{1mm}   \vspace{1mm}

[Almeida(minimal) 1985] Sebastiao Almeida, Minimal Hypersurfaces of a Positive Scalar Curvature. Math. Z. 190, 73-82 (1985).

\vspace{1mm}   \vspace{1mm}

[AndDahl(asymptotically
locally hyperbolic) 1998]  Lars Andersson and Mattias Dahl, {\sl Scalar curvature rigidity for asymptotically
locally hyperbolic manifolds,}   Ann. Global Anal. Geom. 16
(1998), no. 1, 1-27.
    \vspace{1mm}  \vspace{1mm}

   [AndMinGal(asymptotically hyperbolic) 2007] Lars Andersson, Mingliang Cai, and Gregory J. Galloway, {\sl Rigidity and positivity of mass for asymptotically hyperbolic manifolds}, Ann. Henri Poincar\'e  (2008), no. 1, 1-33.  \vspace{1mm}  \vspace{1mm}
   
    [AS(index) 1971] M.F.Atiyah, I.M.Singer, The index of elliptic operators IV, V, Ann. of Math. 93(1971), pp. 119-149.  \vspace{1mm}  \vspace{1mm}

 [Atiyah(global)  1969]   M. F. Atiyah. Global theory of elliptic s. InProc. of the Int. Symposium onFunctional Analysis, 1969, pages 21-30, Tokyo. University of Tokyo Press.
   
    \vspace{1mm}  \vspace{1mm}

[Atiyah (L$_2$)  1976]  M. F. Atiyah. Elliptic s, discrete groups and von Neumann
algebras. pages 43-72. Asterisque, No. 32{33, 1976.
   
     \vspace{1mm}  \vspace{1mm}

[Atiyah(eigenvalues} 1984] M Atiyah, {\sl 
Eigenvalues of the Dirac operator}
Mathematical Institute. Oxford OXI 3LB. 

   \vspace{1mm}  \vspace{1mm}

   [Atiyah-Schmid(discrete
series)   1977] M. Atiyah and W. Schmid. {\sl A geometric construction of the discrete
series for semisimple Lie groups}. Invent. Math., 42:1, 62, 1977.
   
     \vspace{1mm}  \vspace{1mm}

  [BaDoSo(sewing Riemannian manifolds) 2018] J. Basilio ,J. Dodziuk, C. Sormani, Sewing Riemannian Manifolds with Positive Scalar Curvature, The Journal of Geometric AnalysisDecember 2018, Volume 28, Issue 4, pp 3553-3602.

 \vspace{1mm}  \vspace{1mm}

[Baer-Hanke(local flexibility) 2020]	C. Baer, B. Hanke,  {\sl Local Flexibility for Open Partial Differential Relations,} arXiv:1809.05703. \vspace{1mm}  \vspace{1mm}

 [Ballmann(lectures) 2006] W. Ballmann {\sl Lectures on K\"ahler Manifolds,}
A publication of the European Mathematical Society.     \vspace{1mm}  \vspace{1mm}

    [Bamler(Ricci flow proof) 2016] R. [Bamler(Ricci flow proof), A Ricci flow proof of a result by Gromov on lower bounds for scalar curvature  arXiv:1505.00088v1c  \vspace{1mm}   \vspace{1mm}

[Barboza-Conrado](disks)   E. Barbosa, F. Conrado, {\sl  Disks area-minimizing in mean convex Riemannian n-manifolds},
	: 	arXiv:2006.11788 [math.DG]

   \vspace{1mm}   \vspace{1mm}

   [B\"ar-Bleecker(deformed algebraic) 1999]     C.B\"ar-D. Bleecker
    {\sl The Dirac  Operator and the Scalar Curvature of Continuously Deformed Algebraic Varieties} 	Zeitschrift: Contemp. Math.Verlag: American Mathematical SocietySeiten: 3-24Band: 242.
   \vspace{1mm}   \vspace{1mm}

   [B\"ar-Hanke(bounary) 2021]        Christian B\"ar, Bernhard Hanke, {\sl Boundary conditions for scalar curvature}
arXiv:2012.09127   \vspace{1mm} \vspace{1mm}

[Bartnik(prescribed scalar) 1993]] R.Bartnik,{\sl Quasi-spherical metrics and prescribed 
scalar curvature, J. Differential Geom.} 37 (1993) 31-71.   \vspace{1mm} \vspace{1mm}

 [Bartnik(asymptotically flat) 1986]) R.Bartnik .{\sl The Mass of an Asymptotically Flat Manifold}, 
   \url {http://www.math.jhu.edu/~js/Math646/bartnik.mass.pdf}   \vspace{1mm} \vspace{1mm}

  [Bartnik(mass) 2003]Mass and 3-metrics of Non-negative
Scalar Curvature, 	arXiv:math/0304259

   \vspace{1mm} \vspace{1mm}

[BaSo(sequences) 2019] J. Basilio, C. Sormani, {\sl Sequences of three dimensional manifolds with positive scalar curvature}, 	arXiv:1911.02152

   \vspace{1mm}  \vspace{1mm}

 [BD(totally non-spin) 2015]    D. Bolotov, A. Dranishnikov {\sl On Gromov's conjecture for totally non-spin manifolds}
arXiv:1402.4510
  \vspace{1mm} \vspace{1mm}
[Be-Be-Ma(Ricci flow) 2011] Laurent Bess\`eeres, G\'erard Besson, and Sylvain Maillot. {\sf Ricci flow on open 3-manifolds
and positive scalar curvature}. Geom. Topol., 15(2):927?975, 2011.
\vspace{1mm}  \vspace{1mm}

[Be-Be-Ma-Ma(deforming 3-manifolds)  2017] Laurent Bessi\`eres (IMB), G\'erard Besson (IF), Sylvain Maillot (IMAG), Fernando Coda Marques,{\sl  Deforming 3-manifolds of bounded geometry and uniformly positive scalar curvature,} arXiv:1711.02457 [math.DG]
  	\vspace{1mm}  \vspace{1mm}

[Bengtsson(trapped surfaces) 2011] I. Bengtsso,   {\sl Some Examples of Trapped Surfaces} arXiv:1112.5318

 \vspace{1mm} \vspace{1mm}

 [Bern-Heit(enlargeability-foliations) 2018] M.-T. Bernameur and J. L. Heitsch, {\sl Enlargeability, foliations, and positive scalar curvature.}

Preprint, arXiv: 1703.02684.\vspace{1mm}  \vspace{1mm}

[BoEW(infinite loop spaces) 2014] {\sl Infinite loop spaces and positive scalar curvature,} B. Botvinnik, J.  Ebert, O. Randal-Williams,

 	arXiv:1411.7408 \vspace{1mm}  \vspace{1mm}

[Boileau(lectures)  2005]   M Boileau {\sl Lectures on Cheeger-Gromov Theory of Riemannian manifolds}, Summer School and Conference on Geometry and Topology of 3-Manifolds |, Trieste, Italy 2005

\url{https://pdfs.semanticscholar.org/9db2/2df12b52f2d1ce3f4c3434731b37bb69d4e6.pdf}\vspace{1mm}  \vspace{1mm}

[Bot-Ros)elementary abelian) 2004]  B. Botvinnik and J.Rosenberg {\sl  Positive scalar curvature for manifolds with elementary abelian fundamental group,}  arXiv:1808.06007.
\vspace{1mm}  \vspace{1mm}

[Bray(Penrose inequality) 2009]  H.L Bray, {\sl The Penrose inequality in general relativity and volume comparison theorems involving scalar curvature (thesis)}  arXiv:0902.3241v1  \vspace{1mm}  \vspace{1mm}

[Bray-Stern 2019(harmonic forms)] H. Bray D. Stern, Scalar curvature and harmonic one- forms on three-manifolds with boundary, arXiv:1911.06803v1. \vspace{1mm}  \vspace{1mm}

[Bre-Mar-Nev(hemisphere) 2011] S. Brendle, F. Marques, and A. Neves, {\sl  Deformations of the hemisphere that increase
scalar curvature}, Invent. Math. 185, 175?197 (2011)
\vspace{1mm}  \vspace{1mm}

[Brendle(rigidity2010] S. Brendle {\sl Rigidity phenomena involving scalar curvature}, arXiv:1008.3097v3
\vspace{1mm}  \vspace{1mm}

 [Brendle-Marques(balls in $S^n$( 2011] S. Brendle, F. Marques, {\sl Scalar curvature rigidity of geodesic balls in  $S^n$}, J. Diff.
Geom. 88, 379?394, (2011).\vspace{1mm}  \vspace{1mm}

[Breuillard-Gelander(non-amenable) 2005] E. Breuillard, T. Gelander, {\sl Cheeger constant and algebraic entropy of linear groups}
arXiv:math/0507441

\vspace{1mm}  \vspace{1mm}

[Brewin(ADM) 2006] Leo Brewin, {\sl A simple expression for the ADM mass.}
arXiv:gr-qc/0609079

[Brun-Han(large and small) 2009] M. Brunnbauer, B. Hanke, {\sl Large and small group homology}, J.Topology 3 (2010) 463-486. 

 \vspace{1mm}  \vspace{1mm}
 [Buhovsky-Opshtein($C^0$-symplectic)2014]  Lev Buhovsky, Emmanuel Opshtein {\sl   Some quantitative results in $C^0$-symplectic 
 geometry} arXiv:1404.0875v2.
\vspace{1mm}  \vspace{1mm}

[Bu-Hu-Sey($C^0$ counterexample)  2016]
L. Buhovsky, V. Humili\`ere, S. Seyfaddini,  {\sl  A $C^0$-counterexample to the Arnold conjecture},
arXiv:1609.09192 [math.SG]\vspace{1mm}  \vspace{1mm}

[Bu-Hu-Sey($C^0$ symplectic)  2020]
L. Buhovsky, V. Humili\`ere, S. Seyfaddini,
The action spectrum and $C^0$-symplectic topology arXiv:1808.09790v2 [math.SG].
  \vspace{1mm}   \vspace{1mm}

[Bunke(relative index) 1992]  Ulrich Bunke  {\sl Relative index theory}. J. Funct. Anal.,105 (1992), 63-76

\vspace{1mm}  \vspace{1mm}

[Bunke(orbifolds) 2007] U.Bunke  {\sl Orbifold index and equivariant K-homology},
Ulrich Bunke	arXiv:math/0701768 [math.KT]
 	\vspace{1mm}  \vspace{1mm}

[Burkhart-Guim(regularizing Ricci flow)   2019]  P. Burkhardt-Guim,  {\sl Pointwise lower scalar curvature bounds for $C^0$-metrics via regularizing Ricci flow}, arXiv:1907.13116 [math.DG]
   
     \vspace{1mm}   \vspace{1mm}

[BurTop(curvature bounded  above)1973] Yu.  Burago and V.  Toponogov, {\sl On 3-dimensional Riemannian spaces with curvature bounded
above}. Math. Zametki 13 (1973), 881-887.  \vspace{1mm}  \vspace{1mm}

[Calegari(lectures on minimal)2019] D. Galegary, {\sl Lectures on Minimal Surfaces} 
\url {https://math.uchicago.edu/~dannyc/courses/minimal_surfaces_2014/minimal_surfaces_notes.pdf}
     \vspace{1mm}   \vspace{1mm}

 \vspace{1mm}  \vspace{1mm}

 [Cecchini(Callias) 2018] S. Cecchini, {\sl   
Callias-type operators in C*-algebras and positive scalar curvature on noncompact manifolds}, Journal of Topology and Analysis on Line.

\url{ https://doi.org/10.1142/S1793525319500687}

 \vspace{1mm}  \vspace{1mm}

[Cecchini(long neck) 2020] S. Cecchini, {\sl A long neck principle for Riemannian spin manifolds with positive scalar curvature}
arXiv:2002.07131v1

\vspace{1mm}  \vspace{1mm}

[Cecchini-Zeidler(generalized Callias) 2021]  S. Cecchini,  R. Zeidler  {\sl Scalar curvature and generalized Callias s}
\vspace {2mm}

[Cecchini-Zeidler(Scalar\&mean) 2021]  S. Cecchini,  R. Zeidler {\sl  Scalar and mean curvature comparison
via the Dirac operator},
arXiv:2103.06833

[CGM(Lipschitz control) 1993] A. Connes, M. Gromov, M, H,  Moscovici, H. {\sl Group cohomology with Lipschitz control and higher signatures}. Geometric and Functional Analysis, 3(1), 1-78. (1993).
 \vspace{1mm}  \vspace{1mm}

 [Cheeger(singular) 1983] J. Cheeger,  {\sl Spectral geometry of singular Riemannian spaces.} J. Differential Geom. 18 (1983), no. 4, 575-657.\vspace{1mm}  \vspace{1mm}

 [Chodosh-Li(bubbles) 2020]  O. Chodosh, C. Li {\sl Generalized soap bubbles and the topology of manifolds with positive scalar curvature},  arXiv:2008.11888v3\vspace{1mm}  \vspace{1mm}

{\color {red} new}[Chodosh-Li-Liokumovich
(classification) 2021]  Otis Chodosh, Chao Li, Yevgeny Liokumovich,  {\sl Classifying sufficiently connected PSC manifolds in 4 and 5 dimensions}, arXiv:2105.07306. \vspace {2mm}

[Chrusciel-Delay(hyperbolic positive energy) 2019] P. Chru\`spiel,  Delay, {\sl The hyperbolic positive energy theorem}
arXiv:1901.05263\vspace{1mm}  \vspace{1mm}

[Chrusciel-Herzlich [asymptotically hyperbolic)  2003] P. Chru\`spiel, M. Herzlich  {\sl The mass of asymptotically hyperbolic Riemannian manifolds}, arXiv:math/0110035v2
\vspace{1mm}  \vspace{1mm}

[Connes(survey of foliations) 1983] A. Connes   {\sl A survey of foliations and  algebras.}
 
\url {Con  www.alainconnes.org/docs/foliationsfine.ps}

\vspace{1mm}  \vspace{1mm}

  [Connes(cyclic cohomology-foliation) 1986]  A. Connes, {\sl Cyclic cohomology and the transverse fundamental class of a foliation}. Geometric methods in  algebras (Kyoto, 1983), 52-144, PitmanRes. Notes Math. Ser., 123,Longman Sci. Tech., Harlow, 1986
  
  \vspace{1mm}  \vspace{1mm}
 
 [Connes(book) 1994]  A. Connes,{
 \sl  Non-commutative geometry} Academic Press.  
 \url{ https://www.alainconnes.org/docs/book94bigpdf.pdf}
 
 \vspace{1mm}  \vspace{1mm}

 [Connes-Moscovici($L_2-index$ for homogeneous)  1982] A. Connes and H. Moscovici. {\sl  The L2-index theorem for homogeneous
spaces of Lie groups}. Ann. of Math. (2), 115(2):291-330, 1982. \vspace{1mm}  \vspace{1mm}

 [Darmos(Gravitation einsteinienne) 1927]    Georges Darmois { Les \'equations de la Gravitation einsteinienne (M\'emorial des Sciences math\'ematiques dirig\'e par Henri Villat; fasc. XXV).Edit\'e par Gauthier-Villars 1927 
 \vspace{1mm}  \vspace{1mm}

   [Davaux(spectrum) 2003] H\'el\`ene Davaux, {\sl An optimal inequality between scalar curvature and spectrum of the Laplacian}  Mathematische Annalen,
 Volume 327, Issue 2, pp 271-292  (2003)   \vspace{1mm}  \vspace{1mm}

 [Davis(orbifolds) 2008] M. Davis, {\sl Lectures on orbifolds and reflection groups.} 
 
 \hspace{-6mm} \url {https://math.osu.edu/sites/math.osu.edu/files/08-05-MRI-preprint.pdf} \vspace{1mm}
    \vspace{1mm}

[Debarre(lectures) 2003]  O. Debarre,  {\sl Fano varieties},

\url{https://webusers.imj-prg.fr/uploads/olivier.debarre/budapest.pdf}\vspace{1mm}  \vspace{1mm}

 [DFW(flexible) 2003] A. Dranishnikov, S. Ferry, S. Weinberger, 
  {\sl  Large Riemannian manifolds which are flexible}.	Ann.Math  157(3), Pages 919-938.
	  \vspace{1mm}  \vspace{1mm}

[Diek-Kantowski (Clifford History)1995]  A. Diek,  R. Kantowski, {\sl  Some Clifford Algebra History} \url {file:///Users/misha/Downloads/Diek-Kantowski1995_Chapter_SomeCliffordAlgebraHistory.pdf} 
\vspace{1mm}  \vspace{1mm}
 
  [Dra-Kee-Usp( Higson corona) 1998]  A. Dranishnikov, J.Keesling, V.Uspenskij. {\sl On the Higson corona of uniformly
 contractible spaces}, Topology 37 (1998), 791-803. \vspace{1mm}  \vspace{1mm}

  [Dranishnikov(asymptotic) 2000]     A.  Dranishnikov, {\sl Asymptotic topology}, Uspekhi Mat. Nauk, 2000, Volume 55
 
\vspace{1mm}  \vspace{1mm}

 [Dranishnikov( hypereuclidean) 2006] {\sl On Hypereuclidean Manifolds}  Geometriae Dedicata (2006)
  117: 215-231 \vspace{1mm}  \vspace{1mm}

   [Dranishnikov (large scale) 1999]   A.Dranishnikov, {\sl On large scale properties of manifolds}, \url{https://arxiv.org/pdf/math/9912062.pdf}\vspace{1mm}  \vspace{1mm}

 [Dranishnikov(macroscopic) 2010] A. Dranishnikov {\sl On macroscopic dimension of rationally essential manifolds}
arXiv:1005.0424

  \vspace{1mm}  \vspace{1mm}
  
  [Dra-Kee-Usp( Higson corona) 1998]  A. Dranishnikov, J.Keesling, V.Uspenskij. {\sl On the Higson corona of uniformly
 contractible spaces}, Topology 37 (1998), 791-803

    \vspace{1mm}  \vspace{1mm}

 [Ebert-Williams(infinite loop spaces) 2017] Johannes Ebert, Oscar Randal-Williams, {\sl 
Infinite loop spaces and positive scalar curvature in the presence of a fundamental group.}

arXiv:1711.11363v1.
  \vspace{1mm} \vspace{1mm}

   [Ebert-Williams(cobordism category) 2019] Johannes Ebert, Oscar Randal-Williams,
   {\sl The positive scalar curvature cobordism category} arXiv:1904.12951v1

    \vspace{1mm}
   \vspace{1mm}
 
[EM(wrinkling) 1998] Ya.   Eliashberg, N. Mishachev {\sl Wrinkling of smooth mappings III. Foliations of codimension greater  than one.}
Topological Methods in Nonlinear Analysis Journal of the Juliusz Schauder Center Volume 11, 1998, 321-350. 
 
   \vspace{1mm}  \vspace{1mm}

 [EMW(boundary) 2009] Eichmair, P. Miao and X. Wang  {\sl Boundary effect on compact manifolds with nonnegative scalar curvature - a generalization of a theorem of Shi and Tam}.  Calc. Var. Partial Differential Equations., 43 (1-2): 45-56, 2012. [arXiv:0911.0377] 
  
      \vspace{1mm}  \vspace{1mm}

 [Entov(Hofer metric) 2001]   M. Entov,  {\sl K-area, Hofer metric and geometry of conjugacy classes in Lie groups}, Invent.Math., 146 (2001), pp. 93-141.

   \vspace{1mm}  \vspace{1mm}

[Federer(singular) 1970] H. Federer,  {\sl The singular sets of area minimizing rectifiable currents with
codimension one and of area minimizing flat chains modulo two with arbitrary
codimension}, Bull. Amer. Math. Soc, 76 (1970), 767-771.   \vspace{1mm}
    \vspace{1mm}

 [Fisher-C\&S(stable minimal) 1980] D. Fischer-Colbrie, R. Schoen The structure of complete stable
minimal surfaces in 3-manifolds of non-negative scalar curvature, Comm. Pure
Appl. Math., 33 (1980) 199-211.   \vspace{1mm}  \vspace{1mm}

  [Frigero(Bounded Cohomology) 2016]  Roberto Frigerio. {\sl Bounded Cohomology of Discrete Groups}, 
  arXiv:1610.08339.

   \vspace{1mm}
    \vspace{1mm}

 [Futaki(scalar-flat) 1993] A. Futaki {\sl Scalar-flat closed manifolds not admitting positive
scalar curvature metrics},
   Invent. math. 112, 23-29 (1993)

  [F-W(zero-in-the-spectrum) 1999]  M. Farber, S. Weinberger, {\sl On the zero-in-the-spectrum conjecture.}
	arXiv:math/9911077 
  \vspace{1mm}
    \vspace{1mm}

   [G(filling) 1983])  [\sl Filling Riemannian manifolds}, J. Differential Geom. 18 (1983), no. 1, 1-147.     \vspace{1mm}  \vspace{1mm}
   
   [G(infinite) 1983]) M. Gromov {\sl Infinite groups as geometric objects}, ICM,  Warsaw, 1984, pp  385-392. \vspace{1mm} \vspace{1mm}

    [G(foliated) 1991])  M. Gromov, {\sl The foliated plateau problem}, Part I: Minimal varieties, Geometric and Functional Analysis (GAFA) 1:1 (1991), 14-79.   \vspace{1mm}  \vspace{1mm}

    [G(large)  1986]  M. Gromov,  {\sl Large Riemannian manifolds.}  In: Shiohama K., Sakai T., Sunada T. (eds) Curvature and Topology of Riemannian Manifolds. Lecture Notes in Mathematics, vol 1201 (1986), 108-122.  \vspace{1mm}  \vspace{1mm}
   
   [G(positive) 1996] M. Gromov, {\sl Positive curvature, macroscopic dimension, spectral
gaps and higher signatures.} In Functional analysis on the eve of the 21st century,
Vol. II (New Brunswick, NJ, 1993) ,volume 132 of Progr. Math., pages 1-213,
Birkh\"auser, 1996.
   
   \vspace{1mm}  \vspace{1mm}
   
[G(waists) 2003] M. Gromov, {\sl Isoperimetry of waists and concentration of maps},  GAFA, Geom. funct. anal.,13 (2003), 178-215.

  \vspace{1mm}  \vspace{1mm}

   [G(singularities)  2009]  {\sl Singularities, expanders and topology of maps. Part 1 : Homology versus volume in the spaces of cycles}, GAFA,  19 (2009) n 3, 743-841.
     \vspace{1mm}  \vspace{1mm}

   [G(Hilbert)  2012] M. Gromov, {\sl Hilbert volume in metric spaces. Part 1}, Open Mathematics
(formerly Central European Journal of Mathematics)

\url {https://doi.org/10.2478/s11533-011-0143-7}

    \vspace{1mm}  \vspace{1mm}

[G(Plateau-Stein) 2014]  M. Gromov, {\sf Plateau-Stein manifolds,} Central European Journal of Mathematics,
 Volume 12, Issue 7, pp 923-95.    \vspace{1mm}  \vspace{1mm}
  
   [G(billiards) 2014] M. Gromov, {\sl Dirac and Plateau billiards in domains with corners}, Central European Journal of Mathematics, Volume 12, Issue 8,   2014, pp 1109-1156. \vspace{1mm}
    \vspace{1mm}

 [G (hyperbolic) 2016] M. Gromov, {\sl   Hyperbolic dynamics, Markov partitions and
Symbolic Categories, Chapters 1 and 2.}
\url{https://www.ihes.fr/~gromov/wp-content/uploads/2018/08/SymbolicDynamicalCategories.pdf}\vspace{1mm}
    \vspace{1mm}

   [G(101) 2017] M. Gromov 101 Questions, Problems and {\color {red!50!black} Conjecture}s around Scalar Curvature. (Incomplete and Unedited Version)
   
\url {https://www.ihes.fr/~gromov/wp-content/uploads/2018/08/101-problemsOct1-2017.pdf}\vspace{1mm}  \vspace{1mm}

 [G(Morse Spectra) 2017]  M.Gromov {\sl Morse Spectra, Homology Measures and Parametric Packing Problems,}
	arXiv:1710.03616 [math.MG]
  	
    \vspace{1mm}  \vspace{1mm}

   [G(inequalities) 2018] 
{\sl Metric Inequalities with Scalar Curvature} Geometric and Functional Analysis Volume 28, Issue 3, pp 645-726.       \vspace{1mm}  \vspace{1mm}
  
     [G(boundary) 2019]  M. Gromov {\sl Scalar Curvature of Manifolds with Boundaries: Natural Questions and  Artificial Constructions.} 

https://arxiv.org/pdf/1811.04311

   \vspace{1mm}  \vspace{1mm}
  
  [G(mean) 2019] M.Gromov {\sl Mean Curvature in the Light of Scalar Curvature}, arXiv:1812.09731v2
   \vspace{1mm} \vspace{1mm}

  [G(aspherical) 2020]     M.Gromov   {\sl  No metrics with Positive Scalar Curvatures on Aspherical 5-Manifolds,}   	arXiv:2009.05332 c\vspace{1mm} \vspace{1mm}

{\color {red} new} [Gadgil-Seshadri(isotropic)2008] S. Gadgil, H. Seshadri, {\sl On the topology of manifolds with positive isotropic curvature}, arXiv:0801.2221 \vspace {2mm}

  [Ge(linking) 2021]  J. Ge,{\sl
  Gehring Link Problem, Focal Radius and Over-torical width}, arXiv:2102.0590
 
 \vspace{1mm} \vspace{1mm}

  [Gerhard(capillarity) 1976] C. Gerhard, {\sl Global regularity of the solutions to the capillarity problem,}
Annali della Scuola Normale Superiore di Pisa, Classe di Scienze 4e s\'erie, tome 3, no 1,
(1976), p. 157-175.    \vspace {2mm}
  
 [Geroch (extraction) 1973] R Geroch,  {\sl Energy Extraction}. Publication: Sixth Texas Symposium on Relativistic Astrophysics. 1973. {\url{https://nyaspubs.onlinelibrary.wiley.com/doi/abs/10.1111/j.1749-6632.1973.tb41445.x}

  [Geroch(relativity)  1975] R. Geroch {\sl General Relativity} Proc. of Symp. in Pure Math., 27, Amer. Math. Soc., 1975, pp.401-414. 
    \vspace{1mm}  \vspace{1mm}

  [Ginoux(3d clasification) 2013]  N.Ginoux.{\sl  The classification of 3-manifolds admitting positive scslr curvature}, 
  \url{ http://www.iecl.univ-lorraine.fr/~Nicolas.Ginoux/classif3dimpsc.pdf} \vspace{1mm} \vspace{1mm}

   [GL(classification) 1980] M.Gromov, B Lawson,     "The classification of simply connected manifolds of positive scalar curvature" Ann. of Math. , 111 (1980) pp. 423-434. \vspace{1mm}  \vspace{1mm}

  [GL(spin) 1980] M.Gromov, B Lawson,  {\sl Spin and Scalar Curvature in the Presence of a Fundamental Group I}
Annals of Mathematics, 111 (1980), 209-230. 
    \vspace{1mm}

  \vspace{1mm}
  
   [GL(complete) 1983]  M.Gromov, B Lawson,  
{\sl Positive scalar curvature and the Dirac  operator on complete Riemannian
manifolds,} Inst. Hautes Etudes Sci. Publ. Math.58 (1983), 83-196.  \vspace{1mm}\vspace{1mm}

 [Goette(alternating torsion)2007] S. Goette,
  {\sl Scalar Curvature Estimates by Parallel Alternating Torsion}, arXiv:0709.4586\vspace{1mm}\vspace{1mm}

 [Goette-Semmelmann(Hermitian) 1999] S. Goette, U. Semmelmann, \{sl  Spin$^c$ Structures and Scalar Curvature Estimates},
arXiv:math/9905089

 \vspace{1mm}\vspace{1mm}

 [Goette-Semmelmann(symmetric) 2002]   S. Goette and U. Semmelmann, {\sl  Scalar curvature estimates for compact symmetric spaces.}
Differential Geom. Appl.  16(1):65-78, 2002.  \vspace{1mm}\vspace{1mm}

  [GSh(Riemann-Roch) 1993] M. Gromov M. Shubin, {\sf The Riemann-Roch theorem for elliptic s}, I. M. Gelfand Seminar, Adv. Soviet Math., vol. 16, Amer. Math. Soc., Providence, RI, 1993, pp. 211-241. MR 1237831 \vspace{1mm}  \vspace{1mm}

 [Guijarro-Wilhelm(focal radius) 2017]  
 L. Guijarro, F. Wilhelm, {\sl Focal Radius, Rigidity, and Lower Curvature Bounds},
arXiv:1603.04050 [math.DG].{\color {red}(new !!!!!!!!)}

 \vspace{1mm}

\vspace{1mm}

  [Guo-Xie-Yu(quantitative K-theory) 2020] H. Guo, Z. Xie, G. Yu
  {\sl Quantitative K-theory, positive scalar curvature, and band width},
	arXiv:2010.01749
\vspace{1mm}

\vspace{1mm}

[Guth(Steenrod) 2007] L. Guth, {\sl  Minimax problems related to cup powers and Steenrod squares}
 	arXiv:math/0702066

 \vspace{1mm}

\vspace{1mm}

 [Guth(metaphors) 2010] L. Guth,  {\sl  Metaphors in systolic geometry}, 	arXiv:1003.4247. 

\vspace{1mm}

\vspace{1mm}

[Guth(volumes of balls-large}) 2011] L.Guth, {\sl Volumes of balls in large Riemannian manifolds}.    Annals of Mathematics173(2011), 51-76.  \vspace{1mm}  \vspace{1mm}
  
  [Guth(volumes of balls-width) 2011]L. Guth,  Volumes of balls in Riemannian manifolds and Uryson width.
Journal of Topology and Analysis,
Vol. 09, No. 02, pp. 195-219 (2017) \vspace{1mm}  \vspace{1mm}

  [Guth (waist) 2014] L Guth,  {\sf The waist inequality in Gromov's work,}  MIT Mathematics
\url {math.mit.edu/~lguth/Exposition/waist.pdf}
  
   \vspace{1mm}  \vspace{1mm}

  [GWY (Novikov) 2019] Sherry Gong, Jianchao Wu, Guoliang Yu
 {\sl The Novikov conjecture, the group of volume preserving diffeomorphisms and Hilbert-Hadamard spaces}
arXiv:1811.02086v3 [math.KT]

      \vspace{1mm}  \vspace{1mm}

   [G-Z(area) 2021] M. Gromov, J. Zhu,  {\sl Area and Gauss-Bonnet inequalities with scalar curvature}, in preparation.
   
     \vspace{1mm}  \vspace{1mm}
   
 [Hanke(Abelian)2020] B. Hanke, {\sl Positive scalar curvature on manifolds with
odd order abelian fundamental groups}  arXiv:1908.00944v3
   \vspace{1mm}  
  \vspace{1mm}

    [HaSchSt(space of metrics)2014]    Bernhard Hanke,Thomas Schick , Wolfgang Steimle, 
{\sl The space of metrics of positive scalar curvature}. Publications math\'ematiques de l'IHES
 Volume 120, Issue 1, pp 335-367 (2014) \vspace{1mm}  \vspace{1mm}  
 
  [Hawking(black holes)  1972] S. Hawking {\sl Black Holes in General Relativity},
    Commun. math. Phys. 25, 152-166 (1972). 
     \vspace{1mm}  
  \vspace{1mm}

[Hess-Schrader-Uhlenbrock(Kato) 1980] H. Hess, R.  Schrader,  A. Uhlenbrock, {\sl Kato's inequality and the spectral
distribution of Laplacians on compact Riemannian manifolds}, J. Diff. Geom.
15 (1980), 27-38.

   \vspace{1mm} \vspace{1mm}
   
    [Herzlich(mass) 2021]    Marc Herzlich {\sl  Scalar curvature, mass  and other asymptotic invariants,} this volume.
     
     \vspace{1mm}  \vspace{1mm}

           [Higson(cobordism invariance) 1991]  N. Higson,  {\sl A note on the cobordism invariance of the index,} Topology 30:3 (1991), 439-443.
\vspace{1mm}\vspace{1mm}

             [Hitchin(spinors)1974]  N. Hitchin, Harmonic spinors, Advances in Math. 14 (1974), 1-55.  \vspace{1mm}  \vspace{1mm}

 [Heier-Wong(uniruled) 2012]  Gordon Heier and Bun Wong,
 {\sl Scalar curvature and uniruledness on projective manifolds},  	 arXiv:1206.2576v1

    \vspace {1mm}  \vspace {1mm}

 [Huang-Jang-Martin(hyperbolic mass rigidity) 2019]  L-H. Huang, H.Jang, D.Martin
 {\sl Mass rigidity for hyperbolic manifolds},arXiv:1904.12010.

     \vspace {1mm}  \vspace {1mm}

    [Hu-Shi(Bartnik cobordism) 1920] X. Hu  Y. Shi, {\sl  NNSC-Cobordism of Bartnik Data in High Dimensions} 	arXiv:2001.05607\vspace {2mm}

    [Ivanov(Lipschitz)  2008]   S. Ivanov, {\sl  Volumes and areas of Lipschitz metrics} Jour. Algebra i Analiz  
    vol 20:3, pp 74-111.

  \vspace {1mm}  \vspace {1mm}

	 [Jang-Miao(hyperbolic mass) 2021] H. Jang, P. Miao {\sl Hyperbolic mass via horospheres}
	arXiv:2102.01036.
	
	\vspace {2mm}

[Jauregui(fill-in)2013]  J. Jauregui, {\sl Fill-ins of nonnegative scalar curvature, static metrics and quasi-local mass},  arXiv:1106.4339v2. \vspace {1mm}  \vspace {1mm}

[Jauregui-Miao-Tam(extensions and fill-ins) 2013]  J. Jauregui, P. Miao, L-F Tam{\sl Extensions and fill-ins with nonnegative scalar curvature},  arXiv:1304.0721

     [JW(exotic) 2008]  M. Joahim, D. J. Wraith, {\sl  Exotic spheres and curvature},  Bull. Amer. Math. Soc.45  no. 4 (2008), 595-616.

  \vspace{1mm}  \vspace{1mm}

[Karami-Zadeh-Sadegh(relative-partitioned) 2018]  M. Karami, M. Zadeh, A. Sadegh  {\sl Relative-partitioned index theorem} arXiv:1411.6090

  \vspace{1mm}  \vspace{1mm}

[Kasparov(index)  1973] G. Kasparov, {\sl The generalized index of elliptic s}, Funkc. Anal i Prilozhen.7no. 3 (1973), 82-83; English translation, Funct. Anal. Appl.7(1973), 238?240. \vspace{1mm}  \vspace{1mm}

[Kasparov(elliptic)   1975] G. Kasparov, {\sl Topological invariants of elliptic s,} I:K-homology,Izv. Akad. Nauk SSSR, Ser. Mat.39(1975), 796-838; English translation,Math. USSR?Izv.9(1975), 751?792. \vspace{1mm}  \vspace{1mm}

[Kasp-Scan (Novikov) 1991] G. Kasparov, G. Skandalis  {\sl Groups Acting on Buildings,  K-Theory,
and Novikov's Conjecture},}
 K-Theory 4: 303-337, 1991. 303
 
 
   \vspace{1mm}\vspace{1mm}
 
 [Katz(systolic geometry)  2017]         M. Katz                        {\sl  Systolic Geometry and Topology, With an Appendix by Jake P. Solomon}. Mathematical. Surveys and. Monographs. Volume 137
 American Mathematical Society, Providence, RI,
 
   \vspace{1mm}\vspace{1mm}
 
  [Kaz(4-manifolds) 2019] D. Kazaras, {\sl 
Desingularizing positive scalar curvature 4-manifolds}
arXiv:1905.05306 [math.DG] \vspace{1mm}  \vspace{1mm}

[Kazdan(complete) 1982] J. Kazdan {\sl Deformation to Positive Scalar Curvature on Complete Manifolds.}  Mathematische Annalen 261 (1982): 227-234.

\vspace{1mm}  \vspace{1mm}

 [Kazdan-Warner(conformal) 1975] by J. Kazdan, F.Warner 
{\sl  Scalar curvature and conformal deformation of Riemannian structure}. J. Differential Geom. 10 (1975), no. 1, 113--134.
  \vspace{1mm}\vspace{1mm}

 [Klartag(waists) 2017] B. Klartag,  {\sl Convex geometry and waist inequalities},
arXiv:1608.04121v2.

 \vspace{2mm}

[Kleiner-Lott(on Perelman's) 2008] B.Kleiner,  J. Lott,{\sl  Notes on Perelman's papers},
Bruce Kleiner, John Lott
arXiv:math/0605667v5 \vspace {2mm}

{\color  {red}( new!!!)}
[Kreck-Su](5-manifolds) 2017] Matthias Kreck, Yang Su, {\sl  On 5-manifolds with free fundamental group and simple boundary links in $S^5$,}  arXiv:1602.01943.  \vspace {2mm}

[Lang-Schroeder(Kirszbraun)1997] U. Lang,  V. Schroeder, {\sc Kirszbraun's Theorem and Metric Spaces of Bounded Curvature} 
GAFA volume 7, pages 535-560.
 \vspace{1mm}  \vspace{1mm}

[Lawson\&Michelsohn(spin geometry) 1989]  B. Lawson, M.-L. Michelsohn, {\sl  Spin geometry.} Princeton University Press 1989.
 
 \url{https://zulfahmed.files.wordpress.com/2014/04/spingeometry.pdf}

 \vspace{2mm}

     [LeBrun(Yamabe)  1999] C. LeBrun,  {\it Kodaira Dimension and the Yamabe Problem}, 
Communications in Analysis and Geometry,Volume7, Number1,133-156 (1999).
  \vspace{1mm}\vspace{1mm}

      [Lesourd-Unger-Yau(arbitrary ends) 2021] M. Lesourd, R.Unger, S-T. Yau. {\sl The Positive Mass Theorem with Arbitrary Ends.} arXiv:2103.02744 [math.DG]\vspace{1mm}\vspace{1mm}
     
 \vspace{1mm}\vspace{1mm}

       [Li(comparison) 2019] C. Li,    {\sl A polyhedron comparison theorem for 3-manifolds with positive scalar curvature} arXiv:1710.08067.[v2] Tue, 25 Jun 2019 13:11:17 UTC (57 KB)
         \vspace{1mm}   \vspace{1mm}
        
          [Li(rigidity) 2019] C. Li,  {\sl The dihedral rigidity conjecture for n-cubes},
arXiv:1907.03855 [math.DG]\vspace{1mm}   \vspace{1mm}
  	
        [Li(parabolic) 2020] C. Li,   
         {\sl Dihedral rigidity of parabolic polyhedrons in hyperbolic spaces}.
arXiv:2007.12563 [math.DG]
 	          \vspace{1mm}   \vspace{1mm}

[Liang(capillarity) 2005] F-T Liang   {\sl Boundary Regularity for Capillary Surfaces},    
Georgian Mathematical Journal 12(2)(2005)

         \vspace{1mm}   \vspace{1mm}
       
        [Lichnerowitz(spineurs harmoniques) 1963]   A. Lichnerowicz,
{\sl Spineurs harmoniques.}
C. R. Acad. Sci. Paris, S\'erie A, 257 (1963), 7-9.
      
         \vspace{1mm}  \vspace{1mm}

     [Listing(symmetric  spaces) 2010] M. Listing,  {\sl Scalar  curvature  on  compact  symmetric  spaces.} arXiv:1007.1832, 2010. 
        \vspace{1mm}
      
        \vspace{1mm}
       
        [Llarull(sharp estimates)  1998]  M. Llarull {\sl Sharp 
        estimates and the Dirac operator,} Mathematische Annalen January 1998, Volume 310, Issue 1, pp 55-71.
 
 \vspace{1mm}  \vspace{1mm}

 [Lohkamp(negative  Ricci curvature) 1994] J. Lohkamp,
{\sl  Metrics of negative Ricci curvature}, Annals of Mathematics, 140 (1994), 655-683.
  \vspace{1mm}

\vspace{1mm}

  [Lohkamp(hammocks)  1999] \label {111}    J. Lohkamp, Scalar curvature and hammocks, Math. Ann. 313, 385-407, 1999.\vspace{1mm}
  \vspace{1mm}

  [Lohkamp(smoothing) 2018] J. Lohkamp, {\sl Minimal Smoothings of Area Minimizing Cones}, https://arxiv.org/abs/1810.03157\vspace{1mm}  \vspace{1mm}

[Lott(boundary) 2020] J. Lott,
{\sl Index theory for scalar curvature on manifolds with boundary},
arXiv:2009.07256

\vspace{1mm}  \vspace{1mm}

[Lafont-Schimidt(simplicial volume) 2017] J.-F. Lafont, B. Schmidt, 
 {\sl Simplicial volume of closed locally symmetric spaces of non-compact type.}
   	arXiv:math/0504338

\vspace{1mm}  \vspace{1mm}

[Lee-Naber-Neumayer](convergence) 2019]     M. Lee, A. Naber  R. Neumayer.

{\sl $d_p$-Convergence and epsilon regularity theorems for entropy and scalar curvature lower bounds}
arxiv.org/pdf/1808.07841

\vspace{1mm}  \vspace{1mm}

[Lieberman(capillary) 1987] G. Lieberman, a  {\sl Holder continuity of the gradient at a corner 
for the capillary problem 
and related results,} 
Pacific J. Math. 133 (1988), no. 1, 115-135.

\vspace{1mm}  \vspace{1mm}

{Lio-Max (waist inequality) 2020]   Y.  Liokumovich, D. Maximo, {\sl  Waist inequality for 3-manifolds with positive scalar curvature.} This volume.

\vspace{1mm}  \vspace{1mm}

[Lio-Li-Na-Ro(filling) 2019]  Y. Liokumovich, B. Lishak, A. Nabutovsky, R.Rotman       {\sl  Filling metric spaces},
	arXiv:1905.06522 [math.DG].

\vspace{1mm}  \vspace{1mm}

[Lusztig(Novikov) 1972]  G. Lusztig,  {\sl Novikov's higher signature and families of elliptic s}, J.Diff.Geom. 7(1972), 229-256. 
\vspace{1mm}  \vspace{1mm}

[Lusztig(cohomology) 1996]  G. Lusztig,  {\sl Cohomology of Classifying Spaces and Hermitian Representations,} 
Representation Theory, An Electronic Journal of the American Mathematical SocietyVolume 1, Pages 31-36 (November 4, 1996),

\hspace {-6mm}\url {https://pdfs.semanticscholar.org/026c/e6b1c8f754143f6ed72008fb8f044af7d835.pdf}

\vspace{1mm}  \vspace{1mm}

[MarMin(global riemannian) 2012] S. Markvorsen, M. Min-Oo, Global Riemannian Geometry: Curvature
and Topology, 2012 Birkh\"auser.
   
   \url{ https://pdfs.semanticscholar.org/4890/0527441badea97c64130819fb338daa5f864.pdf}

    \vspace{1mm}    \vspace{1mm}

[Marques(deforming $Sc>0$)2012]   F. Marques, {\sl  Deforming three-manifolds with positive scalar curvature},
 arXiv:0907.2444v2.\vspace{1mm}  \vspace{1mm}

[Marques-Neves(min-max spheres in  3d)  2011] F. Marques, A. Neves {\sl Rigidity of min-max minimal spheres in three manifolds},
https://arxiv.org/pdf/1105.4632.pdf
  \vspace{1mm}  \vspace{1mm}

  [Matheu(Dirac) 2012]  A. Matheu, {\sl The Dirac Operator}.  
  \url{math.uchicago.edu/~amathew/dirac.pdf}
    \vspace{1mm}  \vspace{1mm}
  
  [Meinrenken(lectures) 2017]  Eckhard Meinrenken, {\sl 
Lie Groupoids and Lie Algebroids,}   Lecture notes , Fall 2017,  

\url{ http://www.math.toronto.edu/mein/teaching/MAT1341_LieGroupoids/Groupoids2.pdf}
\vspace{1mm}  

\vspace{1mm}

 [Miao(corners) 2002]  P. Miao, {\sl Positive mass theorem on manifolds admitting corners along a hypersurface}, Adv. Theor. Math. Phys., 6 (6): 1163--1182, 2002.

   \vspace{1mm}  

\vspace{1mm}

   [Miao(nonexistence of fill-ins) 2020] P. Miao, {\sl  Nonexistence 
of NNSC fill-ins with large mean curvature}
	arXiv:2009.04976

   \vspace{1mm}  

\vspace{1mm}

   [Min-Oo(hyperbolic)  1989] M. Min-Oo, {\sl Scalar curvature rigidity of asymptotically hyperbolic spin manifolds},
Math. Ann. 285 (1989), 527?539.\vspace{1mm}  \vspace{1mm}

[Min-Oo(Hermitian) 1998] M. Min-Oo, {\sl Scalar Curvature Rigidity of Certain Symmetric Spaces,}  Geometry, Topology and Dynamics
(Montreal, PQ, 1995), CRM Proc. Lecture Notes, 15, Amer. Math. Soc., Providence, RI, 1998, pp. 127-136.\vspace{1mm}  \vspace{1mm}
 
[Min-Oo(K-Area) 2002] M. Min-Oo,  {\sl K-Area mass and asymptotic geometry}
  \url{https://pdfs.semanticscholar.org/4890/0527441badea97c64130819fb338daa5f864.pdf}
\vspace{1mm}  \vspace{1mm}

[Min-Oo(scalar) 2020] {\sl The Lichnerowicz formula 
and lower bounds for the scalar curvature} these proceedings"

 \vspace{1mm}  \vspace{1mm}

[Mishchenko(infinite-dimensional) 1974] A. Mishchenko, {\sl  Infinite-dimensional representations of discrete groups, and higher signatures,}			  	
Izv. Akad. Nauk SSSR Ser. Mat., 38:1 (1974), 81-106; Math. USSR-Izv., 8:1 (1974), 85-111

c\vspace{1mm}\vspace{1mm}

[Morgan(isoperimetric) 2003] F. Morgan, Regularity of isoperimetric hypersurfaces in Riemannian manifolds, Trans. AMS 355 (2003),
5041?5052.
\vspace{1mm}\vspace{1mm}

[MW(mapping classes) 2018]  F. Manin, S. Weinberger {\sl Integral and rational mapping classes}
arXiv preprint arXiv:1802.0578

\vspace{1mm}\vspace{1mm}

  [NaSchSt(localization)   2001] V. Nazaikinskii, B-F. Schulze, B. Sternin, {\sl  The Localization Problem in Index Theory of Elliptic s,} 

\url {https://publishup.uni-potsdam.de/opus4-ubp/frontdoor/deliver/index/docId/2411/file/2001_34.pdf}

\vspace{1mm}\vspace{1mm}

%[NT(gauge-invariant) 2004] V. Nistor, E. Troitsky {\sf  An index for gauge-invariant operators and the Dixmier-Douady invariant.}  Trans. Amer. Math. Soc. 356 (2004), 185-218. \vspace{1mm}\vspace{1mm}

  [NS(Lipschitz) 2007] J. Naumann, C. Simader,   {\sf  Measure and Integration on Lipschitz Manifolds.} 

\url{http://www2.mathematik.hu-berlin.de/publ/pre/2007/P-07-15.pdf}
  \vspace{1mm}\vspace{1mm}
  
 [Nuchi(cube) 2018]  H.  Nuchi, {\sl Doubling a cube with positive scalar curvature},
  	arXiv:1810.01579
  
    \vspace{1mm}\vspace{1mm}

  [Ono(spectrum) 1988]  K.Ono, {\sl The scalar curvature and the spectrum of the Laplacian of spin manifolds}. Math. Ann. 281, 163-168 (1988).\vspace{1mm}  \vspace{1mm}

[Pa-Ke-Pe(graphical tori) 2020]   A. Pacheco,  C.  Ketterer,  R.Perales, {\sl Stability of graphical tori with almost nonnegative scalar curvature}  	arXiv:1902.03458
 \vspace {2mm}

[Papasoglu(width) 2019]   P. Papasoglu {\sl Uryson width and volume,}
	arXiv:1909.03738 [math.DG]

 \vspace {2mm}

 [Penrose(gravitational collapse) 1965] R. Penrose, {\sl Gravitational Collapse 
 and Space-Time Singularities}
Phys. Rev. Lett. 14, 57.\vspace {2mm}

 [Penrose(naked singularities) 1973] R. Penrose, Naked singularities, Ann. New York Acad.Sci.224(1973), 125-34.
 \vspace{1mm}\vspace{1mm}

 [Perelman(width) 1995] G. Perelman, {\sl Widths of non-negatively curved spaces} 
 Geometric And Functional Analysis, 
Vol. 5, No. 2 (1995) ? 1995  
 
  \vspace{1mm}\vspace{1mm}

 [Petrunin(convex) 2003] A.Petrunin, {\sl Polyhedral approximations of Riemannian manifolds}, Proceedings of 9th G\"okova
Geometry-Topology Conference,
pp. 1-15. \vspace{1mm}\vspace{1mm}

 [Petrunin(upper bound) 2008)] A.Petrunin,  {\sl An upper bound for curvature integral} St. Petersburg Math. J.
Tom 20 (2008), pp. 255-265.

  \vspace{1mm}\vspace{1mm}
 
[Philippis-Maggi(capillary) 2015]
{\sl Regularity of free boundaries in anisotropic capillarity 
problems and the validity of Young's law,}
arXiv:1402.0549\vspace {2mm}

 [Polterovich(rigidity) 1996] L. Polterovich, {\sl Gromov's K-area and symplectic rigidity}, Geom. an  Funct. Analysis 6(1996), 726-739.

 \vspace{1mm}\vspace{1mm}

[Richard(2-systoles) 2020]  T. Richard   {\sl On the 2-systole of stretched enough positive scalar 
curvature metrics on $S^2\times S^2$,}   arXiv:2007.02705v2. \vspace{1mm}\vspace{1mm}

[Roe(coarse geometry 1996]) John Roe, {\sl Index Theory, Coarse Geometry, and Topology of Manifolds.}
Regional Conference Series in Mathematics
Number 90:
CBMS Conference on Index Theory, Coarse Geometry,
and Topology of Manifolds
held at the University of Colorado,
August 1995.\vspace{1mm}\vspace{1mm}

[Roe (partial vanishing) 2012] John Roe, {\sl  Positive curvature, partial vanishing theorems, and coarse indices.}
arXiv:1210.6100 [math.KT] \vspace{1mm}\vspace{1mm}

[Rosenberg($C^\ast$-algebras - positive  scalar) 1984] J. Rosenberg,  {\sl $C^\ast$-algebras,    positive   scalar   curvature,    and   the   Novikov   conjecture. }  Inst.   Hautes?Etudes   Sci.   Publ.   Math.58,    197-212 
  \vspace{1mm}  \vspace{1mm}

[Ros(isoperimetric) 2001]  A. Ros {\sl The isoperimetric problem}  \url{http://www.ugr.es/~aros/isoper.pdf}

 [Ros-Souam(capillary) 1997]   A. Ros and R. Souam, {\sl   On stability of capillary surfaces in a ball}, 
 pacific journal of mathematics,
Vol. 178, No. 2, 1997.
    \vspace{1mm}  \vspace{1mm}

  [Salamon(lectures) 1999] Dietmar Salamon,
 {\sl Spin geometry and Seiberg-Witten invariants } \url {https://people.math.ethz.ch/~salamon/PREPRINTS/witsei.pdf}
 \vspace{1mm}  \vspace{1mm}

[Savelyev(jumping) 2012]  Ya. Savelyev, {\sf Gromov K-area and jumping curves in $\mathbb CP^n$},  	arXiv:1006.4383 [math.SG] 2012.

 \vspace{1mm}  \vspace{1mm}

[Schick(counterexample)1998]
T. Schick, A counterexample to the (unstable) Gromov-Lawson-Rosenberg conjecture,
Topology 37 (1998), no. 6, 1165-1168. \vspace{1mm}  \vspace{1mm}

[Schick-Zadeh(multi-partitioned) 2015] T. Schick M. Zadeh, {\sl Large scale index of multi-partitioned manifolds}, arXiv:1308.0742.
\vspace{1mm}  \vspace{1mm}

[Schoen(topics-scalar) 2017] R. Schoen, {\sl Topics in Scalar Curvature}

 \url {http://www.homepages.ucl.ac.uk/~ucahjdl/Schoen_Topics_in_scalar_curvature_2017.pdf}
  \vspace{1mm}  \vspace{1mm}

 [Shi-Tam(positive mass)  2002]    Yuguang Shi and Luen-Fai Tam J. {\sl Positive Mass Theorem and the y Behaviors of Compact Manifolds with Nonnegative Scalar Curvature} J.
  Differential Geom.
    Volume 62, Number 1 (2002), 79-125.
    \vspace{1mm}  \vspace{1mm}

[Simon(capillary) 1980] L. Simon, {\sl Regularity of capillary surfaces over domains with corners.}
Leon Simon
Pacific J. Math. 88(2): 363-377. 
\vspace{1mm}  \vspace{1mm}

[Simons(minimal) 1968] J. Simons, {\sl Minimal varieties in riemannian manifolds}, Ann. of Math., 88 (1968) pp. 62-105.
 
  \vspace{1mm}  \vspace{1mm}

  [Simon-Spruck(capillary) 1976]  L. Simon \& J. Spruck {\sl Existence and regularity of a capillary surface with prescribed contact angle}, Archive for Rational Mechanics and Analysis volume 61, pages19-34(1976)

    \vspace {2mm}

  [Smale(generic regularity) 2003]   N. Smale,
{\sl Generic regularity of homologically area minimizing hyper
surfaces in eight-dimensional mani-
folds}, Comm. Anal. Geom. 1, no. 2 (1993), 217-228.\vspace {1mm}
  \vspace{1mm}

[Sormani(scalar curvature-convergence)  2016] C. Sormani   {\sl Scalar Curvature and Intrinsic Flat Convergence,}

arXiv:1606.08949 [math.MG]
  	\vspace {1mm}
  \vspace{1mm}

[Sormani(conjectures on convergence) 2021] C. Sormani, {\sl Conjectures on Convergence and Scalar Curvature}, arXiv:2103.10093\vspace {2mm}

 [Sormani-Wenger(intrinsic flat) 2011]  C. Sormani and S. Wenger, {Intrinsic Flat Distance between Riemannian Manifolds and other Integral Current Spaces}  Journal of Differential Geometry, Vol 87, 2011, 117-199
 
 \vspace {1mm}
  \vspace{1mm}

 [Souam (Schl\"afli) 2004] Rabah Souam. {\sl The Schl\"afli formula for polyhedra and piecewise smooth hypersurfaces}Differential Geometry and its Applications, Volume 20, Issue 1, 2004, Pages 31-45.\vspace {1mm}
  \vspace{1mm}

 [Stern(harmonic) 2019] D. Stern. Scalar curvature and harmonic maps to S1. arXiv:1908.09754\vspace {1mm}
  \vspace{1mm}

  [Stolz(classification) 1995] S.  Stolz, {\sl Positive Scalar Curvature Metrics - Existence
and Classification Questions,}

    \url {https://link.springer.com/chapter/10.1007/978-3-0348-9078-6_56z}
   \vspace{1mm} 	 \vspace{1mm}

  [Stolz(simply connected) 1992] S. Stolz, {\sl Simply   connected  manifolds   of positive  scalar  curvature},
  Ann.  of  Math.  (2) 136
  (1992),  511-540.

 \vspace{1mm} 	 \vspace{1mm}

   [Stolz(survey) 2001] S.  Stolz   {\sl Manifolds of Positive Scalar Curvature}
    \url {http://users.ictp.it/~pub_off/lectures/lns009/Stolz/Stolz.pdf}
   
  \vspace{1mm} 	 \vspace{1mm}

 [Su(foliations) 2018]  G.Su,  {\sl Lower bounds of Lipschitz constants on foliations.} 
arXiv:1801.06967 [math.DG]
  	
 \vspace{1mm} 	 \vspace{1mm}
 
 [Su-Wang-Zhang(area decreasing foliations) 2021] G. Su, X. Wang, W. Zhang,
{\sl Nonnegative scalar curvature and area decreasing maps on complete foliated manifolds}.arXiv:2104.03472
\vspace {2mm}

[Sun-Dai(bi-invariant)2020] {\sl Gromov Rigidity of Bi-Invariant Metrics on Lie Groups and Homogeneous Spaces}
arXiv:2005.00161{\color {red}(new !!!!!!!!)}  \vspace {2mm}

  [SWWZ(fill-in) 2019] Y.Shi, W. Wang, G. Wei, J.  Zhu, 
  {\sl  On the Fill-in of Nonnegative Scalar Curvature Metrics},
	arXiv:1907.12173 [math.DG]    \vspace{1mm}  \vspace{1mm}

   \vspace{1mm}

 \vspace{1mm}

[SWW(total mean)   2020]   Y. Shi, W. Wang, G. Wei, {\sl Total mean curvature of the boundary and nonnegative scalar curvature fill-ins}	 arXiv:2007.06756v2.

     \vspace{1mm} 	 \vspace{1mm}

[SY(incompressible) 1979] R.Schoen, S.T Yau {\sl  Existence of incompressible minimal surfaces and the topology of three dimensional
manifolds of non-negative scalar curvature,} Ann. of Math. 110 (1979), 127-142.\vspace{1mm}  \vspace{1mm}

    [SY(positive mass) 1979] R.Schoen, S.T Yau,
     {\sl On the proof of the positive mass conjecture in general relativity,} Commun. Math. Phys. 65,  (1979). 45-76.
    
    \vspace{1mm}  \vspace{1mm}

   [SY(structure) 1979]   R. Schoen and S. T. Yau, On the structure of manifolds with positive scalar
curvature, Manuscripta Math. 28 (1979), 159-183. \vspace{1mm}  \vspace{1mm}
  
[SY(singularities) 2017] R. Schoen and S. T. Yau {\sl Positive Scalar Curvature and Minimal Hypersurface Singularities.}
arXiv:1704.05490
  \vspace{1mm} 	 \vspace{1mm}

 [Vafa-Witten(fermions) 1984] C. Vafa, E. Witten, {\sl  Eigenvalue inequalities for fermions in gauge theories},
 Comm. Math. Phys. 95(3): 257-276 (1984).

    \vspace{1mm}  \vspace{1mm}

 [Wang(Contractible)  2019] J. Wang, {\sl Contractible 3-manifold and Positive scalar curvature (I),(2)},

arXiv:1901.04605,  	arXiv:1906.04128.
  	    \vspace{1mm} 	 \vspace{1mm}
   
    [Wang(topological characterization)  2021] J. Wang, {\sl Topological Characterization of Contractible 3-manifold and Positive Scalar Curvature)}, These Proceedings.
    \vspace{1mm} 	 \vspace{1mm}

   [Wang-Xie-Yu(decay) 2021] J. Wang, Z.  Xie, G. Yu {\sl Decay of scalar curvature on uniformly contractible manifolds with finite asymptotic dimension,}
Jinmin Wang, Zhizhang Xie, Guoliang Yu
 arXiv:2101.11584 [math.DG] 
     \vspace{1mm} 	 \vspace{1mm}
  	 
     [Wang-Xie-Yu(cube inequality) 2021] J. Wang, Z.  Xie, G. Yu {\sl An index theoretic proof of Gromov's cube inequality on scalar curvature}
	arXiv:2105.12054 [math.DG]

     \vspace{1mm} 	 \vspace{1mm}

   [Weinstein(Positively curved) 1970] A. Weinstein, {\sl Positively curved n
-manifolds in $R^{n+2}$}.
   J. Differential Geom.Volume 4, Number 1 (1970), 1-4.
    
         \vspace{1mm} 	 \vspace{1mm}

 [Wenger(filling) 2007]  S. Wenger
 {\sl A short proof of Gromov's filling inequality,} arXiv:math/070388.

 \vspace{1mm} 	 \vspace{1mm}

 [White(minimal) 2016] B. White,  {\sl Lectures on Minimal Surface Theory}, arXiv:1308.3325

 \vspace{1mm} 	 \vspace{1mm}

     [Witten(Positive Energy) 1981] E. Witten, A New Proof  of the Positive Energy Theorem. Communications in   Math.  Phys.  80, 381-402  (1981)                          
  \vspace{1mm} 	 \vspace{1mm}

[Zeidler(bands)  2019]  R. Zeidler,
{\sl Band width estimates via the Dirac operator}, 
arXiv:1905.08520v2  \vspace{1mm} 	 \vspace{1mm}

[Zeidler(width)  2020]  R. Zeidler, {\sl Width, largeness and index theory.}
arXiv:2008.13754  \vspace{1mm} 	 \vspace{1mm}

[Zhang(foliations) 2016]   W. Zhang {\sl Positive scalar curvature on foliations.}
arXiv:1508.04503 [math.DG]  \vspace{1mm} 	 \vspace{1mm}

[Zhang(foliations:enlargeability)  2018]   W. Zhang {\sl Positive scalar curvature on foliations: the enlargeability} arXiv:1703.04313v2\vspace{2mm}

[Zhang(area decreasing)  2020] W. Zhang  {\sl  Nonnegative Scalar Curvature and Area Decreasing Maps},
arXiv:1912.03649v3.

\vspace{2mm}

[Zhang(deformed Dirac) 2021]    W.Zhang  {\sl  Deformed Dirac Operators and the Scalar curvature }, These Proceedings.\vspace{2mm}

[Zhu(rigidity)  2019]  J. Zhu, {\sl  Rigidity of Area-Minimizing 2-Spheres in n-Manifolds with Positive Scalar Curvature}   	arXiv:1903.0578c\vspace{2mm}

[Zhu(rigidity)  2020]  J. Zhu, {\sl  Rigidity results for complete manifolds with nonnegative scalar curvature}

	arXiv:2008.07028 [math.DG].

\end{document}